\documentclass[12pt]{amsbook}
\usepackage{amssymb}
\usepackage[all]{xy}

\usepackage[colorlinks=true,linkcolor=magenta,citecolor=magenta]{hyperref}
\makeindex

\newtheorem{theorem}{Theorem}[chapter]
\newtheorem{addendum}[theorem]{Addendum}

\newtheorem{cat}[theorem]{Cat}
\newtheorem{comment}[theorem]{Comment}
\newtheorem{conclusion}[theorem]{Conclusion}
\newtheorem{conjecture}[theorem]{Conjecture}
\newtheorem{conventions}[theorem]{Conventions}
\newtheorem{definition}[theorem]{Definition}

\newtheorem{exercise}[theorem]{Exercise}
\newtheorem{fact}[theorem]{Fact}
\newtheorem{plan}[theorem]{Plan}
\newtheorem{principle}[theorem]{Principle}
\newtheorem{problem}[theorem]{Problem}
\newtheorem{proposition}[theorem]{Proposition}
\newtheorem{question}[theorem]{Question}
\newtheorem{remark}[theorem]{Remark}
\newtheorem{speculation}[theorem]{Speculation}

\newtheorem{warning}[theorem]{Warning}

\begin{document}

\title{Methods of free probability}

\author{Teo Banica}
\address{Department of Mathematics, University of Cergy-Pontoise, F-95000 Cergy-Pontoise, France. {\tt teo.banica@gmail.com}}


\begin{abstract}
This is an introduction to free probability, with all the needed preliminaries included. We first discuss the most basic noncommuting variables, namely the usual matrices, followed by the basic random matrices, and with a look into Lie groups too. Then we present, following Voiculescu, the main concepts and tools of free probability, notably with the various limiting theorems that can be obtained. We then go on a more technical discussion, featuring cumulants, functional transforms and group theory, on the correspondence between classical and free. Finally, we discuss matrix modeling questions for the various more specialized free laws and objects that can be imagined.
\end{abstract}

\maketitle

\chapter*{Preface}

Probability theory, and the probabilistic way of thinking, have seen a considerable surge in the last years, with virtually every single branch of mathematics being affected. It goes without saying that everything mathematics coming from quantum mechanics, which actually accounts for a big part of pure mathematics as we know it, has some probability behind, and this has become more and more visible during recent years. The same goes of course for statistical mechanics, once again somehow by definition. As in what regards classical mechanics, randomness of the initial data is certainly a very fruitful idea too. Finally, old branches of pure mathematics, such as number theory, are increasingly becoming more analytic, and more probabilistic too.

\bigskip

At the technical level, probability theory comes in many flavors. However, if there is one thing to be known, having interesting mathematics and physics behind, this is the fact that classical probability theory has a ``twin sister'', namely free probability.

\bigskip

Free probability was introduced by Voiculescu in the mid 1980s, with motivation coming from general quantum mechanics, and more specifically with a number of operator algebra questions in mind. Among the main discoveries of Voiculescu was the fact that Wigner's semicircle law, coming from advanced quantum physics and random matrices, appears as the ``free analogue'' of the normal law. This has led to a lot of interest in free probability, with the subject having now deep ties to operator algebras, random matrices, quantum groups, noncommutative geometry, and virtually any other branch of mathematics coming from quantum mechanics, or statistical mechanics.

\bigskip

This book is an introduction to free probability, with the aim of keeping things as simple and concrete as possible, while still being relatively complete. Our goals will be on one hand that of explaining the definition and main properties of free probability, in analogy with the definition and main properties of classical probability, and by keeping the presentation as elementary as possible, and on the other hand to go, at least a little bit, into each of the above-mentioned classes of examples and applications, namely operator algebras, random matrices, quantum groups and noncommutative geometry.

\bigskip

As mentioned, our presentation is meant to be elementary, based on calculus and combinatorics. However, with our free random variables being of matrix type, some familiarity with linear algebra and its ramifications can be of help. As for the classical probability theory needed, in order to get started, that basically reduces to the knowledge of the normal and Poisson laws, and we will review this in the opening chapter.

\bigskip

The book is organized in 4 parts, as follows:

\bigskip

I. We first discuss the most basic noncommuting variables, namely the usual matrices, followed by the basic random matrices, and with a look into Lie groups too.

\bigskip

II. Then we present, following Voiculescu, the main concepts and tools of free probability, notably with the various limiting theorems that can be obtained.

\bigskip

III. We then go on a more technical discussion, featuring cumulants, functional transforms and group theory, on the correspondence between classical and free. 

\bigskip

IV. Finally, we discuss matrix modeling questions for the various more specialized free laws and objects that can be imagined.

\bigskip

All in all, many things to be discussed. As a complement to what we will be doing here, for advanced combinatorics and operator algebra aspects you have \cite{hpe}, \cite{nsp}, \cite{vdn}, and for advanced random matrix theory you have \cite{agz}, \cite{bos}, \cite{msp}. So, in the hope that you will like free probability, and end up learning everything, from here and from \cite{agz}, \cite{bos}, \cite{hpe}, \cite{msp}, \cite{nsp}, \cite{vdn}, with the precise order being more a matter of taste.

\bigskip

I learned myself free probability long ago, as a graduate student, from \cite{vdn}, with my first research paper being a 1996 note on the circular variables \cite{ba1}. Later I started doing quantum groups, and some random matrices too, with free probability always in mind. I am grateful to Mireille Capitaine, Beno\^it Collins, Steve Curran, Ion Nechita, Roland Speicher and the others, for substantial joint work on the subject. Many thanks go as well to my cats. No serious science can be done without advice from a cat or two.

\bigskip

\

{\em Cergy, September 2026}

\smallskip

{\em Teo Banica}

\baselineskip=15.95pt
\tableofcontents
\baselineskip=14pt

\part{Matrices, laws}

\ \vskip50mm

\begin{center}
{\em The Magical Mystery Tour

Is coming to take you away

Coming to take you away

Take you today}
\end{center}

\chapter{Normal variables}

\section*{1a. Random variables}

Generally speaking, probability theory is best learned by flipping coins and rolling dice. At a more advanced level, which is playing cards, we have:

\index{poker}

\begin{theorem}
The probabilities at poker are as follows:
\begin{enumerate}
\item One pair: $0.533$.

\item Two pairs: $0.120$.

\item Three of a kind: $0.053$.

\item Straight: $0.020$.

\item Full house: $0.006$.

\item Four of a kind: $0.001$.

\item Flush: $0.001$.

\item Straight flush: $0.000$.
\end{enumerate}
\end{theorem}

\begin{proof}
Let us consider indeed our deck of 32 cards, $7,8,9,10,J,Q,K,A$. The total number of possibilities for a poker hand is:
$$\binom{32}{5}
=\frac{32\cdot 31\cdot 30\cdot 29\cdot 28}{2\cdot 3\cdot 4\cdot 5}
=32\cdot 31\cdot 29\cdot 7$$

(1) For having a pair, the number of possibilities is:
$$N
=\binom{8}{1}\binom{4}{2}\times\binom{7}{3}\binom{4}{1}^3
=8\cdot 6\cdot 35\cdot 64$$

Thus, the probability of having a pair is:
$$P
=\frac{8\cdot 6\cdot 35\cdot 64}{32\cdot 31\cdot 29\cdot 7}
=\frac{6\cdot 5\cdot 16}{31\cdot 29}
=\frac{480}{899}
=0.533$$

(2) For having two pairs, the number of possibilities is:
$$N
=\binom{8}{2}\binom{4}{2}^2\times\binom{24}{1}
=28\cdot 36\cdot 24$$

Thus, the probability of having two pairs is:
$$P
=\frac{28\cdot 36\cdot 24}{32\cdot 31\cdot 29\cdot 7}
=\frac{36\cdot 3}{31\cdot 29}
=\frac{108}{899}
=0.120$$

(3) For having three of a kind, the number of possibilities is:
$$N
=\binom{8}{1}\binom{4}{3}\times\binom{7}{2}\binom{4}{1}^2
=8\cdot 4\cdot 21\cdot 16$$

Thus, the probability of having three of a kind is:
$$P
=\frac{8\cdot 4\cdot 21\cdot 16}{32\cdot 31\cdot 29\cdot 7}
=\frac{3\cdot 16}{31\cdot 29}
=\frac{48}{899}
=0.053$$

(4) For having a straight, the number of possibilities is:
$$N
=4\left[\binom{4}{1}^5-4\right]
=16\cdot 255$$

Thus, the probability of having a straight is:
$$P
=\frac{16\cdot 255}{32\cdot 31\cdot 29\cdot 7}
=\frac{255}{2\cdot 31\cdot 29\cdot 7}
=\frac{255}{12586}
=0.020$$

(5) For having full house, the number of possibilities is:
$$N
=\binom{8}{1}\binom{4}{3}\times\binom{7}{1}\binom{4}{2}
=8\cdot 4\cdot 7\cdot 6$$

Thus, the probability of having full house is:
$$P
=\frac{8\cdot 4\cdot 7\cdot 6}{32\cdot 31\cdot 29\cdot 7}
=\frac{6}{31\cdot 29}
=\frac{6}{899}
=0.006$$

(6) For having four of a kind, the number of possibilities is:
$$N
=\binom{8}{1}\binom{4}{4}\times\binom{7}{1}\binom{4}{1}
=8\cdot 7\cdot 4$$

Thus, the probability of having four of a kind is:
$$P
=\frac{8\cdot 7\cdot 4}{32\cdot 31\cdot 29\cdot 7}
=\frac{1}{31\cdot 29}
=\frac{1}{899}
=0.001$$

(7) For having a flush, the number of possibilities is:
$$N
=4\left[\binom{8}{5}-4\right]
=4\cdot 52$$

Thus, the probability of having a flush is:
$$P
=\frac{4\cdot 52}{32\cdot 31\cdot 29\cdot 7}
=\frac{13}{2\cdot 31\cdot 29\cdot 7}
=\frac{13}{12586}
=0.001$$

(8) For having a straight flush, the number of possibilities is:
$$N=4\cdot 4$$

Thus, the probability of having a straight flush is:
$$P
=\frac{4\cdot 4}{32\cdot 31\cdot 29\cdot 7}
=\frac{1}{2\cdot 31\cdot 29\cdot 7}
=\frac{1}{12586}
=0.000$$

Thus, we have obtained the numbers in the statement.
\end{proof}

Summarizing, probability is about binomials and factorials. Getting now to more formal work, with the idea in mind of doing things a bit abstractly, we have:

\index{probability space}
\index{random variable}
\index{moments}
\index{law}
\index{distribution}

\begin{definition}
Let $X$ be a probability space, that is, a space with a probability measure, and with the corresponding integration denoted $E$, and called expectation.
\begin{enumerate}
\item The random variables are the real functions $f\in L^\infty(X)$.

\item The moments of such a variable are the numbers $M_k(f)=E(f^k)$.

\item The law of such a variable is the measure given by $M_k(f)=\int_\mathbb Rx^kd\mu_f(x)$.
\end{enumerate}
\end{definition}

Here the fact that $\mu_f$ exists indeed is well-known. By linearity, we would like to have a real probability measure making hold the following formula, for any $P\in\mathbb R[X]$:
$$E(P(f))=\int_\mathbb RP(x)d\mu_f(x)$$

By using a standard continuity argument, it is enough to have this formula for the characteristic functions $\chi_I$ of the measurable sets of real numbers $I\subset\mathbb R$:
$$E(\chi_I(f))=\int_\mathbb R\chi_I(x)d\mu_f(x)$$

But this latter formula, which reads $P(f\in I)=\mu_f(I)$, can serve as a definition for $\mu_f$, and we are done. Next, we can talk as well about independence, as follows:

\index{independence}

\begin{definition}
Two variables $f,g\in L^\infty(X)$ are called independent when
$$E(f^kg^l)=E(f^k)\,E(g^l)$$
happens, for any $k,l\in\mathbb N$.
\end{definition}

Again, this definition hides some non-trivial things. Indeed, by linearity, we would like to have a formula as follows, valid for any polynomials $P,Q\in\mathbb R[X]$:
$$E[P(f)Q(g)]=E[P(f)]\,E[Q(g)]$$

By using a continuity argument, it is enough to have this formula for characteristic functions $\chi_I,\chi_J$ of the measurable sets of real numbers $I,J\subset\mathbb R$:
$$E[\chi_I(f)\chi_J(g)]=E[\chi_I(f)]\,E[\chi_J(g)]$$

Thus, we are led to the usual definition of independence, namely:
$$P(f\in I,g\in J)=P(f\in I)\,P(g\in J)$$

All this might seem a bit abstract, but in practice, the idea is of course that $f,g$ must be independent, in an intuitive, real-life sense. As a first result now, we have:

\begin{proposition}
Assuming that $f,g\in L^\infty(X)$ are independent, we have
$$\mu_{f+g}=\mu_f*\mu_g$$
where $*$ is the convolution of real probability measures.
\end{proposition}

\begin{proof}
We have the following computation, using the independence of $f,g$:
\begin{eqnarray*}
M_k(f+g)
&=&E((f+g)^k)\\
&=&\sum_r\binom{k}{r}E(f^rg^{k-r})\\
&=&\sum_r\binom{k}{r}M_r(f)M_{k-r}(g)
\end{eqnarray*}

On the other hand, by using the Fubini theorem, we have as well:
\begin{eqnarray*}
\int_\mathbb Rx^kd(\mu_f*\mu_g)(x)
&=&\int_{\mathbb R\times\mathbb R}(x+y)^kd\mu_f(x)d\mu_g(y)\\
&=&\sum_r\binom{k}{r}\int_\mathbb Rx^rd\mu_f(x)\int_\mathbb Ry^{k-r}d\mu_g(y)\\
&=&\sum_r\binom{k}{r}M_r(f)M_{k-r}(g)
\end{eqnarray*}

Thus $\mu_{f+g}$ and $\mu_f*\mu_g$ have the same moments, so they coincide, as desired.
\end{proof}

Here is now a second result on independence, which is something more advanced:

\index{independence}
\index{Fourier transform}

\begin{theorem}
Assuming that $f,g\in L^\infty(X)$ are independent, we have
$$F_{f+g}=F_fF_g$$
where $F_f(x)=E(e^{ixf})$ is the Fourier transform.
\end{theorem}

\begin{proof}
We have the following computation, using Proposition 1.4 and Fubini:
\begin{eqnarray*}
F_{f+g}(x)
&=&\int_\mathbb Re^{ixz}d\mu_{f+g}(z)\\
&=&\int_\mathbb Re^{ixz}d(\mu_f*\mu_g)(z)\\
&=&\int_{\mathbb R\times\mathbb R}e^{ix(z+t)}d\mu_f(z)d\mu_g(t)\\
&=&\int_\mathbb Re^{ixz}d\mu_f(z)\int_\mathbb Re^{ixt}d\mu_g(t)\\
&=&F_f(x)F_g(x)
\end{eqnarray*}

Thus, we are led to the conclusion in the statement.
\end{proof}

\section*{1b. Normal laws}

The main result in classical probability is the Central Limit Theorem (CLT), that we would like to explain now. We will need the following standard fact:

\index{polar coordinates}
\index{Jacobian}
\index{Gauss integral}

\begin{theorem}
We have the following formula,
$$\int_\mathbb Re^{-x^2}dx=\sqrt{\pi}$$
called Gauss integral formula.
\end{theorem}

\begin{proof}
Let $I$ be the integral in the statement. By using polar coordinates, namely $x=r\cos t$, $y=r\sin t$, with the corresponding Jacobian being $r$, we have:
$$I^2
=\int_\mathbb R\int_\mathbb Re^{-x^2-y^2}dxdy
=\int_0^{2\pi}\int_0^\infty e^{-r^2}r\,drdt
=2\pi\int_0^\infty\left(-\frac{e^{-r^2}}{2}\right)'dr
=\pi$$

Thus, we are led to the formula in the statement.
\end{proof}

We can now introduce the normal distributions, as follows:

\index{normal law}
\index{Gaussian law}

\begin{definition}
The normal law of parameter $t>0$ is the following measure:
$$g_t=\frac{1}{\sqrt{2\pi t}}\,e^{-x^2/2t}dx$$
This is also called Gaussian distribution, with ``g'' standing for Gauss.
\end{definition}

Observe that the above law has indeed mass 1, as it should, due to the Gauss formula. Many things can be said about the normal laws, with a useful result being:

\index{Fourier transform}
\index{convolution semigroup}

\begin{proposition}
We have the following formula, valid for any $t>0$:
$$F_{g_t}(x)=e^{-tx^2/2}$$
In particular, the normal laws satisfy $g_s*g_t=g_{s+t}$, for any $s,t>0$.
\end{proposition}

\begin{proof}
The Fourier transform formula can be established as follows:
\begin{eqnarray*}
F_{g_t}(x)
&=&\frac{1}{\sqrt{2\pi t}}\int_\mathbb Re^{-z^2/2t+ixz}dz\\
&=&\frac{1}{\sqrt{2\pi t}}\int_\mathbb Re^{-(z/\sqrt{2t}-\sqrt{t/2}\,iz)^2-tx^2/2}dz\\
&=&\frac{1}{\sqrt{2\pi t}}\int_\mathbb Re^{-y^2-tx^2/2}\sqrt{2t}\,dy\\
&=&\frac{1}{\sqrt{\pi}}\,e^{-tx^2/2}\int_\mathbb Re^{-y^2}dy\\
&=&e^{-tx^2/2}
\end{eqnarray*}

As for $g_s*g_t=g_{s+t}$, this comes from Theorem 1.5, $\log F_{g_t}$ being linear in $t$.
\end{proof}

We can now state and prove the Central Limit Theorem, as follows:

\index{CLT}
\index{Central Limit Theorem}

\begin{theorem}[CLT]
Given real variables $f_1,f_2,f_3,\ldots\in L^\infty(X)$ which are independent identically distributed (i.i.d.), centered, with common variance $t>0$, we have
$$\frac{1}{\sqrt{n}}\sum_{i=1}^nf_i\sim g_t$$
with $n\to\infty$, in moments.
\end{theorem}

\begin{proof}
In terms of moments, the Fourier transform $F_f(x)=E(e^{ixf})$ is given by:
$$F_f(x)
=E\left(\sum_{r=0}^\infty\frac{(ixf)^r}{r!}\right)
=\sum_{r=0}^\infty\frac{(ix)^rE(f^r)}{r!}
=\sum_{r=0}^\infty\frac{i^rM_r(f)}{r!}\,x^r$$

Thus, the Fourier transform of the variable in the statement is:
$$F(x)
=\left[F_f\left(\frac{x}{\sqrt{n}}\right)\right]^n
=\left[1-\frac{tx^2}{2n}+O(n^{-2})\right]^n
\simeq e^{-tx^2/2}$$

But this function being the Fourier transform of $g_t$, we obtain the result.
\end{proof}

Let us discuss now some further properties of the normal law. We first have:

\index{moments}

\begin{theorem}
The even moments of the normal law are the numbers
$$M_k(g_t)=t^{k/2}\times k!!$$
where $k!!=(k-1)(k-3)(k-5)\ldots\,$, and the odd moments vanish. 
\end{theorem}

\begin{proof}
We have the following computation, valid for any integer $k\in\mathbb N$:
\begin{eqnarray*}
M_k
&=&\frac{1}{\sqrt{2\pi t}}\int_\mathbb Ry^ke^{-x^2/2t}dx\\
&=&\frac{1}{\sqrt{2\pi t}}\int_\mathbb R(tx^{k-1})\left(-e^{-x^2/2t}\right)'dx\\
&=&\frac{1}{\sqrt{2\pi t}}\int_\mathbb Rt(k-1)x^{k-2}e^{-x^2/2t}dx\\
&=&t(k-1)\times\frac{1}{\sqrt{2\pi t}}\int_\mathbb Rx^{k-2}e^{-x^2/2t}dx\\
&=&t(k-1)M_{k-2}
\end{eqnarray*}

But $M_0=1$, $M_1=0$, so we are led to the formula in the statement.
\end{proof}

Next, we have the following reformulation of the above result, whose raison d'etre we will understand later, when investigating some other probability measures:

\begin{theorem}
The moments of the normal law are the numbers
$$M_k(g_t)=\sum_{\pi\in P_2(k)}t^{|\pi|}$$
where $P_2(k)$ is the set of pairings of $\{1,\ldots,k\}$, and $|.|$ is the number of blocks.
\end{theorem}

\begin{proof}
Since the number of blocks of a pairing of $\{1,\ldots,k\}$ is trivially $k/2$, independently of the pairing, our formula says in practice that we have:
$$M_k(g_t)=t^{k/2}|P_2(k)|$$

In order to prove this, let us count the pairings of $\{1,\ldots,k\}$. In order to have such a pairing, we must pair $1$ with one of the numbers $2,\ldots,k$, and then use a pairing of the remaining $k-2$ numbers. Thus, we have the following recurrence formula:
$$|P_2(k)|=(k-1)|P_2(k-2)|$$

With the initial data being $|P_2(1)|=0$, $|P_2(2)|=1$, we are led to the result.
\end{proof}

In a purely mathematical context, the normal laws appear by looking at the coordinates over the unit spheres $S^{N-1}_\mathbb R\subset\mathbb R^N$, in the $N\to\infty$ limit. We will need:

\index{Wallis formula}
\index{trigonometric integral}

\begin{theorem}[Wallis]
We have the following formula,
$$\int_0^{\pi/2}\cos^rt\sin^st\,dt=\left(\frac{\pi}{2}\right)^{\varepsilon(r)\varepsilon(s)}\frac{r!!s!!}{(r+s+1)!!}$$
where $\varepsilon(r)=1$ if $r$ is even, and $\varepsilon(r)=0$ if $r$ is odd.
\end{theorem}

\begin{proof}
Let us call $I_{rs}$ the integral in the statement. In order to first compute $I_{r0}$, we can use partial integration. We have the following formula:
\begin{eqnarray*}
(\cos^rt\sin t)'
&=&r\cos^{r-1}t(-\sin t)\sin t+\cos^rt\cos t\\
&=&(r+1)\cos^{r+1}t-k\cos^{r-1}t
\end{eqnarray*}

By integrating between $0$ and $\pi/2$, we obtain the following formula:
$$(r+1)I_{r+1,0}=rI_{r-1,0}$$

Thus we can compute $I_{r0}$ by recurrence, and we obtain in this way:
$$I_{r0}
=\frac{r-1}{r}\,I_{r-2,0}
=\frac{r-1}{r}\cdot\frac{r-3}{r-2}\,I_{r-4,0}
=\ldots\ldots
=\frac{r!!}{(r+1)!!}\,I_{1-\varepsilon(r),0}$$

The initial data being $I_{00}=\pi/2$ and $I_{10}=1$, we obtain the result. Next, observe that the integrals $I_{r1}$ can be computed as well, this time trivially, as follows:
$$I_{r1}
=\int_0^{\pi/2}\cos^rt\sin t\,dt
=-\frac{1}{r+1}\int_0^{\pi/2}(\cos^{r+1}t)'\,dt
=\frac{1}{r+1}$$

Getting now to the general case, in order to do the partial integration, we have:
\begin{eqnarray*}
(\cos^rt\sin^st)'
&=&r\cos^{r-1}t(-\sin t)\sin^st+\cos^rt\cdot s\sin^{s-1}t\cos t\\
&=&-r\cos^{r-1}t\sin^{s+1}t+s\cos^{r+1}t\sin^{s-1}t
\end{eqnarray*}

By integrating between $0$ and $\pi/2$, we obtain, for $r,s>0$:
$$rI_{r-1,s+1}=sI_{r+1,s-1}$$

Thus, we can compute $I_{rs}$ by recurrence. When $s$ is even we obtain:
$$I_{rs}
=\frac{s-1}{r+1}\,I_{r+2,s-2}
=\ldots\ldots
=\frac{r!!s!!}{(r+s)!!}\,I_{r+s,0}
=\left(\frac{\pi}{2}\right)^{\varepsilon(r)}\frac{r!!s!!}{(r+s+1)!!}$$

As for the case where $s$ is odd, here the recurrence goes as follows:
$$I_{rs}
=\frac{s-1}{r+1}\,I_{r+2,s-2}
=\ldots\ldots
=\frac{r!!s!!}{(r+s-1)!!}\,I_{r+s-1,1}
=\frac{r!!s!!}{(r+s+1)!!}$$

Thus, we are led to the formula in the statement.
\end{proof}

More generally now, over the spheres in arbitrary dimensions, we have:

\index{spherical integral}

\begin{theorem}
We have the following integration formula, for $k_i$ even,
$$\int_{S^{N-1}_\mathbb R}x_1^{k_1}\ldots x_N^{k_N}\,dx=\frac{(N-1)!!k_1!!\ldots k_N!!}{(N+\sum k_i-1)!!}$$
with respect to the mass $1$ uniform measure on the unit sphere $S^{N-1}_\mathbb R\subset\mathbb R^N$.
\end{theorem}

\begin{proof}
We can use spherical coordinates in $N$ dimensions, namely:
$$\begin{cases}
x_1\!\!\!&=\ r\cos t_1\\
x_2\!\!\!&=\ r\sin t_1\cos t_2\\
\vdots\\
x_{N-1}\!\!\!&=\ r\sin t_1\sin t_2\ldots\sin t_{N-2}\cos t_{N-1}\\
x_N\!\!\!&=\ r\sin t_1\sin t_2\ldots\sin t_{N-2}\sin t_{N-1}
\end{cases}$$

Regarding the Jacobian, by developing over the last column, we have:
\begin{eqnarray*}
J_N
&=&r\sin t_1\ldots\sin t_{N-2}\sin t_{N-1}\times \sin t_{N-1}J_{N-1}\\
&+&r\sin t_1\ldots \sin t_{N-2}\cos t_{N-1}\times\cos t_{N-1}J_{N-1}\\
&=&r\sin t_1\ldots\sin t_{N-2}(\sin^2 t_{N-1}+\cos^2 t_{N-1})J_{N-1}\\
&=&r\sin t_1\ldots\sin t_{N-2}J_{N-1}
\end{eqnarray*}

Thus, by recurrence, the Jacobian is given by the following formula:
$$J_N=r^{N-1}\sin^{N-2}t_1\sin^{N-3}t_2\,\ldots\,\sin^2t_{N-3}\sin t_{N-2}$$

As a first application, we can compute the volume of the sphere, as follows:
\begin{eqnarray*}
\frac{V}{2^N}
&=&\int_0^1\int_0^{\pi/2}\ldots\int_0^{\pi/2}r^{N-1}\sin^{N-2}t_1\ldots\sin t_{N-2}\,drdt_1\ldots dt_{N-1}\\
&=&\int_0^1r^{N-1}\,dr\int_0^{\pi/2}\sin^{N-2}t_1\,dt_1\ldots\int_0^{\pi/2}\sin t_{N-2}dt_{N-2}\int_0^{\pi/2}1dt_{N-1}\\
&=&\frac{1}{N}\times\left(\frac{\pi}{2}\right)^{[N/2]}\times\frac{(N-2)!!}{(N-1)!!}\cdot\frac{(N-3)!!}{(N-2)!!}\ldots\frac{2!!}{3!!}\cdot\frac{1!!}{2!!}\cdot1\\
&=&\left(\frac{\pi}{2}\right)^{[N/2]}\frac{1}{(N+1)!!}
\end{eqnarray*}

Next, let us write the integral in the statement in spherical coordinates:
$$I=\frac{2^N}{A}\int_0^{\pi/2}\ldots\int_0^{\pi/2}x_1^{k_1}\ldots x_N^{k_N}J_N\,dt_1\ldots dt_{N-1}$$

Here $A=NV$ is the area of the sphere, and the unnormalized integral is:
\begin{eqnarray*}
I'=\int_0^{\pi/2}\ldots\int_0^{\pi/2}
&&(\cos t_1)^{k_1}
(\sin t_1\cos t_2)^{k_2}\\
&&\vdots\\
&&(\sin t_1\sin t_2\ldots\sin t_{N-2}\cos t_{N-1})^{k_{N-1}}\\
&&(\sin t_1\sin t_2\ldots\sin t_{N-2}\sin t_{N-1})^{k_N}\\
&&\sin^{N-2}t_1\sin^{N-3}t_2\ldots\sin^2t_{N-3}\sin t_{N-2}\\
&&dt_1\ldots dt_{N-1}
\end{eqnarray*}

By rearranging the terms, we obtain the following formula:
\begin{eqnarray*}
I'
&=&\int_0^{\pi/2}\cos^{k_1}t_1\sin^{k_2+\ldots+k_N+N-2}t_1\,dt_1\\
&&\int_0^{\pi/2}\cos^{k_2}t_2\sin^{k_3+\ldots+k_N+N-3}t_2\,dt_2\\
&&\vdots\\
&&\int_0^{\pi/2}\cos^{k_{N-2}}t_{N-2}\sin^{k_{N-1}+k_N+1}t_{N-2}\,dt_{N-2}\\
&&\int_0^{\pi/2}\cos^{k_{N-1}}t_{N-1}\sin^{k_N}t_{N-1}\,dt_{N-1}
\end{eqnarray*}

Now by using the formula in Theorem 1.12, this gives:
\begin{eqnarray*}
I'
&=&\left(\frac{\pi}{2}\right)^{\varepsilon(N-2)}\frac{k_1!!(k_2+\ldots+k_N+N-2)!!}{(k_1+\ldots+k_N+N-1)!!}\\
&&\left(\frac{\pi}{2}\right)^{\varepsilon(N-3)}\frac{k_2!!(k_3+\ldots+k_N+N-3)!!}{(k_2+\ldots+k_N+N-2)!!}\\
&&\vdots\\
&&\left(\frac{\pi}{2}\right)^{\varepsilon(1)}\frac{k_{N-2}!!(k_{N-1}+k_N+1)!!}{(k_{N-2}+k_{N-1}+k_N+2)!!}\\
&&\left(\frac{\pi}{2}\right)^{\varepsilon(0)}\frac{k_{N-1}!!k_N!!}{(k_{N-1}+k_N+1)!!}
\end{eqnarray*}

By simplifying this, and putting everything together, we obtain:
$$I=\left(\frac{2}{\pi}\right)^{[N/2]}(N-1)!!\cdot\left(\frac{\pi}{2}\right)^{[N/2]}
\frac{k_1!!\ldots k_N!!}{(N+\sum k_i-1)!!}
=\frac{(N-1)!!k_1!!\ldots k_N!!}{(N+\sum k_i-1)!!}$$

Thus, we have obtained the formula in the statement.
\end{proof}

We can now recover the normal laws, geometrically, as follows:

\index{hyperspherical law}
\index{asymptotic independence}
\index{asymptotic law}

\begin{theorem}
The even moments of the hyperspherical variables are
$$\int_{S^{N-1}_\mathbb R}x_i^kdx=\frac{(N-1)!!k!!}{(N+k-1)!!}$$
and the rescaled variables $y_i=\sqrt{N}x_i$ become normal and independent with $N\to\infty$.
\end{theorem}

\begin{proof}
The moment formula in the statement follows from Theorem 1.13. As a consequence, with $N\to\infty$ we have the following estimate:
$$\int_{S^{N-1}_\mathbb R}x_i^kdx
\simeq N^{-k/2}\times k!!
=N^{-k/2}M_k(g_1)$$

Thus, the rescaled variables $\sqrt{N}x_i$ become normal with $N\to\infty$, as claimed. As for the proof of the asymptotic independence, this is standard too, once again by using the formula in Theorem 1.13. Indeed, the joint moments of $x_1,\ldots,x_N$ are given by:
$$\int_{S^{N-1}_\mathbb R}x_1^{k_1}\ldots x_N^{k_N}\,dx
=\frac{(N-1)!!k_1!!\ldots k_N!!}{(N+\sum k_i-1)!!}
\simeq N^{-\sum k_i/2}\times k_1!!\ldots k_N!!$$

By rescaling, the joint moments of the variables $y_i=\sqrt{N}x_i$ are given by:
$$\int_{S^{N-1}_\mathbb R}y_1^{k_1}\ldots y_N^{k_N}\,dx\simeq k_1!!\ldots k_N!!$$

Thus, we have multiplicativity, and so independence with $N\to\infty$, as claimed.
\end{proof}

\section*{1c. Complex variables} 

As a next topic for this chapter, let us discuss now the complex analogues of all the above. To start with, we have the following definition:

\index{complex Gaussian law}
\index{complex normal law}

\begin{definition}
The complex Gaussian law of parameter $t>0$ is
$$G_t=law\left(\frac{a+ib}{\sqrt{2}}\right)$$
where $a,b$ are independent, each following the law $g_t$.
\end{definition}

Here the $\sqrt{2}$ factor comes for having the correct variances in the complex CLT, to be discussed in a moment. As in the real case, these laws form convolution semigroups:

\index{convolution semigroup}

\begin{proposition}
The complex Gaussian laws have the property
$$G_s*G_t=G_{s+t}$$
for any $s,t>0$, and so they form a convolution semigroup.
\end{proposition}

\begin{proof}
This follows indeed from the real result, namely $g_s*g_t=g_{s+t}$, established in Proposition 1.8, simply by taking real and imaginary parts.
\end{proof}

We have as well the following complex analogue of the CLT:

\index{CCLT}
\index{Complex CLT}

\begin{theorem}[CCLT]
Given complex variables $f_1,f_2,f_3,\ldots\in L^\infty(X)$ which are i.i.d., centered, and with common variance $t>0$, we have
$$\frac{1}{\sqrt{n}}\sum_{i=1}^nf_i\sim G_t$$
with $n\to\infty$, in moments.
\end{theorem}

\begin{proof}
This follows indeed from the real CLT, established in Theorem 1.9, simply by taking the real and imaginary parts of all variables involved.
\end{proof}

Regarding now the moments, the situation is more complicated than in the real case, because in order to have good results, we have to take into account both the complex variables, and their conjugates. Let us formulate the following definition:

\index{colored integers}
\index{colored moments}

\begin{definition}
The moments a complex variable $f\in L^\infty(X)$ are the numbers
$$M_k=E(f^k)$$
depending on colored integers $k=\circ\bullet\bullet\circ\ldots\,$, with the conventions
$$f^\emptyset=1\quad,\quad f^\circ=f\quad,\quad f^\bullet=\bar{f}$$
and multiplicativity, in order to define the colored powers $f^k$.
\end{definition}

Observe that, since the variables $f,\bar{f}$ commute, we can permute terms, and restrict the attention to exponents of type $k=\circ\ldots\circ\bullet\ldots\bullet\,$, if we want to. However, our result about the complex Gaussian laws, and other complex laws, later on, will actually look better without doing this, so we will use Definition 1.18 as stated. We first have:

\begin{theorem}
The moments of the complex normal law are given by
$$M_k(G_t)=\begin{cases}
t^pp!&(k\ {\rm uniform, of\ length}\ 2p)\\
0&(k\ {\rm not\ uniform})
\end{cases}$$
where $k=\circ\bullet\bullet\circ\ldots$ is called uniform when it contains the same number of $\circ$ and $\bullet$.
\end{theorem}

\begin{proof}
We must compute the moments, with respect to colored integer exponents $k=\circ\bullet\bullet\circ\ldots$\,, of the variable from Definition 1.15, namely:
$$f=\frac{a+ib}{\sqrt{2}}$$

We can assume that we are in the case $t=1$, and the proof here goes as follows:

\medskip

(1) As a first observation, in the case where our exponent $k=\circ\bullet\bullet\circ\ldots$ is not uniform, a standard rotation argument shows that the corresponding moment of $f$ vanishes. To be more precise, the variable $f'=wf$ is complex Gaussian too, for any complex number $w\in\mathbb T$, and from $M_k(f)=M_k(f')$ we obtain $M_k(f)=0$, in this case.

\medskip

(2) In the uniform case now, where the exponent $k=\circ\bullet\bullet\circ\ldots$ consists of $p$ copies of $\circ$ and $p$ copies of $\bullet$\,, the corresponding moment can be computed as follows:
\begin{eqnarray*}
M_k
&=&\frac{1}{2^p}E[(a^2+b^2)^p]\\
&=&\frac{1}{2^p}\sum_r\binom{p}{r}(2r)!!(2p-2r)!!\\
&=&\frac{p!}{4^p}\sum_r\binom{2r}{r}\binom{2p-2r}{p-r}
\end{eqnarray*}

(3) In order to finish now the computation, let us recall that we have the following formula, coming from the generalized binomial formula, at exponent $-1/2$:
$$\frac{1}{\sqrt{1+t}}=\sum_q\binom{2q}{q}\left(\frac{-t}{4}\right)^q$$

By taking the square of this series, we obtain the following formula:
$$\frac{1}{1+t}
=\sum_{qr}\binom{2q}{q}\binom{2r}{r}\left(\frac{-t}{4}\right)^{q+r}
=\sum_p\left(\frac{-t}{4}\right)^p\sum_r\binom{2r}{r}\binom{2p-2r}{p-r}$$

Now by looking at the coefficient of $t^p$ on both sides, we conclude that the sum on the right equals $4^p$. Thus, we can finish the moment computation in (2), as follows:
$$M_k=\frac{p!}{4^p}\times 4^p=p!$$

We are therefore led to the conclusion in the statement.
\end{proof}

As before with the real Gaussian laws, a better-looking statement  is in terms of partitions. Given a colored integer $k=\circ\bullet\bullet\circ\ldots\,$, we say that a pairing $\pi\in P_2(k)$ is matching when it pairs $\circ-\bullet$ symbols. With this convention, we have the following result:

\index{matching pairings}

\begin{theorem}
The moments of the complex normal law are the numbers
$$M_k(G_t)=\sum_{\pi\in\mathcal P_2(k)}t^{|\pi|}$$
where $\mathcal P_2(k)$ are the matching pairings of $\{1,\ldots,k\}$, and $|.|$ is the number of blocks.
\end{theorem}

\begin{proof}
This is a reformulation of Theorem 1.19. Indeed, we can assume that we are in the case $t=1$, and here we know from Theorem 1.19 that the moments are:
$$M_k=\begin{cases}
(|k|/2)!&(k\ {\rm uniform})\\
0&(k\ {\rm not\ uniform})
\end{cases}$$

On the other hand, the numbers $|\mathcal P_2(k)|$ are given by exactly the same formula. Indeed, in order to have a matching pairing of $k$, our exponent $k=\circ\bullet\bullet\circ\ldots$ must be uniform, consisting of $p$ copies of $\circ$ and $p$ copies of $\bullet$, with $p=|k|/2$. But then the matching pairings of $k$ correspond to the permutations of the $\bullet$ symbols, as to be matched with $\circ$ symbols, and so we have $p!$ such pairings. Thus, we have the same formula as for the moments of $f$, and we are led to the conclusion in the statement.
\end{proof}

In practice, we will also need to know how to compute joint moments of independent normal variables. We have here the following result, to be heavily used later on:

\index{complex normal law}
\index{colored moments}
\index{matching pairings}
\index{Wick formula}

\begin{theorem}[Wick formula]
Given independent variables $f_i$, each following the complex normal law $G_t$, with $t>0$ being a fixed parameter, we have the formula
$$E\left(f_{i_1}^{k_1}\ldots f_{i_s}^{k_s}\right)=t^{s/2}\#\left\{\pi\in\mathcal P_2(k)\Big|\pi\leq\ker i\right\}$$
where $k=k_1\ldots k_s$ and $i=i_1\ldots i_s$, for the joint moments of these variables, where $\pi\leq\ker i$ means that the indices of $i$ must fit into the blocks of $\pi$, in the obvious way.
\end{theorem}

\begin{proof}
This is something well-known, which can be proved as follows:

\medskip

(1) Let us first discuss the case where we have a single variable $f$, which amounts in taking $f_i=f$ for any $i$ in the formula in the statement. What we have to compute here are the moments of $f$, with respect to colored integer exponents $k=\circ\bullet\bullet\circ\ldots\,$, and the formula in the statement tells us that these moments must be:
$$E(f^k)=t^{|k|/2}|\mathcal P_2(k)|$$

But this is the formula in Theorem 1.20, so we are done with this case.

\medskip

(2) In general now, when expanding the product $f_{i_1}^{k_1}\ldots f_{i_s}^{k_s}$ and rearranging the terms, we are left with doing a number of computations as in (1), and then making the product of the expectations that we found. But this amounts in counting the partitions in the statement, with the condition $\pi\leq\ker i$ there standing for the fact that we are doing the various type (1) computations independently, and then making the product.
\end{proof}

Getting now to geometric aspects, involving spheres, we first have:

\index{complex sphere}

\begin{theorem}
We have the following integration formula over the complex unit sphere $S^{N-1}_\mathbb C\subset\mathbb C^N$, with respect to the mass $1$ uniform measure, 
$$\int_{S^{N-1}_\mathbb C}|z_1|^{2k_1}\ldots|z_N|^{2k_N}\,dz=\frac{(N-1)!k_1!\ldots k_n!}{(N+\sum k_i-1)!}$$
valid for any exponents $k_i\in\mathbb N$. As for the other polynomial integrals in $z_1,\ldots,z_N$ and their conjugates $\bar{z}_1,\ldots,\bar{z}_N$, these all vanish.
\end{theorem}

\begin{proof}
By using transformations of type $z\to\lambda z$ with $|\lambda|=1$, we can see that our integral $I$ vanishes, unless each $z_a$ appears as many times as $\bar{z}_a$ does, and this gives the last assertion. So, assume now that we are in the non-vanishing case. Then the $k_a$ copies of $z_a$ and the $k_a$ copies of $\bar{z}_a$ produce by multiplication a factor $|z_a|^{2k_a}$, so we have:
$$I=\int_{S^{N-1}_\mathbb C}|z_1|^{2k_1}\ldots|z_N|^{2k_N}\,dz$$

Now by using the standard identification $S^{N-1}_\mathbb C\simeq S^{2N-1}_\mathbb R$, we obtain:
\begin{eqnarray*}
I
&=&\int_{S^{2N-1}_\mathbb R}(x_1^2+y_1^2)^{k_1}\ldots(x_N^2+y_N^2)^{k_N}\,d(x,y)\\
&=&\sum_{r_1\ldots r_N}\binom{k_1}{r_1}\ldots\binom{k_N}{r_N}\int_{S^{2N-1}_\mathbb R}x_1^{2k_1-2r_1}y_1^{2r_1}\ldots x_N^{2k_N-2r_N}y_N^{2r_N}\,d(x,y)\\
&=&\sum_{r_1\ldots r_N}\binom{k_1}{r_1}\ldots\binom{k_N}{r_N}\frac{(2N-1)!!(2r_1)!!\ldots(2r_N)!!(2k_1-2r_1)!!\ldots (2k_N-2r_N)!!}{(2N+2\sum k_i-1)!!}\\
&=&\sum_{r_1\ldots r_N}\binom{k_1}{r_1}\ldots\binom{k_N}{r_N}
\frac{(N-1)!(2r_1)!\ldots (2r_N)!(2k_1-2r_1)!\ldots (2k_N-2r_N)!}{4^{\sum k_i}(N+\sum k_i-1)!r_1!\ldots r_N!(k_1-r_1)!\ldots (k_N-r_N)!}
\end{eqnarray*}

Observe now that can rewrite this quantity in the following way:
\begin{eqnarray*}
I
&=&\sum_{r_1\ldots r_N}\frac{k_1!\ldots k_N!(N-1)!(2r_1)!\ldots (2r_N)!(2k_1-2r_1)!\ldots (2k_N-2r_N)!}{4^{\sum k_i}(N+\sum k_i-1)!(r_1!\ldots r_N!(k_1-r_1)!\ldots (k_N-r_N)!)^2}\\
&=&\sum_{r_1}\binom{2r_1}{r_1}\binom{2k_1-2r_1}{k_1-r_1}\ldots\sum_{r_N}\binom{2r_N}{r_N}\binom{2k_N-2r_N}{k_N-r_N}\frac{(N-1)!k_1!\ldots k_N!}{4^{\sum k_i}(N+\sum k_i-1)!}\\
&=&4^{k_1}\times\ldots\times 4^{k_N}\times\frac{(N-1)!k_1!\ldots k_N!}{4^{\sum k_i}(N+\sum k_i-1)!}\\
&=&\frac{(N-1)!k_1!\ldots k_N!}{(N+\sum k_i-1)!}
\end{eqnarray*}

Thus, we obtain the formula in the statement.
\end{proof}

We can now formulate the following complex analogue of Theorem 1.14:

\begin{theorem}
The moments of the complex hyperspherical variables are
$$\int_{S^{N-1}_\mathbb C}|z_i|^{2k}dz=\frac{(N-1)!k!}{(N+k-1)!}$$
and the variables $y_i=\sqrt{N}z_i$ become complex normal and independent with $N\to\infty$.
\end{theorem}

\begin{proof}
According to the formula in Theorem 1.22, the polynomials integrals in $z_i,\bar{z}_i$ vanish, unless the number of $z_i,\bar{z}_i$ is the same, and in this latter case, we have:
$$\int_{S^{N-1}_\mathbb C}|z_i|^{2k}\,dz=\frac{(N-1)!k!}{(N+k-1)!}
\simeq N^{-k}\times k!$$

Thus, the rescaled variables $\sqrt{N}z_i$ become normal with $N\to\infty$, as claimed. Next, again by using Theorem 1.22, the joint moments of $z_1,\ldots,z_N$ are given by:
$$\int_{S^{N-1}_\mathbb R}|z_1|^{2k_1}\ldots|z_N|^{2k_N}\,dx
=\frac{(N-1)!k_1!\ldots k_n!}{(N+\sum k_i-1)!}
\simeq N^{-\sum k_i}\times k_1!\ldots k_N!$$

By rescaling, the joint moments of the variables $y_i=\sqrt{N}z_i$ are given by:
$$\int_{S^{N-1}_\mathbb R}|y_1|^{2k_1}\ldots|y_N|^{2k_N}\,dx\simeq k_1!\ldots k_N!$$

Thus, we have multiplicativity, and so independence with $N\to\infty$, as claimed.
\end{proof}

\section*{1d. Poisson laws}

Let us discuss now the discrete counterpart of the above results. The mathematics here will involve the Poisson laws $p_t$, which are constructed as follows:

\index{Poisson law}

\begin{definition}
The Poisson law of parameter $t>0$ is the following measure,
$$p_t=e^{-t}\sum_{k\in\mathbb N}\frac{t^k}{k!}\,\delta_k$$
with the letter ``p'' standing for Poisson.
\end{definition}

Many interesting things can be said about the Poisson laws, and going now directly to the essentials, here is the technical result that we will need, for our purposes:

\index{convolution semigroup}
\index{Fourier transform}

\begin{proposition}
The Fourier transform of $p_t$ is given by
$$F_{p_t}(x)=\exp\left((e^{ix}-1)t\right)$$
for any $t>0$. In particular we have $p_s*p_t=p_{s+t}$, for any $s,t>0$.
\end{proposition}

\begin{proof}
We have indeed the following computation, as desired:
$$F_{p_t}(x)
=e^{-t}\sum_k\frac{t^k}{k!}\,e^{ikx}
=e^{-t}\exp(e^{ix}t)
=\exp\left((e^{ix}-1)t\right)$$

As for $p_s*p_t=p_{s+t}$, this comes from this, $\log F_{p_t}$ being linear in $t$.
\end{proof}

We can now establish the Poisson Limit Theorem, as follows:

\index{PLT}
\index{Poisson Limit Theorem}
\index{Bernoulli laws}

\begin{theorem}[PLT]
We have the following convergence, in moments,
$$\left(\left(1-\frac{t}{n}\right)\delta_0+\frac{t}{n}\,\delta_1\right)^{*n}\to p_t$$
for any $t>0$.
\end{theorem}

\begin{proof}
Let us denote by $\nu_n$ the measure under the convolution sign. We have the following computation, for the Fourier transform of the limit: 
\begin{eqnarray*}
F_{\nu_n}(x)=\left(1-\frac{t}{n}\right)+\frac{t}{n}\,e^{ix}
&\implies&F_{\nu_n^{*n}}(x)=\left(\left(1-\frac{t}{n}\right)+\frac{t}{n}\,e^{ix}\right)^n\\
&\implies&F_{\nu_n^{*n}}(x)=\left(1+\frac{(e^{ix}-1)t}{n}\right)^n\\
&\implies&F(x)=\exp\left((e^{ix}-1)t\right)
\end{eqnarray*}

Thus, we obtain indeed the Fourier transform of $p_t$, as desired.
\end{proof}

At the level of moments now, things are quite tricky for the Poisson laws:

\index{Bell numbers}
\index{Stirling numbers}

\begin{theorem}
The moments of $p_t$ with $t>0$ are given by:
$$M_k(p_t)=\sum_{\pi\in P(k)}t^{|\pi|}$$
In particular at $t=1$ we have $M_k(p_1)=|P(k)|$, the Bell numbers.
\end{theorem}

\begin{proof}
We have the following recurrence formula for the moments:
\begin{eqnarray*}
M_{k+1}
&=&e^{-t}\sum_{m\geq0}\frac{t^{m+1}(m+1)^k}{m!}\\
&=&e^{-t}\sum_{m\geq0}\frac{t^{m+1}m^k}{m!}\left(1+\frac{1}{m}\right)^k\\
&=&e^{-t}\sum_{m\geq0}\frac{t^{m+1}m^k}{m!}\sum_{s=0}^k\binom{k}{s}m^{-s}\\
&=&\sum_{s=0}^k\binom{k}{s}\cdot e^{-t}\sum_{m\geq0}\frac{t^{m+1}m^{k-s}}{m!}\\
&=&t\sum_{s=0}^k\binom{k}{s}M_{k-s}
\end{eqnarray*}

But $N_k=\sum_{\pi\in P(k)}t^{|\pi|}$ satisfy the same recurrence, and we get the result.
\end{proof}

In order to discuss the occurrence of Poisson laws in geometry, let us start with:

\index{character}
\index{main character}
\index{truncated character}

\begin{definition}
Given a closed subgroup $G\subset U_N$, the function
$$\chi_t:G\to\mathbb C\quad,\quad 
\chi_t(g)=\sum_{i=1}^{[tN]}g_{ii}$$
is called main truncated character of $G$, of parameter $t\in(0,1]$.
\end{definition}

As a first observation, we actually already met such truncated characters, but in a disguised form, when talking about spheres, with our results there implying:

\index{normal law}
\index{orthogonal group}
\index{unitary group}

\begin{proposition}
For the orthogonal and unitary groups $O_N,U_N$, the rescalings
$$\rho_{1/N}=\sqrt{N}\cdot\chi_{1/N}$$
become respectively real and complex Gaussian, in the $N\to\infty$ limit.
\end{proposition}

\begin{proof}
According to our conventions, given a closed subgroup $G\subset U_N$, the main character truncated at $t=1/N$ is simply the first coordinate:
$$\chi_{1/N}(g)=g_{11}$$

On the other hand, the standard action $O_N\curvearrowright S^{N-1}_\mathbb R$ provides us with an embedding as follows, which makes correspond the respective integration functionals:
$$C(S^{N-1}_\mathbb R)\subset C(O_N)\quad,\quad 
x_i\to g_{1i}$$

But with this identification made, the result for $O_N$ follows from Theorem 1.14. As for discussion for $U_N$, this is similar, using this time Theorem 1.23.
\end{proof}

We will see later in this book that for $O_N$ we have $\chi_t\sim g_t$, and that for $U_N$ we have $\chi_t\sim G_t$, in the $N\to\infty$ limit. Now back to Poisson, we have the following result:

\index{Poisson law}

\begin{theorem}
Consider the symmetric group $S_N$, regarded as the permutation group, $S_N\subset O_N$, of the $N$ coordinate axes of $\mathbb R^N$.
\begin{enumerate}
\item The variable $\chi_t$ counts the number of fixed points among $1,\ldots,[tN]$.

\item The law of this variable $\chi_t$ becomes Poisson $(t)$, in the $N\to\infty$ limit.
\end{enumerate}
\end{theorem}

\begin{proof}
We have the following computation, coming from definitions:
$$\chi_t(g)
=\sum_{i=1}^{[tN]}\delta_{g(i)i}
=\#\left\{i\in\{1,\ldots,[tN]\}\Big|g(i)=i\right\}$$

As for the second assertion, there are several proofs here, all instructive:

\medskip

(1) Inclusion-exclusion. Let us first compute $P(\chi_1=0)$, which is the probability for a random $g\in S_N$ to be a derangement. With $S_N^i=\left\{g\in S_N|g(i)=i\right\}$, we have:
\begin{eqnarray*}
P(\chi_1=0)
&=&\frac{1}{N!}\left(|S_N|-\sum_i|S_N^i|+\sum_{i<j}|S_N^i\cap S_N^j|-\sum_{i<j<k}|S_N^i\cap S_N^j\cap S_N^k|+\ldots\right)\\
&=&\frac{1}{N!}\sum_{r=0}^N(-1)^r\binom{N}{r}(N-r)!\\
&=&\sum_{r=0}^N\frac{(-1)^r}{r!}
\end{eqnarray*}

Thus we have $P(\chi_1=0)\simeq 1/e$, and it is then straightforward to generalize this into $P(\chi_1=k)\simeq 1/(k!e)$, and then into $P(\chi_t=k)\simeq 1/(k!e^t)$, leading to the result.

\medskip

(2) Integration. The polynomial integrals over $S_N\subset O_N$ are given by:
\begin{eqnarray*}
\int_{S_N}g_{i_1j_1}\ldots g_{i_kj_k}\,dg
&=&\frac{1}{N!}\#\left\{g\in S_N\Big|g(j_1)=i_1,\ldots,g(j_k)=i_k\right\}\\
&=&\delta_{\ker i=\ker j}\,\frac{(N-|\ker i|)!}{N!}
\end{eqnarray*}

But with this in hand, we have the following formula for the moments:
\begin{eqnarray*}
\int_{S_N}\chi_t^k
&=&\sum_{i_1,\ldots,i_k=1}^{[tN]}\int_{S_N}g_{i_1i_1}\ldots g_{i_ki_k}\,dg\\
&=&\sum_{\pi\in P(k)}\sum_{i_1,\ldots,i_k=1}^{[tN]}\delta_{\ker i=\pi}\,\frac{(N-|\pi|)!}{N!}
\end{eqnarray*}

With $S_{kb}$ being the Stirling numbers, counting the partitions having $b$ blocks, we get:
\begin{eqnarray*}
\int_{S_N}\chi_t^k
&=&\sum_{\pi\in P(k)}\frac{[tN]!}{([tN]-|\pi|!)}\cdot\frac{(N-|\pi|!)}{N!}\\
&=&\sum_{b=1}^{[tN]}\frac{[tN]!}{([tN]-b)!}\cdot\frac{(N-b)!}{N!}\cdot S_{kb}\\
&\simeq&\sum_{b=1}^kS_{kb}t^b
\end{eqnarray*}

But this is the $k$-th moment of the Poisson law $p_t$, so we are done.

\medskip

(3) Exact count. A variation of the count in (1), that we will leave as an exercise, leads to the following exact formula, for the law that we are interested in:
$$law(\chi_t)=\frac{[tN]!}{N!}\sum_{p=0}^{[tN]}\frac{(N-p)!}{([tN]-p)!}
\cdot\frac{\left(\delta_1-\delta_0\right)^{*p}}{p!}$$ 

The coefficients in this law formula can be estimated via Stirling, as follows:
$$c_p
\simeq\frac{(tN)^{tN}}{N^N}\cdot\frac{(N-p)^{N-p}}{(tN-p)^{tN-p}}
\simeq t^p$$

We conclude that the asymptotic law is given by the following formula:
$$law(\chi_t)\simeq \sum_{p=0}^st^p\cdot\frac{\left(\delta_1-\delta_0\right)^{*p}}{p!}$$

But on the right we have $p_t$, say via $F(x)=\exp(t(e^{ix}-1))$, and we are done.

\medskip

(4) Easiness. It is not hard to prove that $S_N\subset_u U_N$ is easy, coming from the category of partitions $P$, in the sense that we have the following formula, for any $k\in\mathbb N$, with $T_\pi=\sum_i\delta_\pi(i)e_{i_1}\otimes\ldots\otimes e_{i_k}$, where $\delta_\pi=0,1$, depending on whether indices fit or not:
$$Fix(u^{\otimes k})=span\left(T_\pi\Big|\pi\in P(k)\right)$$

Now since the vectors $T_\pi$ are well-known to be linearly independent with $N\to\infty$, this gives right away $\int\chi_1^k=|P(k)|$, and so $\chi_1\sim p_1$, in the $N\to\infty$ limit. And with a bit more know-how we obtain, more generally, $\int\chi_t^k=\sum_{\pi\in P(k)}t^{|\pi|}$, and so $\chi_t\sim p_t$, in the $N\to\infty$ limit, for any $t\in(0,1]$. And more on this, easiness, later in this book.
\end{proof}

As a last question that we would like to investigate, in this chapter, what is the complex analogue of the Poisson law? In answer, things here are quite tricky, and we first have the following notion, extending the Poisson limit theory developed above:

\index{compound Poisson law}

\begin{definition}
Associated to any compactly supported positive measure $\nu$ on $\mathbb C$, not necessarily of mass $1$, is the probability measure
$$p_\nu=\lim_{n\to\infty}\left(\left(1-\frac{t}{n}\right)\delta_0+\frac{1}{n}\,\nu\right)^{*n}$$
where $t=mass(\nu)$, called compound Poisson law.
\end{definition}

In what follows we will be mainly interested in the case where the measure $\nu$ is discrete, as is for instance the case for $\nu=t\delta_1$ with $t>0$, which produces the Poisson laws. The following standard result allows one to detect compound Poisson laws:

\index{Fourier transform}

\begin{proposition}
For $\nu=\sum_{i=1}^st_i\delta_{z_i}$ with $t_i>0$ and $z_i\in\mathbb C$, we have
$$F_{p_\nu}(x)=\exp\left(\sum_{i=1}^st_i(e^{ixz_i}-1)\right)$$
where $F$ denotes as usual the Fourier transform.
\end{proposition}

\begin{proof}
Let $\eta_n$ be the measure in Definition 1.31, under the $*$ sign. We have:
\begin{eqnarray*}
F_{\eta_n}(x)=\left(1-\frac{t}{n}\right)+\frac{1}{n}\sum_{i=1}^st_ie^{ixz_i}
&\implies&F_{\eta_n^{*n}}(x)=\left(\left(1-\frac{t}{n}\right)+\frac{1}{n}\sum_{i=1}^st_ie^{ixz_i}\right)^n\\
&\implies&F_{p_\nu}(x)=\exp\left(\sum_{i=1}^st_i(e^{ixz_i}-1)\right)
\end{eqnarray*}

Thus, we have obtained the formula in the statement.
\end{proof}

Next, we have the following result, providing an alternative to Definition 1.31, and which will be our formulation here of the Compound Poisson Limit Theorem:

\index{CPLT}
\index{Compound Poisson Limit Theorem}

\begin{theorem}[CPLT]
For $\nu=\sum_{i=1}^st_i\delta_{z_i}$ with $t_i>0$ and $z_i\in\mathbb C$, we have
$$p_\nu=law\left(\sum_{i=1}^sz_i\alpha_i\right)$$
with the variables $\alpha_i$ being Poisson $(t_i)$, independent.
\end{theorem}

\begin{proof}
Let $\alpha$ be the sum of Poisson variables in the statement. We have:
\begin{eqnarray*}
F_{\alpha_i}(x)=\exp(t_i(e^{ix}-1))
&\implies&F_{z_i\alpha_i}(x)=\exp(t_i(e^{ixz_i}-1))\\
&\implies&F_\alpha(x)=\exp\left(\sum_{i=1}^st_i(e^{ixz_i}-1)\right)
\end{eqnarray*}

Thus we have indeed the same formula as in Proposition 1.32, as desired.
\end{proof}

At the level of main examples now, following \cite{bb+}, we have the following result:

\index{Bessel law}
\index{Bessel function}
\index{complex reflection}

\begin{theorem}
The Bessel laws, $p^s_t=p_{t\varepsilon_s}$ with $\varepsilon_s$ being the uniform measure on the $s$-th roots of unity, appear as follows, with $\alpha_i\sim p_{t/s}$ independent, and $w=e^{2\pi i/s}$:
$$p^s_t=law\left(\sum_{i=1}^sw^i\alpha_i\right)$$
The Fourier transform is as follows, with the convention $\exp_sx=\sum_{k=0}^\infty x^{sk}/(sk)!$,
$$F_{p^s_t}(x)=\exp\left(t\left(\exp_s(ix)-1\right)\right)$$
and in particular we have $p^s_t*p^s_{t'}=p^s_{t+t'}$. The density is given by the formula
$$p^s_t=e^{-t}\sum_{c_1=0}^\infty\ldots\sum_{c_s=0}^\infty\frac{1}{c_1!\ldots c_s!}\,\left(\frac{t}{s}\right)^{c_1+\ldots+c_s}\delta\left(\sum_{k=1}^sw^kc_k\right)$$
with at $s=2$ the density function being the Bessel function of the first kind:
$$p^2_t=e^{-t}\sum_{r=-\infty}^\infty\delta_r\sum_{p=0}^\infty \frac{(t/2)^{|r|+2p}}{(|r|+p)!p!}$$
The moments are as follows, $P^s$ being the partitions satisfying $\#\circ=\#\bullet(s)$ in each block:
$$M_k(p_t^s)=\sum_{\pi\in P^s(k)}t^{|\pi|}$$
Also, for the complex reflection group $H_N^s=S_N\wr\mathbb Z_s$, we have $\chi_t\sim p^s_t$, with $N\to\infty$.
\end{theorem}

\begin{proof}
This is something quite routine, and fundamental too, generalizing what we know about the Poisson laws $p_t$, in the complex setting, the idea being as follows:

\medskip

(1) To start with, we can talk about $p^s_t=p_{t\varepsilon_s}$ as above, with $s\in\{1,2,\ldots,\infty\}$, with the convention that at $s=\infty$, the $s$-th roots of unity are the unit circle $\mathbb T$.

\medskip

(2) In the discrete case, $s<\infty$, we have the Poisson sum formula, and the Fourier transform formula too, coming from Proposition 1.32 and Theorem 1.33.

\medskip

(3) Since $\log F_{p^s_t}$ is linear in $t$, we obtain $p^s_t*p^s_{t'}=p^s_{t+t'}$, at $s<\infty$. And this formula holds as $s=\infty$ too, say appearing as the $s\to\infty$ limit of the $s<\infty$ ones.

\medskip

(4) Regarding now the formula of the density, at $s<\infty$, the Fourier transform of the measure on the right in the statement is given by the following formula:
\begin{eqnarray*}
F(x)
&=&e^{-t}\sum_{p_1=0}^\infty\ldots\sum_{p_s=0}^\infty\frac{1}{p_1!\ldots p_s!}\left(\frac{t}{s}\right)^{p_1+\ldots+p_s}\exp\left(\sum_{k=1}^siw^kp_kx\right)\\
&=&e^{-t}\sum_{r=0}^\infty\left(\frac{t}{s}\right)^r\sum_{\sum p_i=r}\frac{\exp\left(\sum_{k=1}^siw^kp_kx\right)}{p_1!\ldots p_s!}
\end{eqnarray*}

If we multiply by $e^t$, and compute the derivative with respect to $t$, we get:
\begin{eqnarray*}
(e^tF(x))'
&=&\frac{1}{s}\sum_{r=1}^\infty\left(\frac{t}{s}\right)^{r-1}\sum_{\sum p_i=r}\sum_{l=1}^s\frac{\exp\left(\sum_{k=1}^siw^kp_kx\right)}{p_1!\ldots p_{l-1}!(p_l-1)!p_{l+1}!\ldots p_s!}\\
&=&\frac{1}{s}\sum_{u=0}^\infty\left(\frac{t}{s}\right)^u\sum_{\sum q_i=u}\sum_{l=1}^s\frac{\exp\left(iw^lx+\sum_{k=1}^siw^kq_kx\right)}{q_1!\ldots q_s!}\\
&=&\exp_s(ix)\cdot e^tF(x)
\end{eqnarray*}

But $\Phi(t)=\exp(t\exp_s(ix))$ satisfies as well $\Phi'(t)=(\exp_s(ix))\Phi(t)$, so we have:
$$e^tF(x)=\exp(t\exp_s(ix))$$

Thus $F(x)=\exp(t(\exp_s(ix)-1))$, so the measure in statement is indeed $p^s_t$. And in particular, we can see that at $s=2$ we have a Bessel function density, as stated.

\medskip

(5) Next, the asymptotic truncated character formula $\chi_t\sim p^s_t$ for $H_N^s=S_N\wr\mathbb Z_s$ can be established, by using $\chi_t\sim p_t$ for $S_N$, via a direct density computation. 

\medskip

(6) Finally, the moments are best found via $H_N^s$, the idea being that this group is easy, coming from the category $P^s$. And more on this, easiness, later in this book.
\end{proof}

\section*{1e. Exercises}

Welcome to probability, are here are a few exercises, in relation with the above:

\begin{exercise}
Learn about the binomial laws, produced by coin tosses.
\end{exercise}

\begin{exercise}
What is the law of the sum, when rolling a die $n$ times?
\end{exercise}

\begin{exercise}
Learn about the inequalities of Markov and Chebycheff.
\end{exercise}

\begin{exercise}
Learn about the law of large numbers, weak and strong.
\end{exercise}

\begin{exercise}
Meditate on that $\sqrt{2}$ factor, in the definition of $G_t$.
\end{exercise}

\begin{exercise}
Assuming $f\sim g_t$ or $f\sim G_t$, compute the laws of $|f|$ and $|f|^2$.
\end{exercise}

\begin{exercise}
Have all $4$ proofs of $\chi_t\sim p_t$ with $N\to\infty$ for $S_N$ understood.
\end{exercise}

\begin{exercise}
Work out the details, in our last theorem, regarding $p^s_t$.
\end{exercise}

As bonus exercise, for the axiomatics, have a look at Durrett \cite{dur} or Feller \cite{fel}.

\chapter{Matrices, laws}

\section*{2a. Linear algebra}

Free probability, that we would like to learn in this book, is about random variables $f,g:X\to\mathbb C$ which are free, in the sense that there is no algebraic relation between them. In particular, these variables must not commute, $fg\neq gf$, and this is already a contradiction, because we perfectly know that functions commute, $fg=gf$.

\bigskip

In answer, free probability is in fact about certain ``generalized random variables'' $f,g\in B$ which are free, with the algebras $B$ that we can use being something more general than the usual random variable algebras $B=L^\infty(X)$ that we are used to.

\bigskip

So, what are the basic examples? Looking around, the simplest algebras where we can find noncommutativity are the matrix algebras $M_N(\mathbb C)$. But, can we do probability over such algebras? And in answer, yes, as shown by the following speculation:

\begin{speculation}
We can do probability, in a generalized sense, over the matrix algebras $M_N(\mathbb C)$, by declaring that the expectation is the normalized trace
$$E(A)=\frac{Tr(A)}{N}$$
which is positive and unital, as it should. In this framework any matrix $A\in M_N(\mathbb C)$ can be regarded as a generalized random variable $A:M_N\to\mathbb C$, has moments
$$M_k=E(A^k)$$
and has even a law, encoding the collection of moments $\{M_k\}$. Moreover, this can be generalized to $M_N(L^\infty(X))$, and even to $B(H)$, leading to free probability.
\end{speculation}

To be more precise here, the functional $E=Tr/N$ is notoriously positive and unital, and more on this in a moment. Regarding the convention $A:M_N\to\mathbb C$, this is something intuitive, and we will make later something rigorous out of this. Next, we can certainly talk about moments $M_k=E(A^k)$, and by being a bit abstract, and in the lack of a good idea on the spot, we can simply declare that the law is the collection of moments $\{M_k\}$. Finally, the assertions at the end, regarding random matrices, linear operators and freeness, are something a bit more specialized, that we can leave for later. 

\bigskip

So, can we make something rigorous out of Speculation 2.1? In answer, yes, and we will do this here. Then in chapters 3-4 we will study some basic examples of random matrices, and with a look into Lie groups as well, as to get introduced to freeness. And then, in chapter 5 and afterwards we will axiomatize and study free probability.

\bigskip

Getting started now, we must first talk about positivity, in order to understand why the expectation $E(A)=Tr(A)/N$ is indeed positive. What we would need here is:
$$A>0\implies Tr(A)>0$$

Now in order to discuss this, we first need to know about self-adjointness:
$$A=A^*\implies Tr(A)\in\mathbb R$$

So, let us have as starting point the construction of the adjoint matrix, $A\to A^*$. You certainly know about this, with the basics of this operation being as follows:

\index{adjoint operator}
\index{adjoint matrix}
\index{scalar product}
\index{unitary}
\index{projection}

\begin{theorem}
Given $A\in M_N(\mathbb C)$, define its adjoint matrix by $(A^*)_{ij}=\bar{A}_{ji}$.
\begin{enumerate}
\item $<Ax,y>=<x,A^*y>$, for any two vectors $x,y\in\mathbb C^N$.

\item $x\to Ux$ with $U\in M_N(\mathbb C)$ is an isometry precisely when $U^*=U^{-1}$.

\item $x\to Px$ with $P\in M_N(\mathbb C)$ is a projection precisely when $P^2=P^*=P$.
\end{enumerate}
\end{theorem}

\begin{proof}
This is something very standard, the idea being as follows:

\medskip

(1) The formula $<Ax,y>=<x,A^*y>$ is indeed clear from definitions on the basis vectors, $x,y\in\{e_1,\ldots,e_N\}$, and then by linearity, it must hold in general, as stated.

\medskip

(2) We recall that the norm can be recovered from the scalar product, thanks to $||x||=\sqrt{<x,x>}$, and vice versa too, thanks to the polarization identity, namely:
\begin{eqnarray*}
&&||x+y||^2-||x-y||^2+i||x+iy||^2-i||x-iy||^2\\
&=&||x||^2+||y||^2-||x||^2-||y||^2+i||x||^2+i||y||^2-i||x||^2-i||y||^2\\
&&+2Re(<x,y>)+2Re(<x,y>)+2iIm(<x,y>)+2iIm(<x,y>)\\
&=&4<x,y>
\end{eqnarray*}

Now given a matrix $U\in M_N(\mathbb C)$, we have indeed the following equivalences, with the first one coming from the polarization identity, and the other ones being clear:
\begin{eqnarray*}
||Ux||=||x||
&\iff&<Ux,Uy>=<x,y>\\
&\iff&<x,U^*Uy>=<x,y>\\
&\iff&U^*Uy=y\\
&\iff&U^*U=1\\
&\iff&U^*=U^{-1}
\end{eqnarray*}

(3) Given a matrix $P\in M_N(\mathbb C)$, in order for $x\to Px$ to be an oblique projection, we must have $P^2=P$. Now observe that this projection is orthogonal when:
\begin{eqnarray*}
<Px-x,Py>=0
&\iff&<P^*Px-P^*x,y>=0\\
&\iff&P^*Px-P^*x=0\\
&\iff&P^*P-P^*=0\\
&\iff&P^*P=P^*
\end{eqnarray*}

The point now is that by conjugating the last formula, we get $P^*P=P$. Now by comparing with $P^*P=P^*$, we obtain $P=P^*$, and this gives the result. 
\end{proof}

With this understood, we can talk about positivity, as follows:

\begin{theorem}
In analogy with what happens for the real numbers, where $x>0$ means $x=z\bar{z}$ with $z\neq0$, let us call $A\in M_N(\mathbb C)$ positive, and write $A>0$, when:
$$A=BB^*\quad,\quad B\neq0$$
Then the functional $E:M_N(\mathbb C)\to\mathbb C$ given by $E(A)=Tr(A)/N$ satisfies $E(1)=1$ and
$$A>0\implies Tr(A)>0$$
that is, is unital and positive, exactly as the probabilistic expectations are.
\end{theorem}

\begin{proof}
We have several things going on here, the idea being as follows:

\medskip

(1) The first assertion stands as a definition for the positive matrices, and of course in the hope that the argument there is convincive. As a further piece of evidence, observe that at $N=1$ we recover the usual positivity notion for the real numbers. On top of this, we will learn in a moment that a matrix $A\in M_N(\mathbb C)$ is positive precisely when it is diagonalizable via a unitary, with positive eigenvalues $\lambda_i>0$, which sounds good too. Finally, we will learn as well that positivity also means $A=C^2$ with $C=C^*$, $C\neq0$, which again sounds good, in analogy with what we know about the real numbers.

\medskip

(2) The second assertion is trivial, because with $A=BB^*$ with $B\neq0$, we have:
$$E(A)=\frac{Tr(BB^*)}{N}=\frac{1}{N}\sum_{ij}|B_{ij}|^2>0$$

(3) Finally, the conclusion at the end is something subjective, I mean I challenge you to find a general property of the probability expectations $E:L^\infty(X)\to\mathbb C$ that our expectation $E:M_N(\mathbb C)\to\mathbb C$ does not have. Unitality, positivity, trace property $E(AB)=E(BA)$ too, everything is there, like a gift coming from the outer space.
\end{proof}

Getting back to Speculation 2.1, the first claim there, regarding the expectation, is now understood. In what regards the second claim, the fact that we can view any matrix $A\in M_N(\mathbb C)$ as some sort of function $A:M_N\to\mathbb C$, here is about that:

\begin{fact}
If we define the quantum space $M_N$ by the formula $C(M_N)=M_N(\mathbb C)$, then the matrices $A\in M_N(\mathbb C)$ become elements $A\in C(M_N)$, and so functions $A:M_N\to\mathbb C$, in some generalized sense. These functions generically don't commute, $AB\neq BA$.
\end{fact}

And up to you to see if you are satisfied by this. If you are a mathematician, you would probably say that this is some physics garbage, with no mathematical value. On the opposite, if you are a physicist, you would probably say that this is some mathematical garbage, using that fancy term ``quantum'', without knowing what we're talking about. Well, in answer to you guys, I noted your criticism, and we will be back to this, later. 

\bigskip

Getting now to the third claim in Speculation 2.1, moments and law, there is certainly some mathematics to be done there. Here is a precise question in this sense:

\begin{question}
For a diagonal matrix $A=diag(\lambda_1,\ldots,\lambda_N)$ we have
$$E(P(A))=\int_\mathbb CP(z)d\mu_A(z)$$
for any polynomial $P\in\mathbb C<X,X^*>$, with $\mu_A$ being the complex discrete measure
$$\mu_A=\frac{1}{N}\sum_i\delta_{\lambda_i}$$
which can be thought of as being the law of $A$. How far can we generalize this?
\end{question}

To be more precise, the above formula for $E(P(A))$ can be deduced as follows:
$$E(P(A))
=\frac{1}{N}\,Tr\begin{pmatrix}
P(\lambda_1)&&0\\
&\ddots\\
0&&P(\lambda_N)
\end{pmatrix}\\
=\frac{1}{N}\sum_iP(\lambda_i)
=\int_\mathbb CP(z)d\mu_A(z)$$

And the problem is now, how far can we generalize this. In answer, here is a first theorem, called spectral theorem for the self-adjoint matrices, that can be of help:

\begin{theorem}
Any matrix $A\in M_N(\mathbb C)$ which is self-adjoint, $A=A^*$, is diagonalizable, with the diagonalization being of the following type,
$$A=UDU^*$$
with $U\in U_N$, and with $D\in M_N(\mathbb R)$ diagonal. The converse holds too.
\end{theorem}

\begin{proof}
As a first remark, the converse trivially holds, because if we take a matrix of the form $A=UDU^*$, with $U$ unitary and $D$ diagonal and real, then we have:
$$A^*
=(UDU^*)^*
=UD^*U^*
=UDU^*
=A$$

In the other sense now, assume that $A$ is self-adjoint, $A=A^*$.  Our first claim is that the eigenvalues are real. Indeed, assuming $Av=\lambda v$, we have:
\begin{eqnarray*}
\lambda<v,v>
&=&<Av,v>\\
&=&<v,Av>\\
&=&\bar{\lambda}<v,v>
\end{eqnarray*}

Thus we obtain $\lambda\in\mathbb R$, as claimed. Our next claim now is that the eigenspaces corresponding to different eigenvalues are pairwise orthogonal. Assume indeed that:
$$Av=\lambda v\quad,\quad 
Aw=\mu w$$

We have then the following computation, by using the fact that we have $\mu\in\mathbb R$:
\begin{eqnarray*}
\lambda<v,w>
&=&<Av,w>\\
&=&<v,Aw>\\
&=&\mu<v,w>
\end{eqnarray*}

Thus $\lambda\neq\mu$ implies $v\perp w$, as claimed. In order now to finish the proof, it remains to prove that the eigenspaces of $A$ span the whole space $\mathbb C^N$. For this purpose, we will use a recurrence method. Let us pick an eigenvector of our matrix:
$$Av=\lambda v$$

Assuming now that we have a vector $w$ orthogonal to it, $v\perp w$, we have:
\begin{eqnarray*}
<Aw,v>
&=&<w,Av>\\
&=&<w,\lambda v>\\
&=&\lambda<w,v>\\
&=&0
\end{eqnarray*}

Thus, if $v$ is an eigenvector, then the vector space $v^\perp$ is invariant under $A$. Moreover, since a matrix $A$ is self-adjoint precisely when $<Av,v>\in\mathbb R$ for any vector $v\in\mathbb C^N$, as one can see by expanding the scalar product, the restriction of $A$ to the subspace $v^\perp$ is self-adjoint. Thus, we can proceed by recurrence, and we obtain the result.
\end{proof}

As a first piece of answer now to Question 2.5, we have:

\begin{theorem}
For a self-adjoint matrix $A=A^*$ we have the formula
$$E(P(A))=\int_\mathbb RP(z)d\mu_A(z)$$
for any polynomial $P\in\mathbb C<X,X^*>$, with $\mu_A$ being the real discrete measure
$$\mu_A=\frac{1}{N}\sum_i\delta_{\lambda_i}$$
which can be therefore thought of as being the law of $A$.
\end{theorem}

\begin{proof}
If we write $A=UDU^*$ as in Theorem 2.6, with $U\in U_N$ and with $D\in M_N(\mathbb R)$ diagonal, then we have the following computation, using $E(BC)=E(CB)$:
\begin{eqnarray*}
E(P(A))
&=&E(UP(D)U^*)\\
&=&E(U^*UP(D))\\
&=&E(P(D))
\end{eqnarray*}

But with this formula in hand, our result follows from our previous computation for the diagonal matrices, which was something trivial, made after Question 2.5.
\end{proof}

Observe that Theorem 2.7, while being certainly something non-trivial, based on the spectral theorem for self-adjoints, is not exactly an answer to Question 2.5, because what we have is a generalization of the result for the diagonal matrices $D$, but only in the real case, $D\in M_N(\mathbb R)$. Thus, in order to fully answer Question 2.5, still long way to go.

\bigskip

Before getting further into this, as a byproduct of our study, let us record some useful consequences of Theorem 2.6 in relation with matrix positivity, fully clarifying our considerations from Theorem 2.3 and its proof, and putting an end to that discussion:

\begin{theorem}
For a matrix $A\in M_N(\mathbb C)$ the following conditions are equivalent, and if they are satisfied, we say that $A$ is positive:
\begin{enumerate}
\item $A=B^2$, with $B=B^*$.

\item $A=CC^*$, for some $C\in M_N(\mathbb C)$.

\item $<Ax,x>\geq0$, for any vector $x\in\mathbb C^N$.

\item $A=A^*$, and the eigenvalues are positive, $\lambda_i\geq0$.

\item $A=UDU^*$, with $U\in U_N$ and with $D\in M_N(\mathbb R_+)$ diagonal.
\end{enumerate}
\end{theorem}

\begin{proof}
This is something very standard, the idea being as follows:

\medskip

$(1)\implies(2)$ This is clear, because we can take $C=B$.

\medskip

$(2)\implies(3)$ This follows from the following computation:
$$<Ax,x>
=<CC^*x,x>
=<C^*x,C^*x>
\geq0$$

$(3)\implies(4)$ By using the fact that $<Ax,x>$ is real, we have:
$$<Ax,x>
=<x,A^*x>
=<A^*x,x>$$

Thus we have $A=A^*$, and the remaining assertion, regarding the eigenvalues, follows from the following computation, assuming $Ax=\lambda x$:
$$<Ax,x>
=<\lambda x,x>
=\lambda<x,x>
\geq0$$

$(4)\implies(5)$ This follows indeed by using Theorem 2.6.

\medskip

$(5)\implies(1)$ Assuming $A=UDU^*$ with $U\in U_N$, and with $D\in M_N(\mathbb R_+)$ diagonal, we can set $B=U\sqrt{D}U^*$. Then $B$ is self-adjoint, and we have $B^2=A$, as desired.
\end{proof}

Let us record as well the following technical version of the above result:

\begin{theorem}
For a matrix $A\in M_N(\mathbb C)$ the following conditions are equivalent, and if they are satisfied, we say that $A$ is strictly positive, and write $A>0$:
\begin{enumerate}
\item $A=B^2$, with $B=B^*$, invertible.

\item $A=CC^*$, for some $C\in M_N(\mathbb C)$ invertible.

\item $<Ax,x>>0$, for any nonzero vector $x\in\mathbb C^N$.

\item $A=A^*$, and the eigenvalues are strictly positive, $\lambda_i>0$.

\item $A=UDU^*$, with $U\in U_N$ and with $D\in M_N(\mathbb R_+^*)$ diagonal.
\end{enumerate}
\end{theorem}

\begin{proof}
This follows either from Theorem 2.8, by adding the extra assumptions in the statement, or from the proof of Theorem 2.8, by modifying where needed.
\end{proof}

\section*{2b. Matrix laws}

What is next? Looking for non-self-adjoint matrices that we can diagonalize, and here some natural candidates are the unitaries, as shown by the following key formula:
$$\begin{pmatrix}\cos t&-\sin t\\ \sin t&\cos t\end{pmatrix}
=\frac{1}{2}\begin{pmatrix}1&1\\i&-i\end{pmatrix}
\begin{pmatrix}e^{-it}&0\\0&e^{it}\end{pmatrix}
\begin{pmatrix}1&-i\\1&i\end{pmatrix}$$

So, let us have a look now into this, unitaries. Here is our general result:

\begin{theorem}
Any matrix $U\in M_N(\mathbb C)$ which is unitary, $U^*=U^{-1}$, is diagonalizable, with the eigenvalues on the unit circle $\mathbb T$. More precisely we have
$$U=VDV^*$$
with $V\in U_N$, and with $D\in M_N(\mathbb T)$ diagonal. The converse holds too.
\end{theorem}

\begin{proof}
As a first remark, the converse trivially holds, because given a matrix of type $U=VDV^*$, with $V\in U_N$, and with $D\in M_N(\mathbb T)$ being diagonal, we have:
$$U^*
=VD^*V^*
=(V^*)^{-1}D^{-1}V^{-1}
=U^{-1}$$

Let us prove now the first assertion, stating that the eigenvalues of a unitary matrix $U\in U_N$ belong to $\mathbb T$. Indeed, assuming $Uv=\lambda v$, we have:
\begin{eqnarray*}
<v,v>
&=&<U^*Uv,v>\\
&=&<Uv,Uv>\\
&=&<\lambda v,\lambda v>\\
&=&|\lambda|^2<v,v>
\end{eqnarray*}

Thus we obtain $\lambda\in\mathbb T$, as claimed. Our next claim now is that the eigenspaces corresponding to different eigenvalues are pairwise orthogonal. Assume indeed that:
$$Uv=\lambda v\quad,\quad 
Uw=\mu w$$

We have then the following computation, using $U^*=U^{-1}$ and $\lambda,\mu\in\mathbb T$:
\begin{eqnarray*}
\lambda<v,w>
&=&<\lambda v,w>\\
&=&<Uv,w>\\
&=&<v,U^*w>\\
&=&<v,U^{-1}w>\\
&=&<v,\mu^{-1}w>\\
&=&\mu<v,w>
\end{eqnarray*}

Thus $\lambda\neq\mu$ implies $v\perp w$, as claimed. In order now to finish the proof, it remains to prove that the eigenspaces of $U$ span the whole space $\mathbb C^N$. For this purpose, we will use a recurrence method. Let us pick an eigenvector of our matrix:
$$Uv=\lambda v$$

Assuming that we have a vector $w$ orthogonal to it, $v\perp w$, we have:
\begin{eqnarray*}
<Uw,v>
&=&<w,U^*v>\\
&=&<w,U^{-1}v>\\
&=&<w,\lambda^{-1}v>\\
&=&\lambda<w,v>\\
&=&0
\end{eqnarray*}

Thus, if $v$ is an eigenvector, then the vector space $v^\perp$ is invariant under $U$. Now since $U$ is an isometry, so is its restriction to this space $v^\perp$. Thus this restriction is a unitary, and so we can proceed by recurrence, and we obtain the result.
\end{proof}

Well done all this, but we have now to unify Theorem 2.6, dealing with the self-adjoints, with Theorem 2.10, dealing with the unitaries. And fortunately, we have here:

\begin{theorem}
Any matrix $A\in M_N(\mathbb C)$ which is normal, $AA^*=A^*A$, is diagonalizable, with the diagonalization being of the following type,
$$A=UDU^*$$
with $U\in U_N$, and with $D\in M_N(\mathbb C)$ diagonal. The converse holds too.
\end{theorem}

\begin{proof}
As a first remark, the converse trivially holds, because if we take a matrix of the form $A=UDU^*$, with $U$ unitary and $D$ diagonal, then we have:
$$AA^*
=UDD^*U^*
=UD^*DU^*
=A^*A$$

In the other sense now, this is something more technical. Our first claim is that a matrix $A$ is normal precisely when the following happens, for any vector $v$:
$$||Av||=||A^*v||$$

Indeed, the above equality can be written as follows:
$$<AA^*v,v>=<A^*Av,v>$$

But this is equivalent to $AA^*=A^*A$, by expanding the scalar products. Our claim now is that $A,A^*$ have the same eigenvectors, with conjugate eigenvalues:
$$Av=\lambda v\implies A^*v=\bar{\lambda}v$$

Indeed, this follows from the following computation, and from the trivial fact that if $A$ is normal, then so is any matrix of type $A-\lambda 1_N$:
$$||(A^*-\bar{\lambda}1_N)v||
=||(A-\lambda 1_N)^*v||
=||(A-\lambda 1_N)v||
=0$$

Let us prove now, by using this, that the eigenspaces of $A$ are pairwise orthogonal. Assume that we have two eigenvectors, corresponding to different eigenvalues, $\lambda\neq\mu$:
$$Av=\lambda v\quad,\quad 
Aw=\mu w$$

We have the following computation, which shows that $\lambda\neq\mu$ implies $v\perp w$:
\begin{eqnarray*}
\lambda<v,w>
&=&<\lambda v,w>\\
&=&<Av,w>\\
&=&<v,A^*w>\\
&=&<v,\bar{\mu}w>\\
&=&\mu<v,w>
\end{eqnarray*}

In order to finish, it remains to prove that the eigenspaces of $A$ span the whole $\mathbb C^N$. This is something that we have already seen for the self-adjoint matrices, and for unitaries, and we will use here these results, in order to deal with the general normal case. As a first observation, given an arbitrary matrix $A$, the matrix $AA^*$ is self-adjoint:
$$(AA^*)^*=AA^*$$

Thus, we can diagonalize this matrix $AA^*$, as follows, with the passage matrix being a unitary, $V\in U_N$, and with the diagonal form being real, $E\in M_N(\mathbb R)$:
$$AA^*=VEV^*$$

Now observe that, for matrices of type $A=UDU^*$, which are those that we supposed to deal with, we have the following formulae:
$$V=U\quad,\quad 
E=D\bar{D}$$

In particular, the matrices $A$ and $AA^*$ have the same eigenspaces. So, this will be our idea, proving that the eigenspaces of $AA^*$ are eigenspaces of $A$. In order to do so, let us pick two eigenvectors $v,w$ of the matrix $AA^*$, corresponding to different eigenvalues, $\lambda\neq\mu$. The eigenvalue equations are then as follows:
$$AA^*v=\lambda v\quad,\quad 
AA^*w=\mu w$$

We have the following computation, using the normality condition $AA^*=A^*A$, and the fact that the eigenvalues of $AA^*$, and in particular $\mu$, are real:
\begin{eqnarray*}
\lambda<Av,w>
&=&<\lambda Av,w>\\
&=&<A\lambda v,w>\\
&=&<AAA^*v,w>\\
&=&<AA^*Av,w>\\
&=&<Av,AA^*w>\\
&=&<Av,\mu w>\\
&=&\mu<Av,w>
\end{eqnarray*}

We conclude that we have $<Av,w>=0$. But this reformulates as follows:
$$\lambda\neq\mu\implies A(E_\lambda)\perp E_\mu$$

Now since the eigenspaces of $AA^*$ are pairwise orthogonal, and span the whole $\mathbb C^N$, we deduce from this that these eigenspaces are invariant under $A$:
$$A(E_\lambda)\subset E_\lambda$$

But with this result in hand, we can finish the proof of the theorem. Indeed, we can decompose the problem, and the matrix $A$ itself, following these eigenspaces of $AA^*$, which in practice amounts in saying that we can assume that we only have 1 eigenspace. By rescaling, this is the same as assuming that we have $AA^*=1$, and so we are now into the unitary case, that we know how to solve, as explained in Theorem 2.10.
\end{proof}

Getting back now to our questions regarding moments and laws, we have:

\index{normal matrix}

\begin{theorem}
For a normal matrix, $AA^*=A^*A$, we have the formula
$$E(P(A))=\int_\mathbb CP(z)d\mu_A(z)$$
for any polynomial $P\in\mathbb C<X,X^*>$, with $\mu_A$ being the complex discrete measure
$$\mu_A=\frac{1}{N}\sum_i\delta_{\lambda_i}$$
which can be therefore thought of as being the law of $A$.
\end{theorem}

\begin{proof}
If we write our matrix as $A=UDU^*$, as in Theorem 2.11, we have:
\begin{eqnarray*}
E(P(A))
&=&E(UP(D)U^*)\\
&=&E(U^*UP(D))\\
&=&E(P(D))
\end{eqnarray*}

But with this formula in hand, our result follows from our previous computation for the diagonal matrices, which was something trivial, made after Question 2.5.
\end{proof}

With this discussed, what is next? A no-go theorem, closing the discussion:

\begin{theorem}
With the law of a complex matrix $A\in M_N(\mathbb C)$ being the following functional, on the algebra of polynomials in two noncommuting variables $X,X^*$,
$$\mu_A:\mathbb C<X,X^*>\to\mathbb C\quad,\quad P\to E(P(A))$$
this law is a usual complex probability measure $\mu_A\in\mathcal P(\mathbb C)$, in the sense that
$$E(P(A))=\int_\mathbb CP(z)\,d\mu_A(z)$$
precisely when $A$ is normal. In this case $\mu_A=\frac{1}{N}\sum_i\delta_{\lambda_i}$, with $\lambda_i\in\mathbb C$ being the eigenvalues.
\end{theorem}

\begin{proof}
This comes from Theorem 2.12, with only the word ``precisely'' being new. So, let us prove that a non-normal matrix cannot have a law $\mu_A\in\mathcal P(\mathbb C)$.

\medskip

(1) Assuming $AA^*\neq A^*A$, in order to show that $\mu_A$ is not a measure, we can use:
\begin{eqnarray*}
AA^*-A^*A\neq0
&\implies&(AA^*-A^*A)^2>0\\
&\implies&AA^*AA^*-AA^*A^*A-A^*AAA^*+A^*AA^*A>0\\
&\implies&E(AA^*AA^*-AA^*A^*A-A^*AAA^*+A^*AA^*A)>0\\
&\implies&E(AA^*AA^*+A^*AA^*A)>E(AA^*A^*A+A^*AAA^*)\\
&\implies&E(AA^*AA^*)>E(AAA^*A^*)
\end{eqnarray*}

(2) Indeed, assume by contradiction that the law of $A$ is a usual probability measure $\mu_A\in\mathcal P(\mathbb C)$. In particular, we should have the following formulae:
$$E(AA^*AA^*)=\int_\mathbb Cz\bar{z}z\bar{z}\,d\mu_A(z)
\quad,\quad E(AAA^*A^*)=\int_\mathbb Czz\bar{z}\bar{z}\,d\mu_A(z)$$

But the integrands in both formulae are equal to $|z|^2$, so the two integrals are equal, and this is in contradiction with the inequality found in (1), as desired.
\end{proof}

Now with the axiomatic discussion closed, time perhaps to meditate a bit, on what we found. Looking at Theorem 2.13, the ``good'' case, where the matrix is normal, is that where the algebra $<A>\subset M_N(\mathbb C)$ is commutative. And the ``bad'' case, where the matrix is non-normal, is that where the algebra $<A>\subset M_N(\mathbb C)$ is not commutative. So, the dichotomy in Theorem 2.13 is something of algebraic nature, and we have:

\begin{comment}
There are two types of matrices $A\in M_N(\mathbb C)$:
\begin{enumerate}
\item The normal ones, for which $<A>\subset M_N(\mathbb C)$ is commutative. For these we have $<A>\simeq L^\infty(X)$, with $X=\{\lambda_1,\ldots,\lambda_N\}$ being the eigenvalues, and the relevant random variable is $id:X\to\mathbb C$, having law  $\mu=\frac{1}{N}\sum_i\delta_{\lambda_i}$.

\item The non-normal ones, for which $<A>\subset M_N(\mathbb C)$ is not commutative. For these we cannot have an interpretation as classical variables, and we must define the law as being the functional $\mu:\mathbb C<X,X^*>\to\mathbb C$, $\mu(P)=E(P(A))$.
\end{enumerate}
\end{comment}

And we will leave some further thinking at this as an instructive exercise. And with the comment that we will be regularly doing such yoga type sessions in this book, in more and more general frameworks, until reaching to peace and understanding.

\bigskip

Excuse me, but the cat, who has been monitoring my work since we started this chapter on free probability, is here, meowing something. What is it, cat?

\begin{cat}
You're fooling the reader, what you talked about so far are variables of type $id:X\to X$, with $X\subset\mathbb C$ finite. There is no free probability in this.
\end{cat}

Okay, thanks cat, good point indeed, but as just mentioned, this was a matter of doing some yoga, not all of us humans are born cats, takes some time to become one.

\bigskip

So, following now what cat says, namely burying the normal matrices, on the grounds that these are trivial objects, and getting into the unknown, the non-normal matrices, as a first thing that can be said here, which is quite interesting, we have:

\index{Jordan block}

\begin{theorem}
For the simplest non-normal matrix, namely the Jordan block 
$$J=\begin{pmatrix}0&1\\0&0\end{pmatrix}$$
the corresponding law is given by the following moment formula, and linearity
$$M_k(J)=\begin{cases}
\frac{1}{2}&{\rm if}\ k\ {\rm  is\ alternating,\ of\ even\ length}\\
0&{\rm otherwise}\end{cases}$$
and with this obviously not coming from a complex probability measure.
\end{theorem}

\begin{proof}
Skipping some discussion about Jordan form and Jordan blocks, that we will technically not need at this stage, let us just look for the simplest possible non-normal matrix. And here, the answer is a no-brainer, with the culprits being as follows:
$$J=\begin{pmatrix}0&1\\0&0\end{pmatrix}
\quad,\quad J^*=\begin{pmatrix}0&0\\1&0\end{pmatrix}$$

Now observe that the length 2 products between $J$ and $J^*$ are as follows:
$$JJ=J^*J^*=\begin{pmatrix}0&0\\0&0\end{pmatrix}
\quad,\quad JJ^*=\begin{pmatrix}1&0\\0&0\end{pmatrix}
\quad,\quad J^*J=\begin{pmatrix}0&0\\0&1\end{pmatrix}$$

Next, observe that, due to $JJ=J^*J^*=0$, only the alternating products between $J,J^*$ will be nonzero. Thus, we only have to compute two length 3 products, namely:
$$JJ^*J=\begin{pmatrix}1&0\\0&0\end{pmatrix}\begin{pmatrix}0&1\\0&0\end{pmatrix}=\begin{pmatrix}0&1\\0&0\end{pmatrix}=J$$
$$J^*JJ^*=\begin{pmatrix}0&0\\0&1\end{pmatrix}\begin{pmatrix}0&0\\1&0\end{pmatrix}=\begin{pmatrix}0&0\\1&0\end{pmatrix}=J^*$$

But with this done, because at higher length, since $JJ^*$ and $J^*J$ are both idempotents, the quantity of these idempotents will not matter. Thus, we have proved: 
$$JJ^*\ldots JJ^*=\begin{pmatrix}1&0\\0&0\end{pmatrix}\quad,\quad 
JJ^*\ldots JJ^*J=\begin{pmatrix}0&1\\0&0\end{pmatrix}$$
$$J^*J\ldots J^*J=\begin{pmatrix}0&0\\0&1\end{pmatrix}\quad,\quad 
J^*J\ldots J^*JJ^*=\begin{pmatrix}0&0\\1&0\end{pmatrix}$$

Now by applying $E$, we are led to the conclusion in the statement.
\end{proof}

For more complicated $2\times2$ matrices, such as Jordan blocks with eigenvalues, or arbitrary upper triangular matrices, things get worse. As for higher dimensions, here things are certainly quite complicated. With exercise for you, to explore a bit all this.

\bigskip

In what follows we will not further go into this direction, which is itself something very interesting, and this in order not to get hit by a further demolishing comment of cat. The problem indeed with the non-normal matrices $A\in M_N(\mathbb C)$ is that, while the algebras $<A>\subset M_N(\mathbb C)$ are not commutative, and with this being very good news, these algebras remain finite dimensional, and this prevents freeness to happen.

\section*{2c. Infinite matrices}

With the above discussed, let us turn now to the case of infinite matrices, which will, hopefully, lead to free probability. The infinite matrices appear as $A\in M_\infty(\mathbb C)$, or equivalently as linear maps $T:\mathbb C^\infty\to\mathbb C^\infty$, with of course some due attention to convergence matters. However, for various reasons, in order to get started it is technically convenient not to use a basis for the ambient space, and to rely instead on:

\index{Hilbert space}
\index{scalar product}

\begin{definition}
A complex Hilbert space is a complex vector space $H$ with a scalar product $<x,y>$, taken linear at left and antilinear at right,
$$<\lambda x,y>=\lambda<x,y>\quad,\quad <x,\lambda y>=\bar{\lambda}<x,y>$$
which is complete with respect to corresponding norm
$$||x||=\sqrt{<x,x>}$$
in the sense that any sequence $\{x_n\}$ which is a Cauchy sequence, having the property $||x_n-x_m||\to0$ with $n,m\to\infty$, has a limit, $x_n\to x$.
\end{definition}

Here the fact that $||x||=\sqrt{<x,x>}$ is indeed a norm, satisfying $||x+y||\leq||x||+||y||$, follows from Cauchy-Schwarz, $|\!<x,y>\!|\leq||x||\cdot||y||$, which itself follows as in the $H=\mathbb C^N$ case, by looking at the discriminant of the following degree 2 polynomial:
$$f(t)=||wx+ty||^2$$

As for our convention for scalar products, written $<x,y>$ and taken linear at left, this is the so-called feline convention, which can be a source of many good things. 

\bigskip

Getting now to theorems about such spaces, we first have something standard:

\begin{theorem}
The following are Hilbert spaces:
\begin{enumerate}
\item The space $H=\mathbb C^N$, with scalar product $<x,y>=\sum_ix_i\bar{y}_i$.

\item More generally, $H=l^2(I)$, with scalar product $<x,y>=\sum_ix_i\bar{y}_i$.

\item Even more generally, $H=L^2(X)$, with $<f,g>=\int_Xf(x)\overline{g(x)}dx$.
\end{enumerate}
\end{theorem}

\begin{proof}
Here all the assertions are very standard, with (1) being clear, then (2) being an easy remake of (1), and then (3) being an easy remake of (2). 
\end{proof}

Less trivial now is the following result, getting us into advanced mathematics:

\index{separable space}
\index{Gram-Schmidt}

\begin{theorem}
There is only one separable Hilbert space, that is, a space having a countable orthonormal basis, and this unique space can be, upon needs and taste:
\begin{enumerate}
\item $H=l^2(\mathbb N)$. Alternatively, $H=l^2(\mathbb Z)$.

\item $H=L^2[0,1]$. Alternatively, $H=L^2[-1,1]$, or $H=L^2[0,\infty)$, or $H=L^2(\mathbb R)$.

\item $H=L^2([0,1]^n)$. Or $H=L^2([-1,1]^n)$, $H=L^2([0,\infty)^n)$, $H=L^2(\mathbb R^n)$.

\item $H=L^2(\mathbb T)$. Or more generally, $H=L^2(\mathbb T^n)$, or $H=L^2(S^{n-1}_\mathbb R)$.
\end{enumerate}
In fact, what happens is that $L^2(X)$ is separable, provided that $X$ itself is separable. 
\end{theorem}

\begin{proof}
The first assertion is clear, based on the fact that any Hilbert space has an orthonormal basis, which itself comes from Gram-Schmidt, with some care in relation with infinite dimensionality, when doing the recurrence. Indeed, $H$ being separable means to have a countable orthonormal basis $\{e_i\}_{i\in I}$, so $H\simeq l^2(I)\simeq l^2(\mathbb N)$. As for the rest:

\medskip

(1) This is trivial, coming from what we just found, namely $H\simeq l^2(I)$, with $I=\mathbb N$ being the simplest countable set, and with $I=\mathbb Z$ being the next simplest one.

\medskip

(2) This comes from the Weierstrass approximation theorem, stating that $\{x^n\}_{n\in\mathbb N}$ is an algebraic basis of $L^2[0,1]$, that we can orthogonalize afterwards via Gram-Schmidt. Next, the passage to $H=L^2[-1,1]$ or $H=L^2[0,\infty)$ or $H=L^2(\mathbb R)$ is standard. With the remark however that, if you want explicit orthogonal bases for these spaces, you are led into quite subtle mathematics, involving lots of interesting orthogonal polynomials, due to Legendre, Chebycheff, Gegenbauer, Jacobi, Laguerre and Hermite.

\medskip

(3) This is a routine generalization of (2), but again coming with the remark that, if you want explicit bases for these spaces, you are led into orthogonal polynomials.

\medskip

(4) This is something of different nature, coming from the Fourier approximation theorem, stating that $\{z^n\}_{n\in\mathbb Z}$ is an orthonormal basis of $L^2(\mathbb T)$. Then, the passage to higher dimensions ramifies, leading to both $H=L^2(\mathbb T^n)$ and  $H=L^2(S^{n-1}_\mathbb R)$, and this because $\mathbb T\subset\mathbb R^2$ is at the same time the unit torus, and the unit sphere. 

\medskip

(5) Finally, the last assertion comes from measure theory, the idea being that we can short-circuit all the above interesting mathematics, with an abstract argument. 
\end{proof}

Quite interesting the above result, and you get the point I guess, at the beginner level, which is ours for the moment, in order to avoid advanced mathematics, we must assume that our Hilbert spaces come as in Definition 2.17, as abstract spaces, without bases.

\bigskip

Moving ahead, now that we know what our vector spaces are, we can talk about matrices with respect to them. And the situation here is as follows:

\index{infinite matrix}
\index{linear operator}
\index{bounded operator}

\begin{theorem}
Given a Hilbert space $H$, consider the linear operators $T:H\to H$, and for each such operator define its norm by the following formula:
$$||T||=\sup_{||x||=1}||Tx||$$
The operators which are bounded, $||T||<\infty$, form then a complex algebra $B(H)$, which is complete with respect to $||.||$. When $H$ comes with a basis $\{e_i\}_{i\in I}$, we have
$$B(H)\subset M_I(\mathbb C)$$
with the correspondence $T\to M$ coming via the usual linear algebra formulae, namely:
$$T(x)=Mx\quad,\quad M_{ij}=<Te_j,e_i>$$
In infinite dimensions, the inclusion $B(H)\subset M_I(\mathbb C)$ is not an equality.
\end{theorem}

\begin{proof}
This is something straightforward, the idea being as follows:

\medskip

(1) The fact that we have indeed an algebra, satisfying the product condition in the statement, follows from the following estimates, which are all elementary:
$$||S+T||\leq||S||+||T||\quad,\quad 
||\lambda T||=|\lambda|\cdot||T||\quad,\quad 
||ST||\leq||S||\cdot||T||$$

(2) Regarding now the completness assertion, if $\{T_n\}\subset B(H)$ is Cauchy then $\{T_nx\}$ is Cauchy for any $x\in H$, so we can define the limit $T=\lim_{n\to\infty}T_n$ by setting:
$$Tx=\lim_{n\to\infty}T_nx$$

(3) It is then clear that $x\to Tx$ is linear, and it remains to prove that we have $T\in B(H)$, and that $T_n\to T$ in norm. For this purpose, observe that we have:
\begin{eqnarray*}
||T_n-T_m||\leq\varepsilon\ ,\ \forall n,m\geq N
&\implies&||T_nx-T_mx||\leq\varepsilon\ ,\ \forall||x||=1\ ,\ \forall n,m\geq N\\
&\implies&||T_nx-Tx||\leq\varepsilon\ ,\ \forall||x||=1\ ,\ \forall n\geq N\\
&\implies&||T_Nx-Tx||\leq\varepsilon\ ,\ \forall||x||=1\\
&\implies&||T_N-T||\leq\varepsilon
\end{eqnarray*}

But this gives both $T\in B(H)$, and $T_N\to T$ in norm, and we are done.

\medskip

(4) Finally, the correspondence $T\to M$ in the statement is linear, and its kernel is $\{0\}$, so we have an embedding $B(H)\subset M_I(\mathbb C)$. In finite dimensions this is an isomorphism. However, in infinite dimensions we have matrices not producing operators, as for instance the all-one matrix, so this embedding $B(H)\subset M_I(\mathbb C)$ is not an isomorphism.
\end{proof}

As a second and last basic result regarding the operators, we will need:

\index{adjoint operator}
\index{adjoint matrix}

\begin{theorem}
Each operator $T\in B(H)$ has an adjoint $T^*\in B(H)$, given by: 
$$<Tx,y>=<x,T^*y>$$
The operation $T\to T^*$ is antilinear, antimultiplicative, involutive, and satisfies:
$$||T||=||T^*||\quad,\quad ||TT^*||=||T||^2$$
When $H$ comes with a basis $\{e_i\}_{i\in I}$, the operation $T\to T^*$ corresponds to
$$(M^*)_{ij}=\overline{M}_{ji}$$ 
at the level of the associated matrices $M\in M_I(\mathbb C)$.
\end{theorem}

\begin{proof}
This is standard too, and can be proved in 3 steps, as follows:

\medskip

(1) The existence of the adjoint operator $T^*$, given by the formula in the statement, comes from the fact that the function $\varphi(x)=<Tx,y>$ being a linear map $H\to\mathbb C$, we must have a formula as follows, for a certain vector $T^*y\in H$:
$$\varphi(x)=<x,T^*y>$$

Moreover, since this vector is unique, $T^*$ is unique too, and we have as well:
$$(S+T)^*=S^*+T^*\quad,\quad
(\lambda T)^*=\bar{\lambda}T^*\quad,\quad 
(ST)^*=T^*S^*\quad,\quad 
(T^*)^*=T$$

Observe that we have indeed $T^*\in B(H)$, as shown by the following computation:
\begin{eqnarray*}
||T||
&=&\sup_{||x||=1}\sup_{||y||=1}<Tx,y>\\
&=&\sup_{||y||=1}\sup_{||x||=1}<x,T^*y>\\
&=&||T^*||
\end{eqnarray*}

(2) Next, $||TT^*||\leq||T||\cdot||T^*||=||T||^2$. In the other sense now, observe that:
\begin{eqnarray*}
||T||^2
&=&\sup_{||x||=1}|<Tx,Tx>|\\
&=&\sup_{||x||=1}|<x,T^*Tx>|\\
&\leq&||T^*T||
\end{eqnarray*}

But by replacing $T\to T^*$ we obtain from this $||T||^2\leq||TT^*||$, as desired.

\medskip

(3) Finally, when $H$ comes with a basis, the formula $<Tx,y>=<x,T^*y>$ applied with $x=e_i$, $y=e_j$ translates into the formula $(M^*)_{ij}=\overline{M}_{ji}$, as desired.
\end{proof}

Next, we need to talk about operator algebras. Which can be something quite tricky, but fortunately, we have here the following clever definition, due to Gelfand:

\begin{definition}
An abstract operator algebra, or $C^*$-algebra, is a complex algebra $A$ having a norm $||.||$ and an involution $*$, subject to the following conditions:
\begin{enumerate}
\item $A$ is closed with respect to the norm.

\item We have $||aa^*||=||a||^2$, for any $a\in A$.
\end{enumerate}
\end{definition}

In other words, what we did here is to axiomatize the abstract properties of the operator algebras $A\subset B(H)$, chosen closed under the norm and under the adjoint operation, both very natural conditions, without reference to the ambient Hilbert space $H$.

\bigskip

As a first observation, which is of particular importance for the considerations to follow, the random matrix algebras are covered by our formalism, as shown by:

\index{random matrix algebra}

\begin{theorem}
The random matrix algebras $A=M_N(L^\infty(X))$ are $C^*$-algebras, with their usual norm and involution, given by:
$$||Z||=\sup_{x\in X}||Z_x||\quad,\quad 
(Z^*)_{ij}=\overline{Z}_{ij}$$
These algebras generalize both the algebras $M_N(\mathbb C)$, and the algebras $L^\infty(X)$.
\end{theorem}

\begin{proof}
The fact that the $C^*$-algebra axioms are satisfied is clear from definitions. As for the last assertion, this follows by taking $X=\{.\}$ and $N=1$, respectively.
\end{proof}

We can in fact say more about the random matrix algebras, as follows:

\begin{theorem}
Any random variable algebra is an operator algebra, as follows:
$$L^\infty(X)\subset B(L^2(X))\quad,\quad 
f\to(g\to fg)$$
More generally, any random matrix algebra is an operator algebra, as follows,
$$M_N(L^\infty(X))\subset B\left(\mathbb C^N\otimes L^2(X)\right)$$
with the embedding being the above one, tensored with the identity.
\end{theorem}

\begin{proof}
We have two assertions to be proved, the idea being as follows:

\medskip

(1) Given $f\in L^\infty(X)$, consider the following operator, acting on $H=L^2(X)$:
$$T_f(g)=fg$$

Observe that $T_f$ is indeed well-defined, and bounded as well, because:
$$||fg||_2
=\sqrt{\int_X|f(x)|^2|g(x)|^2d\mu(x)}
\leq||f||_\infty||g||_2$$

The application $f\to T_f$ being linear, involutive, continuous, and injective as well, we obtain in this way a $C^*$-algebra embedding $L^\infty(X)\subset B(H)$, as desired.

\medskip

(2) Regarding the second assertion, this is best viewed in the following way:
\begin{eqnarray*}
M_N(L^\infty(X))
&=&M_N(\mathbb C)\otimes L^\infty(X)\\
&\subset&M_N(\mathbb C)\otimes B(L^2(X))\\
&=&B\left(\mathbb C^N\otimes L^2(X)\right)
\end{eqnarray*}

Here we have used (1), and some standard tensor product identifications.
\end{proof}

Our purpose in what follows will be to develop the spectral theory of the $C^*$-algebras, and in particular that of the random matrix algebras $M_N(L^\infty(X))$ that we are mainly interested in, our main objective being that of talking about spectral measures, in the normal case, in analogy with what we know about the usual matrices. We first have:

\index{spectrum}
\index{functional calculus}
\index{rational function}
\index{self-adjoint}
\index{unitary}
\index{spectral radius}
\index{normal element}

\begin{theorem}
Given an element $a\in A$ of a $C^*$-algebra, define its spectrum as:
$$\sigma(a)=\left\{\lambda\in\mathbb C\Big|a-\lambda\notin A^{-1}\right\}$$
The following spectral theory results hold, exactly as in the $A=B(H)$ case:
\begin{enumerate}
\item We have $\sigma(ab)\cup\{0\}=\sigma(ba)\cup\{0\}$.

\item We have $\sigma(f(a))=f(\sigma(a))$, for any $f\in\mathbb C(X)$ having poles outside $\sigma(a)$.

\item The spectrum $\sigma(a)$ is compact, non-empty, and contained in $D_0(||a||)$.

\item The spectra of unitaries $(u^*=u^{-1})$ and self-adjoints $(a=a^*)$ are on $\mathbb T,\mathbb R$.

\item The spectral radius of normal elements $(aa^*=a^*a)$ is given by $\rho(a)=||a||$.
\end{enumerate}
In addition, assuming $a\in A\subset B$, the spectra of $a$ with respect to $A$ and to $B$ coincide.
\end{theorem}

\begin{proof}
Here the assertions (1-5), which are of course formulated a bit informally, are well-known for the full operator algebra $A=B(H)$, and the proof in general is similar:

\medskip

(1) Assuming that $1-ab$ is invertible, with inverse $c$, we have $abc=cab=c-1$, and it follows that $1-ba$ is invertible too, with inverse $1+bca$. Thus $\sigma(ab),\sigma(ba)$ agree on $1\in\mathbb C$, and by linearity, it follows that $\sigma(ab),\sigma(ba)$ agree on any point $\lambda\in\mathbb C^*$.

\medskip

(2) The formula $\sigma(f(a))=f(\sigma(a))$ is clear for polynomials, $f\in\mathbb C[X]$, by factorizing $f-\lambda$, with $\lambda\in\mathbb C$. Then, the extension to the rational functions is straightforward, because $P(a)/Q(a)-\lambda$ is invertible precisely when $P(a)-\lambda Q(a)$ is.

\medskip

(3) By using $1/(1-b)=1+b+b^2+\ldots$ for $||b||<1$ we obtain that $a-\lambda$ is invertible for $|\lambda|>||a||$, and so $\sigma(a)\subset D_0(||a||)$. It is also clear that $\sigma(a)$ is closed, so what we have is a compact set. Finally, assuming $\sigma(a)=\emptyset$ the function $f(\lambda)=\varphi((a-\lambda)^{-1})$ is well-defined, for any $\varphi\in A^*$, and by Liouville we get $f=0$, contradiction.

\medskip

(4) Assuming $u^*=u^{-1}$ we have $||u||=1$, and so $\sigma(u)\subset D_0(1)$. But with $f(z)=z^{-1}$ we obtain via (2) that we have as well $\sigma(u)\subset f(D_0(1))$, and this gives $\sigma(u)\subset\mathbb T$. As for the result regarding the self-adjoints, this can be obtained from the result for the unitaries, by using (2) with functions of type $f(z)=(z+it)/(z-it)$, with $t\in\mathbb R$.

\medskip

(5) It is routine to check, by integrating quantities of type $z^n/(z-a)$ over circles centered at the origin, and estimating, that the spectral radius is given by $\rho(a)=\lim||a^n||^{1/n}$. But in the self-adjoint case, $a=a^*$, this gives $\rho(a)=||a||$, by using exponents of type $n=2^k$, and then the extension to the general normal case is straightforward.

\medskip 

(6) Regarding now the last assertion, the inclusion $\sigma_B(a)\subset\sigma_A(a)$ is clear. For the converse, assume $a-\lambda\in B^{-1}$, and set $b=(a-\lambda )^*(a-\lambda )$. We have then:
$$\sigma_A(b)-\sigma_B(b)=\left\{\mu\in\mathbb C-\sigma_B(b)\Big|(b-\mu)^{-1}\in B-A\right\}$$

Thus this difference in an open subset of $\mathbb C$. On the other hand $b$ being self-adjoint, its two spectra are both real, and so is their difference. Thus the two spectra of $b$ are equal, and in particular $b$ is invertible in $A$, and so $a-\lambda\in A^{-1}$, as desired.
\end{proof}

We can now prove a key result about the operator algebras, as follows:

\index{Gelfand theorem}
\index{commutative algebra}

\begin{theorem}[Gelfand]
If $X$ is a compact space,  the algebra $C(X)$ of continuous functions on it $f:X\to\mathbb C$ is a $C^*$-algebra, with usual norm and involution, namely:
$$||f||=\sup_{x\in X}|f(x)|\quad,\quad 
f^*(x)=\overline{f(x)}$$
Conversely, any commutative $C^*$-algebra is of this form, $A=C(X)$, with 
$$X=\Big\{\chi:A\to\mathbb C\ ,\ {\rm normed\ algebra\ character}\Big\}$$
with topology making continuous the evaluation maps $ev_a:\chi\to\chi(a)$.
\end{theorem}

\begin{proof}
Given a commutative $C^*$-algebra $A$, let us define $X$ as in the statement. Then $X$ is compact, and $a\to ev_a$ is a morphism of algebras, as follows:
$$ev:A\to C(X)$$

We first prove that $ev$ is involutive. We use the following formula, which is similar to the $z=Re(z)+iIm(z)$ decomposition formula for usual complex numbers:
$$a=\frac{a+a^*}{2}+i\cdot\frac{a-a^*}{2i}$$

Thus it is enough to prove $ev_{a^*}=ev_a^*$ for the self-adjoint elements $a$. But this is the same as proving that $a=a^*$ implies that $ev_a$ is a real function, which is in turn true, by Theorem 2.25, because $ev_a(\chi)=\chi(a)$ is an element of $\sigma(a)$, contained in $\mathbb R$. Next, since $A$ is commutative, each element is normal, so $ev$ is isometric, as shown by:
$$||ev_a||
=\rho(a)
=||a||$$

It remains to prove that $ev$ is surjective. But this follows from the Stone-Weierstrass theorem, because $ev(A)$ is a closed subalgebra of $C(X)$, which separates the points.
\end{proof}

As a first remark, the Gelfand theorem puts on firm ground Fact 2.4, because in view of it, we can think of any $C^*$-algebra as being of the form $A=C(X)$, with $X$ being a compact quantum space. Of course, this is just a beginning, and more on this, later.

\bigskip

As a more concrete consequence now of the Gelfand theorem, we have:

\begin{theorem}
For any normal element $a\in A$ we have $<a>=C(\sigma(a))$. In addition, given a function $f\in C(\sigma(a))$, we can apply it to $a$, and we have
$$\sigma(f(a))=f(\sigma(a))$$
which generalizes the previous rational calculus formula, in the normal case.
\end{theorem}

\begin{proof}
Since $a$ is normal, the $C^*$-algebra $<a>$ that is generates is commutative, so if we denote by $X$ the space of the characters $\chi:<a>\to\mathbb C$, we have $<a>=C(X)$. Now since the map $X\to\sigma(a)$ given by evaluation at $a$ is bijective, we obtain $<a>=C(\sigma(a))$. Thus, we are dealing here with usual functions, and this gives all the assertions.
\end{proof}

So long for the basic theory of $C^*$-algebras. This has obviously many similarities with the theory of operator algebras $A\subset B(H)$, chosen closed under the usual operator involution and norm, and in fact, coming as some sort of half-surprise, we have:

\begin{fact}
It can be shown, by reasoning along the lines of the proof of Theorem 2.24, that any $C^*$-algebra must appear as a closed $*$-subalgebra $A\subset B(H)$. 
\end{fact}

However, this is something quite advanced, whose full understanding involves as well the spectral theorem for normal operators, and the related notion of von Neumann algebra, which are both technical subjects. We will not need this right next, our plan for what follows being to keep enjoying Definition 2.22 as such. More on all this later.

\section*{2d. Spectral measures}

In order to get now towards noncommutative probability, we first have to develop the theory of positive elements, and linear forms. We first have the following result:

\index{positive element}

\begin{theorem}
For an element $a\in A$ the following conditions are equivalent, and if they are satisfied, we call $a$ positive, and we write $a\geq0$:
\begin{enumerate}
\item $a=a^*$ and $\sigma(a)\subset[0,\infty)$.

\item $a=b^2$, for some $b\in A$ satisfying $b=b^*$.

\item $a=cc^*$, for some $c\in A$.
\end{enumerate}
\end{theorem}

\begin{proof}
This is something very standard, as follows:

\medskip

$(1)\implies(2)$ Since $a=a^*$, with this being in fact implied by $\sigma(a)\subset[0,\infty)$, the algebra $<a>$ is commutative, and by using Theorem 2.27, we can set $b=\sqrt{a}$.

\medskip

$(2)\implies(3)$ This is trivial, because we can simply set $c=b$. 

\medskip

$(3)\implies(1)$ We can prove $\sigma(cc^*)\subset[0,\infty)$ by contradiction. Indeed, by multiplying $c$ by a suitable element of the algebra $<cc^*>$, we are led to the existence of an element $d\neq0$ satisfying $-dd^*\geq0$. By writing now $d=x+iy$ with $x=x^*,y=y^*$ we have:
$$dd^*+d^*d
=2(x^2+y^2)
\geq0$$

Thus $d^*d\geq0$, which is easily seen to contradict the condition $-dd^*\geq0$.
\end{proof}

Next, and again in analogy with what we know about the usual matrices, we have as well a useful strict version of the above result, as follows:

\begin{theorem}
For an element $a\in A$ the following conditions are equivalent, and if they are satisfied, we call $a$ strictly positive, and we write $a>0$:
\begin{enumerate}
\item $a=a^*$ and $\sigma(a)\subset(0,\infty)$.

\item $a=b^2$, for some $b\in A$ invertible satisfying $b=b^*$.

\item $a=cc^*$, for some $c\in A$ invertible.
\end{enumerate}
\end{theorem}

\begin{proof}
This follows either from Theorem 2.29, by adding the extra assumptions in the statement, or from the proof of Theorem 2.29, by modifying where needed.
\end{proof}

As a non-trivial comment now, coming in relation with Fact 2.28, we have:

\begin{comment}
For a true operator, $T\in B(H)$, the following happen:
\begin{enumerate}
\item $T\geq0$ precisely when $<Tx,x>\geq0$, for any $x\in H$.

\item $T>0$ precisely when $T$ is invertible, and $T\geq0$.

\item $T>0$ implies $<Tx,x>>0$ for $x\neq0$, but the converse fails.
\end{enumerate}
\end{comment}

To be more precise here, (1) is something standard, coming a bit as for the usual matrices, and (2) is something trivial, coming from definitions. The subtlety is at the level of (3), with $\implies$ coming from $<Tx,x>=||\sqrt{T}x||^2>0$, and with the failure of the converse coming from $T=diag(1,\frac{1}{2},\frac{1}{3},\ldots)$ on $H=l^2(\mathbb N)$, which is not invertible.

\bigskip

Next, we can talk as well about positive linear forms, as follows:

\index{positive linear form}

\begin{definition}
Consider a continuous linear map $\varphi:A\to\mathbb C$.
\begin{enumerate}
\item $\varphi$ is called positive when $a\geq0\implies\varphi(a)\geq0$.

\item $\varphi$ is called faithful and positive when $a>0\implies\varphi(a)>0$.
\end{enumerate}
\end{definition}

In the commutative case, $A=C(X)$, the positive linear forms appear as follows, with $\mu$ being positive, and strictly positive if we want $\varphi$ to be faithful and positive:
$$\varphi(f)=\int_Xf(x)d\mu(x)$$

In general, the positive linear forms can be thought of as being integration functionals with respect to some underlying positive measures. As a first application, we have:

\begin{fact}
In order to prove the Gelfand-Naimark-Segal embedding theorem,
$$A\subset B(H)$$
what you need is a faithful positive form $\varphi:A\to\mathbb C$, whose existence is quite standard.
\end{fact}

To be more precise here, let us get back to Fact 2.28. As mentioned there, the GNS theorem is a generalization of Theorem 2.24, and looking at what we did in the proof of Theorem 2.24, that was mostly algebra, the only analytic ingredient being the construction of the Hilbert space $H=L^2(A)$, using the expectation $E:A\to\mathbb C$. But this makes it clear that we can prove the GNS theorem in general, provided that we have available a faithful positive form $\varphi:A\to\mathbb C$, and the existence of such a form is quite standard.

\bigskip

But probably enough of this, the plan was that of enjoying Definition 2.22 as such, enjoyable pure algebra, but we ended up with a string of analytic facts and comments coming in parallel to this learning, putting us in a situation similar to schizophrenia. More on such things, ramified way of thinking at all this, later in this book, in chapter 5.

\bigskip

Getting back now to our usual business, probability, we can formulate:

\index{colored moments}

\begin{definition}
Let $A$ be a $C^*$-algebra, given with a positive trace $tr:A\to\mathbb C$.
\begin{enumerate}
\item The elements $a\in A$ are called random variables.

\item The moments of such a variable are the numbers $M_k(a)=tr(a^k)$.

\item The law of such a variable is the functional $\mu_a:P\to tr(P(a))$.
\end{enumerate}
\end{definition}

To be more precise, here the exponent $k=\circ\bullet\bullet\circ\ldots$ is as usual a colored integer, and the powers $a^k$ are defined by the following formulae, and multiplicativity: 
$$a^\emptyset=1\quad,\quad
a^\circ=a\quad,\quad
a^\bullet=a^*$$ 

As for $P$, this is a noncommuting $*$-polynomial in one variable, $P\in\mathbb C<X,X^*>$. At the level of the general theory, we have the following key result:

\index{normal element}
\index{spectral measure}

\begin{theorem}
Assuming that $a\in A$ is normal, $aa^*=a^*a$, the following happen:
\begin{enumerate}
\item $\mu_a$ is a complex probability measure, satisfying $supp(\mu_a)\subset\sigma(a)$.

\item In the self-adjoint case, $a=a^*$, this measure $\mu_a$ is real.

\item Assuming that $tr$ is faithful, we have $supp(\mu_a)=\sigma(a)$.
\end{enumerate}
Also, when $a\in A$ is not normal, $aa^*\neq a^*a$, and $tr$ is faithful, all this fails.
\end{theorem}

\begin{proof}
According to Theorem 2.27, we have an identification as follows: 
$$<a>=C(\sigma(a))$$

Thus the functional $f(a)\to tr(f(a))$ can be regarded as an integration functional on the algebra $C(\sigma(a))$, and by the Riesz theorem, this latter functional must come from a probability measure $\mu_a$ on the spectrum $\sigma(a)$, in the sense that we must have:
$$tr(f(a))=\int_{\sigma(a)}f(x)d\mu_a(x)$$

We are therefore led to the various conclusions in the statement. Finally, the last assertion comes from the positivity computation in the proof of Theorem 2.13.
\end{proof}

As a main application of this, to the case of the random matrices, we have:

\begin{theorem}
Given a random matrix $Z\in M_N(L^\infty(X))$, assumed to be normal, $ZZ^*=Z^*Z$, its law, when restricted to the usual polynomials in two variables,
$$\mu_Z:\mathbb C[X,X^*]\to\mathbb C\quad,\quad 
P\to\frac{1}{N}\int_XTr(P(Z))$$
must come from a probability measure on the spectrum $\sigma(Z)\subset\mathbb C$, as follows:
$$\mu_Z(P)=\int_{\sigma(T)}P(x)d\mu_Z(x)$$
We agree to use the symbol $\mu_Z$ for all these notions.
\end{theorem}

\begin{proof}
This follows indeed from what we know from Theorem 2.35, applied to the normal element $a=Z$, belonging to the $C^*$-algebra $A=M_N(L^\infty(X))$. 
\end{proof}

All this sounds very good, and we have some beginning of theory here, that we will further explore in chapter 3, with some explicit computations of spectral measures, for some basic examples of normal random matrices, due to Wigner and Wishart.

\bigskip

As a last task for this chapter, in analogy with what we previously did for the usual matrices, let us quickly explore now the non-normal case. The most basic example of a non-normal operator is the shift on the space $l^2(\mathbb N)$, and regarding it, we have:

\index{shift}

\begin{theorem}
The moments of the shift $S(e_i)=e_{i+1}$ on $H=l^2(\mathbb N)$, with respect to the positive linear form $E(T)=<Te_0,e_0>$, which is actually non tracial, are
$$M_k(S)=\begin{cases}
1&{\rm if}\ S^k=1\\
0&{\rm otherwise}
\end{cases}$$
and with $S^k=1$ appearing by iterating the rule $S^*S=1$.
\end{theorem}

\begin{proof}
The shift operator and its adjoint are given by the following formulae:
$$S(e_i)=e_{i+1}\quad,\quad S^*(e_i)
=\begin{cases}
e_{i-1}&(i>0)\\
0&(i=0)
\end{cases}$$

Next, the shift is a non-unitary isometry, as shown by the following formulae:
$$S^*S=1\quad,\quad SS^*=Proj(e_0^\perp)$$

Now observe that, when looking at an arbitrary colored power $S^k$, after iterating the rule $S^*S=1$, we are left with a reduced word of the following type:
$$S^k=S\ldots SS^*\ldots S^*$$

But, with respect to $E(T)=<Te_0,e_0>$, the expectation of this word is $0$, if there is at least one $S^*$ present. And then, assuming that there is no $S^*$ present, the expectation of this word remains $0$ too, assuming that at least one $S$ is present. Thus, as a conclusion, the expectation is 0, unless our word is the empty word, $S^k=1$, as stated.
\end{proof}

The above result is quite nice, and we will heavily use it later, when doing free probability, where $S$ plays a key role. However, there is something frustrating about it, coming from the fact that the expectation $E(T)=<Te_0,e_0>$ is not a trace, as it should.

\bigskip

In answer, we can look instead at the group algebras $C^*_{red}(G)\subset B(l^2(G))$, appearing as norm completions of the usual group algebras $\mathbb C[G]$, acting via $g(\delta_h)=\delta_{gh}$, which do have traces, namely $tr(T)=<T\delta_1,\delta_1>$. Here the standard generators $g\in G$ are unitaries, and so normal, but their sums $\sum_ig_i$ are generically not normal, and this can lead to some interesting mathematics. And more on this, later in this book.

\section*{2e. Exercises}

Welcome to noncommutative probability, and as exercises, we have:

\begin{exercise}
Compute unitaries and projections in low dimensions.
\end{exercise}

\begin{exercise}
Learn more about positive matrices, from Sylvester and others.
\end{exercise}

\begin{exercise}
Learn about polar decomposition and singular values.
\end{exercise}

\begin{exercise}
Do some computations, for more complicated Jordan blocks.
\end{exercise}

\begin{exercise}
Learn about Weierstrass, Fourier and orthogonal polynomials.
\end{exercise}

\begin{exercise}
Clarify what we said, regarding spectral theory in $C^*$-algebras.
\end{exercise}

\begin{exercise}
Learn about the GNS theorem, and about finite dimensions too.
\end{exercise}

\begin{exercise}
Compute some laws, in group $C^*$-algebras of your choice.
\end{exercise}

As bonus exercise, and no surprise here, read an operator theory book.

\chapter{Random matrices}

\section*{3a. Gaussian matrices}

We recall that a random matrix algebra is an algebra of type $A=M_N(L^\infty(X))$, and that we are interested in the computation of the laws of the operators $Z\in A$, called random matrices. Regarding the precise classes of random matrices that we are interested in, first we have the complex Gaussian matrices, which are constructed as follows:

\index{Gaussian matrix}

\begin{definition}
A complex Gaussian matrix is a random matrix of type
$$Z\in M_N(L^\infty(X))$$
which has i.i.d. complex normal entries.
\end{definition}

Next, and coming as a somewhat surprising remark, using real normal variables in the above definition, instead of complex ones, leads nowhere. The correct real versions of the Gaussian matrices are the Wigner random matrices, constructed as follows: 

\index{Wigner matrix}

\begin{definition}
A Wigner matrix is a random matrix of type
$$Z\in M_N(L^\infty(X))$$
which has i.i.d. complex normal entries, up to the constraint $Z=Z^*$.
\end{definition}

In other words, a Wigner matrix must be as follows, with the diagonal entries being real normal variables, $a_i\sim g_t$, for some $t>0$, the upper diagonal entries being complex normal variables, $b_{ij}\sim G_t$, the lower diagonal entries being the conjugates of the upper diagonal entries, as indicated, and with all the variables $a_i,b_{ij}$ being independent: 
$$Z=\begin{pmatrix}
a_1&b_{12}&\ldots&\ldots&b_{1N}\\
\bar{b}_{12}&a_2&\ddots&&\vdots\\
\vdots&\ddots&\ddots&\ddots&\vdots\\
\vdots&&\ddots&a_{N-1}&b_{N-1,N}\\
\bar{b}_{1N}&\ldots&\ldots&\bar{b}_{N-1,N}&a_N
\end{pmatrix}$$

As a comment here, for many concrete applications the Wigner matrices are in fact the central objects in random matrix theory, and in particular, they are often more important than the Gaussian matrices. In fact, these are the random matrices which were first considered and investigated, a long time ago, by Wigner himself \cite{wig}.

\bigskip

We will be interested as well in the complex Wishart matrices, which are the positive versions of the above random matrices, constructed as follows: 

\index{Wishart matrix}

\begin{definition}
A complex Wishart matrix is a random matrix of type
$$Z=YY^*\in M_N(L^\infty(X))$$
with $Y$ being a complex Gaussian matrix.
\end{definition}

As before with the Wigner matrices, there are several comments that can be made here. First, using in the above real Gaussian variables instead of complex ones, as Wishart originally did, leads to a less interesting combinatorics. Next, the complex Wishart matrices were systematically studied by Marchenko-Pastur not long after Wigner, and so historically came second \cite{mpa}. Finally, in what regards their combinatorics and applications, these matrices quite often come first, before both the Gaussian and the Wigner ones, with all this being of course a matter of knowledge, needs and taste.

\bigskip

Summarizing, we have three main types of random matrices, which can be somehow designated as ``complex'', ``real'' and ``positive'', and that we will study in what follows. Let us also mention that there are many other interesting classes of random matrices, usually appearing as modifications of the above, notably with some more general types of Wishart matrices, $Z=YY^*$ with $Y$ rectangular Gaussian. More on these later.

\bigskip

As a last thing to be discussed, before getting into computations, what about the precise motivations, for the study of these matrices? In answer, we have:

\begin{fact}
The philosophy of the random matrices $Z\in M_N(L^\infty(X))$ is as follows:
\begin{enumerate}
\item When the usual random variables $f\in L^\infty(X)$ are not enough, for the modeling of your problem, the random matrices with $N\geq2$ are the answer.

\item As we know from chapter 2, each random matrix has a law, and in relation with your initial question, the main problem is that of computing this law.

\item Most questions involve random matrices coming in series, one for each $N\in\mathbb N$, and the main problem is that of computing the law in the $N\to\infty$ limit.
\end{enumerate}
\end{fact}

So, this was for the story and motivations, and lacking of course remains an explicit example of an application. However, and here comes the difficulty, although applications abound at the advanced level, to questions ranging from quantum physics to stock markets, at the beginner level, which is ours for the moment, there is none.

\bigskip

So, we will have to live with this, I mean just trust me, and we will learn about applications later. At a philosophical level, however, in the lack of a concrete illustration for Fact 3.4, let us point out the obvious similarity with what we did in chapter 1 for the reflection groups. Indeed, what we have is, a bit as there, certain beasts depending on $N\in\mathbb N$, constructed out of nothing or almost, using symmetry constraints, and having a law, that we are interested to compute in the $N\to\infty$ limit.

\bigskip

Getting to work now, let us first study the Gaussian matrices. These matrices are not normal, so we will be interested in computing the moments, and depending on what we find, speculate on the corresponding asymptotic law $\mu:\mathbb C<X,X^*>\to\mathbb C$. Now thinking a bit, for the moment computation we will certainly need the Wick formula, that we learned in chapter 1. So, let us first recall that formula:

\begin{theorem}[Wick formula]
Given independent variables $f_i$, each following the complex normal law $G_t$, with $t>0$ being a fixed parameter, we have the formula
$$E\left(f_{i_1}^{k_1}\ldots f_{i_s}^{k_s}\right)=t^{s/2}\#\left\{\pi\in\mathcal P_2(k)\Big|\pi\leq\ker i\right\}$$
with $\mathcal P_2(k)$ standing for the pairings of $k$ which are matching, pairing $\circ-\bullet$ symbols, and where $\pi\leq\ker i$ means that the indices of $i$ must fit into the blocks of $\pi$.
\end{theorem}

\begin{proof}
This is something from chapter 1, the idea being that for single variables the Wick formula is the moment formula for $G_t$, which was as follows:
$$M_k(G_t)=\sum_{\pi\in\mathcal P_2(k)}t^{|\pi|}$$

As for the general case, that of several variables, this follows from this. 
\end{proof}

Next, it is also pretty much clear that for doing our computations, we will need a good knowledge of the partitions. In practice, here is what we need to know: 

\index{cycling inside blocks}

\begin{proposition}
For a permutation $\sigma\in S_k$, let $|\sigma|\in\{1,\ldots,k\}$ be its number of cycles, and $l(\sigma)=k-|\sigma|\in\{0,\ldots,k-1\}$ be its length, that is, the minimal number of transpositions needed for decomposing $\sigma$ as a product of transpositions.
\begin{enumerate}
\item $d(\sigma,\pi)=l(\sigma^{-1}\pi)$ is a distance function on $S_k$.

\item If $\gamma\in S_k$ is a full cycle, $l(\sigma)+l(\sigma^{-1}\gamma)\geq k-1$.

\item Equivalently, we have $|\sigma|+|\sigma^{-1}\gamma|\leq k+1$.

\item Equality happens when $\sigma\in NC(k)\subset S_k$, by cycling inside each block.
\end{enumerate} 
\end{proposition}

\begin{proof}
To start with, given a permutation $\sigma\in S_k$, we can certainly talk about the number of cycles $|\sigma|\in\{1,\ldots,k\}$, and about the length $l(\sigma)\in\{0,\ldots,k-1\}$ too, with the convention $l(id)=0$. In order to prove now $l(\sigma)=k-|\sigma|$, let us decompose $\sigma$ as a product of distinct cycles, having lengths $k_1,\ldots,k_r\geq1$, as follows:
$$\sigma=\gamma_{k_1}\ldots\gamma_{k_r}$$

We can then compute the length of $\sigma$, in the obvious way, and we get, as desired:
\begin{eqnarray*}
l(\sigma)
&=&(k_1-1)+\ldots+(k_r-1)\\
&=&(k_1+\ldots+k_r)-r\\
&=&k-|\sigma|
\end{eqnarray*}

Thus we have the basics, and regarding the proof of the rest, this goes as follows:

\medskip

(1) First, we have $d(\sigma,\sigma)=l(id)=0$. Next, the symmetry axiom comes from:
$$d(\sigma,\pi)=l(\sigma^{-1}\pi)=l\left((\sigma^{-1}\pi)^{-1}\right)=l(\pi^{-1}\sigma)=d(\pi,\sigma)$$

Finally, the triangle inequality $d(\sigma,\pi)+d(\pi,\nu)\geq d(\sigma,\nu)$ amounts in saying that we have $l(\sigma^{-1}\pi)+l(\pi^{-1}\nu)\geq l(\sigma^{-1}\nu)$. But this is clear, coming from $\sigma^{-1}\pi\cdot \pi^{-1}\nu=\sigma^{-1}\nu$.

\medskip

(2) This is the triangle inequality, applied to the triangle $id-\sigma-\gamma$.

\medskip

(3) This comes indeed from (2), via the formula $l(\sigma)=k-|\sigma|$.

\medskip

(4) Consider the embedding $P(k)\subset S_k$, obtained by cycling inside each block, in the obvious way. As a first observation, for $\sigma\in P(k)$ our previous notation $|\sigma|$ for the number of blocks matches with our new notation $|\sigma|$ for the number of cycles. Thus, our formalism is fine, and what happens now is that equality happens in (2,3), meaning that $id-\sigma-\gamma$ must be aligned, precisely when $\sigma\in NC(k)\subset P(k)\subset S_k$. We will leave the details as an exercice, and for this, and more, you can have a look at Biane \cite{bia}.
\end{proof}

With this discussed, here is the moment computation for the Gaussian matrices:

\begin{theorem}
Given a sequence of Gaussian random matrices $Z_N\in M_N(L^\infty(X))$ having independent $G_t$ variables as entries, for some fixed $t>0$, we have
$$M_k\left(\frac{Z_N}{\sqrt{N}}\right)\simeq t^{|k|/2}|\mathcal{NC}_2(k)|$$
in the $N\to\infty$ limit, with $\mathcal{NC}_2\subset\mathcal P_2$ being the noncrossing matching pairings.
\end{theorem}

\begin{proof}
This is something very standard, which can be done as follows:

\medskip

(1) We fix $N\in\mathbb N$, and we let $Z=Z_N$. Let us first compute the trace of $Z^k$. With $k=k_1\ldots k_s$, and with the convention $(ij)^\circ=ij,(ij)^\bullet=ji$, we have:
\begin{eqnarray*}
Tr(Z^k)
&=&Tr(Z^{k_1}\ldots Z^{k_s})\\
&=&\sum_{i_1=1}^N\ldots\sum_{i_s=1}^N(Z^{k_1})_{i_1i_2}(Z^{k_2})_{i_2i_3}\ldots(Z^{k_s})_{i_si_1}\\
&=&\sum_{i_1=1}^N\ldots\sum_{i_s=1}^N(Z_{(i_1i_2)^{k_1}})^{k_1}(Z_{(i_2i_3)^{k_2}})^{k_2}\ldots(Z_{(i_si_1)^{k_s}})^{k_s}
\end{eqnarray*}

(2) Next, we rescale our variable $Z$ by a $\sqrt{N}$ factor, as in the statement, and we also replace the usual trace by its normalized version, $tr=Tr/N$. Our formula becomes:
$$tr\left(\left(\frac{Z}{\sqrt{N}}\right)^k\right)=\frac{1}{N^{s/2+1}}\sum_{i_1=1}^N\ldots\sum_{i_s=1}^N(Z_{(i_1i_2)^{k_1}})^{k_1}(Z_{(i_2i_3)^{k_2}})^{k_2}\ldots(Z_{(i_si_1)^{k_s}})^{k_s}$$

Thus, the moment that we are interested in is given by:
$$M_k\left(\frac{Z}{\sqrt{N}}\right)=\frac{1}{N^{s/2+1}}\sum_{i_1=1}^N\ldots\sum_{i_s=1}^N\int_X(Z_{(i_1i_2)^{k_1}})^{k_1}(Z_{(i_2i_3)^{k_2}})^{k_2}\ldots(Z_{(i_si_1)^{k_s}})^{k_s}$$

(3) Let us apply now the Wick formula, from Theorem 3.5. We conclude that the moment that we are interested in is given by the following formula:
\begin{eqnarray*}
&&M_k\left(\frac{Z}{\sqrt{N}}\right)\\
&=&\frac{t^{s/2}}{N^{s/2+1}}\sum_{i_1=1}^N\ldots\sum_{i_s=1}^N\#\left\{\pi\in\mathcal P_2(k)\Big|\pi\leq\ker\left((i_1i_2)^{k_1},(i_2i_3)^{k_2},\ldots,(i_si_1)^{k_s}\right)\right\}\\
&=&t^{s/2}\sum_{\pi\in\mathcal P_2(k)}\frac{1}{N^{s/2+1}}\#\left\{i\in\{1,\ldots,N\}^s\Big|\pi\leq\ker\left((i_1i_2)^{k_1},(i_2i_3)^{k_2},\ldots,(i_si_1)^{k_s}\right)\right\}
\end{eqnarray*}

(4) Our claim now is that in the $N\to\infty$ limit the combinatorics of the above sum simplifies, with only the noncrossing partitions contributing to the sum, and with each of them contributing precisely with a 1 factor, so that we will have, as desired:
\begin{eqnarray*}
M_k\left(\frac{Z}{\sqrt{N}}\right)
&=&t^{s/2}\sum_{\pi\in\mathcal P_2(k)}\Big(\delta_{\pi\in NC_2(k)}+O(N^{-1})\Big)\\
&\simeq&t^{s/2}\sum_{\pi\in\mathcal P_2(k)}\delta_{\pi\in NC_2(k)}\\
&=&t^{s/2}|\mathcal{NC}_2(k)|
\end{eqnarray*}

(5) In order to prove this, the first observation is that when $k$ is not uniform, in the sense that it contains a different number of $\circ$, $\bullet$ symbols, we have $\mathcal P_2(k)=\emptyset$, and so:
$$M_k\left(\frac{Z}{\sqrt{N}}\right)=t^{s/2}|\mathcal{NC}_2(k)|=0$$

(6) Thus, we are left with the case where $k$ is uniform. Let us examine first the case where $k$ consists of an alternating sequence of $\circ$ and $\bullet$ symbols, as follows:
$$k=\underbrace{\circ\bullet\circ\bullet\ldots\ldots\circ\bullet}_{2p}$$

In this case it is convenient to relabel our multi-index $i=(i_1,\ldots,i_s)$, with $s=2p$, in the form $(j_1,l_1,j_2,l_2,\ldots,j_p,l_p)$. With this done, our moment formula becomes:
$$M_k\left(\frac{Z}{\sqrt{N}}\right)
=t^p\sum_{\pi\in\mathcal P_2(k)}\frac{1}{N^{p+1}}\#\left\{j,l\in\{1,\ldots,N\}^p\Big|\pi\leq\ker\left(j_1l_1,j_2l_1,j_2l_2,\ldots,j_1l_p\right)\right\}$$

Now observe that, with $k$ being as above, we have an identification $\mathcal P_2(k)\simeq S_p$, obtained in the obvious way. With this done too, our moment formula becomes:
$$M_k\left(\frac{Z}{\sqrt{N}}\right)
=t^p\sum_{\pi\in S_p}\frac{1}{N^{p+1}}\#\left\{j,l\in\{1,\ldots,N\}^p\Big|j_r=j_{\pi(r)+1},l_r=l_{\pi(r)},\forall r\right\}$$

(7) We are now ready to do our asymptotic study, and prove the claim in (4). Let indeed $\gamma\in S_p$ be the full cycle, which is by definition the following permutation:
$$\gamma=(1 \, 2 \, \ldots \, p)$$

In terms of $\gamma$, the conditions $j_r=j_{\pi(r)+1}$ and $l_r=l_{\pi(r)}$ found above read:
$$\gamma\pi\leq\ker j\quad,\quad 
\pi\leq\ker l$$

Counting the number of free parameters in our moment formula, we obtain:
$$M_k\left(\frac{Z}{\sqrt{N}}\right)
=\frac{t^p}{N^{p+1}}\sum_{\pi\in S_p}N^{|\pi|+|\gamma\pi|}
=t^p\sum_{\pi\in S_p}N^{|\pi|+|\gamma\pi|-p-1}$$

(8) Now by Proposition 3.6 the last exponent is $\leq 0$, with equality precisely when $\pi\in S_p$ comes from a noncrossing partition. Thus, with $N\to\infty$, we obtain:
$$M_k\left(\frac{Z}{\sqrt{N}}\right)\simeq t^p|\mathcal{NC}_2(k)|$$

But this finishes the proof for the exponents $k$ which are alternating, and the case where $k$ is an arbitrary uniform exponent is similar, by permuting everything.
\end{proof}

We can reformulate what we found in Theorem 3.7, in the following way:

\begin{theorem}
Given a sequence of complex Gaussian matrices $Z_N\in M_N(L^\infty(X))$ having independent $G_t$ variables as entries, with $t>0$, we have
$$\frac{Z_N}{\sqrt{N}}\sim\Gamma_t$$
in the $N\to\infty$ limit, with the limiting law being defined via its colored moments by
$$M_k(\Gamma_t)=\sum_{\pi\in\mathcal{NC}_2(k)}t^{|\pi|}$$
and with this being a law which is not a usual measure.
\end{theorem}

\begin{proof}
According to what we found in Theorem 3.7, the limiting law $\Gamma_t$, that of the sequence of matrices $(Z_N/\sqrt{N})_{N\in\mathbb N}$, has moments as follows:
$$M_k(\Gamma_t)=t^{|k|/2}|\mathcal{NC}_2(k)|$$

Now since the number of blocks of a pairing $\pi\in\mathcal P_2(k)$ is always $|k|/2$, independently of the pairing, we are led to the moment formula for $\Gamma_t$ in the statement.
\end{proof}

Next, we can make the following speculation on the law $\Gamma_t$ that we found:

\begin{speculation}
In view of the moment formulae for $G_t$ and $\Gamma_t$, namely
$$M_k(G_t)=\sum_{\pi\in\mathcal P_2(k)}t^{|\pi|}\quad,\quad 
M_k(\Gamma_t)=\sum_{\pi\in\mathcal{NC}_2(k)}t^{|\pi|}$$
and assuming that the crossings correspond somehow to commutation
$$\slash\hskip-2.1mm\backslash\quad\longleftrightarrow\quad ab=ba$$
we can say that $\Gamma_t$ is the free analogue of the complex normal law $G_t$.
\end{speculation}

And wild speculation this is, but if you have something better, at this stage of things, let me know. As some sort of rationale for this, we can regard $\slash\hskip-2.1mm\backslash$ as a Feynman diagram, which will make the incoming particles $a,b$ switch position at the arrival:
$$\xymatrix@R=35pt@C=25pt{
a\ar[dr]&b\ar[dl]\\
b&a}$$

Now since obviously nothing happened, I mean there was no collision or scattering or something, we conclude that having $\slash\hskip-2.1mm\backslash$ in our theory means to assume $ab=ba$. Which sounds sweet, in the hope that I convinced you guys of the validity of my theory. 

\bigskip

We will see later that free probability provides a conceptual approach to all this, with the free analogue of $G_t$ being indeed the above law $\Gamma_t$, called Voiculescu circular law. As for our $\slash\hskip-2.1mm\backslash\longleftrightarrow ab=ba$ speculation, this will be not forgotten either, and we will make some mathematics out of it, using quantum groups and easiness. More later.

\bigskip

As a second comment now on Theorem 3.8, which is advanced too, we have:

\begin{comment}
Given a family of sequences of complex Gaussian matrices, 
$$Z^i_N\in M_N(L^\infty(X))\quad,\quad i\in I$$
having independent $G_t$ variables as entries, with $t>0$, the rescaled sequences
$$\frac{Z^i_N}{\sqrt{N}}\in M_N(L^\infty(X))\quad,\quad i\in I$$
become with $N\to\infty$ circular, each following the Voiculescu law $\Gamma_t$, and free.
\end{comment}

To be more precise, what is new here with respect to Theorem 3.8, reformulated as in Speculation 3.9, taking into account the fact that we are now dealing with a family of Gaussian matrices, is the word ``free'' at the end, and with this coming from a famous theorem of Voiculescu \cite{vo4}. So, this is how freeness basically appears, in the random matrix context, and more on this later in this book, when doing free probability.

\section*{3b. Wigner matrices}

Getting now to the Wigner matrices, these are meant to be more understandable, because they are self-adjoint. We first have the following result, about them:

\begin{theorem}
Given a sequence of Wigner random matrices $Z_N\in M_N(L^\infty(X))$ having independent $G_t$ variables as entries, with $t>0$, up to $Z_N=Z_N^*$, we have
$$M_k\left(\frac{Z_N}{\sqrt{N}}\right)\simeq t^{k/2}|NC_2(k)|$$
for any $k\in\mathbb N$, in the $N\to\infty$ limit, with $NC_2=P_2\cap NC$, the noncrossing pairings.
\end{theorem}

\begin{proof}
We have two possible proofs here, as follows:

\medskip

(1) The formula in the statement can be certainly established via a direct computation based on the Wick formula, similar to that from the proof of Theorem 3.7. We will leave the calculations here, featuring no particular difficulties, as an exercise.

\medskip

(2) However, the best is to deduce this result from Theorem 3.7 itself. Indeed, we know from there that for Gaussian matrices $Y_N\in M_N(L^\infty(X))$ we have the following formula, valid for any colored integer $K=\circ\bullet\bullet\circ\ldots\,$, in the $N\to\infty$ limit:
$$M_K\left(\frac{Y_N}{\sqrt{N}}\right)\simeq t^{|K|/2}|\mathcal{NC}_2(K)|$$

We deduce from this that we have the following formula for the moments of the matrices $Re(Y_N)$, with respect to usual exponents, $k\in\mathbb N$:
\begin{eqnarray*}
M_k\left(\frac{Re(Y_N)}{\sqrt{N}}\right)
&=&2^{-k}\cdot M_k\left(\frac{Y_N}{\sqrt{N}}+\frac{Y_N^*}{\sqrt{N}}\right)\\
&=&2^{-k}\sum_{|K|=k}M_K\left(\frac{Y_N}{\sqrt{N}}\right)\\
&\simeq&2^{-k}\sum_{|K|=k}t^{k/2}|\mathcal{NC}_2(K)|\\
&=&2^{-k}\cdot t^{k/2}\cdot 2^{k/2}|\mathcal{NC}_2(k)|\\
&=&2^{-k/2}\cdot t^{k/2}|NC_2(k)|
\end{eqnarray*}

In terms of the matrices $Z_N=\sqrt{2}Re(Y_N)$, which are of Wigner type, this reads:
$$M_k\left(\frac{Z_N}{\sqrt{2N}}\right)\simeq 2^{-k/2}\cdot t^{k/2}|NC_2(k)|$$

Thus, we are led to the moment formula in the statement.
\end{proof}

We can reformulate what we found in Theorem 3.11, in the following way:

\begin{theorem}
Given a sequence of Wigner random matrices $Z_N\in M_N(L^\infty(X))$ having independent $G_t$ variables as entries, with $t>0$, up to $Z_N=Z_N^*$, we have
$$\frac{Z_N}{\sqrt{N}}\sim\gamma_t$$
in the $N\to\infty$ limit, with the limiting law being defined via its moments by
$$M_k(\gamma_t)=\sum_{\pi\in NC_2(k)}t^{|\pi|}$$
and with this law being this time a usual real probability measure.
\end{theorem}

\begin{proof}
According to what we found in Theorem 3.11, the limiting law $\gamma_t$, that of the sequence of matrices $(Z_N/\sqrt{N})_{N\in\mathbb N}$, has moments as follows:
$$M_k(\gamma_t)=t^{k/2}|NC_2(k)|$$

Now since the number of blocks of a pairing $\pi\in P_2(k)$ is always $k/2$, independently of the pairing, we are led to the moment formula for $\gamma_t$ in the statement.
\end{proof}

Next, in analogy with our previous Speculation 3.9, we can formulate:

\begin{speculation}
In view of the moment formulae for $g_t$ and $\gamma_t$, namely
$$M_k(g_t)=\sum_{\pi\in P_2(k)}t^{|\pi|}\quad,\quad 
M_k(\gamma_t)=\sum_{\pi\in NC_2(k)}t^{|\pi|}$$
we can say that $\gamma_t$ is the free analogue of the normal law $g_t$.
\end{speculation}

To be more precise, this is something quite self-explanatory, and more of course on such things later, when doing free probability. And with the further comment that we have as well asymptotic freeness in this setting, along the lines of Comment 3.10. 

\bigskip

With this discussed, what is next? More computations I guess, the plan being:

\begin{plan}
Leaving aside everything advanced, what we have to do is:
\begin{enumerate}
\item Compute the numbers $C_k=|NC_2(2k)|$.

\item Find the measure $\gamma_t$ having $M_{2k}=t^kC_k$ as even moments.
\end{enumerate}
\end{plan}

In what regards the first question, the numbers $C_k=|NC_2(2k)|$ that we came upon are the well-known Catalan numbers from combinatorics, counting many things:

\index{Catalan numbers}
\index{Dyck paths}

\begin{theorem}
The Catalan numbers $C_k=\frac{1}{k+1}\binom{2k}{k}$ count:
\begin{enumerate}
\item The length $2k$ loops on $\mathbb N$, based at $0$.

\item The noncrossing pairings of $1,\ldots,2k$.

\item The noncrossing partitions of $1,\ldots,k$.

\item The length $2k$ Dyck paths in the plane.
\end{enumerate}
\end{theorem}

\begin{proof}
All this is standard combinatorics, the idea being as follows:

\medskip

(1) To start with, in what regards the various objects involved, the length $2k$ loops on $\mathbb N$ are something self-explanatory, and the same goes for the noncrossing pairings of $1,\ldots,2k$, that we already met, and for the noncrossing partitions of $1,\ldots,k$ too, the idea here being that you must be able to draw the partition in a noncrossing way. 

\medskip

(2) Regarding now the length $2k$ Dyck paths in the plane, these are by definition the paths from $(0,0)$ to $(k,k)$, marching North-East over the integer lattice $\mathbb Z^2\subset\mathbb R^2$, by staying inside the square $[0,k]\times[0,k]$, and staying as well under the diagonal of this square. As an example, here are the 5 possible Dyck paths at $n=3$:
$$\xymatrix@R=4pt@C=4pt
{\circ&\circ&\circ&\circ\\
\circ&\circ&\circ&\circ\ar@{-}[u]\\
\circ&\circ&\circ&\circ\ar@{-}[u]\\
\circ\ar@{-}[r]&\circ\ar@{-}[r]&\circ\ar@{-}[r]&\circ\ar@{-}[u]}
\qquad
\xymatrix@R=4pt@C=4pt
{\circ&\circ&\circ&\circ\\
\circ&\circ&\circ&\circ\ar@{-}[u]\\
\circ&\circ&\circ\ar@{-}[r]&\circ\ar@{-}[u]\\
\circ\ar@{-}[r]&\circ\ar@{-}[r]&\circ\ar@{-}[u]&\circ}
\qquad
\xymatrix@R=4pt@C=4pt
{\circ&\circ&\circ&\circ\\
\circ&\circ&\circ\ar@{-}[r]&\circ\ar@{-}[u]\\
\circ&\circ&\circ\ar@{-}[u]&\circ\\
\circ\ar@{-}[r]&\circ\ar@{-}[r]&\circ\ar@{-}[u]&\circ}
\qquad
\xymatrix@R=4pt@C=4pt
{\circ&\circ&\circ&\circ\\
\circ&\circ&\circ&\circ\ar@{-}[u]\\
\circ&\circ\ar@{-}[r]&\circ\ar@{-}[r]&\circ\ar@{-}[u]\\
\circ\ar@{-}[r]&\circ\ar@{-}[u]&\circ&\circ}
\qquad
\xymatrix@R=4pt@C=4pt
{\circ&\circ&\circ&\circ\\
\circ&\circ&\circ\ar@{-}[r]&\circ\ar@{-}[u]\\
\circ&\circ\ar@{-}[r]&\circ\ar@{-}[u]&\circ\\
\circ\ar@{-}[r]&\circ\ar@{-}[u]&\circ&\circ}
$$

(3) Thus, we have definitions for all the objects involved, and in each case, if you start counting them, you always end up with the same sequence of numbers:
$$1,2,5,14,42,132,429,1430,4862,16796,58786,\ldots$$

(4) In order to prove now that (1-4) produce indeed the same numbers, many things can be said. The idea is that, leaving aside mathematical brevity, and more specifically abstract reasonings of type $a=b,b=c\implies a=c$, what we have to do, in order to fully understand what is going on, is to etablish $\binom{4}{2}=6$ equalities, via bijective proofs.

\medskip

(5) But this can be done, indeed. As an example here, the noncrossing pairings of $1,\ldots,2k$ from (2) are in bijection with the noncrossing partitions of $1,\ldots,k$ from (3), via  fattening the pairings and shrinking the partitions. We will leave the details here as an instructive exercise, and exercise as well, to add (1) and (4) to the picture.

\medskip

(6) However, as a matter of having our theorem formally proved, here is quick argument, which does the job. The point is that, in each of the cases (1-4) under consideration, the numbers $C_k$ that we get are easily seen to be subject to the following recurrence:
$$C_{k+1}=\sum_{a+b=k}C_aC_b$$ 

Now the initial data being the same, namely $C_1=1$ and $C_2=2$, in each of the cases (1-4) under consideration, we get indeed the same numbers.

\medskip

(7) Next, let us count the Dyck paths in the plane. For this purpose, we use a trick. Indeed, if we ignore the assumption that our path must stay under the diagonal of the square, we have $\binom{2k}{k}$ such paths. And among these, we have the ``good'' ones, those that we want to count, and then the ``bad'' ones, those that we want to ignore.

\medskip

(8) So, let us count the bad paths, those crossing the diagonal of the square, and reaching the higher diagonal next to it, the one joining $(0,1)$ and $(k,k+1)$. In order to count these, the trick is to ``flip'' their bad part over that higher diagonal, as follows:
$$\xymatrix@R=6pt@C=6pt
{\cdot&\cdot&\cdot&\cdot&\cdot&\cdot\\
\circ&\circ&\circ&\circ\ar@{-}[r]&\circ\ar@{-}[r]\ar@{.}[u]&\circ\\
\circ&\circ\ar@{.}[r]&\circ\ar@{.}[r]&\circ\ar@{.}[r]\ar@{-}[u]&\circ\ar@{.}[u]&\circ\\
\circ&\circ\ar@{.}[u]&\circ&\circ\ar@{-}[u]&\circ&\circ\\
\circ&\circ\ar@{-}[r]\ar@{.}[u]&\circ\ar@{-}[r]&\circ\ar@{-}[u]&\circ&\circ\\
\circ&\circ\ar@{-}[u]&\circ&\circ&\circ&\circ\\
\circ\ar@{-}[r]&\circ\ar@{-}[u]&\circ&\circ&\circ&\circ}$$

(9) Now observe that, as it is obvious on the above picture, due to the flipping, the flipped bad path will no longer end in $(k,k)$, but rather in $(k-1,k+1)$. Moreover, more is true, in the sense that, by thinking a bit, we see that the flipped bad paths are precisely those ending in $(k-1,k+1)$. Thus, we can count these flipped bad paths, and so the bad paths, and so the good paths too, and so good news, we are done.

\medskip

(10) To finish now, by putting everything together, we have:
$$C_k
=\binom{2k}{k}-\binom{2k}{k-1}
=\binom{2k}{k}-\frac{k}{k+1}\binom{2k}{k}
=\frac{1}{k+1}\binom{2k}{k}$$

Thus, we are led to the formula in the statement.
\end{proof}

At a more advanced level, as an analytic alternative to Theorem 3.15, we have:

\begin{theorem}
The Catalan numbers have the following properties:
\begin{enumerate}
\item They satisfy $C_{k+1}=\sum_{a+b=k}C_aC_b$.

\item The series $f(z)=\sum_{k\geq0}C_kz^k$ satisfies $zf^2-f+1=0$.

\item This series is given by $f(z)=\frac{1-\sqrt{1-4z}}{2z}$.

\item We have the formula $C_k=\frac{1}{k+1}\binom{2k}{k}$.
\end{enumerate}
\end{theorem}

\begin{proof}
This is something very standard, the idea being as follows: 

\medskip

(1) Let us count the noncrossing pairings of $\{1,\ldots,2k+2\}$. Such a pairing appears by pairing 1 to an odd number, $2a+1$, and then inserting a noncrossing pairing of $\{2,\ldots,2a\}$, and a noncrossing pairing of $\{2a+2,\ldots,2k+2\}$. Thus we have, as claimed:
$$C_{k+1}=\sum_{a+b=k}C_aC_b$$ 

(2) Consider now the generating series of the Catalan numbers, $f(z)=\sum_{k\geq0}C_kz^k$. In terms of this generating series, the above recurrence gives, as desired:
$$zf^2
=\sum_{a,b\geq0}C_aC_bz^{a+b+1}
=\sum_{k\geq1}\sum_{a+b=k-1}C_aC_bz^k
=\sum_{k\geq1}C_kz^k
=f-1$$

(3) By solving the equation $zf^2-f+1=0$ found above, and choosing the solution which is bounded at $z=0$, we obtain the following formula, as claimed:
$$f(z)=\frac{1-\sqrt{1-4z}}{2z}$$ 

(4) In order to compute this function, we use the generalized binomial formula:
$$(1+t)^p=\sum_{k=0}^\infty\binom{p}{k}t^k$$

For the exponent $p=1/2$, the generalized binomial coefficients are:
\begin{eqnarray*}
\binom{1/2}{k}
&=&\frac{1/2(-1/2)(-3/2)\ldots(3/2-k)}{k!}\\
&=&(-1)^{k-1}\frac{1\cdot 3\cdot 5\ldots(2k-3)}{2^kk!}\\
&=&(-1)^{k-1}\frac{(2k-2)!}{2^{k-1}(k-1)!2^kk!}\\
&=&\frac{(-1)^{k-1}}{2^{2k-1}}\cdot\frac{1}{k}\binom{2k-2}{k-1}\\
&=&-2\left(\frac{-1}{4}\right)^k\cdot\frac{1}{k}\binom{2k-2}{k-1}
\end{eqnarray*}

Thus the generalized binomial formula at exponent $p=1/2$ reads:
$$\sqrt{1+t}=1-2\sum_{k=1}^\infty\frac{1}{k}\binom{2k-2}{k-1}\left(\frac{-t}{4}\right)^k$$

Now with $t=-4z$ we obtain from this the following formula:
$$\sqrt{1-4z}=1-2\sum_{k=1}^\infty\frac{1}{k}\binom{2k-2}{k-1}z^k$$

Getting back to our series $f$, we obtain the following formula for it:
$$f(z)
=\frac{1-\sqrt{1-4z}}{2z}
=\sum_{k=1}^\infty\frac{1}{k}\binom{2k-2}{k-1}z^{k-1}
=\sum_{k=0}^\infty\frac{1}{k+1}\binom{2k}{k}z^k$$

Thus, the Catalan numbers are indeed given by $C_k=\frac{1}{k+1}\binom{2k}{k}$, as claimed.
\end{proof}

Good to have this done, fully answering question (1) in our Plan 3.14. In order to deal now with question (2) there, we will need a non-trivial result, as follows:

\index{Cauchy transform}
\index{Stieltjes inversion}
\index{density of measure}
\index{sequence of moments}
\index{moment problem}

\begin{theorem}
The density of a real probability measure $\mu$ can be recaptured from the sequence of moments $(M_k)_{k\geq0}$ via the Stieltjes inversion formula
$$d\mu (x)=\lim_{s\searrow 0}-\frac{1}{\pi}\,Im\left(G(x+is)\right)\cdot dx$$
with the function on the right, given in terms of moments by
$$G(\xi)=\xi^{-1}+M_1\xi^{-2}+M_2\xi^{-3}+\ldots$$
being the Cauchy transform of the measure $\mu$.
\end{theorem}

\begin{proof}
The Cauchy transform of our measure $\mu$ is given by:
\begin{eqnarray*}
G(\xi)
&=&\xi^{-1}\sum_{k=0}^\infty M_k\xi^{-k}\\\
&=&\int_\mathbb R\frac{\xi^{-1}}{1-\xi^{-1}y}\,d\mu(y)\\
&=&\int_\mathbb R\frac{1}{\xi-y}\,d\mu(y)
\end{eqnarray*}

Now with $\xi=x+is$, we obtain the following formula:
\begin{eqnarray*}
Im(G(x+is))
&=&\int_\mathbb RIm\left(\frac{1}{x-y+is}\right)d\mu(y)\\
&=&\int_\mathbb R\frac{1}{2i}\left(\frac{1}{x-y+is}-\frac{1}{x-y-is}\right)d\mu(y)\\
&=&-\int_\mathbb R\frac{s}{(x-y)^2+s^2}\,d\mu(y)
\end{eqnarray*}

By integrating over $[a,b]$ we obtain, with the change of variables $x=y+sz$:
\begin{eqnarray*}
\int_a^bIm(G(x+is))dx
&=&-\int_\mathbb R\int_a^b\frac{s}{(x-y)^2+s^2}\,dx\,d\mu(y)\\
&=&-\int_\mathbb R\int_{(a-y)/s}^{(b-y)/s}\frac{s}{(sz)^2+s^2}\,s\,dz\,d\mu(y)\\
&=&-\int_\mathbb R\int_{(a-y)/s}^{(b-y)/s}\frac{1}{1+z^2}\,dz\,d\mu(y)\\
&=&-\int_\mathbb R\left(\arctan\frac{b-y}{s}-\arctan\frac{a-y}{s}\right)d\mu(y)
\end{eqnarray*}

Now observe that with $s\searrow0$ we have:
$$\lim_{s\searrow0}\left(\arctan\frac{b-y}{s}-\arctan\frac{a-y}{s}\right)
=\begin{cases}
\frac{\pi}{2}-\frac{\pi}{2}=0& (y<a)\\
\frac{\pi}{2}-0=\frac{\pi}{2}& (y=a)\\
\frac{\pi}{2}-(-\frac{\pi}{2})=\pi& (a<y<b)\\
0-(-\frac{\pi}{2})=\frac{\pi}{2}& (y=b)\\
-\frac{\pi}{2}-(-\frac{\pi}{2})=0& (y>b)
\end{cases}$$

We therefore obtain the following formula:
$$\lim_{s\searrow0}\int_a^bIm(G(x+is))dx=-\pi\left(\mu(a,b)+\frac{\mu(a)+\mu(b)}{2}\right)$$

Now when $\mu$ is continuous, this leads to the formula in the statement.
\end{proof}

As a comment on this, Theorem 3.17 assumes of course that the measure $\mu$ to be found has indeed a density, the point being that, in the opposite case, the density found there will not integrate up to 1. However, we can deal with atoms too, as follows:
$$\mu=\nu+\sum_ic_i\delta_{a_i}\implies G_\mu(\xi)=G_{\nu}(\xi)+\sum_i\frac{c_i}{\xi-a_i}$$

In other words, the poles of the Cauchy transform $G_\mu$ correspond to the atoms of $\mu$. We can now solve our moment problem for the Catalan numbers, as follows:

\begin{theorem}
The centered measure having $M_{2k}=t^kC_k$ as even moments is
$$\gamma_t=\frac{1}{2\pi t}\sqrt{4t-x^2}\,dx$$
the Wigner semicircle law on $[-2\sqrt{t},2\sqrt{t}]$.
\end{theorem}

\begin{proof}
We have two possible proofs for this, as follows:

\medskip

(1) Stieltjes inversion. With $M_{2k}=t^kC_k$, the Cauchy transform is: 
$$G(\xi)
=\xi^{-1}\sum_{k=0}^\infty\frac{1}{k+1}\binom{2k}{k}(t\xi^{-2})^k
=\frac{\xi-\sqrt{\xi^2-4t}}{2t}$$

Now let us apply Theorem 3.17. The study here goes as follows:

\medskip

-- According to the general philosophy of the Stieltjes formula, the first term, namely $\xi/2t$, which is ``trivial'', will not contribute to the density. 

\medskip

-- As for the second term, which is something non-trivial, this will contribute to the density, the rule here being that the square root $\sqrt{\xi^2-4t}$ will be replaced by the ``dual'' square root $\sqrt{4t-x^2}\,dx$, and that we have to multiply everything by $-1/\pi$. 

\medskip

As a conclusion, by Stieltjes inversion we obtain the following density, as claimed:
$$d\mu(x)
=-\frac{1}{\pi}\cdot-\frac{\sqrt{4t-x^2}}{2t}\,dx
=\frac{1}{2\pi t}\sqrt{4t-x^2}dx$$

(2) Cheating. Assuming that a cat or a good physicist suggested you to look at the semicircle law $\gamma_t$, you can prove that this law $\gamma_t$ does indeed the job, with the moment computation being as follows, using the change of variables $x=2\sqrt{t}\cos s$:
\begin{eqnarray*}
M_{2k}
&=&\frac{1}{\pi t}\int_0^{2\sqrt{t}}x^{2k}\sqrt{4t-x^2}\,dx\\
&=&\frac{1}{\pi t}\int_0^{\pi/2}(2\sqrt{t}\cos s)^{2k}\sqrt{4t-4t\cos^2s}\cdot2\sqrt{t}\sin s\,ds\\
&=&\frac{4^{k+1}t^k}{\pi}\int_0^{\pi/2}\cos^{2s}t\sin^2s\,ds\\
&=&\frac{4^{k+1}t^k}{\pi}\cdot\frac{\pi}{2}\cdot\frac{(2k)!!2!!}{(2k+3)!!}\\
&=&2\cdot 4^kt^k\cdot\frac{(2k)!/2^kk!}{2^{k+1}(k+1)!}\\
&=&t^kC_k
\end{eqnarray*}

As for the odd moments, these all vanish, because $\gamma_t$ is an even function. 
\end{proof}

Good news, we can now formulate the Wigner theorem, as follows:

\begin{theorem}
Given a sequence of Wigner random matrices $Z_N\in M_N(L^\infty(X))$
having independent $G_t$ variables as entries, with $t>0$, up to $Z_N=Z_N^*$, we have
$$\frac{Z_N}{\sqrt{N}}\sim\gamma_t$$
in the $N\to\infty$ limit, where $\gamma_t$ is the Wigner semicircle law.
\end{theorem}

\begin{proof}
This comes indeed from Theorems 3.12 and 3.18.
\end{proof}

So long for the Wigner theorem. Of course, this was just a beginning, and far more things can be said, regarding the precise nature of the convergence, and the fluctuations, with the problematics here being highly non-trivial, and still subject to research. For more on all this, you can consult the various random matrix books and articles referenced at the end of this book. And by the way, do not forget to learn, either from Wigner himself \cite{wig}, or from a book like Mehta \cite{meh}, how Wigner came upon these matrices, in relation with the physics of the heavy atoms. First class science, all this.

\section*{3c. Wishart matrices}

Let us discuss now the Wishart matrices, which are the positive analogues of the Wigner matrices. Quite surprisingly, the computation here leads to the Catalan numbers, but not in the same way as for the Wigner matrices, the result being as follows:

\begin{theorem}
Given a sequence of complex Wishart matrices
$$Z_N=Y_NY_N^*\in M_N(L^\infty(X))$$
with $Y_N$ being $N\times N$ complex Gaussian of parameter $t>0$, we have
$$M_k\left(\frac{Z_N}{N}\right)\simeq t^kC_k$$
for any exponent $k\in\mathbb N$, in the $N\to\infty$ limit.
\end{theorem}

\begin{proof}
There are several possible proofs for this result, as follows:

\medskip

(1) A first method is by using the formula that we have in Theorem 3.7, for the Gaussian matrices $Y_N$. Indeed, we know from there that we have the following formula, valid for any colored integer $K=\circ\bullet\bullet\circ\ldots\,$, in the $N\to\infty$ limit:
$$M_K\left(\frac{Y_N}{\sqrt{N}}\right)\simeq t^{|K|/2}|\mathcal{NC}_2(K)|$$

With $K=\circ\bullet\circ\bullet\ldots\,$, alternating word of length $2k$, with $k\in\mathbb N$, this gives:
$$M_k\left(\frac{Y_NY_N^*}{N}\right)\simeq t^k|\mathcal{NC}_2(K)|$$

Thus, in terms of the Wishart matrix $Z_N=Y_NY_N^*$ we have, for any $k\in\mathbb N$:
$$M_k\left(\frac{Z_N}{N}\right)\simeq t^k|\mathcal{NC}_2(K)|$$

The point now is that, by doing some combinatorics, we have:
$$|\mathcal{NC}_2(K)|=|NC_2(2k)|=C_k$$

Thus, we are led to the formula in the statement.

\medskip

(2) A second method, that we will explain now as well, is by proving the result directly, starting from definitions. The matrix entries of our matrix $Z=Z_N$ are given by:
$$Z_{ij}=\sum_{r=1}^NY_{ir}\bar{Y}_{jr}$$

Thus, the normalized traces of powers of $Z$ are given by the following formula:
\begin{eqnarray*}
tr(Z^k)
&=&\frac{1}{N}\sum_{i_1=1}^N\ldots\sum_{i_k=1}^NZ_{i_1i_2}Z_{i_2i_3}\ldots Z_{i_ki_1}\\
&=&\frac{1}{N}\sum_{i_1=1}^N\ldots\sum_{i_k=1}^N\sum_{r_1=1}^N\ldots\sum_{r_k=1}^NY_{i_1r_1}\bar{Y}_{i_2r_1}Y_{i_2r_2}\bar{Y}_{i_3r_2}\ldots Y_{i_kr_k}\bar{Y}_{i_1r_k}
\end{eqnarray*}

By rescaling now $Z$ by a $1/N$ factor, as in the statement, we obtain:
$$tr\left(\left(\frac{Z}{N}\right)^k\right)=\frac{1}{N^{k+1}}\sum_{i_1=1}^N\ldots\sum_{i_k=1}^N\sum_{r_1=1}^N\ldots\sum_{r_k=1}^NY_{i_1r_1}\bar{Y}_{i_2r_1}Y_{i_2r_2}\bar{Y}_{i_3r_2}\ldots Y_{i_kr_k}\bar{Y}_{i_1r_k}$$

By using now the Wick rule, we obtain the following formula for the moments, with $K=\circ\bullet\circ\bullet\ldots\,$, alternating word of lenght $2k$, and with $I=(i_1r_1,i_2r_1,\ldots,i_kr_k,i_1r_k)$:
\begin{eqnarray*}
M_k\left(\frac{Z}{N}\right)
&=&\frac{t^k}{N^{k+1}}\sum_{i_1=1}^N\ldots\sum_{i_k=1}^N\sum_{r_1=1}^N\ldots\sum_{r_k=1}^N\#\left\{\pi\in\mathcal P_2(K)\Big|\pi\leq\ker(I)\right\}\\
&=&\frac{t^k}{N^{k+1}}\sum_{\pi\in\mathcal P_2(K)}\#\left\{i,r\in\{1,\ldots,N\}^k\Big|\pi\leq\ker(I)\right\}
\end{eqnarray*}

In order to compute this quantity, we use the standard bijection $\mathcal P_2(K)\simeq S_k$. By identifying the pairings $\pi\in\mathcal P_2(K)$ with their counterparts $\pi\in S_k$, we obtain:
\begin{eqnarray*}
M_k\left(\frac{Z}{N}\right)
&=&\frac{t^k}{N^{k+1}}\sum_{\pi\in S_k}\#\left\{i,r\in\{1,\ldots,N\}^k\Big|i_s=i_{\pi(s)+1},r_s=r_{\pi(s)},\forall s\right\}
\end{eqnarray*}

As before in our study of the Gaussian matrices, let $\gamma=(1 \, 2 \, \ldots \, k)\in S_k$ be the full cycle. By counting the number of free parameters in our moment formula, we get:
$$M_k\left(\frac{Z}{N}\right)
=t^k\sum_{\pi\in S_k}N^{|\pi|+|\gamma\pi|-k-1}$$

The point now is that the last exponent is $\leq 0$, with equality precisely when $\pi$ comes from a noncrossing partition. Thus we obtain, in the $N\to\infty$ limit:
$$M_k\left(\frac{Z}{N}\right)\simeq t^kC_k$$

We are therefore led to the conclusion in the statement.
\end{proof}

We can reformulate what we found in Theorem 3.20, in the following way:

\begin{theorem}
Given a sequence of complex Wishart matrices
$$Z_N=Y_NY_N^*\in M_N(L^\infty(X))$$
with $Y_N$ being $N\times N$ complex Gaussian of parameter $t>0$, we have
$$\frac{Z_N}{tN}\sim\pi_1$$
in the $N\to\infty$ limit, with the limiting law being defined via its moments by
$$M_k(\pi_1)=|NC(k)|$$
and with this law being a usual real positive probability measure.
\end{theorem}

\begin{proof}
This follows indeed from our formula from Theorem 3.20.
\end{proof}

Next, in analogy with our previous study of the Gaussian and Wigner matrices, we have the following speculation, and more on this later, when doing free probability:

\begin{speculation}
In view of the moment formulae for $p_1$ and $\pi_1$, namely
$$M_k(p_1)=|P(k)|\quad,\quad 
M_k(\pi_1)=|NC(k)|$$
we can say that $\pi_1$ is the free analogue of the Poisson law $p_1$.
\end{speculation}

As yet another thing to be said about all this, this time of pure combinatorial nature, regarding the Catalan numbers, further clarifying the above, we have:

\index{fattening of partitions}
\index{shrinking partitions}
\index{noncrossing partitions}
\index{noncrossing pairings}

\begin{proposition}
We have a bijection between noncrossing partitions and pairings, $NC(k)\simeq NC_2(2k)$, constructed as follows:
\begin{enumerate}
\item The application $NC(k)\to NC_2(2k)$ is the ``fattening'' one, obtained by doubling all the legs, and doubling all the strings as well.

\item Its inverse $NC_2(2k)\to NC(k)$ is the ``shrinking'' application, obtained by collapsing pairs of consecutive neighbors.
\end{enumerate}
\end{proposition}

\begin{proof}
This is something self-explanatory, and in order to see how this works, let us discuss an example. Consider a noncrossing partition, say the following one:
$$\xymatrix@R=10pt@C=10pt{
\ar@{-}[rrrrrrr]&&&&&&&\\
&\ar@{-}[rr]&&&&\ar@{-}[r]&&\\
1\ar@{-}[uu]&2\ar@{-}[u]&3\ar@{-}[u]&4\ar@{-}[u]&5\ar@{-}[uu]&6\ar@{-}[u]&7\ar@{-}[u]&8\ar@{-}[uu]
}$$

Now let us ``fatten'' this partition, by doubling everything, as follows:
$$\xymatrix@R=10pt@C=10pt{
\ar@{=}[rrrrrrr]&&&&&&&\\
&\ar@{=}[rr]&&&&\ar@{=}[r]&&\\
11'\ar@{=}[uu]&22'\ar@{=}[u]&33'\ar@{=}[u]&44'\ar@{=}[u]&55'\ar@{=}[uu]&66'\ar@{=}[u]&77'\ar@{=}[u]&88'\ar@{=}[uu]
}$$

We can see emerging here a noncrossing pairing, and by relabeling the points $1,\ldots,16$, and properly redrawing the picture, what we have is indeed a noncrossing pairing:
$$\xymatrix@R=10pt@C=5pt{
\ar@{-}[rrrrrrrrrrrrrrr]&&&&&&&&&&&&&&&\\
&\ar@{-}[rrrrrrr]&&&&&&&&\ar@{-}[rrrrr]&&&&&&\\
&&\ar@{-}[rrrrr]&&&&&&&&\ar@{-}[rrr]&&&&&\\
&&&\ar@{-}[r]&&\ar@{-}[r]&&&&&&\ar@{-}[r]&&&&\\
1\ar@{-}[uuuu]&2\ar@{-}[uuu]&3\ar@{-}[uu]&4\ar@{-}[u]&5\ar@{-}[u]&6\ar@{-}[u]&7\ar@{-}[u]&8\ar@{-}[uu]&9\ar@{-}[uuu]&10\ar@{-}[uuu]&11\ar@{-}[uu]&12\ar@{-}[u]&13\ar@{-}[u]&14\ar@{-}[uu]&15\ar@{-}[uuu]&16\ar@{-}[uuuu]
}$$

As for the reverse operation, that is obviously obtained by ``shrinking'' our pairing, by collapsing pairs of consecutive neighbors, that is, by identifying $1=2$, then $3=4$, then $5=6$, and so on. Thus, we are led to the conclusion in the statement.
\end{proof}

Now back to random matrices, by solving the moment problem, we are led to: 

\begin{theorem}
Given a sequence of complex Wishart matrices
$$Z_N=Y_NY_N^*\in M_N(L^\infty(X))$$
with $Y_N$ being $N\times N$ complex Gaussian of parameter $t>0$, we have
$$\frac{Z_N}{tN}\sim\frac{1}{2\pi}\sqrt{4x^{-1}-1}\,dx$$
with $N\to\infty$, with the limiting measure being the Marchenko-Pastur law $\pi_1$.
\end{theorem}

\begin{proof}
We have three possible methods here, all working:

\medskip

(1) Squaring. Indeed, we can capture the density from that of $\gamma_1$, by using:
$$f\sim\gamma_1\ \implies\ f^2\sim\pi_1$$

To be more precise, the computation here, which is straightforward, is as follows:
\begin{eqnarray*}
E(\varphi(f^2))
&=&\frac{1}{2\pi}\int_{-2}^2\varphi(x^2)\sqrt{4-x^2}\,dx\\
&=&\frac{1}{\pi}\int_0^2\varphi(x^2)\sqrt{4-x^2}\,dx\\
&=&\frac{1}{\pi}\int_0^4\varphi(y)\sqrt{4-y}\cdot\frac{dy}{2\sqrt{y}}\\
&=&\frac{1}{2\pi}\int_0^4\varphi(y)\sqrt{4y^{-1}-1}\,dy
\end{eqnarray*}

(2) Stieltjes inversion. This is standard too, because we have:
$$M(z)=\frac{1-\sqrt{1-4z}}{2z}\quad ,\quad G(\xi)=\frac{1-\sqrt{1-4\xi^{-1}}}{2}$$

Thus, by Stieltjes inversion we obtain the following density, as stated:
$$d\mu(x)
=-\frac{1}{\pi}\cdot-\frac{\sqrt{4x^{-1}-1}}{2}\,dx
=\frac{1}{2\pi}\sqrt{4x^{-1}-1}\,dx$$

(3) Cheating. Again this is standard, because the moments of the law $\pi_1$ in the statement can be computed with the change of variable $x=4\cos^2t$, as follows:
\begin{eqnarray*}
M_k
&=&\frac{1}{2\pi}\int_0^4x^k\sqrt{4x^{-1}-1}dx\\
&=&\frac{1}{2\pi}\int_0^{\pi/2}(4\cos^2t)^k\cdot\frac{\sin t}{\cos t}\cdot2\cos t\sin t\,dt\\
&=&\frac{4^{k+1}}{\pi}\int_0^{\pi/2}\cos^{2k}t\sin^2t\,dt\\
&=&\frac{4^{k+1}}{\pi}\cdot\frac{\pi}{2}\cdot\frac{(2k)!!2!!}{(2k+3)!!}\\
&=&2\cdot 4^k\cdot\frac{(2k)!/2^kk!}{2^{k+1}(k+1)!}\\
&=&C_k
\end{eqnarray*}

Thus, one way or another, we are led to the density in the statement. 
\end{proof}

\section*{3d. Marchenko-Pastur}

Looking at what we did in the last section, at $t=1$ that is definitely quality mathematics, but what we got at general $t>0$ is not very satisfactory, because the limiting law $\pi_1$ stays the same, modulo a trivial rescaling. Which is, obviously, not good.

\bigskip

The reasons behind this phenomenon are quite subtle, and skipping some discussion here, the point is that Theorem 3.24 is indeed something ``fake'' at general values $t>0$. The correct generalization of the $t=1$ computation, involving more general classes of complex Wishart matrices, is something a bit more complicated, as follows:

\index{Marchenko-Pastur}

\begin{theorem}
Given a sequence of general complex Wishart matrices
$$Z_N=Y_NY_N^*\in M_N(L^\infty(X))$$
with $Y_N$ being $N\times M$ complex Gaussian of parameter $1$, we have
$$M_k\left(\frac{Z}{N}\right)
\simeq\sum_{\pi\in NC(k)}t^{|\pi|}$$
for any exponent $k\in\mathbb N$, in the $M=tN\to\infty$ limit.
\end{theorem}

\begin{proof}
This is again something which is very standard, as follows:

\medskip

(1) In order to prove the formula in the statement, we can proceed as usual, by using the Wick formula. The matrix entries of our Wishart matrix $Z=Z_N$ are given by:
$$Z_{ij}=\sum_{r=1}^MY_{ir}\bar{Y}_{jr}$$

Thus, the normalized traces of powers of $Z$ are given by the following formula:
\begin{eqnarray*}
tr(Z^k)
&=&\frac{1}{N}\sum_{i_1=1}^N\ldots\sum_{i_k=1}^NZ_{i_1i_2}Z_{i_2i_3}\ldots Z_{i_ki_1}\\
&=&\frac{1}{N}\sum_{i_1=1}^N\ldots\sum_{i_k=1}^N\sum_{r_1=1}^M\ldots\sum_{r_k=1}^MY_{i_1r_1}\bar{Y}_{i_2r_1}Y_{i_2r_2}\bar{Y}_{i_3r_2}\ldots Y_{i_kr_k}\bar{Y}_{i_1r_k}
\end{eqnarray*}

By rescaling now $Z$ by a $1/N$ factor, as in the statement, we obtain:
$$tr\left(\left(\frac{Z}{N}\right)^k\right)=\frac{1}{N^{k+1}}\sum_{i_1=1}^N\ldots\sum_{i_k=1}^N\sum_{r_1=1}^M\ldots\sum_{r_k=1}^MY_{i_1r_1}\bar{Y}_{i_2r_1}Y_{i_2r_2}\bar{Y}_{i_3r_2}\ldots Y_{i_kr_k}\bar{Y}_{i_1r_k}$$

(2) By using now the Wick rule, we obtain the following formula for the moments, with $K=\circ\bullet\circ\bullet\ldots\,$, alternating word of lenght $2k$, and $I=(i_1r_1,i_2r_1,\ldots,i_kr_k,i_1r_k)$:
\begin{eqnarray*}
M_k\left(\frac{Z}{N}\right)
&=&\frac{1}{N^{k+1}}\sum_{i_1=1}^N\ldots\sum_{i_k=1}^N\sum_{r_1=1}^M\ldots\sum_{r_k=1}^M\#\left\{\pi\in\mathcal P_2(K)\Big|\pi\leq\ker I\right\}\\
&=&\frac{1}{N^{k+1}}\sum_{\pi\in\mathcal P_2(K)}\#\left\{i\in\{1,\ldots,N\}^k,r\in\{1,\ldots,M\}^k\Big|\pi\leq\ker I\right\}
\end{eqnarray*}

(3) In order to compute this quantity, we use the standard bijection $\mathcal P_2(K)\simeq S_k$. By identifying the pairings $\pi\in\mathcal P_2(K)$ with their counterparts $\pi\in S_k$, we obtain:
\begin{eqnarray*}
M_k\left(\frac{Z}{N}\right)
&=&\frac{1}{N^{k+1}}\sum_{\pi\in S_k}\#\left\{i\in\{1,\ldots,N\}^k,r\in\{1,\ldots,M\}^k\Big|i_s=i_{\pi(s)+1},r_s=r_{\pi(s)}\right\}
\end{eqnarray*}

With $\gamma=(1 \, 2 \, \ldots \, k)\in S_k$ being the full cycle, by counting the number of free parameters in our expectation formula, as before for other random matrices, we obtain:
$$M_k\left(\frac{Z}{N}\right)
=\frac{1}{N^{k+1}}\sum_{\pi\in S_k}N^{|\gamma\pi|}M^{|\pi|}
=\sum_{\pi\in S_k}N^{|\gamma\pi|-k-1}M^{|\pi|}$$

(4) Now by using the same arguments as in the case $M=N$, from the proof of Theorem 3.20, we conclude that in the $M=tN\to\infty$ limit the permutations $\pi\in S_k$ which matter are those coming from noncrossing partitions, and so that we have:
$$M_k\left(\frac{Z}{N}\right)
\simeq\sum_{\pi\in NC(k)}N^{-|\pi|}M^{|\pi|}
=\sum_{\pi\in NC(k)}t^{|\pi|}$$

Thus, we are led to the conclusion in the statement.
\end{proof}

We can reformulate what we found in Theorem 3.25, in the following way:

\begin{theorem}
Given a sequence of general complex Wishart matrices
$$Z_N=Y_NY_N^*\in M_N(L^\infty(X))$$
with $Y_N$ being $N\times M$ complex Gaussian of parameter $1$, we have
$$\frac{Z_N}{N}\sim\pi_t$$
in the $M=tN\to\infty$ limit, with the limiting law being given by
$$M_k(\pi_t)=\sum_{\pi\in NC(k)}t^{|\pi|}$$
and with this law being a usual real positive probability measure.
\end{theorem}

\begin{proof}
This follows indeed from our formula from Theorem 3.25.
\end{proof}

Next, in analogy with our previous speculations, we can formulate:

\begin{speculation}
In view of the moment formulae for $p_t$ and $\pi_t$, namely
$$M_k(p_t)=\sum_{\pi\in P(k)}t^{|\pi|}\quad,\quad 
M_k(\pi_t)=\sum_{\pi\in NC(k)}t^{|\pi|}$$
we can say that $\pi_t$ is the free analogue of the Poisson law $p_t$.
\end{speculation}

To be more precise, this is something quite self-explanatory, in tune with our previous speculations, and more of course on such things later, when doing free probability. And with the further comment that we have as well asymptotic freeness in this setting. 

\bigskip

Finally, by solving the moment problem, we are led to the following result: 

\begin{theorem}
Given a sequence of general complex Wishart matrices
$$Z_N=Y_NY_N^*\in M_N(L^\infty(X))$$
with $Y_N$ being $N\times M$ complex Gaussian of parameter $1$, we have
$$\frac{Z_N}{N}\sim\max(1-t,0)\delta_0+\frac{\sqrt{4t-(x-1-t)^2}}{2\pi x}\,dx$$
with $M=tN\to\infty$, with the limiting measure being the Marchenko-Pastur law $\pi_t$.
\end{theorem}

\begin{proof}
This is something quite tricky, the idea being as follows:

\medskip

(1) Let us first review the $t=1$ case. Here the moment generating function is:
$$M(z)=\sum_{k\geq0}C_kz^k$$

For $\pi\in NC(k+1)$ we can consider the last leg of the first block, and make a cut at right of it. This gives a decomposition of $\pi$ into $2$ partitions in $NC$, and we get:
$$C_{k+1}=\sum_{k=k_0+k_1}C_{k_0}C_{k_1}$$

By multiplying by $z^{k+1}$, then summing over $k$, we obtain the following equation:
$$M-1=zM^2$$

But this gives $M$, then $G$, then $\pi_1$, as explained in the proof of Theorem 3.24.

\medskip

(2) At $t>0$ now, let us denote by $F_{kb}$ the number of partitions in $NC(k)$ having $b$ blocks, and we set $F_{kb}=0$ for other integer values of $k,b$. Also, let us make the convention that, in what follows, all sums will be over integer indices $\geq 0$. With these conventions, the moment generating function for our measure is given by the following formula:
$$M(z)=\sum_{kb}F_{kb}z^kt^b$$

Now observe that the recurrence formula for the numbers $C_k$ in (1) becomes:
$$\sum_bF_{k+1,b}=\sum_{k=k_0+k_1}\sum_{b_0,b_1}F_{k_0b_0}F_{k_1b_1}$$

In this formula, each term contributes to $F_{k+1,b}$ with a contribution $b=\sum b_i$, except for those of the form $F_{00}F_{k_1b_1}$, which contribute to $F_{k+1,b+1}$. Thus, we have:
$$F_{k+1,b}
=\sum_{k=k_0+k_1}\sum_{b=b_0+b_1}F_{k_0b_0}F_{k_1b_1}
+\sum_{k=k_0+k_1}F_{k_1,b-1}
-\sum_{k=k_0+k_1}F_{k_1,b-1}$$

But this gives the following formula for the polynomials $P_k=\sum_bF_{kb}t^b$:
$$P_{k+1}=\sum_{k=k_0+k_1}P_{k_0}P_{k_1}+(t-1)\sum_{k=k_0+k_1}P_{k_1}$$

In terms of $M(z)=\sum_kP_kz^k$, we obtain the following equation:
$$M-1=zM^2+(t-1)zM$$

(3) But with this, we can go ahead like we did before at $t=1$. Our equation is:
$$zM^2+((t-1)z-1)M+1=0$$

Thus, the moment generating series is given by the following formula:
$$M(z)=\frac{(1-t)z+1-\sqrt{((t-1)z-1)^2-4z}}{2z}$$

Equivalently, the Cauchy transform is given by the following formula:
\begin{eqnarray*}
G(\xi)
&=&\frac{1-t+\xi-\sqrt{(t-1-\xi)^2-4\xi}}{2\xi}\\
&=&\frac{1-t+\xi-\sqrt{(\xi-1-t)^2-4t}}{2\xi}
\end{eqnarray*}

Now by Stieltjes inversion we obtain the following density:
$$d\mu(x)
=\max(1-t,0)\delta_0+\frac{\sqrt{4t-(x-1-t)^2}}{2\pi x}\,dx$$

Summarizing, we have obtained the density in the statement of $\pi_t$.
\end{proof}

So long for the basic theory of the random matrices. We will be back later to this, with further explanations, and some new results too, when doing free probability.

\section*{3e. Exercises}

Welcome to random matrices, and as exercises on them, we have:

\begin{exercise}
Clarify what we said, regarding the partitions $\pi\in NC(k)\subset S_k$.
\end{exercise}

\begin{exercise}
Have some thinking at the circular law $\Gamma_t$ that we found.
\end{exercise}

\begin{exercise}
Establish the combinatorial bijections for the Catalan numbers.
\end{exercise}

\begin{exercise}
Learn more about moments, including Hankel determinants.
\end{exercise}

\begin{exercise}
Learn more about the fattening and shrinking of partitions.
\end{exercise}

\begin{exercise}
Learn about the beta laws, generalizing both $\gamma_1$ and $\pi_1$.
\end{exercise}

\begin{exercise}
Work out the moment details, for the general Wishart matrices.
\end{exercise}

\begin{exercise}
Can you recover $\pi_t$, without using Stieltjes inversion?
\end{exercise}

As bonus exercise, already mentioned, read Wigner's paper \cite{wig}.

\chapter{Rotation groups}

\section*{4a. Characters, laws}

With the basics of random matrix theory learned, the natural continuation would be getting into free probability. However, before doing that, let us go back to the subgroups $G\subset U_N$, and study them more in detail, along the lines of what we did in chapter 1, for the symmetric group $S_N$ and other groups. We have 3 motivations, as follows:

\bigskip

(1) To start with, we have unfinished business with $O_N,U_N$. So, we will do this, integrtation over $O_N,U_N$, called Weingarten formula, and application to characters. And with this, Weingarten formula, being something extremely useful, for both random matrices and free probability, in potentially any question having $O_N,U_N$ invariance.

\bigskip

(2) As another motivation, once we will be reasonably advanced in free probability, we will start looking for free manifolds $X^+$, in order to do probability over them. And, these free manifolds will appear via a scheme of type $X=G/H\to X^+=G^+/H^+$, with $H\subset G\subset U_N$ being certain subgroups, whose mathematics we need to master well.

\bigskip

(3) As yet another motivation, we will see that for the groups $SU_2,SO_3$, playing a key role in particle physics, the main character $\chi_1$ follows the Wigner and Marchenko-Pastur laws $\gamma_1,\pi_1$. Thus, as strange as this might seem, we have freeness inside $SU_2,SO_3$. And with understanding this freeness being in fact a key problem in particle physics.

\bigskip

Getting to work now, we talked in chapter 1 about the main characters $\chi_1$ of certain special compact subgroups $G\subset U_N$. In general, let us have as starting point:

\index{representation}
\index{character}

\begin{definition}
A unitary representation of a compact group $G$ is a morphism
$$u:G\to U_N\quad,\quad g\to u_g$$
which can be faithful or not. The character of such a representation is the function
$$\chi:G\to\mathbb C\quad,\quad g\to Tr(u_g)$$
where $Tr$ is the usual, unnormalized trace of the $N\times N$ matrices.
\end{definition}

At the level of examples, most of the compact groups that we met so far, finite or continuous, naturally appear as closed subgroups $G\subset U_N$. In this case, the embedding $u:G\subset U_N$ is of course a representation, called fundamental representation.

\bigskip

In general now, let us first discuss the various operations on the representations. We have here the following elementary result, coming from definitions:

\begin{proposition}
The representations of a compact group $G$ are subject to:
\begin{enumerate}
\item Making sums. Given representations $u,v$, of dimensions $N,M$, 
their sum is the $N+M$-dimensional representation $u+v=diag(u,v)$.

\item Making products. Given representations $u,v$, of dimensions $N,M$, their product is the $NM$-dimensional representation $(u\otimes v)_{ia,jb}=u_{ij}v_{ab}$.

\item Taking conjugates. Given a $N$-dimensional representation $u$, its conjugate is the $N$-dimensional representation $(\bar{u})_{ij}=\bar{u}_{ij}$.

\item Spinning by unitaries. Given a $N$-dimensional representation $u$, and a unitary $V\in U_N$, we can spin $u$ by this unitary, $u\to VuV^*$.
\end{enumerate}
\end{proposition}

\begin{proof}
The fact that the operations in the statement are indeed well-defined, among morphisms from $G$ to unitary groups, is indeed clear from definitions.
\end{proof}

In relation now with characters, we have the following result:

\begin{proposition}
We have the following formulae, regarding characters
$$\chi_{u+v}=\chi_u+\chi_v\quad,\quad 
\chi_{u\otimes v}=\chi_u\chi_v\quad,\quad 
\chi_{\bar{u}}=\bar{\chi}_u\quad,\quad
\chi_{VuV^*}=\chi_u$$
in relation with the basic operations for the representations.
\end{proposition}

\begin{proof}
All these assertions are elementary, coming from the following formulae:
$$Tr(diag(U,V))=Tr(U)+Tr(V)\quad,\quad 
Tr(U\otimes V)=Tr(U)Tr(V)$$
$$Tr(\bar{U})=\overline{Tr(U)}\quad,\quad 
Tr(VUV^*)=Tr(U)$$

Thus, we are led to the various character formulae in the statement.
\end{proof}

Assume now that we are given a closed subgroup $G\subset U_N$. By using the above operations, we can construct a whole family of representations of $G$, as follows:

\index{Peter-Weyl representations}

\begin{definition}
Given a closed subgroup $G\subset U_N$, its Peter-Weyl representations are the tensor products between the fundamental representation and its conjugate:
$$u:G\subset U_N\quad,\quad 
\bar{u}:G\subset U_N$$ 
We denote these tensor products $u^{\otimes k}$, with $k=\circ\bullet\bullet\circ\ldots$ being a colored integer, with the colored tensor powers being defined according to the rules 
$$u^{\otimes\circ}=u\quad,\quad
u^{\otimes\bullet}=\bar{u}\quad,\quad
u^{\otimes kl}=u^{\otimes k}\otimes u^{\otimes l}$$
and with the convention that $u^{\otimes\emptyset}$ is the trivial representation $1:G\to U_1$.
\end{definition}

Here are a few examples of such representations, namely those coming from the colored integers of length 2, which will often appear in what follows:
$$u^{\otimes\circ\circ}=u\otimes u\quad,\quad 
u^{\otimes\circ\bullet}=u\otimes\bar{u}$$
$$u^{\otimes\bullet\circ}=\bar{u}\otimes u\quad,\quad
u^{\otimes\bullet\bullet}=\bar{u}\otimes\bar{u}$$

In order to advance, we must develop some general theory. Let us start with:

\index{intertwiners}

\begin{definition}
Given a compact group $G$, and two of its representations,
$$u:G\to U_N\quad,\quad 
v:G\to U_M$$
we define the space of intertwiners between these representations as being 
$$Hom(u,v)=\left\{T\in M_{M\times N}(\mathbb C)\Big|Tu_g=v_gT,\forall g\in G\right\}$$
and we use the following conventions:
\begin{enumerate}
\item We use the notations $Fix(u)=Hom(1,u)$, and $End(u)=Hom(u,u)$.

\item We write $u\sim v$ when $Hom(u,v)$ contains an invertible element.

\item We say that $u$ is irreducible, and write $u\in Irr(G)$, when $End(u)=\mathbb C1$.
\end{enumerate}
\end{definition}

To be more precise, the terminology here is very standard, with Fix, Hom, End standing respectively for fixed points, homomorphisms and endomorphisms. We will see a bit later that irreducible means indecomposable, in a suitable sense. 

\bigskip

Here are now a few basic results, regarding the above linear spaces:

\begin{proposition}
The spaces of intertwiners have the following properties:
\begin{enumerate}
\item $T\in Hom(u,v),S\in Hom(v,w)\implies ST\in Hom(u,w)$.

\item $S\in Hom(u,v),T\in Hom(w,z)\implies S\otimes T\in Hom(u\otimes w,v\otimes z)$.

\item $T\in Hom(u,v)\implies T^*\in Hom(v,u)$.
\end{enumerate}
In abstract terms, we say that the Hom spaces form a tensor $*$-category.
\end{proposition}

\begin{proof}
All the formulae in the statement are clear from definitions, via elementary computations. As for the last assertion, this is something coming from (1,2,3). We will be back to tensor categories later on, with more details on this latter fact.
\end{proof}

In order to advance, we will need the following standard linear algebra fact, which stands as well as a useful complement to the $C^*$-algebra theory from chapter 3:

\begin{proposition}
Let $A\subset M_N(\mathbb C)$ be a $*$-algebra.
\begin{enumerate}
\item We have $1=p_1+\ldots+p_k$, with $p_i\in A$ being central minimal projections.

\item Each of the spaces $A_i=p_iAp_i$ is a non-unital $*$-subalgebra of $A$.

\item We have a non-unital $*$-algebra sum decomposition $A=A_1\oplus\ldots\oplus A_k$.

\item We have unital $*$-algebra isomorphisms $A_i\simeq M_{n_i}(\mathbb C)$, with $n_i=rank(p_i)$.

\item Thus, we have a $*$-algebra isomorphism $A\simeq M_{n_1}(\mathbb C)\oplus\ldots\oplus M_{n_k}(\mathbb C)$.
\end{enumerate}
\end{proposition}

\begin{proof}
Consider indeed an arbitrary $*$-algebra of the $N\times N$ matrices, $A\subset M_N(\mathbb C)$. Let us first look at the center of this algebra, $Z(A)=A\cap A'$. It is elementary to prove that this center, as an algebra, is of the following form:
$$Z(A)\simeq\mathbb C^k$$

Consider now the standard basis $e_1,\ldots,e_k\in\mathbb C^k$, and let  $p_1,\ldots,p_k\in Z(A)$ be the images of these vectors via the above identification. In other words, these elements $p_1,\ldots,p_k\in A$ are central minimal projections, summing up to 1:
$$p_1+\ldots+p_k=1$$

The idea is then that this partition of the unity eventually leads to the block decomposition of $A$, as in the statement, and we leave the details here as an exercise.
\end{proof}

We can now formulate our first Peter-Weyl theorem, as follows:

\index{Peter-Weyl}

\begin{theorem}[PW1]
Let $u:G\to U_N$ be a representation, consider the algebra $A=End(u)$, and write its unit $1=p_1+\ldots+p_k$ as above. We have then 
$$u=v_1+\ldots+v_k$$
with each $v_i$ being an irreducible representation, obtained by restricting $u$ to $Im(p_i)$.
\end{theorem}

\begin{proof}
This follows indeed from Proposition 4.6 and Proposition 4.7:

\medskip

(1) We first associate to our representation $u:G\to U_N$ the corresponding action map on $\mathbb C^N$. If a linear subspace $V\subset\mathbb C^N$ is invariant, the restriction of the action map to $V$ is an action map too, which must come from a subrepresentation $v\subset u$.

\medskip

(2) Consider now a projection $p\in End(u)$. From $pu=up$ we obtain that the linear space $V=Im(p)$ is invariant under $u$, and so this space must come from a subrepresentation $v\subset u$. It is routine to check that the operation $p\to v$ maps subprojections to subrepresentations, and minimal projections to irreducible representations.

\medskip

(3) With these preliminaries in hand, let us decompose the algebra $End(u)$ as above, by using the decomposition $1=p_1+\ldots+p_k$ into central minimal projections. If we denote by $v_i\subset u$ the subrepresentation coming from the vector space $V_i=Im(p_i)$, then we obtain in this way a decomposition $u=v_1+\ldots+v_k$, as in the statement.
\end{proof}

Here is now our second Peter-Weyl theorem, complementing Theorem 4.8:

\begin{theorem}[PW2]
Given a closed subgroup $G\subset_uU_N$, any of its irreducible smooth representations 
$$v:G\to U_M$$
appears inside a tensor product of the fundamental representation $u$ and its adjoint $\bar{u}$.
\end{theorem}

\begin{proof}
This basically follows from Theorem 4.8, by reasoning as follows:

\medskip

(1) Given $v:G\to U_M$, consider its space of coefficients $C_v\subset C(G)$. The operation $w\to C_w$ is then functorial, mapping subrepresentations into linear subspaces.

\medskip

(2) A closed subgroup $G\subset_uU_N$ is a Lie group, and a representation $v:G\to U_M$ is smooth when we have an inclusion $C_v\subset<C_u>$. This is indeed well-known.

\medskip

(3) By definition of the Peter-Weyl representations, as arbitrary tensor products between the fundamental representation $u$ and its conjugate $\bar{u}$, we have:
$$<C_u>=\sum_kC_{u^{\otimes k}}$$

(4) Now by putting together the above observations (2,3) we conclude that we must have an inclusion as follows, for certain exponents $k_1,\ldots,k_p\in\mathbb N$:
$$C_v\subset C_{u^{\otimes k_1}\oplus\ldots\oplus u^{\otimes k_p}}$$

(5) By using now (1), we deduce that we have an inclusion $v\subset u^{\otimes k_1}\oplus\ldots\oplus u^{\otimes k_p}$, and by applying Theorem 4.8, this leads to the conclusion in the statement.
\end{proof}

In order to further advance, we will need the following standard fact:

\index{Haar measure}
\index{Ces\`aro limit}

\begin{theorem}
The integration over a compact group $G\subset U_N$ can be constructed by starting with any faithful positive unital linear form $\varphi\in C(G)^*$, and setting:
$$\int_G=\lim_{n\to\infty}\frac{1}{n}\sum_{k=1}^n\varphi^{*k}$$
Moreover, for any representation $v:G\to U_n$ we have the following formula,
$$\left(\int_Gv_{ij}\right)_{ij}=P$$
where $P$ is the orthogonal projection onto the linear space $Fix(v)$.
\end{theorem}

\begin{proof}
To start with, given a positive unital linear form $\varphi\in C(G)^*$, our claim is that the following Ces\`aro limit converges, for any $f\in C(G)$:
$$\int_\varphi f=\lim_{n\to\infty}\frac{1}{n}\sum_{k=1}^n\varphi^{*k}(f)$$

Indeed, by linearity we can assume that $f$ is the coefficient of certain representation $v:G\to U_n$. But in this case, an elementary computation gives the following formula, with $P_\varphi\geq P$ being the orthogonal projection onto the $1$-eigenspace of $[\varphi(v_{ij})]_{ij}$:
$$\left(\int_\varphi v_{ij}\right)_{ij}=P_\varphi$$

The point now is that when $\varphi\in C(G)^*$ is faithful, by using a standard positivity trick, we can prove that we have $P_\varphi=P$. Assume indeed $P_\varphi\xi=\xi$, and let us set:
$$f=\sum_i\left(\sum_jv_{ij}\xi_j-\xi_i\right)\overline{\left(\sum_kv_{ik}\xi_k-\xi_i\right)}$$

A straightforward computation shows then that $\varphi(f)=0$, and so $f=0$, as desired. Thus, we proved our claim, and with this in hand, the left and right invariance of $\int_G=\int_\varphi$ is clear on coefficients, and so in general, and this gives all the assertions.
\end{proof}

We will need as well an algebraic ingredient for our study, as follows:

\index{Frobenius isomorphism}

\begin{proposition}
We have a Frobenius type isomorphism
$$Hom(v,w)\simeq Fix(v\otimes\bar{w})$$
valid for any two representations $v,w$.
\end{proposition}

\begin{proof}
According to the definitions, we have the following equivalences:
\begin{eqnarray*}
T\in Hom(v,w)
&\iff&Tv=wT\\
&\iff&\sum_jT_{aj}v_{ji}=\sum_bw_{ab}T_{bi},\forall a,i
\end{eqnarray*}

On the other hand, we have as well the following equivalences:
\begin{eqnarray*}
T\in Fix(v\otimes\bar{w})
&\iff&(v\otimes\bar{w})T=\xi\\
&\iff&\sum_{jb}v_{ij}w_{ab}^*T_{bj}=T_{ai}\forall a,i
\end{eqnarray*}

But with this in hand, both inclusions follow from the unitarity of $v,w$.
\end{proof}

Good news, we can now formulate our third Peter-Weyl theorem, as follows:

\index{Peter-Weyl}

\begin{theorem}[PW3]
The dense subalgebra $\mathcal C(G)\subset C(G)$ generated by the coefficients of the fundamental representation decomposes as a direct sum 
$$\mathcal C(G)=\bigoplus_{v\in Irr(G)}M_{\dim(v)}(\mathbb C)$$
with the summands being pairwise orthogonal with respect to $<f,g>=\int_Gf\bar{g}$.
\end{theorem}

\begin{proof}
By combining the previous two Peter-Weyl results, we deduce that we have a linear space decomposition as follows:
$$\mathcal C(G)
=\sum_{v\in Irr(G)}C_v
=\sum_{v\in Irr(G)}M_{\dim(v)}(\mathbb C)$$

Thus, in order to conclude, it is enough to prove that for any two irreducible representations $v,w\in Irr(G)$, the corresponding spaces of coefficients are orthogonal:
$$v\not\sim w\implies C_v\perp C_w$$ 

But this follows by Frobenius duality, by integrating. Let us set indeed:
$$P_{ia,jb}=\int_Gv_{ij}\bar{w}_{ab}$$

Then $P$ is the orthogonal projection onto the following vector space:
$$Fix(v\otimes\bar{w})
\simeq Hom(v,w)
=\{0\}$$

Thus we have $P=0$, and this gives the result.
\end{proof}

Finally, we have the following result, completing the Peter-Weyl theory:

\index{Peter-Weyl}

\begin{theorem}[PW4]
The characters of irreducible representations belong to
$$\mathcal C(G)_{central}=\left\{f\in\mathcal C(G)\Big|f(gh)=f(hg),\forall g,h\in G\right\}$$
called algebra of central functions on $G$, and form an orthonormal basis of it.
\end{theorem}

\begin{proof}
Observe first that $\mathcal C(G)_{central}$ is indeed an algebra, which contains all the characters. Conversely, consider a function $f\in\mathcal C(G)$, written as follows:
$$f=\sum_{v\in Irr(G)}f_v$$

The condition $f\in\mathcal C(G)_{central}$ states then that for any $v\in Irr(G)$, we must have:
$$f_v\in\mathcal C(G)_{central}$$

But this means that $f_v$ must be a scalar multiple of $\chi_v$, so the characters form a basis of $\mathcal C(G)_{central}$, as stated. Also, the fact that we have an orthogonal basis follows from Theorem 4.12. As for the fact that the characters have norm 1, this follows from:
$$\int_G\chi_v\bar{\chi}_v
=\sum_{ij}\int_Gv_{ii}\bar{v}_{jj}
=1$$

Here we have used the fact, coming from Frobenius duality, that the various integrals $\int_Gv_{ij}\bar{v}_{kl}$ form altogether the orthogonal projection onto the following vector space:
$$Fix(v\otimes\bar{v})\simeq End(v)=\mathbb C1$$

Thus, good news, we just proved all the Peter-Weyl theorems.
\end{proof}

Getting back now to our usual business, probability, we have, as an application:

\begin{theorem}
Given $G\subset_uU_N$, the moments of the main character are
$$\int_G\chi^k=\dim\left(Fix(u^{\otimes k})\right)$$
whose computation needs a basis for $Fix(u^{\otimes k})$. Alternatively, we have
$$\int_G\chi^k=\#\left\{1\in u^{\otimes k}\right\}$$
whose computation needs decomposing $u^{\otimes k}$ into irreducibles.
\end{theorem}

\begin{proof}
This is indeed something self-explanatory, with the first assertion coming from the formula $(\int_Gv_{ij})_{ij}=P$ from Theorem 4.10, with $v=u^{\otimes k}$, by applying the trace, and with the second assertion coming from the Peter-Weyl theory.
\end{proof}

As a comment now, among the above two methods for computing the law of $\chi$, the first one is in fact better, because a suitable modification of the argument leading to it will allow us to compute, more generally, the laws of the truncated characters $\chi_t$. We will discuss this in a moment, directly in the case which is of main interest, the easy one. 

\section*{4b. Tannaka, easiness}

Motivated by Theorem 4.14, given a closed subgroup $G\subset_uU_N$, we would like to compute the spaces $F_k=Fix(u^{\otimes k})$. And here, it is convenient to enlarge the attention to the collection of spaces $C_{kl}=Hom(u^{\otimes k},u^{\otimes l})$, which form a tensor category. Observe by the way that when doing this, our question does not become harder, because by Frobenius duality the collection of spaces $\{C_{kl}\}$ coincides with the collection of spaces $\{F_k\}$. 

\bigskip

Let us start our work with a straightforward categorical definition, as follows:

\index{tensor category}

\begin{definition}
A tensor category over $H=\mathbb C^N$ is a collection $C=(C_{kl})$ of linear spaces $C_{kl}\subset\mathcal L(H^{\otimes k},H^{\otimes l})$ satisfying the following conditions:
\begin{enumerate}
\item $S,T\in C$ implies $S\otimes T\in C$.

\item If $S,T\in C$ are composable, then $ST\in C$.

\item $T\in C$ implies $T^*\in C$.

\item Each $C_{kk}$ contains the identity operator.

\item $C_{\emptyset k}$ with $k=\circ\bullet,\bullet\circ$ contain the operator $R:1\to\sum_ie_i\otimes e_i$.

\item $C_{kl,lk}$ with $k,l=\circ,\bullet$ contain the flip operator $\Sigma:a\otimes b\to b\otimes a$.
\end{enumerate}
\end{definition}

Here, as usual, the tensor powers $H^{\otimes k}$, which are Hilbert spaces depending on a colored integer $k=\circ\bullet\bullet\circ\ldots\,$, are defined by the following formulae, and multiplicativity:
$$H^{\otimes\emptyset}=\mathbb C\quad,\quad 
H^{\otimes\circ}=H\quad,\quad
H^{\otimes\bullet}=\bar{H}\simeq H$$

We have already met such categories, when dealing with the Tannakian categories of the closed subgroups $G\subset U_N$, and our knowledge can be summarized as follows:

\index{Tannakian category}

\begin{proposition}
Given a closed subgroup $G\subset_uU_N$, its Tannakian category
$$C=\Big(Hom(u^{\otimes k},u^{\otimes l})\Big)_{kl}$$
is a tensor category over $H=\mathbb C^N$. Conversely, given a tensor category $C$ over $\mathbb C^N$,
$$G=\left\{g\in U_N\Big|Tg^{\otimes k}=g^{\otimes l}T,\forall k,l,\forall T\in C_{kl}\right\}$$
is a closed subgroup of $U_N$.
\end{proposition}

\begin{proof}
This is something that we basically know, the idea being as follows:

\medskip

(1) Regarding the first assertion, we have to check here the axioms (1-6) in Definition 4.15. The axioms (1-4) being all clear from definitions, let us establish (5). But this follows from the fact that each element $g\in G$ is a unitary, which can be reformulated as follows, with $R:1\to\sum_ie_i\otimes e_i$ being the map in Definition 4.15:
$$R\in Hom(1,g\otimes\bar{g})\quad,\quad 
R\in Hom(1,\bar{g}\otimes g)$$

Regarding now the condition in Definition 4.15 (6), this comes from the fact that the matrix coefficients $g\to g_{ij}$ and their conjugates $g\to\bar{g}_{ij}$ commute with each other.

\medskip

(2) Regarding the second assertion, we have to check that the subset $G\subset U_N$ constructed in the statement is a closed subgroup. But this is clear from definitions.
\end{proof}

Summarizing, we have so far precise axioms for the tensor categories $C=(C_{kl})$, given in Definition 4.15, as well as correspondences as follows:
$$G\to C_G\quad,\quad 
C\to G_C$$

We will prove in what follows that these correspondences are inverse to each other. In order to get started, we first have the following technical result:

\begin{proposition}
Consider the following conditions:
\begin{enumerate}
\item $C=C_{G_C}$, for any tensor category $C$.

\item $G=G_{C_G}$, for any closed subgroup $G\subset U_N$.
\end{enumerate}
We have then $(1)\implies(2)$. Also, $C\subset C_{G_C}$ is automatic.
\end{proposition}

\begin{proof}
Given $G\subset U_N$, we have $G\subset G_{C_G}$. On the other hand, by using (1) we have $C_G=C_{G_{C_G}}$. Thus, we have an inclusion of closed subgroups of $U_N$, which becomes an isomorphism at the level of the associated Tannakian categories, so $G=G_{C_G}$. Finally, the fact that we have an inclusion $C\subset C_{G_C}$ is clear from definitions.
\end{proof}

The point now is that it is possible to prove that we have $C_{G_C}\subset C$, by doing some abstract algebra, and we are led in this way to the following conclusion:

\begin{theorem}
The Tannakian duality constructions 
$$C\to G_C\quad,\quad 
G\to C_G$$
are inverse to each other.
\end{theorem}

\begin{proof}
This is something quite tricky, the idea being as follows:

\medskip

(1) According to Proposition 4.17, we must prove $C_{G_C}\subset C$. For this purpose, given a tensor category $C=(C_{kl})$ over a Hilbert space $H$, consider the following $*$-algebra:
$$E_C
=\bigoplus_{k,l}C_{kl}
\subset\bigoplus_{k,l}B(H^{\otimes k},H^{\otimes l})
\subset B\left(\bigoplus_kH^{\otimes k}\right)$$

Consider also, inside this $*$-algebra, the following $*$-subalgebra:
$$E_C^{(s)}
=\bigoplus_{|k|,|l|\leq s}C_{kl}
\subset\bigoplus_{|k|,|l|\leq s}B(H^{\otimes k},H^{\otimes l})
=B\left(\bigoplus_{|k|\leq s}H^{\otimes k}\right)$$

(2) It is then routine to check that we have equivalences as follows:
\begin{eqnarray*}
C_{G_C}\subset C
&\iff&E_{C_{G_C}}\subset E_C\\
&\iff&E_{C_{G_C}}^{(s)}\subset E_C^{(s)},\forall s\\
&\iff&E_{C_{G_C}}^{(s)'}\supset E_C^{(s)'},\forall s
\end{eqnarray*}

(3) Summarizing, we would like to prove that we have inclusions $E_C^{(s)'}\subset E_{C_{G_C}}^{(s)'}$. But this can be done by doing some algebra, as explained by Malacarne in \cite{mal}. For an introduction to all this, you can have a look at any of my algebra books.
\end{proof}

With this discussed, let us develop now a second idea which is needed, namely delinearization and Brauer theorems. Following \cite{bsp}, we have the following key definition:

\index{category of partitions}

\begin{definition}
A category of partitions is a collection $D=\bigsqcup_{k,l}D(k,l)$ of subsets $D(k,l)\subset P(k,l)$, having the following properties:
\begin{enumerate}
\item Stability under the horizontal concatenation, $(\pi,\sigma)\to[\pi\sigma]$.

\item Stability under vertical concatenation $(\pi,\sigma)\to[^\sigma_\pi]$, with matching middle symbols.

\item Stability under the upside-down turning $*$, with switching of colors, $\circ\leftrightarrow\bullet$.

\item Each set $P(k,k)$ contains the identity partition $||\ldots||$.

\item The sets $P(\emptyset,\circ\bullet)$ and $P(\emptyset,\bullet\circ)$ both contain the semicircle $\cap$.

\item The sets $P(k,\bar{k})$ with $|k|=2$ contain the crossing partition $\slash\hskip-2.0mm\backslash$.
\end{enumerate}
\end{definition} 

Observe the similarity with Definition 4.15, and more on this in a moment. Also, at the level of the examples, there are plenty of them, that we already met, namely $P$ itself, then $NC$, then $P_2$ and $NC_2$, and then $\mathcal P_2,\mathcal{NC}_2$ and more. We will be back to this, examples, a bit later, when discussing classification results for such categories.

\bigskip

In order to construct a Tannakian category out of such a category, we will need:

\begin{proposition}
Each partition $\pi\in P(k,l)$ produces a linear map
$$T_\pi:(\mathbb C^N)^{\otimes k}\to(\mathbb C^N)^{\otimes l}$$
given by the following formula, with $e_1,\ldots,e_N$ being the standard basis of $\mathbb C^N$,
$$T_\pi(e_{i_1}\otimes\ldots\otimes e_{i_k})=\sum_{j_1\ldots j_l}\delta_\pi\begin{pmatrix}i_1&\ldots&i_k\\ j_1&\ldots&j_l\end{pmatrix}e_{j_1}\otimes\ldots\otimes e_{j_l}$$
and with the Kronecker type symbols $\delta_\pi\in\{0,1\}$ depending on whether the indices fit or not. The assignement $\pi\to T_\pi$ is categorical, in the sense that we have
$$T_\pi\otimes T_\sigma=T_{[\pi\sigma]}\quad,\quad 
T_\pi T_\sigma=N^{c(\pi,\sigma)}T_{[^\sigma_\pi]}\quad,\quad 
T_\pi^*=T_{\pi^*}$$
where $c(\pi,\sigma)$ are certain integers, coming from the erased components in the middle.
\end{proposition}

\begin{proof}
The concatenation axiom comes from the following computation:
\begin{eqnarray*}
&&(T_\pi\otimes T_\sigma)(e_{i_1}\otimes\ldots\otimes e_{i_p}\otimes e_{k_1}\otimes\ldots\otimes e_{k_r})\\
&=&\sum_{j_1\ldots j_q}\sum_{l_1\ldots l_s}\delta_\pi\begin{pmatrix}i_1&\ldots&i_p\\j_1&\ldots&j_q\end{pmatrix}\delta_\sigma\begin{pmatrix}k_1&\ldots&k_r\\l_1&\ldots&l_s\end{pmatrix}e_{j_1}\otimes\ldots\otimes e_{j_q}\otimes e_{l_1}\otimes\ldots\otimes e_{l_s}\\
&=&\sum_{j_1\ldots j_q}\sum_{l_1\ldots l_s}\delta_{[\pi\sigma]}\begin{pmatrix}i_1&\ldots&i_p&k_1&\ldots&k_r\\j_1&\ldots&j_q&l_1&\ldots&l_s\end{pmatrix}e_{j_1}\otimes\ldots\otimes e_{j_q}\otimes e_{l_1}\otimes\ldots\otimes e_{l_s}\\
&=&T_{[\pi\sigma]}(e_{i_1}\otimes\ldots\otimes e_{i_p}\otimes e_{k_1}\otimes\ldots\otimes e_{k_r})
\end{eqnarray*}

The composition axiom comes from the following computation:
\begin{eqnarray*}
T_\pi T_\sigma(e_{i_1}\otimes\ldots\otimes e_{i_p})
&=&\sum_{j_1\ldots j_q}\delta_\sigma\begin{pmatrix}i_1&\ldots&i_p\\j_1&\ldots&j_q\end{pmatrix}
\sum_{k_1\ldots k_r}\delta_\pi\begin{pmatrix}j_1&\ldots&j_q\\k_1&\ldots&k_r\end{pmatrix}e_{k_1}\otimes\ldots\otimes e_{k_r}\\
&=&\sum_{k_1\ldots k_r}N^{c(\pi,\sigma)}\delta_{[^\sigma_\pi]}\begin{pmatrix}i_1&\ldots&i_p\\k_1&\ldots&k_r\end{pmatrix}e_{k_1}\otimes\ldots\otimes e_{k_r}\\
&=&N^{c(\pi,\sigma)}T_{[^\sigma_\pi]}(e_{i_1}\otimes\ldots\otimes e_{i_p})
\end{eqnarray*}

Finally, the involution axiom comes from the following computation:
\begin{eqnarray*}
T_\pi^*(e_{j_1}\otimes\ldots\otimes e_{j_q})
&=&\sum_{i_1\ldots i_p}<T_\pi^*(e_{j_1}\otimes\ldots\otimes e_{j_q}),e_{i_1}\otimes\ldots\otimes e_{i_p}>e_{i_1}\otimes\ldots\otimes e_{i_p}\\
&=&\sum_{i_1\ldots i_p}\delta_\pi\begin{pmatrix}i_1&\ldots&i_p\\ j_1&\ldots& j_q\end{pmatrix}e_{i_1}\otimes\ldots\otimes e_{i_p}\\
&=&T_{\pi^*}(e_{j_1}\otimes\ldots\otimes e_{j_q})
\end{eqnarray*}

Summarizing, our correspondence is indeed categorical.
\end{proof}

Following \cite{bsp}, we can now formulate a key theoretical result, as follows:

\index{easy group}

\begin{theorem}
Any category of partitions $D\subset P$ produces a series of compact groups $G=(G_N)$, with $G_N\subset U_N$ for any $N\in\mathbb N$, via the formula
$$Hom(u^{\otimes k},u^{\otimes l})=span\left(T_\pi\Big|\pi\in D(k,l)\right)$$
for any $k,l$, and Tannakian duality. We call such groups easy.
\end{theorem}

\begin{proof}
Indeed, once we fix an integer $N\in\mathbb N$, the various axioms in Definition 4.19 show, via Proposition 4.20, that the following spaces form a Tannakian category:
$$C_{kl}=span\left(T_\pi\Big|\pi\in D(k,l)\right)$$

Thus, Tannakian duality applies, and provides us with a closed subgroup $G_N\subset U_N$ such that the following equalities are satisfied, for any colored integers $k,l$:
$$Hom(u^{\otimes k},u^{\otimes l})=C_{kl}$$

We are therefore led to the conclusion in the statement.
\end{proof}

And with this, good news, done with the general theory. At the level of examples, and generalizing what we previously knew from chapter 1, we have:

\index{bistochastic group}

\begin{theorem}
We have easy groups as follows, with $H_N=H_N^2,K_N=H_N^\infty$, and $B_N\subset O_N,C_N\subset U_N$ being the bistochastic groups, with sum $1$ on each row and column:
$$\xymatrix@R=18pt@C=14pt{
&K_N\ar[rr]&&U_N\\
H_N\ar[rr]\ar[ur]&&O_N\ar[ur]\\
&S_N\ar[rr]\ar[uu]&&C_N\ar[uu]\\
S_N\ar[uu]\ar[ur]\ar[rr]&&B_N\ar[uu]\ar[ur]
}\qquad\xymatrix@R=18pt@C=14pt{\\ \\ :}
\qquad\xymatrix@R=18pt@C=10pt{
&\mathcal P_{even}\ar[dl]\ar[dd]&&\mathcal P_2\ar[dl]\ar[ll]\ar[dd]\\
P_{even}\ar[dd]&&P_2\ar[dd]\ar[ll]\\
&P\ar[dl]&&\mathcal P_{12}\ar[dl]\ar[ll]\\
P&&P_{12}\ar[ll]}$$
Moreover, the reflection groups $H_N^s$ with $s\in\{1,2,\ldots,\infty\}$ fit on the diagonal of the left face of the left cube, corresponding to the categories $P^s$, fitting on the right cube. 
\end{theorem}

\begin{proof}
This is something very standard, based on Theorem 4.21, as follows:

\medskip

(1) We certainly have the groups on the left, and we must prove that these are easy, coming from the categories on the right, with ``even'' standing for having blocks of even size, $12$ standing for singletons and pairings, calligraphic standing for matching, meaning satisfying $\#\circ=\#\bullet$, as a weighted equality, in each block, and $P^s$ being the category of partitions satisfying $\#\circ=\#\bullet(s)$, as a weighted equality, in each block.

\medskip

(2) Regarding now the proof, as mentioned in chapter 1, for the symmetric group $S_N$, and more generally for the reflection groups $H_N^s$, this is something elementary, coming from a direct computation of the spaces $Hom(u^{\otimes k},u^{\otimes l})$. However, for $O_N,U_N$ and their bistochastic versions $B_N,C_N$, the result, due to Brauer, is something more subtle, requiring the use of Theorem 4.21, or something equivalent, meaning advanced technology. 

\medskip

(3) In short, we have a mixture of trivial and non-trivial facts, and the best in practice is to prove everything, both the trivial and non-trivial results, using Theorem 4.21. Getting started now, the unitary group $U_N$ being defined via the relations $u^*=u^{-1}$, $u^t=\bar{u}^{-1}$, the associated Tannakian category is $C=span(T_\pi|\pi\in D)$, with:
$$D
=<{\ }^{\,\cap}_{\circ\bullet}\,\,,{\ }^{\,\cap}_{\bullet\circ}>
=\mathcal P_2$$

(4) The orthogonal group $O_N\subset U_N$ being defined by imposing the relations $u_{ij}=\bar{u}_{ij}$, the associated Tannakian category is $C=span(T_\pi|\pi\in D)$, with:
$$D
=<\mathcal P_2,|^{\hskip-1.32mm\circ}_{\hskip-1.32mm\bullet},|_{\hskip-1.32mm\circ}^{\hskip-1.32mm\bullet}>
=P_2$$

(5) The bistochastic group $C_N\subset U_N$ being defined by imposing the relations $u\xi=\xi$, $\bar{u}\xi=\xi$, with $\xi=(1)_i$, the associated Tannakian category is $C=span(T_\pi|\pi\in D)$, with:
$$D=<\mathcal P_2,|_{\hskip-1.32mm\circ},|_{\hskip-1.32mm\bullet}>
=\mathcal P_{12}$$
  
(6) The orthogonal bistochastic group appears as $B_N=C_N\cap O_N$, and we conclude that the associated Tannakian category is $C=span(T_\pi|\pi\in D)$, with:
$$D=<\mathcal P_{12},P_2>=P_{12}$$

(7) In order to discuss now $S_N$, consider the one-block ``fork'' partition, namely:
$$\xymatrix@R=1mm@C=2mm{\\ \\ \mu\ \ =\\ \\ }\ \ \ 
\xymatrix@R=2mm@C=3mm{
\circ\ar@/_/@{-}[dr]&&\circ\\
&\ar@/_/@{-}[ur]\ar@{-}[dd]\\
&&&\\
&\circ}$$

We have then $T_\mu(e_i\otimes e_j)=\delta_{ij}e_i$, and by using this formula, we obtain:
$$T_\mu\in Hom(u^{\otimes 2},u)\iff u_{ij}u_{ik}=\delta_{jk}u_{ij},\forall i,j,k$$

But on the right we have the relations defining $S_N\subset O_N$, and we conclude that $S_N$ is indeed easy, coming from the following category of partitions:
$$D=<\mu>=P$$

(8) In order to discuss $H_N$, consider the following one-block partition:
$$\xymatrix@R=0.5mm@C=2mm{\\ \\ \\ \chi\ \ =\\ \\ }\ \ \ 
\xymatrix@R=2mm@C=3mm{
\circ\ar@/_/@{-}[dr]&&\circ\\
&\ar@/_/@{-}[ur]\ar@{-}[dd]\\
&&&\\
&\ar@/^/@{-}[dr]\ar@/_/@{-}[dl]\\
\circ&&\circ}$$

We have then $T_\chi(e_i\otimes e_j)=\delta_{ij}e_i\otimes e_i$, and by using this formula, we obtain:
$$T_\chi\in End(u^{\otimes 2})\iff \delta_{ik}u_{ia}u_{ib}=\delta_{ab}u_{ia}u_{ka},\forall i,k,a,b$$

But on the right we have the relations defining $H_N\subset O_N$, and we conclude that $H_N$ is indeed easy, coming from the following category of partitions:
$$D=<\chi>=P_{even}$$

(9) In order to discuss now $K_N\subset U_N$, consider the following colored partition:
$$\xymatrix@R=0.5mm@C=2mm{\\ \\ \\ \chi'\ \ =\\ \\ }\ \ \ 
\xymatrix@R=2mm@C=3mm{
\circ\ar@/_/@{-}[dr]&&\bullet\\
&\ar@/_/@{-}[ur]\ar@{-}[dd]\\
&&&\\
&\ar@/^/@{-}[dr]\ar@/_/@{-}[dl]\\
\bullet&&\circ}$$

Our computations from the previous proof, for the group $H_N$, modify into:
$$T_{\chi'}\in Hom(u\otimes\bar{u},\bar{u}\otimes u)\iff \delta_{ik}u_{ia}\bar{u}_{ib}=\delta_{ab}\bar{u}_{ia}u_{ka},\forall i,k,a,b$$

But on the right we have the relations defining $K_N\subset U_N$, and we conclude that $K_N$ is indeed easy, coming from the following category of partitions:
$$D=<\chi'>=\mathcal P_{even}$$

(10) Finally, regarding $H_N^s\subset K_N$, consider the following partition, with $s+2$ legs:
$$\xymatrix@R=6pt@C=12pt{
&&\ar@{-}[dd]\ar@{-}[rrrrr]&\ar@{-}[dd]&&\ar@{-}[dd]&\ar@{-}[dd]&\ar@{-}[dd]\\
\xi&=&&&\ldots\\
&&\circ&\circ&&\circ&\circ&\bullet}$$

We have then $T_\xi=\sum_je_j^{\otimes s+2}$, and by using this formula, inside $K_N$, we obtain:
$$T_\xi\in Fix(u^{\otimes s+1}\otimes\bar{u})\iff\sum_ju_{ij}^{s+1}\bar{u}_{ij}=1$$

But on the right we have the relations defining $H_N^s\subset K_N$, and we conclude that $H_N^s$ is indeed easy, coming from the following category of partitions:
$$D=<\mathcal P_{even},\xi>=P^s$$

Summarizing, theorem proved, in exactly 2 pages. Tannakian duality rules.
\end{proof}

Now back to probability, in order to make use of Theorem 4.14, we will need:

\index{Lindst\"om formula}
\index{Gram determinant}

\begin{theorem}[Lindst\"om]
The Gram matrix of the vectors $\{T_\pi|\pi\in P(k)\}$ is
$$G_{kN}(\pi,\nu)=N^{|\pi\vee\nu|}$$
where $\vee$ is the superposition operation for the partitions, and we have the formula
$$\det(G_{kN})=\prod_{\pi\in P(k)}\frac{N!}{(N-|\pi|)!}$$
at $N\geq k$. In particular, the vectors $\{T_\pi|\pi\in P(k)\}$ are linearly independent at $N\geq k$.
\end{theorem}

\begin{proof}
Regarding the first formula, this comes from the following computation:
\begin{eqnarray*}
G_{kN}(\pi,\nu)
&=&\sum_{j_1\ldots j_l}
\delta_\pi\begin{pmatrix}j_1&\ldots& j_l\end{pmatrix}
\delta_\nu\begin{pmatrix}j_1&\ldots& j_l\end{pmatrix}\\
&=&\sum_{j_1\ldots j_l}
\delta_{\pi\vee\nu}\begin{pmatrix}j_1&\ldots& j_l\end{pmatrix}\\
&=&N^{|\pi\vee\nu|}
\end{eqnarray*}

Next, in order to compute $\det(G_{kN})$, observe that we have the following formula:
\begin{eqnarray*}
N^{|\pi\vee\nu|}
&=&\#\left\{i_1,\ldots,i_k\in\{1,\ldots,N\}\Big|\ker i\geq\pi\vee\nu\right\}\\
&=&\sum_{\tau\geq\pi\vee\nu}\#\left\{i_1,\ldots,i_k\in\{1,\ldots,N\}\Big|\ker i=\tau\right\}\\
&=&\sum_{\tau\geq\pi\vee\nu}N(N-1)\ldots(N-|\tau|+1)
\end{eqnarray*}

Now in view of this latter formula, consider the following matrix:
$$L_{kN}(\pi,\nu)=
\begin{cases}
N(N-1)\ldots(N-|\pi|+1)&{\rm if}\ \nu\leq\pi\\
0&{\rm otherwise}
\end{cases}$$

If we denote by $A_k$ the adjacency matrix of $P(k)$, we conclude that we have:
\begin{eqnarray*}
G_{kN}(\pi,\nu)
&=&\sum_{\tau\geq\pi}L_{kN}(\tau,\nu)\\
&=&\sum_\tau A_k(\pi,\tau)L_{kN}(\tau,\nu)\\
&=&(A_kL_{kN})(\pi,\nu)
\end{eqnarray*}

Summarizing, we have decomposed our Gram matrix as $G_{kN}=A_kL_{kN}$. Now if we order $P(k)$ as usual, with respect to the number of blocks, and then lexicographically, the matrix $A_k$ is upper triangular, and the matrix $L_{kN}$ is lower triangular, and so:
\begin{eqnarray*}
\det(G_{kN})
&=&\det(A_k)\det(L_{kN})\\
&=&\det(L_{kN})\\
&=&\prod_{\pi\in P(k)}\frac{N!}{(N-|\pi|)!}
\end{eqnarray*}

And with this done, we are led to the conclusions in the statement. 
\end{proof}

Now back to our main examples of easy groups, we have:

\begin{theorem}
For the main easy groups the asymptotic law of $\chi_1$ is
$$\xymatrix@R=18pt@C=14pt{
&K_N\ar[rr]&&U_N\\
H_N\ar[rr]\ar[ur]&&O_N\ar[ur]\\
&S_N\ar[rr]\ar[uu]&&C_N\ar[uu]\\
S_N\ar[uu]\ar[ur]\ar[rr]&&B_N\ar[uu]\ar[ur]
}\qquad\xymatrix@R=18pt@C=14pt{\\ \\ :}
\qquad\xymatrix@R=16.5pt@C=20pt{
&P_1\ar@{-}[rr]\ar@{-}[dd]&&G_1\ar@{-}[dd]\\
p_1^2\ar@{-}[rr]\ar@{-}[dd]\ar@{-}[ur]&&g_1\ar@{-}[dd]\ar@{-}[ur]\\
&p_1\ar@{-}[rr]\ar@{-}[uu]&&G_1^1\ar@{.}[uu]\\
p_1\ar@{-}[uu]\ar@{-}[ur]\ar@{-}[rr]&&g_1^1\ar@{-}[uu]\ar@{-}[ur]
}$$
and on the left we can insert $H_N^s$ with $s\in\{1,2,\ldots,\infty\}$, corresponding to $p^s_1$.
\end{theorem}

\begin{proof}
According to Theorems 4.14, 4.21 and 4.23, for an easy group $G=(G_N)$ coming from a category of partitions $D=(D(k,l))$, we have:
\begin{eqnarray*}
\lim_{N\to\infty}\int_{G_N}\chi^k
&=&\lim_{N\to\infty}\dim\left(Fix(u^{\otimes k})\right)\\
&=&\lim_{N\to\infty}\dim\left(span\left(T_\pi\Big|\pi\in D(k)\right)\right)\\
&=&|D(k)|
\end{eqnarray*}

Now for the groups in Theorem 4.22, this leads to the measures on the right, via the moment formulae for $g_1,G_1,p^s_1$ from chapter 1, with the notation $P_1=p^\infty_1$, and with the remaining results, involving the shifts $g_1^1,G_1^1$, being an easy exercise for you.
\end{proof}

\section*{4c. Weingarten formula}

Getting now towards truncated characters, things are a bit more tricky here, requiring some further general easiness theory. We first have the following key result:

\index{Weingarten formula}
\index{Gram matrix}
\index{Weingarten matrix}

\begin{theorem}
For an easy group $G\subset U_N$, coming from a category of partitions $D=(D(k,l))$, we have the Weingarten integration formula
$$\int_Gu_{i_1j_1}^{e_1}\ldots u_{i_kj_k}^{e_k}=\sum_{\pi,\nu\in D(k)}\delta_\pi(i)\delta_\nu(j)W_{kN}(\pi,\nu)$$
for any multi-indices $i,j$ and any exponent $k=e_1\ldots e_k$, where $D(k)=D(\emptyset,k)$, the $\delta$ numbers are the usual Kronecker type symbols, and $W_{kN}=G_{kN}^{-1}$, with 
$$G_{kN}(\pi,\nu)=N^{|\pi\vee\nu|}$$
where $|.|$ is the number of blocks. 
\end{theorem}

\begin{proof}
We know from Theorems 4.10 and 4.21 that the integrals in the statement form altogether the orthogonal projection $P$ onto the following space:
$$Fix(u^{\otimes k})=span\left(T_\pi\Big|\pi\in D(k)\right)$$

In order to compute this projection $P$, consider the following linear map:
$$E(x)=\sum_{\pi\in D(k)}<x,T_\pi>T_\pi$$

By a standard linear algebra computation, it follows that we have the following formula, where $W$ is the inverse of the restriction of $E$ to the space $Fix(u^{\otimes k})$: 
$$P=WE$$

But this latter restriction is the linear map given by the Gram matrix $G$ of the vectors $T_\pi$, given by $G(\pi,\nu)=N^{|\pi\vee\nu|}$, so we have $W=G^{-1}$, which gives the result.
\end{proof}

Next, in order to obtain good results from this, we will need to assume uniformity:

\index{uniform group}
\index{removing blocks}

\begin{proposition}
For an easy group $G=(G_N)$, coming from a category of partitions $D\subset P$, the following conditions are equivalent:
\begin{enumerate}
\item $G_{N-1}=G_N\cap U_{N-1}$, via the embedding $U_{N-1}\subset U_N$ given by $u\to diag(u,1)$.

\item $G_{N-1}=G_N\cap U_{N-1}$, via the $N$ possible diagonal embeddings $U_{N-1}\subset U_N$.

\item $D$ is stable under the operation which consists in removing blocks.
\end{enumerate}
If these conditions are satisfied, we say that $G=(G_N)$ is uniform.
\end{proposition}

\begin{proof}
We use the general easiness theory developed in the above:

\medskip

$(1)\iff(2)$ This is something standard, coming from the inclusion $S_N\subset G_N$, which makes everything $S_N$-invariant. The result follows as well from the proof of $(1)\iff(3)$ below, which can be converted into a proof of $(2)\iff(3)$, in the obvious way.

\medskip

$(1)\iff(3)$ Given a subgroup $K\subset U_{N-1}$, with fundamental representation $u$, consider the $N\times N$ matrix $v=diag(u,1)$. Our claim is that for $\pi\in P(k)$ we have:
$$\xi_\pi\in Fix(v^{\otimes k})\iff\xi_{\pi'}\in Fix(v^{\otimes k'}),\,\forall\pi'\in P(k'),\pi'\subset\pi$$

In order to prove this, we must study the condition on the left. We have:
\begin{eqnarray*}
\xi_\pi\in Fix(v^{\otimes k})
&\iff&(v^{\otimes k}\xi_\pi)_{i_1\ldots i_k}=(\xi_\pi)_{i_1\ldots i_k},\forall i\\
&\iff&\sum_j(v^{\otimes k})_{i_1\ldots i_k,j_1\ldots j_k}(\xi_\pi)_{j_1\ldots j_k}=(\xi_\pi)_{i_1\ldots i_k},\forall i\\
&\iff&\sum_j\delta_\pi(j_1,\ldots,j_k)v_{i_1j_1}\ldots v_{i_kj_k}=\delta_\pi(i_1,\ldots,i_k),\forall i
\end{eqnarray*}

Now let us recall that our representation has the special form $v=diag(u,1)$. We conclude from this that for any index $a\in\{1,\ldots,k\}$, we must have:
$$i_a=N\implies j_a=N$$

With this observation in hand, if we denote by $i',j'$ the multi-indices obtained from $i,j$ obtained by erasing all the above $i_a=j_a=N$ values, and by $k'\leq k$ the common length of these new multi-indices, our condition becomes:
$$\sum_{j'}\delta_\pi(j_1,\ldots,j_k)(v^{\otimes k'})_{i'j'}=\delta_\pi(i_1,\ldots,i_k),\forall i$$

Here the index $j$ is by definition obtained from $j'$ by filling with $N$ values. In order to finish now, we have two cases, depending on $i$, as follows:

\medskip

\underline{Case 1}. Assume that the index set $\{a|i_a=N\}$ corresponds to a certain subpartition $\pi'\subset\pi$. In this case, the $N$ values will not matter, and our formula becomes:
$$\sum_{j'}\delta_\pi(j'_1,\ldots,j'_{k'})(v^{\otimes k'})_{i'j'}=\delta_\pi(i'_1,\ldots,i'_{k'})$$

\underline{Case 2}. Assume now the opposite, namely that the set $\{a|i_a=N\}$ does not correspond to a subpartition $\pi'\subset\pi$. In this case the indices mix, and our formula reads:
$$0=0$$

Thus, we are led to $\xi_{\pi'}\in Fix(v^{\otimes k'})$, for any subpartition $\pi'\subset\pi$, as claimed. Thus our claim is proved, and with this in hand, the result follows from Tannakian duality.
\end{proof}

We can now formulate a nice probabilistic result, as follows:

\begin{theorem}
For a uniform easy group $G=(G_N)$, we have the formula
$$\lim_{N\to\infty}\int_{G_N}\chi_t^k=\sum_{\pi\in D(k)}t^{|\pi|}$$
with $D\subset P$ being the corresponding category of partitions.
\end{theorem}

\begin{proof}
By using the Weingarten formula from Theorem 4.25, we obtain:
\begin{eqnarray*}
\int_{G_N}(u_{11}+\ldots +u_{ss})^k
&=&\sum_{i_1=1}^{s}\ldots\sum_{i_k=1}^s\int_{G_N}u_{i_1i_1}\ldots u_{i_ki_k}\\
&=&\sum_{\pi,\nu\in D(k)}W_{kN}(\pi,\nu)\sum_{i_1=1}^{s}\ldots\sum_{i_k=1}^s\delta_\pi(i)\delta_\nu(i)\\
&=&\sum_{\pi,\nu\in D(k)}W_{kN}(\pi,\nu)G_{ks}(\nu,\pi)\\
&=&Tr(W_{kN}G_{ks})
\end{eqnarray*}

The point now is that in the uniform case the Gram matrix is asymptotically diagonal, $G_{kN}\simeq diag(N^{|\pi|})$, with this being something clear for all the examples, and true in general too, either via an abstract argument, or via the classification results in \cite{twe}. Now the inverse of $G_{kN}$ being asymptotically diagonal too, $W_{kN}\simeq diag(N^{-|\pi|})$, we obtain:
$$\lim_{N\to\infty}\int_{G_N}\chi_t^k=\sum_{\pi\in D(k)}t^{|\pi|}$$

Thus, we have indeed the asymptotic moment formula in the statement.
\end{proof}

At the level of the examples, upgrading our previous knowledge, we have:

\begin{theorem}
For the main easy groups the asymptotic law of $\chi_t$ is
$$\xymatrix@R=18pt@C=14pt{
&K_N\ar[rr]&&U_N\\
H_N\ar[rr]\ar[ur]&&O_N\ar[ur]\\
&S_N\ar[rr]\ar[uu]&&C_N\ar[uu]\\
S_N\ar[uu]\ar[ur]\ar[rr]&&B_N\ar[uu]\ar[ur]
}\qquad\xymatrix@R=18pt@C=14pt{\\ \\ :}
\qquad\xymatrix@R=16.5pt@C=20pt{
&P_t\ar@{-}[rr]\ar@{-}[dd]&&G_t\ar@{-}[dd]\\
p_t^2\ar@{-}[rr]\ar@{-}[dd]\ar@{-}[ur]&&g_t\ar@{-}[dd]\ar@{-}[ur]\\
&p_t\ar@{-}[rr]\ar@{-}[uu]&&G_t^t\ar@{.}[uu]\\
p_t\ar@{-}[uu]\ar@{-}[ur]\ar@{-}[rr]&&g_t^t\ar@{-}[uu]\ar@{-}[ur]
}$$
and on the left we can insert $H_N^s$ with $s\in\{1,2,\ldots,\infty\}$, corresponding to $p^s_t$.
\end{theorem}

\begin{proof}
According to Theorem 4.27, for the easy groups in Theorem 4.22, which are all uniform, the asymptotic laws of truncated characters $\chi_t$ are given by:
$$\lim_{N\to\infty}\int_{G_N}\chi_t^k=\sum_{\pi\in D(k)}t^{|\pi|}$$

But this leads to the laws on the right, via the moment formulae for $g_t,G_t,p^s_t$ from chapter 1, with the notation $P_t=p^\infty_t$, and with the remaining results, involving the shifts $g_t^t,G_t^t$, being an easy exercise. To be more precise, regarding these latter exercises:

\medskip

(1) According to Theorem 4.22, for the unitary bistochastic group $C_N$, we have:
$$\lim_{N\to\infty}\int_{C_N}\chi_t^k=\sum_{\pi\in\mathcal P_{12}(k)}t^{|\pi|}
=\sum_{l\subset_u k}t^{|k|-|l|/2}(|l|/2)!$$

To be more precise here, in all this $k=\circ\bullet\bullet\circ\ldots$ is a colored integer, and on the right $l\subset_uk$ stands for the choice of a colored subinteger, which is uniform in the sense of chapter 5, meaning containing the same number of $\circ$ and $\bullet$ symbols. This uniform subinteger $l\subset_uk$ corresponds to the matching pairings inside $\pi$, and there are $(|l|/2)!$ choices of such matching pairings, and the $|k|-|l|$ points left are singletons.

\medskip

(2) On the other hand, according to the formulae in chapter 1, we have:
\begin{eqnarray*}
M_k(G_t^a)
&=&\frac{1}{\pi t}\int_\mathbb Cz^ke^{-|z-a|^2/t}dz\\
&=&\frac{1}{\pi t}\int_\mathbb C(z+a)^ke^{-|z|^2/t}dz\\
&=&\frac{1}{\pi t}\int_\mathbb C\sum_{l\subset k}z^la^{k-l}e^{-|z|^2/t}dz\\
&=&\sum_{l\subset k}a^{k-l}\frac{1}{\pi t}\int_\mathbb Cz^le^{-|z|^2/t}dz\\
&=&\sum_{l\subset k}a^{k-l}M_l(g_t)
\end{eqnarray*}

But we know from chapter 1 that $M_l(g_t)$ vanishes when $l$ is not uniform, and is given by $M_l(g_t)=t^{|l|/2}(|l|/2)!$ in the uniform case. Thus, our moment formula reads:
$$M_k(G_t^a)=\sum_{l\subset_u k}a^{k-l}t^{|l|/2}(|l|/2)!$$

(3) Now in the case where $a$ is a real number, the quantity $a^{k-l}$ is a usual power of $a$, the exponent being the length $|k-l|=|k|-|l|$, and our formula becomes:
$$M_k(G_t^a)=\sum_{l\subset_u k}a^{|k|-|l|}t^{|l|/2}(|l|/2)!$$

But with $a=t$ we obtain precisely the moment formula in (1), and we are done with $C_N$. As for the group $B_N$, here things are even simpler, say easy exercise for you.
\end{proof}

\section*{4d. Clebsch-Gordan}

Time for some physics, along the lines evoked in the beginning of this chapter. We will be talking about the Clebsch-Gordan rules, which are essential to the functioning of both our beloved pure mathematics, and our beloved Standard Model. We first have:

\begin{theorem}
In the standard picture of the group $SU_2$, which is as follows,
$$SU_2=\left\{\begin{pmatrix}x+iy&z+it\\ -z+it&x-iy\end{pmatrix}\Big|(x,y,z,t)\in S^3_\mathbb R\right\}$$ 
the main character is $\chi=2x$, following the Wigner semicircle law $\gamma_1$.
\end{theorem}

\begin{proof}
The starting point is the defining formula for the group $SU_2$, as being the group of $2\times2$ unitary matrices having determinant one:
$$SU_2=\left\{U\in M_2(\mathbb C)\Big|U^*=U^{-1},\,\det U=1\right\}$$

In practice, if we pick $U=\binom{a\ b}{c\ d}$ of determinant 1, the equation $U^*=U^{-1}$ reads:
$$\begin{pmatrix}\bar{a}&\bar{c}\\ \bar{b}&\bar{d}\end{pmatrix}
=\begin{pmatrix}d&-b\\ -c&a\end{pmatrix}$$

Thus we must have $d=\bar{a}$, $c=-\bar{b}$, and since with these conditions imposed, the condition $\det U=1$ reads $|a|^2+|b|^2=1$, we are led to the following conclusion:
$$SU_2=\left\{\begin{pmatrix}a&b\\ -\bar{b}&\bar{a}\end{pmatrix}\Big|\ |a|^2+|b|^2=1\right\}$$

But with $a=x+iy$, $b=z+it$ this gives the formula in the statement. Now observe that in that picture, the main character is given by the following formula:
$$\chi\begin{pmatrix}x+iy&z+it\\ -z+it&x-iy\end{pmatrix}=2x$$

Moreover, the uniform measure on $SU_2$ corresponds to the uniform measure on $S^3_\mathbb R$, so we can finish by using our spherical integral knowledge from chapter 1, as follows:
$$\int_{S^3_\mathbb R}x^{2k}
=\frac{3!!(2k)!!}{(2k+3)!!}
=\frac{1}{4^k}\cdot\frac{1}{k+1}\binom{2k}{k}
=\frac{C_k}{4^k}$$

Thus the variable $2x\in C(S^3_\mathbb R)$ follows the Wigner semicircle law $\gamma_1$, as claimed.
\end{proof}

Regarding the group $SO_3$, which is the close friend of $SU_2$, we have:

\index{Euler-Rodrigues}

\begin{theorem}
In the standard picture of the elements $U\in SO_3$, namely
$$U=\begin{pmatrix}
x^2+y^2-z^2-t^2&2(yz-xt)&2(xz+yt)\\
2(xt+yz)&x^2+z^2-y^2-t^2&2(zt-xy)\\
2(yt-xz)&2(xy+zt)&x^2+t^2-y^2-z^2
\end{pmatrix}$$
with $(x,y,z,t)\in S^3_\mathbb R$, the variable $\chi+1=4x^2$ follows the Marchenko-Pastur law $\pi_1$.
\end{theorem}

\begin{proof}
This is something quite tricky, the idea being as follows:

\medskip

(1) To start with, following Euler, any rotation $U\in SO_3$ has a rotation axis. Indeed, we have the following computation, using some linear algebra magic:
\begin{eqnarray*}
\det(U-1)
&=&\det(U^t-1)\\
&=&\det(U^t(1-U))\\
&=&\det(U^t)\det(1-U)\\
&=&\det(1-U)
\end{eqnarray*}

Thus $\det(U-1)=0$, which tells us that $U$ must have a $1$-eigenvector, as desired. 

\medskip

(2) Next, following Euler-Rodrigues, we have the formula in the statement, with the parameters $(x,y,z,t)\in S^3_\mathbb R$ encapsulating the direction of the rotation axis, and the rotation angle. Alternatively, if you want to trade heavy trigonometry for heavy algebra, the Euler-Rodrigues formula comes from the double cover map $SU_2\to SO_3$, obtained via the adjoint action of $SU_2$ on $\mathbb C^2\simeq\mathbb R^4$, by removing the trivial fixed direction.

\medskip

(3) In short, quite complicated all this, and with the complete picture, needed for properly understanding the double cover map $SU_2\to SO_3$, involving as well the 4 Pauli spin matrices, appearing as matrix coefficients of $x,y,z,t$ in the formula from Theorem 4.29. Have a look at any of my algebra books, where all this is explained.

\medskip

(4) With this discussed, in the Euler-Rodrigues picture, the main character is:
\begin{eqnarray*}
\chi
&=&3x^2-y^2-z^2-t^2\\
&=&4x^2-(x^2+y^2+z^2+t^2)\\
&=&4x^2-1
\end{eqnarray*}

Thus $\chi+1=4x^2$, which by Theorem 4.29 follows the Marchenko-Pastur law $\pi_1$.
\end{proof}

As a continuation of the above, we know from Peter-Weyl that the law of the main character ultimately comes by decomposing the representations $u^{\otimes k}$ into irreducibles. So, let us try to understand this. For the group $SU_2$, we have the following result:

\index{Clebsch-Gordan}
\index{fusion rules}

\begin{theorem}
The irreducible representations of $SU_2$ are all self-adjoint, and can be labeled by positive integers, with their fusion rules being as follows,
$$r_k\otimes r_l=r_{|k-l|}+r_{|k-l|+2}+\ldots+r_{k+l}$$
called Clebsch-Gordan rules. The corresponding dimensions are $\dim r_k=k+1$.
\end{theorem}

\begin{proof}
With the knowledge that we have, we can prove this result as follows:

\medskip

(1) According to Theorem 4.29, the moments of the main character are:
$$\int_{SU_2}\chi^{2k}=C_k$$

(2) Our claim is that we can construct, by recurrence on $k\in\mathbb N$, a sequence $r_k$ of irreducible, self-adjoint and distinct representations of $SU_2$, satisfying:
$$r_0=1\quad,\quad
r_1=u\quad,\quad 
r_k+r_{k-2}=r_{k-1}\otimes r_1$$

Indeed, assume that $r_0,\ldots,r_{k-1}$ are constructed, and let us construct $r_k$. We have:
$$r_{k-1}+r_{k-3}=r_{k-2}\otimes r_1$$

Thus $r_{k-1}\subset r_{k-2}\otimes r_1$, and since $r_{k-2}$ is irreducible, by Frobenius we have:
$$r_{k-2}\subset r_{k-1}\otimes r_1$$

We conclude there exists a certain representation $r_k$ such that:
$$r_k+r_{k-2}=r_{k-1}\otimes r_1$$

(3) By recurrence, $r_k$ is self-adjoint. Now observe that according to our recurrence formula, we can split $u^{\otimes k}$ as a sum of the following type, with positive coefficients:  
$$u^{\otimes k}=c_kr_k+c_{k-2}r_{k-2}+\ldots$$

We conclude by Peter-Weyl that we have an inequality as follows, with equality precisely when $r_k$ is irreducible, and non-equivalent to the other summands $r_i$:
$$\sum_ic_i^2\leq\dim(End(u^{\otimes k}))$$

(4) But according to (1) the number on the right is $C_k$, and some straightforward combinatorics, based on the fusion rules, shows that the number on the left is $C_k$ as well. Thus we have equality in our estimate, so our representation $r_k$ is irreducible, and non-equivalent to $r_{k-2},r_{k-4},\ldots$ Moreover, this representation $r_k$ is not equivalent to $r_{k-1},r_{k-3},\ldots$ either, with this coming from $r_p\subset u^{\otimes p}$ for any $p$, and from:
$$\dim(Fix(u^{\otimes 2s+1}))=\int_{SU_2}\chi^{2s+1}=0$$

(5) Thus, we proved our claim. Now since each irreducible representation of $SU_2$ appears into some $u^{\otimes k}$, and we know how to decompose each $u^{\otimes k}$ into sums of representations $r_k$, these representations $r_k$ are all the irreducible representations of $SU_2$, and we are done with the main assertion. As for the dimension formula, this is clear.
\end{proof}

Regarding now $SO_3$, we have here a similar result, as follows:

\begin{theorem}
The irreducible representations of $SO_3$ are all self-adjoint, and can be labeled by positive integers, with their fusion rules being as follows,
$$r_k\otimes r_l=r_{|k-l|}+r_{|k-l|+1}+\ldots+r_{k+l}$$ 
also called Clebsch-Gordan rules. The corresponding dimensions are $\dim r_k=2k+1$.
\end{theorem}

\begin{proof}
As before with $SU_2$, there are many possible proofs here, which are all instructive. Here is our take on the subject, in the spirit of our proof for $SU_2$:

\medskip

(1) According to Theorem 4.30, the moments of the main character are:
$$\int_{SO_3}\chi^k=C_k$$

(2) Our claim now is that we can construct, by recurrence on $k\in\mathbb N$, a sequence $r_k$ of irreducible, self-adjoint and distinct representations of $SO_3$, satisfying:
$$r_0=1\quad,\quad
r_1=u-1\quad,\quad 
r_k+r_{k-1}+r_{k-2}=r_{k-1}\otimes r_1$$

Indeed, assume that $r_0,\ldots,r_{k-1}$ are constructed, and let us construct $r_k$. The Frobenius trick from the proof for $SU_2$ will no longer work, due to some technical reasons, so we have to invoke (1). To be more precise, by integrating characters we obtain:
$$r_{k-1},r_{k-2}\subset r_{k-1}\otimes r_1$$

We conclude that there exists a representation $r_k$ such that:
$$r_{k-1}\otimes r_1=r_k+r_{k-1}+r_{k-2}$$

(3) Once again by integrating characters, we conclude that $r_k$ is irreducible, and non-equivalent to $r_1,\ldots,r_{k-1}$, and this proves our claim. Also, since any irreducible representation of $SO_3$ must appear in some tensor power of $u$, and we can decompose each $u^{\otimes k}$ into sums of representations $r_p$, we conclude that these representations $r_p$ are all the irreducible representations of $SO_3$. Finally, the dimension formula is clear.
\end{proof}

As a conclusion to this, which is good to know, the advanced combinatorics of the Wigner and Marchenko-Pastur laws $\gamma_1$ and $\pi_1$ comes from the Clebsch-Gordan rules for $SU_2$ and $SO_3$. There are of course far more things that can be said, with for instance an alternative proof for this coming from a certain ``super-easiness'' property of $SU_2$ and $SO_3$. And, for more on this, you can have a look at the quantum algebra literature.

\section*{4e. Exercises}

Welcome to easiness, and as exercises about it, all easy I guess, we have:

\begin{exercise}
Learn about compact Lie groups, and their embeddings $G\subset U_N$.
\end{exercise}

\begin{exercise}
Learn more about the Haar measure, on various types of groups.
\end{exercise}

\begin{exercise}
Learn the various possible versions of Tannakian duality.
\end{exercise}

\begin{exercise}
Learn about the classification results for the easy groups.
\end{exercise}

\begin{exercise}
Clarify everything that we said, in relation with uniformity.
\end{exercise}

\begin{exercise}
Further study the shifted normal laws $g_t^t$ and $G_t^t$.
\end{exercise}

\begin{exercise}
Learn more about $SU_2$, as a continuation of the above.
\end{exercise}

\begin{exercise}
Learn more about $SO_3$, as a continuation of the above.
\end{exercise}

As bonus exercise, and no surprise here, read an easy quantum group book.

\part{Free probability}

\ \vskip50mm

\begin{center}
{\em Life is a mystery

Everyone must stand alone

I hear you call my name

And it feels like home}
\end{center}

\chapter{Operator algebras}

\section*{5a. Spectral theorems}

Welcome to free probability, and many things to be done. Our plan will be to axiomatize freeness in the noncommutative probability setting $(A,tr)$ developed in chapter 2, then work out an analogue of the logarithm of the Fourier transform, linearizing the free convolution, then talk about limiting theorems, applications to random matrices, easiness and free geometry, and many more. Excited about this I hope, and so am I.

\bigskip

Getting to work now.. but wait, cat is back here, meowing something, now that I started to talk again about free probability, let's see what the expert has to say:

\begin{cat}
There is a bug in your algebra formalism, apply the Gelfand theorem to $L^\infty(X)$, and you will understand.
\end{cat}

Damn. Technically, what cat says leads indeed to unspeakable horrors, I mean imagine a bit how the compactification given by $L^\infty(X)=C(\widehat{X})$ should look like, that is not mathematics, it's something else. In short, good point, we will have to postpone free probability, and first invent another class of operator algebras, ideally with the commutative ones being, by a rival to the Gelfand theorem, those of the form $A=L^\infty(X)$.

\bigskip

So, thanks cat, and getting now to work, since we are very interested in the commutativity of the algebras $A\subset B(H)$, our first goal will be that of having a spectral theorem for the normal operators, those for which $<T>$ is commutative. Let us start with:

\index{spectrum of operator}
\index{eigenvalues of operator}

\begin{definition}
The spectrum of an operator $T\in B(H)$ is the set
$$\sigma(T)=\left\{\lambda\in\mathbb C\Big|T-\lambda\not\in B(H)^{-1}\right\}$$
with $B(H)^{-1}\subset B(H)$ being the set of invertible operators.
\end{definition}

As a first observation, the spectrum contains the eigenvalue set, namely:
$$\varepsilon(T)=\left\{\lambda\in\mathbb C\Big|\ker(T-\lambda)\neq\{0\}\right\}$$

In finite dimensions $\varepsilon(T)\subset\sigma(T)$ is of course an equality. In infinite dimensions, however, this is wrong, the standard counterexample coming from the shift $S(e_i)=e_{i+1}$ on $H=l^2(\mathbb N)$, which is injective but not surjective, and so has $0\in\sigma(S)-\varepsilon(S)$.

\bigskip

The basic properties of the spectra can be summarized as follows:

\begin{theorem}
The following spectral theory results hold:
\begin{enumerate}
\item We have $\sigma(ST)\cup\{0\}=\sigma(TS)\cup\{0\}$.

\item In infinite dimensions $\sigma(ST)=\sigma(TS)$ can fail.

\item We have $\sigma(P(T))=P(\sigma(T))$, for any polynomial $P\in\mathbb C[X]$.

\item In fact, $\sigma(f(T))=f(\sigma(T))$ for any $f\in\mathbb C(X)$ having poles outside $\sigma(T)$.

\item The spectrum $\sigma(T)$ is compact, non-empty, and contained in $D_0(||T||)$.

\item For a unitary operator, $U^*=U^{-1}$, we have $\sigma(U)\subset\mathbb T$.

\item For a self-adjoint operator, $T=T^*$, we have $\sigma(T)\subset\mathbb R$.
\end{enumerate}
\end{theorem}

\begin{proof}
We already talked about this in chapter 2, in the $C^*$-algebra framework, but without much details, so here is the proof, with everything or almost included:

\medskip

(1) Assuming that $1-ST$ is invertible, with inverse $R$, we have $STR=RST=R-1$, and it follows that $1-TS$ is invertible too, with inverse $1+TRS$. Thus $\sigma(ST),\sigma(TS)$ agree on $1\in\mathbb C$, and by linearity, $\sigma(ST),\sigma(TS)$ agree on any point $\lambda\in\mathbb C^*$.

\medskip

(2) The standard counterexample here is provided by the shift $S(e_i)=e_{i+1}$ and its adjoint $S^*(e_i)=\delta_{i0}e_{i-1}$ on $H=l^2(\mathbb N)$. Indeed, $0\in\sigma(SS^*)-\sigma(S^*S)$.

\medskip

(3) If we pick $\lambda\in\mathbb C$, and decompose $P(X)-\lambda=c(X-r_1)\ldots(X-r_n)$, then we have the following equivalences, which are all clear, and which give the result:
\begin{eqnarray*}
\lambda\notin\sigma(P(T))
&\iff&P(T)-\lambda\in B(H)^{-1}\\
&\iff&c(T-r_1)\ldots(T-r_n)\in B(H)^{-1}\\
&\iff&T-r_1,\ldots,T-r_n\in B(H)^{-1}\\
&\iff&r_1,\ldots,r_n\notin\sigma(T)\\
&\iff&\lambda\notin P(\sigma(T))
\end{eqnarray*}

(4) Indeed, we can pick a scalar $\lambda\in\mathbb C$, write $f=P/Q$, and set $F=P-\lambda Q$. By using what we found in (3) for this polynomial $F\in\mathbb C[X]$ we obtain, as desired:
\begin{eqnarray*}
\lambda\in\sigma(f(T))
&\iff&F(T)\notin B(H)^{-1}\\
&\iff&0\in\sigma(F(T))\\
&\iff&0\in F(\sigma(T))\\
&\iff&\exists\mu\in\sigma(T),F(\mu)=0\\
&\iff&\lambda\in f(\sigma(T))
\end{eqnarray*}

(5) By using $(1-S)^{-1}=1+S+S^2+\ldots$ for $||S||<1$ we obtain that $\lambda-T$ is invertible for $||T||<|\lambda|$, and so $\sigma(T)\subset D_0(||T||)$. It is also clear that $\sigma(T)$ is closed, so what we have is a compact set. Finally, assuming $\sigma(T)=\emptyset$ the function $f(\lambda)=\varphi((T-\lambda)^{-1})$ is well-defined, for any $\varphi\in B(H)^*$, and by Liouville we get $f=0$, contradiction.

\medskip

(6) Given a unitary, $U^*=U^{-1}$, we have the following norm computation:
$$||U||
=\sqrt{||UU^*||}
=\sqrt{1}
=1$$

Now if we denote by $D$ the unit disk, we obtain from this, by using (5), that we have $\sigma(U)\subset D$. On the other hand, once again by using $U^*=U^{-1}$, we have as well:
$$||U^{-1}||
=||U^*||
=||U||
=1$$

Thus by (5) we have $\sigma(U^{-1})\subset D$, and then by (4) this gives $\sigma(U)\subset D^{-1}$. We conclude that we have $\sigma(U)\subset D\cap D^{-1}=\mathbb T$, as desired.

\medskip

(7) Consider the following rational function, depending on a parameter $r\in\mathbb R$:
$$f(z)=\frac{z+ir}{z-ir}$$

Then for $r>>0$ the operator $f(T)$ is well-defined, and if $T=T^*$, we have:
$$\left(\frac{T+ir}{T-ir}\right)^*
=\frac{T-ir}{T+ir}
=\left(\frac{T+ir}{T-ir}\right)^{-1}$$

Thus $f(T)$ is unitary, and by (6) we have $\sigma(T)\subset f^{-1}(\mathbb T)=\mathbb R$, as desired.
\end{proof}

Getting now to the spectral radius $\rho(T)=\sup_{\lambda\in\sigma(T)}|\lambda|$, we first have:

\index{spectral radius}

\begin{theorem}
The spectral radius of an operator $T\in B(H)$ is given by
$$\rho(T)=\lim_{n\to\infty}||T^n||^{1/n}$$
and in this formula, we can replace the limit by an inf.
\end{theorem}

\begin{proof}
This is something quite technical, the idea being as follows:

\medskip

(1) Our first claim is that the numbers $u_n=||T^n||^{1/n}$ are subject to:
$$(n+m)u_{n+m}\leq nu_n+mu_m$$

Indeed, we have the following estimate, using the Young inequality $ab\leq a^p/p+b^q/q$ from classical analysis, with exponents $p=(n+m)/n$ and $q=(n+m)/m$:
\begin{eqnarray*}
u_{n+m}
&=&||T^{n+m}||^{1/(n+m)}\\
&\leq&||T^n||^{1/(n+m)}||T^m||^{1/(n+m)}\\
&\leq&||T^n||^{1/n}\cdot\frac{n}{n+m}+||T^m||^{1/m}\cdot\frac{m}{n+m}\\
&=&\frac{nu_n+mu_m}{n+m}
\end{eqnarray*}

(2) Our second claim is that the last assertion of the theorem holds, namely:
$$\lim_{n\to\infty}||T^n||^{1/n}=\inf_n||T^n||^{1/n}$$

For this purpose, we can use the inequality found in (1). Indeed, fix $m\geq1$, let $n\geq1$, and write $n=lm+r$ with $0\leq r\leq m-1$. By using twice $u_{ab}\leq u_b$, we get:
$$u_n
\leq\frac{1}{n}( lmu_{lm}+ru_r)
\leq u_{m}+\frac{r}{n}\,u_1$$

It follows that we have $\lim\sup_nu_n\leq u_m$, which proves our claim.

\medskip

(3) Summarizing, we are left with proving the main formula, which is as follows, and with the remark that we already know that the sequence on the right converges:
$$\rho(T)=\lim_{n\to\infty}||T^n||^{1/n}$$

In one sense, we can use the polynomial calculus formula $\sigma(T^n)=\sigma(T)^n$, coming from Theorem 5.3. Indeed, this gives the following estimate, valid for any $n$, as desired:
$$\rho(T)
=\sup_{\lambda\in\sigma(T)}|\lambda|
=\sup_{\rho\in\sigma(T)^n}|\rho|^{1/n}
=\sup_{\rho\in\sigma(T^n)}|\rho|^{1/n}
=\rho(T^n)^{1/n}
\leq||T^n||^{1/n}$$

(4) For the reverse inequality, we fix a number $\rho>\rho(T)$, and we want to prove that we have $\rho\geq\lim_{n\to\infty}||T^n||^{1/n}$. By using the Cauchy formula, we have:
\begin{eqnarray*}
\frac{1}{2\pi i}\int_{|z|=\rho}\frac{z^n}{z-T}\,dz
&=&\frac{1}{2\pi i}\int_{|z|=\rho}\sum_{k=0}^\infty z^{n-k-1}T^k\,dz\\
&=&\sum_{k=0}^\infty\frac{1}{2\pi i}\left(\int_{|z|=\rho}z^{n-k-1}dz\right)T^k\\
&=&\sum_{k=0}^\infty\delta_{n,k+1}T^k\\
&=&T^{n-1}
\end{eqnarray*}

(5) By applying the norm we obtain from this formula the following estimate:
$$||T^{n-1}||
\leq\frac{1}{2\pi}\int_{|z|=\rho}\left|\left|\frac{z^n}{z-T}\right|\right|\,dz
\leq\rho^n\cdot\sup_{|z|=\rho}\left|\left|\frac{1}{z-T}\right|\right|$$

Since the sup does not depend on $n$, by taking $n$-th roots, we obtain in the limit:
$$\rho\geq\lim_{n\to\infty}||T^n||^{1/n}$$

(6) Now recall that $\rho$ was by definition an arbitrary number satisfying $\rho>\rho(T)$. Thus, we have obtained the following estimate, valid for any $T\in B(H)$:
$$\rho(T)\geq\lim_{n\to\infty}||T^n||^{1/n}$$

Thus, we are led to the conclusion in the statement.
\end{proof}

In the case of the normal elements, we have the following finer result:

\index{normal operator}

\begin{theorem}
The spectral radius of a normal element,
$$TT^*=T^*T$$
is equal to its norm.
\end{theorem}

\begin{proof}
We can proceed in two steps, as follows:

\medskip

\underline{Step 1}. In the case $T=T^*$ we have $||T^n||=||T||^n$ for any exponent of the form $n=2^k$, by using the formula $||TT^*||=||T||^2$, and by taking $n$-th roots we get:
$$\rho(T)\geq||T||$$

Thus, we are done with the self-adjoint case, with the result $\rho(T)=||T||$.

\medskip

\underline{Step 2}. In the general normal case $TT^*=T^*T$ we have $T^n(T^n)^*=(TT^*)^n$, and by using this, along with the result from Step 1, applied to $TT^*$, we obtain:
\begin{eqnarray*}
\rho(T)
&=&\lim_{n\to\infty}||T^n||^{1/n}\\
&=&\sqrt{\lim_{n\to\infty}||T^n(T^n)^*||^{1/n}}\\
&=&\sqrt{\lim_{n\to\infty}||(TT^*)^n||^{1/n}}\\
&=&\sqrt{\rho(TT^*)}\\
&=&\sqrt{||T||^2}\\
&=&||T||
\end{eqnarray*}

Thus, we are led to the conclusion in the statement.
\end{proof}

By using Theorem 5.5 we can say a number of non-trivial things about the normal operators, known as ``spectral theorem for normal operators''. We first have:

\begin{theorem}
Given $T\in B(H)$ normal, we have a morphism of algebras
$$\mathbb C[X]\to B(H)\quad,\quad 
P\to P(T)$$
having the properties $||P(T)||=||P_{|\sigma(T)}||$, and $\sigma(P(T))=P(\sigma(T))$.
\end{theorem}

\begin{proof}
This is an improvement of the polynomial calculus, in the normal case, with the extra assertion being the norm estimate. But the element $P(T)$ being normal, we can apply to it the spectral radius formula for normal elements, and we obtain:
$$||P(T)||
=\rho(P(T))
=\sup_{\lambda\in\sigma(P(T))}|\lambda|
=\sup_{\lambda\in P(\sigma(T))}|\lambda|
=||P_{|\sigma(T)}||$$

Thus, we are led to the conclusions in the statement.
\end{proof}

At a more advanced level now, we have the following result:

\begin{theorem}
Given $T\in B(H)$ normal, we have a morphism of algebras
$$C(\sigma(T))\to B(H)\quad,\quad 
f\to f(T)$$
which is isometric, $||f(T)||=||f||$, and has the property $\sigma(f(T))=f(\sigma(T))$.
\end{theorem}

\begin{proof}
The idea here is to ``complete'' the morphism in Theorem 5.6. Indeed, by Stone-Weierstrass, that morphism has a unique isometric extension, as follows:
$$C(\sigma(T))\to B(H)\quad,\quad  
f\to f(T)$$

It remains to prove $\sigma(f(T))=f(\sigma(T))$, and we can do this by double inclusion:

\medskip

``$\subset$'' Given a continuous function $f\in C(\sigma(T))$, we must prove that we have:
$$\lambda\notin f(\sigma(T))\implies\lambda\notin\sigma(f(T))$$

For this purpose, we can apply a suitable function to $T$, as follows:
$$\left(\frac{1}{f-\lambda}\right)T=\frac{1}{f(T)-\lambda}$$

In particular $f(T)-\lambda$ is invertible, so $\lambda\notin\sigma(f(T))$, as desired.

\medskip

``$\supset$'' Given a continuous function $f\in C(\sigma(T))$, we must prove that we have: 
$$\lambda\in f(\sigma(T))\implies\lambda\in\sigma(f(T))$$

But this is the same as proving that we have the following implication:
$$\mu\in\sigma(T)\implies f(\mu)\in\sigma(f(T))$$

For this purpose, we approximate our function by polynomials, $P_n\to f$, and we examine the following convergence, which follows from $P_n\to f$:
$$P_n(T)-P_n(\mu)\to f(T)-f(\mu)$$

We know from polynomial functional calculus that we have:
$$P_n(\mu)
\in P_n(\sigma(T))
=\sigma(P_n(T))$$

Thus, the operators $P_n(T)-P_n(\mu)$ are not invertible. On the other hand, we know that the set formed by the invertible operators is open, so its complement is closed. Thus the limit $f(T)-f(\mu)$ is not invertible either, and so $f(\mu)\in\sigma(f(T))$, as desired.
\end{proof}

Even more generally now, we have the following result:

\index{Borel functions}
\index{measurable calculus}

\begin{theorem}
Given $T\in B(H)$ normal, we have a morphism of algebras as follows, with $L^\infty$ standing for abstract measurable functions, or Borel functions,
$$L^\infty(\sigma(T))\to B(H)\quad,\quad 
f\to f(T)$$
which is isometric, $||f(T)||=||f||$, and has the property $\sigma(f(T))=f(\sigma(T))$.
\end{theorem}

\begin{proof}
As before, the idea will be that of completing what we have:

\medskip

(1) Given a vector $x\in H$, consider the following functional:
$$C(\sigma(T))\to\mathbb C\quad,\quad 
g\to<g(T)x,x>$$

By the Riesz theorem, this functional must be the integration with respect to a certain measure $\mu$ on the space $\sigma(T)$. Thus, we have a formula as follows:
$$<g(T)x,x>=\int_{\sigma(T)}g(z)d\mu(z)$$

Now given an arbitrary Borel function $f\in L^\infty(\sigma(T))$, as in the statement, we can define a number $<f(T)x,x>\in\mathbb C$, by using exactly the same formula, namely:
$$<f(T)x,x>=\int_{\sigma(T)}f(z)d\mu(z)$$

Thus, we have managed to define numbers $<f(T)x,x>\in\mathbb C$, for all vectors $x\in H$, and in addition we can recover these numbers as follows, with $g_n\in C(\sigma(T))$:
$$<f(T)x,x>=\lim_{g_n\to f}<g_n(T)x,x>$$ 

(2) In order to define now numbers $<f(T)x,y>\in\mathbb C$, for all vectors $x,y\in H$, we can use a polarization trick. Indeed, for any operator $S\in B(H)$ we have:
$$<S(x+y),x+y>=<Sx,x>+<Sy,y>+<Sx,y>+<Sy,x>$$

By replacing $y\to iy$, we have as well the following formula:
$$<S(x+iy),x+iy>=<Sx,x>+<Sy,y>-i<Sx,y>+i<Sy,x>$$

By multiplying this formula by $i$, and summing with the first one, we obtain:
\begin{eqnarray*}
<S(x+y),x+y>+i<S(x+iy),x+iy>
&=&(1+i)[<Sx,x>+<Sy,y>]\\
&+&2<Sx,y>
\end{eqnarray*}

(3) But with this, we can finish. Indeed, by combining (1,2), given a Borel function $f\in L^\infty(\sigma(T))$, we can define numbers $<f(T)x,y>\in\mathbb C$ for any $x,y\in H$, and we obtain in this way a certain operator $f(T)\in B(H)$, having all the desired properties.
\end{proof}

\index{spectral theorem}
\index{self-adjoint operator}

Good news, we can now diagonalize the normal operators. Let us start with:

\begin{theorem}
Any self-adjoint operator $T\in B(H)$ can be diagonalized,
$$T=U^*M_fU$$
with $U:H\to L^2(X)$ being a unitary operator from $H$ to a certain $L^2$ space associated to $T$, with $f:X\to\mathbb R$ being a certain function, once again associated to $T$, and with
$$M_f(g)=fg$$
being the usual multiplication operator by $f$, on the Hilbert space $L^2(X)$.
\end{theorem}

\begin{proof}
The construction of $U,f$ can be done in several steps, as follows:

\medskip

(1) We first prove the result in the special case where our operator $T$ has a cyclic vector $x\in H$, with this meaning that the following holds:
$$\overline{span\left(T^kx\Big|n\in\mathbb N\right)}=H$$

For this purpose, let us go back to the proof of Theorem 5.8. We will use the following formula from there, with $\mu$ being the measure on $X=\sigma(T)$ associated to $x$:
$$<g(T)x,x>=\int_{\sigma(T)}g(z)d\mu(z)$$

Our claim is that we can define a unitary $U:H\to L^2(X)$, first on the dense part spanned by the vectors $T^kx$, by the following formula, and then by continuity:
$$U[g(T)x]=g$$

Indeed, the following computation shows that $U$ is well-defined, and isometric:
$$||g(T)x||^2
=<|g|^2(T)x,x>
=\int_{\sigma(T)}|g(z)|^2d\mu(z)
=||g||_2^2$$

We can then extend this isometry $U$ by continuity into a unitary $U:H\to L^2(X)$, as claimed. Now observe that we have the following formula:
$$UTU^*g
=U[Tg(T)x]
=U[(zg)(T)x]
=zg$$

Thus our result is proved in the present case, with $U$ as above, and with $f(z)=z$.

\medskip

(2) We discuss now the general case. By using the Zorn lemma, $H$ has a decomposition as follows, with each $H_i$ being invariant under $T$, and admitting a cyclic vector $x_i$:
$$H=\bigoplus_iH_i$$

But with this, we can make a direct sum of the diagonalizations obtained in (1), for each of the restrictions $T_{|H_i}$, and we obtain the formula in the statement.
\end{proof}

We have the following technical generalization of the above result:

\index{family of operators}

\begin{theorem}
Any family of commuting self-adjoint operators $T_i\in B(H)$ can be jointly diagonalized,
$$T_i=U^*M_{f_i}U$$
with $U:H\to L^2(X)$ being a unitary operator from $H$ to a certain $L^2$ space associated to $\{T_i\}$, with $f_i:X\to\mathbb R$ being certain functions, once again associated to $T_i$, and with
$$M_{f_i}(g)=f_ig$$
being the usual multiplication operator by $f_i$, on the Hilbert space $L^2(X)$.
\end{theorem}

\begin{proof}
This is very similar to the proof of Theorem 5.9, say exercise for you.
\end{proof}

We can discuss now the case of the arbitrary normal operators, as follows:

\begin{theorem}
Any normal operator $T\in B(H)$ can be diagonalized,
$$T=U^*M_fU$$
with $U:H\to L^2(X)$ being a unitary operator from $H$ to a certain $L^2$ space associated to $T$, with $f:X\to\mathbb C$ being a certain function, once again associated to $T$, and with
$$M_f(g)=fg$$
being the usual multiplication operator by $f$, on the Hilbert space $L^2(X)$.
\end{theorem}

\begin{proof}
Consider the decomposition of $T$ into its real and imaginary parts:
$$T=\frac{T+T^*}{2}+i\cdot\frac{T-T^*}{2i}$$

Since our operator $T$ was assumed to be normal, these real and imaginary parts commute, so Theorem 5.10 applies to them, and gives the result.
\end{proof}

As our last result, regarding the families of commuting normal operators, we have:

\index{commuting normal operators}
\index{diagonalization}

\begin{theorem}
Any family of commuting normal operators $T_i\in B(H)$ can be jointly diagonalized,
$$T_i=U^*M_{f_i}U$$
with $U:H\to L^2(X)$ being a unitary operator from $H$ to a certain $L^2$ space associated to $\{T_i\}$, with $f_i:X\to\mathbb C$ being certain functions, once again associated to $T_i$, and with
$$M_{f_i}(g)=f_ig$$
being the usual multiplication operator by $f_i$, on the Hilbert space $L^2(X)$.
\end{theorem}

\begin{proof}
This is indeed straightforward, obtained by suitably modifying the above.
\end{proof}

And with this, end of our learning of spectral theorems. Good work that we did, and in what regards the applications of this technology, sky is the limit.

\section*{5b. Operator algebras}

In order to further build on the above spectral theorems, we need to know about the various topologies on $B(H)$, and fine-tune our operator theory knowledege, based on them. The result that we will need, which is something very basic, is as follows: 

\index{weak topology}
\index{weak operator topology}
\index{weakly closed}
\index{strongly closed}

\begin{proposition}
For a subalgebra $A\subset B(H)$, the following are equivalent:
\begin{enumerate}
\item $A$ is closed under the weak operator topology, making each of the linear maps $T\to<Tx,y>$ continuous.

\item $A$ is closed under the strong operator topology, making each of the linear maps $T\to Tx$ continuous.
\end{enumerate}
In the case where these conditions are satisfied, $A$ is closed under the norm topology.
\end{proposition}

\begin{proof}
It is clear that the norm topology is stronger than the strong operator topology, which is in turn stronger than the weak operator topology. Thus, we are left with proving that for any algebra $A\subset B(H)$, strongly closed implies weakly closed. Now for this purpose, consider the Hilbert space obtained by summing $n$ times $H$ with itself:
$$K=H\oplus\ldots\oplus H$$

The operators over $K$ can be regarded as being square matrices with entries in $B(H)$, and in particular, we have a representation $\pi:B(H)\to B(K)$, as follows:
$$\pi(T)=\begin{pmatrix}
T\\
&\ddots\\
&&T
\end{pmatrix}$$

Assume now that we are given an operator $T\in\bar{A}$, with the bar denoting the weak closure. We have then, by using the Hahn-Banach theorem, for any $x\in K$:
\begin{eqnarray*}
T\in\bar{A}
&\implies&\pi(T)\in\overline{\pi(A)}\\
&\implies&\pi(T)x\in\overline{\pi(A)x}\\
&\implies&\pi(T)x\in\overline{\pi(A)x}^{\,||.||}
\end{eqnarray*}

Now observe that the last formula tells us that for any $x=(x_1,\ldots,x_n)$, and any $\varepsilon>0$, we can find $S\in A$ such that the following holds, for any $i$:
$$||Sx_i-Tx_i||<\varepsilon$$

Thus $T$ belongs to the strong operator closure of $A$, as desired.
\end{proof}

Following von Neumann \cite{von}, we can now formulate a key definition, as follows:

\index{von Neumann algebra}

\begin{definition}
A von Neumann algebra is an operator $*$-algebra
$$A\subset B(H)$$
which is closed under the weak operator topology.
\end{definition}

As basic examples, we have the algebra $B(H)$ itself, then the singly generated algebras, $A=<T>$ with $T\in B(H)$, and then the multiply generated algebras, $A=<T_i>$ with $T_i\in B(H)$. There are of course many other examples, and more on this later.

\bigskip

But for the moment, let us keep things simple, and build directly on Definition 5.14, by using basic functional analysis methods. We will need the following key result:

\index{bicommutant theorem}
\index{amplification trick}

\begin{theorem}
For an operator algebra $A\subset B(H)$, we have
$$A''=\bar{A}$$
with $A''$ being the bicommutant inside $B(H)$, and $\bar{A}$ being the weak closure.
\end{theorem}

\begin{proof}
We can prove this by double inclusion, as follows:

\medskip

``$\supset$'' Since any operator commutes with the operators that it commutes with, we have a trivial inclusion $S\subset S''$, valid for any set $S\subset B(H)$. In particular, we have:
$$A\subset A''$$

Our claim now is that the algebra $A''$ is closed, with respect to the strong operator topology. Indeed, assuming that we have $T_i\to T$ in this topology, we have:
\begin{eqnarray*}
T_i\in A''
&\implies&ST_i=T_iS,\ \forall S\in A'\\
&\implies&ST=TS,\ \forall S\in A'\\
&\implies&T\in A
\end{eqnarray*}

Thus our claim is proved, and together with Proposition 5.13, which allows us to pass from the strong to the weak operator topology, this gives $\bar{A}\subset A''$, as desired.

\medskip

``$\subset$'' Here we must prove that we have the following implication, valid for any $T\in B(H)$, with the bar denoting as usual the weak operator closure:
$$T\in A''\implies T\in\bar{A}$$

For this purpose, we use the amplification trick from the proof of Proposition 5.13. With $\pi:B(H)\to B(K)$ being as there, the image of our algebra $A\subset B(H)$ is:
$$\pi(A)=\left\{\begin{pmatrix}
T\\
&\ddots\\
&&T
\end{pmatrix}\Big|T\in A\right\}$$

The commutant of this image is then given by the following formula:
$$\pi(A)'=\left\{\begin{pmatrix}
S_{11}&\ldots&S_{1n}\\
\vdots&&\vdots\\
S_{n1}&\ldots&S_{nn}
\end{pmatrix}\Big|S_{ij}\in A'\right\}$$

We conclude from this that, given an operator $T\in A''$ as above, we have:
$$\begin{pmatrix}
T\\
&\ddots\\
&&T
\end{pmatrix}\in\pi(A)''$$

In other words, the conclusion is that we have the following implication:
$$T\in A''\implies \pi(T)\in\pi(A)''$$

Now given a vector $x\in K$, consider the orthogonal projection $P\in B(K)$ on the norm closure of the vector space $\pi(A)x\subset K$. Since the subspace $\pi(A)x\subset K$ is invariant under the action of $\pi(A)$, so is its norm closure inside $K$, and we obtain from this:
$$P\in\pi(A)'$$

By combining this with what we found above, we conclude that we have:
$$T\in A''\implies \pi(T)P=P\pi(T)$$

Since this holds for any $x\in K$, we conclude that any operator $T\in A''$ belongs to the strong operator closure of $A$. By using now Proposition 5.13, which allows us to pass from the strong to the weak operator closure, we conclude that we have:
$$A''\subset\bar{A}$$

Thus, we have the desired reverse inclusion, and this finishes the proof.
\end{proof}

Now by getting back to the von Neumann algebras, from Definition 5.14, we have the following result, which is a reformulation of Theorem 5.15, by using this notion:

\begin{theorem}
For an operator algebra $A\subset B(H)$, the following are equivalent:
\begin{enumerate}
\item $A$ is weakly closed, so it is a von Neumann algebra.

\item $A$ equals its algebraic bicommutant $A''$, taken inside $B(H)$.
\end{enumerate}
\end{theorem}

\begin{proof}
This follows from the formula $A''=\bar{A}$ from Theorem 5.15, along with the trivial fact that the commutants are automatically weakly closed.
\end{proof}

The above statement, called bicommutant theorem, and due to von Neumann \cite{von}, is quite interesting, philosophically speaking. Among others, it shows that the von Neumann algebras are exactly the commutants of the self-adjoint sets of operators:

\index{commutant}

\begin{theorem}
Given a subset $S\subset B(H)$ which is closed under $*$, the commutant
$$A=S'$$
is a von Neumann algebra. Any von Neumann algebra appears in this way.
\end{theorem}

\begin{proof}
We have two assertions here, the idea being as follows:

\medskip

(1) Given $S\subset B(H)$ satisfying $S=S^*$, the commutant $A=S'$ satisfies $A=A^*$, and is also weakly closed. Thus, $A$ is a von Neumann algebra. Note that this follows as well from the following ``tricommutant formula'', which follows from Theorem 5.16:
$$S'''=S'$$

(2) Given a von Neumann algebra $A\subset B(H)$, we can take $S=A'$. Then $S$ is closed under the involution, and we have $S'=A$, as desired.
\end{proof}

Observe that Theorem 5.17 can be regarded as yet another alternative definition for the von Neumann algebras, and with this definition being probably the best one when talking about quantum mechanics, where the self-adjoint operators $T:H\to H$ can be though of as being ``observables'' of the system, and with the commutants $A=S'$ of the sets of such observables $S=\{T_i\}$ being the algebras $A\subset B(H)$ that we are interested in. And with all this actually needing some discussion about self-adjointness, and about  boundedness too, but let us not get into this here, and stay mathematical, as before.

\bigskip

As an interesting consequence of Theorem 5.17, let us record:

\index{center of algebra}

\begin{proposition}
Given a von Neumann algebra $A\subset B(H)$, its center
$$Z(A)=A\cap A'$$
regarded as an algebra $Z(A)\subset B(H)$, is a von Neumann algebra too.
\end{proposition}

\begin{proof}
We know from Theorem 5.17 that $A'\subset B(H)$ is a von Neumann algebra. Thus, the intersection $Z(A)=A\cap A'$ must be a von Neumann algebra too.
\end{proof}

Moving on, in finite dimensions things are very simple, as follows:

\index{matrix algebra}
\index{multimatrix algebra}
\index{sum of matrix algebras}
\index{finite dimensional algebra}

\begin{theorem}
The von Neumann algebras $A\subset M_N(\mathbb C)$ are the algebras
$$A=M_{n_1}(\mathbb C)\oplus\ldots\oplus M_{n_k}(\mathbb C)$$
depending on parameters $k\in\mathbb N$ and $n_1,\ldots,n_k\in\mathbb N$ satisfying
$$n_1+\ldots+n_k=N$$
embedded into $M_N(\mathbb C)$ via the obvious block embedding, twisted by a unitary $U\in U_N$.
\end{theorem}

\begin{proof}
This is something that we already know, from chapter 4 when doing Peter-Weyl, obtained by writing the center as $Z(A)=\mathbb C^k$, and analyzing the situation.
\end{proof}

As an interesting application now of Theorem 5.19, clarifying the relation with operator theory, in finite dimensions, we have the following result:

\begin{theorem}
Given an operator $T\in B(H)$ in finite dimensions, $H=\mathbb C^N$, the von Neumann algebra $A=<T>$ that it generates inside $B(H)=M_N(\mathbb C)$ is
$$A=M_{n_1}(\mathbb C)\oplus\ldots\oplus M_{n_k}(\mathbb C)$$
with the sizes of the blocks $n_1,\ldots,n_k\in\mathbb N$ coming from the spectral theory of the associated matrix $M\in M_N(\mathbb C)$. In the normal case $TT^*=T^*T$, this decomposition comes from
$$T=UDU^*$$
with $D\in M_N(\mathbb C)$ diagonal, and with $U\in U_N$ unitary.
\end{theorem}

\begin{proof}
This is something which is routine, by using the linear algebra and spectral theory developed in chapter 2, for the matrices $M\in M_N(\mathbb C)$. To be more precise:

\medskip

(1) The fact that $A=<T>$ decomposes into a direct sum of matrix algebras is something that we already know, coming from Theorem 5.19.

\medskip

(2) By using standard linear algebra, we can compute the block sizes $n_1,\ldots,n_k\in\mathbb N$, from the knowledge of the spectral theory of the associated matrix $M\in M_N(\mathbb C)$.

\medskip

(3) In the normal case, $TT^*=T^*T$, we can simply invoke the spectral theorem, and by suitably changing the basis, we are led to the conclusion in the statement.
\end{proof}

In infinite dimensions now, by using our spectral theorems, we first have:

\begin{theorem}
Given a bounded operator $T\in B(H)$ which is normal, $TT^*=T^*T$, the von Neumann algebra $A=<T>$ that it generates inside $B(H)$ is
$$<T>=L^\infty(\sigma(T))$$
with $\sigma(T)\subset\mathbb C$ being as usual its spectrum.
\end{theorem}

\begin{proof}
The measurable functional calculus theorem for the normal operators tells us that we have a weakly continuous morphism of $*$-algebras, as follows:
$$L^\infty(\sigma(T))\to B(H)\quad,\quad 
f\to f(T)$$

Moreover, we also know that this morphism is injective, and that its image is the weakly closed algebra $<T>$ generated by $T,T^*$. Thus, we obtain the result.
\end{proof}

More generally now, along the same lines, we have the following result:

\begin{theorem}
Given operators $T_i\in B(H)$ which are normal, and which commute, the von Neumann algebra $A=<T_i>$ that these operators generates inside $B(H)$ is
$$<T_i>=L^\infty(X)$$
with $X$ being a certain measured space, associated to the family $\{T_i\}$.
\end{theorem}

\begin{proof}
This is once again routine, by using the spectral theory for the families of commuting normal operators $T_i\in B(H)$, developed earlier in this chapter. 
\end{proof}

As a fundamental consequence now of the above results, we have:

\index{commutative algebra}

\begin{theorem}
The commutative von Neumann algebras are the algebras
$$A=L^\infty(X)$$
with $X$ being a measured space.
\end{theorem}

\begin{proof}
We have two assertions to be proved, the idea being as follows:

\medskip

(1) In one sense, we must prove that given a measured space $X$, we can realize the $A=L^\infty(X)$ as a von Neumann algebra, on a certain Hilbert space $H$. But this is something that we know since chapter 2, the representation being as follows: 
$$L^\infty(X)\subset B(L^2(X))\quad,\quad f\to(g\to fg)$$

(2) In the other sense, given a commutative von Neumann algebra $A\subset B(H)$, we must construct a certain measured space $X$, and an identification $A=L^\infty(X)$. But this follows from Theorem 5.22, because we can write our algebra as follows:
$$A=<T_i>$$

To be more precise, $A$ being commutative, any element $T\in A$ is normal, so we can pick a basis $\{T_i\}\subset A$, and then we have $A=<T_i>$ as above, with $T_i\in B(H)$ being commuting normal operators. Thus Theorem 5.22 applies, and gives the result.

\medskip

(3) Alternatively, and more explicitly, we can deduce this from Theorem 5.21, applied with $T=T^*$. Indeed, by using $T=Re(T)+iIm(T)$, we conclude that any von Neumann algebra $A\subset B(H)$ is generated by its self-adjoint elements $T\in A$. Moreover, by using measurable functional calculus, we conclude that $A$ is linearly generated by its projections. But then, assuming $A=\overline{span}\{p_i\}$, with $p_i$ being projections, we can set:
$$T=\sum_{i=0}^\infty\frac{p_i}{3^i}$$

Then $T=T^*$, and by functional calculus we have $p_0\in<T>$, then $p_1\in <T>$, and so on. Thus $A=<T>$, and $A=L^\infty(X)$ comes now via Theorem 5.21, as claimed.
\end{proof}

\section*{5c. Projections, order}

In order to further investigate the von Neumann algebras, the key idea will be that of looking at the projections. Let us start with some generalities:

\begin{definition}
Associated to any two projections $P,Q\in B(H)$ are:
\begin{enumerate}
\item The projection $P\wedge Q$, projecting on the common range.

\item The projection $P\vee Q$, projecting on the span of the ranges.
\end{enumerate}
\end{definition} 

Abstractly speaking, the above two operations can be thought of as being inf and sup type operations, and all the known algebraic formulae for inf and sup hold in this setting. Let us record here the following basic formula, which is something very useful:

\begin{proposition}
We have the following formula,
$$P+Q=P\wedge Q+P\vee Q$$
valid for any two projections $P,Q\in B(H)$.
\end{proposition}

\begin{proof}
This is clear from definitions, because when computing $P+Q$ we obtain the projection $P\vee Q$ on the span on the ranges, modulo the fact that the vectors in the common range are obtained twice, which amounts in saying that we must add $P\wedge Q$.
\end{proof}

Less trivial now is the following result, involving some functional analysis:

\begin{theorem}
Consider two projections $P,Q\in B(H)$.
\begin{enumerate}
\item In finite dimensions we have $(PQ)^n\to P\wedge Q$, in norm.

\item In general, we have $(PQ)^nx\to (P\wedge Q)x$, for any $x\in H$.

\item The operators $(PQ)^n$ do not necessarily converge in norm.
\end{enumerate}
\end{theorem}

\begin{proof}
We have several assertions here, the proof being as follows:

\medskip

(1) Assume that we are in the case $P,Q\in M_N(\mathbb C)$. By substracting $P\wedge Q$ from both $P,Q$, we can assume $P\wedge Q=0$, and we must prove that we have:
$$P\wedge Q=0\implies(PQ)^n\to 0$$

Our claim is that we have $||PQ||<1$. Indeed, we know that we have:
$$||PQ||\leq||P||\cdot||Q||=1$$

Assuming now by contradiction that we have $||PQ||=1$, since we are in finite dimensions, we must have, for a certain norm one vector, $||x||=1$:
$$||PQx||=1$$

Thus, we must have equalities in the following estimate:
$$||PQx||\leq||Qx||\leq||x||$$

But the second equality tells us that we must have $x\in Im(Q)$, and with this in hand, the first equality tells us that we must have $x\in Im(P)$. But this contradicts $P\wedge Q=0$, so we have proved our claim, and the convergence $(PQ)^n\to 0$ follows.

\medskip

(2) In infinite dimensions now, as before by substracting $P\wedge Q$ from both $P,Q$, we can assume $P\wedge Q=0$, and we must prove that we have, for any $x\in H$:
$$P\wedge Q=0\implies(PQ)^nx\to 0$$

For this purpose, we use a trick. Consider the following operator:
$$R=PQP$$

This operator is positive, because we have $R=(PQ)(PQ)^*$, and we have:
$$||R||\leq||P||\cdot||Q||\cdot||P||=1$$

Our claim, which will finish the proof, is that for any $x\in H$ we have:
$$R^nx\to0$$

In order to prove this claim, let us diagonalize $R$, by using the spectral theorem for the self-adjoint operators. If all the eigenvalues are $<1$ then we are done. If not, this means that we can find a nonzero vector $x\in H$ such that:
$$||Rx||=||x||$$

But this condition means that we must have equalities in the following estimate:
$$||PQPx||\leq||QPx||\leq||Px||\leq||x||$$

The point now is that this is impossible, due to our assumption $P\wedge Q=0$. Indeed, the last equality tells us that we must have $x\in Im(P)$, and with this in hand, the middle equality tells us that we must have $x\in Im(Q)$. But this contradicts $P\wedge Q=0$, so we have proved our claim, and the convergence $(PQ)^nx\to 0$ follows.

\medskip

(3) Finally, for a counterexample to $(PQ)^n\to0$, in infinite dimensions, we can take $H=l^2(\mathbb N)$, and then find projections $P,Q$ such that $(PQ)^ne_k\to 0$ for any $k$, but with the convergence arbitrarily slowing down with $k\to\infty$. Thus, $(PQ)^n\not\to0$.
\end{proof}

As a consequence, in connection with the von Neumann algebras, we have:

\begin{theorem}
Given two projections $P,Q\in B(H)$, the projections
$$P\wedge Q\quad,\quad P\vee Q$$
both belong to the von Neumann algebra generated by $P,Q$.
\end{theorem}

\begin{proof}
This comes from the above. Indeed, in what regards $P\wedge Q$, this is something that follows from Theorem 5.26. As for $P\vee Q$, here the result follows from the result for $P\wedge Q$, and from the formula $P+Q=P\wedge Q+P\vee Q$ from Proposition 5.25.
\end{proof}

The idea now will be that of studying the von Neumann algebras $A\subset B(H)$ by using their projections, $p\in A$. Let us start with the following result:

\begin{theorem}
Any von Neumann algebra is generated by its projections.
\end{theorem}

\begin{proof}
This is something from the proof of Theorem 5.23, coming from the measurable functional calculus, which can cut any normal operator into projections.
\end{proof}

As a key observation now, providing us with some geometric intuition, we have:

\index{equivalence of projections}

\begin{proposition}
Given a von Neumann algebra $A\subset B(H)$, the formula
$$p\simeq q\iff\exists u,\ \begin{cases}
uu^*=p\\
u^*u=q
\end{cases}$$
defines an equivalence relation for the projections $p\in A$.
\end{proposition}

\begin{proof}
This is something elementary, which follows from definitions, with the transitivity coming by composing the corresponding partial isometries.
\end{proof}

At a more advanced level, along the same lines, following von Neumann, we have:

\index{order on projections}

\begin{theorem}
Given a von Neumann algebra $A\subset B(H)$, we have a partial order on the projections $p\in A$, constructed as follows, with $u$ being a partial isometry,
$$p\preceq q\iff\exists u,\ \begin{cases}
uu^*=p\\
u^*u\leq q
\end{cases}$$
which is related to the equivalence relation $\simeq$ constructed above by:
$$p\simeq q\iff p\preceq q,\ q\preceq p$$
Thus, $\preceq$ is a partial order on the equivalence classes of projections $p\in A$.
\end{theorem}

\begin{proof}
We have several assertions here, the idea being as follows:

\medskip

(1) The fact that we have indeed a partial order is clear, with the transitivity coming, as before, by composing the corresponding partial isometries.

\medskip

(2) Regarding now the relation with $\simeq$, via the equivalence in the statement, the implication $\implies$ is clear. Thus, we are left with proving $\Longleftarrow$, which reads:
$$p\preceq q,\ q\preceq p\implies p\simeq q$$

Our assumption is that we have partial isometries $u,v$ such that:
$$uu^*=p\quad,\quad u^*u\leq q$$
$$v^*v\leq p\quad,\quad vv^*=q$$

We can construct then two sequences of decreasing projections, as follows:
$$p_0=p\quad,\quad p_{n+1}=v^*q_nv$$
$$q_0=q\quad,\quad q_{n+1}=u^*p_nu$$

Consider now the limits of these two sequences of projections, namely:
$$p_\infty=\bigwedge_ip_i\quad,\quad q_\infty=\bigwedge_iq_i$$

In terms of all these projections that we constructed, we have the following decomposition formulae for the original projections $p,q$:
$$p=(p-p_1)+(p_1-p_2)+\ldots+p_\infty$$
$$q=(q-q_1)+(q_1-q_2)+\ldots+q_\infty$$

Now observe that the summands are equivalent, with this being clear from the  definition of $p_n,q_n$ at the finite indices $n<\infty$, and with $p_\infty\simeq q_\infty$ coming from:
$$v^*q_\infty v=p_\infty\quad,\quad q_\infty vv^*q_\infty=q_\infty$$

Thus we obtain that we have $p\simeq q$, as desired, by summing.

\medskip

(3) Finally, the fact that the order relation $\preceq$ factorizes indeed to the equivalence classes under $\simeq$ follows from the equivalence established in (2).
\end{proof}

As a conclusion now, in view of Theorems 5.28 and 5.30, we can formulate:

\begin{fact}
We can think of any von Neumann algebra $A\subset B(H)$ as being a kind of object belonging to ``mathematical logic'', consisting of equivalence classes of projections $p\in A$, ordered via the relation $\preceq$, and producing $A$ itself via transport by partial isometries, and then linear combinations, and weak limits.
\end{fact}

Which is something quite remarkable, who on Earth could have guessed, when we were struggling with the von Neumann algebra basics, that we will end up with something that luminous. Well, that person on Earth who found this was von Neumann himself, back in the 1930s. And the above Fact 5.31, called ``von Neumann vision'' of the operator algebras, has been extremely useful ever since, and is still largely used nowadays.

\bigskip

From a probability perspective this is something very interesting too, reminding the hardline viewpoint on measure and integration, namely everything coming from the measurable sets, via step functions and limits. Powerful technology that we have here.

\section*{5d. Freeness, factors}

Among the von Neumann algebras $A\subset B(H)$, of particular interest are the ``free'' ones, having trivial center, $Z(A)=\mathbb C$. These algebras are called factors:

\index{factor}

\begin{definition}
A factor is a von Neumann algebra $A\subset B(H)$ whose center
$$Z(A)=A\cap A'$$
which is a commutative von Neumann algebra, reduces to the scalars, $Z(A)=\mathbb C$.
\end{definition}

This notion is something of key importance. The idea is that there are two main reasons for the study of factors. First, at the intuitive level, we have:

\begin{principle}[Freeness]
The following happen:
\begin{enumerate}
\item The condition $Z(A)=\mathbb C$ defining the factors is, obviously, opposite to the condition $Z(A)=A$ defining the commutative von Neumann algebras. 

\item Therefore, the factors are the von Neumann algebras which are ``free'', meaning as far as possible from the commutative ones.

\item Equivalently, with $A=L^\infty(X)$, the quantum spaces $X$ coming from factors are those which are ``free'', as far as possible from the classical spaces.
\end{enumerate}
\end{principle}

At a more advanced level, another motivation for the study of factors, which among others justifies the name ``factors'' for them, comes from the reduction theory of von Neumann, which is something non-trivial, that can be summarized as follows:

\index{reduction theory}
\index{integral of algebras}

\begin{principle}[Reduction theory]
Given a von Neumann algebra $A\subset B(H)$, if we write its center $Z(A)\subset A$, which is a commutative von Neumann algebra, as 
$$Z(A)=L^\infty(X)$$
with $X$ being a measured space, then the whole algebra decomposes as
$$A=\int_XA_x\,dx$$
with the fibers $A_x$ being factors, that is, satisfying $Z(A_x)=\mathbb C$.
\end{principle}

As a first comment, we have already seen an instance of such decomposition results before, when talking about finite dimensional algebras. Indeed, such algebras decompose, in agreement with the above, as direct sums of matrix algebras, as follows:
$$A=\bigoplus_xM_{n_x}(\mathbb C)$$

In general, however, things are more complicated than this, and technically speaking, and as opposed to Principle 5.33, which was more of a triviality, Principle 5.34 is a tough theorem, due to von Neumann. We will not get here further into this.

\bigskip

This was for the story, and let us close this philosophical discussion with:

\begin{conclusion}
Regardless of the approach and technical level, be that beginner or advanced, the von Neumann factors are the algebras that matter.
\end{conclusion}

Moving ahead now, remember that as probabilists, we are mainly interested in the case where our algebras have traces. And in relation with the factors, by leaving aside the rather trivial case of the matrix algebras $A=M_N(\mathbb C)$, we are led in this way to:

\begin{definition}
A ${\rm II}_1$ factor is a von Neumann algebra $A\subset B(H)$ which:
\begin{enumerate}
\item Is infinite dimensional, $\dim A=\infty$.

\item Has trivial center, $Z(A)=\mathbb C$.

\item Has a trace $tr:A\to\mathbb C$.
\end{enumerate}
\end{definition}

Here the order of the axioms is a bit random, with any of the possible $3!=6$ choices making sense, and corresponding to a slightly different vision on what the ${\rm II}_1$ factors truly are. With the above order, with (1) we are making it clear, right from the beginning, that we are not here for revolutionizing linear algebra. Then with (2) we adhere to Definition 5.32, and to what was said next about it, on freeness and reduction. And finally with (3) we adhere to our usual principle, that von Neumann algebras must have traces.

\bigskip

More technically now, the above definition is motivated by some heavy classification work of Murray-von Neumann and Connes \cite{co1}, \cite{mvo}, coming on top of the reduction theory of von Neumann, whose conclusion is more or less that everything in von Neumann algebras reduces, via certain product type procedures, to the study of ${\rm II}_1$ factors:

\begin{fact}
The ${\rm II}_1$ factors are the building blocks of the whole von Neumann algebra theory.
\end{fact}

Getting to work now, how to study the ${\rm II}_!$ factors? And here, leaving aside all sorts of advanced things that we said in the above, in practice, it is not even clear that such beasts exist. Fortunately, group theory comes to the rescue, as follows:

\index{group algebra}
\index{ICC}
\index{conjugacy classes}

\begin{theorem}
Given a discrete group $G$, consider its von Neumann algebra:
$$L(G)=\left\{\sum_g\lambda_gg\Big|\lambda_g\in\mathbb C\right\}''\subset B(l^2(G))\quad,\quad g(\delta_h)=\delta_{gh}$$
The center of this group von Neumann algebra is then given by
$$Z(L(G))=\left\{\sum_g\lambda_gg\Big|\lambda_{gh}=\lambda_{hg}\right\}''$$
and if $G\neq\{1\}$ has infinite conjugacy classes, in the sense that 
$$\Big|\{ghg^{-1}|g\in G\}\Big|=\infty\quad,\quad\forall h\neq1$$
with this being called ICC property, the algebra $L(G)$ is a ${\rm II}_1$ factor.
\end{theorem}

\begin{proof}
There are two assertions here, the idea being as follows:

\medskip

(1) The formula of the center comes indeed from the following computation:
\begin{eqnarray*}
\sum_g\lambda_gg\in Z(A)
&\iff&\sum_g\lambda_ggh=\sum_g\lambda_ghg\\
&\iff&\sum_k\lambda_{kh^{-1}}k=\sum_k\lambda_{h^{-1}k}k\\
&\iff&\lambda_{kh^{-1}}=\lambda_{h^{-1}k}
\end{eqnarray*}

(2) Regarding the last assertion, it is clear that $L(G)$ has a trace, given by $tr(g)=\delta_{g1}$. Concerning now the center, the condition $\lambda_{gh}=\lambda_{hg}$ that we found is equivalent to the fact that $g\to\lambda_g$ is constant on the conjugacy classes, and we obtain:
$$Z(L(G))=\mathbb C\iff G={\rm ICC}$$

Finally, assuming that this ICC condition is satisfied, with $G\neq\{1\}$, our group $G$ must be infinite, so the algebra $L(G)$ is infinite dimensional, as desired.
\end{proof}

In general now, following Murray and von Neumann \cite{mvo}, we have: 

\index{continuous dimension}
\index{unique trace}

\begin{theorem}
The trace of a ${\rm II}_1$ factor $tr:A\to\mathbb C$ is unique, and we have
$$\left\{tr(p)\Big|p\in A\ {\rm projection}\right\}=[0,1]$$
showing that, geometrically, we are in the context of ``continuous dimension''.
\end{theorem}

\begin{proof}
This is something quite technical, the idea being as follows:

\medskip

(1) Consider first an arbitrary factor $A$, and two nonzero projections $p,q\in A$. Since by factoriality we have $\vee_{u\in U_A}u^*pu=1$, we can find a unitary $u\in U_A$ such that $puq\neq0$. Now with $v\neq0$ being the polar part of $puq$, we have $vv^*\leq p$, $v^*v\leq q$. 

\medskip

(2) But with this, we can see that the partial order $\preceq$ constructed in Theorem 5.30 is a total order. Observe by the way that for the ${\rm II}_1$ factors the overall result is quite clear, because the key ingredient, $p\preceq q,\,q\preceq p\implies p\simeq q$, follows by using the trace.

\medskip

(3) Next, assume that a factor $A$ has a minimal projection $p\in A$. It is then routine to show, by using (1), that we have $A\simeq B(K)$, for some Hilbert space $K$. 

\medskip

(4) With this discussed, and getting now to what our statement says, let $A$ be a ${\rm II}_1$ factor, consider the set $S\subset[0,1]$ formed by the traces of the projections, and set:
$$c=\inf(S-\{0\})$$

(5) Our claim is that $c=0$. In order to prove this, assume by contradiction $c>0$, pick $\varepsilon>0$, and pick a projection $p\in A$ such that the following condition is satisfied:
$$tr(p)<c+\varepsilon$$

(6) Since we are in a ${\rm II}_1$ factor, this projection $p\in A$ cannot be minimal, so we can find a projection $q\in A$ satisfying $q<p$. But this leads to a contradiction, due to:
$$tr(p-q)
=tr(p)-tr(q)
\leq tr(p)-c
<\varepsilon$$

(7) Summarizing, the traces of projections can be arbitrarily small. But then $S=[0,1]$ follows from this, by using standard cut and paste arguments, based on Fact 5.31.

\medskip

(8) Finally, in what regards the uniqueness of the trace, this comes again from Fact 5.31, the idea being that two traces must agree on projections, and then in general.
\end{proof}

And with this, end of our study of the von Neumann algebras, good work that we did, we majestically answered the criticism from Cat 5.1. So, what do you think, cat?

\begin{cat}
There is a bug in your latest formalism too, because in quantum mechanics the operator algebra $A$ came first, and the Hilbert space $H$ came after.
\end{cat}

Good Lord. So, all our work in this chapter was in vain? Someone please axiomatize quantum mechanics, so that we mathematicians should know, what algebras to use.

\bigskip

In practice however, thinking some more, all this is not a disaster. We can use for our future work both $C^*$-algebras and von Neumann algebras, depending on needs.

\section*{5e. Exercises}

Welcome to von Neumann algebras, and as exercises on them, we have:

\begin{exercise}
Learn about the holomorphic functional calculus.
\end{exercise}

\begin{exercise}
Learn about the more specialized topologies available on $B(H)$.
\end{exercise}

\begin{exercise}
Learn about basic quantum mechanics, and unbounded operators.
\end{exercise}

\begin{exercise}
Meditate on the von Neumann vision of operator algebras.
\end{exercise}

\begin{exercise}
Have a look at the reduction theorem of von Neumann.
\end{exercise}

\begin{exercise}
Learn, with full details, the basic theory of the ${\rm II}_1$ factors.
\end{exercise}

\begin{exercise}
Read the classification work of Murray-von Neumann and Connes.
\end{exercise}

\begin{exercise}
Get to know about the hyperfinite ${\rm II}_1$ factor $R$.
\end{exercise}

As bonus exercise, if stuck with this, read the originals, \cite{co1}, \cite{mvo}, \cite{von}.

\chapter{Free probability}

\section*{6a. Quantum spaces}

Welcome to free probability, take two. Our plan will be to talk a bit about axiomatics, then talk about freeness, then work out an analogue of the logarithm of the Fourier transform, linearizing the free convolution, and then talk about limiting theorems, applications to random matrices, easiness and free geometry, and many more.

\bigskip

Let us start with something that we know well since chapter 2, namely:

\begin{definition}
A $C^*$-algebra is a complex normed algebra $A$, complete with respect to the norm, given with an involution $*$, such that $||aa^*||=||a||^2$, for any $a\in A$.
\end{definition}

We are already quite familiar with this notion, and our knowledge on the subject, coming from various results from chapters 2 and 5, can be summarized as follows:

\bigskip

(1) As basic examples, we have the matrix algebras $M_N(\mathbb C)$, and more generally, the algebras $B(H)$ of bounded operators on a Hilbert space $H$. Even more generally, any closed $*$-algebra $A\subset B(H)$ is a $C^*$-algebra. The GNS representation theorem, that we will prove in a moment, states that any $C^*$-algebra appears in this way, $A\subset B(H)$.

\bigskip

(2) The Gelfand theorem, based on spectral theory developed in the $C^*$-algebra framework, states that the commutative $C^*$-algebras are precisely those of the form $A=C(X)$, with $X$ being a compact space. And with this being something non-trivial, related in the case $A=<T>\subset B(H)$, with $T$ normal, to the spectral theorem for $T$. 

\bigskip

(3) The Gelfand theorem, while extremely useful, has however a number of limitations. Indeed, for the algebras of type $A=L^\infty(X)$, that we are particularly interested in, as probabilists, this theorem provides us with an isomorphism of type $A=C(\widehat{X})$, with $\widehat{X}$ being a certain complicated, and essentially useless, compactification of $X$.

\bigskip

(4) As a solution to this, we can talk instead about von Neumann algebras, which are the weakly closed $*$-algebras $A\subset B(H)$. These are of course $C^*$-algebras, of a special type, and do well the job with respect to probability requirements, the commutative ones being precisely those of the form $A=L^\infty(X)$, with $X$ being a measured space.

\bigskip

As a conclusion to this, if we want to do some serious mathematics, we need to master both theories. Indeed, looking now at the classical case, frankly, do you think that someone very specialized in measure theory, who has never heard about continuous functions, is a true mathematician? Or that someone very specialized in topology, who has never heard about measurable functions, is a true mathematician? Certainly not, and let's not do like them, we need to master both continuity and measurability.

\bigskip

Moving on, now that we better know our goals, we have the following question:

\begin{question}
Can we have the compact quantum measured spaces $X$ axiomatized, via certain $C^*$-algebras $A=C(X)$, and von Neumann algebras $A''=L^\infty(X)$?
\end{question}

And good question this is. In answer now, if we start with the von Neumann algebra viewpoint, we are a bit in the dark, because in the commutative case already, the isomorphism $A\simeq L^\infty(X)$, coming from the spectral theorem, is a bit hard to decipher, with several choices of $X$ being possible, some topological and compact, and some other not. So, in the hope that we agree on this, wrong way, we won't get far with this.

\bigskip

Summarizing, the question that makes sense is that of axiomatizing the compact quantum measured spaces $X$ via the corresponding $C^*$-algebras $A=C(X)$. And then talking, later, about the corresponding von Neumann algebras $A''=L^\infty(X)$ too.

\bigskip

Now in order to deal with this latter question, recall the Riesz theorem, telling us that the positive measures on the compact spaces $X$ are in correspondence with the associated integration functionals $\int:C(X)\to\mathbb C$. Thus, in order to answer Question 6.2, we need to talk about integration functionals $\int:A\to\mathbb C$, and we can do this as follows:

\index{positive element}
\index{positive linear form}
\index{faithful form}

\begin{definition}
Consider a continuous linear form on a $C^*$-algebra, $\varphi:A\to\mathbb C$.
\begin{enumerate}
\item $\varphi$ is called positive when $a\geq0\implies\varphi(a)\geq0$.

\item $\varphi$ is called faithful and positive when $a>0\implies\varphi(a)>0$.
\end{enumerate}
\end{definition}

In the commutative case, $A=C(X)$, the positive linear forms $\varphi$ appear as follows, with $\mu$ being positive, and strictly positive if we want $\varphi$ to be faithful and positive:
$$\varphi(f)=\int_Xf(x)d\mu(x)$$

In general, the positive linear forms can be thought of as being integration functionals with respect to some underlying ``positive measures''. We can use them as follows:

\index{GNS construction}

\begin{theorem}[GNS construction]
Let $\varphi:A\to\mathbb C$ be a positive linear form.
\begin{enumerate}
\item $<a,b>=\varphi(b^*a)$ defines a generalized scalar product on $A$.

\item Separation and completion gives a Hilbert space $H$.

\item $\pi(a):\widehat{b}\to\widehat{ab}$ defines a representation $\pi:A\to B(H)$.

\item If the linear form $\varphi$ is faithful, then $\pi$ is faithful.
\end{enumerate}
\end{theorem}

\begin{proof}
All this is very straightforward, the idea being as follows:

\medskip

(1) The pairing $<a,b>=\varphi(b^*a)$ is indeed linear at left, antilinear at right, satisfies $<a,b>=\overline{<b,a>}$, and is positive too, $<a,a>\geq0$. Thus, what we have a generalized scalar product, with the only missing axiom being $<a,a>=0\implies a=0$.

\medskip

(2) This is a standard procedure, which works for any generalized scalar product, the idea being that of dividing by the vectors satisfying $<x,x>=0$, then completing. Observe that we have a linear map $A\to H$, that we will denote $a\to\widehat{a}$.

\medskip

(3) There are several things to be proved here. To start with, given $a\in A$, for any $b\in A$ we have $b^*a^*ab\leq||a||^2b^*b$, and by applying $\varphi$ we have the following inequality:
$$||\widehat{ab}||\leq||a||\cdot||\widehat{b}||$$

Thus the linear map $\widehat{b}\to\widehat{ab}$, defined a priori only on the dense subspace $\widehat{A}\subset H$, is continuous on this subspace. By density we get a unique extension to $H$, that we may call $\pi(a)$. We prove now that $a\to\pi(a)$ is a representation. Associativity comes from:
$$\pi(ab)\widehat{c}=\widehat{abc}=\pi(a)\widehat{bc}=\pi(a)\pi(b)\widehat{c}$$

Also, the linearity and unitality are clear, and the involutivity comes from:
$$<\pi(a)\widehat{b},\widehat{c}>=<\widehat{ab},\widehat{c}>
=\varphi(c^*ab)=<\widehat{b},\widehat{a^*c}>
=<\widehat{b},\pi(a^*)\widehat{c}>$$

(4) In order to prove this, observe that with $\xi=\widehat{1}$ we have the following formula:
$$\varphi(a)
=<\widehat{a},\widehat{1}>
=<\pi(a)\widehat{1},\widehat{1}>
=<\pi(a)\xi,\xi>$$

Now assuming that $\varphi$ is faithful, let $b\in A$, $b\neq 0$. With $a=bb^*$ we have $\varphi(a)>0$, so the above formula gives $\pi(a)\neq 0$. Thus $\pi(b)\neq 0$, and so $\pi$ is faithful, as claimed.
\end{proof}

As a comment, observe the special role played by $\xi=\widehat{1}$, which is cyclic for $\pi$, in the sense that $\overline{\pi(A)\xi}=H$. Thus, the GNS construction produces from the pair $(A,\varphi)$ a certain triple $(H,\pi,\xi)$, with $H$ being a Hilbert space, $\pi:A\to B(H)$ being a representation, and $\xi\in H$ being a cyclic vector such that $\varphi(a)=<\pi(a)\xi,\xi>$. In fact, the GNS construction can be abstractly characterized in this way, say easy exercise for you.

\bigskip

Still following Gelfand, Naimark and Segal, we can now prove a key result:

\index{GNS theorem}
\index{Gelfand-Naimark-Segal}

\begin{theorem}[GNS]
Any $C^*$-algebra appears as closed $*$-algebra of operators
$$A\subset B(H)$$
over a certain Hilbert space $H$. When $A$ is separable, $H$ can be taken to be separable.
\end{theorem}

\begin{proof}
In view of Theorem 6.4, it is enough to prove that any $C^*$-algebra has a faithful positive linear form $\varphi:A\to\mathbb C$, and with this generalizing the fact that any compact space $X$ has a positive measure $\mu$. But this can be done as follows: 

\medskip

(1) Our first claim is that any positive linear form $\varphi:A\to\mathbb C$ is continuous. But this follows indeed from Theorem 6.4, via the following estimate:
$$|\varphi(a)|
\leq||\pi(a)||\varphi(1)
\leq||a||\varphi(1)$$

(2) Next, our claim is that a linear form $\varphi:A\to\mathbb C$ is positive precisely when there is a norm one positive element $h\in A_+$ such that $||\varphi||=\varphi(h)$. Indeed, in one sense we can take $h=1$, and we are done. Conversely, let $a\in A_+$, $||a||\leq1$. We have:
$$|\varphi(h)-\varphi(a)|
\leq||\varphi||\cdot||h-a||
\leq\varphi(h)$$

Thus we have $Re(\varphi(a))\geq0$, and with $a=1-h$ we obtain in this way:
$$Re(\varphi(1-h))\geq0$$

Thus $Re(\varphi(1))\geq||\varphi||$, and so $\varphi(1)=||\varphi||$, so we can assume $h=1$. Now observe that for any self-adjoint element $a$, and any $t\in\mathbb R$ we have, with $\varphi(a)=x+iy$:
\begin{eqnarray*}
\varphi(1)^2(1+t^2||a||^2)
&\geq&\varphi(1)^2||1+t^2a^2||\\
&=&||\varphi||^2\cdot||1+ita||^2\\
&\geq&|\varphi(1+ita)|^2\\
&=&|\varphi(1)-ty+itx|\\
&\geq&(\varphi(1)-ty)^2
\end{eqnarray*}

Thus we have $y=0$, and this finishes the proof of our remaining claim.

\medskip

(3) Next, for any $a\in A$ there exists a positive norm one form $\varphi:A\to\mathbb C$ such that $\varphi(aa^*)=||a||^2$. Indeed, we can set $\varphi(\lambda aa^*)=\lambda||a||^2$ on $span(aa^*)$, then extend this functional by Hahn-Banach, to the whole $A$, with the positivity coming from (2).

\medskip

(4) Our next claim is that if $A$ is separable, meaning that it has a dense sequence $\{a_n\}\subset A$, there is a faithful positive form $\varphi:A\to\mathbb C$. Indeed, for any $n$ we can construct as in (3) a positive form satisfying $\varphi(a_na_n^*)=||a_n||^2$, then define $\varphi$ as follows:
$$\varphi=\sum_{n=1}^\infty\frac{\varphi_n}{2^n}$$

In order to prove that $\varphi$ is faithful, let $a\in A$, $a\neq0$. If we pick $a_n$ close to $a$, and consider the GNS triple $(H,\pi,\xi)$ associated to the pair $(A,\varphi_n)$, we have:
\begin{eqnarray*}
\sqrt{\varphi_n(aa^*)}
&=&||\pi(a)\xi||\\
&\geq&||\pi(a_n)\xi||-||a-a_n||\\
&=&||a_n||-||a-a_n||\\
&>&0
\end{eqnarray*}

Thus $\varphi_n(aa^*)>0$, and so $\varphi(aa^*)>0$ and we have proved our claim.

\medskip

(5) Getting now to what our statement says, given an arbitrary $C^*$-algebra $A$, for any nonzero element $x\in A$ let $\varphi_x$ be a positive form such that $\varphi_x(xx^*)=||x||^2$. The GNS construction gives a certain triple $(H_x,\pi_x,\xi_x)$, and we have the following estimate:
$$||\pi_x(x)||^2
\geq||\pi_x(x)\xi_x||^2
=\varphi_x(xx^*)
=||x||^2$$

Thus we have $||\pi_x(x)||=||x||$, and if we take now $\pi$ to be the sum of all the representations $\pi_x$, we obtain in this way a certain embedding $\pi:A\subset H$, as desired.

\medskip

(6) Finally, in relation with separability, let $\varphi:A\to\mathbb C$ be a faithful positive form, as found in (4), and consider the GNS construction $(H,\pi ,\xi )$. We know that $\pi$ is faithful, and since $\pi (A)\xi$ is dense and $A$ is separable, $H$ follows to be separable, as desired.
\end{proof}

Good news, we can answer now Question 6.2, as follows:

\index{compact quantum space}
\index{compact quantum measured space}
\index{quantum measured space}

\begin{theorem}
We can talk about compact quantum measured spaces, as follows:
\begin{enumerate}
\item The category of compact quantum measured spaces $(X,\mu)$ is the category of $C^*$-algebras with faithful positive forms $(A,\varphi)$, with the arrows reversed. 

\item In the case where we have a non-faithful positive form $\varphi$, we can still talk about the corresponding space $(X,\mu)$, by performing the GNS construction.

\item In either case, by taking the weak closure in the GNS representation, we obtain a von Neumann algebra, that can be denoted $A''=L^\infty(X)$.
\end{enumerate}
\end{theorem}

\begin{proof}
All this follows from Theorems 6.4 and 6.5, and from other things that we already know, with the whole result itself being something rather philosophical.
\end{proof}

Moving on, what about imposing the extra axiom $\varphi(ab)=\varphi(ba)$? This certainly makes sense, for a myriad reasons, with for instance an advanced piece of motivation here coming from what we learned in chapter 5, namely that the ${\rm II}_1$ factors, which do have traces $tr:A\to\mathbb C$, are the building blocks of the whole von Neumann algebra theory.

\bigskip

So, let us complement Theorem 6.6 with the following addendum:

\begin{addendum}
For most common questions, we can restrict the attention to the tracial compact quantum measured spaces, those coming from a trace $tr:A\to\mathbb C$.
\end{addendum}

However, the story is not over with this, because with free probability ideas in mind, we can go nuclear, ignore everything classical, and focus on the ${\rm II}_1$ factor case:

\begin{addendum}
With free probability ideas in mind, we can restrict the attention to the compact quantum measured spaces coming from ${\rm II}_1$ factors.
\end{addendum}

And there is even more, along these lines, because among the ${\rm II}_1$ factors, just by looking at the examples of type $L(G)$, we can see that some of them, say with $G$ being free, are most likely more interesting for us than those with $G$ amenable.

\bigskip

In short, take it easy. What we have in Theorem 6.6, of perhaps in Addendum 6.7, is good work, but that is not the end of the story. And more on all this, later.

\section*{6b. Free variables}

With the above discussed, time now to get back to probability, and talk about freeness. As starting point, let us have our good old formalism from chapter 2, namely:

\index{random variable}
\index{moments}
\index{colored moments}
\index{law}
\index{distribution}

\begin{definition}
Let $A$ be a $C^*$-algebra, given with a positive trace $tr:A\to\mathbb C$.
\begin{enumerate}
\item The elements $a\in A$ are called random variables.

\item The moments of such a variable are the numbers $M_k(a)=tr(a^k)$.

\item The law of such a variable is the functional $\mu:P\to tr(P(a))$.
\end{enumerate}
\end{definition}

To be more precise, here $k=\circ\bullet\bullet\circ\ldots$ is a colored integer, and the corresponding powers $a^k$ are defined by the following formulae, and multiplicativity: 
$$a^\emptyset=1\quad,\quad
a^\circ=a\quad,\quad
a^\bullet=a^*$$

As for the polynomial $P$, this is a noncommuting $*$-polynomial in one variable:
$$P\in\mathbb C<X,X^*>$$

Observe that the law is uniquely determined by the moments, because we have:
$$P(X)=\sum_k\lambda_kX^k\implies\mu(P)=\sum_k\lambda_kM_k(a)$$

Generally speaking, the above definition is something quite abstract, but there is no other way of doing things, at least at this level of generality. However, in certain special cases, the formalism simplifies, and we recover more familiar objects, as follows:

\index{normal element}
\index{spectral measure}

\begin{theorem}
Assuming that $a\in A$ is normal, $aa^*=a^*a$, the following happen:
\begin{enumerate}
\item $\mu_a$ is a complex probability measure, satisfying $supp(\mu_a)\subset\sigma(a)$.

\item In the self-adjoint case, $a=a^*$, this measure $\mu_a$ is real.

\item Assuming that $tr$ is faithful, we have $supp(\mu_a)=\sigma(a)$.
\end{enumerate}
Also, when $a\in A$ is not normal, $aa^*\neq a^*a$, and $tr$ is faithful, all this fails.
\end{theorem}

\begin{proof}
This is something that we already know, from chapter 2. Indeed, according to the Gelfand theorem, we have an identification of algebras as follows: 
$$<a>=C(\sigma(a))$$

Thus the functional $f(a)\to tr(f(a))$ can be regarded as an integration functional on the algebra $C(\sigma(a))$, and by the Riesz theorem, this latter functional must come from a probability measure $\mu_a$ on the spectrum $\sigma(a)$, in the sense that we must have:
$$tr(f(a))=\int_{\sigma(a)}f(x)d\mu_a(x)$$

We are therefore led to the conclusions in the statement. Finally, the last assertion comes from $aa^*\neq a^*a\implies tr(aa^*aa^*)>tr(aaa^*a^*)$, as explained in chapter 2.
\end{proof}

Before going further, you might perhaps wonder where all our previous axiomatic discussion has gone, in relation with all this. In answer, the point comes from:

\begin{comment}
We reserve the right to slightly modify our formalism, upon needs:
\begin{enumerate}
\item If for certain questions we need our trace $tr:A\to\mathbb C$ to be faithful, we will assume so, by indicating that we are using a faithful trace.

\item If for certain questions we need to use positive forms $\varphi:A\to\mathbb C$ which are not tracial, we will do so, by indicating that we are exiting the tracial case.

\item If for certain questions we need to use von Neumann algebras $A$, or even more specialized objects, like ${\rm II}_1$ factors $A$, we will again do so.
\end{enumerate}
\end{comment}

To be more precise here, (1) is something that we already met, in Theorem 6.10, and we will keep doing that, namely assuming from time to time that $tr:A\to\mathbb C$ is faithful, as the usual probabilistic expectations $E:L^\infty(X)\to\mathbb C$ are, by mentioning this.

\bigskip

Regarding (2), this is something more subtle, and we will see a first instance of this very soon, when discussing models for freeness, with the simplest such models, leading to concrete results, being non-tracial. We will see illustrations of this later too.

\bigskip

Finally, regarding (3), the von Neumann algebra assumption is something quite natural, with the usual probabilistic expectations $E:L^\infty(X)\to\mathbb C$ being of this type. Observe that the random matrix expectations $E\circ tr:M_N(L^\infty(X))\to\mathbb C$ are of this type too.

\bigskip

Getting to work now, let us first discuss independence, in the framework of Definition 6.9. We have here the following notion, generalizing the one from chapter 1:

\index{independence}

\begin{definition}
Two subalgebras $A,B\subset C$ are called independent when the following condition is satisfied, for any $a\in A$ and $b\in B$: 
$$tr(ab)=tr(a)tr(b)$$
Equivalently, the following condition must be satisfied, for any $a\in A$ and $b\in B$: 
$$tr(a)=tr(b)=0\implies tr(ab)=0$$
Also, two variables $a,b\in C$ are called independent when the algebras that they generate, 
$$A=<a>\quad,\quad 
B=<b>$$
are independent inside $C$, in the above sense.
\end{definition}

Observe that the above two independence conditions are indeed equivalent, with this following from the following computation, with the convention $a'=a-tr(a)$:
\begin{eqnarray*}
tr(ab)
&=&tr[(a'+tr(a))(b'+tr(b))]\\
&=&tr(a'b')+t(a')tr(b)+tr(a)tr(b')+tr(a)tr(b)\\
&=&tr(a'b')+tr(a)tr(b)\\
&=&tr(a)tr(b)
\end{eqnarray*}

The other remark is that the above notion generalizes indeed the usual notion of independence, from the classical case, the precise result here being as follows:

\index{Fubini theorem}

\begin{theorem}
Given two compact measured spaces $X,Y$, the algebras
$$C(X)\subset C(X\times Y)\quad,\quad 
C(Y)\subset C(X\times Y)$$
are independent in the above sense, and a converse of this fact holds too.
\end{theorem}

\begin{proof}
We have two assertions here, the idea being as follows:

\medskip

(1) First of all, given two abstract compact spaces $X,Y$, we have embeddings of algebras as in the statement, defined by the following formulae:
$$f\to[(x,y)\to f(x)]\quad,\quad 
g\to[(x,y)\to g(y)]$$

In the measured space case now, the Fubini theorem tells us that we have:
$$\int_{X\times Y}f(x)g(y)=\int_Xf(x)\int_Yg(y)$$

Thus, the algebras $C(X),C(Y)$ are independent in the sense of Definition 6.12.

\medskip

(2) Conversely, assume that $A,B\subset C$ are independent, with $C$ being commutative. Let us write our algebras as follows, with $X,Y,Z$ being certain compact spaces:
$$A=C(X)\quad,\quad 
B=C(Y)\quad,\quad
C=C(Z)$$ 

In this picture, the inclusions $A,B\subset C$ must come from quotient maps, as follows:
$$p:Z\to X\quad,\quad 
q:Z\to Y$$

Regarding now the independence condition from Definition 6.12, in the above picture, this tells us that the following equality must happen:
$$\int_Zf(p(z))g(q(z))=\int_Zf(p(z))\int_Xg(q(z))$$

Thus we are in a Fubini type situation, and we obtain from this:
$$X\times Y\subset Z$$

Thus, the independence of the algebras $A,B\subset C$ appears as in (1) above.
\end{proof}

As an alternative notion to independence, we have Voiculescu's freeness \cite{vo1}:

\index{freeness}
\index{free variables}

\begin{definition}
Two subalgebras $A,B\subset C$ are called free when the following condition is satisfied, for any $a_i\in A$ and $b_i\in B$:
$$tr(a_i)=tr(b_i)=0\implies tr(a_1b_1a_2b_2\ldots)=0$$
Also, two variables $a,b\in C$ are called free when the algebras that they generate,
$$A=<a>\quad,\quad 
B=<b>$$
are free inside $C$, in the above sense.
\end{definition}

In short, freeness appears by definition as a kind of ``free analogue'' of usual independence, taking into account the fact that the variables do not necessarily commute. 

\bigskip

As a first observation, of theoretical nature, there is actually a certain lack of symmetry between Definition 6.12 and Definition 6.14, because in contrast to the former, the latter does not include an explicit formula for the quantities of the following type:
$$tr(a_1b_1a_2b_2\ldots)$$

However, this is not an issue, and is simply due to the fact that the formula in the free case is something more complicated, the precise result being as follows:

\begin{proposition}
Assuming that $A,B\subset C$ are free, the restriction of $tr$ to $<A,B>$ can be computed in terms of the restrictions of $tr$ to $A,B$. To be more precise,
$$tr(a_1b_1a_2b_2\ldots)=P\Big(\{tr(a_{i_1}a_{i_2}\ldots)\}_i,\{tr(b_{j_1}b_{j_2}\ldots)\}_j\Big)$$
where $P$ is certain polynomial in several variables, depending on the length of the word $a_1b_1a_2b_2\ldots$, and having as variables the traces of products of type
$$a_{i_1}a_{i_2}\ldots\quad,\quad 
b_{j_1}b_{j_2}\ldots$$
with the indices being chosen increasing, $i_1<i_2<\ldots$ and $j_1<j_2<\ldots$
\end{proposition}

\begin{proof}
This is something a bit theoretical, so let us begin with an example. Our claim is that if $a,b$ are free then, exactly as in the case where we have independence:
$$tr(ab)=tr(a)tr(b)$$

Indeed, let us go back to the independence computation performed after Definition 6.12, which was as follows, with the convention $a'=a-tr(a)$:
\begin{eqnarray*}
tr(ab)
&=&tr[(a'+tr(a))(b'+tr(b))]\\
&=&tr(a'b')+t(a')tr(b)+tr(a)tr(b')+tr(a)tr(b)\\
&=&tr(a'b')+tr(a)tr(b)\\
&=&tr(a)tr(b)
\end{eqnarray*}

But this computation perfectly works under the sole freeness assumption. Indeed, the only non-trivial equality is the last one, and this follows from:
$$tr(a')=tr(b')=0\implies tr(a'b')=0$$

In general now, the situation is of course more complicated than this, but the same trick applies. To be more precise, we can start our computation as follows:
\begin{eqnarray*}
tr(a_1b_1a_2b_2\ldots)
&=&tr\big[(a_1'+tr(a_1))(b_1'+tr(b_1))(a_2'+tr(a_2))(b_2'+tr(b_2))\ldots\ldots\big]\\
&=&tr(a_1'b_1'a_2'b_2'\ldots)+{\rm other\ terms}\\
&=&{\rm other\ terms}
\end{eqnarray*}

Observe that we have used here the freeness condition, in the following form:
$$tr(a_i')=tr(b_i')=0\implies tr(a_1'b_1'a_2'b_2'\ldots)=0$$

Now regarding the ``other terms'', those which are left, each of them will consist of a product of traces of type $tr(a_i)$ and $tr(b_i)$, and then a trace of a product still remaining to be computed, which is of the following form, for some elements $\alpha_i\in A$ and $\beta_i\in B$:
$$tr(\alpha_1\beta_1\alpha_2\beta_2\ldots)$$

To be more precise, the variables $\alpha_i\in A$ appear as ordered products of those $a_i\in A$ not getting into individual traces $tr(a_i)$, and the variables $\beta_i\in B$ appear as ordered products of those $b_i\in B$ not getting into individual traces $tr(b_i)$. Now since the length of each such alternating product $\alpha_1\beta_1\alpha_2\beta_2\ldots$ is smaller than the length of the original product $a_1b_1a_2b_2\ldots$, we are led into of recurrence, and this gives the result.
\end{proof}

Let us discuss now some models for independence and freeness. We have the following result, from \cite{vo1}, which clarifies the analogy between independence and freeness:

\index{tensor product}
\index{free product}
\index{free product trace}

\begin{theorem}
Given two algebras $(A,tr)$ and $(B,tr)$, the following hold:
\begin{enumerate}
\item $A,B$ are independent inside their tensor product $A\otimes B$, endowed with its canonical tensor product trace, given by $tr(a\otimes b)=tr(a)tr(b)$.

\item $A,B$ are free inside their free product $A*B$, endowed with its canonical free product trace, given by the formulae in Proposition 6.15.
\end{enumerate}
\end{theorem}

\begin{proof}
Both the above assertions are clear from definitions, as follows:

\medskip

(1) This is clear with either of the definitions of the independence, from Definition 6.12, because we have by construction of the product trace:
$$tr(ab)
=tr[(a\otimes1)(1\otimes b)]
=tr(a\otimes b)
=tr(a)tr(b)$$

Observe that there is a relation here with Theorem 6.13 as well, due to the following formula for compact spaces, with $\otimes$ being a topological tensor product:
$$C(X\times Y)=C(X)\otimes C(Y)$$

To be more precise, the present statement generalizes the first assertion in Theorem 6.13, and the second assertion tells us that this generalization is more or less the same thing as the original statement. All this comes of course from basic measure theory.

\medskip

(2) This is clear too from definitions, the only point being that of showing that the notion of freeness, or the recurrence formulae in Proposition 6.15, can be used in order to construct a canonical free product trace, on the free product of the algebras involved:
$$tr:A*B\to\mathbb C$$

But this can be checked for instance by using a GNS construction. Indeed, consider the GNS constructions for the algebras $(A,tr)$ and $(B,tr)$:
$$A\to B(L^2(A))\quad,\quad
B\to B(L^2(B))$$

By taking the free product of these representations, we obtain a representation as follows, with the $*$ on the right being a free product of pointed Hilbert spaces:
$$A*B\to B(L^2(A)*L^2(B))$$

Now by composing with the linear form $T\to<T\xi,\xi>$, where $\xi=1_A=1_B$ is the common distinguished vector of $L^2(A)$, $L^2(B)$, we obtain a linear form, as follows:
$$tr:A*B\to\mathbb C$$

It is routine then to check that $tr$ is indeed a trace, and this is the ``canonical free product trace'' from the statement. Then, an elementary computation shows that $A,B$ are free inside $A*B$, with respect to this trace, and this finishes the proof. See \cite{vo1}.
\end{proof}

More concretely now, we have the following result, also from \cite{vo1}:

\index{group algebras}

\begin{theorem}
We have the following results, valid for group algebras:
\begin{enumerate}
\item $C^*(G),C^*(H)$ are independent inside $C^*(G\times H)$.

\item $C^*(G),C^*(H)$ are free inside $C^*(G*H)$.
\end{enumerate}
\end{theorem}

\begin{proof}
This something very standard, the idea being as follows:

\medskip

(1) To start with, recall that associated to a discrete group $G$ is its group algebra $\mathbb C[G]$, formal linear span of the group elements $g\in G$, with involution $g^*=g^{-1}$. We can then talk about $C^*$-norms on this algebra, subject to $||aa^*||=||a||^2$, and define $C^*(G)$ as being the completion of  $\mathbb C[G]$ with respect to the biggest such norm.

\medskip

(2) And this, coming with the remark that such $C^*$-norms exist indeed, a basic example being the one coming from the embedding $\mathbb C[G]\subset B(l^2(G))$ given by $g(\delta_h)=\delta_{gh}$. And with the extra remark that the biggest $C^*$-norm is bounded, due to the fact that the standard generators $g\in G$ being unitaries, $g^*=g^{-1}$, they must satisfy $||g||=1$.

\medskip

(3) Getting now to what the statement says, we can use the general results in Theorem 6.16, along with the following two isomorphisms, which are both standard:
$$C^*(G\times H)=C^*(G)\otimes C^*(H)\quad,\quad
C^*(G*H)=C^*(G)*C^*(H)$$

Alternatively, we can check the independence and freeness formulae on group elements, which is something trivial, and then conclude by linearity.
\end{proof}

All the above was quite theoretical, and as a concrete application of the above results, bringing us into probability, we have the following result, from \cite{vo1}:

\index{free convolution}

\begin{theorem}
We have a free convolution operation $\boxplus$ for the distributions
$$\mu:\mathbb C<X,X^*>\to\mathbb C$$
which is well-defined by the following formula, with $a,b$ taken to be free:
$$\mu_a\boxplus\mu_b=\mu_{a+b}$$
This restricts to an operation, still denoted $\boxplus$, on the real probability measures.
\end{theorem}

\begin{proof}
We have several verifications to be performed here, as follows:

\medskip

(1) We first have to check that given two variables $a,b$ which live respectively in certain $C^*$-algebras $A,B$, we can recover inside some $C^*$-algebra $C$, with exactly the same distributions $\mu_a,\mu_b$, as to be able to sum them and talk about $\mu_{a+b}$. But this comes from Theorem 6.16, because we can set $C=A*B$, as explained there.

\medskip

(2) The other verification which is needed is that of the fact that if two variables $a,b$ are free, then the distribution $\mu_{a+b}$ depends only on the distributions $\mu_a,\mu_b$. But for this purpose, we can use the general formula from Proposition 6.15, namely:
$$tr(a_1b_1a_2b_2\ldots)=P\Big(\{tr(a_{i_1}a_{i_2}\ldots)\}_i,\{tr(b_{j_1}b_{j_2}\ldots)\}_j\Big)$$

Now by plugging in arbitrary powers of $a,b$ as variables $a_i,b_j$, we obtain a family of formulae of the following type, with $Q$ being certain polyomials:
$$tr(a^{k_1}b^{l_1}a^{k_2}b^{l_2}\ldots)=Q\Big(\{tr(a^k)\}_k,\{tr(b^l)\}_l\Big)$$

Thus the moments of $a+b$ depend only on the moments of $a,b$, with of course colored exponents in all this, according to our moment conventions, and this gives the result.

\medskip

(3) Finally, in what regards the last assertion, regarding the real measures, this is clear from the fact that if the variables $a,b$ are self-adjoint, then so is their sum $a+b$.
\end{proof}

Along the same lines, but with some technical subtleties this time, we can talk as well about multiplicative free convolution, following \cite{vo2}, as follows:

\index{multiplicative free convolution}

\begin{theorem}
We have a free convolution operation $\boxtimes$ for the distributions
$$\mu:\mathbb C<X,X^*>\to\mathbb C$$
which is well-defined by the following formula, with $a,b$ taken to be free:
$$\mu_a\boxtimes\mu_b=\mu_{ab}$$
In the case of positive variables, we can equally set, with the same outcome,
$$\mu_a\boxtimes\mu_b=\mu_{\sqrt{a}\,b\,\sqrt{a}}$$
so we have an operation, still denoted $\boxtimes$, on the real positive probability measures.
\end{theorem}

\begin{proof}
We have two statements here, the idea being as follows:

\medskip

(1) The verifications for the fact that $\boxtimes$ as above is indeed well-defined at the general distribution level are identical to those done before for $\boxplus$, with the result basically coming from the formula in Proposition 6.15, and with Theorem 6.16 invoked as well, in order to say that we have a model, and so we can indeed use this formula.

\medskip

(2) Regarding the last assertion, real measures, this was something trivial for $\boxplus$, but is something trickier for $\boxtimes$, because if we take $a,b$ to be self-adjoint, their product $ab$ will be not self-adjoint, in general. However, when $a>0$ we can trick, by setting:
$$c=\sqrt{a}\,b\,\sqrt{a}$$

Indeed, this new variable is then self-adjoint, and its moments are given by:
\begin{eqnarray*}
tr(c^k)
&=&tr[(\sqrt{a}\,b\,\sqrt{a})^k]\\
&=&tr[\sqrt{a}\,baba\ldots bab\,\sqrt{a}]\\
&=&tr[\sqrt{a}\cdot \sqrt{a}\,baba\ldots bab]\\
&=&tr[(ab)^k]
\end{eqnarray*}

We conclude that the $\boxtimes$ operation restricts into an operation as follows:
$$\boxtimes:\mathcal P(\mathbb R_+)\times\mathcal P(\mathbb R)\to\mathcal P(\mathbb R)$$

Moreover, since $b>0\implies c>0$, we can make things symmetric, as stated.
\end{proof}

\section*{6c. Linearization}

We have seen so far the foundations of free probability, in analogy with those of classical probability, taken with a functional analysis touch. The idea now is that with a bit of luck, the basic theory from the classical case, namely the Fourier transform, and then the CLT, should have free extensions. Let us being our discussion with the following definition, from \cite{vo1}, coming from the theory developed in the above:

\index{free convolution}

\begin{definition}
The real probability measures are subject to operations $*$ and $\boxplus$, called classical and free convolution, given by the formulae
$$\mu_a*\mu_b=\mu_{a+b}\quad,\quad 
\mu_\alpha\boxplus\mu_\beta=\mu_{\alpha+\beta}$$
with $a,b$ being independent, and $\alpha,\beta$ being free, and all variables being self-adjoint.
\end{definition}

The problem now is that of linearizing these operations $*$ and $\boxplus$. In what regards $*$, we know from chapter 1 that this operation is linearized by the logarithm $\log F$ of the Fourier transform, which in the present setting, where $E=tr$, is given by:
$$F_a(x)=tr(e^{ixa})$$

In order to find a similar result for $\boxplus$, we need some efficient models for the pairs of free random variables $(a,b)$. This is a priori not a problem, because once we have $a\in A$ and $b\in B$, we can form the free product $A*B$, which contains $a,b$ as free variables. 

\bigskip

However, the initial choice, that of the variables $a\in A$, $b\in B$ modeling some given laws $\mu,\nu\in\mathcal P(\mathbb R)$, matters a lot. Indeed, any kind of abstract choice here would lead us into an abstract algebra $A*B$, and so into the abstract combinatorics of the free convolution, that cannot be solved with bare hands, and that we want to avoid.

\bigskip

In short, we must be tricky, at least in what concerns the beginning of our computation. Following \cite{vo1}, the idea will be that of temporarily lifting the self-adjointness assumption on our variables $a,b$, and looking instead for random variables $\alpha,\beta$, not necessarily self-adjoint, modeling in integer moments our given laws $\mu,\nu\in\mathcal P(\mathbb R)$, as follows:
$$tr(\alpha^k)=M_k(\mu)\quad,\quad tr(\beta^k)=M_k(\nu)$$

To be more precise, assuming that $\alpha,\beta$ are indeed not self-adjoint, the above formulae are not the general formulae for $\alpha,\beta$, simply because these latter formulae involve colored integers $k=\circ\bullet\bullet\circ\ldots$ as exponents. Thus, in the context of the above formulae, $\mu,\nu$ are not the distributions of $\alpha,\beta$, but just some ``parts'' of these distributions.

\bigskip

Now with this idea in mind, due to Voiculescu and quite tricky, the solution to the law modeling problem comes in a quite straightforward way, involving the good old Hilbert space $H=l^2(\mathbb N)$ and the good old shift operator $S\in B(H)$, as follows:

\index{shift}

\begin{theorem}
Consider the shift operator on $H=l^2(\mathbb N)$, given by $S(e_i)=e_{i+1}$. The operators of the following type, with $f\in\mathbb C[X]$ being a polynomial, 
$$T=S^*+f(S)$$
model then in uncolored moments, up to finite order, all the distributions $\mu:\mathbb C[X]\to\mathbb C$, when regarded as random variables with respect to the state $\varphi(T)=<Te_0,e_0>$.
\end{theorem}

\begin{proof}
This is something quite tricky, the idea being as follows:

\medskip

(1) Regarding the statement, as just mentioned, this is something subtle, exiting the framework of Definition 6.9. That definition was something clean, with the expectation being a trace, $tr(ab)=tr(ba)$, and with us being interested in the colored moments $M_k(a)$. What we are doing here is to use the state $\varphi(T)=<Te_0,e_0>$, which is not a trace, and get interested in uncolored moments only, that is, $M_k(a)=\varphi(a^k)$ with $k\in\mathbb N$.

\medskip

(2) Good, so with the statement agreed upon, let us turn now to the proof. We know that the shift and its adjoint are given by the following formulae:
$$S(e_i)=e_{i+1}\quad,\quad S^*(e_i)=
\begin{cases}
e_{i-1}&(i>0)\\
0&(i=0)
\end{cases}$$

Consider now a variable as in the statement, written as follows:
$$T=S^*+a_01+a_1S+a_2S^2+\ldots+a_nS^n$$

The computation of the moments of $T$ amounts then in expanding $T^k$, and looking for terms which simplify up to 1, according to the rule $S^*S=1$. In practice, we get:
$$\varphi(T)=\varphi(a_01)=a_0$$
$$\varphi(T^2)=\varphi(a_0^211+a_1S^*S)=a_0^2+a_1$$
$$\varphi(T^3)=\varphi(a_0^3111+a_0a_11S^*S+a_0a_1S^*1S+a_0a_1S^*S1)=a_0^2+3a_0a_1$$

At order 4 the computation is similar, but a bit more complex, as follows:
\begin{eqnarray*}
\varphi(T^4)
&=&\varphi\big[a_0^41111+a_0^2a_1(11S^*S+1S^*1S+1S^*S1+S^*11S+S^*1S1+S^*S11)\\
&&+a_1^2(S^*S^*SS+S^*SS^*S)+a_2S^*S^*S^2\big]\\
&=&a_0^4+6a_0^2a_1+2a_1^2+a_2
\end{eqnarray*}

And so on, you get the point, the idea being that we are led into a certain recurrence, with the $2k$-th moment of our distribution $\varphi(T^{2k})$ providing us with the needed value of the coefficient $a_k$. Thus, we are led to the conclusion in the statement.
\end{proof}

Before getting further, with free products of such models, let us work out a very basic example, which is something fundamental, that we will need in what follows:

\index{random walk}
\index{Dyck paths}
\index{Catalan numbers}

\begin{proposition}
In the context of the above correspondence, the variable
$$T=S+S^*$$
follows the Wigner semicircle law $\gamma_1=\frac{1}{2\pi}\sqrt{4-x^2}\,dx$.
\end{proposition}

\begin{proof}
In order to compute the law of variable $T$ in the statement, we can use the moment method. Indeed, the moments of this variable are as follows:
$$M_k
=\varphi(T^k)
=\varphi((S+S^*)^k)
=\#(1\in(S+S^*)^k)$$

Now since $S$ shifts to the right on $\mathbb N$, and $S^*$ shifts to the left, while remaining positive, we are left with counting the length $k$ paths on $\mathbb N$ starting and ending at 0. But this is something that we are familiar with, the number of such paths being $L_{2k}=C_k$, the Catalan numbers, and with the corresponding measure being the Wigner law $\gamma_1$.
\end{proof}

Getting back now to our linearization program for $\boxplus$, the next step is that of taking a free product of the model found in Theorem 6.21 with itself.  There are two approaches here, one being a bit abstract, and the other one being more concrete. We first have:

\index{semigroup algebra}
\index{shift}

\begin{proposition}
We can talk about semigroup algebras $C^*_{red}(G)\subset B(l^2(G))$, a bit as we did for the group algebras, and at the level of examples:
\begin{enumerate}
\item With $G=\mathbb N$ we recover the shift algebra $A=<S>$ on $H=l^2(\mathbb N)$.

\item With $G=\mathbb N*\mathbb N$, we obtain the algebra $A=<S_1,S_2>$ on $H=l^2(\mathbb N*\mathbb N)$.
\end{enumerate}
\end{proposition}

\begin{proof}
We can talk indeed about semigroup algebras $C^*_{red}(G)\subset B(l^2(G))$, with the semigroup elements $g\in G$ being now isometries, and with this in hand:

\medskip

(1) With $G=\mathbb N$ we recover the shift algebra $A=<S>$ on the Hilbert space $H=l^2(\mathbb N)$, the shift $S$ itself being the isometry associated to the element $1\in\mathbb N$.

\medskip

(2) With $G=\mathbb N*\mathbb N$ we recover the double shift algebra $A=<S_1,S_2>$ on the Hilbert space $H=l^2(\mathbb N*\mathbb N)$, the two shifts $S_1,S_2$ themselves being the isometries associated to two copies of the element $1\in\mathbb N$, one for each of the two copies of $\mathbb N$ which are present.
\end{proof}

In what follows we will rather use an equivalent, second approach to our problem, which is exactly the same thing, but formulated in a less abstract way, as follows:

\index{free Fock space}
\index{creation operator}

\begin{proposition}
We can talk about the algebra of creation operators $S_x:v\to x\otimes v$ on the free Fock space associated to a real Hilbert space $H$, given by 
$$F(H)=\mathbb C\Omega\oplus H\oplus H^{\otimes2}\oplus\ldots$$
and at the level of examples, we have:
\begin{enumerate}
\item With $H=\mathbb C$ we recover the shift algebra $A=<S>$ on $H=l^2(\mathbb N)$.

\item With $H=\mathbb C^2$, we obtain the algebra $A=<S_1,S_2>$ on $H=l^2(\mathbb N*\mathbb N)$.
\end{enumerate}
\end{proposition}

\begin{proof}
We can talk indeed about the algebra $A(H)$ of creation operators on the free Fock space $F(H)$ associated to a real Hilbert space $H$, with the remark that, in terms of the abstract semigroup notions from Proposition 6.23, we have:
$$A(\mathbb C^k)=C^*_{red}(\mathbb N^{*k})\quad,\quad 
F(\mathbb C^k)=l^2(\mathbb N^{*k})$$

As for the assertions (1,2) in the statement, these are both clear, either directly, or by passing via (1,2) from Proposition 6.23, which were both clear as well.
\end{proof}

The advantage with this latter model comes from the following result, from \cite{vo1}, which has a very simple formulation, without linear combinations or anything:

\index{creation operator}
\index{annihilation operator}
\index{vacuum vector}
\index{free Fock space}

\begin{proposition}
Given a real Hilbert space $H$, and two orthogonal vectors $x\perp y$, the corresponding creation operators $S_x$ and $S_y$ are free with respect to
$$tr(T)=<T\Omega,\Omega>$$
called state associated to the vacuum vector.
\end{proposition}

\begin{proof}
In standard tensor product notation for the elements of the free Fock space $F(H)$, the formula of a creation operator associated to a vector $x\in H$ is as follows:
$$S_x(y_1\otimes\ldots\otimes y_n)=x\otimes y_1\otimes\ldots\otimes y_n$$

As for the formula of the adjoint of this creation operator, called annihilation operator associated to the vector $x\in H$, this is as follows: 
$$S_x^*(y_1\otimes\ldots\otimes y_n)=<x,y_1>\otimes y_2\otimes\ldots\otimes y_n$$

We obtain from this the following formula, for any two vectors $x,y\in H$:
$$S_x^*S_y=<x,y>id$$

But with this formula in hand, the result follows by doing some elementary computations, in the spirit of those done for the group algebras, in the above.
\end{proof}

With this technology in hand, let us go back to our linearization program for $\boxplus$. We know from Theorem 6.21 how to model the individual distributions $\mu\in\mathcal P(\mathbb R)$. Thus, it remains to freely sum in the model, by using Proposition 6.25, and then compute the distribution. And the point here is that, still following Voiculescu \cite{vo1}, we have:

\index{freeness}
\index{rotation trick}

\begin{theorem}
Given two polynomials $f,g\in\mathbb C[X]$, consider the variables 
$$S^*+f(S)\quad,\quad 
T^*+g(T)$$
where $S,T$ are two creation operators, or shifts, associated to a pair of  orthogonal norm $1$ vectors. These variables are then free, and their sum has the same law as
$$R^*+(f+g)(R)$$
with $R$ being the usual shift on $l^2(\mathbb N)$.
\end{theorem}

\begin{proof}
We have two assertions here, the idea being as follows:

\medskip

(1) The freeness assertion comes from the general freeness result from Proposition 6.25, via the various identifications coming from the previous results.

\medskip

(2) Regarding the second assertion, the idea is that this comes from a $45^\circ$ rotation trick. Let us write indeed the two variables in the statement as follows:
$$X=S^*+a_0+a_1S+a_2S^2+\ldots$$
$$Y=T^*+b_0+b_1T+a_2T^2+\ldots$$

Now let us perform the following $45^\circ$ base change, on the real span of the vectors $s,t\in H$ producing our two shifts $S,T$, as follows:
$$r=\frac{s+t}{\sqrt{2}}\quad,\quad 
u=\frac{s-t}{\sqrt{2}}$$

The new shifts, associated to these vectors $r,u\in H$, are then given by:
$$R=\frac{S+T}{\sqrt{2}}\quad,\quad
U=\frac{S-T}{\sqrt{2}}$$

By using now these two new shifts, which are free according to Proposition 6.25, we obtain the following equality of distributions:
\begin{eqnarray*}
X+Y
&=&S^*+T^*+\sum_ka_kS^k+b_kT^k\\
&=&\sqrt{2}R^*+\sum_ka_k\left(\frac{R+U}{\sqrt{2}}\right)^k+b_k\left(\frac{R-U}{\sqrt{2}}\right)^k\\
&\sim&\sqrt{2}R^*+\sum_ka_k\left(\frac{R}{\sqrt{2}}\right)^k+b_k\left(\frac{R}{\sqrt{2}}\right)^k\\
&\sim&R^*+\sum_ka_kR^k+b_kR^k
\end{eqnarray*}

To be more precise, here at the end we have used the freeness property of $R,U$ in order to cut $U$ from the computation, as it cannot bring anything, and then we did a basic rescaling at the end. Thus, we are led to the conclusion in the statement.
\end{proof}

As a conclusion, the operation $\mu\to f$ from Theorem 6.21 linearizes $\boxplus$. In order to reach now to something concrete, we are left with a computation inside $C^*(\mathbb N)$, which is elementary, and whose conclusion is that $R_\mu=f$ can be recaptured from $\mu$ via the Cauchy transform $G_\mu$. The precise result here, due to Voiculescu \cite{vo1}, is as follows:

\index{Cauchy transform}
\index{R-transform}
\index{free convolution}
\index{free Fourier transform}

\begin{theorem}
Given a real probability measure $\mu$, define its $R$-transform as follows:
$$G_\mu(\xi)=\int_\mathbb R\frac{d\mu(x)}{\xi-x}\implies G_\mu\left(R_\mu(z)+\frac{1}{z}\right)=z$$
The free convolution operation is then linearized by this $R$-transform.
\end{theorem}

\begin{proof}
This can be done by using the above results, in several steps, as follows:

\medskip

(1) According to Theorem 6.26, the operation $\mu\to f$ in Theorem 6.21 linearizes the free convolution operation $\boxplus$. We are therefore left with a computation inside $C^*(\mathbb N)$. To be more precise, consider a variable as in Theorem 6.21, as follows, with $f\in\mathbb C[X]$:
$$X=S^*+f(S)$$

In order to establish the result, we must prove that the $R$-transform of $X$, constructed according to the procedure in the statement, is the polynomial $f$ itself.

\medskip

(2) In order to do so, we fix $|z|<1$ in the complex plane, and we set:
$$q_z=\delta_0+\sum_{k=1}^\infty z_k\delta_k$$

The shift operator and its adjoint act then on this vector as follows:
$$Sq_z=z^{-1}(q_z-\delta_0)\quad,\quad 
S^*q_z=zq_z$$

It follows that the adjoint of our operator $X$ acts on this vector as follows:
\begin{eqnarray*}
X^*q_z
&=&(S+f(S^*))q_z\\
&=&z^{-1}(q_z-\delta_0)+f(z)q_z\\
&=&(z^{-1}+f(z))q_z-z^{-1}\delta_0
\end{eqnarray*}

Now observe that the above formula can be written in the following way:
$$z^{-1}\delta_0=(z^{-1}+f(z)-X^*)q_z$$

The point now is that when $|z|$ is small, the operator appearing on the right is invertible. Thus, with this assumption, we can rewrite the above formula as follows:
$$(z^{-1}+f(z)-X^*)^{-1}\delta_0=zq_z$$

Now by applying the trace, we are led to the following formula:
\begin{eqnarray*}
tr\left[(z^{-1}+f(z)-X^*)^{-1}\right]
&=&\left<(z^{-1}+f(z)-X^*)^{-1}\delta_0,\delta_0\right>\\
&=&<zq_z,\delta_0>\\
&=&z
\end{eqnarray*}

(3) Let us apply now the procedure in the statement to the real probability measure $\mu$ modeled by our variable $X$. The Cauchy transform $G_\mu$ is given by:
$$G_\mu(\xi)
=tr((\xi-X)^{-1})
=\overline{tr\Big((\bar{\xi}-X^*)^{-1}\Big)}
=tr((\xi-X^*)^{-1})$$

Now observe that, with the choice $\xi=z^{-1}+f(z)$ for our complex variable, the trace formula found in (2) above tells us that we have:
$$G_\mu\big(z^{-1}+f(z)\big)=z$$

Thus, by definition of the $R$-transform, we have the following formula:
$$R_\mu(z)=f(z)$$

But this finishes the proof of the result, as explained before in step (1).
\end{proof}

Summarizing, the situation in free probability is quite similar to the one in classical probability, the product spaces needed for the basic properties of the Fourier transform being replaced by something ``noncommutative'', namely the free Fock space models.

\section*{6d. Central limits}

With the above linearization technology in hand, we can do many things. First, we have the following free analogue of the CLT, due to Voiculescu \cite{vo1}:

\index{FCLT}
\index{Free CLT}
\index{R-transform}

\begin{theorem}[Free CLT]
Given self-adjoint variables $f_1,f_2,f_3,\ldots$ which are f.i.d., centered, with variance $t>0$, we have, with $n\to\infty$, in moments,
$$\frac{f_1+\ldots+f_n}{\sqrt{n}}\sim\gamma_t$$
with the limiting measure being the Wigner semicircle law $\gamma_t$.
\end{theorem}

\begin{proof}
The $R$-transform of the variable on the left in the statement can be computed by using the linearization property from Theorem 6.27, and is given by:
$$R(z)
=nR_f\left(\frac{z}{\sqrt{n}}\right)
\simeq tz$$

Regarding now the right term, as explained in chapter 3, the Cauchy transform is:
$$G_{\gamma_t}(\xi)=\frac{\xi-\sqrt{\xi^2-4t}}{2t}$$

But this gives $R_{\gamma_t}(z)=tz$, according to the following computation:
\begin{eqnarray*}
G_{\gamma_t}(tz+z^{-1})
&=&\frac{tz+z^{-1}-\sqrt{(tz+z^{-1})^2-4t}}{2t}\\
&=&\frac{tz+z^{-1}+tz-z^{-1}}{2t}\\
&=&z
\end{eqnarray*}

Observe that $R_{\gamma_t}(z)=tz$ follows in fact as well from the following formula, coming from Proposition 6.22, and from the technical details of the $R$-transform:
$$S+S^*\sim\gamma_1$$

Thus, the laws in the statement have the same $R$-transforms, so they are equal.
\end{proof}

Regarding the limiting measures $\gamma_t$, that we already met in the previous chapters, in relation with the Wigner matrices, one problem that we were having was that of understanding how $\gamma_t$ exactly appears, out of $\gamma_1$. We can now solve this question:

\index{Wigner law}
\index{free convolution semigroup}

\begin{theorem}
The Wigner semicircle laws have the property
$$\gamma_s\boxplus\gamma_t=\gamma_{s+t}$$
so they form a $1$-parameter semigroup with respect to free convolution.
\end{theorem}

\begin{proof}
This follows either from Theorem 6.28, or from the fact, that we know from its proof, that the $R$-transform of $\gamma_t$, which is given by $R_{\gamma_t}(\xi)=t\xi$, is linear in $t$.
\end{proof}

As a continuation of this, let us discuss now the complex analogues of the above results. Following Voiculescu \cite{vo1}, we first have the following definition:

\index{circular law}
\index{Voiculescu law}

\begin{definition}
The Voiculescu circular law of parameter $t>0$ is given by
$$\Gamma_t=law\left(\frac{a+ib}{\sqrt{2}}\right)$$
where $a,b$ are free, each following the Wigner semicircle law $\gamma_t$.
\end{definition}

In other words, the passage $\gamma_t\to\Gamma_t$ is by definition entirely similar to the passage $g_t\to G_t$ from the classical case, by taking real and imaginary parts. As before in other similar situations, the fact that $\Gamma_t$ is indeed well-defined is clear from definitions. 

\bigskip

As a first result now on these circular laws, we have:

\index{free convolution semigroup}

\begin{theorem}
The Voiculescu circular laws have the property
$$\Gamma_s\boxplus\Gamma_t=\Gamma_{s+t}$$
so they form a $1$-parameter semigroup with respect to free convolution.
\end{theorem}

\begin{proof}
This comes indeed from the semigroup result from Theorem 6.29, by taking the real and imaginary parts of the various variables involved.
\end{proof}

Next, we have the following free analogue of the CCLT, also from \cite{vo1}:

\index{FCCLT}
\index{Free CCLT}

\begin{theorem}[Free CCLT]
Given random variables $f_1,f_2,f_3,\ldots$ which are f.i.d., centered, with variance $t>0$, we have, with $n\to\infty$, in moments,
$$\frac{f_1+\ldots+f_n}{\sqrt{n}}\sim\Gamma_t$$
where $\Gamma_t$ is the Voiculescu circular law of parameter $t$.
\end{theorem}

\begin{proof}
Let us decompose our variables into real and imaginary parts, as follows:
$$f_i=\frac{h_i+ik_i}{\sqrt{2}}$$

The variables $h_i$ and $k_i$ satisfy then the assumptions of the free CLT, so we have:
$$\frac{h_1+\ldots+h_n}{\sqrt{n}}\sim\gamma_t\quad,\quad 
\frac{k_1+\ldots+k_n}{\sqrt{n}}\sim\gamma_t$$

Now since the two limiting semicircle laws that we obtain in this way are free, their rescaled sum is circular, in the above sense, and this gives the result.
\end{proof}

Going ahead with more study of the Voiculescu circular variables, less trivial now is the computation of their moments. For our computations, we will need explicit models for these variables. Following \cite{vo1}, let us start with the following key result:

\index{shift}

\begin{theorem}
Let $H$ be the complex Hilbert space having as basis the colored integers $k=\circ\bullet\bullet\circ\ldots$\,, and consider the shift operators on this space: 
$$S:k\to\circ k\quad,\quad 
T:k\to\bullet k$$
We have then the following equalities of distributions,
$$S+S^*\sim\gamma_1\quad,\quad 
S+T^*\sim\Gamma_1$$
with respect to the state $\varphi(T)=<Te,e>$, with $e$ being the empty word.
\end{theorem}

\begin{proof}
This is standard combinatorics, the idea being as follows:

\medskip

(1) The first formula, namely $S+S^*\sim\gamma_1$, is something that we already know, from Proposition 6.22, and which played a key role in the proof of the CLT.

\medskip

(2) Regarding the second formula, $S+T^*\sim\Gamma_1$, the idea is that we can deduce this from $S+S^*\sim\gamma_1$, by using the freeness result from Proposition 6.25, and the rotation tricks from the proof of Theorem 6.26. To start with, due to $S+S^*\sim\gamma_1$, we have:
$$\frac{(S+S^*)+i(T+T^*)}{\sqrt{2}}\sim\Gamma_1$$

Next, we can rewrite this formula in the following more convenient way:
$$\frac{(S+iT)+(S^*+iT^*)}{\sqrt{2}}\sim\Gamma_1$$

In order to further process this latter formula, we can use a rotation trick, as in the proof of Theorem 6.26. Indeed, with the notations there, let us perform the following base change, on the span of the vectors $s,t\in H$ producing our two shifts $S,T$:
$$r=\frac{s+it}{\sqrt{2}}\quad,\quad 
u=\frac{s-it}{\sqrt{2}}$$

The new shifts, associated to these two vectors $r,u\in H$, are then given by:
$$R=\frac{S+iT}{\sqrt{2}}\quad,\quad
U=\frac{S-iT}{\sqrt{2}}$$

Observe now that in terms of $R,U$, our previous modeling formula for $\Gamma_1$ reads:
$$R+U^*\sim\Gamma_1$$

Now since $R,U$ are free shifts, their joint distribution is identical to the one of the original shifts $S,T$. We conclude that we have indeed $S+T^*\sim\Gamma_1$, as stated.
\end{proof}

We can now answer our moment question for the circular variables, as follows:

\index{circular variable}
\index{noncrossing pairings}

\begin{theorem}
A variable $a\in A$ follows the law $\Gamma_1$ precisely when
$$tr(a^k)=|\mathcal{NC}_2(k)|$$
for any colored integer $k=\circ\bullet\bullet\circ\ldots$
\end{theorem}

\begin{proof}
This is standard combinatorics, the idea being as follows:

\medskip

(1) By using Theorem 6.33, it is enough to do the computation in the model there. To be more precise, we can use the following explicit formulae for $S,T$:
$$S:k\to\circ k\quad,\quad 
T:k\to\bullet k$$

With these formulae in hand, our claim is that we have the following formula:
$$<(S+T^*)^ke,e>=|\mathcal{NC}_2(k)|$$

(2) In order to prove this formula, we can proceed as before, for the semicircle laws. Indeed, let us expand the quantity $(S+T^*)^k$, and then apply the state $\varphi$. With respect to the previous computation, for the semicircle laws, what happens now is that the contributions will come this time via the following formulae, which must succesively apply, as to collapse the whole product of $S,S^*,T,T^*$ variables into a 1 quantity:
$$S^*S=1\quad,\quad
T^*T=1$$

(3) As before, in our computation for the semicircle laws, these applications of the rules $S^*S=1$, $T^*T=1$ must appear in a noncrossing manner. But, in contrast with our previous computation for the semicircle laws, where $S+S^*$ was self-adjoint, what happens now is that at each point where the exponent $k$ has a $\circ$ entry we must use $T^*T=1$, and at each point where the exponent $k$ has a $\bullet$ entry we must use $S^*S=1$. 

\medskip

(4) Thus the contributions, which are each worth 1, are parametrized by the partitions $\pi\in\mathcal{NC}_2(k)$. We are therefore led to the moment formula in (1), as desired.
\end{proof}

More generally now, by rescaling, we have the following result:

\begin{theorem}
A variable $a\in A$ follows the law $\Gamma_t$ precisely when
$$tr(a^k)=t^{|k|/2}|\mathcal{NC}_2(k)|$$
for any colored integer $k=\circ\bullet\bullet\circ\ldots$
\end{theorem}

\begin{proof}
This follows indeed from Theorem 6.34, by rescaling. Alternatively, we can get this as well directly, by suitably modifying Theorem 6.33 first.
\end{proof}

Summarizing, we have good moment formulae for both the semicircular and circular variables. As a continuation of this, regarding the free families of semicircular and circular variables, called semicircular and circular systems, their combinatorics is something quite subtle, and there has been a lot of work developed in this direction. For more on all this, we refer to \cite{nsp}, \cite{vdn}. In what concerns us, we will be back to this on several occasions, and notably in chapter 12 below, when talking about free cumulants.

\bigskip

Getting back now to the case of the single variables, as in Theorem 6.35, the formula there has the following more conceptual interpretation:

\begin{theorem}
The moments of the Voiculescu laws are the numbers
$$M_k(\Gamma_t)=\sum_{\pi\in\mathcal{NC}_2(k)}t^{|\pi|}$$
with ``$\mathcal{NC}_2$'' standing for the noncrossing matching pairings.
\end{theorem}

\begin{proof}
This follows from the moment formula in Theorem 6.35. Indeed, we know from there that a variable $a\in A$ is circular, of parameter $t>0$, precisely when we have the following formula, for any colored integer $k=\circ\bullet\bullet\circ\ldots\,$:
$$tr(a^k)=t^{|k|/2}|\mathcal{NC}_2(k)|$$

Now since the number of blocks of a pairing $\pi\in\mathcal{NC}_2(k)$ is given by $|\pi|=|k|/2$, independently of the pairing, this formula can be written in the following way:
$$tr(a^k)=\sum_{\pi\in\mathcal{NC}_2(k)}t^{|\pi|}$$

Thus, we are led to the conclusion in the statement.
\end{proof}

All the above is quite nice, when compared with the similar results from the classical case. As a conclusion to these considerations, we can now formulate a global result regarding the classical and free, real and complex Gaussian laws, as follows:

\begin{theorem}
We have classical and free limiting theorems, as follows,
$$\xymatrix@R=48pt@C=30pt{
FCLT\ar@{-}[r]\ar@{-}[d]&FCCLT\ar@{-}[d]\\
CLT\ar@{-}[r]&CCLT
}\qquad\xymatrix@R=25pt@C=25pt{\\ :\\}
\qquad
\xymatrix@R=45pt@C=55pt{
\gamma_t\ar@{-}[r]\ar@{-}[d]&\Gamma_t\ar@{-}[d]\\
g_t\ar@{-}[r]&G_t
}$$
and the limiting laws $g_t,G_t,\gamma_t,\Gamma_t$ have the following properties:
\begin{enumerate}
\item They form semigroups with respect to the operations $*$ and $\boxplus$.

\item Their moments are $M_k=\sum_{\pi\in D(k)}t^{|\pi|}$, with $D=P_2,\mathcal P_2,NC_2,\mathcal{NC}_2$.
\end{enumerate}
\end{theorem}

\begin{proof}
The first assertion comes from our various limiting results, in classical and free probability, from chapter 1 and from here. As for (1) and (2), these come from our various semigroup and moment results, again from chapter 1 and from here.
\end{proof}

We will see later that in this book that the above squares can be substantially enlarged, with Poisson and Bessel laws and their free analogues, and with some shifts too.

\section*{6e. Exercises}

Welcome to free probability, and as exercises on all this, we have:

\begin{exercise}
Clarify the characterization of the GNS construction $(H,\pi,\xi)$.
\end{exercise}

\begin{exercise}
Further meditate on the various formalisms for quantum spaces.
\end{exercise}

\begin{exercise}
Work out all the details, for the $A*B$ models for freeness. 
\end{exercise}

\begin{exercise}
Work out all the details, for the $C^*(G*H)$ models for freeness. 
\end{exercise}

\begin{exercise}
Compute some further moments of the operators $T=S^*+f(S)$.
\end{exercise}

\begin{exercise}
Learn more about free Fock spaces, and their various uses.
\end{exercise}

\begin{exercise}
Fill in the details, for the moments of circular variables.
\end{exercise}

\begin{exercise}
Experiment a bit with the joint moments of circular variables.
\end{exercise}

As bonus exercise, have a look at Voiculescu's original paper \cite{vo1}.

\chapter{Circular variables}

\section*{7a. Gaussian matrices}

As an application of the semicircular and circular variable theory developed so far, and of free probability in general, let us go back now to the random matrices. Following Voiculescu's paper \cite{vo4}, we will prove a number of key freeness results for them, complementing the basic random matrix theory developed in chapter 3. As a first result, completing our asymptotic law study for the Gaussian matrices, we have:

\index{circular law}

\begin{theorem}
Given a sequence of complex Gaussian matrices, $Z_N\in M_N(L^\infty(X))$ having independent $G_t$ variables as entries, with $t>0$, we have
$$\frac{Z_N}{\sqrt{N}}\sim\Gamma_t$$
in the $N\to\infty$ limit, with the limiting measure being Voiculescu's circular law.
\end{theorem}

\begin{proof}
We know from chapter 3 that the asymptotic moments of the complex Gaussian matrices are given by the following formula:
$$M_k\left(\frac{Z_N}{\sqrt{N}}\right)\simeq t^{|k|/2}|\mathcal{NC}_2(k)|$$

But the numbers on the right being the moments of $\Gamma_t$, we obtain the result.
\end{proof}

The above result is of course something very nice, but remains of theoretical nature. Getting now to truly new results regarding the random matrices, by using the free probability technology, following Voiculescu's paper \cite{vo4}, we first have:

\index{Wigner matrix}
\index{asymptotic freeness}

\begin{theorem}
Given a family of sequences of Wigner matrices, $Z^i_N\in M_N(L^\infty(X))$, $i\in I$ with pairwise independent entries, following $G_t$, up to $Z_N^i=(Z_N^i)^*$, the matrices
$$\frac{Z^i_N}{\sqrt{N}}\in M_N(L^\infty(X))\quad,\quad i\in I$$
become with $N\to\infty$ semicircular, each following the Wigner law $\gamma_t$, and free.
\end{theorem}

\begin{proof}
This is something quite subtle, the idea being as follows:

\medskip

(1) First of all, we know from chapter 3 that for any $i\in I$ the corresponding sequence of rescaled Wigner matrices becomes semicircular in the $N\to\infty$ limit:
$$\frac{Z_N^i}{\sqrt{N}}\simeq\gamma_t$$ 

(2) Thus, what is new here, and that we have to prove, is the asymptotic freeness assertion. For this purpose we can assume that we are dealing with the case of 2 sequences of matrices, $|I|=2$. So, assume that we have Wigner matrices as follows:
$$Z_N,Z_N'\in M_N(L^\infty(X))$$

We have to prove that these matrices become asymptotically free, with $N\to\infty$.

\medskip

(3) Normally this is something that can be proved directly, via various routine computations with partitions, as we did before in this book, which simplify as usual in the $N\to\infty$ limit, and bring freeness. However, we can prove this as well by using a trick, based on the result in Theorem 7.1. Consider indeed the following random matrix:
$$Y_N=\frac{Z_N+iZ_N'}{\sqrt{2}}$$

This matrix is then complex Gaussian, and so by using Theorem 7.1, we obtain that in the limit $N\to\infty$, we have the following distribution estimate:
$$\frac{Y_N}{\sqrt{N}}\simeq\Gamma_t$$

(4) Now recall that the circular law $\Gamma_t$ was by definition the law of the following variable, with $a,b$ being semicircular, each following the law $\gamma_t$, and free:
$$c=\frac{a+ib}{\sqrt{2}}$$

We are therefore in the situation where the variable $(Z_N+iZ_N')/\sqrt{N}$, which has asymptotically semicircular real and imaginary parts, converges to the distribution of $a+ib$, equally having semicircular real and imaginary parts, and with these real and imaginary parts being free. Thus $Z_N,Z_N'$ become asymptotically free, as desired.
\end{proof}

Getting now to the complex case, we have a similar result here, as follows:

\index{Gaussian matrix}
\index{asymptotic freeness}

\begin{theorem}
Given a family of sequences of Gaussian matrices, $Z^i_N\in M_N(L^\infty(X))$, $i\in I$ with pairwise independent entries, each following $G_t$, with $t>0$, the matrices
$$\frac{Z^i_N}{\sqrt{N}}\in M_N(L^\infty(X))\quad,\quad i\in I$$
become with $N\to\infty$ circular, each following the Voiculescu law $\Gamma_t$, and free.
\end{theorem}

\begin{proof}
This follows indeed from Theorem 7.2, which applies to the real and imaginary parts of our complex Gaussian matrices, and gives the result.
\end{proof}

Still following the paper of Voiculescu \cite{vo4}, we have the following important addition to Theorems 7.2 and 7.3, which is something essential for applications:

\begin{theorem}
In the context of families of Gaussian or Wigner matrices
$$Z^i_N\in M_N(L^\infty(X))\quad,\quad i\in I$$
the rescalings $Z^i_N/\sqrt{N}$ become asymptotically free from any scalar diagonal matrices
$$D_N\in M_N(\mathbb C)$$
whose law converges with $N\to\infty$.
\end{theorem}

\begin{proof}
There are several proofs for this fact \cite{agz}, \cite{msp}, \cite{vo4}, but none is trivial. Here is the argument from Mingo-Speicher \cite{msp}, most in tune with the present book:

\medskip

(1) We can assume that we are dealing with Gaussian matrices, because this will give the result for the Wigner matrices too, by looking at the real and imaginary parts. Also, we can assume that we are dealing with a single sequence of such matrices. 

\medskip

(2) Thus, what we want to prove is that a sequence of rescaled Gaussian matrices $Y_N=Z_N/\sqrt{N}\in M_N(L^\infty(X))$ becomes asymptotically free from a sequence of scalar diagonal matrices $D_N\in M_N(\mathbb C)$, whose law converges, as in the statement.

\medskip

(3) For this purpose, let us look at a joint moment of $Y_N,D_N$. It is actually most convenient for the computation to drop the $N$ subscripts, and with this convention, we can write this joint moment as follows, with $a_1,\ldots,a_s$ being certain colored integers, of arbitrary lengths, and with $k=(k_1\ldots k_s)$ being a colored integer too:
$$M=\frac{1}{N}\cdot E\Big(Tr\big(D^{a_1}Y^{k_1}D^{a_2}Y^{k_2}\ldots D^{a_s}Y^{k_s}\big)\Big)$$

(4) So, let us compute this joint moment. To start with, we have:
\begin{eqnarray*}
M
&=&\frac{1}{N^{s/2+1}}\cdot E\Big(\sum_{i_1}\big(D^{a_1}Z^{k_1}D^{a_2}Z^{k_2}\ldots D^{a_s}Z^{k_s}\big)_{i_1i_1}\Big)\\
&=&\frac{1}{N^{s/2}+1}\cdot E\Big(\sum_{i_1\ldots i_s}D_{i_1i_1}^{a_1}(Z^{k_1})_{i_1i_2}D_{i_2i_2}^{a_2}(Z^{k_2})_{i_2i_3}\ldots D_{i_si_s}^{a_s}(Z^{k_s})_{i_si_1}\Big)\\
&=&\frac{1}{N^{s/2+1}}\sum_{i_1\ldots i_s}D_{i_1i_1}^{a_1}D_{i_2i_2}^{a_2}\ldots D_{i_si_s}^{a_s}E\Big((Z^{k_1})_{i_1i_2}(Z^{k_2})_{i_2i_3}\ldots(Z^{k_s})_{i_si_1}\Big)
\end{eqnarray*}

(5) Next, let us look at the expectation $E$ on the right. With the convention for indices $(ij)^\circ=ij,(ij)^\bullet=ji$, as in chapter 3, and by applying the Wick rule, we have:
\begin{eqnarray*}
E
&=&E\left[(Z_{(i_1i_2)^{k_1}})^{k_1}(Z_{(i_2i_3)^{k_2}})^{k_2}\ldots(Z_{(i_si_1)^{k_s}})^{k_s}\right]\\
&=&t^{s/2}\#\left\{\pi\in\mathcal P_2(k)\Big|\pi\leq\ker\left((i_1i_2)^{k_1},(i_2i_3)^{k_2},\ldots,(i_si_1)^{k_s}\right)\right\}
\end{eqnarray*}

(6) Now let us call $I$ the multi-index appearing on the right, so that the condition on the right reads $\pi\leq\ker I$. With this convention, and getting back to $M$, we have:
\begin{eqnarray*}
M
&=&\frac{t^{s/2}}{N^{s/2+1}}\sum_{i_1\ldots i_s}D_{i_1i_1}^{a_1}D_{i_2i_2}^{a_2}\ldots D_{i_si_s}^{a_s}\#\left\{\pi\in\mathcal P_2(k)\Big|\pi\leq\ker I\right\}\\
&=&\frac{t^{s/2}}{N^{s/2+1}}\sum_{\pi\in\mathcal P_2(k)}\sum_{\pi\leq\ker I}D_{i_1i_1}^{a_1}D_{i_2i_2}^{a_2}\ldots D_{i_si_s}^{a_s}
\end{eqnarray*}

(7) In order to understand what is going on, let us work out an explicit example, namely $k=\circ\bullet\bullet\,\circ$. Here we have two matching partitions, namely $\cap\cap$ and $\cap\hskip-1.5mm\cap$. For the partition $\pi=\cap\cap$ the above condition $\pi\leq\ker I$ can be processed as follows:
\begin{eqnarray*}
\cap\cap\leq\ker I
&\iff&\cap\cap\leq\ker((i_1i_2)^\circ,(i_2i_3)^\bullet,(i_3i_4)^\bullet,(i_4i_1)^\circ)\\
&\iff&\cap\cap\leq\ker((i_1i_2),(i_3i_2),(i_4i_3),(i_4i_1))\\
&\iff&(i_1i_2)=(i_3i_2),(i_4i_3)=(i_4i_1)\\
&\iff&i_1=i_3
\end{eqnarray*}

We conclude that the contribution of $\pi=\cap\cap$ to the moment $M$ is given by:
\begin{eqnarray*}
M_{\cap\cap}
&=&\frac{t^2}{N^3}\sum_{i_1i_2i_4}D_{i_1i_1}^{a_1}D_{i_2i_2}^{a_2}D_{i_1i_1}^{a_3}D_{i_4i_4}^{a_4}\\
&=&\frac{t^2}{N^3}\sum_{i_1}D_{i_1i_1}^{a_1}D_{i_1i_1}^{a_3}\sum_{i_2}D_{i_2i_2}^{a_2}\sum_{i_4}D_{i_4i_4}^{a_4}\\
&=&t^2\cdot tr(D^{a_1}D^{a_3})tr(D^{a_2})tr(D^{a_4})
\end{eqnarray*}

(8) As for the other partition, $\pi=\cap\hskip-1.5mm\cap$, here the condition $\pi\leq\ker I$ reads:
\begin{eqnarray*}
\cap\hskip-1.5mm\cap\leq\ker I
&\iff&\cap\hskip-1.5mm\cap\leq\ker((i_1i_2),(i_3i_2),(i_4i_3),(i_4i_1))\\
&\iff&(i_1i_2)=(i_4i_3),(i_3i_2)=(i_4i_1)\\
&\iff&i_1=i_2=i_3=i_4
\end{eqnarray*}

We conclude that the contribution of $\pi=\cap\hskip-1.5mm\cap$ to the moment $M$ is given by:
\begin{eqnarray*}
M_{\cap\hskip-1.2mm\cap}
&=&\frac{t^2}{N^3}\sum_{i_1}D_{i_1i_1}^{a_1}D_{i_1i_1}^{a_2}D_{i_1i_1}^{a_3}D_{i_1i_1}^{a_4}\\
&=&\frac{t^2}{N^2}\cdot tr(D^{a_1}D^{a_2}D^{a_3}D^{a_4})
\end{eqnarray*}

(9) Summarizing, some partitions, like $\cap\cap$, contribute to the asymptotic moment, while some other, like $\cap\hskip-1.5mm\cap$, do not. In general now, by performing an asymptotic analysis as in chapter 3, based on the fact that only the partitions $\pi\in\mathcal{NC}_2(k)$ are geodesic for $1-\gamma$, with $\gamma$ being the full cycle, we are led to the conclusion that only these latter partitions contribute, in the $N\to\infty$ limit, so that we have the following formula:
$$M
\simeq\frac{t^{s/2}}{N^{s/2+1}}\sum_{\pi\in\mathcal{NC}_2(k)}\sum_{\pi\leq\ker I}D_{i_1i_1}^{a_1}D_{i_2i_2}^{a_2}\ldots D_{i_si_s}^{a_s}$$

(10) Next, with our computation from (7) in mind, let us try now to find the exact contributions of the partitions $\pi\in\mathcal{NC}_2(k)$. And here, we are led to the following formula, with $tr_{\gamma\pi}$ standing for the operation consisting in tracing inside each cycle of $\gamma\pi$:
$$M_\pi=t^{s/2}\cdot tr_{\gamma\pi}\big(D^{a_1},\ldots,D^{a_s}\big)$$

To be more precise, recall from chapter 3 that we have an embedding $P(k)\subset S_k$, obtained by cycling inside each block. Thus, we can form the product $\gamma\pi\in S_k$, and then $tr_{\gamma\pi}$ means to take the product of traces of products, following the cycles of $\gamma\pi$.

\medskip

(11) As an illustration for this, let us go back to the partition $\pi=\cap\cap$, studied in (7). Via the standard embedding $P(4)\subset S_4$, obtained by cycling inside each block, this partition corresponds to the product of transpositions $\pi=(12)(34)$, so we have:
$$\xymatrix@R=10pt@C=15pt{
&&1\ar[dr]&2\ar[dl]&3\ar[dr]&4\ar[dl]&&1\ar[ddrr]&2\ar[dd]&3\ar[ddll]&4\ar[dd]\\
\gamma\pi&=&1\ar[dr]&2\ar[dr]&3\ar[dr]&4\ar[dlll]&=\\
&&1&2&3&4&&1&2&3&4}$$

Thus, the contribution from (10) coincides with the one found in (7), namely:
$$M_{\cap\cap}=t^2\cdot tr(D^{a_1}D^{a_3})tr(D^{a_2})tr(D^{a_4})$$

(12) Now by putting everything together, namely formula from (9), and individual contributions of the partitions computed in (10), we conclude that we have:
$$M\simeq t^{s/2}\sum_{\pi\in\mathcal{NC}_2(k)}tr_{\gamma\pi}\big(D^{a_1},\ldots,D^{a_s}\big)$$

Equivalently, with $d$ being a variable following the asymptotic law of the matrices $(D_N)$, we have the following exact formula, for our asymptotic moment:
$$M=t^{s/2}\sum_{\pi\in\mathcal{NC}_2(k)}tr_{\gamma\pi}\big(d^{a_1},\ldots,d^{a_s}\big)$$

(13) As a next computation now, which will finish the proof, with $d$ being as above, following the asymptotic law of the matrices $(D_N)$, consider as well an abstract circular variable $c\sim\Gamma_t$ free from it, and let us look at the joint moments of $c,d$:
$$M'=tr\Big(d^{a_1}c^{k_1}\,d^{a_2}c^{a_2}\ldots d^{a_s}c^{k_s}\Big)$$

(14) But computing this abstract joint moment $M'$ is a basic free probability question, only needing the basic freeness rule from chapter 6. And when doing the computation, we can see appearing the noncrossing partitions $\pi\in\mathcal{NC}_2(k)$, needed for having the non-vanishing of the components. In the end, we are led to the following formula, with $tr_{\gamma\pi}$ standing as before for the operation consisting in tracing inside each cycle of $\gamma\pi$:
$$M'=t^{s/2}\sum_{\pi\in\mathcal{NC}_2(k)}tr_{\gamma\pi}\big(d^{a_1},
\ldots,d^{a_s}\big)$$

(15) As an illustration, let us take as before $k=\circ\bullet\bullet\,\circ$. With the notation $x_i=d^{a_i}$, and with our usual convention $x'=x-tr(x)$, we have the following formula:
\begin{eqnarray*}
M'
&=&tr(d^{a_1}c^\circ d^{a_2}c^\bullet d^{a_3}c^\bullet d^{a_4}c^\circ)\\
&=&tr(x_1cx_2c^*x_3c^*x_4c)\\
&=&tr\left([x_1'+tr(x_1)]c[x_2'+tr(x_2)]c^*[x_3'+tr(x_3)]c^*[x_4'+tr(x_4)]c\right)
\end{eqnarray*}

Now when expanding into a sum of 16 terms, and regrouping, we obtain:
\begin{eqnarray*}
M'
&=&tr(x_1)tr(cx_2c^*x_3c^*x_4c)\quad\longleftarrow\ 8\ terms\\
&&+tr(x_2)tr(x_1'cc^*x_3c^*x_4c)\quad\longleftarrow\ 4\ terms\\
&&+tr(x_3)tr(x_1'cx_2'c^*c^*x_4c)\quad\longleftarrow\ 2\ terms\\
&&+tr(x_4)tr(x_1'cx_2'c^*x_3'c^*c)\quad\longleftarrow\ 1\ term\\
&&+tr(x_1'cx_2'c^*x_3'c^*x_4c)\quad\longleftarrow\ 1\ term
\end{eqnarray*}

(16) Let us analyse now the 5 components of our sum. The last one vanishes by freeness. The second to last vanishes by freeness too, by writing $c^*c=(c^*c)'+tr(c^*c)$. The third to last vanishes by freeness too, by writing $x_4=x_4'+tr(x_4)$. For the fourth to last component we must again expand and regroup, and we are led to:
\begin{eqnarray*}
tr(x_2)tr(x_1'cc^*x_3c^*x_4c)
&=&tr(x_2)tr(x_1'cc^*[x_3'+tr(x_3)]c^*[x_4'+tr(x_4)]c)\\
&=&tr(x_2)tr(x_3)tr(x_1'cc^*c^*x_4c)\quad\longleftarrow\ 2\ terms\\
&&+tr(x_2)tr(x_4)tr(x_1'cc^*x_3'c^*c)\quad\longleftarrow\ 1\ term\\
&&+tr(x_2)tr(x_1'cc^*x_3'c^*x_4'c)\quad\longleftarrow\ 1\ term\\
&=&0+tr(x_2)tr(x_4)tr(x_1'cc^*x_3'c^*c)+0\\
&=&tr(x_2)tr(x_4)tr(x_1'x_3')tr(cc^*)^2
\end{eqnarray*}

(17) Finally, for the first component of our sum in (15), we must again expand and regroup, by using the usual $x'=x-tr(x)$ trick, and we are led to:
\begin{eqnarray*}
tr(x_1)tr(cx_2c^*x_3c^*x_4c)
&=&tr(x_1)tr(x_2)tr(c^*x_3c^*x_4c^2)\quad\longleftarrow\ 4\ terms\\
&&+tr(x_1)tr(x_3)tr(x_2'c^*c^*x_4c^2)\quad\longleftarrow\ 2\ terms\\
&&+tr(x_1)tr(x_4)tr(x_2'c^*x_3'c^*c^2)\quad\longleftarrow\ 1\ term\\
&&+tr(x_1)tr(x_2'c^*x_3'c^*x_4'c^2)\quad\longleftarrow\ 1\ term
\end{eqnarray*}

(18) Now observe that in this sum, the last 3 components vanish by freeness. As for the first component, here we must again expand and regroup, and we are led to:
\begin{eqnarray*}
tr(x_1)tr(x_2)tr(c^*x_3c^*x_4c^2)
&=&tr(x_1)tr(x_2)tr(c^*[x_3'+tr(x_3)]c^*[x_4'+tr(x_4)]c^2)\\
&=&tr(x_1)tr(x_2)tr(x_3)tr(c^*c^*x_4c^2)\quad\longleftarrow\ 2\ terms\\
&&+tr(x_1)tr(x_2)tr(x_4)tr(c^*x_3'c^*c^2)\quad\longleftarrow\ 1\ term\\
&&+tr(x_1)tr(x_2)tr(c^*x_3'c^*x_4'c^2)\quad\longleftarrow\ 1\ term\\
&=&tr(x_1)tr(x_2)tr(x_3)tr(x_4)tr(c^*c^*c^2)
\end{eqnarray*}

(19) Now by putting everything together, the moment that we were computing is:
\begin{eqnarray*}
M'
&=&tr(x_2)tr(x_4)tr(x_1'x_3')tr(cc^*)^2+tr(x_1)tr(x_2)tr(x_3)tr(x_4)tr(c^*c^*c^2)\\
&=&tr(x_2)tr(x_4)[t^2tr(x_1'x_3')+t^2tr(x_1)tr(x_3)]\\
&=&t^2[tr(x_1'x_3')+tr(x_1)tr(x_3)]tr(x_2)tr(x_4)\\
&=&t^2[tr((x_1'+tr(x_1))(x_3'+tr(x_3))]tr(x_2)tr(x_4)\\
&=&t^2\cdot tr(x_1x_3)tr(x_2)tr(x_4)
\end{eqnarray*}

But, as explained in (11), this matches with the general formula in (14).

\medskip

(20) As a conclusion to all this, we have proved $M\simeq M'$, which reads:
$$\frac{1}{N}\cdot E\Big(Tr\big(D^{a_1}Y^{k_1}D^{a_2}Y^{k_2}\ldots D^{a_s}Y^{k_s}\big)\Big)\simeq tr\Big(d^{a_1}c^{k_1}\,d^{a_2}c^{a_2}\ldots d^{a_s}c^{k_s}\Big)$$

Now since we know that the input distributions match, $Y\sim c$ and $D\sim d$ in the $N\to\infty$ limit, we conclude that we have asymptotic freenesss for $Y,D$, as stated.
\end{proof}

Still following Voiculescu \cite{vo3}, as an illustration for the above random matrix technology, we have the following application, to basic questions in free probability:

\index{semicircular system}
\index{circular system}
\index{compression by projection}

\begin{theorem}
Assuming that $\{s_1,\ldots,s_n\}$ is a semicircular system, in the sense that $s_i\sim\gamma_1$ are free, and $p$ is a projection free from $s_1,\ldots,s_n$, the family
$$\left\{\frac{ps_1p}{\sqrt{tr(p)}},\ldots\ldots,\frac{ps_np}{\sqrt{tr(p)}}\right\}$$
is a semicircular system too. The same happens for the circular systems.
\end{theorem}

\begin{proof}
We can model indeed $s_1,\ldots,s_n$ by large Wigner matrices $Z_N^i$, and our projection $p$ by a sequence of diagonal matrices of type $P_N=diag(1,\ldots,1,0,\ldots,0)$, with the rank being chosen as for the trace to match asymptotically $tr(p)$. Now the point is that in the model, $P_NZ_N^iP_N$ are rescaled Wigner matrices, and this gives the result. As for the proof for the circular systems, this is similar, using Gaussian matrices.
\end{proof}

Finally, still following Voiculescu \cite{vo3}, we have yet another key application of the above random matrix technology, to basic questions in free probability, as follows:

\index{Haar unitary}

\begin{theorem}
Assume that $c$ is circular, $u$ is a Haar unitary, in the sense that
$$tr(u^k)=\delta_{k0}\quad,\quad\forall k\in\mathbb Z$$
and $c,u$ are free. Then the variable $cu$ is circular.
\end{theorem}

\begin{proof}
This is something which is not obvious to prove with bare hands, I mean try a bit and you will understand, but we can get it easily with our matrix models:

\medskip

(1) We can certainly model $c$ by large Gaussian matrices, as in Theorem 7.1. Regarding now $u$, our claim is that this can be modeled by the following matrices:
$$D_N=diag(w^i)\quad,\quad w=e^{2\pi i/N}$$

Indeed, we have the following moment computation, for these matrices:
$$M_k(D_N)=\frac{1}{N}\sum_iw^{ik}=\delta_{N|k}\simeq\delta_{k0}$$

(2) Thus we have our models for $c,u$, and by Theorem 7.4 the variables modeling $c,u$ are asymptotically free. Now observe that when making the product, we have:
$$(cu)_{ij}=c_{ij}u_{jj}=w^jc_{ij}$$

But these are Gaussian matrices, so we have $cu\sim\Gamma_1$ with $N\to\infty$, as desired.
\end{proof}

\section*{7b. Product questions}

Getting back now to the free probability basics, we would like to discuss a number of more advanced combinatorial results regarding the semicircular and circular laws, which are of multiplicative nature, and quite often have no classical counterpart.

\bigskip

In general, in order to deal with multiplicative questions for the free random variables, we are in need of results regarding the multiplicative free convolution operation $\boxtimes$. And here, in the classical case, to start with, the situation is very simple, as follows:

\begin{proposition}
Assuming that $fg=h$, with $f,g$ being independent, we have:
$$E(f^{ix})E(g^{ix})=E(h^{ix})$$
Thus the multiplicative convolution operation $\times$ is linearized by $\log E(f^{ix})$.
\end{proposition} 

\begin{proof}
We have indeed the following computation, using the Fourier transform:
\begin{eqnarray*}
fg=h
&\implies&\log f+\log g=\log h\\
&\implies&F_{\log f}(x)F_{\log g}(x)=F_{\log h}(x)\\
&\implies&E(e^{ix\log f})E(e^{ix\log g})=E(e^{ix\log h})\\
&\implies&E(f^{ix})E(g^{ix})=E(h^{ix})
\end{eqnarray*}

Alternatively, the result follows as well from a direct Fubini computation.
\end{proof}

In the free case now, let us recall from chapter 6 that we have the following result:

\index{multiplicative free convolution}

\begin{definition}
We have a free convolution operation $\boxtimes$, constructed as follows:
\begin{enumerate}
\item For abstract distributions, via $\mu_a\boxtimes\mu_b=\mu_{ab}$, with $a,b$ free.

\item For real positive measures, via $\mu_a\boxtimes\mu_b=\mu_{\sqrt{a}\,b\,\sqrt{a}}$, with $a,b$ positive and free.
\end{enumerate}
\end{definition}

All this is quite tricky, explained in chapter 6, the idea being that, while (1) is something straightforward, (2) is not, and comes by considering the variable $c=\sqrt{a}\,b\,\sqrt{a}$, which unlike $ab$ is always self-adjoint, and whose moments are given by:
\begin{eqnarray*}
tr(c^k)
&=&tr[(\sqrt{a}\,b\,\sqrt{a})^k]\\
&=&tr[\sqrt{a}\,baba\ldots bab\,\sqrt{a}]\\
&=&tr[\sqrt{a}\cdot \sqrt{a}\,baba\ldots bab]\\
&=&tr[(ab)^k]
\end{eqnarray*}

Quite remarkably, the above free multiplicative convolution operation $\boxtimes$ can be linearized too, with the result here, due to Voiculescu \cite{vo2}, being as follows:

\index{S-transform}
\index{multiplicative free convolution}

\begin{theorem}
The operation $\boxtimes$ can be linearized as follows:
\begin{enumerate}
\item Start with $\psi(z)=M_1z+M_2z^2+M_3z^3+\ldots$

\item Compute its inverse, $\psi(\chi(z))=z$.

\item Compute $S(z)=(1+z^{-1})\chi(z)$.

\item Then $\log S$ linearizes the free multiplicative convolution, $S_{\mu\boxtimes\nu}=S_\mu S_\nu$.
\end{enumerate}
\end{theorem}

\begin{proof}
There are several possible proofs here, with the one of Haagerup \cite{haa}, which is the most in tune with the present book, being as follows:

\medskip

(1) According to our conventions from chapter 6, we want to prove that, given noncommutative variables $a,b$ which are free, we have the following formula:
$$S_{\mu_{ab}}(z)=S_{\mu_a}(z)S_{\mu_b}(z)$$

(2) For this purpose, consider the orthogonal shifts $S,T$ on the free Fock space, from chapter 6. By using the algebraic arguments there, from the proof of the $R$-transform theorem, we can assume, a bit as there, that our variables are as follows:
$$a=(1+S)f(S^*)\quad,\quad 
b=(1+T)g(T^*)$$

Our claim, which will prove the theorem, is that we have the following formulae, for the $S$-transforms of the various variables involved:
$$S_{\mu_a}(z)=\frac{1}{f(z)}\quad,\quad 
S_{\mu_b}(z)=\frac{1}{g(z)}\quad,\quad 
S_{\mu_{ab}}(z)=\frac{1}{f(z)g(z)}$$

(3) Let us first compute $S_{\mu_a}$. We know that we have $a=(1+S)f(S^*)$, with $S$ being the shift on $l^2(\mathbb N)$. Given $|z|<1$, consider the following vector:
$$p=\sum_{k\geq0}z^ke_k$$ 

The shift and its adjoint act on this vector in the following way:
$$Sp=\sum_{k\geq0}z^ke_{k+1}=\frac{p-e_0}{z}\quad,\quad 
S^*p=\sum_{k\geq1}z^ke_{k-1}=zp$$

Thus $f(S^*)p=f(z)p$, and we deduce from this that we have:
$$ap
=(1+S)f(z)p
=f(z)\left(p+\frac{p-e_0}{z}\right)
=\left(1+\frac{1}{z}\right)f(z)p-\frac{f(z)}{z}\,e_0$$

By dividing everything by $(1+1/z)f(z)$, this formula becomes:
$$\frac{z}{1+z}\cdot\frac{1}{f(z)}\,ap=p-\frac{e_0}{1+z}$$

We can write this latter formula in the following way:
$$\left(1-\frac{z}{1+z}\cdot\frac{1}{f(z)}\,a\right)p=\frac{e_0}{1+z}$$

Now by inverting, we obtain from this the following formula:
$$\left(1-\frac{z}{1+z}\cdot\frac{1}{f(z)}\,a\right)^{-1}e_0=(1+z)p$$

(4) But this gives us the formula of $S_{\mu_a}$. Indeed, consider the following function:
$$\rho(z)=\frac{z}{1+z}\cdot\frac{1}{f(z)}$$

With this notation, the formula that we found in (3) becomes:
$$(1-\rho(z)a)^{-1}e_0=(1+z)p$$

By using this, in terms of $\varphi(T)=<Te_0,e_0>$, we obtain:
\begin{eqnarray*}
\varphi\left((1-\rho(z)a)^{-1}\right)
&=&<(1-\rho(z)a)^{-1}e_0,e_0>\\
&=&<(1+z)p,e_0>\\
&=&1+z
\end{eqnarray*}

Thus the above function $\rho$ is the inverse of the following function:
$$\psi(z)=\varphi\left(\frac{1}{1-za}\right)-1$$

But this latter function is the $\psi$ function from the statement, and so $\rho$ is the function $\chi$ from the statement, and we can finish our computation, as follows:
$$S_{\mu_a}(z)
=\frac{1+z}{z}\cdot\rho(z)
=\frac{1+z}{z}\cdot\frac{z}{1+z}\cdot\frac{1}{f(z)}
=\frac{1}{f(z)}$$

(5) A similar computation, or just a symmetry argument, gives $S_{\mu_b}(z)=1/g(z)$. In order to compute now $S_{\mu_{ab}}(z)$, we use a similar trick. Consider the following vector of $l^2(\mathbb N*\mathbb N)$, with the primes and double primes referring to the two copies of $\mathbb N$:
$$q=e_0+\sum_{k\geq1}(e_1'+e_1''+e_1'\otimes e_1'')^{\otimes k}$$

The adjoints of the shifts $S,T$ act as follows on this vector:
$$S^*q=z(1+T)q\quad,\quad T^*q=zq$$

By using these formulae, we have the following computation:
\begin{eqnarray*}
abq
&=&(1+S)f(S^*)(1+T)g(T^*)q\\
&=&(1+S)f(S^*)(1+T)g(z)q\\
&=&g(z)(1+S)f(S^*)(1+T)q
\end{eqnarray*}

In order to compute the last term, observe that we have:
$$S^*(1+T)q
=(S^*+S^*T)q
=S^*q
=z(1+T)q$$

We conclude that we have the following formula:
$$f(S^*)(1+T)q=f(z)(1+T)q$$

Now back to our computation, we have the following formula:
\begin{eqnarray*}
abq
&=&g(z)(1+S)f(z)(1+T)q\\
&=&f(z)g(z)(1+S)(1+T)q\\
&=&f(z)g(z)\left(\frac{1+z}{z}\cdot q-\frac{e_0}{z}\right)
\end{eqnarray*}

Now observe that we can write this formula as follows:
$$\left(1-\frac{z}{1+z}\cdot\frac{1}{f(z)g(z)}\cdot ab\right)q=\frac{e_0}{1+z}$$

By inverting, we obtain from this the following formula:
$$\left(1-\frac{z}{1+z}\cdot\frac{1}{f(z)g(z)}\cdot ab\right)^{-1}e_0=(1+z)q$$

(6) But this formula that we obtained is similar to the formula that we obtained at the end of (3) above. Thus, we can use the same argument as in (4), and we obtain:
$$S_{\mu_{ab}}(z)=\frac{1}{f(z)g(z)}$$

We are therefore done with the computations, and this finishes the proof.
\end{proof}

As an illustration for the method in Theorem 7.9, here are a few computations of $S$-transforms, along with the formulae of the corresponding $R$-transforms too:

\begin{theorem}
We have the following formulae of $R$ and $S$ transforms:
\begin{enumerate}
\item For $\delta_t$ we have $R=t$, $S=1/t$.

\item For $\gamma_t$ we have $R=tz$, $S=1/\sqrt{tz}$.

\item For $\gamma_t^t$ we have $R=t+tz$, $S=(-t+\sqrt{t^2+4tz})/2tz$.

\item For $b_t$ we have $R=(z-1+\sqrt{(1-z)^2+4tz})/2z$, $S=(1+z)/(t+z)$.
\end{enumerate}
\end{theorem}

\begin{proof}
This is very standard, by using the formulae in Theorem 7.9:

\medskip

(1) For $\delta_t$ the $R$-transform can be computed as follows:
$$G=\frac{1}{\xi-t}\ \implies\ z=\frac{1}{R+z^{-1}-t}\ \implies\ R=t$$

As for the corresponding $S$-transform, this can be computed as follows:
\begin{eqnarray*}
\psi=tz/(1-tz)
&\implies&z=t\chi/(1-t\chi)\\
&\implies&z^{-1}=(t\chi)^{-1}-1\\
&\implies&t\chi=1/(1+z^{-1})\\
&\implies&S=1/t
\end{eqnarray*}

(2) We know from chapter 6 that the $R$-transform of $\gamma_t$ is $R=tz$. As for the $S$-transform of $\gamma_t$, this is standard to compute too, as follows:
\begin{eqnarray*}
tG^2=\xi G-1
&\implies&tz^2(1+\psi)^2=\psi\\
&\implies&t\chi^2(1+z)^2=z\\
&\implies&tz^2S^2=z\\
&\implies&S=1/\sqrt{tz}
\end{eqnarray*}

(3) Regarding now $\gamma_t^t$, that we included here in view of the results from chapter 4 regarding $g_t^t$, the $R$-transform can be computed as follows:
\begin{eqnarray*}
tG^2=(\xi-t)G-1
&\implies&tz^2=(R+z^{-1}-t)z-1\\
&\implies&tz^2=Rz-tz\\
&\implies&R=t+tz
\end{eqnarray*}

As for the corresponding $S$-transform, this can be computed as follows:
\begin{eqnarray*}
tG^2=(\xi-t)G-1
&\implies&tz^2(1+\psi)^2=(1-zt)(1+\psi)-1\\
&\implies&t\chi^2(1+z)^2=(1-t\chi)(1+z)-1\\
&\implies&tz^2S^2=z(1-tS)\\
&\implies&S=\frac{-t+\sqrt{t^2+4tz}}{2tz}
\end{eqnarray*}

(4) The $R$-transform of $b_t=(1-t)\delta_0+t\delta_1$ can be computed as follows:
\begin{eqnarray*}
M=1+tz/(1-z)
&\implies&G=(t+\xi-1)/(\xi^2-\xi)\\
&\implies&(Rz+1)(Rz+1-z)=tz+Rz+1-z\\
&\implies&R(Rz+1-z)=t\\
&\implies&R=\frac{z-1+\sqrt{(1-z)^2+4tz}}{2z}
\end{eqnarray*}

As for the corresponding $S$-transform, this can be computed as follows:
\begin{eqnarray*}
M=1+tz/(1-z)
&\implies&\psi=tz/(1-z)\\
&\implies&z=t\chi/(1-\chi)\\
&\implies&\chi=z/(t+z)\\
&\implies&S=(1+z)/(t+z)
\end{eqnarray*}

Thus, we are led to the formulae in the statement.
\end{proof}

As an application now of the above $S$-transform technology, to some basic questions from operator theory, following the original paper of Voiculescu \cite{vo2}, we have the following result, solving the spectral measure question for a product of two free projections:

\index{free projections}

\begin{theorem}
Given free projections $p,q$, of traces $s,t$, we have
$$law(pqp)=c_0\delta_0+c_1\delta_1+\frac{\sqrt{(x-a)(b-x)}}{2\pi x(x-1)}\,dx$$
with $c_0=1-\min(s,t)$, $c_1=\max(s+t-1,0)$ and $a,b=s+t-2st\pm2\sqrt{st(1-s)(1-t)}$.
\end{theorem}

\begin{proof}
According to our formulae in Theorem 7.10, the $S$-transforms of $p,q$ are:
$$S_p(z)=\frac{1+z}{s+z}\quad,\quad S_q(z)=\frac{1+z}{t+z}$$

By using Theorem 7.9, we conclude that the $S$-transform of $pqp$ is given by:
$$S=\frac{(1+z)^2}{(s+z)(t+z)}$$

In terms of $\chi(z)=S(z)/(1+z^{-1})$, this formula takes the following form:
$$\chi=\frac{z(1+z)}{(s+z)(t+z)}$$

Next, by replacing $z\to\psi(z)$, so that $\chi\to z$, this formula takes the following form:
$$z=\frac{\psi(1+\psi)}{(s+\psi)(t+\psi)}$$

Now with $z=\xi^{-1}$, and by using $\psi(\xi^{-1})=\xi G(\xi)-1$, we obtain from this:
$$\xi^{-1}=\frac{(\xi G-1)\xi G}{(s+\xi G-1)(t+\xi G-1)}$$

We conclude that the Cauchy transform satisfies the following equation:
$$(\xi G+s-1)(\xi G+t-1)=\xi^2 G(\xi G-1)$$

In other words, the Cauchy transform satisfies the following equation:
$$\xi^2(1-\xi)G^2+\xi(\xi+s+t-2)G+(s-1)(t-1)=0$$

We are led in this way to the following formula, for the Cauchy transform:
$$G=\frac{2-\xi-s-t+\sqrt{\xi^2+(4st-2s-2t)\xi+(s-t)^2}}{2\xi(1-\xi)}$$

Equivalently, with $a,b=s+t-2st\pm2\sqrt{st(1-s)(1-t)}$ as in the statement, we have:
$$G=\frac{2-\xi-s-t+\sqrt{(\xi-a)(\xi-b)}}{2\xi(1-\xi)}$$

Now by Stieltjes inversion, we are led to the formula in the statement.
\end{proof}

As a final comment here, the technology in Theorem 7.9 is of genuine real nature, and cannot be applied as such to complex variables, like the circular ones. However, it is possible to further build on Theorem 7.9, with technical complex analogues. See \cite{nsp}.

\section*{7c. Semigroup models}

We would like to discuss now, following \cite{ba1}, some further models for the circular variables, using semigroups, with these generalizing the main algebraic models that we have, namely the group algebras and the free Fock spaces. Let us start with:

\index{semigroup}

\begin{definition}
We call semigroup a unital semigroup, embeddable into a group, $M\subset G$. For such a semigroup $M$, we use the notation 
$$M^{-1}=\left\{m^{-1}\Big|m\in M\right\}$$
regarded as a subset of some group $G$ containing $M$.
\end{definition}

As a first observation, the above embeddability assumption $M\subset G$ tells us that the usual group cancellation rules hold in $M$, namely:
$$ab=ac\implies b=c\quad,\quad 
ba=ca\implies b=c$$

Regarding the precise relation between $M$ and the various groups $G$ containing it, it is possible to talk here about the Grothendieck group $G$ associated to such a semigroup $M$. However, we will not need this in what follows, and use Definition 7.12 as such.

\bigskip

With the above definition in hand, we have the following construction, which unifies the main models that we have, namely the group algebras and the free Fock spaces:

\index{semigroup algebra}
\index{free Fock space}

\begin{proposition}
Given a semigroup $M$, we have an embedding of semigroups as follows, with the operators $\lambda(m)$ being isometries, but not necessarily unitaries:
$$\lambda:M\subset B(l^2(M))\quad,\quad \lambda(m)=[\delta_n\to\delta_{mn}]$$
We denote by $C_{red}^*(M)\subset B(l^2(M))$ the $C^*$-algebra generated by the image $\lambda(M)$ of this embedding, and we endow it with the state $\varphi(T)=<T\delta_1,\delta_1>$.
\end{proposition}

\begin{proof}
We have indeed the following computation, for any $m\in M$:
\begin{eqnarray*}
\lambda(m)^*(\delta_n)
&=&\sum_{x\in M}<\lambda(m)^*\delta_n,\delta_x>\delta_x\\
&=&\sum_{x\in M}<\delta_n,\lambda(m)\delta_x>\delta_x\\
&=&\sum_{x\in M}\delta_{n ,mx}\delta_x
\end{eqnarray*}

Thus we have the isometry property for the operators $\lambda(m)$, as shown by:
$$\lambda(m)^*\lambda(m)(\delta_p)
=\lambda(m)^*\delta_{mp}
=\sum_{x\in M}\delta_{mp,mx}\delta_x
=\delta_p$$

As for unitarity, this holds for groups, $M=G$, but not in general. Indeed, for $M=\mathbb N$ the operator $S=\lambda(1)$ is the shift, which is an isometry which is not a unitary.
\end{proof}

At the level of examples now, as announced above, we have:

\begin{proposition}
The construction $M\to C_{red}^*(M)$ is as follows:
\begin{enumerate}
\item For the discrete groups, $M=G$, we obtain in this way the usual reduced discrete group algebras $C_{red}^*(G)\subset B(l^2(G))$, with their canonical traces. 

\item For a free semigroup, $M=\mathbb N^{*I}$, we obtain the algebra of creation operators over the full Fock space over $\mathbb R^I$, with the state associated to the vacuum vector.
\end{enumerate}
\end{proposition}

\begin{proof}
All this is clear from definitions, with (1) being obvious, and (2) coming via our usual identifications for the free Fock spaces and related algebras.
\end{proof}

As a key observation now, enabling us to do some probability, we have:

\begin{proposition}
If $M\subset N$ are semigroups satisfying the condition 
$$M(N-M)=N-M$$
then for every family $\{a_i\}_{i\in I}$ of elements in $M$, we have the formula
$$\{\lambda_N(a_i)\}_{i\in I}\sim\{\lambda_M(a_i)\}_{i\in I}$$
as an equality of joint distributions, with respect to the canonical states.
\end{proposition}

\begin{proof}
Assuming $M\subset N$ we have $l^2(M)\subset l^2(N)$, and for $m,m'\in M$ we have:
$$\lambda_M(m)\delta_{m'}=\lambda_N(m)\delta_{m'}$$

Thus if we suppose $M(N-M)=N-M$, as in the statement, then we have:
$$\lambda_M(m)^*\delta_{m'}
=\sum_{x\in M}\delta_{m',mx}\delta_x
=\sum_{x\in N}\delta_{m',mx}\delta_x
=\lambda_N(m)^*\delta_{m'}$$

In particular, if $m_1,\ldots,m_k\in M$, and $\alpha_1,\ldots,\alpha_k$ are exponents in $ \{1,*\}$, then:
$$\lambda_M(m_1)^{\alpha_1}\ldots\lambda_M(m_k)^{\alpha_k}\delta_1
=\lambda_N(m_1)^{\alpha_1}\ldots\lambda_N(m_k)^{\alpha_k}\delta_1$$

We are therefore led to the conclusion in the statement.
\end{proof}

In order to further build on this, we will need the following technical notion:

\begin{definition}
Let $N$ be a semigroup, and consider the following order on it:
$$a\preceq b\iff b\in aN$$
We say that $N$ is in the class $E$ if it satisfies one of the following equivalent conditions:
\begin{enumerate}
\item For $\preceq$ every bounded subset is totally ordered.

\item $a\preceq c,b\preceq c\implies a\preceq b$ or $b\preceq a$.

\item $aN\cap bN\neq\emptyset\implies aN\subset bN$ or $bN\subset aN$.

\item $NN^{-1}\cap N^{-1}N=N\cup N^{-1}$.
\end{enumerate}
\end{definition}

In relation with this, let us introduce as well the following notion, which is something standard in the combinatorial theory of semigroups:

\index{prefix}
\index{code}

\begin{definition}
Let $(a_i)_{i\in I}$ be a family of elements in a semigroup $N$. 
\begin{enumerate}
\item We say that $(a_i)_{i\in I}$ is a code if the semigroup $M\subset N$ generated by the $a_i$ is isomorphic to $\mathbb N^{*I}$, via $a_i\to e_i$, and satisfies $M(N-M)=N-M$.

\item We say that $(a_i)_{i\in I}$ is a prefix if $a_i\in a_jN\implies i=j$, which means that the elements $a_i$ are not comparable via the order relation $\preceq$.
\end{enumerate}
\end{definition}

In the probabilistic setting, the above notion of code is of interest, due to the following result, which is to be related to the key formula $S+T^*\sim\Gamma_1$ from chapter 6:

\begin{proposition}
Assuming that $(a_i,b_i)_{i\in I}$ is a code, the family
$$\big(\lambda(a_i)+\lambda(b_i)^*\big)_{i\in I}$$ 
is a circular family, in the sense of free probability.
\end{proposition}

\begin{proof}
Let $(a_i,b_i)_{i\in I}$ be a code, and consider the following family:
$$\Big(\lambda_N(a_i),\lambda_N(b_i)\Big)_{i\in I}\in B(l^2(N))$$

By using Proposition 7.15, this family has the same distribution as a family of creation operators associated to a family of $2I$ orthonormal vectors, on the free Fock space:
$$\Big(\lambda_{\mathbb N^{*I}}(e_i),\lambda_{\mathbb N^{*I}}(f_i)\Big)_{i\in I}\in B(l^2(\mathbb N^{*I}))$$

Thus we obtain the result, via the standard Fock space theory from chapter 6.
\end{proof}

Summarizing, we are in need of codes. And here, the class $E$ introduced in Definition 7.16 and the notion of prefix from Definition 7.17 come into play, as follows:

\begin{proposition}
For a semigroup $N\in E$, a family 
$$(a_i)_{i\in I}\subset N$$
having at least two elements is a prefix if and only if it is a code.
\end{proposition}

\begin{proof}
We have two implications to be proved, as follows:

\medskip

(1) Let $(a_i)_{i\in I}$ be a code which is not a prefix, for instance because we have $a_i=a_jn$ with $i\neq j,n\in N$. Then $n$ is in the semigroup $M$ generated by the $a_k$ and $a_i=a_jn$ with $i\neq j$, so $M$ cannot be free, and this is a contradiction, as desired. 

\medskip

(2) Conversely, assume that $(a_i)_{i\in I}$ is a prefix and let, with $m\in N$:
$$A=a_{i_1}^{{\alpha}_1}\ldots a_{i_n}^{{\alpha}_n}m=a_{j_ 1}^{{\beta}_1}\ldots a_{j_s}^{ {\beta}_s}$$ 

We have then $a_{i_1}\preceq A$, $a _{j_1}\preceq A$, and so $i_1=j_1$. We can therefore simplify $A$ to the left by $a_{i_1}$. A reccurence on $\sum\alpha_ i$ shows then that we have $n\leq s$ and:
$$a_{i_k}=a_{j_k}\ ,\ \forall k\leq n\quad,\quad 
{\alpha}_k={\beta}_k\ ,\ \forall k<n\quad,\quad 
{\alpha}_n\leq{\beta}_n$$
$$m=a_{j_n}^{{\beta}_n-{\alpha}_n}a_{j_{n+1}}^{{\beta}_{n+1}}\ldots a_{j_s
}^{{\beta}_s}$$

We know that $m$ is in the semigroup $M=<a_i>$, so $M(N-M)=N-M$ holds. Moreover, for $m=1$ we obtain that we have $n=s$, $a_{j_k}=a_{i_k}$ and ${\alpha}_k={\beta}_k$ for any $k\leq n$. Thus the variables $a_i$ freely generate $M$, and so $(a_i)_{i\in I}$ is a code, as desired.
\end{proof}

Summarizing, we have some good freeness results, for our semigroups. Before getting into applications, let us discuss the examples. We have here the following result:

\begin{proposition}
The class $E$ has the following properties:
\begin{enumerate}
\item All groups are in $E$.

\item The positive parts of totally ordered abelian groups are in $E$.

\item If $G$ is a group and $M\in E$, then $M\times G\in E$.

\item If $A_1,A_2$ are in $E$, then the free product $A_1*A_2$ is in $E$.
\end{enumerate}
\end{proposition}

\begin{proof}
This is something elementary, whose proof goes as follows:

\medskip

(1) This is obvious, coming from definitions.

\medskip

(2) This is obvious too, because $M$ is here totally ordered by $\preceq$. 

\medskip

(3) Let $G$ be a group and $M\in E$. We have then, as desired:
\begin{eqnarray*}
&&(M\times G)(M\times G)^{-1}\cap (M\times G)^{-1}(M\times G)\\
&=&(M\times G)(M^{-1}\times G)\cap (M^{-1}\times G)(M\times G)\\
&=&(MM^{-1}\times G)\cap (M^{-1}M\times G)\\
&=&(MM^{-1}\cap M^{-1}M)\times G\\
&=&(M\cup M^{-1})\times G\\
&=&(M\times G)\cup (M^{-1}\times G)\\
&=&(M\times G)\cup (M\times G)^{-1}
\end{eqnarray*}

(4) Let $a,b,c\in A_1*A_2$ such that $ab=c$. We write, as reduced words:
$$a=x_1\ldots x_n\quad,\quad 
b=y_1\ldots y_m\quad,\quad
c=z_1\ldots z_p$$

Now let $s$ be such that the following equalities happen:
$$x_ny_1=1\quad,\quad \ldots\quad,\quad 
x_{n-s+1}y_s=1\quad,\quad
x_{n-s}y_{s+1}\neq 1$$

With this done, consider now the following element:
$$u=x_{n-s+1}\ldots x_n=(y_1\ldots y_s)^{-1}$$

We have then the following computation, based on the above formulae:
$$c=ab=x_1\ldots x_{n-s}y_{s+1}\ldots y_m$$

Now let $i\in\{ 1,2\}$ be such that $z_{n-s}\in A_i$. There are two cases:

\medskip

-- If $x_{n-s}\in A_1$ and $y_{s+1}\in A_2$ or if $x_{n-s}\in A_2$ and $y_{s+1}\in A_1$, then $x_1\ldots x_{n-s}y_{s+1}\ldots y_m$ is a reduced word. In particular, we have $x_1=z_1$, $x_2=z_2$, and so on up to $x_{n-s}=z_{n-s}$. Thus we have $a=z_1\ldots z_{n-s}u$, with $u$ invertible.

\medskip

-- If $x_{n-s},y_{s+1}\in A_i$ then $x_1=z_1$ and so on, up to $x_{n-s-1}=z_{n-s-1}$ and $x_{n-s}y_{s+1}=z_{n-s}$. In this case we have $a=z_1\ldots z_{n-s-1}x_{n-s}u$, with $u$ invertible.

\medskip

Now observe that in both cases we obtained that $a$ is of the form $z_1\ldots z_fxu$ for some $f$, with $u$ invertible and such that if $z_{f+1}\in A_i$, then there exists $y\in A_i$ such that:
$$xy=z_{f+1}$$

Indeed, we can take $f=n-s-1$ and $x=z_{n-s},y=1$ in the first case, and $x=x_{n-s},y=y_{s+1}$ in the second one. Suppose now that $A_1,A_2\in E$ and let $a,b,a',b'\in A_1*A_2$ such that $ab=a'b'$. Let $z_1\ldots z_p$ be the decomposition of $ab=a'b'$ as a reduced word. Then we can decompose our words, as above, in the following way:
$$a=z_1\ldots z_fxu\quad,\quad 
a'=z_1\ldots z_{f'}x'u'$$

We have to show that $a=a'm$ or that $a'=am$ for some $m\in A_1*A_2$. But this is clear in all three cases that can appear, namely $f<f'$, $f'<f$, $f=f'$.
\end{proof}

We can now formulate a main result about semigroup freeness, as follows:

\begin{theorem}
The following happen:
\begin{enumerate}
\item Given $M\subset N$, both in the class $E$, satisfying $M(N-M)=N-M$, any $x$ in the $*$-algebra generated by $\lambda (M)$ can be written as follows, with $p_i,q_i\in M$:
$$x=\sum_ia_i\lambda_N(p_i)\lambda_N(q_i)^*$$

\item Given $A,B\in E$, let $x$ be in the $*$-algebra generated by $\lambda_{A*B}(A)$, such that $\varphi(x)=0$. If $W_A,W_B$ are the words beginning by elements of $A,B$, then $x$ acts as:
$$l^2(W_B\cup\{1\})\to l^2(W_A)$$

\item Let $A,B\in E$. Then $\lambda_{A*B}(A)$ and $\lambda_{A*B}(B)$ are free.
\end{enumerate}
\end{theorem}

\begin{proof}
This comes from our various results above, the idea being as follows:

\medskip

(1) It is enough to prove this for elements of the form $x=\lambda(m)^*\lambda(n)$ with $m,n\in M$, because the general case will follow from this. In order to do so, observe that $x=\lambda(m)^*\lambda(n)$ is different from $0$ precisely when there exist $a,b\in N$ such that:
$$<\lambda(m)^*\lambda(n)\delta_a,\delta_b>\neq 0$$

Thus $na=mb$. But we know that there exists $c\in N$ with $n=mc$ or with $m=nc$. Moreover, as $M(N-M)=N-M$, it follows that $c\in M$. Thus $x=\lambda(m)^*\lambda (n)\neq 0$ implies that $x=\lambda(c)$ or $x=\lambda(c)^*$ with $c\in M$, and this finishes the proof.

\medskip

(2) We apply (1) with $M=A$ and $N=A*B$ for writing, with $p_i,q_i\in A$:
$$x=\sum_ia_i\lambda(p_i)\lambda (q_i)^*$$

In order to prove the result, consider the following quantity:
$$\varphi(\lambda(p_i)\lambda(q_i)^*)=\sum_x\delta_{1,p_ix}\delta_{1,q_ix}$$

This quantity is nonzero precisely when $p_i=q_i$ is invertible, and in this case:
$$\lambda(p_i)\lambda(q_i)^*=1$$

Now since we assumed $\varphi(x)=0$, it follows that we can write:
$$x=\sum a_i\lambda(p_i)\lambda (q_i)^*\quad,\quad 
\varphi(\lambda(p_i)\lambda (q_i)^*)=0$$

By linearity, it is enough to prove the result for $x=\lambda(p_i)\lambda (q_i)^*$. Let $m\in  W_B\cup\{1\}$ and assume that $x\delta_m\neq 0$. Then $\lambda(q_i)^*\delta_m\neq0$ implies that $m=q_ic$ for some word $c\in A*B$. As $q_i\in A$ and $m\in W_B\cup\{1\}$, it follows that $q_i$ is invertible. Now observe that:
$$p_iq_i^{-1}=1\implies\varphi(x)=1$$

We conclude that we have $x\delta_m=\delta_{p_iq_i^{-1}m}\in l^2(W_A)$, as desired.

\medskip

(3) This follows from (2). Indeed, let $P=x_n\ldots x_1$ be a product of elements in $ker(\varphi)$, such that $x_{2k}$ is in the $*$-algebra generated by $\lambda(B)$, and $x_{2k+1}$ is in the $*$-algebra generated by $\lambda(A)$. Then $x_1\delta_1\in l^2(W_A)$. Thus $x_2x_1\delta_1\in l^2(W_B)$, and so on, and by a reccurence, $P\delta_1$ is in $l^2(W_A)$ or in $l^2(W_B)$. But this implies $\varphi(P)=0$, as desired.
\end{proof}

As a main application of the above semigroup technology, we have:

\begin{theorem}
Consider a Haar unitary $u$, free from a semicircular $s$. Then 
$$c=us$$
is a circular variable.
\end{theorem}

\begin{proof}
Let us denote by $z$ the image of $1\in\mathbb Z$ and by $n$ the image of $1\in \mathbb N$ by the canonical embeddings into $\mathbb Z*\mathbb N$. Let $\lambda ={\lambda} _{\mathbb Z *\mathbb N}$. We know that $\mathbb Z *\mathbb N\in E$. Also $(zn,nz^{-1})$ is obviously a prefix, so it is a code. Thus, the following variable is circular:
$$c=\lambda (zn)+\lambda (nz^{-1})^*$$

The point now is that we have the following formula:
$$\lambda (zn)+\lambda (nz^{-1})^*=us$$

But this gives the result, in our model and so in general as well, because $u={\lambda}(z)$ is a Haar unitary, $s=\lambda(n)+\lambda (n)^*$ is semicircular, and $u$ and $s$ are  free. 
\end{proof}

\section*{7d. Polar decomposition}

Following Voiculescu \cite{vo3}, let us discuss now the $\chi$ variable problematics in the context of free probability. We will need the following standard fact:

\index{polar decomposition}
\index{partial isometry}

\begin{theorem}
Each $T\in B(H)$ decomposes as follows, $U$ being a partial isometry:
$$T=U|T|\quad,\quad |T|=\sqrt{T^*T}$$.
Also, assuming $T\in A\subset B(H)$, we have $U,|T|\in\bar{A}$, weak closure.
\end{theorem}

\begin{proof}
According to our definition of the modulus, $|T|=\sqrt{T^*T}$, we have:
$$<|T|x,|T|y>
=<x,|T|^2y>
=<x,T^*Ty>
=<Tx,Ty>$$

We conclude that the following linear application is well-defined, and isometric:
$$U:Im|T|\to Im(T)\quad,\quad 
|T|x\to Tx$$

Now by continuity we can extend this isometry $U$ into an isometry between certain Hilbert subspaces of $H$, as follows:
$$U:\overline{Im|T|}\to\overline{Im(T)}\quad,\quad 
|T|x\to Tx$$

Moreover, we can further extend $U$ into a partial isometry $U:H\to H$, by setting $Ux=0$, for any $x\in\overline{Im|T|}^{\,\perp}$, and with this convention, the result follows. As for the last assertion, regarding the weak closure, this is something which is clear too.
\end{proof}

We will be interested in the polar decomposition of normal variables $a\in(A,tr)$, with $A$ being a von Neumann algebra, and $tr$ being faithful. In the $g_1$ case, we have:

\begin{theorem}
The polar decomposition of standard normal variables is $n=eh$, with the variables $e,h$ being as follows:
\begin{enumerate}
\item $e$ is a symmetry, having moments $1,0,1,0,1,\ldots$

\item $h\sim\sqrt{2/\pi}\,e^{-x^2/2}dx$ on $[0,\infty)$.

\item $e,h$ are independent.
\end{enumerate}
\end{theorem}

\begin{proof}
It is enough to prove the result in a model of our choice, and the best choice here is the most straightforward model for the normal variables, namely:
$$n=x\in L^\infty(\mathbb R,g_1)$$

Indeed, in this model, the polar decomposition of the our normal variable $x=(x\to x)$ is then $n=eh$, with the variables $e,h$ being as follows:
$$e=sgn(x)\quad,\quad 
h=|x|$$

Now since $e$ is a symmetry having moments $1,0,1,0,1,\ldots\,$, and also $h$ follows the law in the statement, and finally $e,h$ are independent, this gives the result.
\end{proof}

Regarding now the semicircular variables, the situation is very similar:

\index{quarter circular}

\begin{theorem}
The polar decomposition of semicircular variables is $s=eq$, with the variables $e,q$ being as follows:
\begin{enumerate}
\item $e$ is a symmetry, having moments $1,0,1,0,1,\ldots$

\item $q\sim\frac{1}{\pi}\sqrt{4-x^2}\,dx$ on $[0,2]$.

\item $e,q$ are independent.
\end{enumerate}
\end{theorem}

\begin{proof}
This is similar to the proof of Theorem 7.24, by using the following model:
$$s=x\in L^\infty\Big([-2,2]\,,\,\gamma_1\Big)$$

Thus $e=sgn(x)$ and $h=|x|$, leading to the conclusions in the statement.
\end{proof}

Before going further with a study in the complex case, let us point out the fact that the laws of the variables $h,q$ above are not as trivial as they seem, as shown by:

\begin{theorem}
The moments of the chi law $\sqrt{2/\pi}\,e^{-x^2/2}dx$ are given by
$$M_{2k}=(2k)!!\quad,\quad M_{2k+1}=2^k\sqrt{\frac{2}{\pi}}\cdot k!$$
with $E=\sqrt{2/\pi}$ and $V=1-2/\pi$. As for the quarter-circular law $\frac{1}{\pi}\sqrt{4-x^2}\,dx$, here
$$M_{2k}=\frac{1}{k+1}\binom{2k}{k}
\quad,\quad M_{2k-1}=\frac{16^k}{\pi}\cdot\frac{(k-1)!k!}{(2k+1)!}$$
the expectation being $E=8/(3\pi)$, and the variance being $V=1-(8/(3\pi))^2$.
\end{theorem}

\begin{proof}
We know the even moments in both cases. Regarding now the odd moments, for the chi law, we have the following computation:
\begin{eqnarray*}
M_{2k+1}
&=&\sqrt{\frac{2}{\pi}}\int_0^\infty x^{2k+1}e^{-x^2/2}dx\\
&=&\sqrt{\frac{2}{\pi}}\int_0^\infty(\sqrt{2y})^{2k+1}e^{-y}\,\frac{dy}{\sqrt{2y}}\\
&=&\sqrt{\frac{2^{2k+1}}{\pi}}\int_0^\infty y^ke^{-y}dy\\
&=&2^k\sqrt{\frac{2}{\pi}}\cdot k!
\end{eqnarray*}

Observe that this computation works in fact for the even moments too, that is, works for any $k\in\mathbb Z/2$, with at the end $k!$ being replaced by the gamma function value $\Gamma(k+1)$. Regarding now the quarter-circular law, the computation here is as follows:
\begin{eqnarray*}
M_{2k-1}
&=&\frac{1}{\pi}\int_0^2x^{2k-1}\sqrt{4-x^2}\,dx\\
&=&\frac{1}{\pi}\int_0^{\pi/2}(2\cos t)^{2k-1}2\sin t\cdot 2\sin tdt\\
&=&\frac{2^{2k+1}}{\pi}\int_0^{\pi/2}\cos^{2k-1}t\sin^2 tdt\\
&=&\frac{2^{2k+1}}{\pi}\cdot\frac{(2k-1)!!2!!}{(2k+2)!!}\\
&=&\frac{2^{2k+1}}{\pi}\cdot\frac{2^{k-1}(k-1)!}{(2k+2)!/2^{k+1}(k+1)!}\\
&=&\frac{2^{4k+1}}{\pi}\cdot\frac{(k-1)!(k+1)!}{(2k+2)!}\\
&=&\frac{16^k}{\pi}\cdot\frac{(k-1)!k!}{(2k+1)!}
\end{eqnarray*}

Thus, we are led to the conclusions in the statement.
\end{proof}

Getting now to the complex normal variables, the result here is as follows:

\index{polar decomposition}
\index{circular variable}
\index{Haar unitary}
\index{quarter-circular}
\index{Rayleigh law}

\begin{theorem}
The polar decomposition of complex normal variables is $N=uH$, with the variables $u,H$ being as follows:
\begin{enumerate}
\item $u$ is a Haar unitary.

\item $H\sim 2xe^{-x^2}dx$.

\item $u,H$ are independent.
\end{enumerate}
\end{theorem}

\begin{proof}
This is something very classical, due to Lord Rayleigh, as follows:

\medskip

(1) Our first claim is that for a complex normal variable $N\sim G_1$, the density of the variable $H^2=|N|^2$ is standard exponential, given by the following formula:
$$e_1=e^{-x}dx$$

Indeed, by partial integration, the moments of this predicted density are given by:
$$M_k=\int_0^\infty x^ke^{-x}dx=k!$$

But these are the correct moments, matching with those computed in chapter 1. 

\medskip

(2) Next, for any function $\varphi:[0,\infty)\to\mathbb R$, we have the following computation:
\begin{eqnarray*}
E(\varphi(H))
&=&E\big(\varphi\big(\sqrt{H^2}\,\big)\big)\\
&=&\int_0^\infty\varphi(\sqrt{x})e^{-x}\,dx\\
&=&\int_0^\infty\varphi(y)e^{-y^2}\,2ydy\\
&=&2\int_0^\infty\varphi(y)ye^{-y^2}\,dy
\end{eqnarray*}

Thus the density of $H$ is the one in the statement, and the rest is clear.

\medskip

(3) Finally, let us mention that, in analogy with what we saw in Theorem 7.26 for the 1D chi law, the Rayleigh law, which is the 2D chi law, has many interesting properties. Next, in 3D we have the Maxwell-Boltzmann law, which is a key object too. For more on all this, all classical mathematics and physics, you can have a look at my book \cite{ba2}. 
\end{proof}

Regarding now the case of the circular variables, things here are quite tricky too, with the result, due to Voiculescu \cite{vo3}, being as follows:

\index{polar decomposition}
\index{circular variable}
\index{Haar unitary}
\index{quarter-circular}

\begin{theorem}
The polar decomposition of circular variables is $c=uq$, with the variables $u,q$ being as follows:
\begin{enumerate}
\item $u$ is a Haar unitary.

\item $q$ is quarter-circular.

\item $u,q$ are free.
\end{enumerate}
\end{theorem}

\begin{proof}
This is something that can be proved via various techniques, as follows:

\medskip

(1) The original proof, by Voiculescu in \cite{vo3}, used Gaussian random matrix models for the circular variables, based on the asymptotic freeness result from Theorem 7.4. More specifically, the result follows from the formula $cu$ = circular, from Theorem 7.6.

\medskip

(2) A second proof, obtained by pure combinatorics, in the spirit of Theorem 7.9, regarding the $S$-transform, is the one in given by Nica-Speicher in \cite{nsp}.

\medskip

(3) The third proof, from \cite{ba1}, is by suitably manipulating Theorem 7.22, as to replace the semicircular there by a quarter-circular. Consider indeed the following group:
$$G=\mathbb Z*(\mathbb Z\times\mathbb Z /2\mathbb Z)$$

Let $z,t,a$ be the images of the following elements, into this group $G$:
$$1\in\mathbb Z\quad,\quad
(1,\hat{0})\in\mathbb Z\times(\mathbb Z /2\mathbb Z)\quad,\quad 
(0,\hat{1})\in\mathbb Z\times(\mathbb Z /2\mathbb Z)$$ 

Let $u=\lambda_G (z)$, $d=\lambda_G (a)$ and choose a quarter-circular $q\in <\lambda_G(t)>$. Then $(q,d)$ are independent, so $dq$ is semicircular, and so $c=udq$ is circular, and:

\medskip

-- The modulus of $c$ is $q$, which is a quarter-circular.

\medskip

-- The polar part of $c$ is $ud$, which is obviously a Haar unitary.

\medskip

-- Consider the automorphism of $G$ which is the identity on $\mathbb Z\times\mathbb Z /2\mathbb Z $ and maps $z\to za$. This extends to a trace-preserving automorphism of $C_{red}^*(G)$ which maps:
$$u\to ud\quad,\quad q\to q$$

Since $u,q$ are free, it follows that $ud,q$ are free too, finishing the proof.
\end{proof}

And we will end our discussion with this. More on such things, and more specifically on the related $\chi^2$ problematics for the laws $g_1,G_1,\gamma_1,\Gamma_1$, in the next chapter.

\section*{7e. Exercises}

This was a quite technical chapter, and as exercises here, we have:

\begin{exercise}
Work out all details, for the asymptotic freeness of $Y_N,D_N$.
\end{exercise}

\begin{exercise}
Learn more about the semicircular and circular systems.
\end{exercise}

\begin{exercise}
Compute further $R$ and $S$ transforms, for laws of your choice.
\end{exercise}

\begin{exercise}
Further study the laws of products $pq$ of free projections.
\end{exercise}

\begin{exercise}
Review what we did in chapter 6, using semigroup models.
\end{exercise}

\begin{exercise}
Find generalizations of the $c=us$ result, using semigroups.
\end{exercise}

\begin{exercise}
Learn more about polar decomposition, and singular values too.
\end{exercise}

\begin{exercise}
Futher explore the $\chi$ problematics for the laws $g_1,G_1,\gamma_1,\Gamma_1$.
\end{exercise}

As bonus exercise, more difficult, explore the multivariable $S$-transforms.

\chapter{Poisson limits}

\section*{8a. Poisson limits}

As a direct continuation of what we did in chapter 6, namely free probability basics and limiting theorems, we would like to have as well a free analogue of the Poisson Limit Theorem (PLT). Although elementary from what we have, this was something not done by Voiculescu himself, and not appearing in the foundational book  \cite{vdn}, and only explained later, in the book of Hiai and Petz \cite{hpe}. The statement is as follows:

\index{FPLT}
\index{free PLT}
\index{Marchenko-Pastur law}
\index{free Poisson law}

\begin{theorem}[Free PLT]
The following limit converges, for any $t>0$,
$$\lim_{n\to\infty}\left(\left(1-\frac{t}{n}\right)\delta_0+\frac{t}{n}\,\delta_1\right)^{\boxplus n}$$
and we obtain the Marchenko-Pastur law of parameter $t$, 
$$\pi_t=\max(1-t,0)\delta_0+\frac{\sqrt{4t-(x-1-t)^2}}{2\pi x}\,dx$$
also called free Poisson law of parameter $t$.
\end{theorem}

\begin{proof}
This is very standard, by using the technology from chapter 6:

\medskip

(1) Consider the measure in the statement, under the convolution sign:
$$\eta=\left(1-\frac{t}{n}\right)\delta_0+\frac{t}{n}\,\delta_1$$

The Cauchy transform of this measure is easy to compute, given by:
$$G_\eta(\xi)=\left(1-\frac{t}{n}\right)\frac{1}{\xi}+\frac{t}{n}\cdot\frac{1}{\xi-1}$$

In order to prove the result, we want to compute the following $R$-transform:
$$R
=R_{\eta^{\boxplus n}}(z)
=nR_\eta(z)$$

According to the formula of $G_\eta$, the equation for this function $R$ is as follows:
$$\left(1-\frac{t}{n}\right)\frac{1}{1/z+R/n}+\frac{t}{n}\cdot\frac{1}{1/z+R/n-1}=z$$

By multiplying both sides by $n/z$, this equation can be written as:
$$\frac{t+zR}{1+zR/n}=\frac{t}{1+zR/n-z}$$

Now with $n\to\infty$ things simplify, and we obtain the following formula:
$$t+zR=\frac{t}{1-z}$$

Thus we have the following formula, for the $R$-transform that we are interested in:
$$R=\frac{t}{1-z}$$

(2) On the other hand, let us compute the $R$-transform of the Marchenko-Pastur law $\pi_t$. As explained in chapter 3, this law $\pi_t$ is something quite tricky, at arbitrary $t>0$, naturally appearing via Stieltjes inversion, from the following Cauchy transform:
\begin{eqnarray*}
G(\xi)
&=&\frac{1-t+\xi-\sqrt{(t-1-\xi)^2-4\xi}}{2\xi}\\
&=&\frac{1-t+\xi-\sqrt{(\xi-1-t)^2-4t}}{2\xi}
\end{eqnarray*}

To be more precise, the first formula is the one coming from random matrix computations, and the second formula, equivalent to it, leads via Stieltjes inversion to the density in the statement. For computing the $R$-transform, observe that $G=G(\xi)$ satisfies:
$$(2\xi G+t-1-\xi)^2=(t-1-\xi)^2-4\xi$$

By developing the square on the left, we are led to the following equation:
$$4\xi^2 G^2+4\xi G(t-1-\xi)=-4\xi$$

Now by simplifying, and then setting $\xi=R+z^{-1}$, we obtain in this way:
\begin{eqnarray*}
G^2\xi=G(\xi+1-t)-1
&\implies&z^2(R+z^{-1})=z(R+z^{-1}+1-t)-1\\
&\implies&z^2R+z=zR+z-tz\\
&\implies&zR=R-t\\
&\implies&R=t/(1-z)
\end{eqnarray*}

Thus, we are led to the conclusion in the statement.
\end{proof}

In practice now, most of the basic theory of $\pi_t$ was already done before, in chapter 3, when talking random matrices. One thing which was missing there, however, was that of understanding how the law $\pi_t$, with parameter $t>0$, exactly appears, out of $\pi_1$. We can now solve this question, with a simple answer to it, as follows:

\index{free convolution semigroup}

\begin{theorem}
The Marchenko-Pastur laws have the semigroup property
$$\pi_s\boxplus\pi_t=\pi_{s+t}$$
so in particular, they appear via the formula $\pi_t=\pi_1^{\boxplus t}$ from their $t=1$ instance,
$$\pi_1=\frac{1}{2\pi}\sqrt{4x^{-1}-1}\,dx$$
which $\pi_1$ itself is something quite conceptual, coming via $f\sim\gamma_1\implies f^2\sim\pi_1$.
\end{theorem}

\begin{proof}
Many things can be said here, the idea being as follows:

\medskip

(1) To start with, the semigroup formula $\pi_s\boxplus\pi_t=\pi_{s+t}$, which implies $\pi_t=\pi_1^{\boxplus t}$, follows either from Theorem 8.1, or from the fact that the $R$-transform of $\pi_t$, computed in the proof of Theorem 8.1, and given by $R(z)=t/(1-z)$, is linear in $t>0$.

\medskip

(2) Next, regarding the law $\pi_1$, the formula in the statement, $f\sim\gamma_1\implies f^2\sim\pi_1$, is something that we know from chapter 3, coming from the fact that the moments coincide, $M_{2k}(\gamma_1)=M_k(\pi_1)=C_k$, or from a density computation, or from random matrices too.

\medskip

(3) Still talking $\pi_1$, observe that this is something quite simple from the abstract calculus viewpoint too, as being a basic beta law. In fact, up to rescalings, the most basic beta laws are $\{u_1,e_1,\gamma_1,\pi_1\}$, with $u_1,e_1$ being the standard uniform and exponential laws, and with $\gamma_1,\pi_1$ being the standard Wigner and Marchenko-Pastur laws. We will be back to such things, numerology and mysticism of $\gamma_1,\pi_1$, later in this book.

\medskip

(4) As another comment on $\pi_1$, in view of the formula $f\sim\gamma_1\implies f^2\sim\pi_1$, this can be thought of as being the free analogue of the 1-dimensional chi squared law $\chi_1^2$, given by $f\sim g_1\implies f^2\sim\chi_1^2$. As for higher dimensions, here we can say that the free analogue of $\chi_n^2$ should be the law of $f_1^2+\ldots+f_n^2$, with $f_i\sim\gamma_1$, free, which is given by:
$$law(f_1^2+\ldots+f_n^2)=\underbrace{\pi_1\boxplus\ldots\boxplus\pi_1}_{n\ times}=\pi_n$$

Which again is an interesting piece of numerology, and more on this later.

\medskip

(5) Next, observe that, while the first formulae in the statement are analogous to the formulae $p_s*p_t=p_{s+t}$ and $p_t=p_1^{*t}$ for the usual Poisson laws, in what concerns $f\sim\gamma_1\implies f^2\sim\pi_1$ and its consequences, these are something totally new, of genuine free nature, and this because $p_1,g_1$ have nothing much to do with each other.

\medskip

(6) Still talking philosophy, probably most striking is the conclusion regarding the Wishart matrices, namely ``the asymptotic law of the Wishart matrices is a free Poisson law''. Who could have guessed, that Poisson has something to do, with all this.

\medskip

(7) Finally, you might perhaps wonder, how are we supposed to live and think, with all these miracles around. In answer, my personal take is that everything comes from partitions, and with $f\sim\gamma_1\implies f^2\sim\pi_1$ coming from $NC_2(2k)\simeq NC(k)$. I mean, having this latter isomorphism always in mind leads to peace and understanding.
\end{proof}

In view of what was said above, the reasonable thing to be done next would be to go back to the Wishart matrices, and restudy them, with Theorem 8.1 and the various considerations from chapter 7 in mind. To be more precise, we would be interested here in understanding how asymptotic freeness works, in the Wishart matrix context.

\bigskip

Before this, however, with numerology being something addictive, let us do some more, in the case where the parameter $t>0$ is arbitrary. And here, as a first observation, the passage $p_t\to\pi_t$, which can be understood by replacing a classical Poisson limit by a free Poisson limit, can be better understood in terms of moments, as follows:

\begin{theorem}
The Poisson laws $p_t$ and Marchenko-Pastur laws $\pi_t$, given by
$$p_t=e^{-t}\sum_k\frac{t^k}{k!}\,\delta_k\quad,\quad 
\pi_t=\max(1-t,0)\delta_0+\frac{\sqrt{4t-(x-1-t)^2}}{2\pi x}\,dx$$
have the following properties:
\begin{enumerate}
\item They appear via the PLT, and the free PLT.

\item They form semigroups with respect to $*$ and $\boxplus$.

\item Their transforms are $\log F_{p_t}(x)=t(e^{ix}-1)$, $R_{\pi_t}(z)=t/(1-z)$.

\item Their moments are $M_k=\sum_{\pi\in D(k)}t^{|\pi|}$, with $D=P,NC$.
\end{enumerate}
\end{theorem}

\begin{proof}
These are all results that we already know, the idea being as follows:

\medskip

(1) This comes from the PLT from chapter 1, and the FPLT from here.

\medskip

(2) This comes again from our results from chapter 1 and from here.

\medskip

(3) Again, this comes from our results from chapter 1, and from here. Observe that, by adding the laws $g_t,\gamma_t$ to the picture, we are led to a puzzle as follows:
$$\log F_{p_t}(x)=t(e^{ix}-1)\quad\longleftrightarrow\quad R_{\pi_t}(z)=t/(1-z)$$
$$\log F_{g_t}(x)=-tx^2/2\quad\longleftrightarrow\quad R_{\gamma_t}(z)=tz$$

To be more precise, the puzzle is that of connecting the above formulae, into something of type $\log F_{m_t}(x)=\Phi(R_{\mu_t}(z))$, as to reach to a bijection $m_t\leftrightarrow\mu_t$ between classical and free probability measures. Which is doable, enjoy, and we will be back to this later.

\medskip

(4) Here the moment formula for $p_t$ is something that we already know, from chapter 1, and the formula for $\pi_t$ is something that we know too, from chapter 3. Observe that, by adding $g_t,\gamma_t$ to the picture, this leads us again into some sort of puzzle, namely:
$$P\quad\longleftrightarrow\quad NC$$
$$P_2\quad\longleftrightarrow\quad NC_2$$

But the solution to this latter puzzle is plainly obvious, namely $NC_x=P_x\cap NC$, and $P_x=<NC_x,\slash\hskip-2.1mm\backslash>$. So, this moment assertion (4) will be our favorite one, among the assertions in the statement, and this latter puzzle will be our favorite puzzle too.
\end{proof}

Coming next, and focusing now on the moment formula (4) from Theorem 8.3, we can merge this with our previous continuous results, and we obtain in this way:

\index{noncrossing partitions}
\index{noncrossing pairings}

\begin{theorem}
The moments of the various central limiting measures, namely
$$\xymatrix@R=46.5pt@C=57pt{
\pi_t\ar@{-}[r]\ar@{-}[d]&\gamma_t\ar@{-}[r]\ar@{-}[d]&\Gamma_t\ar@{-}[d]\\
p_t\ar@{-}[r]&g_t\ar@{-}[r]&G_t
}$$
are always given by the same formula, involving partitions, namely
$$M_k=\sum_{\pi\in D(k)}t^{|\pi|}$$
with the sets of partitions $D(k)$ in question being respectively
$$\xymatrix@R=45.7pt@C=45pt{
NC\ar@{-}[r]\ar@{-}[d]&NC_2\ar@{-}[r]\ar@{-}[d]&\mathcal{NC}_2\ar@{-}[d]\\
P\ar@{-}[r]&P_2\ar@{-}[r]&\mathcal P_2
}$$
and where $|.|$ is the number of blocks. 
\end{theorem}

\begin{proof}
This follows indeed by putting together the various moment results that we have, obtained throughout this book, in chapters 1, 3, 6 and here.
\end{proof}

As a comment now, as explained in chapters 1 and 4, with a bit of group theory inspiration in mind, there are a number of further measures that can be added to the picture, namely the Bessel laws $p^s_t$ and the shifts $g_t^t,G_t^t$. But then, we can conjecturally say that these should have free analogues $\pi^s_t$ and $\gamma_t^t,\Gamma_t^t$, fitting into the picture. 

\bigskip

More on this later, and for closing this numerology discussion, let us record:

\begin{question}
How far can we go with $M_k=\sum_{\pi\in D(k)}t^{|\pi|}$, and at the level of the main examples, can we reach to an almighty cube? 
\end{question}

And we will be back to this later, with an answer involving classical and free cumulants, classical and free quantum groups, and many more. Not to forget of course the functional transform speculations from the proof of Theorem 8.3, which are something quite general, potentially going beyong the partition setting, that we will talk about too.

\bigskip

Back to Earth now, as a most pressing task, let us rediscuss the Wishart matrices, as a continuation of what we did in chapter 3. We first have the following result:

\begin{theorem}
Given a family of sequences of complex Wishart matrices, 
$$Z^i_N=Y^i_N(Y^i_N)^*\in M_N(L^\infty(X))\quad,\quad i\in I$$
with each $Y^i_N$ being a $N\times N$ matrix, with entries following the normal law $G_t$, and with all these entries being pairwise independent, the rescaled sequences of matrices
$$\frac{Z^i_N}{tN}\in M_N(L^\infty(X))\quad,\quad i\in I$$
become with $N\to\infty$ Marchenko-Pastur, each following the law $\pi_1$, and free.
\end{theorem}

\begin{proof}
Here the first assertion is the Marchenko-Pastur theorem at $t=1$, from chapter 3, and the second assertion follows from the asymptotic freeness result from chapter 7, for the complex Gaussian matrices, in the obvious way.
\end{proof}

We have as well asymptotic freeness with the scalar diagonal matrices:

\begin{theorem}
In the context of families of complex Wishart matrices
$$Z^i_N\in M_N(L^\infty(X))\quad,\quad i\in I$$
the matrices $Z_N^i/tN$ become asymptotically free from any scalar diagonal matrices
$$D_N\in M_N(\mathbb C)$$
whose law converges with $N\to\infty$.
\end{theorem}

\begin{proof}
Again, this follows from the asymptotic freeness results from chapter 7, for the Gaussian  matrices, with respect to the scalar matrices, in the obvious way.
\end{proof}

In order to discuss now the parametric case, we will need a slight extension of our result from chapter 7 regarding the compressions of circular systems, as follows:

\begin{theorem}
Assuming that $\{s_1,\ldots,s_n\}$ is a circular system, free from an algebra $D$, and that $p\in D$ is a projection, the family 
$$\left\{\frac{ps_1p}{\sqrt{tr(p)}},\ldots\ldots,\frac{ps_np}{\sqrt{tr(p)}}\right\}$$
is a circular system, free from $pDp$. The same happens for the semicircular systems.
\end{theorem}

\begin{proof}
This is again a result from \cite{vo3}, the idea being as follows:

\medskip

(1) We know from chapter 7 that the result holds for $D=<p>$, as a consequence of our asymptotic freeness results there, the idea being that we can model $s_1,\ldots,s_n$ by large Gaussian matrices $Z_N^i$, and $p$ by matrices of type $P_N=diag(1,\ldots,1,0,\ldots,0)$, with the rank being chosen as for the trace to match asymptotically $tr(p)$. 

\medskip

(2) But with this in hand, the extension to the case of arbitrary algebras $D\supset<p>$ is straightforward. Indeed, given $a\in D$, we can come up with an extra Gaussian matrix $Z_N$, and with functions $f_N$ such that the matrices $f_N(Z_N)$ model the law of $a$ with $N\to\infty$, and it follows that our circular system is free from $pap$, as desired.

\medskip

(3) Finally, the assertion regarding semicircular systems follows from this, by taking real and imaginary parts, or can be proved in a similar way, using Wigner matrices.
\end{proof}

We can now generalize Theorems 8.6 and 8.7 to the parametric case, as follows:

\begin{theorem}
Given $t\geq1$ and a family of sequences of complex Wishart matrices, 
$$Z^i_N=Y^i_N(Y^i_N)^*\in M_N(L^\infty(X))\quad,\quad i\in I$$
with each $Y^i_N$ being a $N\times M$ matrix, with entries following the normal law $G_1$, and with all these entries being pairwise independent, the rescaled sequences of matrices
$$\frac{Z^i_N}{tN}\in M_N(L^\infty(X))\quad,\quad i\in I$$
become with $M=tN\to\infty$ Marchenko-Pastur, following the law $\pi_t$, and free. Moreover, the matrices $Z_N^i/tN$ become asymptotically free from any scalar diagonal matrices
$$D_N\in M_N(\mathbb C)$$
whose law converges with $N\to\infty$.
\end{theorem}

\begin{proof}
Here the first assertion is the Marchenko-Pastur theorem at $t\geq1$, from chapter 3, and the other assertions, regarding freeness, can be proved as follows:

\medskip

(1) Since we are in the case $t\geq1$, and so $M\geq N$, we can complete the $N\times M$ Gaussian matrices $Y_N^i$ into square matrices, that is, into $M\times M$ Gaussian matrices $T_M^i$, and our Wishart $N\times N$ matrix will not change, and more specifically, will appear as the upper left $N\times N$ block of the big $M\times M$ Wishart matrix $W_M^i=T_M^i(T_M^i)^*$.

\medskip

(2) But this means $Z_N^i=P_MW_M^iP_M$, with $P_M=diag(1,\ldots,1,0,\ldots,0)$ being the projection on the first $N$ coordinates. Now since these projections $P_M$ converge in law, to the law of a projection having trace $1/t$, by our results in chapter 7, we conclude that the Gaussian matrices $T_M^i$ are asymptotically free from these projections $P_M$.

\medskip

(3) Now with this freeness property in hand, the compression result for circular systems from chapter 7 applies, and shows that the matrices $\sqrt{t}\cdot P_MT_M^iP_M$ form in the $M\to\infty$ limit a circular system, and in particular, are asymptotically free. 

\medskip

(4) But the original Wishart matrices appear as $Z_N^i=(P_MT_M^iP_M)(P_MT_M^iP_M)^*$, and we conclude that these Wishart matrices are asymptotically free, as stated.

\medskip

(5) Finally, the last assertion, regarding the asymptotic freeness from the diagonal scalar matrices, follows as well, as an application of Theorem 8.8.
\end{proof}

\section*{8b. Compound Poisson}

Moving ahead now, let us try to find a free analogue of the CPLT. We will follow the CPLT material from chapter 1, by performing modifications where needed, as to replace everywhere classical probability with free probability. Let us start with the following straightforward definition, similar to the one from the classical case:

\index{compound Poisson law}

\begin{definition}
Associated to any compactly supported positive measure $\nu$ on $\mathbb R$ is the probability measure
$$\pi_\nu=\lim_{n\to\infty}\left(\left(1-\frac{c}{n}\right)\delta_0+\frac{1}{n}\,\nu\right)^{\boxplus n}$$
where $c=mass(\nu)$, called compound free Poisson law.
\end{definition}

In what follows we will be mostly interested in the case where $\nu$ is discrete, as is for instance the case for the measure $\nu=t\delta_1$ with $t>0$, which produces the free Poisson laws. However, we will sometimes use continuous measures too, with the convergence of the above limit coming from a standard approximation argument, based on our formulae below. Finally, we will talk later in this chapter about the complex case too.

\bigskip

The following result allows one to detect the compound free Poisson laws:

\index{R-transform}

\begin{proposition}
For a discrete measure, written as $\nu=\sum_{i=1}^sc_i\delta_{x_i}$ with $c_i>0$ and $x_i\in\mathbb R$, we have the following formula,
$$R_{\pi_\nu}(z)=\sum_{i=1}^s\frac{c_ix_i}{1-x_iz}$$
where $R$ denotes as usual the Voiculescu $R$-transform.
\end{proposition}

\begin{proof}
In order to prove this result, let $\eta_n$ be the measure appearing in Definition 8.10, under the free convolution sign. The Cauchy transform of $\eta_n$ is given by:
$$G_{\eta_n}(\xi)=\left(1-\frac{c}{n}\right)\frac{1}{\xi}+\frac{1}{n}\sum_{i=1}^s\frac{c_i}{\xi-x_i}$$

Consider now the $R$-transform of the measure $\eta_n^{\boxplus n}$, which is given by:
$$R_{\eta_n^{\boxplus n}}(z)=nR_{\eta_n}(z)$$

By using the general theory of the $R$-transform, from chapter 6, the above formula of $G_{\eta_n}$ shows that the equation for $R=R_{\eta_n^{\boxplus n}}$ is as follows:
\begin{eqnarray*}
&&\left(1-\frac{c}{n}\right)\frac{1}{1/z+R/n}+\frac{1}{n}\sum_{i=1}^s\frac{c_i}{1/z+R/n-x_i}=z\\
&\implies&\left(1-\frac{c}{n}\right)\frac{1}{1+zR/n}+\frac{1}{n}\sum_{i=1}^s\frac{c_i}{1+zR/n-x_iz}=1
\end{eqnarray*}

Now multiplying by $n$, then rearranging the terms, and letting $n\to\infty$, we get:
\begin{eqnarray*}
\frac{c+zR}{1+zR/n}=\sum_{i=1}^s\frac{c_i}{1+zR/n-x_iz}
&\implies&c+zR_{\pi_\nu}(z)=\sum_{i=1}^s\frac{c_i}{1-x_iz}\\
&\implies&R_{\pi_\nu}(z)=\sum_{i=1}^s\frac{c_ix_i}{1-x_iz}
\end{eqnarray*}

Thus, we are led to the conclusion in the statement.
\end{proof}

We have as well the following result, providing an alternative to Definition 8.10, and which, together with Definition 8.10, can be thought of as being the free CPLT:

\index{CFPLT}
\index{Compound FPLT}

\begin{theorem}
For a discrete measure, written as $\nu=\sum_{i=1}^sc_i\delta_{x_i}$ with $c_i>0$ and $x_i\in\mathbb R$, we have the formula
$$\pi_\nu=law\left(\sum_{i=1}^sx_i\alpha_i\right)$$
where the variables $\alpha_i$ are free Poisson $(c_i)$, free.
\end{theorem}

\begin{proof}
Let $\alpha$ be the sum of free Poisson variables in the statement. We can compute the $R$-transform of this variable as follows:
$$R_{\alpha_i}(z)=\frac{c_i}{1-z}
\implies R_{x_i\alpha_i}(z)=\frac{c_ix_i}{1-x_iz}
\implies R_\alpha(z)=\sum_{i=1}^s\frac{c_ix_i}{1-x_iz}$$

Thus we have the same formula as in Proposition 8.11, and we are done.
\end{proof}

As an illustration for this, following \cite{bn1}, let us study now in detail the compound free Poisson measures appearing from Bernoulli laws at $\pm1$. Let us start with:

\index{free difference}

\begin{theorem}
Given $m>0,n\geq 1$, the compound free Poisson law
$$\mu_{mn}=\pi_\nu\quad,\quad\nu=\frac{m(n-1)}{2}\,\delta_{-1}+\frac{m(n+1)}{2}\,\delta_1$$
coming from a generic Bernoulli law at $\pm1$, appears more concretely as follows,
$$a\sim\pi_{m(n+1)/2}\ ,\ b\sim\pi_{m(n-1)/2}\ ,\ {\rm free}\quad\implies\quad a-b\sim\pi_{mn}$$
that is, as the free difference of free Poisson laws of parameters $m(n\pm 1)/2$.
\end{theorem}

\begin{proof}
Consider two free variables $a,b$, following free Poisson laws of parameters $s,t$. According to our formulae before, the corresponding $R$-transforms are:
$$R_a=\frac{s}{1-z}\quad,\quad R_b=\frac{t}{1-z}$$

By using the general dilation formula $R_{qb}(z)=qR_b(qz)$  at $q=-1$, we deduce that we have $R_{-b}=-t/(1+z)$. Thus, our two formulae above can be written as follows:
$$R_a=\frac{s}{1-z}\quad,\quad R_{-b}=-\frac{t}{1+z}$$

Now when making the sum, the corresponding $R$-transform will be given by:
$$R_{a-b}=\frac{s}{1-z}-\frac{t}{1+z}$$

But with $s=m(n+1)/2$ and $t=m(n-1)/2$, this latter formula reads:
$$R_{a-b}=\frac{m}{2}\left(\frac{n+1}{1-z}-\frac{n-1}{1+z}\right)$$

On the other hand, for the measure $\nu$ in the statement, we have:
$$R_{\pi_\nu}(z)=\frac{m(n-1)}{2}\cdot\frac{-1}{1+z}+\frac{m(n+1)}{2}\cdot\frac{1}{1-z}
=\frac{m}{2}\left(\frac{n+1}{1-z}-\frac{n-1}{1+z}\right)$$

Thus, we are led to the conclusion in the statement. 
\end{proof}

As a comment now on this, you might wonder what that parameters $m>0,n\geq 1$ stand for. In answer, the laws $\mu_{mn}$ naturally appear in the random matrix context, with $m,n$ having a certain significance there, and we will discuss this later in this book. In the meantime, let us do some mathematics for $\mu_{mn}$. Following \cite{bn1}, we first have:

\index{blocks of even size}

\begin{theorem}
The measure $\mu_{mn}$ has the following properties, 
\begin{enumerate}
\item $R(z)=m(1+nz)/(1-z^2)$,

\item $(\xi G-1)(1-G^2)=mG(1+nG)$,

\item $(M-1)(1-z^2M^2)=mzM(1+nzM)$,

\item The moments are $M_k=\sum_{\pi\in NC(k)}m^{|\pi|}n^{||\pi||}$,

\item The mean is $E=m$, the variance is $V=mn$,
\end{enumerate}
with $|.|$ and $||.||$ being the number of blocks, and the number of blocks of even size.
\end{theorem}

\begin{proof}
This is something quite tricky, the idea being as follows:

\medskip

(1) To start with, the formula from the proof of Theorem 8.13 gives, as desired:
$$R=\frac{m}{2}\left(\frac{n+1}{1-z}-\frac{n-1}{1+z}\right)=\frac{m(1+nz)}{1-z^2}$$

(2) Next, observe that the above $R$-transform satisfies the following equation:
$$zR(1-z^2)=mz(1+nz)$$

In terms of $K=R+z^{-1}$, we have $zK=zR+1$, and our equation becomes:
$$(zK-1)(1-z^2)=mz(1+nz)$$

Now with $K\to\xi$, so that $z\to G$, according to $G=K^{-1}$, we obtain, as desired:
$$(\xi G-1)(1-G^2)=mG(1+nG)$$

(3) Next, with $\xi=z^{-1}$, so that $G=zM$, we obtain, as desired:
$$(M-1)(1-z^2M^2)=mzM(1+nzM)$$

(4) In order to prove this, let us compute the generating function $M=\sum_kM_kz^k$ of the numbers $M_k=\sum_{\pi\in NC(k)}m^{|\pi|}n^{||\pi||}$ in the statement. If we denote by $N(k,b,e)$ the number of partitions in $NC(k)$ having $b$ blocks and $e$ even blocks, we have:
$$M
=1+\sum_{k=1}^\infty\sum_{\pi\in NC(k)}z^km^{|\pi|}n^{||\pi||}
=1+\sum_{k=1}^\infty\sum_{b=0}^\infty\sum_{e=0}^\infty z^km^bn^eN(k,b,e)$$

Let us try now to find a recurrence formula for the numbers $N(k,b,e)$. If we look at the block containing $1$, this block must have $r\geq 0$ other legs, and we get:
\begin{eqnarray*}
N(k,b,e) 
&=&\sum_{r\in 2\mathbb N}\sum_{k=\Sigma k_i+r+1}\sum_{b=\Sigma b_i+1}\sum_{e=\Sigma e_i}N(k_1,b_1,e_1)\ldots N(k_{r+1},b_{r+1},e_{r+1})\\
&+&\sum_{r\in 2\mathbb N+1}\sum_{k=\Sigma k_i+r+1}\sum_{b=\Sigma b_i+1}\sum_{e=\Sigma e_i+1}N(k_1,b_1,e_1)\ldots N(k_{r+1},b_{r+1},e_{r+1})
\end{eqnarray*}

Here $k_1,\ldots,k_{r+1}$ are the number of points between the legs of the block containing 1, so that we have $k=(k_1+\ldots+k_{r+1})+r+1$, and the whole sum is split over two cases, $r$ even or odd, because the parity of $r$ affects the number of even blocks of our partition. Now by multiplying everything by a $z^km^bn^e$ factor, and by carefully distributing the various powers of $z,m,b$ on the right, we obtain the following formula:
\begin{eqnarray*}
z^km^bn^eN(k,b,e)
&=&m\sum_{r\in 2\mathbb N}z^{r+1}\sum_{k=\Sigma k_i+r+1}\sum_{b=\Sigma b_i+1}\sum_{e=\Sigma e_i}\prod_{i=1}^{r+1}z^{k_i}m^{b_i}n^{e_i}N(k_i,b_i,e_i)\\
&+&mn\sum_{r\in 2\mathbb N+1}z^{r+1}\sum_{k=\Sigma k_i+r+1}\sum_{b=\Sigma b_i+1}\sum_{e=\Sigma e_i+1}\prod_{i=1}^{r+1}z^{k_i}m^{b_i}n^{e_i}N(k_i,b_i,e_i)
\end{eqnarray*}

Let us sum now all these equalities, over all $k\geq 1$ and over all $b,e\geq 0$. According to the definition of $M$, at left we obtain $M-1$. As for the two sums appearing on the right, that is, at right of the two $z^{r+1}$ factors, when summing them over all $k\geq 1$ and over all $b,e\geq 0$, we obtain in both cases $M^{r+1}$. Thus, we have the following formula:
\begin{eqnarray*}
M-1
&=&m\sum_{r\in 2\mathbb N}(zM)^{r+1}+mn\sum_{r\in 2\mathbb N+1}(zM)^{r+1}\\
&=&m\,\frac{zM}{1-z^2M^2}+mn\,\frac{z^2M^2}{1-z^2M^2}\\
&=&mzM\,\frac{1+nzM}{1-z^2M^2}
\end{eqnarray*}

But this is exactly the equation found in (3), so our moment formula is proved.

\medskip

(5) According to our moment formula in (4) we have $M_1=m$ and $M_2=m^2+mn$, so the expectation is $E=m$ and the variance is $V=mn$, as stated.
\end{proof}

Still following \cite{bn1}, we have the following key supplementary result, adding to what we have in Theorem 8.14, regarding the densities of the measures $\mu_{mn}$:

\index{Bercovici-Voiculescu}
\index{atoms of convolution}

\begin{theorem}
The density of $\mu_{mn}$ has the following properties:
\begin{enumerate}
\item It has at most one atom, located at $0$, of mass $\max(1-mn,0)$.

\item The support is positive precisely when $n\leq m/4+1/m$ and $m \geq 2$.
\end{enumerate}
\end{theorem}

\begin{proof}
We have two assertions here, the idea being as follows:

\medskip

(1) The first assertion follows from the general characterization of atoms of a free additive convolution, given by  Bercovici and Voiculescu in \cite{bvo}, namely $x$ is an atom for the free additive convolution of two measures $\nu_1$ and $\nu_2$ precisely when:
$$x=x_1+x_2\quad,\quad\nu_1(\{x_1\})+\nu_2(\{x_2\})>1$$

In addition, if this is the case, we have the following formula, also from \cite{bvo}:
$$[\nu_1\boxplus\nu_2](\{x\})=\nu_1(\{x_1\})+\nu_2(\{x_2\})-1$$

In our situation, we can apply this with $\nu_1$ being a free Poisson law of parameter $m(n+1)/2$, and $\nu_2$ being the image of a free Poisson law of parameter $m(n-1)/2$ through the negation map, and we are led to the formula in the statement.

\medskip

(2) This is something more complicated, based on a detailed study of the degree 3 equation  for the Cauchy transform, from Theorem 8.14, namely:
$$(\xi G-1)(1-G^2)=mG(1+nG)$$

We use the standard fact that the support of the absolutely continuous part of a probability measure is a union of intervals defined by the points where the analyticity of the Cauchy transform $G$ breaks. In our case, for the measure $\mu_{mn}$, these points are the roots of the discriminant of the above equation defining $G$, which is as follows:
\begin{eqnarray*}
\Delta(\xi)
&=&4\xi^4-12m\xi^3+(n^2m^2+12m^2-20nm-8)\xi^2\\
&&\ +(-2n^2m^3-4m^3+22nm^2-20m)\xi\\
&&\ +n^2m^4-4n^3m^3-2nm^3+12n^2m^2+m^2-12nm+4
\end{eqnarray*}

Now observe that $\Delta(\xi)=0$ is a degree 4 equation in $\xi$ that has either 2 or 4 real solutions, with 4 complex solutions being ruled out, since the support must be non-empty. In order to decide between these cases, we compute the discriminant of $\Delta$:
$$\Delta_2(m,n)=-256m^2(n-1)(n+1)\left(m^3n^3+15m^2n^2-27m^2+48mn-64\right)^3$$

Indeed, recall that the discriminant is negative when the quartic has two real and two complex roots. In our case, due to our assumption $n\geq1$, and with the case $n=1$ being trivial, with $\mu_{m1}=\pi_m$ in this case, the sign of $\Delta_2$ is opposite to the sign of: 
$$P(m,n) = m^3 n^3+15 m^2 n^2+48 m n-27 m^2-64$$

But this leads, via some calculus, partly computer aided, to the conclusion in the statement. For details and pictures here, taking about 3 pages, we refer to \cite{bn1}.
\end{proof}

\section*{8c. Bessel laws}

Following the classical Poisson law material from chapter 1, and the paper \cite{bb+}, we would like now to use as ``input'' for the CPLT the measures of type $\nu=t\varepsilon_s$, with $t>0$, and with $\varepsilon_s$ being the uniform measure on the $s$-th roots of unity.

\bigskip

As a first observation, we would need a complex extension of the CPLT theory developed in the previous section, and this is something quite tricky, as shown by:

\begin{warning}
It is tempting to use complex measures $\nu\in\mathcal P(\mathbb C)$ in the free CPLT, but some of the formulae will hold, with suitable definitions, while some other not.
\end{warning}

As an illustration for this, let us try for instance to talk about the $R$-transform, in the complex setting. A good question here is that of defining the $R$-transform of a circular variable, $c=(a+ib)/\sqrt{2}$, and intuition tells us to set $R_c=R_{a/\sqrt{2}}+R_{ib/\sqrt{2}}$. But, will this work. By using the dilation formula $R_{qf}(z)=qR_f(qz)$, we obtain in this way:
\begin{eqnarray*}
R_c(z)
&=&\frac{1}{\sqrt{2}}\,R_a\left(\frac{z}{\sqrt{2}}\right)
+\frac{i}{\sqrt{2}}\,R_b\left(\frac{iz}{\sqrt{2}}\right)\\
&=&\frac{1}{\sqrt{2}}\cdot\frac{tz}{\sqrt{2}}
+\frac{i}{\sqrt{2}}\cdot\frac{tiz}{\sqrt{2}}\\
&=&0
\end{eqnarray*}

Which does not look very good. Of course, by thinking some more it is possible to ultimately fix this. But you get the point I guess, Warning 8.16 is something real.

\bigskip

In view of this, and with the idea in mind of avoiding difficulties, even if this will eventually lead us to less theorems, the most convenient is to proceed as follows:

\index{compound free Poisson law}

\begin{definition}
Given a complex discrete measure $\nu\in\mathcal P(\mathbb C)$, not necessarily of mass $1$, written as $\nu=\sum_{r=1}^sc_r\delta_{z_r}$ with $c_r>0$ and $z_r\in\mathbb C$, we set
$$\pi_\nu=law\left(\sum_{r=1}^sz_r\alpha_r\right)$$
with $\alpha_r\sim\pi_{c_r}$ being free, and call this compound free Poisson law.
\end{definition}

Observe that this fits both with what we did in chapter 1, in the classical complex case, and in the previous section, in the free real case. Thus, we have some compound free Poisson law theory up and working, save for the CPLT itself, and for the related $R$-transform formulae too, that we have sacrificed, in the name of mathematical rigor.

\bigskip

Getting now to what we wanted to do, free Bessel laws, we can formulate:

\index{free Bessel law}
\index{Bessel law}
\index{real Bessel law}
\index{complex Bessel law}

\begin{definition}
The free Bessel laws, depending on parameters $s\in\mathbb N\cup\{\infty\}$ and $t>0$, are the following compound free Poisson laws,
$$\pi^s_t=\pi_{t\varepsilon_s}$$
with $\varepsilon_s$ being the uniform measure on the $s$-th roots of unity. In particular:
\begin{enumerate}
\item At $s=1$ we recover the free Poisson law $\pi_t$.

\item At $s=2$ we have the real free Bessel law $\pi_t^2$.

\item At $s=\infty$ we have the purely complex free Bessel law $\Pi_t$.
\end{enumerate}
\end{definition}

To be more precise, what we have here is a rigorous definition for the free Bessel laws $\pi^s_t$, based on Definition 8.17, which in practice means the following:
$$\pi^s_t=law\left(\sum_{r=1}^sw^r\alpha_r\right)\quad,\quad w=e^{2\pi i/s}\quad ,\quad  \alpha_r\sim\pi_{t/s}\ ,\ {\rm free}$$

Next, (1) is a remark, coming from definitions. Then (2) is a definition, and more on the laws $\pi^2_t$ in a moment. And then (3) is again a definition, actually coming a bit in advance, and this because Definition 8.17 as stated does not really allow us to talk about $s=\infty$, where $\varepsilon_\infty$ is the uniform measure on the unit circle $\mathbb T$. More on this later.

\bigskip

As a first result about the free Bessel laws, which is straightforward, we have:

\index{convolution semigroup}

\begin{theorem}
The free Bessel laws have the property
$$\pi^s_t\boxplus\pi^s_{t'}=\pi^s_{t+t'}$$
so they form a $1$-parameter semigroup with respect to free convolution.
\end{theorem}

\begin{proof}
This is clear indeed from our definition of the free Bessel laws.
\end{proof}

Next, we can say many things at $s=2$, where our measure is real, as follows:

\begin{theorem}
The free Bessel law $\pi_t^2$ of parameter $t>0$, appearing as
$$\pi_t^2=\pi_{t\varepsilon}\quad,\quad\varepsilon=\frac{\delta_{-1}+\delta_1}{2}$$
equals $\mu_{mn}$ in the $m=t/n\to0$ limit, and has the following properties:
\begin{enumerate}
\item The $R$-transform is $R(z)=tz/(1-z^2)$.

\item $G=G(\xi)$ satisfies $(\xi G-1)(1-G^2)=tG^2$.

\item $M=M(z)$ satisfies $(M-1)(1-z^2M^2)=tz^2M^2$.

\item The moments are $M_k=\sum_{\pi\in NC_{even}(k)}t^{|\pi|}$.

\item The mean is $E=0$, the variance is $V=t$.

\item The density is analytic, with an atom at $0$ when $t<1$.
\end{enumerate}
\end{theorem}

\begin{proof}
In what regards the first assertion, this comes from Theorem 8.13, because with $m=t/n\to 0$, the measure $\nu$ appearing there, producing $\mu_{mn}=\pi_\nu$, is given by: 
$$\nu
=\frac{m(n-1)}{2}\,\delta_{-1}+\frac{m(n+1)}{2}\,\delta_1
=\frac{t-m}{2}\,\delta_{-1}+\frac{t+m}{2}\,\delta_1
\simeq t\cdot\frac{\delta_{-1}+\delta_1}{2}$$

As for the rest, this follows from Theorems 8.14 and 8.15, with $m=t/n\to0$:

\medskip

(1) $R(z)=m(1+nz)/(1-z^2)$ gives $R(z)=tz/(1-z^2)$ with $m=t/n\to0$.

\medskip

(2) $(\xi G-1)(1-G^2)=mG(1+nG)$ gives $(\xi G-1)(1-G^2)=tG^2$.

\medskip

(3) $(M-1)(1-z^2M^2)=mzM(1+nzM)$ gives $(M-1)(1-z^2M^2)=tz^2M^2$.

\medskip

(4) $M_k=\sum_{\pi\in NC(k)}m^{|\pi|}n^{||\pi||}$ gives $M_k=\sum_{\pi\in NC_{even}(k)}t^{|\pi|}$.

\medskip

(5) $E=m$, $V=mn$ gives $E=0$, $V=t$.

\medskip

(6) This follows from Theorem 8.15, with $m=t/n\to0$.
\end{proof}

Getting now to the general case, as a main result here, following \cite{bb+}, we have:

\begin{theorem}
The moments of the free Bessel law $\pi^s_t$ are the numbers
$$M_k=\sum_{\pi\in NC^s(k)}t^{|\pi|}$$
where $NC^s$ are the noncrossing partitions satisfying $\#\circ=\#\bullet(s)$ in each block.
\end{theorem}

\begin{proof}
According to what was said in Definition 8.18 and afterwards, the free Bessel law $\pi^s_t$ with $s\in\mathbb N$ and $t>0$ appears by definition as follows:
$$\pi^s_t=law\left(\sum_{r=1}^sw^r\alpha_r\right)\quad,\quad w=e^{2\pi i/s}\quad ,\quad  \alpha_r\sim\pi_{t/s}\ ,\ {\rm free}$$

We must compute the moments of the above variable $\alpha=\sum_rw^r\alpha_r$, and this can be done in several possible ways, which are all instructive, as follows:

\medskip

(1) Combinatorics. Given a colored integer exponent $k=(k_1\ldots k_n)$, we have:
\begin{eqnarray*}
M_k
&=&tr(\alpha^{k_1}\ldots\alpha^{k_n})\\
&=&tr\left(\sum_{r_1=1}^sw^{k_1r_1}\alpha_{r_1}\ldots\sum_{r_n=1}^sw^{k_nr_n}\alpha_{r_n}\right)\\
&=&\sum_{r_1\ldots r_n}w^{<k,r>}tr(\alpha_{r_1}\ldots\alpha_{r_n})
\end{eqnarray*}

But the trace on the right can be computed by using the basic freeness rule, from chapter 6, and we are led to the following formula for our moment, as desired:
\begin{eqnarray*}
M_k
&=&\sum_{r_1\ldots r_n}w^{<k,r>}\sum_{\pi\in NC(n)}\left(\frac{t}{s}\right)^{|\pi|}\delta_{\pi\leq\ker r}\\
&=&\sum_{\pi\in NC(n)}\left(\frac{t}{s}\right)^{|\pi|}\sum_{\pi\leq\ker r}w^{<k,r>}\\
&=&\sum_{\pi\in NC(n)}\left(\frac{t}{s}\right)^{|\pi|}\cdot\,s^{|\pi|}\,\delta_{\pi\in NC^s(k)}\\
&=&\sum_{\pi\in NC^s(k)}t^{|\pi|}
\end{eqnarray*}

(2) Gaussian matrices. At $t=1$, consider a family of $sN\times N$ Gaussian matrices $Z_1,\ldots,Z_s$, with independent $G_1$ entries. By rotating, we obtain the following family of matrices $Y_1,\ldots,Y_s$, which are also Gaussian, with independent $G_1$ entries:
$$Y_r=\frac{1}{\sqrt{s}}\sum_{p=1}^sw^{rp}Z_p$$

Thus, the following matrix provides an asymptotic model for $\alpha=\sum_rw^r\alpha_r$:
\begin{eqnarray*}
T
&=&\frac{1}{sN}\sum_{r=1}^sw^rY_rY_r^*\\
&=&\frac{1}{sN}\cdot\frac{1}{s}\sum_{r=1}^s\sum_{p=1}^s\sum_{q=1}^sw^{r+rp-rq}Z_pZ_q^*\\
&=&\frac{1}{sN}\sum_{p=1}^s\sum_{q=1}^s\left(\frac{1}{s}\sum_{r=1}^s\left(w^{1+p-q}\right)^r\right)Z_pZ_q^*\\
&=&\frac{1}{sN}(Z_1Z_2^*+Z_2Z_3^*+\ldots+Z_{s-1}Z_s^*+Z_sZ_1^*)
\end{eqnarray*}

Now observe that this latter matrix can be written as follows:
\begin{eqnarray*}
T
&=&\frac{1}{sN}\begin{pmatrix}Z_1&Z_2&\ldots&Z_{s-1}&Z_s\end{pmatrix}
\begin{pmatrix}Z_2^*\\ Z_3^*\\\vdots\\ Z_s^*\\ Z_1^*\end{pmatrix}\\
&=&\frac{1}{sN}\begin{pmatrix}Z_1&Z_2&\ldots&Z_{s-1}&Z_s\end{pmatrix}
\begin{pmatrix}
0&1_N&0&\ldots&0\\
0&0&1_N&\ldots&0\\
&&&\ddots&&\\
0&0&0&\ldots&1_N\\
1_N&0&0&\ldots&0
\end{pmatrix}
\begin{pmatrix}Z_1^*\\ Z_2^*\\\vdots\\ Z_{s-1}^*\\ Z_s^*\end{pmatrix}\\
&=&\frac{1}{sN}\cdot ZOZ^*
\end{eqnarray*}

To be more precise, in the last formula $Z=(Z_1\ \ldots\  Z_s)$ is the $sN\times sN$ Gaussian matrix obtained by concatenating $Z_1,\ldots,Z_s$, and $O$ is the matrix in the middle. But this latter matrix $O$ is of the form $O=UDU^*$, with $U$ unitary, and with:
$$D=\begin{pmatrix}
1_N&0&&0\\
0&w1_N&&0\\
&&\ddots&\\
0&0&&w^{s-1}1_N
\end{pmatrix}$$

Thus $T=ZUDU^*Z^*/(sN)$, and since $ZU$ is a Gaussian matrix, we conclude that:
$$law(T)=law\left(\frac{ZDZ^*}{sN}\right)$$

Getting now to moments, given a colored integer $k=(k_1\ldots k_n)$, we have:
\begin{eqnarray*}
M_k(T)
&=&\frac{1}{(sN)^n}\,M_k(ZDZ^*)\\
&=&\frac{1}{(sN)^n}\,(E\circ tr)[(ZDZ^*)^{k_1}(ZDZ^*)^{k_2}\ldots(ZDZ^*)^{k_n}]\\
&=&\frac{1}{(sN)^n}\,(E\circ tr)[ZD^{k_1}Z^*ZD^{k_2}Z^*\ldots ZD^{k_n}Z^*]\\
&=&\frac{1}{(sN)^n}\,(E\circ tr)[D^{k_1}Z^*ZD^{k_2}Z^*Z\ldots D^{k_n}Z^*Z]
\end{eqnarray*}

By using now the formula from chapter 7 for such quantities we obtain, with $d\sim\varepsilon_s$, $\pi$ being the partition obtained by shrinking $\nu$, and $b\in\pi$ being the blocks of $\pi$:
\begin{eqnarray*}
M_k(T)
&\simeq&\frac{1}{(sN)^n}\,(sN)^n\sum_{\nu\in\mathcal{NC}_2((\bullet\circ)^n)}tr_{\gamma\nu}(d^{k_1},1,d^{k_2},1\ldots,d^{k_n},1)\\
&=&\sum_{\nu\in\mathcal{NC}_2((\bullet\circ)^n)}tr_{\gamma\nu}(d^{k_1},1,d^{k_2},1\ldots,d^{k_n},1)\\
&=&\sum_{\pi\in NC(k)}\prod_{b\in\pi}tr(d^{\#b})\\
&=&\sum_{\pi\in NC(k)}\prod_{b\in\pi}\delta_{s|\#b}\\
&=&|NC^s(k)|
\end{eqnarray*}

As for the case $t>0$, this is more technical, and we will leave it as an exercise.

\medskip

(3) Wishart matrices. This is a variation of the above proof, taking advantage of some previously known formulae for the Wishart matrices. Let $W$ be a complex Wishart matrix of parameters $(sN,N)$, and consider as well the following diagonal matrix:
$$D=\begin{pmatrix}
1_N&0&&0\\
0&w1_N&&0\\
&&\ddots&\\
0&0&&w^{s-1}1_N
\end{pmatrix}$$

By asymptotic freeness we have the following formula, in the $N\to\infty$ limit:
$$DW\sim\pi^s_1$$ 

Now in order to compute the moments, we can use the following general formula of Graczyk, Letac and Massam \cite{glm}, whose proof is via standard combinatorics:
$$E(Tr(DW)^k)=\sum_{\sigma\in S_k}M^{|\sigma^{-1}\gamma|-k}\,Tr_\sigma(D)$$

Here $W$ is by definition a complex Wishart matrix of parameters $(M,N)$, and $D$ is a deterministic $M\times M$ matrix. As for the right term, here $\gamma$ is the cycle $(1,\ldots,k)$, then $|\sigma|$ is the number of disjoint cycles of $\sigma$, and if we denote by $C(\sigma)$ the set of such cycles, and we let $c\to l(c)$ be the cycle length, then the function on the right is given by:
$$Tr_\sigma(D)=\prod_{c\in C(\sigma)}Tr(D^{l(c)})$$

In our situation, since $D$ is uniformly formed by $s$-roots of unity, we have:
$$Tr(D^p)=
\begin{cases}
sN\mbox{ if }s|p\\
0\ \ \,\mbox{ if }s\!\!\not|p
\end{cases}$$

Thus if we denote by $S_k^s$ the set of permutations $\sigma\in S_k$ having the property that all cycles of $\sigma$ have as length a multiple of $s$, the above formula gives:
$$E(Tr(DW)^k)=\sum_{\sigma\in S_k^s}(sN)^{|\sigma^{-1}\gamma|-k}\,(sN)^{|\sigma|}$$

In terms of the normalized trace $tr$, we obtain the following formula:
$$E(tr(DW)^k)=\sum_{\sigma\in S_k^s}(sN)^{|\sigma^{-1}\gamma|+|\sigma|-k-1}$$

Now when taking the $N\to\infty$ limit, only the permutations $\sigma\in S_k^s$ which are geodesic for $id-\gamma$ will contribute, each with a 1 value, and we therefore obtain:
$$\lim_{N\to\infty}E(tr(DW)^k)
=\#\left\{\sigma\in S_k^s\big|\sigma\in [id,\gamma]\right\}
=|NC^s(k)|$$

As for the case $t>0$, this is more technical, and we will leave it as an exercise.
\end{proof}

As a consequence of Theorem 8.21, we can now talk about $\Pi_t=\pi^\infty_t$, by saying that this is by definition the law having the following moments, with $\mathcal{NC}_{even}=NC^\infty$:
$$M_k=\sum_{\pi\in\mathcal{NC}_{even}(k)}t^{|\pi|}$$

At the combinatorial level, all this is quite interesting, and we have:

\index{standard cube}
\index{limiting measures}

\begin{theorem}
For the main limiting measures we have $M_k=\sum_{\pi\in D(k)}t^{|\pi|}$,
$$\xymatrix@R=19pt@C=22pt{
&\Pi_t\ar@{-}[rr]\ar@{-}[dd]&&\Gamma_t\ar@{-}[dd]\\
\pi^2_t\ar@{-}[rr]\ar@{-}[dd]\ar@{-}[ur]&&\gamma_t\ar@{-}[dd]\ar@{-}[ur]\\
&P_t\ar@{-}[rr]\ar@{-}[uu]&&G_t\ar@{-}[uu]\\
p^2_t\ar@{-}[uu]\ar@{-}[ur]\ar@{-}[rr]&&g_t\ar@{-}[uu]\ar@{-}[ur]
}\qquad\xymatrix@R=20pt@C=22pt{\\ \\ :\\}\qquad 
\xymatrix@R=20pt@C=0pt{
&\mathcal{NC}_{even}\ar[dl]\ar[dd]&&\ \ \ \mathcal{NC}_2\ \ \ \ar[ll]\ar[dd]\ar[dl]\\
NC_{even}\ar[dd]&&NC_2\ar[ll]\ar[dd]\\
&\mathcal P_{even}\ar[dl]&&\mathcal P_2\ar[ll]\ar[dl]\\
P_{even}&&P_2\ar[ll]
}$$
with the categories of partitions $D\subset P$ involved being those on the right. 
\end{theorem}

\begin{proof}
This follows indeed by putting together the various moment results that we have, with in the case of the Bessel laws, by restricting the attention to the cases $s=2,\infty$, which can be thought of as being ``fully real'' and ``purely complex''.
\end{proof}

\section*{8d. Real versions}

Let us discuss now, still following \cite{bb+}, some further aspects of the free Bessel laws. We have the following key definition, coming from a chi type idea:

\index{modified free Bessel law}

\begin{definition}
The modified free Bessel law $\tilde{\pi}^s_t$ is the real measure given by
$$M_k(\tilde{\pi}^s_t)=M_{sk}(\pi^s_t)$$
which must hold for any $k\in\mathbb N$.
\end{definition}

This is something quite tricky. To start with, at $s=2$, the free Bessel law $\pi^2_t$ being a real probability measure, we have available a squared version of it, given by:
$$f\sim\pi^2_t\ \implies\ f^2\sim\tilde{\pi}^2_t$$

Equivalently, in terms of moments, this modified free Bessel law $\tilde{\pi}^2_t$ is given by:
$$M_k(\tilde{\pi}^2_t)=M_{2k}(\pi^2_t)$$

In view of this, the problem appears, what is the correct generalization of the above law $\tilde{\pi}^2_t$, in the general case $s\in\mathbb N\cup\{\infty\}$? And here, at $s=1$ to start with, since the Marchenko-Pastur law $\pi_t=\pi_t^1$ will obviously not simplify when squaring, or doing other types of similar manipulations on it, we must leave this law untouched:
$$\tilde{\pi}^1_t=\pi^1_t$$

But then, some further thinking along these lines, and some computations too, suggest a manipulation as follows, which is in tune with what we have above at $s=1,2$:
$$f\sim\pi^s_t\ \implies\ f^s\sim\tilde{\pi}^s_t$$

However, passed $s=\infty$ that we obviously have to leave aside, this won't work well at $s\geq3$, because we will not obtain in this way a real measure. So, the trick is the one in Definition 8.23, namely restricting the attention to the uncolored moments of $f^s$, which amounts in defining $\tilde{\pi}^s_t$ by the following formula, which must hold for any $k\in\mathbb N$:
$$M_k(\tilde{\pi}^s_t)=M_{sk}(\pi^s_t)$$

Of course, it still remains to be proved that $\tilde{\pi}^s_t$ exists indeed, as a measure. But this will follow from our various results below, providing alternative definitions for $\tilde{\pi}^s_t$.

\bigskip

With this discussed, following \cite{bb+}, let us get now into the mathematics of $\tilde{\pi}^s_t$. We first have the following result, regarding the moment generating function:

\begin{theorem}
The moment generating function of $\tilde{\pi}^s_t$ satisfies
$$M-1=zM^s(M+t-1)$$
and in particular at $t=1$ we have the formula $M-1=zM^{s+1}$.
\end{theorem}

\begin{proof}
This can be done via some combinatorics, in two steps, as follows:

\medskip

(1) At $t=1$, with the notation $C_k=|NC^s(sk)|$, the moment generating function is:
$$M(z)=\sum_kC_{k}z^k$$

For a given partition $p\in NC^s(sk+s)$ we can consider the last $s$ legs of the first block, and make cuts at right of them. This gives a decomposition of $p$ into $s+1$ partitions in $NC^s$, and we obtain the following recurrence formula for the numbers $C_k$:
$$C_{k+1}=\sum_{\Sigma k_i=k}C_{k_0}\ldots C_{k_s}$$

By multiplying now by $z^{k+1}$, and then summing over $k$, we obtain that the generating series of these numbers $C_k$ satisfies the following equation, as desired:
$$M-1=zM^{s+1}$$

(2) At $t>0$, with the above notations from (1), let $F_{kb}$ be the number of partitions in $NC^s(sk)$ having $b$ blocks, and set $F_{kb}=0$ for other integer values of $k,b$. All sums will be over integer indices $\geq 0$. The moment generating function is then:
$$M(z)=\sum_{kb}F_{kb}z^kt^b$$

In terms of $F_{kb}$, the recurrence formula for the numbers $C_k$ in (1) becomes:
$$\sum_bF_{k+1,b}=\sum_{\Sigma k_i=k}\sum_{b_i}F_{k_0b_0}\ldots F_{k_sb_s}$$

In this formula, each term contributes to $F_{k+1,b}$ with $b=\Sigma b_i$, except for those of the form $F_{00}F_{k_1b_1}\ldots F_{k_sb_s}$, which contribute to $F_{k+1,b+1}$. We get:
\begin{eqnarray*}
F_{k+1,b}&=&\sum_{\Sigma k_i=k}\sum_{\Sigma b_i=b}F_{k_0b_0}\ldots F_{k_sb_s}\cr
&+&\sum_{\Sigma k_i=k}\sum_{\Sigma b_i=b-1}F_{k_1b_1}\ldots F_{k_sb_s}\cr
&-&\sum_{\Sigma k_i=k}\sum_{\Sigma b_i=b}F_{k_1b_1}\ldots F_{k_sb_s}
\end{eqnarray*}

This gives the following formula for the polynomials $P_k=\sum_bF_{kb}t^b$:
$$P_{k+1}=\sum_{\Sigma k_i=k}P_{k_0}\ldots  P_{k_s}+(t-1)\sum_{\Sigma k_i=k}P_{k_1}\ldots P_{k_s}$$

In terms of $M=\sum_kP_kz^k$, we obtain the following equation:
$$M-1=zM^{s+1}+(t-1)zM^s$$

Thus, we are led to the conclusion in the statement.
\end{proof}

Regarding the moments themselves, we have here the following result:

\index{Fuss-Catalan numbers}
\index{Fuss-Narayana numbers}

\begin{theorem}
The moments of $\tilde{\pi}^s_1$ and $\tilde{\pi}^s_t$ are given by
$$M_k(\tilde{\pi}^s_1)=\frac{1}{sk+1}\binom{sk+k}{k}\quad,\quad 
M_k(\tilde{\pi}^s_t)=\sum_{b=1}^k\frac{1}{b}\binom{k-1}{b-1}\binom{sk}{b-1}t^b$$
which are the Fuss-Catalan numbers, and the Fuss-Narayana numbers.
\end{theorem}

\begin{proof}
This is standard from what we have, the idea being as follows:

\medskip

(1) At $t=1$ we know that we have $M_k=|NC_s(sk)|$, and the formula in the statement, which is something very classical, follows then by counting such partitions.

\medskip

(2) At $t>0$, we know from the proof of Theorem 8.24 that we have the following formula, where $F_{kb}$ is the number of partitions in $NC_s(sk)$ having $b$ blocks:
$$M_k=\sum_bF_{kb}t^b$$

With this observation in hand, the formula in the statement follows by counting such partitions, with this count being again something well-known.

\medskip

(3) Finally, observe that at $s=1$ we get the moments of $\pi_1,\pi_t$. At $s=2$ we obtain some new formulae for the moments $\pi^2_1,\pi^2_t$, adding to what we have in Theorem 8.20.
\end{proof}

Still following \cite{bb+}, we have the following key result, regarding the laws $\tilde{\pi}^s_t$:

\begin{theorem}
$\tilde{\pi}^s_t$ is subject to the following formulae, involving $\pi=\pi_1$ only,
$$\tilde{\pi}^s_t=\pi^{\boxtimes s-1}\boxtimes\pi^{\boxplus t}\quad,\quad 
\tilde{\pi}^s_t=((1-t)\delta_0+t\delta_1)\boxtimes\pi^{\boxtimes s}$$
and based on these, we can talk about $\tilde{\pi}^s_t$ with $s\geq1,t>0$, or $s>0,t\in(0,1]$.
\end{theorem}

\begin{proof}
This is something quite tricky, the idea being as follows:

\medskip

(1) We first have the following computation, for the $S$-transform of $\pi_t=\pi^{\boxplus t}$:
\begin{eqnarray*}
\xi G^2+1=(\xi+1-t)G
&\implies&z(\psi+1)^2+1=(1+z-zt)(\psi+1)\\
&\implies&\chi(z+1)^2+1=(1+\chi-\chi t)(z+1)\\
&\implies&\chi(z+1)(t+z)=z\\
&\implies&S=1/(t+z)
\end{eqnarray*}

We conclude that the $S$-transform of the measure $\pi^{\boxtimes s-1}\boxtimes\pi^{\boxplus t}$ is given by:
$$S(z)=\frac{1}{(1+z)^{s-1}}\cdot\frac{1}{t+z}$$

(2) Next, we know from chapter 7 that the $S$-transform of $(1-t)\delta_0+t\delta_1$ is:
$$S=\frac{1+z}{t+z}$$

We conclude that the $S$-transform of the measure $((1-t)\delta_0+t\delta_1)\boxtimes\pi^{\boxtimes s}$ is given by the following formula, which coincides with the one found in (1):
$$S(z)=\frac{1}{(1+z)^s}\cdot\frac{1+z}{t+z}$$

(3) Summarizing, done with two measures. In order to deal now with $\tilde{\pi}^s_t$, we can convert the equation of $M$ found in Theorem 8.24 into a formula for $S$, as follows:
\begin{eqnarray*}
M-1=zM^s(M+t-1)
&\implies&\psi=z(1+\psi)^s(t+\psi)\\
&\implies&z=\chi(1+z)^s(t+z)\\
&\implies&(1+z^{-1})z=S(1+z)^s(t+z)\\
&\implies&S=(1+z)^{1-s}/(t+z)
\end{eqnarray*}

Thus the $S$-transform of $\tilde{\pi}^s_t$ coincides with the $S$-transforms computed in (1,2), and we have proved in this way both the formulae in the statement.

\medskip

(4) As a first consequence, and answering a remark made after Definition 8.23, our formulae show that $\tilde{\pi}^s_t$ exists indeed as a probability measure, at any $s\in\mathbb N$ and $t>0$.

\medskip

(5) Next, a careful examination of what we did shows that we can define $\tilde{\pi}^s_t$ at any $s\geq1,t>0$, by setting $\tilde{\pi}^s_t=\pi^{\boxtimes s-1}\boxtimes\pi^{\boxplus t}$, and at any $s>0,t\in(0,1]$ too, by setting $\tilde{\pi}^s_t=((1-t)\delta_0+t\delta_1)\boxtimes\pi^{\boxtimes s}$. Thus, we are led to the conclusions in the statement.
\end{proof}

Still following \cite{bb+}, we can somehow go back from $\tilde{\pi}^s_t$ to $\pi^s_t$, as follows:

\begin{theorem}
The formal measure $\pi^{s\prime}_t$ having moment generating function $M(z)=\tilde{M}(z^s)$, where $\tilde{M}$ is the moment generating function of $\tilde{\pi}^s_t$, has $R(z)=tz^{s-1}/(1-z^s)$.
\end{theorem} 

\begin{proof}
We have indeed the following computation, based on Theorem 8.24:
\begin{eqnarray*}
\tilde{M}=1+z\tilde{M}^s(\tilde{M}+t-1)
&\implies&M=1+(zM)^s(M+t-1)\\
&\implies&\xi G=1+G^s(\xi G+t-1)\\
&\implies&Kz=1+z^s(Kz+t-1)\\
&\implies&Rz+1=1+z^s(Rz+t)\\
&\implies&R(z)=tz^{s-1}/(1-z^s)
\end{eqnarray*}

And for more on all this, $\pi^s_t$, $\tilde{\pi}^s_t$ and $\pi^{s\prime}_t$, we refer to \cite{bb+}. 
\end{proof}

Finally, still following \cite{bb+}, we have the following matrix modeling result:

\index{Wishart matrix}

\begin{theorem}
Let $Z_1,\ldots,Z_s$ be a family of $N\times N$ independent matrices formed by independent centered Gaussian variables, of variance $1/N$. Then with 
$$T=Z_1\ldots Z_s$$
the moments of $TT^*$ converge to the corresponding moments of $\tilde{\pi}^s_1$, as $N\to\infty$.
\end{theorem} 

\begin{proof}
We prove this by recurrence. At $s=1$ it is well-known that $TT^*$ is a model for $\tilde{\pi}^1_1=\pi$. So, assume that the result holds for $s-1\geq 1$. We have:
\begin{eqnarray*}
tr(TT^*)^k
&=&tr(Z_1\ldots Z_sZ_s^*\ldots Z_1^*)^k\\
&=&tr\big(Z_1(Z_2\ldots Z_sZ_s^*\ldots Z_1^*Z_1)^{k-1}Z_2\ldots Z_sZ_s^*\ldots Z_1^*\big)
\end{eqnarray*}

We can pass the first $Z_1$ matrix to the right, and we get:
\begin{eqnarray*}
tr(TT^*)^k
&=&tr\big((Z_2\ldots Z_sZ_s^*\ldots Z_1^*Z_1)^{k-1}Z_2\ldots Z_sZ_s^*\ldots Z_1^*Z_1\big)\\
&=&tr(Z_2\ldots Z_sZ_s^*\ldots Z_1^*Z_1)^k\\
&=&tr((Z_2\ldots Z_sZ_s^*\ldots Z_2^*)(Z_1^*Z_1))^k
\end{eqnarray*}

We know that $Z_1^*Z_1$ is a Wishart matrix, hence is a model for $\pi$:
$$Z_1^*Z_1\sim\pi$$

Also, we know by recurrence that $Z_2\ldots Z_sZ_s^*\ldots Z_2^*$ gives a matrix model for $\tilde{\pi}^{s-1}_1$: 
$$Z_2\ldots Z_sZ_s^*\ldots Z_2^*\sim\pi^{s-1}_1$$

Now since the matrices $Z_1^*Z_1$ and $Z_2\ldots Z_sZ_s^*\ldots Z_2^*$ are asymptotically free, their product gives a matrix model for $\tilde{\pi}^{s-1}_1\boxtimes\pi=\tilde{\pi}^s_1$, and we are done.
\end{proof}

\section*{8e. Exercises}

This was an exciting combinatorial chapter, and as exercises, we have:

\begin{exercise}
Further meditate at the phenomenology of $\gamma_1,\pi_1$.
\end{exercise}

\begin{exercise}
Solve the functional transform puzzle $\log F_m\leftrightarrow R_\mu$.
\end{exercise}

\begin{exercise}
Try finding a formula for the moments of $\mu_{mn}$.
\end{exercise}

\begin{exercise}
Identify among $\mu_{mn}$ measures that you previously knew.
\end{exercise}

\begin{exercise}
Meditate at the FCPLT, in its complex formulation.
\end{exercise}

\begin{exercise}
Meditate as well at what exactly happens at $s=\infty$.
\end{exercise}

\begin{exercise}
Learn more about the modified free Bessel laws $\tilde{\pi}^s_t$.
\end{exercise}

\begin{exercise}
Have a look as well at the modified Bessel laws $\tilde{p}^s_t$.
\end{exercise}

As bonus exercise, further experiment with the FCPLT, in general.

\part{The bijection}

\ \vskip50mm

\begin{center}
{\em Winterlude, Winterlude, my little daisy

Winterlude by the telephone wire

Winterlude, it's making me lazy

Come on, sit by the logs in the fire}
\end{center}

\chapter{The bijection}

\section*{9a. Classical cumulants}

We have seen so far the basics of free probability, with this being an exciting mix of operator theory, random matrices, classical analysis and combinatorics. So, what is next? Normally the same thing, operators, matrices, analysis and combinatorics. All this is quite addictive, and the temptation is high, to keep exploring freeness in this way.

\bigskip

However, it is perhaps a good moment to calm down, and be a bit philosophers. After all, we have enough theory, for a start, and more theory can wait. So, as a most pressing task, based on the things that we already have, what about trying to organize them, as to abstractly understand the correspondence between classical and free?

\bigskip

In answer, let us have a closer look at what we did in Part II, all the formulae found there, regarding $\gamma_t,\Gamma_t,\pi_t,\pi_\nu$, and compare them with our previous formulae from Part I, regarding their classical counterparts $g_t,G_t,p_t,p_\nu$. We are led in this way to:

\begin{speculation}
The correspondence $m_t\leftrightarrow\mu_t$ between the classical and free can be understood in several possible ways, as follows,
\begin{enumerate}
\item Basic guideline: ``independence gets replaced by freeness''.

\item Same thing, sounding better: ``the operation $*$ gets replaced by $\boxplus$''.

\item Moments, we have here the following connecting formula, with the categories $D\subset P$ and $D'\subset NC$ being related by $D'=D\cap NC$ and $D=<D',\slash\hskip-2.1mm\backslash>$:
$$M_k(m_t)=\sum_{\pi\in D(k)}t^{|\pi|}\quad \longleftrightarrow\quad M_k(\mu_t)=\sum_{\pi\in D'(k)}t^{|\pi|}$$

\item Functional transforms, we have here the following connecting formula:
$$\log F_{m_t}(-iz)=\sum_{n=1}^\infty c_n\cdot\frac{z^n}{n!}
\quad \longleftrightarrow\quad zR_{\mu_t}(z)=\sum_{n=1}^\infty c_n\cdot z^n$$
\end{enumerate}
and there are certainly more ways, by looking at orthogonal polynomials, Gram or Hankel determinants, or Weingarten functions. The only wrong way is via densities.
\end{speculation}

To be more precise here, (1) and (2) are generalities, which sound like Cold War propaganda, are these quotes from Stalin or Truman, I don't remember. Then (3) is something that we know well for $g_t\leftrightarrow\gamma_t$, $G_t\leftrightarrow\Gamma_t$, $p_t\leftrightarrow\pi_t$, $p^s_t\leftrightarrow\pi^s_t$, and more on this soon. Next, (4) comes from an exercise that I left to you some time ago, and more on this in a second. As for the assertions at the end, these are about a technological escalation of the conflict, with no holds barred, I would say Brezhnev, or most likely Reagan.

\bigskip

Getting now to the details of (4), this comes by closely examining the various $R$-transform formulae from Part II. Indeed, for the Poisson laws, this certainly holds:
$$\log F_{p_t}(-iz)=(e^z-1)t=\sum_{n=1}^\infty t\cdot\frac{z^n}{n!}$$
$$zR_{\pi_t}(z)=\frac{zt}{1-z}=\sum_{n=1}^\infty t\cdot z^n$$

Observe that by linearity this holds, more generally, for the compound Poisson laws. Next, for the Gaussian laws this holds too, due to more trivial reasons, as follows:
$$\log F_{g_t}(-iz)=\frac{z^2t}{2}=t\cdot\frac{z^2}{2!}$$
$$zR_{\gamma_t}(z)=z^2t=t\cdot z^2$$

As yet another piece of evidence, for the Dirac mass $\delta_t$, which reigns over probability theory at large, being at the same time classical and free, this works too:
$$\log F_{\delta_t}(-iz)=zt=t\cdot\frac{z}{1!}$$
$$zR_{\delta_t}(z)=zt=t\cdot z$$

Summarizing, we have our speculation up and running, and the problem is now, which of the 2 main ways there shall we follow. In answer, we will follow a middle path:

\begin{principle}
The corrrespondence betwen classical and free is best understood via quantities $k_n,\kappa_n$, called cumulants and free cumulants, appearing as follows:
$$\log F_m(-iz)=\sum_{n=1}^\infty k_n(m)\cdot\frac{z^n}{n!}
\quad,\quad 
zR_\mu(z)=\sum_{n=1}^\infty \kappa_n(\mu)\cdot z^n$$
Indeed, we can say that we have $m\leftrightarrow\mu$ when $k_n(m)=\kappa_n(\mu)$. And with some further combinatorial work, this will clarify our correspondence for the moments too. 
\end{principle} 

So, this is the situation, hope you got it, what we have initially suggests jumping into heavy analysis, by comparing $\log F_m$ and $R_\mu$, and this is indeed what Bercovici and Pata did in \cite{bpa}, when first examining the correspondence question. However, from a more modern perspective, it is actually better to replace analysis by combinatorics, by comparing instead the coefficients of $\log F_m$ and $R_\mu$, with this being equivalent to what Bercovici-Pata did, but coming with a bonus, namely moment results too.

\bigskip

Getting to work now, we first need to talk about cumulants $k_n$. And here, following Rota, let us formulate the following definition, which is something elementary:

\index{cumulant}
\index{cumulant-generating function}
\index{Taylor coefficient}
\index{Fourier transform}
\index{generating series}

\begin{definition}
Associated to any real probability measure $\mu=\mu_f$ is the following modification of the logarithm of the Fourier transform $F_\mu(z)=E(e^{izf})$,
$$K_\mu(z)=\log E(e^{zf})$$
called cumulant-generating function. The Taylor coefficients $k_n(\mu)$ of this series, given by
$$K_\mu(z)=\sum_{n=1}^\infty k_n(\mu)\,\frac{z^n}{n!}$$
are called cumulants of the measure $\mu$. We also use the notations $k_f,K_f$ for these cumulants and their generating series, where $f$ is a variable following the law $\mu$.
\end{definition}

In other words, the cumulants are more or less the coefficients of the logarithm of the Fourier transform $\log F_\mu$, up to some normalizations. To be more precise, we have $K_\mu(z)=\log F_\mu(-iz)$, so the formula relating $\log F_\mu$ to the cumulants $k_n(\mu)$ is:
$$\log F_\mu(-iz)=\sum_{n=1}^\infty k_n(\mu)\,\frac{z^n}{n!}$$

As a first observation, the sequence of cumulants $k_1,k_2,k_3,\ldots$ appears as a modification of the sequence of moments $M_1,M_2,M_3,\ldots\,$, the numerics being as follows:

\index{sequence of moments}
\index{sequence of cumulants}

\begin{proposition}
The sequence of cumulants $k_1,k_2,k_3,\ldots$ appears as a modification of the sequence of moments $M_1,M_2,M_3,\ldots\,$, and uniquely determines $\mu$. We have
$$k_1=M_1$$
$$k_2=-M_1^2+M_2$$
$$k_3=2M_1^3-3M_1M_2+M_3$$
$$k_4=-6M_1^4+12M_1^2M_2-3M_2^2-4M_1M_3+M_4$$
$$\vdots$$
in one sense, and in the other sense we have
$$M_1=k_1$$
$$M_2=k_1^2+k_2$$
$$M_3=k_1^3+3k_1k_2+k_3$$
$$M_4=k_1^4+6k_1^2k_2+3k_2^2+4k_1k_3+k_4$$
$$\vdots$$
with in both cases the correspondence being polynomial, with integer coefficients.
\end{proposition}

\begin{proof}
We know from Definition 9.3 that the cumulants are given by:
$$\log E(e^{zf})=\sum_{s=1}^\infty k_s(f)\,\frac{z^s}{s!}$$

By exponentiating, we obtain from this the following formula:
$$E(e^{zf})=\exp\left(\sum_{s=1}^\infty k_s(f)\,\frac{z^s}{s!}\right)$$

Now by looking at the terms of order $1,2,3,4$, this gives the above formulae.
\end{proof}

The interest in cumulants comes from the fact that $\log F_\mu$, and so the cumulants $k_n(\mu)$ too, linearize the convolution. To be more precise, we have the following result:

\index{linearization of convolution}

\begin{theorem}
The cumulants have the following properties:
\begin{enumerate}
\item $k_n(cf)=c^nk_n(f)$.

\item $k_1(f+d)=k_1(f)+d$, and $k_n(f+d)=k_n(f)$ for $n>1$.

\item $k_n(f+g)=k_n(f)+k_n(g)$, if $f,g$ are independent.
\end{enumerate}
\end{theorem}

\begin{proof}
Here (1) and (2) are both clear from definitions, because we have:
\begin{eqnarray*}
K_{cf+d}(z)
&=&\log E(e^{z(cf+d)})\\
&=&\log[e^{zd}\cdot E(e^{zcf})]\\
&=&zd+K_f(cz)
\end{eqnarray*}

As for (3), this follows from the fact that the Fourier transform $F_f(x)=E(e^{ixf})$ satisfies the following formula, whenever $f,g$ are independent random variables:
$$F_{f+g}(x)=F_f(x)F_g(x)$$

Indeed, by applying the logarithm, we obtain the following formula:
$$\log F_{f+g}(x)=\log F_f(x)+\log F_g(x)$$

With the change of variables $x=-iz$, we obtain the following formula:
$$K_{f+g}(z)=K_f(z)+K_g(z)$$

Thus, at the level of coefficients, we obtain $k_n(f+g)=k_n(f)+k_n(g)$, as claimed.
\end{proof}

At the level of basic examples now, we have the following result:

\begin{theorem}
The sequence of cumulants $k_1,k_2,k_3,k_4,\ldots$ is as follows:
\begin{enumerate}
\item For $\mu=\delta_t$ the cumulants are $t,0,0,0,\ldots$

\item For $\mu=g_t$ the cumulants are $0,t,0,0,\ldots$

\item For $\mu=g_t^t$ the cumulants are $t,t,0,0,\ldots$

\item For $\mu=p_t$ the cumulants are $t,t,t,t,\ldots$

\item For $\mu=p_t^2$ the cumulants are $0,t,0,t,\ldots$
\end{enumerate}
\end{theorem}

\begin{proof}
We have 5 computations to be done, the idea being as follows:

\medskip

(1) For a Dirac mass, $\mu=\delta_t$, this comes from the following computation:
$$K_\mu(z)
=\log E(e^{tz})
=\log(e^{tz})
=tz$$

Indeed, the plain coefficients of this series are the numbers $t,0,0,0,\ldots\,$, so the Taylor coefficients of this series are these same numbers $t,0,0,0,\ldots\,$, as claimed.

\medskip

(2) For a normal law, $\mu=g_t$, this comes from the following computation:
$$K_\mu(z)
=\log F_\mu(-iz)
=\log\exp\left[-t(-iz)^2/2\right]
=tz^2/2$$

Indeed, the plain coefficients of this series are the numbers $0,t/2,0,0,\ldots\,$, so the Taylor coefficients of this series are the numbers $0,t,0,0,\ldots\,$, as claimed.

\medskip

(3) For a shifted normal law, $\mu=g_t^t$, this comes from the following computation:
$$K_\mu(z)
=\log F_\mu(-iz)
=\log\exp\left[it(-iz)-t(-iz)^2/2\right]
=tz+tz^2/2$$

Indeed, the plain coefficients of this series are the numbers $t,t/2,0,0,\ldots\,$, so the Taylor coefficients of this series are the numbers $t,t,0,0,\ldots\,$, as claimed.

\medskip

(4) For a Poisson law, $\mu=p_t$, this comes from the following computation:
$$K_\mu(z)
=\log F_\mu(-iz)
=\log\exp\left[(e^{i(-iz)}-1)t\right]
=(e^z-1)t$$

Indeed, the plain coefficients of this series are the numbers $t/n!$, so the Taylor coefficients of this series are the numbers $t,t,t,t,\ldots\,$, as claimed.

\medskip

(5) For a Bessel law, $\mu=p_t^2$, this comes from the following computation:
$$K_\mu(z)
=\log F_\mu(-iz)
=\log\exp\left[\left(\frac{e^z+e^{-z}}{2}-1\right)t\right]
=\left(\frac{e^z+e^{-z}}{2}-1\right)t$$

Indeed, the plain coefficients of this series are the numbers $(1+(-1)^n)t/n!$, so the Taylor coefficients of this series are the numbers $0,t,0,t,\ldots\,$, as claimed.
\end{proof}

At a more theoretical level, we have the following result, generalizing (4,5) above, and which is something very useful, when dealing with the compound Poisson laws:

\begin{theorem}
For a compound Poisson law $p_\nu$ we have
$$k_n(p_\nu)=M_n(\nu)$$
valid for any integer $n\geq1$.
\end{theorem}

\begin{proof}
We can assume, by using a continuity argument, that our measure $\nu$ is discrete, as follows, with $t_i>0$ and $x_i\in\mathbb R$, and with the sum being finite:
$$\nu=\sum_i t_i\delta_{x_i}$$

By using now the Fourier transform formula for $p_\nu$ from chapter 1, we obtain:
\begin{eqnarray*}
K_{p_\nu}(z)
&=&\log F_{p_\nu}(-iz)\\
&=&\log\exp\left[\sum_it_i(e^{x_iz}-1)\right]\\
&=&\sum_it_i\sum_{n\geq1}\frac{(x_iz)^n}{n!}\\
&=&\sum_{n\geq1}\frac{z^n}{n!}\sum_it_ix_i^n\\
&=&\sum_{n\geq1}\frac{z^n}{n!}\,M_n(\nu)
\end{eqnarray*}

Thus, we are led to the conclusion in the statement.
\end{proof}

Here are as well two supplementary computations, which are more tricky:

\begin{theorem}
For the exponential law, $e_t=te^{-tx}dx$ on $[0,\infty)$, we have
$$k_n=\frac{(n-1)!}{t^n}$$
and for the standard uniform law, $u_t=dx/t$ on $[0,t]$, the cumulants are
$$k_1=\frac{t}{2}\quad,\quad k_n=\frac{t^n}{n}\cdot B_n\ (n>1)$$
with $B_n\in\mathbb Q$ being the $n$-th Bernoulli number.
\end{theorem}

\begin{proof}
We have two assertions to be proved, the idea being as follows:

\medskip

(1) For the exponential law $e_t$ we have $F(x)=t/(t-ix)$, and we get, as desired:
\begin{eqnarray*}
K(z)
&=&\log\left(\frac{t}{t-z}\right)\\
&=&-\log\left(1-\frac{z}{t}\right)\\
&=&\sum_{n=1}^\infty\frac{1}{n}\left(\frac{z}{t}\right)^n\\
&=&\sum_{n=1}^\infty\frac{(n-1)!}{t^n}\cdot\frac{z^n}{n!}
\end{eqnarray*}

\medskip

(2) For the uniform law $u_t$ we have $F(x)=(e^{itx}-1)/(itx)$, and this gives:
$$K(z)=\log F(-iz)=\log\left(\frac{e^{tz}-1}{tz}\right)$$

In order to process this formula, consider the following function:
$$L(y)=\log\left(\frac{e^y-1}{y}\right)$$

By differentiating with respect to $y$, we obtain the following formula:
$$L'(y)=\frac{y}{e^y-1}\cdot\frac{e^yy-e^y+1}{y^2}
=\frac{y-1}{y}+\frac{1}{e^y-1}$$

Now recall that the Bernoulli numbers $B_n\in\mathbb Q$ are given by:
$$\frac{y}{e^y-1}=\sum_{n=0}^\infty B_n\cdot\frac{y^n}{n!}$$

In view of this latter formula, consider the following function:
$$M(y)=\frac{y}{2}+\sum_{n=2}^\infty B_n\cdot\frac{y^n}{n\cdot n!}$$

By differentiating with respect to $y$, and using $B_0=1$, $B_1=-1/2$, we obtain:
\begin{eqnarray*}
M'(y)
&=&\frac{1}{2}+\sum_{n=2}^\infty B_n\cdot\frac{y^{n-1}}{n!}\\
&=&\frac{1}{2}+\frac{1}{y}\sum_{n=2}^\infty B_n\cdot\frac{y^n}{n!}\\
&=&\frac{1}{2}+\frac{1}{y}\left(\frac{y}{e^y-1}-1+\frac{y}{2}\right)\\
&=&\frac{y-1}{y}+\frac{1}{e^y-1}
\end{eqnarray*}

Thus $L'=M'$, so $L=M$, and with $y=tz$, the formula $L(tz)=M(tz)$ reads:
\begin{eqnarray*}
K(z)
&=&\frac{tz}{2}+\sum_{n=2}^\infty B_n\cdot\frac{(tz)^n}{n\cdot n!}\\
&=&\frac{t}{2}\cdot z+\sum_{n=2}^\infty \frac{t^nB_n}{n}\cdot\frac{z^n}{n!}
\end{eqnarray*}

Thus the cumulants are $k_1=t/2$, and $k_n=t^nB_n/n$ at $n>1$, as stated.
\end{proof}

\section*{9b. Inversion formula}

Getting back to theory now, the sequence of cumulants $k_1,k_2,k_3,\ldots$ appears as a modification of the sequence of moments $M_1,M_2,M_3,\ldots\,$, and understanding the relation between moments and cumulants will be our next task. Let us start with:

\index{M\"obius function}
\index{lattice}
\index{lattice of partitions}

\begin{definition}
The M\"obius function of any lattice, and so of $P$, is given by
$$\mu(\pi,\nu)=\begin{cases}
1&{\rm if}\ \pi=\nu\\
-\sum_{\pi\leq\tau<\nu}\mu(\pi,\tau)&{\rm if}\ \pi<\nu\\
0&{\rm if}\ \pi\not\leq\nu
\end{cases}$$
with the construction being performed by recurrence.
\end{definition}

As an illustration, for $P(2)=\{||,\sqcap\}$, we have $\mu(||,||)=\mu(\sqcap,\sqcap)=1$. Also, $||<\sqcap$, with no intermediate partition in between, so $\mu(||,\sqcap)=-\mu(||,||)=-1$. Next, $\sqcap\not\leq||$, which gives $\mu(\sqcap,||)=0$. Thus, the M\"obius matrix $M_{\pi\nu}=\mu(\pi,\nu)$ of $P(2)$ is:
$$M=\begin{pmatrix}1&-1\\ 0&1\end{pmatrix}$$

At $k=3$ now, we have the following formula for the M\"obius matrix $M_{\pi\nu}=\mu(\pi,\nu)$, once again written with the indices picked increasing in $P(3)=\{|||,\sqcap|,\sqcap\hskip-3.2mm{\ }_|\,,|\sqcap,\sqcap\hskip-0.7mm\sqcap\}$:
$$M=\begin{pmatrix}
1&-1&-1&-1&2\\
0&1&0&0&-1\\
0&0&1&0&-1\\
0&0&0&1&-1\\
0&0&0&0&1
\end{pmatrix}$$

Back to the general case now, the main interest in the M\"obius function comes from the M\"obius inversion formula, which can be formulated as follows:

\index{M\"obius inversion}
\index{M\"obius matrix}

\begin{theorem}
We have the following implication,
$$f(\pi)=\sum_{\nu\leq\pi}g(\nu)
\quad\implies\quad
g(\pi)=\sum_{\nu\leq\pi}\mu(\nu,\pi)f(\nu)$$
valid for any two functions $f,g:P(n)\to\mathbb C$.
\end{theorem}

\begin{proof}
Consider the adjacency matrix of $P$, given by the following formula:
$$A_{\pi\nu}=\begin{cases}
1&{\rm if}\ \pi\leq\nu\\
0&{\rm if}\ \pi\not\leq\nu
\end{cases}$$

Our claim is that the inverse of this matrix is the M\"obius matrix of $P$, given by:
$$M_{\pi\nu}=\mu(\pi,\nu)$$

Indeed, the above matrix $A$ is upper triangular, and when trying to invert it, we are led to the recurrence in Definition 9.9, so to the M\"obius matrix $M$. Thus we have:
$$M=A^{-1}$$

Thus, in practice, we are led to the inversion formula in the statement.
\end{proof}

As an illustration here, for $P(2)$ the formula $M=A^{-1}$ appears as follows:
$$\begin{pmatrix}1&-1\\ 0&1\end{pmatrix}=
\begin{pmatrix}1&1\\ 0&1\end{pmatrix}^{-1}$$

At $k=3$ now, the formula $M=A^{-1}$ for $P(3)$ reads:
$$\begin{pmatrix}
1&-1&-1&-1&2\\
0&1&0&0&-1\\
0&0&1&0&-1\\
0&0&0&1&-1\\
0&0&0&0&1
\end{pmatrix}=
\begin{pmatrix}
1&1&1&1&1\\
0&1&0&0&1\\
0&0&1&0&1\\
0&0&0&1&1\\
0&0&0&0&1
\end{pmatrix}^{-1}$$

In general, the formula $M=A^{-1}$ looks quite similar. For more on this we refer to the Lindst\"om determinant formula from chapter 4, which is related to all this.

\bigskip

With these ingredients in hand, let us go back to probability. We first have:

\index{cumulant}
\index{generalized cumulant}
\index{multiplicativity over blocks}

\begin{definition}
We define quantities $M_\pi(f),k_\pi(f)$, depending on partitions 
$$\pi\in P(k)$$
by starting with $M_n(f),k_n(f)$, and using multiplicativity over the blocks. 
\end{definition}

To be more precise, the convention here is that for the one-block partition $1_n\in P(n)$, the corresponding moment and cumulant are the usual ones, namely:
$$M_{1_n}(f)=M_n(f)\quad,\quad k_{1_n}(f)=k_n(f)$$

Then, for an arbitrary partition $\pi\in P(k)$, we decompose this partition into blocks, having sizes $b_1,\ldots,b_s$, and we set, by multiplicativity over blocks:
$$M_\pi(f)=M_{b_1}(f)\ldots M_{b_s}(f)\quad,\quad k_\pi(f)=k_{b_1}(f)\ldots k_{b_s}(f)$$

With this convention, following Rota, we can now formulate a key result, fully clarifying the relation between moments and cumulants, as follows:

\index{moment-cumulant formula}

\begin{theorem}
We have the moment-cumulant formulae
$$M_n(f)=\sum_{\nu\in P(n)}k_\nu(f)\quad,\quad 
k_n(f)=\sum_{\nu\in P(n)}\mu(\nu,1_n)M_\nu(f)$$
or, equivalently, we have the moment-cumulant formulae
$$M_\pi(f)=\sum_{\nu\leq\pi}k_\nu(f)\quad,\quad 
k_\pi(f)=\sum_{\nu\leq\pi}\mu(\nu,\pi)M_\nu(f)$$
where $\mu$ is the M\"obius function of $P(n)$.
\end{theorem}

\begin{proof}
There are several things going on here, the idea being as follows:

\medskip

(1) According to our conventions, and to M\"obius inversion, the 4 formulae in the statement are all equivalent. In what follows we will focus on the first 2 formulae.

\medskip

(2) Before anything, let us first work out some examples. At $n=1,2,3$ the moment formula in the statement gives, in tune with Proposition 9.4:
$$M_1=k_|=k_1$$
$$M_2=k_{|\,|}+k_\sqcap=k_1^2+k_2$$
$$M_3=k_{|\,|\,|}+k_{\sqcap|}+k_{\sqcap\hskip-2.8mm{\ }_|}+k_{|\sqcap}+k_{\sqcap\hskip-0.5mm\sqcap}=k_1^3+3k_1k_2+k_3$$

At $n=4$ now, the computation is as follows, again in tune with Proposition 9.4:
\begin{eqnarray*}
M_4
&=&k_{|\,|\,|}+(\underbrace{k_{\sqcap\,|\,|}+\ldots}_{6\ terms})+(\underbrace{k_{\sqcap\,\sqcap}+\ldots}_{3\ terms})+(\underbrace{k_{\sqcap\hskip-0.5mm\sqcap\,|}+\ldots}_{4\ terms})+k_{\sqcap\hskip-0.5mm\sqcap\hskip-0.5mm\sqcap}\\
&=&k_1^4+6k_1^2k_2+3k_2^2+4k_1k_3+k_4
\end{eqnarray*}

As for the cumulant formula in the statement, at $n=1,2,3$ this gives the following formulae for the low order cumulants, again in tune with Proposition 9.4:
$$k_1=M_|=M_1$$
$$k_2=(-1)M_{|\,|}+M_\sqcap=-M_1^2+M_2$$
$$k_3=2M_{|\,|\,|}+(-1)M_{\sqcap|}+(-1)M_{\sqcap\hskip-2.8mm{\ }_|}+(-1)M_{|\sqcap}+M_{\sqcap\hskip-0.5mm\sqcap}=2M_1^3-3M_1M_2+M_3$$

Finally, at $n=4$, after computing the M\"obius function of $P(4)$, we obtain the following formula for the fourth cumulant, again in tune with Proposition 9.4:
\begin{eqnarray*}
k_4
&=&(-6)M_{|\,|\,|}+2(\underbrace{M_{\sqcap\,|\,|}+\ldots}_{6\ terms})+(-1)(\underbrace{M_{\sqcap\,\sqcap}+\ldots}_{3\ terms})+(-1)(\underbrace{M_{\sqcap\hskip-0.5mm\sqcap\,|}+\ldots}_{4\ terms})+M_{\sqcap\hskip-0.5mm\sqcap\hskip-0.5mm\sqcap}\\
&=&-6M_1^4+12M_1^2M_2-3M_2^2-4M_1M_3+M_4
\end{eqnarray*}

(3) After all these preliminaries, time now to get to work, and prove the result. As mentioned above, our formulae are all equivalent, and it is enough to prove just one of them. We will prove in what follows the first formula, namely:
$$M_n(f)=\sum_{\nu\in P(n)}k_\nu(f)$$

In order to do this, we use the very definition of the cumulants, namely:
$$\log E(e^{zf})=\sum_{s=1}^\infty k_s(f)\,\frac{z^s}{s!}$$

By exponentiating, we obtain from this the following formula:
$$E(e^{zf})=\exp\left(\sum_{s=1}^\infty k_s(f)\,\frac{z^s}{s!}\right)$$

(4) Let us first compute the function on the left. This is easily done, as follows:
$$E(e^{zf})
=E\left(\sum_{n=0}^\infty\frac{(zf)^n}{n!}\right)
=\sum_{n=0}^\infty M_n(f)\,\frac{z^n}{n!}$$

(5) Regarding now the function on the right, this is given by:
\begin{eqnarray*}
\exp\left(\sum_{s=1}^\infty k_s(f)\,\frac{z^s}{s!}\right)
&=&\sum_{p=0}^\infty\frac{\left(\sum_{s=1}^\infty k_s(f)\,\frac{z^s}{s!}\right)^p}{p!}\\
&=&\sum_{p=0}^\infty\frac{1}{p!}\sum_{s_1=1}^\infty k_{s_1}(f)\,\frac{z^{s_1}}{s_1!}\ldots\ldots\sum_{s_p=1}^\infty k_{s_p}(f)\,\frac{z^{s_p}}{s_p!}\\
&=&\sum_{p=0}^\infty\frac{1}{p!}\sum_{s_1=1}^\infty\ldots\sum_{s_p=1}^\infty k_{s_1}(f)\ldots k_{s_p}(f)\,\frac{z^{s_1+\ldots+s_p}}{s_1!\ldots s_p!}
\end{eqnarray*}

(6) But the point now is that all this leads us into partitions. Indeed, we are summing over indices $s_1,\ldots,s_p\in\mathbb N$, which can be thought of as corresponding to a partition of $n=s_1+\ldots+s_p$. So, let us rewrite our sum, as a sum over partitions. For this purpose, recall that the number of partitions $\nu\in P(n)$ having blocks of sizes $s_1,\ldots,s_p$ is:
$$\binom{n}{s_1,\ldots,s_p}=\frac{n!}{p_1!\ldots p_s!}$$

Also, when resumming over partitions, there will be a $p!$ factor as well, coming from the permutations of $s_1,\ldots,s_p$. Thus, our sum can be rewritten as follows:
\begin{eqnarray*}
\exp\left(\sum_{s=1}^\infty k_s(f)\,\frac{z^s}{s!}\right)
&=&\sum_{n=0}^\infty\sum_{p=0}^\infty\frac{1}{p!}\sum_{s_1+\ldots+s_p=n}k_{s_1}(f)\ldots k_{s_p}(f)\,\frac{z^n}{s_1!\ldots s_p!}\\
&=&\sum_{n=0}^\infty\frac{z^n}{n!}\sum_{p=0}^\infty\frac{1}{p!}\sum_{s_1+\ldots+s_p=n}\binom{n}{s_1,\ldots,s_p}k_{s_1}(f)\ldots k_{s_p}(f)\\
&=&\sum_{n=0}^\infty\frac{z^n}{n!}\sum_{\nu\in P(n)}k_\nu(f)
\end{eqnarray*}

(7) We are now in position to conclude. According to the above, we have:
$$\sum_{n=0}^\infty M_n(f)\,\frac{z^n}{n!}=\sum_{n=0}^\infty\frac{z^n}{n!}\sum_{\nu\in P(n)}k_\nu(f)$$

Thus, we have the following formula, valid for any $n\in\mathbb N$:
$$M_n(f)=\sum_{\nu\in P(n)}k_\nu(f)$$

We are therefore led to the conclusions in the statement.
\end{proof}

\section*{9c. Free cumulants}

Getting back now to free probability, following Speicher \cite{sp1}, let us formulate:

\index{cumulant}
\index{free cumulant}

\begin{definition}
The free cumulants $\kappa_n(a)$ of a variable $a\in A$ are defined by
$$R_a(z)=\sum_{n=1}^\infty\kappa_n(a)z^{n-1}$$
with $R$ being as usual the Voiculescu $R$-transform.
\end{definition}

As before with the classical cumulants, we have a number of basic examples and illustrations, and a number of basic general results. Let us start with some numerics:

\begin{proposition}
The free cumulants $\kappa_1,\kappa_2,\kappa_3,\ldots$ appear as a modification of the moments $M_1,M_2,M_3,\ldots\,$, and uniquely determine $\mu$. We have
$$\kappa_1=M_1$$
$$\kappa_2=-M_1^2+M_2$$
$$\kappa_3=2M_1^3-3M_1M_2+M_3$$
$$\kappa_4=-5M_1^4+10M_1^2M_2-2M_2^2-4M_1M_3+M_4$$
$$\vdots$$
in one sense, and in the other sense we have
$$M_1=\kappa_1$$
$$M_2=\kappa_1^2+\kappa_2$$
$$M_3=\kappa_1^3+3\kappa_1\kappa_2+\kappa_3$$
$$M_4=\kappa_1^4+6\kappa_1^2\kappa_2+2\kappa_2^2+4\kappa_1\kappa_3+\kappa_4$$
$$\vdots$$
with in both cases the correspondence being polynomial, with integer coefficients.
\end{proposition}

\begin{proof}
This is something quite instructive, a bit more complicated than the similar result in the classical case, worth discussing in detail, the idea being as follows:

\medskip

(1) We know that the Cauchy transform is the following function:
$$G(\xi)=\sum_{n=0}^\infty\frac{M_n}{\xi^{n+1}}$$

By inverting, $K(G(\xi))=\xi$, the numbers $\kappa_n$ appear as follows, with $G=G(\xi)$:
$$\frac{1}{G}+\sum_{n=1}^\infty\kappa_nG^{n-1}=\xi$$

With $\xi=z^{-1}$ our equation takes the following form, with $G=G(z^{-1})$:
$$\frac{z}{G}+z\sum_{n=1}^\infty\kappa_nG^{n-1}=1$$

With $\psi=\sum_{n=1}^\infty M_nz^n$ we have $G/z=1+\psi$, and our equation becomes:
$$\frac{1}{1+\psi}+\sum_{n=1}^\infty\kappa_nz^n(1+\psi)^{n-1}=1$$

(2) By expanding the fraction on the left, our equation becomes:
$$\sum_{n=0}^\infty(-\psi)^n+\sum_{n=1}^\infty\kappa_nz^n(1+\psi)^{n-1}=1$$

Moreover, we can cancel the 1 term on both sides, and our equation becomes:
$$\sum_{n=1}^\infty(-\psi)^n+\sum_{n=1}^\infty\kappa_nz^n(1+\psi)^{n-1}=0$$

Alternatively, we can write our equation in the following way:
$$\sum_{n=1}^\infty\kappa_nz^n(1+\psi)^{n-1}=-\sum_{n=1}^\infty(-\psi)^n$$

(3) Good news, this latter equation is something that we are eventually happy with. By remembering that we have $\psi=\sum_{n=1}^\infty M_nz^n$, our equation looks as follows:
\begin{eqnarray*}
&&\kappa_1z+\kappa_2z^2(1+M_1z+M_2z^2+\ldots)+\kappa_3z^3(1+M_1z+M_2z^2+\ldots)^2+\ldots\\
&=&(M_1z+M_2z^2+\ldots)-(M_1z+M_2z^2+\ldots)^2+(M_1z+M_2z^2+\ldots)^3-\ldots
\end{eqnarray*}

Which looks nice and exploitable, and with this, we are basically done.

\medskip

(4) Indeed, in practice, by looking at the terms of order $1,2,3,4$ in the above formula, we obtain the following equations, relating the free cumulants and moments:
$$\kappa_1=M_1$$
$$\kappa_2=M_2-M_1^2$$
$$\kappa_2M_1+\kappa_3=M_3-2M_1M_2+M_1^3$$
$$\kappa_4+2\kappa_3M_1+\kappa_2M_2=M_4-2M_1M_3-M_2^2+3M_1^2M_2-M_1^4$$

Thus, we are led to the formulae of $\kappa_1,\kappa_2,\kappa_3,\kappa_4$ in the statement. And then, by inversion, to the formulae of $M_1,M_2,M_3,M_4$ as well, as desired. 
\end{proof}

Observe the similarity with the formulae for classical cumulants. In fact, we have:

\begin{conclusion}
The first three classical and free cumulants coincide, but
$$k_4=-6M_1^4+12M_1^2M_2-3M_2^2-4M_1M_3+M_4$$
$$\kappa_4=-5M_1^4+10M_1^2M_2-2M_2^2-4M_1M_3+M_4$$
are different, and the same happens at higher order as well.
\end{conclusion}

This is something quite interesting, and we will back later with a conceptual explanation for this, via partitions, the idea being that all this comes from:
$$P(n)=NC(n)\iff n\leq 3$$

But more on this later. At the level of basic general results, in analogy with what we know about the classical cumulants, we first have the following result:

\begin{theorem}
The free cumulants have the following properties:
\begin{enumerate}
\item $\kappa_n(\lambda a)=\lambda^n\kappa_n(a)$.

\item $\kappa_n(a+b)=\kappa_n(a)+\kappa_n(b)$, if $a,b$ are free.
\end{enumerate}
\end{theorem}

\begin{proof}
This is something very standard, the idea being as follows:

\medskip

(1) We have the following Cauchy transform computation:
$$G_{\lambda a}(\xi)
=\int_\mathbb R\frac{d\mu_{\lambda a}(x)}{\xi-x}
=\int_\mathbb R\frac{d\mu_a(y)}{\xi-\lambda y}
=\frac{1}{\lambda}\int_\mathbb R\frac{d\mu_a(y)}{\xi/\lambda-y}
=\frac{1}{\lambda}\,G_a\left(\frac{\xi}{\lambda}\right)$$

But this gives the following formula, by using the definition of the $R$-transform:
$$G_{\lambda a}\left(\lambda R_a(\lambda z)+\frac{1}{z}\right)
=\frac{1}{\lambda}\,G_a\left(R_a(\lambda z)+\frac{1}{\lambda z}\right)
=\frac{1}{\lambda}\cdot\lambda z
=z$$

Thus we have the formula $R_{\lambda a}(z)=\lambda R_a(\lambda z)$, which gives (1). 

\medskip

(2) This follows indeed from the fact that the $R$-transform linearizes $\boxplus$.
\end{proof}

Again in analogy with our previous results from the classical case, at the level of the main examples of free cumulants, we have the following result:

\begin{theorem}
The sequence of free cumulants $\kappa_1,\kappa_2,\kappa_3,\kappa_4,\ldots$ is as follows:
\begin{enumerate}
\item For $\mu=\delta_t$ the free cumulants are $t,0,0,0,\ldots$

\item For $\mu=\gamma_t$ the free cumulants are $0,t,0,0,\ldots$

\item For $\mu=\gamma_t^t$ the free cumulants are $t,t,0,0,\ldots$

\item For $\mu=\pi_t$ the free cumulants are $t,t,t,t,\ldots$

\item For $\mu=\pi_t^2$ the free cumulants are $0,t,0,t,\ldots$
\end{enumerate}
\end{theorem}

\begin{proof}
The proofs are analogous to those from the classical case, as follows:

\medskip

(1) For $\mu=\delta_t$ we have $G_\mu(\xi)=1/(\xi-t)$, and so $R_\mu(z)=t$, as desired.

\medskip

(2) For $\mu=\gamma_t$ we have, as computed before, $R_\mu(z)=tz$, as desired.

\medskip

(3) For $\mu=\gamma_t^t$ we have, as computed before, $R_\mu(z)=t+tz$, as desired.

\medskip

(4) For $\mu=\pi_t$ we have, also from before, $R_\mu(z)=t/(1-z)$, as desired.

\medskip

(5) For $\mu=\pi_t^2$ we have also from before, $R_\mu(z)=tz/(1-z^2)$, as desired.
\end{proof}

At a more theoretical level, we have the following result, generalizing (4,5) above, and which is something very useful, when dealing with the compound free Poisson laws:

\begin{theorem}
For a compound free Poisson law $\pi_\nu$ we have
$$\kappa_n(\pi_\nu)=M_n(\nu)$$
valid for any integer $n\geq1$.
\end{theorem}

\begin{proof}
With $\nu=\sum_ic_i\delta_{x_i}$ we have, using the formula of $R$ from chapter 8:
\begin{eqnarray*}
R_{\pi_\nu}(z)
&=&\sum_i\frac{c_ix_i}{1-x_iz}\\
&=&\sum_{n\geq1}z^{n-1}\sum_ic_ix_i^n\\
&=&\sum_{n\geq 1}z^{n-1}\,M_n(\nu)
\end{eqnarray*}

Thus, we are led to the conclusion in the statement.
\end{proof}

Here are as well two supplementary computations, which are more tricky:

\begin{theorem}
For the arcsine law, $\alpha_t=1/(\pi\sqrt{4t-x^2})dx$ on $[-2\sqrt{t},2\sqrt{t}]$,
$$\kappa_{2r+2}=(-1)^r2t^{r+1}\frac{1}{r+1}\binom{2r}{r}$$
and for the standard uniform law, $u_t=dx/t$ on $[0,t]$, the free cumulants are
$$\kappa_n=\frac{(-t)^n}{n!}\cdot B_n$$
with $B_n\in\mathbb Q$ being the $n$-th Bernoulli number.
\end{theorem}

\begin{proof}
We have two assertions to be proved, the idea being as follows:

\medskip

(1) The moments of the arcsine law $\alpha_t$, as parametrized in the statement, are rescaled central binomial coefficients, $M_{2k}=t^k\binom{2k}{k}$, and this gives succesively:
$$M(z)=\frac{1}{\sqrt{1-4tz^2}}\quad,\quad G(\xi)=\frac{1}{\sqrt{\xi^2-4t}}$$

Next, from the formula of $G$ we can obtain an equation for $R$, as follows:
\begin{eqnarray*}
G^2(\xi^2-4t)=1
&\implies&z^2((R+z^{-1})^2-4t)=1\\
&\implies&(zR+1)^2-4tz^2=1\\
&\implies&zR=\sqrt{1+4tz^2}-1
\end{eqnarray*}

Now by extracting the square root, using the generalized binomial formula, we get:
\begin{eqnarray*}
zR
&=&2\sum_{r=0}^\infty(-1)^r\frac{1}{r+1}\binom{2r}{r}(tz^2)^{r+1}\\
&=&\sum_{r=0}^\infty(-1)^r2t^{r+1}\frac{1}{r+1}\binom{2r}{r}\cdot z^{2r+2}
\end{eqnarray*}

Thus the free cumulants are $\kappa_{2r+2}=(-1)^r2t^{r+1}\frac{1}{r+1}\binom{2r}{r}$, as stated.

\medskip

(2) For the uniform law $u_t$ the Cauchy transform can be computed as follows:
$$G(\xi)=\sum_{k=0}^\infty\frac{t^k}{k+1}\,\xi^{-k-1}=\frac{1}{t}\sum_{k=0}^\infty\frac{(t\xi^{-1})^{k+1}}{k+1}=-\frac{\log(1-t\xi^{-1})}{t}$$

Next, we can obtain from this an equation for the $R$-transform, as follows:
\begin{eqnarray*}
-tG=\log(1-t\xi^{1})
&\implies&e^{-tG}=1-t\xi^{-1}\\
&\implies&e^{-tz}=1-\frac{t}{R+z^{-1}}\\
&\implies&R+z^{-1}=\frac{t}{1-e^{-tz}}\\
&\implies&zR=-1+\frac{tz}{1-e^{-tz}}
\end{eqnarray*}

But the fraction on the right is computed by the Bernoulli numbers, and we get:
$$zR=-1+\sum_{n=0}^\infty B_n\cdot\frac{(-tz)^n}{n!}=\sum_{n=1}^\infty \frac{(-t)^n}{n!}\cdot B_n\cdot z^n$$

We conclude that the free cumulants are $\kappa_n=(-t)^nB_n/n!$, as stated.
\end{proof}

As before in the classical case, we can define generalized free cumulants, $\kappa_\pi(a)$ with $\pi\in P(k)$, by starting with the numeric free cumulants $\kappa_n(a)$, as follows:

\begin{definition}
We define free cumulants $\kappa_\pi(a)$, depending on partitions 
$$\pi\in P(k)$$
by starting with $\kappa_n(a)$, and using multiplicativity over the blocks. 
\end{definition}

To be more precise, the convention here is that for the one-block partition $1_n\in P(n)$, the corresponding free cumulant is the usual one, namely:
$$\kappa_{1_n}(a)=\kappa_n(a)$$

Then, for an arbitrary partition $\pi\in P(k)$, we decompose this partition into blocks, having sizes $b_1,\ldots,b_s$, and we set, by multiplicativity over blocks:
$$\kappa_\pi(a)=\kappa_{b_1}(a)\ldots\kappa_{b_s}(a)$$

With this convention, we have the following result, due to Speicher \cite{sp1}:

\begin{theorem}
We have the moment-cumulant formulae
$$M_n(a)=\sum_{\nu\in NC(n)}\kappa_\nu(a)\quad,\quad 
\kappa_n(a)=\sum_{\nu\in NC(n)}\mu(\nu,1_n)M_\nu(a)$$
or, equivalently, we have the moment-cumulant formulae
$$M_\pi(a)=\sum_{\nu\leq\pi}\kappa_\nu(a)\quad,\quad 
\kappa_\pi(a)=\sum_{\nu\leq\pi}\mu(\nu,\pi)M_\nu(a)$$
where $\mu$ is the M\"obius function of $NC(n)$.
\end{theorem}

\begin{proof}
As before in the classical case, the 4 formulae in the statement are equivalent, via M\"obius inversion. Thus, it is enough to prove one of them, and we will discuss the first formula, which in practice is the most useful one. Thus, we must prove that:
$$M_n(a)=\sum_{\nu\in NC(n)}\kappa_\nu(a)$$

(1) In order to prove this, let us get back to the construction of the free cumulants, from Definition 9.13. The Cauchy transform of $a$ is the following function:
$$G_a(\xi)=\sum_{n=0}^\infty\frac{M_n(a)}{\xi^{n+1}}$$

Consider the inverse of this Cauchy transform $G_a$, with respect to composition:
$$K_a(G_a(\xi))=\xi$$

According to Definition 9.13, the free cumulants $\kappa_n(a)$ appear then as follows:
$$K_a(z)=\frac{1}{z}+\sum_{n=1}^\infty\kappa_n(a)z^{n-1}$$

(2) In practice, as seen in the proof of Proposition 9.14, this formula leads to:
\begin{eqnarray*}
&&\kappa_1z+\kappa_2z^2(1+M_1z+M_2z^2+\ldots)+\kappa_3z^3(1+M_1z+M_2z^2+\ldots)^2+\ldots\\
&=&(M_1z+M_2z^2+\ldots)-(M_1z+M_2z^2+\ldots)^2+(M_1z+M_2z^2+\ldots)^3-\ldots
\end{eqnarray*}

But, in case you have followed the proof of Proposition 9.14, you know how to exploit this formula at order $n=1,2,3,4$. The same method works in general, and after some computations, this leads to the formula that we want to establish, namely:
$$M_n(a)=\sum_{\nu\in NC(n)}\kappa_\nu(a)$$

(3) We are therefore led to the conclusions in the statement. All this was of course quite brief, and for details here, we refer for instance to Nica-Speicher \cite{nsp}.
\end{proof}

Observe that the above result leads among others to a more conceptual explanation for Conclusion 9.15, with the equalities and non-equalities there simply coming from:
$$P(n)=NC(n)\iff n\leq3$$

As a conclusion to all this, we have now theories of classical and free cumulants, due to Rota and Speicher, coming from $P(n)$ and $NC(n)$. Which is very nice.

\section*{9d. The bijection}

With the above classical and free cumulant theory in hand, we can now formulate the following simple definition, making the connection between classical and free:

\index{free version}
\index{classical version}

\begin{definition}
We say that two measures are in classical/free correspondence
$$m\in\mathcal P(\mathbb R)\quad\longleftrightarrow\quad\mu\in\mathcal P(\mathbb R)$$
when the classical cumulants of $m$ coincide with the free cumulants of $\mu$.
\end{definition}

As a first observation, this definition fits with all the classical and free probability theory developed so far in this book, in the real case, and we have:

\begin{theorem}
In the standard cube of basic real probability limiting measures
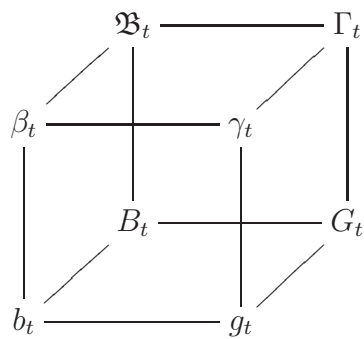
$$\xymatrix@R=16pt@C=20pt{
&\pi_t^2\ar@{-}[rr]\ar@{-}[dd]&&\gamma_t\ar@{-}[dd]\\
\pi_t\ar@{-}[rr]\ar@{-}[dd]\ar@{-}[ur]&&\gamma_t^t\ar@{-}[dd]\ar@{-}[ur]\\
&p_t^2\ar@{-}[rr]\ar@{-}[uu]&&g_t\ar@{.}[uu]\\
p_t\ar@{-}[uu]\ar@{-}[ur]\ar@{-}[rr]&&g_t^t\ar@{-}[uu]\ar@{-}[ur]
}$$
the upper measures appear as free versions of the lower measures.
\end{theorem}

\begin{proof}
This comes indeed from our various cumulant computations, and more specifically from Theorems 9.6 and 9.17, which give exactly the same sequences $k_n,\kappa_n$.
\end{proof} 

Still based on our various cumulant computations above, we have as well:

\begin{theorem}
Given a compactly supported measure $\nu\in\mathcal P(\mathbb R)$, we have
$$p_\nu\quad\longleftrightarrow\quad\pi_\nu$$
classical/free correspondence, between the associated compound Poisson laws.
\end{theorem}

\begin{proof}
When $\nu$ is discrete this follows from Theorems 9.7 and 9.18, which give the same sequences $k_n,\kappa_n$. As for the general case, this follows from this, by density.
\end{proof}

At the theoretical level, we have the following analytic version of Definition 9.22:

\begin{theorem}
The correspondence $m\leftrightarrow\mu$ amounts in having
$$\log F_m(-iz)=\sum_{n=1}^\infty c_n\cdot\frac{z^n}{n!}
\quad \longleftrightarrow\quad zR_\mu(z)=\sum_{n=1}^\infty c_n\cdot z^n$$
at the level of the associated Fourier and $R$-transforms.
\end{theorem}

\begin{proof}
This comes indeed from our definition of the classical and free cumulants $k_n,\kappa_n$, as being the coefficients $c_n$ appearing on the left and on the right.
\end{proof}

Finally, in relation with our few extra computations of cumulants, we have:

\begin{remark}
It is not clear what to do with our computations for $e_t,\alpha_t$, namely
$$k_n(e_t)=\frac{(n-1)!}{t^n}\quad,\quad 
\kappa_{2r+2}(\alpha_t)=(-1)^r2t^{r+1}\frac{1}{r+1}\binom{2r}{r}$$
nor what to do with our computations for the uniform law $u_t$, namely
$$k_1(u_t)=\frac{t}{2}\ ,\  k_n(u_t)=\frac{t^n}{n}\cdot B_n\ (n>1)
\quad,\quad \kappa_n(u_t)=\frac{(-t)^n}{n!}\cdot B_n$$
with $B_n\in\mathbb Q$ being as usual the Bernoulli numbers.
\end{remark}

At a more advanced level now, following Bercovici and Pata \cite{bpa}, we can formulate the following definition, which is something even more conceptual:

\index{convolution semigroup}
\index{free convolution semigroup}
\index{Bercovici-Pata bijection}

\begin{definition}
A convolution semigroup of measures
$$\{m_t\}_{t>0}\quad:\quad m_s*m_t=m_{s+t}$$
is in Bercovici-Pata bijection with a free convolution semigroup of measures
$$\{\mu_t\}_{t>0}\quad:\quad \mu_s\boxplus\mu_t=\mu_{s+t}$$
when the classical cumulants of $m_t$ coincide with the free cumulants of $\mu_t$.
\end{definition}

As a first observation, this new definition kills Remark 9.26, because the measures there do not form convolution semigroups, due to the presence of the powers of $t>0$ in the formulae of the classical and free cumulants. Which is something very nice.

\bigskip

Getting now to the examples, Theorem 9.23 upgrades to this new setting:

\begin{theorem}
In the standard cube of basic real probability limiting measures
$$\xymatrix@R=16pt@C=20pt{
&\pi_t^2\ar@{-}[rr]\ar@{-}[dd]&&\gamma_t\ar@{-}[dd]\\
\pi_t\ar@{-}[rr]\ar@{-}[dd]\ar@{-}[ur]&&\gamma_t^t\ar@{-}[dd]\ar@{-}[ur]\\
&p_t^2\ar@{-}[rr]\ar@{-}[uu]&&g_t\ar@{.}[uu]\\
p_t\ar@{-}[uu]\ar@{-}[ur]\ar@{-}[rr]&&g_t^t\ar@{-}[uu]\ar@{-}[ur]
}$$
the upper semigroups appear as free versions of the lower semigroups.
\end{theorem}

\begin{proof}
This comes indeed from Theorem 9.23, coupled with the fact, that we know well, that we have indeed classical and free convolution semigroups.
\end{proof}

As further good news, Theorem 9.24 upgrades as well to our new setting:

\begin{theorem}
Given a compactly supported measure $\nu\in\mathcal P(\mathbb R)$, we have
$$\{p_{t\nu}\}_{t>0}\quad\longleftrightarrow\quad\{\pi_{t\nu}\}_{t>0}$$
classical/free correspondence, between the associated compound Poisson laws.
\end{theorem}

\begin{proof}
This comes indeed from Theorem 9.24, coupled with the fact, that we know well, that we have indeed classical and free convolution semigroups.
\end{proof}

Next, and with this being something that motivated Bercovici and Pata in their pioneering work \cite{bpa}, we have the following analytic version of Definition 9.27:

\begin{theorem}
The correspondence $\{m_t\}_{t>0}\leftrightarrow\{\mu_t\}_{t>0}$ amounts in having
$$\log F_{m_t}(-iz)=t\sum_{n=1}^\infty c_n\cdot\frac{z^n}{n!}
\quad \longleftrightarrow\quad zR_{\mu_t}(z)=t\sum_{n=1}^\infty c_n\cdot z^n$$
at the level of the associated Fourier and $R$-transforms.
\end{theorem}

\begin{proof}
This comes indeed from our definition of the classical and free cumulants $k_n,\kappa_n$, as being the coefficients $tc_n$ appearing on the left and on the right.
\end{proof}

Finally, in relation with Remark 9.26, as mentioned above, that is essentially buried. Let us nevertheless formulate the following question, regarding the uniform law $u_t$:

\begin{question}
In view of our formulae for $u_t=dx/t$ on $[0,t]$, namely
$$k_1(u_t)=\frac{t}{2}\ ,\  k_n(u_t)=\frac{t^n}{n}\cdot B_n\ (n>1)
\quad,\quad \kappa_n(u_t)=\frac{(-t)^n}{n!}\cdot B_n$$
where does this uniform law fit, in the classical/free probability landscape?
\end{question}

Which sounds quite philosophical, but good question this is, because there is obviously some magic going on in the above cumulant formulae, with the Bernoulli numbers appearing on both sides, and with Bercovici-Pata not providing an answer to this. 

\bigskip

Quite fun all this, hope you agree with me, who could have guessed that the modest uniform law $u_t$ will come to haunt us, as a severely complicated mathematical object. And more on this, momentology and numerology of $u_t$, in the next chapter.

\bigskip

Along the same lines, bizarre puzzles coming from the Bercovici-Pata bijection, and from free probability in general, we have the famous question, from \cite{bbl}:

\begin{question}[Belinschi, Bo\.zejko, Lehner, Speicher]
The normal law $g_1$ being freely infinitely divisible, what is its classical analogue?
\end{question}

And isn't this puzzling, weird question that we have here, obtained by playing with fire, namely applying free transforms to classical measures, although the converse can be certainly attempted too. For more on all this, we refer to their paper \cite{bbl}.

\bigskip

Moving away now from these difficulties, let us get back to the standard cube from Theorem 9.28. The measures there are quite special, and regarding them, we have:

\begin{theorem}
For the main $8$ laws, the moments are $M_k=\sum_{\pi\in D(k)}t^{|\pi|}$, with the corresponding categories of partitions $D\subset P$ being as follows:
$$\xymatrix@R=18pt@C=2pt{
&NC_{even}\ar[dl]\ar[dd]&&NC_2\ar[dl]\ar[ll]\ar[dd]\\
NC\ar[dd]&&NC_{12}\ar[dd]\ar[ll]\\
&P_{even}\ar[dl]&&P_2\ar[dl]\ar[ll]\\
P&&P_{12}\ar[ll]
}$$
Also, the classical cumulants on the bottom are $k_n=t\delta_{n\in L}$, and the free cumulants on top are $\kappa_n=t\delta_{n\in L}$, with $L\subset\mathbb N$ being the set of sizes of blocks of the partitions in $D$.
\end{theorem}

\begin{proof}
This is indeed something self-explanatory, coming from our various moment results for the laws in Theorem 9.28, established earlier in this book.
\end{proof}

Observe that what we have in Theorem 9.33 is much finer than Bercovici-Pata, or than the theory of classical and free cumulants, in general. In fact, in relation with such things, which are quite rigid, we have as well a classification result, as follows:

\begin{theorem}
Under suitable easiness and uniformity assumptions,
$$\xymatrix@R=18pt@C=2pt{
&NC_{even}\ar[dl]\ar[dd]&&NC_2\ar[dl]\ar[ll]\ar[dd]\\
NC\ar[dd]&&NC_{12}\ar[dd]\ar[ll]\\
&P_{even}\ar[dl]&&P_2\ar[dl]\ar[ll]\\
P&&P_{12}\ar[ll]
}$$
are the only classical and free categories of partitions $D\subset P$.
\end{theorem}

\begin{proof}
This is a result from \cite{bsp}, based on the easiness and uniformity notions that we met in chapter 4, in the context of the classical groups, suitably extended to the quantum group setting. We will be back to this, with details, in chapter 11 below.
\end{proof}

And with this, end of our preliminary discussion of Bercovici-Pata. Coming next, in chapter 10 we will keep studying the real case, with calculus methods, and then in chapters 11-12 we will get into the unitary case, using quantum groups and easiness.

\section*{9e. Exercises}

Welcome to liberation, and as exercises on this, we have:

\begin{exercise}
Learn more about the Bernoulli numbers, and their properties.
\end{exercise}

\begin{exercise}
Compute further cumulants, for measures of your choice.
\end{exercise}

\begin{exercise}
Learn more about M\"obius inversion, and its applications.
\end{exercise}

\begin{exercise}
Learn about the infinite divisibility of real measures.
\end{exercise}

\begin{exercise}
Clarify the proof of the free moment-cumulant formula.
\end{exercise}

\begin{exercise}
Compute further free cumulants, for measures of your choice.
\end{exercise}

\begin{exercise}
Learn about the free infinite divisibility of measures.
\end{exercise}

\begin{exercise}
Think a bit at the Bercovici-Pata puzzles formulated above.
\end{exercise}

As bonus exercise, read about freedom. Or, learn this from cats.

\chapter{Advanced calculus}

\section*{10a. The twelve laws}

Time for some more mathematics for the classical/free correspondence $m\leftrightarrow\mu$, this time using tools from advanced classical analysis, such as orthogonal polynomials, and Gram and Hankel determinants. As data for our considerations, we will have the 8 measures from the standard real cube found in the previous chapter, namely:
$$\xymatrix@R=16pt@C=20pt{
&\pi_t^2\ar@{-}[rr]\ar@{-}[dd]&&\gamma_t\ar@{-}[dd]\\
\pi_t\ar@{-}[rr]\ar@{-}[dd]\ar@{-}[ur]&&\gamma_t^t\ar@{-}[dd]\ar@{-}[ur]\\
&p_t^2\ar@{-}[rr]\ar@{-}[uu]&&g_t\ar@{.}[uu]\\
p_t\ar@{-}[uu]\ar@{-}[ur]\ar@{-}[rr]&&g_t^t\ar@{-}[uu]\ar@{-}[ur]
}$$

However, this is not all, and in view of what was said in the previous chapter, and of probability theory in general, it makes sense to add to our study 4 more measures, which are of fundamental nature too. Indeed, we first have the ubiquitous exponential and uniform laws $e_t,u_t$. Then we have the arcsine law $\alpha_t$, which is something fundamental too. And finally, in view of the notorious slippery nature of the arcsine laws, it makes sense to include as well the ``rival'' arcsine law $\alpha_t'$, supported on $[0,4t]$.

\bigskip

So, these will be our laws, and as a plan now for this chapter, we have:

\begin{plan}
We will be interested in the $12$ main laws in probability, namely
\begin{enumerate}
\item Classical laws: $g_t,p_t,g^t_t,p_t^2$,

\item Free laws: $\gamma_t,\pi_t,\gamma_t^t,\pi^2_t$,

\item Special guests: $e_t,u_t,\alpha_t,\alpha_t'$,
\end{enumerate}
and we would like to do extensive classical analysis computations for them.
\end{plan}

As a comment, all this might seem of course a bit subjective. So, you will have to trust me here, the point being that, from a more advanced point of view, involving orthogonal polynomials and related objects, these are, indeed, the 12 main laws. More later.

\bigskip

Getting started now, before doing advanced analysis, it makes sense to review what we already know about our 12 laws, along with a bit more. Here are our conventions:

\index{skewness}
\index{kurtosis}
\index{moment generating function}
\index{Cauchy transform}
\index{Fourier transform}
\index{central moments}
\index{normalized central moments}

\begin{conventions}
Given a real probability measure $\mu$, introduced as it should, with due phenomenology discussion, density and support, we will be interested in:
\begin{enumerate}
\item Moments $M_k=E(f^k)$, with $f\sim\mu$.

\item Mean $E=M_1$, variance $V=M_2-M_1^2$. 

\item Central moments $M_k'=E((f-E)^k)$.

\item Normalized central moments $M_k''=M_k'/\sqrt{V^k}$.

\item Skewness $\gamma=M_3''$, kurtosis $\kappa=M_4''$.

\item Fourier transform $F(x)=E(e^{ifx})$.

\item Moment generating function $M(z)=\sum_{k\geq0}M_kz^k$.

\item Cauchy transform $G(\xi)=\xi^{-1}M(\xi^{-1})$.

\item $R$ transform, given by $G(R(z)+z^{-1})=z$. 

\item $S$ transform, given by $M(S(z)/(1+z^{-1}))=1+z$.

\item Cumulants, given by $\log F(-iz)=\sum_{n\geq1}k_nz^n/n!$. 

\item Free cumulants, given by $zR(z)=\sum_{n\geq1}\kappa_nz^n$.
\end{enumerate}
\end{conventions}

To be more precise, we are already familiar with the above quantities, the only newcomers being the central moments $M_k'$ and normalized central moments $M_k''$, whose interest is self-explanatory, and then the skewness $\gamma=M_3''$ and kurtosis $\kappa=M_4''$, which are quantities coming from statistics, whose abstract interest is self-explanatory too. In relation with the latter, observe that we have the following formulae, with $\sigma=\sqrt{V}$:
$$\gamma=E\left[\left(\frac{f-E}{\sigma}\right)^3\right]\qquad,\qquad 
\kappa=E\left[\left(\frac{f-E}{\sigma}\right)^4\right]$$

Well, a total of $12\times 12=144$ questions to be discussed, that sounds fun. Regarding the normal law $g_t$, skipping the phenomenology that we know well, we have:

\index{normal law}

\begin{theorem}
The normal law $g_t$ of parameter $t>0$, having density
$$g_t=\frac{1}{\sqrt{2\pi t}}\,e^{-x^2/2t}dx$$
on $\mathbb R$, has the following properties:
\begin{enumerate}
\item $M_{2k}=M_{2k}'=t^k(2k)!!$.

\item $M_{2k}''=(2k)!!$.

\item $E=0$, $V=t$, $\gamma=0$, $\kappa=3$.

\item $F(x)=e^{-tx^2/2}$.

\item $k_n=0,t,0,0,\ldots$
\end{enumerate}
\end{theorem}

\begin{proof}
This is something that we know well, from chapters 1 and 9, and with the missing items corresponding to quantities which are not really computable.
\end{proof}

Excuse me, but cat is here, meowing something. What is it, cat?

\begin{cat}
In my opinion, computing the $R$-transform of $g_t$ would be an interesting thing. Why focusing on trivialities only.
\end{cat}

Well, what can I say. Didn't I tell you cat that this is an introductory book, for such damn things the reader can consult afterwards the paper \cite{bbl} of my friends Franz, Marek, Roland and \c Serban, that by the way I already cited in chapter 9. In addition, you're ruining the layout of my book, my carefully prepared 12 theorems were supposed to come without interruption, like a symphony. What about chasing some mice outside.

\bigskip

So, skipping what cat says, and getting now to the Poisson law $p_t$, we have:

\index{Poisson law}

\begin{theorem}
The Poisson law $p_t$ of parameter $t>0$, having density
$$p_t=e^{-t}\sum_{k\geq0}\frac{t^k}{k!}\,\delta_k$$
on $\mathbb N$, has the following properties:
\begin{enumerate}
\item $M_k=\sum_{\pi\in P(k)}t^{|\pi|}$.

\item $E=t$, $V=t$, $\gamma=1/\sqrt{t}$, $\kappa=3+1/t$.

\item $F(x)=\exp((e^{ix}-1)t)$.

\item $k_n=t,t,t,t,\ldots$
\end{enumerate}
\end{theorem}

\begin{proof}
These are again things that we know, from chapters 1 and 9, and with the computation of $\gamma,\kappa$, with things miraculously simplifying, being a must-do for you.
\end{proof}

Next, getting to the shifted normal law $g_t^t$, we have, regarding it:

\index{shifted normal law}

\begin{theorem}
The shifted normal law $g_t^t$ of parameter $t>0$, having density
$$g_t^t=\frac{1}{\sqrt{2\pi t}}\,e^{-(x-t)^2/2t}dx$$
on $\mathbb R$, has the following properties:
\begin{enumerate}
\item $M_k=\sum_{\pi\in P_{12}(k)}t^{|\pi|}$.

\item $M_k=\sum_{r=0}^{[k/2]}\binom{k}{2r}(2r)!!t^{k-r}$.

\item $M_k=tM_{k-1}+(k-1)tM_{k-2}$.

\item $M_{2k}'=t^k(2k)!!$.

\item $M_{2k}''=(2k)!!$.

\item $E=t$, $V=t$, $\gamma=0$, $\kappa=3$.

\item $F(x)=\exp(itx-tx^2/2)$.

\item $k_n=t,t,0,0,\ldots$
\end{enumerate}
\end{theorem}

\begin{proof}
These are again things that we know, from chapters 4 and 9, with exercise for you to establish the equivalence between the above 3 formulae for $M_k$.
\end{proof}

Regarding now our fourth classical measure, the Bessel law $p^2_t$, we have:

\index{Bessel law}

\begin{theorem}
The Bessel law $p_t^2$ of parameter $t>0$, having density
$$p^2_t=e^{-t}\sum_{k=-\infty}^\infty f_k(t/2)\delta_k\quad,\quad 
f_k(t)=\sum_{p=0}^\infty \frac{t^{|k|+2p}}{(|k|+p)!p!}$$
on $\mathbb Z$, has the following properties:
\begin{enumerate}
\item $M_{2k}=\sum_{\pi\in P_{even}(2k)}t^{|\pi|}$.

\item $E=0$, $V=t$, $\gamma=0$, $\kappa=3+1/t$.

\item $F(x)=\exp((\cos x-1)t)$.

\item $k_n=0,t,0,t,\ldots$
\end{enumerate}
\end{theorem}

\begin{proof}
These are again things that we know, from chapters 1 and 9, with exercise for you to compute the skewness $\gamma$ and the kurtosis $\kappa$.
\end{proof}

Next come our four guest measures, namely $e_t,u_t,\alpha_t,\alpha_t'$, with the exponential law $e_t$ being something classical, although lacking a free analogue, then with the uniform law $u_t$ being something quite fascinating, at the same time classical and free, although not self-dual as the Dirac masses $\delta_t$ are, and finally with the arcsine laws $\alpha_t,\alpha_t'$ being free, and combinatorially simpler than $\gamma_t,\pi_t,\gamma_t^t,\pi^2_t$, but lacking classical analogues.

\bigskip

Regarding the first of these beasts, the exponential law $e_t$, we have:

\index{exponential law}

\begin{theorem}
The exponential law $e_t$ of parameter $t>0$, having density
$$e_t=te^{-tx}\,dx$$
on $[0,\infty)$, has the following properties:
\begin{enumerate}
\item $M_k=k!/t^k$.

\item $E=1/t$, $V=1/t^2$, $\gamma=2$, $\kappa=9$.

\item $F(x)=t/(t-ix)$.

\item $k_n=(n-1)!/t^n$.
\end{enumerate}
\end{theorem}

\begin{proof}
We can certainly talk about $e_t$, with the mass one property coming from $\int_0^\infty e^{-x}dx=1$. This law appears for instance in the context of particle decay, because if we denote by $t$ the decay rate, that is, the probability per unit time that the particle will disintegrate, then the probability decay function will be $e_t$. Next, we have:
$$M_k
=t\int_0^\infty x^ke^{-tx}dx
=\frac{k}{t}\cdot t\int_0^\infty x^{k-1}e^{-tx}dx
=\frac{k}{t}\cdot M_{k-1}$$

Thus $M_k=k!/t^k$, which gives the above formulae of $E,V,\gamma,\kappa$. Next, we have:
$$F(x)
=t\int_0^\infty e^{ixy}e^{-ty}dy
=t\left[\frac{e^{(ix-t)y}}{ix-t}\right]_0^\infty
=\frac{t}{t-ix}$$

Finally, the formula of the cumulants is from chapter 9.
\end{proof}

Next, regarding the uniform law $u_t$, which lies somewhere between classical and free, but in a quite complicated way, its basic mathematics is as follows:

\index{uniform law}

\begin{theorem}
The uniform law $u_t$ of parameter $t>0$, having density
$$u_t=\frac{1}{t}\,dx$$
on $[0,t]$, has the following properties:
\begin{enumerate}
\item $M_k=t^k/(k+1)$.

\item $M_{2k}'=(t/2)^{2k}/(2k+1)$.

\item $M_{2k}''=3^k/(2k+1)$.

\item $E=t/2$, $V=t^2/12$, $\gamma=0$, $\kappa=9/5$.

\item $F(x)=(e^{itx}-1)/itx$.

\item $M(z)=-\log(1-tz)/tz$.

\item $G(\xi)=-\log(1-t\xi^{-1})/t$.

\item $R(z)=t/(1-e^{-tz})-1/z$.

\item $\exp(-tzS)=1-tzS/(1+z)$.

\item $k_1=t/2$, and $k_n=t^nB_n/n$ at $n>1$.

\item $\kappa_n=(-t)^nB_n/n!$.
\end{enumerate}
\end{theorem}

\begin{proof}
Here the formula of $M_k$ is clear, and the formula of $M_{2k}'$ comes from: 
\begin{eqnarray*}
M_k'
&=&\frac{1}{t}\int_0^t\left(x-\frac{t}{2}\right)^kdx\\
&=&\frac{1}{t}\left[\frac{1}{k+1}\left(x-\frac{t}{2}\right)^{k+1}\right]_0^t\\
&=&\frac{1}{t}\cdot\frac{1}{k+1}\left[\left(\frac{t}{2}\right)^{k+1}-\left(-\frac{t}{2}\right)^{k+1}\right]
\end{eqnarray*}

Next, the formula of $M_{2k}''$ comes from this, then we get the formulae of $E,V,\gamma,\kappa$, and the formulae of $F,M$ are clear as well. Regarding $G,R$ and $k_n,\kappa_n$, these were discussed in chapter 9. Thus, it remains to study the $S$-transform. We have:
\begin{eqnarray*}
Mtz=-\log(1-tz)
&\implies&e^{-Mtz}=1-tz\\
&\implies&e^{-(1+\psi)tz}=1-tz\\
&\implies&e^{-(1+z)t\chi}=1-t\chi\\
&\implies&\exp\left(-(1+z)t\cdot\frac{S}{1+z^{-1}}\right)=1-t\cdot\frac{S}{1+z^{-1}}\\
&\implies&\exp(-tzS)=1-\frac{tzS}{1+z}
\end{eqnarray*}

Thus, we are led to the equation in the statement for the $S$-transform.
\end{proof}

Next, regarding the arcsine law $\alpha_t$, which is a law of free nature, we have:

\index{arcsine law}

\begin{theorem}
The arcsine law $\alpha_t$ of parameter $t>0$, having density
$$\alpha_t=\frac{1}{\pi\sqrt{4t-x^2}}\,dx$$
on $[-2\sqrt{t},2\sqrt{t}]$, has the following properties:
\begin{enumerate}
\item $M_{2k}=t^k\binom{2k}{k}$.

\item $E=0$, $V=2t$, $\gamma=0$, $\kappa=3/2$.

\item $M(z)=1/\sqrt{1-4tz^2}$.

\item $G(\xi)=1/\sqrt{\xi^2-4t}$.

\item $R(z)=(\sqrt{1+4tz^2}-1)/z$.

\item $S(z)=\sqrt{(z+2)/4tz}$.

\item $\kappa_{2r+2}=(-1)^r2t^{r+1}\frac{1}{r+1}\binom{2r}{r}$.
\end{enumerate}
\end{theorem}

\begin{proof}
To start with, regarding the phenomenology, which is something quite interesting, we can talk about the arcsine law over any interval $[a,b]$, as follows:
$$\alpha_{ab}=\frac{1}{\pi\sqrt{(x-a)(b-x)}}\,dx$$

Indeed, the mass 1 property comes from the following computation, using the linear change of variables $x\to y$ mapping $[a,b]\to[-1,1]$, which makes $a,b$ dissapear:
$$\int_a^b\frac{1}{\sqrt{(x-a)(b-x)}}\,dx=\int_{-1}^1\frac{1}{\sqrt{1-y^2}}\,dy=\Big[\arcsin y\Big]_{-1}^1=\pi$$

In order to choose a winner $\alpha_1$ among these laws, we can look at random walks on $\mathbb Z$. The number of length $2k$ loops based at 0 being the central binomial coefficient $M_{2k}=\binom{2k}{k}$, the winner is the arcsine law on $[-2,2]$, having these numbers as moments. Then, $\alpha_t$ with $t>0$ appears by dilation, which gives $M_{2k}=t^k\binom{2k}{k}$, but with the choice of the scaling being something quite subjective. Next, we talked about $M,G,R$ and $\kappa_{2r+2}$ in chapter 9, so it remains to compute $S$, and this can be done as follows:
\begin{eqnarray*}
M^2(1-4tz^2)=1
&\implies&(1+\psi)^2(1-4tz^2)=1\\
&\implies&(1+z)^2(1-4t\chi^2)=1\\
&\implies&(1+z)^2-4tz^2S^2=1\\
&\implies&4tzS^2=z+2\\
&\implies&S=\sqrt{\frac{z+2}{4tz}}
\end{eqnarray*}

Thus, we are led to the conclusions in the statement.
\end{proof}

Next, we have the rival arcsine law $\alpha_t'$, whose basic mathematics is as follows:

\index{rival arcsine law}

\begin{theorem}
The rival arcsine law $\alpha_t'$ of parameter $t>0$, having density
$$\alpha_t'=\frac{1}{\pi\sqrt{4tx-x^2}}\,dx$$
on $[0,4t]$, has the following properties:
\begin{enumerate}
\item $M_k=t^k\binom{2k}{k}$.

\item $E=2t$, $V=2t^2$, $\gamma=0$, $\kappa=3/2$.

\item $M(z)=1/\sqrt{1-4tz}$.

\item $G(\xi)=1/\sqrt{\xi^2-4t\xi}$.

\item $R(z)=2t+(\sqrt{1+4t^2z^2}-1)/z$.

\item $S(z)=(z+2)/4t(z+1)$.

\item $\kappa_{2r+2}=2t\delta_{r0}+(-1)^r2t^{2r+2}\frac{1}{r+1}\binom{2r}{r}$.
\end{enumerate}
\end{theorem}

\begin{proof}
Regarding the phenomenology, $\alpha_1'$ appears again by looking at random walks on $\mathbb Z$, but this time ignoring the non-existent odd loops, I mean looking for a measure having $M_k=\binom{2k}{k}$ as moments. Then $\alpha_t'$ with $t>0$ appears by dilation, and with the choice of the scaling being a bit subjective. Observe that we have the following formula, showing that the relation between $\alpha_t,\alpha_t'$ is similar to the relation between $\gamma_1,\pi_1$: 
$$f\sim\alpha_t\implies f^2\sim\alpha_t'$$

Concerning now the formulae of $M,G$, these follow from those in Theorem 10.10, so it remains to compute $R,S$ and the free cumulants. And here, we first have:
\begin{eqnarray*}
G^2(\xi^2-4t\xi)=1
&\implies&z^2((R+z^{-1})^2-4t(R+z^{-1}))=1\\
&\implies&(zR+1)^2-4tz(zR+1)=1\\
&\implies&zR+1=2tz+\sqrt{4t^2z^2+1}\\
&\implies&zR=2tz+\sqrt{1+4t^2z^2}-1
\end{eqnarray*}

Now by extracting the square root, using the generalized binomial formula, we get:
$$zR=2tz+\sum_{r=0}^\infty(-1)^r2t^{2r+2}\frac{1}{r+1}\binom{2r}{r}\cdot z^{2r+2}$$

Thus we have $\kappa_{2r+2}=2t\delta_{r0}+(-1)^r2t^{2r+2}\frac{1}{r+1}\binom{2r}{r}$, as stated. Finally, we have:
\begin{eqnarray*}
M^2(1-4tz)=1
&\implies&(1+\psi)^2(1-4tz)=1\\
&\implies&(1+z)^2(1-4t\chi)=1\\
&\implies&(1+z)^2-4t(1+z)zS=1\\
&\implies&4t(1+z)S=z+2\\
&\implies&S=(z+2)/4t(z+1)
\end{eqnarray*}

Thus, we are led to the conclusions in the statement.
\end{proof}

Getting now to the well-established free laws, first we have the Wigner law:

\index{semicircle law}

\begin{theorem}
The semicircle law $\gamma_t$ of parameter $t>0$, having density
$$\gamma_t=\frac{1}{2\pi t}\sqrt{4t-x^2}\,dx$$
on $[-2\sqrt{t},2\sqrt{t}]$, has the following properties:
\begin{enumerate}
\item $M_{2k}=M_{2k}'=\frac{t^k}{k+1}\binom{2k}{k}$.

\item $M_{2k}''=\frac{1}{k+1}\binom{2k}{k}$.

\item $E=0$, $V=t$, $\gamma=0$, $\kappa=2$.

\item $M(z)=(1-\sqrt{1-4tz^2})/2tz^2$.

\item $G(\xi)=(\xi-\sqrt{\xi^2-4t})/2t$.

\item $R(z)=tz$.

\item $S(z)=1/\sqrt{tz}$.

\item $\kappa_n=0,t,0,0,\ldots$
\end{enumerate}
\end{theorem}

\begin{proof}
This is indeed something that we know well, from chapters 3 and 9, and with the comment that, at the level of the phenomenology and scaling, many things can be said, involving random walks on $\mathbb N$, random matrices, free probability and Lie groups. We will be back to this in the next chapter, with one more piece of phenomenology, again leading to the scaling used above, based on character computations for quantum groups.  
\end{proof}

Next, we have the Marchenko-Pastur law, whose mathematics is quite tricky:

\index{Marchenko-Pastur}

\begin{theorem}
The Marchenko-Pastur law $\pi_t$ of parameter $t>0$, having density
$$\pi_t=\max(1-t,0)\delta_0+\frac{\sqrt{4t-(x-1-t)^2}}{2\pi x}\,dx$$
on $[(1-\sqrt{t})^2,(1+\sqrt{t})^2]$, has the following properties:
\begin{enumerate}
\item $M_k=\sum_{\pi\in NC(k)}t^{|\pi|}$.

\item $E=t$, $V=t$, $\gamma=1/\sqrt{t}$, $\kappa=2+1/t$.

\item $M(z)=(1-t)/2+(1-\sqrt{1-2(1+t)z+(1-t)^2z^2})/2z$.

\item $G(\xi)=(1-t+\xi-\sqrt{(\xi-1-t)^2-4t})/2\xi$.

\item $R(z)=t/(1-z)$.

\item $S(z)=1/(t+z)$.

\item $\kappa_n=t,t,t,t,\ldots$
\end{enumerate}
\end{theorem}

\begin{proof}
These are again things that we know well, from chapters 3 and 9, and with our usual comment, as always when talking about the Marchenko-Pastur laws, that the presence of the atom at 0, in the case $t<1$, is something quite subtle.
\end{proof}

Next, we have the shifted semicircle law, whose mathematics is as follows:

\index{shifted semicircle law}

\begin{theorem}
The shifted semicircle law $\gamma_t^t$ of parameter $t>0$, having density
$$\gamma_t^t=\frac{1}{2\pi t}\sqrt{4t-(x-t)^2}\,dx$$
on $[t-2\sqrt{t},t+2\sqrt{t}]$, has the following properties:
\begin{enumerate}
\item $M_k=\sum_{\pi\in NC_{12}(k)}t^{|\pi|}$.

\item $E=t$, $V=t$, $\gamma=0$, $\kappa=2$.

\item $M(z)=(1-tz-\sqrt{(1-tz)^2-4tz^2})/2tz^2$.

\item $G(\xi)=(\xi-t-\sqrt{(\xi-t)^2-4t})/2t$.

\item $R(z)=t+tz$.

\item $S(z)=(-t+\sqrt{t^2+4tz})/2tz$.

\item $\kappa_n=t,t,0,0,\ldots$
\end{enumerate}
\end{theorem}

\begin{proof}
These are again things that we know well, and with the formulae of the $R$-transform and $S$-transform being something that we know from chapter 7.
\end{proof}

Finally, we have the free Bessel law, whose mathematics is as follows:

\index{free Bessel law}

\begin{theorem}
The free Bessel law $\pi_t^2$ of parameter $t>0$, appearing as
$$\pi_t^2=\pi_{t\varepsilon}\quad,\quad\varepsilon=\frac{\delta_{-1}+\delta_1}{2}$$
on $\mathbb R$, has the following properties:
\begin{enumerate}
\item $M_k=\sum_{\pi\in NC_{even}(k)}t^{|\pi|}$.

\item $M_{2k}=\sum_{b=1}^k\frac{1}{b}\binom{k-1}{b-1}\binom{2k}{b-1}t^b$.

\item $E=0$, $V=t$, $\gamma=0$, $\kappa=2+1/t$.

\item $M=1+(zM)^2(M+t-1)$.

\item $\xi G=1+\xi G^3+(t-1)G^2$.

\item $R(z)=tz/(1-z^2)$.

\item $S(z)=1/\sqrt{z(z+t)}$.

\item $\kappa_n=0,t,0,t,\ldots$
\end{enumerate}
\end{theorem}

\begin{proof}
These are again things that we know, from chapters 8 and 9, with the only needed computation being that of the $S$-transform. But this can be done as follows:
\begin{eqnarray*}
M-1=(zM)^2(M+t-1)
&\implies&\psi=z^2(\psi+1)^2(\psi+t)\\
&\implies&z=\chi^2(z+1)^2(z+t)\\
&\implies&1=zS^2(z+t)\\
&\implies&S=\frac{1}{\sqrt{z(z+t)}}
\end{eqnarray*}

Thus, we are led to the conclusions in the statement.
\end{proof}

\section*{10b. Hankel determinants}

Getting now to a more detailed study of our 12 laws, let us begin with the densities. Normally the density way is the wrong way to be followed, for instance because the densities $g_t,p_t$ obviously have nothing to do with each other, and end of the story. 

\bigskip

This being said, we can formulate however a few interesting observations. To start with, we have the following remarkable fact, concerning 6 of our 12 laws:

\index{beta function}
\index{gamma function}

\begin{theorem}
Up to shifts, rescalings and atoms, $u_t,\alpha_t,\alpha_t',\gamma_t,\gamma_t^t,\pi_t$ are beta laws, with these being the laws of the following type, with $a,b>0$, 
$$\nu_{ab}=\frac{x^{a-1}(1-x)^{b-1}}{B(a,b)}\,dx$$
on $[0,1]$, with the normalization factor being the beta function, given by:
$$B(a,b)=\frac{\Gamma(a)\Gamma(b)}{\Gamma(a+b)}$$
In fact, $u_t,\alpha_t,\alpha_t',\gamma_t,\gamma_t^t,\pi_t$ are the simplest beta laws, appearing at $a,b=1/2,1,3/2$.
\end{theorem}

\begin{proof}
There are several things going on here, the idea being as follows:

\medskip

(1) To start with, and standing as a useful complement to what the statement says, let us recall from calculus that $\Gamma(s)$ with $s>0$ is constructed as follows:
$$\Gamma(s)=\int_0^\infty x^{s-1}e^{-x}\,dx$$

By partial integration we have $\Gamma(s+1)=s\Gamma(s)$, and since $\Gamma(1)=1$, trivially, and $\Gamma(1/2)=\sqrt{\pi}$, obtained via $x=y^2$, we have the following formulae, for $N\in\mathbb N$:
$$\Gamma(N)=(N-1)!\quad,\quad\Gamma\left(N+\frac{1}{2}\right)=\frac{(2N)!!}{2^N}\,\sqrt{\pi}$$

Alternatively, for any integer $n\in\mathbb N$ we have, with $c=\sqrt{2},\sqrt{\pi}$ for $n$ even, odd:
$$\Gamma\left(\frac{n}{2}\right)=\frac{(n-1)!!}{2^{(n-1)/2}}\,c$$

Many other things can be said about the gamma function, the idea being that, by some sort of miracle, $n!=\Gamma(n+1)$ works as a definition, for the factorial of real numbers $n>-1$, for most questions. In what concerns us, what we have in the above will do.

\medskip

(2) Getting now to what our statement says, many things there, and temporarily leaving aside our 6 laws, let us first try to understand why $\nu_{ab}$ is a probability measure. For this purpose, we must establish the following formula, for any $a,b>0$:
$$\int_0^1x^{a-1}(1-x)^{b-1}dx=\frac{\Gamma(a)\Gamma(b)}{\Gamma(a+b)}$$

(3) As a first remark, in the case where the exponents are integers, $a,b\in\mathbb N$, the formula in the statement takes the following form, in terms of usual factorials:
$$\int_0^1x^{a-1}(1-x)^{b-1}dx=\frac{(a-1)!(b-1)!}{(a+b-1)!}$$

But with $x=\cos^2t$ this is exactly the Wallis formula from chapter 1, at odd values $p=2a-1$ and $q=2b-1$ of the exponents there. In fact, the Wallis formula from chapter 1 is precisely the formula in (2), at exponents of type $a,b\in\mathbb N/2$.

\medskip

(4) Getting now to the general case, let us attempt to compute $\Gamma(a)\Gamma(b)$. We have:
\begin{eqnarray*}
\Gamma(a)\Gamma(b)
&=&\int_0^\infty x^{a-1}e^{-x}dx\int_0^\infty y^{b-1}e^{-y}dy\\
&=&\int_0^\infty\int_0^\infty x^{a-1}y^{b-1}e^{-x-y}dxdy
\end{eqnarray*}

(5) Now comes the trick. Let us perform the following change of variables:
$$\begin{cases}x=st\\ y=s(1-t)\end{cases}\iff
\begin{cases}s=x+y\\ t=x/(x+y)\end{cases}$$ 

We have then $dxdy=sdsdt$, coming from the following Jacobian computation:
$$\frac{dxdy}{dsdt}
=\begin{vmatrix}dx/ds&dx/dt\\ \\ dy/ds&dy/dt\end{vmatrix}
=\begin{vmatrix}t&s\\ \\1-t&-s\end{vmatrix}=-s$$

(6) We conclude that our integral over $\mathbb R^2$ above takes the following form:
\begin{eqnarray*}
\Gamma(a)\Gamma(b)
&=&\int_0^1\int_0^\infty(st)^{a-1}(s(1-t))^{b-1}e^{-s}sdsdt\\
&=&\int_0^1s^{a+b-1}e^{-s}ds\int_0^\infty t^{a-1}(1-t)^{b-1}dt\\
&=&\Gamma(a+b)\int_0^\infty t^{a-1}(1-t)^{b-1}dt
\end{eqnarray*}

Thus, we are led to the formula in (2), showing that $\nu_{ab}$ has indeed mass 1.

\medskip

(7) Next, what are the simplest beta laws? And here, with $a,b=1/2,1,3/2$, we get, up to symmetry, and up to the beta function constants, the following laws on $[0,1]$:
$$\frac{1}{\sqrt{x(1-x)}}\quad,\quad\frac{1}{\sqrt{x}}
\quad,\quad\sqrt{\frac{1-x}{x}}\quad,\quad 1,\quad,\quad \sqrt{1-x}
\quad,\quad\sqrt{x(1-x)}$$

But these are, up to shifts, rescalings and atoms, $\alpha_t/\alpha_t'$, then the somewhat obscure law $1/\sqrt{x}$, then $\pi_t$, $u_t$, then another obscure law, $1/\sqrt{1-x}$, and then $\gamma_t/\gamma_t^t$, as claimed.
\end{proof}

Let us record as well the following result, generalizing some previous findings:

\index{beta distribution}

\begin{theorem}
The moments of the beta law of parameters $a,b>0$ are
$$M_k=\frac{a(a+1)\ldots(a+k-1)}{(a+b)(a+b+1)\ldots(a+b+k-1)}$$
and in particular, the mean and variance are given by the following formulae,
$$E=\frac{a}{a+b}\quad,\quad V=\frac{ab}{(a+b)^2(a+b+1)}$$
and the skewness and kurtosis are given by the following formulae,
$$\gamma=\frac{2(b-a)\sqrt{a+b+1}}{(a+b+2)\sqrt{ab}}\quad,\quad
\kappa=3+\frac{6[(a-b)^2(a+b+1)-ab(a+b+2)]}{ab(a+b+2)(a+b+3)}$$
and with this generalizing our previous formulae, regarding $u_t,\alpha_t,\alpha_t',\gamma_t,\gamma_t^t,\pi_t$.
\end{theorem}

\begin{proof}
The moments of the beta law can be indeed computed as follows:
\begin{eqnarray*}
M_k
&=&\frac{1}{B(a,b)}\int_0^1x^{a+k-1}(1-x)^{b-1}dx\\
&=&\frac{B(a+k,b)}{B(a,b)}\\
&=&\frac{\Gamma(a+k)\Gamma(b)}{\Gamma(a+b+k)}\cdot\frac{\Gamma(a+b)}{\Gamma(a)\Gamma(b)}\\
&=&\frac{\Gamma(a+k)}{\Gamma(a)}\cdot\frac{\Gamma(a+b)}{\Gamma(a+b+k)}\\
&=&\frac{a(a+1)\ldots(a+k-1)}{(a+b)(a+b+1)\ldots(a+b+k-1)}
\end{eqnarray*}

As for the formulae of $E,V,\gamma,\kappa$, these follow from this, exercise for you.
\end{proof}

We will be back to the beta laws later in this book. Next, still with our 12 main laws in mind, and still talking densities, we can have a look at Stieltjes inversion:

\begin{theorem}
The Stieltjes inversion formula, namely
$$d\mu (x)=\lim_{s\searrow 0}-\frac{1}{\pi}\,Im\left(G(x+is)\right)\cdot dx$$
applies well to $u_t$, to $\alpha_t,\alpha_t'$, and to $\gamma_t,\pi_t,\gamma_t^t$.
\end{theorem}

\begin{proof}
This is obviously something subjective, but at the point where we are, exploring the mysteries of our 12 laws, a bit of subjectivity won't hurt. Here is the idea:

\medskip

(1) To start with, applying Stieltjes inversion to the classical laws $g_t,p_t,g_t^t,p_t^2,e_t$ is not recommended, because we don't have explicit formulae for the corresponding Cauchy transforms. And with this being something particularly annoying in the case of $g_t$, because when performing a central limit, what you get via Fourier are the moments $M_{2k}=t^k(2k)!!$ of the limiting law, and someone has to tell you to look at $g_t$. Thanks Gauss.

\medskip

(2) On the opposite, as we know well from chapter 3, applying Stieltjes inversion to the laws $\gamma_t,\pi_t$ is certainly recommended, because this is how we came upon these laws, via Stieltjes inversion. As for the other free laws, these are $\alpha_t,\alpha_t',\gamma_t^t$, where Stieltjes inversion applies as well without troubles, and then $\pi^2_t$, where things are quite unclear, I mean can the equation $\xi G=1+\xi G^3+(t-1)G^2$ lead to something, good question.

\medskip

(3) We are therefore left with $u_t$. And here, while it looks a bit ridiculous to frantically search, via Stieltjes inversion, for a measure having moments $M_k=t^k/(k+1)$, well, guided by Newton, Einstein, Grothendieck and the others, good puzzle that we have here, so let's solve it, in order to get smarter. Using $M_k=t^k/(k+1)$, the Cauchy transform is:
$$G(\xi)=\sum_{k\geq0}\frac{t^k\xi^{-k-1}}{k+1}=-\frac{\log(1-t\xi^{-1})}{t}$$

Next, let us set $\xi=x+is$, and look at the following quantity:
$$e^{-tG(x+is)}=1-\frac{t}{x+is}$$

The point now is that with $s\searrow 0$, the right term goes to $1-t/x\in\mathbb R$, so we must have $-tIm(G(x+is))\to k\pi$ with $k\in\mathbb Z$, and more specifically, by looking at what happens to the real part too, with $k=1$. Thus, by Stieltjes inversion, the desired density is:
$$d\mu (x)=\lim_{s\searrow 0}-\frac{1}{\pi}\,Im\left(G(x+is)\right)\cdot dx=\frac{1}{t}$$

Which sounds good, puzzle solved, and did we get any smarter? Well, I would say yes, this adds to our previous findings, stating that $u_t$ is something quite interesting.
\end{proof}

Next, still talking Stieltjes inversion, observe that this does not fully solve the moment problem, because we still have the question of understanding when a sequence of numbers $(M_k)_{k\geq0}$ can be the moments of a measure $\mu$. For instance we certainly must have $M_0=1$, but we must have as well the following inequality, corresponding to $V\geq0$:
$$M_2\geq M_1^2$$

In answer now, we have the following result, complementing Stieltjes inversion:

\index{Hankel determinant}
\index{moment problem}

\begin{theorem}
A sequence of numbers $M_0,M_1,M_2,M_3,\ldots\in\mathbb R$, with $M_0=1$, is the series of moments of a real probability measure $\mu$ precisely when:
$$\begin{vmatrix}M_0\end{vmatrix}\geq0\quad,\quad 
\begin{vmatrix}
M_0&M_1\\
M_1&M_2
\end{vmatrix}\geq0\quad,\quad 
\begin{vmatrix}
M_0&M_1&M_2\\
M_1&M_2&M_3\\
M_2&M_3&M_4
\end{vmatrix}\geq0\quad,\quad 
\ldots$$
That is, the associated Hankel determinants must be all positive.
\end{theorem}

\begin{proof}
This is something a bit more advanced, the idea being as follows:

\medskip

(1) Some numerics first. We know that the first two Hankel determinants are $H_1=1$ and $H_2=V$, both positive. Regarding the third Hankel determinant, this is given by:
$$H_3=\begin{vmatrix}
1&M_1&M_2\\
M_1&M_2&M_3\\
M_2&M_3&M_4
\end{vmatrix}=(M_4-M_2^2)V-(M_3-M_1M_2)^2$$

Which does not look obviously positive, so let us do a check for the binomial law $b_{np}=b_p^{*n}$. And here, after 1 hour of computations, we obtain, as desired:
$$H_3(b_{np})=2n^2(n-1)p^3(1-p)^3\geq0$$

Getting now to $H_4$, things get much worse, and as a wise decision, let us assume that our measure is centered, $M_{2k+1}=0$. In this case $H_4$ is easily computable, given by:
$$H_4=\begin{vmatrix}
1&0&M_2&0\\
0&M_2&0&M_4\\
M_2&0&M_4&0\\
0&M_4&0&M_6
\end{vmatrix}=(M_2M_6-M_4^2)(M_4-M_2^2)$$

Now with input from $g_t$, namely $M_2=t,M_4=3t^2,M_6=15t^3$, we obtain:
$$H_1(g_t)=1\quad,\quad H_2(g_t)=t\quad,\quad H_3(g_t)=2t^3\quad,\quad H_4(g_t)=12t^6$$

Which looks good, all positive numbers, but as an overall conclusion to this, the verification of $H_k\geq0$ by brute-force computation of $H_k$ is not a good method.

\medskip

(2) Fortunately, linear algebra comes to the rescue. Recall indeed from there the positivity theorem of Sylvester, stating that a matrix is positive precisely when its principal minors are all positive. In view of this, the conditions $H_k\geq0$, taken altogether, precisely mean that the corresponding matrices, say $A_k$, are all positive, $A_k\geq0$. 

\medskip

(3) Thus, as an overall conclusion, our theorem states that $(M_k)_{k\geq0}$ are the moments of a measure $\mu$ precisely when we have, for any $k\in\mathbb N$, and any $c\in\mathbb C^N$:
$$\sum_{i,j=1}^nM_{i+j}c_i\bar{c}_j\geq0$$

(4) But in one sense this is elementary, coming from the following computation:
$$\sum_{i,j=1}^nM_{i+j}c_i\bar{c}_j
=\int_\mathbb R\sum_{i,j=1}^nc_i\bar{c}_jx^{i+j}d\mu(x)
=\int_\mathbb R\left|\sum_{i=1}^nc_ix^i\right|^2d\mu(x)$$

(5) As for the other sense, this is something more delicate, requiring some functional analysis study. So, exercise for you, to learn more about all this.
\end{proof}

Getting back now to our 12 main probability measures, we have, quite remarkably, the following beautiful and mysterious result, regarding 10 of them:

\begin{theorem}
The Hankel determinants for the main laws at $t=1$ are as follows:
\begin{enumerate}
\item Normal, $M_{2k}=(2k)!!$, here $H_k=\prod_{i=1}^{k-1}i!$.

\item Poisson, $M_k=|P(k)|$, here $H_k=\prod_{i=1}^{k-1}i!$.

\item Shifted normal, $M_k=|P_{12}(k)|$, here $H_k=\prod_{i=1}^{k-1}i!$.

\item Exponential, $M_k=k!$, here $H_k=\prod_{i=1}^{k-1}i!^2$. 

\item Uniform, $M_k=1/(k+1)$, here $H_k^{-1}=k!\prod_{i=1}^{2k-1}\binom{i}{[i/2]}$.

\item Arcsine, $M_{2k}=\binom{2k}{k}$, here $H_k=2^{k-1}$.

\item Rival arcsine, $M_k=\binom{2k}{k}$, here $H_k=2^{k-1}$.

\item Semicircle, $M_{2k}=\frac{1}{k+1}\binom{2k}{k}$, here $H_k=1$.

\item Marchenko-Pastur, $M_k=\frac{1}{k+1}\binom{2k}{k}$, here $H_k=1$.

\item Shifted semicircle, $M_k=|NC_{12}(k)|$, here $H_k=1$.
\end{enumerate}
\end{theorem}

\begin{proof}
This is something tough, that we will not attempt to prove here. This being said, we can do some verifications, and formulate some comments, as follows:

\medskip

(1) Some numerics first. At $k=2$ we obtain indeed the correct variances:
$$H_2\ =\ 1\ ,\ 1\ ,\ 1\ ,\ 1\ ,\ \frac{1}{12}\ ,\ 2\ ,\ 2\ ,\ 1\ ,\ 1\ ,\ 1$$

At $k=3$ we obtain the following numbers, which are the correct ones too, as you can check by using the formula of $H_3$ from the proof of Theorem 10.19:
$$H_3\ =\ 2\ ,\ 2\ ,\ 2\ ,\ 4\ ,\ \frac{1}{2160}\ ,\ 4\ ,\ 4\ ,\ 1 ,\ 1 ,\ 1$$

At $k=4$ we obtain the following numbers, correct too, good exercise for you:
$$H_4\ =\ 12\ ,\ 12\ ,\ 12\ ,\ 144\ ,\ \frac{1}{24192000}\ ,\ 8\ ,\ 8\ ,\ 1\ ,\ 1\ ,\ 1$$

(2) And so on. In general, however, there are several proofs, all complicated. So, come back here after reading the next section, on orthogonal polynomials, which can help.

\medskip

(3) Finally, regarding the parametric versions $g_t,p_t,g_t^t,e_t,u_t,\alpha_t,\alpha_t',\gamma_t,\pi_t,\gamma_t^t$, in some cases we can get away with homogeneity, while in other cases the problem becomes substantially more complicated. With exercise for you, to learn more about all this.
\end{proof}

Regarding now the missing laws $p^2_t,\pi^2_t$, things for them are quite tricky. For the Bessel law $p^2_t$, the even moments at $t=1$, given by $M_{2k}=|P_{even}(2k)|$, are as follows:
$$1\ ,\ 1\ , \ 4\ ,\ 31\ ,\ 379\ ,\ 6556\ ,\ 150349\ ,\ \ldots$$

The first few computations of Hankel determinants are encouraging, as follows:
$$H_1=1\ ,\quad\frac{H_2}{H_1}=1\ ,\quad\frac{H_3}{H_2}=3\ ,\quad\frac{H_4}{H_3}=15\ ,\quad\frac{H_5}{H_4}=120\ ,\quad\frac{H_6}{H_5}=1260$$

However, a bad surprise awaits us at $k=7$, and then at $k=8$ too, as follows:
$$\frac{H_7}{H_6}=2\cdot 9\cdot5\cdot7\cdot29\quad,\quad 
\frac{H_8}{H_7}=4\cdot 27\cdot5\cdot7\cdot83$$

So, where do these numbers, 29 and then 83, come from? Good question, I guess. And the situation for $\pi^2_t$ is quite similar. The problem comes from the fact that the orthogonal polynomials for $p^2_t,\pi^2_t$ are something quite complicated, and more on this later.

\section*{10c. Orthogonal polynomials}

We would like to have now a discussion on orthogonal polynomials, which are advanced objects, related to the above, and whose values for the main 12 laws will be something spectacular again, in the spirit of Theorem 10.20. Let us start with:

\index{orthogonal polynomials}
\index{Gram-Schmidt}

\begin{definition}
The orthogonal polynomials with respect to $d\mu(x)=f(x)dx$ are polynomials $P_k\in\mathbb R[x]$ of degree $k\in\mathbb N$, which are orthogonal inside $H=L^2(\mathbb R,\mu)$:
$$\int_\mathbb RP_k(x)P_l(x)f(x)dx=0\quad,\quad\forall k\neq l$$
Equivalently, these orthogonal polynomials $\{P_k\}_{k\in\mathbb N}$, which are each unique modulo scalars, appear from the Weierstrass basis $\{x^k\}_{k\in\mathbb N}$, by doing Gram-Schmidt.
\end{definition}

As a first observation, the orthogonal polynomials exist indeed for any real measure $d\mu(x)=f(x)dx$, because we can obtain them from the monomials $x^k$ via Gram-Schmidt, as indicated above. It is possible to be a bit more explicit here, as follows:

\begin{theorem}
The orthogonal polynomials with respect to $\mu$ are given by
$$P_k=c_k\begin{vmatrix}
M_0&M_1&\ldots&M_k\\
M_1&M_2&\ldots&M_{k+1}\\
\vdots&\vdots&&\vdots\\
M_{k-1}&M_k&\ldots&M_{2k-1}\\
1&x&\ldots&x^k
\end{vmatrix}$$
where $M_k=\int_\mathbb Rx^kd\mu(x)$ are the moments of $\mu$, and $c_k\in\mathbb R^*$ can be any numbers.
\end{theorem}

\begin{proof}
Let us first see what happens at small values of $k\in\mathbb N$. At $k=0,1$ our formula holds indeed. At $k=2$ things get more complicated, the formula being:
$$P_2=c_2\begin{vmatrix}
M_0&M_1&M_2\\
M_1&M_2&M_3\\
1&x&x^2
\end{vmatrix}$$

However, no need for big computations here, in order to check the orthogonality, because by using the fact that $x^k$ integrates up to $M_k$, we obtain:
$$<1,P_2>=\int_\mathbb RP_2(x)d\mu(x)=c_2\begin{vmatrix}
M_0&M_1&M_2\\
M_1&M_2&M_3\\
M_0&M_1&M_2
\end{vmatrix}=0$$

Similarly, again by using the fact that $x^k$ integrates up to $M_k$, we have as well:
$$<x,P_2>=\int_\mathbb RxP_2(x)d\mu(x)=c_2\begin{vmatrix}
M_0&M_1&M_2\\
M_1&M_2&M_3\\
M_1&M_2&M_3
\end{vmatrix}=0$$

Thus, result proved at $k=0,1,2$, and the proof in general is similar.
\end{proof}

In practice now, for the simplest measured space $X\subset\mathbb R$, which is the interval $[-1,1]$, with its uniform measure, the orthogonal basis problem can be solved as follows:

\index{Legendre polynomials}
\index{Rodrigues formula}
\index{Legendre equation}

\begin{theorem}
The orthogonal polynomials for $L^2[-1,1]$, subject to
$$\int_{-1}^1P_k(x)P_l(x)\,dx=\delta_{kl}$$
and called Legendre polynomials, satisfy the following differential equation,
$$(1-x^2)P_k''(x)-2xP_k'(x)+k(k+1)P_k(x)=0$$
which is the Legendre equation from physics. Moreover, we have the formula
$$(k+1)P_{k+1}(x)=(2k+1)xP_k(x)-kP_{k-1}(x)$$
called Bonnet recurrence formula, as well as the formula
$$P_k(x)=\frac{1}{2^kk!}\cdot\frac{d^k}{dx^k}\left(1-x^2\right)^k$$
called Rodrigues formula for the Legendre polynomials.
\end{theorem}

\begin{proof}
Many things going on here, the idea being as follows:

\medskip

(1) The first assertion is clear, because the Gram-Schmidt procedure applied to the Weierstrass basis $\{x^k\}$ can only lead to a certain family of polynomials $\{P_k\}$, with each $P_k$ being of degree $k$, and also unique, if we assume that it has positive leading coefficient, with this $\pm$ choice being needed, as usual, at each step of Gram-Schmidt.

\medskip

(2) Now thinking about what Gram-Schmidt does, this is certainly something by recurrence. And examining the recurrence leads to the Legendre equation, as stated. As for the Bonnet recurrence formula, the story here is similar.

\medskip

(3) Regarding the Rodrigues formula, by uniqueness no need to try to understand where this formula comes from, and we have two choices here, either by verifying that $\{P_k\}$ is orthonormal, or by verifying the Legendre equation. And both methods work.

\medskip

(4) In relation with the 12 laws, the modified Legendre polynomials $P_k(2x/t-1)$ solve the problem for $u_t$, uniform measure on $[0,t]$. So one done, 11 to go. By the way, the choice $[-1,1]$ in the statement is very standard, coming from physics, the Legendre polynomials defined in this way being able to solve the hydrogen atom. Read about this.
\end{proof}

The above result is just the tip of the iceberg, and as a continuation, we have:

\index{Jacobi polynomials}
\index{Chebycheff polynomials}
\index{Gegenbauer polynomials}

\begin{theorem}
The orthogonal polynomials for $L^2[-1,1]$, with measure
$$d\mu(x)=(1-x)^a(1+x)^bdx$$
called Jacobi polynomials, satisfy as well a degree $2$ equation, namely
$$(1-x^2)P_k''(x)+(b-a-(a+b+2)x)P_k'(x)+k(k+a+b+1)P_k(x)=0$$
as well as an order $2$ recurrence relation, and are given by the following formula:
$$P_k(x)=\frac{(-1)^k}{2^kk!}(1-x)^{-a}(1+x)^{-b}\frac{d^k}{dx^k}\left[(1-x)^a(1+x)^b(1-x^2)^k\right]$$
As main examples, at $a=b=0$ we recover the Legendre polynomials, and at $a,b=\pm\frac{1}{2}$ we recover the four types of Chebycheff polynomials, from trigonometry. 
\end{theorem}

\begin{proof}
Again, many things going on here, the idea being as follows:

\medskip

(1) To start with, in what regards the precise statement, the order 2 recurrence relation mentioned there is something quite complicated, as follows:
\begin{eqnarray*}
&&2k(k+a+b)(2k+a+b-2)P_k(x)\\
&=&(2k+a+b-1)\left[(2k+a+b)(2k+a+b-2)x+a^2-b^2\right]P_{k-1}(x)\\
&-&2(k+a-1)(k+b-1)(2k+a+b)P_{k-2}(x)
\end{eqnarray*}

(2) Regarding now the main particular cases of the Jacobi polynomials, these are the Gegenbauer polynomials, appearing at $a=b$, with as particular cases:

\medskip

-- The Legendre polynomials from Theorem 10.23, appearing at $a=b=0$.

\medskip

-- The Chebycheff polynomials of the first kind, $T_k(\cos t)=\cos(kt)$, at $a=b=-\frac{1}{2}$.

\medskip

-- And of the second kind too, $U_k(\cos t)\sin t=\sin((k+1)t)$, appearing at $a=b=\frac{1}{2}$.

\medskip

(3) Passed Gegenbauer, at $a=-\frac{1}{2},b=\frac{1}{2}$ we have the Chebycheff polynomials of the third kind, $V_k(\cos t)\cos(t/2)=\cos((k+1/2)t)$, and at $a=\frac{1}{2},b=-\frac{1}{2}$ we have the Chebycheff polynomials of the fourth kind, $W_k(\cos t)\sin(t/2)=\sin((k+1/2)t)$.

\medskip

(4) In relation with probability, the Jacobi polynomials solve the orthogonalization problem for the beta distributions from Theorem 10.16, up to some normalizations. In particular, in relation with the 12 laws, $u_t,\alpha_t,\alpha_t',\gamma_t,\gamma_t^t,\pi_t$ done, and 6 more to go.

\medskip

(5) Finally, regarding the proof, the statement itself appears as a generalization of Theorem 10.23, which corresponds to the particular case $a=b=0$, and the proof is quite similar. We will leave learning more about all this as an exercise.
\end{proof}

Getting now to other spaces $X\subset\mathbb R$, we first have the following result, dealing with $e_t$, which complements well Theorem 10.23, for the needs of basic quantum mechanics:

\index{Laguerre polynomials}

\begin{theorem}
The orthogonal polynomials for $L^2[0,\infty)$, with scalar product
$$<f,g>=\int_0^\infty f(x)g(x)e^{-x}\,dx$$
are the Laguerre polynomials $\{P_k\}$, satisfying the following differential equation,
$$xP_k''(x)+(1-x)P_k'(x)+kP_k(x)=0$$
as well as the following order $2$ recurrence relation,
$$(k+1)P_{k+1}(x)=(2k+1-x)P_k(x)-kP_{k-1}(x)$$
and which are given by the following formula,
$$P_k(x)=\frac{e^x}{k!}\cdot\frac{d^k}{dx^k}\left(e^{-x}x^k\right)$$
called Rodrigues formula for the Laguerre polynomials.
\end{theorem}

\begin{proof}
The story here is very similar to that of the Legendre and Jacobi polynomials, with exercise for you to prove this, and learn more about all this. Finally, for the story to be complete, no discussion about the Laguerre polynomials would be complete without a word about their use, in quantum mechanics. So, have a look at a basic quantum physics book, see how these beasts produce the wave functions of hydrogen.
\end{proof}

Finally, regarding the space $X=\mathbb R$ itself, we have here the following result:

\index{Hermite polynomials}

\begin{theorem}
The orthogonal polynomials for $L^2(\mathbb R)$, with scalar product
$$<f,g>=\int_0^\infty f(x)g(x)e^{-x^2}\,dx$$
are the Hermite polynomials $\{P_k\}$, satisfying the following differential equation,
$$P_k''(x)-2xP_k'(x)+P_k(x)=0$$
as well as the following order $2$ recurrence relation,
$$P_{k+1}(x)=2xP_k(x)-2kP_{k-1}(x)$$
and which are given by the following formula,
$$P_k(x)=(-1)^k e^{x^2}\cdot\frac{d^k}{dx^k}\big(e^{-x^2}\big)$$
called Rodrigues formula for the Hermite polynomials.
\end{theorem}

\begin{proof}
As before, the story here is quite similar to that of the Legendre and other orthogonal polynomials, and exercise for you to learn a bit about all this. Finally, do not forget to work out some rescalings too, with the density $e^{-x^2}$ replaced by the more familiar $e^{-x^2/2t}$, corresponding to the normal law $g_t$, and then by the density of $g_t^t$ as well.
\end{proof}

And with this, good news, end of the story with the orthogonal polynomials, at least at the  introductory level, and this due to the following fact, which is something quite technical, and that we will not attempt to prove, or even explain in detail here:

\begin{fact}
From an abstract point of view, coming from degree $2$ equations, and Rodrigues formulae for the solutions, there are only three types of ``classical'' orthogonal polynomials, namely the Jacobi, Laguerre and Hermite ones, discussed above.
\end{fact}

And isn't this amazing, on par with the combinatorial beauties from Theorem 10.20. Simply put, up to some manipulations on our 12 laws, namely merging all the beta ones, and forgetting about $p_t,p_t^2,\pi_t^2$, that is all. And with this, coming by theorem.

\bigskip

The continuation of the story, however, is more complicated, as follows:

\index{Charlier polynomials}
\index{Askey scheme}

\begin{fact}
The orthogonal polynomials for $p_t$ are the Charlier polynomials, which are something quite complicated. The unification with Jacobi is done via the Askey polynomials, advanced level. As for $p^2_t,\pi^2_t$, things are quite complicated here.
\end{fact}

And we will end with this. Take it easy I guess, we learned many things regarding the 12 laws, but this is not the only way of viewing things. Remember indeed that easiness allows us to reconcile $g_t,p_t,g_t^t,p^2_t,\gamma_t,\pi_t,\gamma_t^t,\pi_t^2$, with $e_t,u_t,\alpha_t,\alpha_t'$ being left aside. So, the truth lies somewhere in between easiness, and the Askey scheme. Good to know.

\section*{10d. Gram determinants}

Following \cite{bcu}, we would like to discuss now the Gram determinant question, which is quite central, in all this. Anticipating a bit with quantum groups, we have:

\begin{question}
In the context of the main orthogonal easy quantum groups
$$\xymatrix@R=17pt@C=17pt{
&H_N^+\ar[rr]&&O_N^+\\
S_N^+\ar[rr]\ar[ur]&&B_N^+\ar[ur]\\
&H_N\ar[rr]\ar[uu]&&O_N\ar[uu]\\
S_N\ar[uu]\ar[ur]\ar[rr]&&B_N\ar[uu]\ar[ur]
}$$
how do the Gram determinants behave, on the vertical?
\end{question}

To be more precise, what we call here easy quantum groups are formal objects associated to the main 8 categories of partitions, from chapter 9, namely:
$$\xymatrix@R=18pt@C=2pt{
&NC_{even}\ar[dl]\ar[dd]&&NC_2\ar[dl]\ar[ll]\ar[dd]\\
NC\ar[dd]&&NC_{12}\ar[dd]\ar[ll]\\
&P_{even}\ar[dl]&&P_2\ar[dl]\ar[ll]\\
P&&P_{12}\ar[ll]
}$$

More on quantum groups soon, in chapter 11, and in the meantime, consider Question 10.29 as referring to the above 8 categories. So, our question is well-defined.

\bigskip

In answer now, we first need to know how to compute the above determinants. In the classical discrete case, the answer is something quite simple, as follows:

\begin{theorem}
For the groups $G=S_N,H_N$ we have the Lindst\"om formula
$$\det(G_{kN})=\prod_{\pi\in D(k)}\frac{N!}{(N-|\pi|)!}$$
with $D=P,P_{even}$, and with $|.|$ being the number of blocks.
\end{theorem}

\begin{proof}
This is something that we already know for $G=S_N$, from chapter 4, and the proof for $G=H_N$ is similar, based in the fact that the corresponding categories of partitions $D=P,P_{even}$ have the property of forming semilattices under $\vee$. Consider indeed the following matrix, obtained by making determinant-preserving operations:
$$G_{kN}'(\pi,\sigma)=\sum_{\pi\leq\tau}\mu(\pi,\tau)N^{|\tau\vee\sigma|}$$

It follows then from the M\"obius inversion formula that we have:
$$G_{kN}'(\pi,\sigma)=
\begin{cases}
N(N-1)\ldots(N-|\sigma|+1)&{\rm if}\ \pi\leq\sigma\\
0&{\rm otherwise}
\end{cases}$$

Thus the matrix is upper triangular, and by computing the product on the diagonal we obtain the formula in the statement, exactly as in the $G=S_N$ case before.
\end{proof}

Next, let us discuss the case of the orthogonal group $O_N$. Here things are more complicated, the combinatorics being that of the Young diagrams. So, let us denote by $|.|$ the number of boxes, and use quantity $f^\lambda$, which gives the number of standard Young tableaux of shape $\lambda$. With these conventions, we have the following result:

\index{Young tableaux}

\begin{theorem}
The determinant of the Gram matrix of $O_N$ is given by
$$\det(G_{kN})=\prod_{|\lambda|=k/2}f_N(\lambda)^{f^{2\lambda}}$$
where the quantities on the right are $f_N(\lambda)=\prod_{(i,j)\in\lambda}(N+2j-i-1)$.
\end{theorem}

\begin{proof}
For the group $O_N$ the Gram matrix is diagonalizable, as follows:
$$G_{kN}=\sum_{|\lambda|=k/2}f_N(\lambda)P_{2\lambda}$$

To be more precise, here $1=\sum P_{2\lambda}$ is the standard partition of unity associated to the Young diagrams having $k/2$ boxes, and the coefficients $f_N(\lambda)$ are those in the statement. Now since we have $Tr(P_{2\lambda})=f^{2\lambda}$, this gives the formula in the statement.
\end{proof}

As a variation of the above result, for the bistochastic group $B_N$, we have:

\index{bistochastic group}

\begin{theorem}
For the bistochastic group $B_N$ we have
$$\det(G_{kN})=N^{a_k}\prod_{|\lambda|\leq k/2}f_N(\lambda)^{\binom{k}{2|\lambda|}f^{2\lambda}}$$
where $a_k=\sum_{\pi\in P_{12}(k)}(2|\pi|-k)$, and $f_N(\lambda)=\prod_{(i,j)\in\lambda}(N+2j-i-2)$.
\end{theorem}

\begin{proof}
We can use here the standard fact that we have an isomorphism as follows, where $u,v$ are the fundamental representations of $B_N,O_{N-1}$: 
$$B_N\simeq O_{N-1}\quad,\quad u=v+1$$

In terms of the corresponding Peter-Weyl representations, this gives:
$$Fix(u^{\otimes k})
=Fix\left((v+1)^{\otimes k}\right)
=Fix\left(\sum_{r=0}^k\binom{k}{r}v^{\otimes r}\right)$$

Now if we denote by $\det',f'$ the objects in Theorem 10.31, we obtain:
$$\det(G_{kN})
=N^{a_k}\prod_{r=1}^k{\rm det}'(G_{r,N-1})^{\binom{k}{r}}
=N^{a_k}\prod_{r=1}^k\left(\prod_{|\lambda|=r/2}f'_{N-1}(\lambda)^{f^{2\lambda}}\right)^{\binom{k}{r}}$$

Thus, we are led to the formula in the statement.
\end{proof}

Getting now to the free case, the simplest object here is $O_N^+$, coming from the category of noncrossing pairings $NC_2$. The associated Gram determinant, known as ``meander determinant'', was computed by Di Francesco in \cite{dif}, the result being as follows:

\index{meander determinant}
\index{Gram determinant}
\index{Di Francesco formula}

\begin{theorem}
The determinant of the Gram matrix for $O_N^+$ is given by
$$\det(G_{kN})=\prod_{r=1}^{[k/2]}P_r(N)^{d_{k/2,r}}$$
where $P_r$ are the Chebycheff polynomials, given by
$$P_0=1\quad,\quad 
P_1=X\quad,\quad 
P_{r+1}=XP_r-P_{r-1}$$
and $d_{kr}=f_{kr}-f_{k,r+1}$, with $f_{kr}$ being the following numbers, depending on $k,r\in\mathbb Z$,
$$f_{kr}=\binom{2k}{k-r}-\binom{2k}{k-r-1}$$
with the conventions $f_{kr}=0$ for $k\notin\mathbb Z$, and $\binom{p}{q}=0$ for $q<0$. 
\end{theorem}

\begin{proof}
This is indeed something quite standard, but long and technical. For more on this we refer to \cite{dif}, and to the more recent paper \cite{bcu} too, containing a short proof of this formula, based on advanced planar algebra technology of Jones \cite{jo4}.
\end{proof}

Regarding now $S_N^+$, coming from the category of all noncrossing partitions $NC$, we have here the following formula, also established by Di Francesco in \cite{dif}:

\index{meander determinant}
\index{Gram determinant}
\index{Chebycheff polynomials}
\index{Di Francesco formula}

\begin{theorem}
The determinant of the Gram matrix for $S_N^+$ is given by
$$\det(G_{kN})=(\sqrt{N})^{a_k}\prod_{r=1}^kP_r(\sqrt{N})^{d_{kr}}$$
where $P_r$, $d_{kr}$ are as before, and where $a_k=\sum_{\pi\in NC(k)}(2|\pi|-k)$.
\end{theorem}

\begin{proof}
This comes from Theorem 10.33, via the following standard bijection, obtained via fattening the pairings, and shrinking the partitions:
$$NC(k)\simeq NC_2(2k)$$

Indeed, if we denote by $G'$ the Gram matrix for $O_N^+$, we have the following formula, coming from the above bijection, with $D_{kN}=diag(N^{|\widetilde\pi|/2-k/4})$:
$$G_{kN}=D_{kN}G'_{2k,\sqrt{N}}D_{kN}$$

But with this formula in hand, the result follows from Theorem 10.33.
\end{proof}

As yet another version of Theorem 10.33, we can formulate as well:

\begin{theorem}
For the quantum group $B_N^+$ we have
$$\det(G_{kN})=N^{a_k}\prod_{r=1}^{[k/2]}P_r(N-1)^{\sum_{l=1}^{[k/2]}\binom{k}{2l}d_{lr}}$$
with $P_r$ and $d_{kr}$ being as before, and with $a_k=\sum_{\pi\in NC_{12}(k)}(2|\pi|-k)$.
\end{theorem}

\begin{proof}
The passage $O_N^+\to B_N^+$ is quite similar to the passage $O_N\to B_N$, and by using prime exponents for the various $O_N^+$-related objects, we obtain:
$$\det(G_{kN})
=N^{a_k}\prod_{l=1}^{[k/2]}{\rm det}'(G_{2l,N-1})^{\binom{k}{2l}}
=N^{a_k}\prod_{l=1}^{[k/2]}\left(\prod_{r=1}^lP_r(N-1)^{d_{lr}}\right)^{\binom{k}{2l}}$$

Together with Theorem 10.33, this gives the formula in the statement.
\end{proof}

Finally, in what regards the quantum group $H_N^+$, the result here is as follows:

\begin{theorem}
For the quantum group $H_N^+$ we have the formula
$$\det(G_{kN})=(\sqrt{N})^{a_k}\prod_{r=1}^{[k/2]}P_r(\sqrt{N})^{2d_{k/2,r}'}$$
with $d_{sr}'=f_{sr}'-f_{s,r+1}'$, where 
$f_{sr}'=\binom{3s}{s-r}-\binom{3s}{s-r-1}$ for $s\in\mathbb Z$, $f_{sr}'=0$ for $s\notin\mathbb Z$.
\end{theorem}

\begin{proof}
This can be viewed as yet another variation of Theorem 10.33, involving this time some colors on the strings, and for details here we refer to \cite{dif} and \cite{bcu}.
\end{proof}

Getting now to our original Question 10.29, we have formulae for all determinants there, and the remaining puzzle is that of putting everything together. One of the difficulties here lies in understanding the orthogonal polynomials for $H_N,H_N^+$. See \cite{bcu}.

\section*{10e. Exercises}

Welcome I guess to calculus, if you are an algebraist, and as exercises, we have:

\begin{exercise}
Carefully compute all the skewnesses $\gamma$ and kurtoses $\kappa$.
\end{exercise}

\begin{exercise}
Have a look at the various transforms which were not computed.
\end{exercise}

\begin{exercise}
Have a look at the Hankel determinants for $p^2_1$ and $\pi^2_1$.
\end{exercise}

\begin{exercise}
Have a look at Hankel determinants, in the parametric case.
\end{exercise}

\begin{exercise}
Learn about Charlier polynomials, and the Askey scheme.
\end{exercise}

\begin{exercise}
Have a look at the orthogonal polynomials for $p^2_t$ and $\pi^2_t$.
\end{exercise}

\begin{exercise}
Learn the proof of the Di Francesco determinant formula.
\end{exercise}

\begin{exercise}
Write the Gram determinants in terms of orthogonal polynomials.
\end{exercise}

As bonus exercise, have a look at free quantum groups. But, more on this soon.

\chapter{Quantum groups}

\section*{11a. Quantum groups}

Time to introduce the quantum groups, which already appeared here and there, in the last few chapters, in more or less disguised form. Our official motivation will come from the classical/free correspondence in probability, the point being as follows:

\begin{fact}
While for real measures the correspondence $m_t\leftrightarrow\mu_t$ between classical and free can be understood via transforms and cumulants, according to
$$\log F_{m_t}(-iz)=\sum_{n=1}^\infty c_n\cdot\frac{z^n}{n!}
\quad \longleftrightarrow\quad zR_{\mu_t}(z)=\sum_{n=1}^\infty c_n\cdot z^n$$
for complex measures and abstract distributions this is something quite complicated, both algebrically and analytically. The best is to first understand the easy correspondence
$$M_k(m_t)=\sum_{\pi\in D(k)}t^{|\pi|}\quad \longleftrightarrow\quad M_k(\mu_t)=\sum_{\pi\in D'(k)}t^{|\pi|}$$
with $D\subset P,D'\subset NC$ being related by $D'=D\cap NC$ and $D=<D',\slash\hskip-2.1mm\backslash>$, and with this being itself best understood at the quantum group level, as a correspondence $G\leftrightarrow G^+$.
\end{fact}

Well, not the easiest to decipher mathematical statement ever made, sure I agree with this, and going now more in detail through the above, here is the point:

\bigskip

(1) To start with, for the real measures the correspondence $m_t\leftrightarrow\mu_t$ between classical and free can be understood both analytically via transforms, and combinatorially via cumulants, as explained above. Which is very nice, and we have been doing this, analysis and combinatorics for $m_t\leftrightarrow\mu_t$ in the real case, successfully, in the last 50 pages.

\bigskip

(2) In the complex case the simplest instance of the correspondence is $G_t\leftrightarrow\Gamma_t$, and with $\Gamma_t$ actually exiting the class of complex measures, being an abstract distribution. So, you can see right away the difficulty, shall you choose analysis and transforms, or combinatorics and cumulants, axiomatizing $m_t\leftrightarrow\mu_t$ will be a difficult task.

\bigskip

(3) The solution to this comes from easiness, as indicated, which is of course something quite modest, axiomatizing $G_t\leftrightarrow\Gamma_t$, $G_t^t\leftrightarrow\Gamma_t^t$,  $p^s_t\leftrightarrow\pi^s_t$, and that is pretty much it, but which remains something not bad, at least we know one thing. Afterwards, once easiness mastered, you can of course go back to (2), and fight with the difficulties there.

\bigskip

(4) And with all this, coming with an extra twist, in relation with quantum groups. Indeed, with easiness being meant to be the common ground for both the algebraic and analytic theory of probabilistic liberation, this easiness should be something half-algebraic, half-analytic. And, it is quantum groups that bring the much needed analytic touch.

\bigskip

But probably too much talking, let us get to work, and we'll understand later. Following Woronowicz \cite{wo1}, we can introduce the compact quantum groups as follows:

\begin{definition}
A Woronowicz algebra is a $C^*$-algebra $A$, given with a unitary matrix $u\in M_N(A)$ whose coefficients generate $A$, such that the formulae
$$\Delta(u_{ij})=\sum_ku_{ik}\otimes u_{kj}\quad,\quad
\varepsilon(u_{ij})=\delta_{ij}\quad,\quad 
S(u_{ij})=u_{ji}^*$$
define morphisms of $C^*$-algebras $\Delta:A\to A\otimes A$, $\varepsilon:A\to\mathbb C$ and $S:A\to A^{opp}$, called comultiplication, counit and antipode. 
\end{definition}

To be more precise, here $\otimes$ can be any $C^*$-algebra tensor product, and $A^{opp}$ is the opposite algebra, having multiplication $a\cdot b=ba$. More on this, axioms, in a moment.

\bigskip

So, let us see how this works. With the convention that $A$ is cocommutative when $\Sigma\Delta=\Delta$, where $\Sigma(a\otimes b)=b\otimes a$ is the flip, we have the following key result:

\begin{proposition}
The following are Woronowicz algebras:
\begin{enumerate}
\item $C(G)$, with $G\subset U_N$ compact Lie group, with structural maps as follows:
$$\Delta(f)=\big[(g,h)\to f(gh)\big]\quad,\quad 
\varepsilon(f)=f(1)\quad,\quad
S(f)=\big[g\to f(g^{-1})\big]$$

\item $C^*(\Gamma)$, with $F_N\to\Gamma$ finitely generated group, with structural maps as follows:
$$\Delta(g)=g\otimes g\quad,\quad 
\varepsilon(g)=1\quad,\quad
S(g)=g^{-1}$$
\end{enumerate}
Moreover, we obtain in this way all the commutative/cocommutative algebras.
\end{proposition}

\begin{proof}
For the first assertion we can use the matrix $u=(u_{ij})$ formed by the standard matrix coordinates of $G$, which is by definition given by:
$$g=\begin{pmatrix}
u_{11}(g)&\ldots&u_{1N}(g)\\
\vdots&&\vdots\\
u_{N1}(g)&\ldots&u_{NN}(g)
\end{pmatrix}$$

Indeed, all axioms are clearly satisfied. Next, for the second assertion we can use the diagonal matrix formed by the generators, with again the axioms being satisfied, and with the remark here that, for things to work fine, the target of $S$ must be indeed $A^{opp}$:
$$u=\begin{pmatrix}
g_1&&0\\
&\ddots&\\
0&&g_N
\end{pmatrix}$$

Finally, regarding the last assertion, in the commutative case this follows from the Gelfand theorem, and in the cocommutative case, we will be back to this.
\end{proof}

In order to get now to quantum groups, we will need as well:

\begin{proposition}
Assuming that $G\subset U_N$ is abelian, we have an identification of Woronowicz algebras $C(G)=C^*(\Gamma)$, with $\Gamma$ being the Pontrjagin dual of $G$:
$$\Gamma=\big\{\chi:G\to\mathbb T\big\}$$
Conversely, assuming that $F_N\to\Gamma$ is abelian, we have an identification of Woronowicz algebras $C^*(\Gamma)=C(G)$, with $G$ being the Pontrjagin dual of $\Gamma$:
$$G=\big\{\chi:\Gamma\to\mathbb T\big\}$$
Thus, the Woronowicz algebras which are both commutative and cocommutative are exactly those of type $A=C(G)=C^*(\Gamma)$, with $G,\Gamma$ being abelian, in Pontrjagin duality.
\end{proposition}

\begin{proof}
This follows from the Gelfand theorem applied to $C^*(\Gamma)$, and from the fact that the characters of a group algebra come from the characters of the group.
\end{proof}

The above results are quite interesting, and based on them, and on Gelfand duality in general, we can formulate the following definition, complementing Definition 11.2:

\begin{definition}[continuation]
Given a Woronowicz algebra $(A,u)$, we write it as
$$A=C(G)=C^*(\Gamma)$$
and call $G$ a compact quantum Lie group, and $\Gamma$ a finitely generated discrete quantum group. We say that these quantum groups are dual to each other, and write: 
$$G=\widehat{\Gamma}\quad,\quad\Gamma=\widehat{G}$$
Also, we agree to identify two Woronowicz algebras, $(A,u)=(B,v)$, when we have an isomorphism of $*$-algebras $<u_{ij}>\simeq<v_{ij}>$ mapping $u_{ij}\to v_{ij}$.
\end{definition}

To be more precise, the convention $A=C(G)$ comes from Proposition 11.3 and Gelfand, the convention $A=C^*(\Gamma)$ comes from Proposition 11.3, and the resulting duality between $G,\Gamma$ generalizes the usual Pontrjagin duality, from Proposition 11.4. 

\bigskip

As for the convention at the end, which by the way clarifies our $\otimes$ comment made after Definition 11.2, this is something more technical, coming from the fact that a non-amenable discrete group $\Gamma$ has several possible completions of $\mathbb C[\Gamma]$, that we must identify, in order to have uniqueness of $\Gamma$, as a discrete quantum group. More on this later.

\bigskip

Let us discuss now some tools for studying our beasts. We first have:

\begin{proposition}
Let $(A,u)$ be a Woronowicz algebra.
\begin{enumerate} 
\item $\Delta,\varepsilon$ satisfy the usual axioms for a comultiplication and a counit, namely:
$$(\Delta\otimes id)\Delta=(id\otimes \Delta)\Delta$$
$$(\varepsilon\otimes id)\Delta=(id\otimes\varepsilon)\Delta=id$$

\item $S$ satisfies the antipode axiom, on the $*$-algebra generated by entries of $u$: 
$$m(S\otimes id)\Delta=m(id\otimes S)\Delta=\varepsilon(.)1$$

\item In addition, the square of the antipode is the identity, $S^2=id$.
\end{enumerate}
\end{proposition}

\begin{proof}
As a first observation, the result holds in the commutative case, $A=C(G)$ with $G\subset U_N$. Indeed, here we know from Proposition 11.3 that $\Delta,\varepsilon,S$ appear as functional analytic transposes of the multiplication, unit and inverse maps $m,u,i$:
$$\Delta=m^t\quad,\quad 
\varepsilon=u^t\quad,\quad 
S=i^t$$

Thus, the various conditions in the statement on $\Delta,\varepsilon,S$ simply come from the group axioms satisfied by $m,u,i$. In general now, the first axiom follows from:
$$(\Delta\otimes id)\Delta(u_{ij})=(id\otimes \Delta)\Delta(u_{ij})=\sum_{kl}u_{ik}\otimes u_{kl}\otimes u_{lj}$$

As for the other axioms, the verifications here are similar.
\end{proof}

In order to reach to more advanced results, the idea will be that of doing representation theory. Following Woronowicz \cite{wo1}, let us start with the following definition:

\begin{definition}
Given $(A,u)$, we call corepresentation of it any unitary matrix $v\in M_n(\mathcal A)$, with $\mathcal A=<u_{ij}>$, satisfying the same conditions as $u$, namely:
$$\Delta(v_{ij})=\sum_kv_{ik}\otimes v_{kj}\quad,\quad
\varepsilon(v_{ij})=\delta_{ij}\quad,\quad
S(v_{ij})=v_{ji}^*$$
We also say that $v$ is a representation of the underlying compact quantum group $G$.
\end{definition}

In the commutative case, $A=C(G)$ with $G\subset U_N$, we obtain in this way the finite dimensional unitary smooth representations $v:G\to U_n$, via the following formula:
$$v(g)=\begin{pmatrix}
v_{11}(g)&\ldots&v_{1n}(g)\\
\vdots&&\vdots\\
v_{n1}(g)&\ldots&v_{nn}(g)
\end{pmatrix}$$

With these conventions, we have the following fundamental result, from \cite{wo1}:

\begin{theorem}
Any Woronowicz algebra has a unique Haar integration, 
$$\left(\int_G\otimes\, id\right)\Delta=\left(id\otimes\int_G\right)\Delta=\int_G(.)1$$
and for any corepresentation $v\in M_n(\mathbb C)\otimes A$ we have the formula
$$\left(id\otimes\int_G\right)v=P$$
where $P$ is the orthogonal projection onto $Fix(v)=\{\xi\in\mathbb C^n|v\xi=\xi\}$.
\end{theorem}

\begin{proof}
This is something that we know well in the classical case, for the algebras of type $A=C(G)$ with $G\subset U_N$, from chapter 4, and the proof in general is similar.
\end{proof}

Next, we can develop the Peter-Weyl theory for the corepresentations of $A$. Consider the dense subalgebra $\mathcal A=<u_{ij}>$, and endow it with the scalar product $<a,b>=\int_Gab^*$. With this convention, we have the following key result, also from \cite{wo1}:

\begin{theorem}
We have the following Peter-Weyl type results:
\begin{enumerate}
\item Any corepresentation decomposes as a sum of irreducible corepresentations.

\item Each irreducible corepresentation appears inside a certain $u^{\otimes k}$.

\item $\mathcal A=\bigoplus_{v\in Irr(A)}M_{\dim(v)}(\mathbb C)$, the summands being pairwise orthogonal.

\item The characters of irreducible corepresentations form an orthonormal system.
\end{enumerate}
\end{theorem}

\begin{proof}
This is something that we know well from chapter 4, when $G\subset U_N$ is a compact group. In general, the proof is quite similar, by using Theorem 11.8.
\end{proof}

Finally, again by following \cite{wo1}, and further clarifying the axiomatics, we have:

\begin{theorem}
Let $A_{full}$ be the enveloping $C^*$-algebra of $\mathcal A$, and $A_{red}$ be the quotient of $A$ by the null ideal of the Haar integration. The following are then equivalent:
\begin{enumerate}
\item The Haar functional of $A_{full}$ is faithful.

\item The projection map $A_{full}\to A_{red}$ is an isomorphism.

\item The counit map $\varepsilon:A_{full}\to\mathbb C$ factorizes through $A_{red}$.

\item We have $N\in\sigma(Re(\chi_u))$, the spectrum being taken inside $A_{red}$.
\end{enumerate}
If this is the case, we say that the underlying discrete quantum group $\Gamma$ is amenable.
\end{theorem}

\begin{proof}
This is well-known in the group dual case, $A=C^*(\Gamma)$, with $\Gamma$ being a usual discrete group. In general, the result follows by adapting the group dual case proof:

\medskip

$(1)\iff(2)$ This simply follows from the fact that the GNS construction for the algebra $A_{full}$ with respect to the Haar functional produces the algebra $A_{red}$.

\medskip

$(2)\iff(3)$ Here $\implies$ is trivial, and conversely, a counit $\varepsilon:A_{red}\to\mathbb C$ produces an isomorphism $\Phi:A_{red}\to A_{full}$, by slicing the map $\widetilde{\Delta}:A_{red}\to A_{red}\otimes A_{full}$.

\medskip

$(3)\iff(4)$ Here $\implies$ is clear, coming from $\varepsilon(N-Re(\chi (u)))=0$, and the converse can be proved by doing some functional analysis. See \cite{wo1}.
\end{proof}

\section*{11b. Rotations, reflections} 

The compact quantum groups include the compact Lie groups, $G\subset U_N$, and the abstract duals $G=\widehat{\Gamma}$ of the finitely generated groups $F_N\to\Gamma$. Following Wang \cite{wan}, let us discuss now a number of truly ``new'' quantum groups. We first have:

\index{free orthogonal group}
\index{free unitary group}

\begin{theorem}
The following universal algebras are Woronowicz algebras,
$$C(O_N^+)=C^*\left((u_{ij})_{i,j=1,\ldots,N}\Big|u=\bar{u},u^t=u^{-1}\right)$$
$$C(U_N^+)=C^*\left((u_{ij})_{i,j=1,\ldots,N}\Big|u^*=u^{-1},u^t=\bar{u}^{-1}\right)$$
so the underlying quantum spaces $O_N^+,U_N^+$ are compact quantum groups.
\end{theorem}

\begin{proof}
This comes from the elementary fact that if a matrix $u=(u_{ij})$ is orthogonal or biunitary, then so must be the following matrices:
$$(u^\Delta)_{ij}=\sum_ku_{ik}\otimes u_{kj}\quad,\quad 
(u^\varepsilon)_{ij}=\delta_{ij}\quad,\quad 
(u^S)_{ij}=u_{ji}^*$$

Thus we can define $\Delta,\varepsilon,S$ by using the universal property of $C(O_N^+)$, $C(U_N^+)$.
\end{proof}

The basic properties of these beasts $O_N^+,U_N^+$ can be summarized as follows:

\begin{theorem}
The quantum groups $O_N^+,U_N^+$ have the following properties:
\begin{enumerate}
\item The closed subgroups $G\subset U_N^+$ are exactly the $N\times N$ compact quantum groups. As for the closed subgroups $G\subset O_N^+$, these are those satisfying $u=\bar{u}$.

\item We have liberation embeddings $O_N\subset O_N^+$ and $U_N\subset U_N^+$, obtained by dividing the algebras $C(O_N^+),C(U_N^+)$ by their respective commutator ideals.

\item We have as well embeddings $\widehat{L}_N\subset O_N^+$ and $\widehat{F}_N\subset U_N^+$, where $L_N$ is the free product of $N$ copies of $\mathbb Z_2$, and where $F_N$ is the free group on $N$ generators.
\end{enumerate}
\end{theorem}

\begin{proof}
All these assertions are elementary, as follows:

\medskip

(1) This is clear from definitions, with the remark that, in the context of Definition 11.2, the formula $S(u_{ij})=u_{ji}^*$ shows that the matrix $\bar{u}$ must be unitary too. 

\medskip

(2) This follows from the Gelfand theorem, which shows that we have presentation results for the algebras $C(O_N),C(U_N)$, similar to those in Theorem 11.11, but with the commutativity between the standard coordinates and their adjoints added:
\begin{eqnarray*}
C(O_N)&=&C^*_{comm}\left((u_{ij})_{i,j=1,\ldots,N}\Big|u=\bar{u},u^t=u^{-1}\right)\\
C(U_N)&=&C^*_{comm}\left((u_{ij})_{i,j=1,\ldots,N}\Big|u^*=u^{-1},u^t=\bar{u}^{-1}\right)
\end{eqnarray*}

(3) This follows indeed from (1) and from Theorem 11.11, with the remark that with $u=diag(g_1,\ldots,g_N)$, the condition $u=\bar{u}$ is equivalent to $g_i^2=1$, for any $i$.
\end{proof}

Next, following \cite{bsp}, we have as well bistochastic quantum groups, as follows:

\begin{theorem}
We have bistochastic free quantum groups as follows,
$$B_N^+\subset O_N^+\quad,\quad C_N^+\subset U_N^+$$
obtained via the fact that $u=(u_{ij})$ must have sums $1$ on each row and column.
\end{theorem}

\begin{proof}
This is indeed something self-explanatory, with the quantum group property coming as before for $O_N^+,U_N^+$, by looking at the matrices $u^\Delta,u^\varepsilon,u^S$. 
\end{proof}

Moving on, following Wang \cite{wan}, let us discuss now the construction and basic properties of the quantum permutation group $S_N^+$. The construction is as follows:

\begin{theorem}
The following universal $C^*$-algebra, with magic meaning formed by projections $(p^2=p^*=p)$, summing up to $1$ on each row and each column,
$$C(S_N^+)=C^*\left((u_{ij})_{i,j=1,\ldots,N}\Big|u={\rm magic}\right)$$
is a Woronowicz algebra, and the underlying quantum group $S_N^+$, called quantum permutation group, appears as a liberation of the usual permutation group $S_N$.
\end{theorem}

\begin{proof}
The first assertion comes from the standard fact, say easy exercise for you, that if a matrix $u=(u_{ij})$ is magic, then so must be the following matrices:
$$(u^\Delta)_{ij}=\sum_ku_{ik}\otimes u_{kj}\quad,\quad 
(u^\varepsilon)_{ij}=\delta_{ij}\quad,\quad 
(u^S)_{ij}=u_{ji}^*$$

As for the second assertion, what that precisely says is that we have an inclusion $S_N\subset S_N^+$, with the group $S_N$ being the classical version of $S_N^+$, obtained at the algebra level by dividing $C(S_N^+)$ by its commutator ideal. Which, in practice, means:
$$C(S_N)=C^*_{comm}\left((u_{ij})_{i,j=1,\ldots,N}\Big|u={\rm magic}\right)$$

But this is something coming from Gelfand. Indeed, the algebra on the right must be of the form $A=C(X)$, with $X$ being a certain compact space. Now since we have coordinates $u_{ij}:X\to\mathbb C$, our space must appear as $X\subset M_N(\mathbb C)$, and then the magic condition on the matrix $u=(u_{ij})$ tells us that the matrices $g\in X$ must have 0-1 entries, summing up to 1 on each row and column. Thus $X=S_N$, as desired.
\end{proof}

Our first task now is to make sure that Theorem 11.14 produces indeed a new quantum group, which does not collapse to $S_N$. Still following Wang \cite{wan}, we have:

\begin{theorem}
The embedding $S_N\subset S_N^+$, given at the algebra level by
$$u_{ij}\to\chi\left(\sigma\in S_N\Big|\sigma(j)=i\right)$$
is an isomorphism at $N\leq3$, but not at $N\geq4$, where $S_N^+$ is not classical, nor finite.
\end{theorem} 

\begin{proof}
We can prove this in four steps, as follows:

\medskip

\underline{Case $N=2$}. The fact that $S_2^+$ is indeed classical, and hence collapses to $S_2$, is trivial, because the $2\times2$ magic matrices are as follows, with $p$ being a projection:
$$U=\begin{pmatrix}p&1-p\\1-p&p\end{pmatrix}$$

Thus $C(S_2^+)$ is commutative, and equals its biggest commutative quotient, $C(S_2)$.

\medskip

\underline{Case $N=3$}. It is enough to check that $u_{11},u_{22}$ commute. But this follows from:
\begin{eqnarray*}
u_{11}u_{22}
&=&u_{11}u_{22}(u_{11}+u_{12}+u_{13})\\
&=&u_{11}u_{22}u_{11}+u_{11}u_{22}u_{13}\\
&=&u_{11}u_{22}u_{11}+u_{11}(1-u_{21}-u_{23})u_{13}\\
&=&u_{11}u_{22}u_{11}
\end{eqnarray*}

Indeed, by conjugating, $u_{22}u_{11}=u_{11}u_{22}u_{11}$, so $u_{11}u_{22}=u_{22}u_{11}$, as desired.

\medskip

\underline{Case $N=4$}. Consider the following matrix, with $p,q$ being projections:
$$U=\begin{pmatrix}
p&1-p&0&0\\
1-p&p&0&0\\
0&0&q&1-q\\
0&0&1-q&q
\end{pmatrix}$$ 

This matrix is magic, and we can choose $p,q\in B(H)$ as for the algebra $<p,q>$ to be noncommutative and infinite dimensional. We conclude that $C(S_4^+)$ is noncommutative and infinite dimensional as well, and so $S_4^+$ is non-classical and infinite, as claimed.

\medskip

\underline{Case $N\geq5$}. Here we can use the standard embedding $S_4^+\subset S_N^+$, obtained at the level of the corresponding magic matrices in the following way:
$$u\to\begin{pmatrix}u&0\\ 0&1_{N-4}\end{pmatrix}$$

Indeed, with this in hand, the fact that $S_4^+$ is a non-classical, infinite compact quantum group implies that $S_N^+$ with $N\geq5$ has these two properties as well.
\end{proof}

At a more advanced level, we have the following result, improving Theorem 11.15:

\begin{theorem}
The quantum groups $S_N^+$ have the following properties:
\begin{enumerate}
\item We have $S_N^+\,\hat{*}\,S_M^+\subset S_{N+M}^+$, for any $N,M$.

\item In particular, we have an embedding $\widehat{D_\infty}\subset S_4^+$. 

\item If $\mathbb Z_{N_1}*\ldots*\mathbb Z_{N_k}\to\Gamma$, with $N=\sum N_i$, then $\widehat{\Gamma}\subset S_N^+$.

\item The quantum groups $S_N^+$ with $N\geq5$ are not coamenable.
\end{enumerate}
\end{theorem}

\begin{proof}
This is something quite standard, the idea being as follows:

\medskip

(1) If we denote by $u,v$ the fundamental corepresentations of $C(S_N^+),C(S_M^+)$, the fundamental corepresentation of $C(S_N^+\,\hat{*}\,S_M^+)$ is by definition given by:
$$w=\begin{pmatrix}u&0\\0&v\end{pmatrix}$$

But this matrix is magic, because both $u,v$ are magic, and this gives the result.

\medskip

(2) This result, which refines our $N=4$ trick from the proof of Theorem 11.15, follows from (1) with $N=M=2$. Indeed, we have the following computation:
\begin{eqnarray*}
S_2^+\,\hat{*}\,S_2^+
&=&S_2\,\hat{*}\, S_2
=\mathbb Z_2\,\hat{*}\, \mathbb Z_2\\
&\simeq&\widehat{\mathbb Z_2}\,\hat{*}\, \widehat{\mathbb Z_2}
=\widehat{\mathbb Z_2*\mathbb Z_2}\\
&=&\widehat{D_\infty}
\end{eqnarray*}

(3) This result, which generalizes (2), can be deduced as follows:
\begin{eqnarray*}
\widehat{\Gamma}
&\subset&\widehat{\mathbb Z_{N_1}*\ldots*\mathbb Z_{N_k}}
=\widehat{\mathbb Z_{N_1}}\,\hat{*}\,\ldots\,\hat{*}\,\widehat{\mathbb Z_{N_k}}\\
&\simeq&\mathbb Z_{N_1}\,\hat{*}\,\ldots\,\hat{*}\,\mathbb Z_{N_k}
\subset S_{N_1}\,\hat{*}\,\ldots\,\hat{*}\,S_{N_k}\\
&\subset&S_{N_1}^+\,\hat{*}\,\ldots\,\hat{*}\,S_{N_k}^+
\subset S_N^+
\end{eqnarray*}

(4) This follows from (3), because at $N=5$ the dual of the group $\Gamma=\mathbb Z_2*\mathbb Z_3$, which is well-known not to be amenable, embeds into $S_5^+$. As for the general case, that of $S_N^+$ with $N\geq5$, here the result follows by using the embedding $S_5^+\subset S_N^+$.
\end{proof}

Getting now to quantum reflections, things here are more tricky. Following \cite{bb+}, the free analogues of the reflection groups $H_N^s$ can be constructed as follows:

\begin{definition}
The algebra $C(H_N^{s+})$ is the universal $C^*$-algebra generated by $N^2$ normal elements $u_{ij}$, subject to the following relations,
\begin{enumerate}
\item $u=(u_{ij})$ is unitary,

\item $u^t=(u_{ji})$ is unitary,

\item $p_{ij}=u_{ij}u_{ij}^*$ is a projection,

\item $u_{ij}^s=p_{ij}$,
\end{enumerate}
with Woronowicz algebra maps $\Delta,\varepsilon,S$ constructed by universality.
\end{definition}

Here we allow the value $s=\infty$, with the convention that the last axiom disappears in this case. Observe also that at $s<\infty$ the normality condition is actually redundant, because a partial isometry $a$ subject to the relation $aa^*=a^s$ is normal. Finally, the quantum group property comes as before for $O_N^+,U_N^+,S_N^+$, by looking at $u^\Delta,u^\varepsilon,u^S$. 

\bigskip

As a first result now, making the connection with the groups $H_N^s$, we have:

\begin{theorem}
We have an inclusion of compact quantum groups
$$H_N^s\subset H_N^{s+}$$
which is a liberation, in the sense that $H_N^s$ is the classical version of $H_N^{s+}$.
\end{theorem}

\begin{proof}
This follows as before for $O_N\subset O_N^+$, $U_N\subset U_N^+$, $S_N\subset S_N^+$, by using the Gelfand theorem, applied to the quotient of $C(H_N^{s+})$ by its commutator ideal.
\end{proof}

Let us prove now that $H_N^{s+}$ with $s<\infty$ is a quantum permutation group. For this purpose, we must change the fundamental representation. Let us start with:

\begin{definition}
A $(s,N)$-sudoku matrix is a magic unitary of size $sN$, of the form
$$m=\begin{pmatrix}
a^0&a^1&\ldots&a^{s-1}\\
a^{s-1}&a^0&\ldots&a^{s-2}\\
\vdots&\vdots&&\vdots\\
a^1&a^2&\ldots&a^0
\end{pmatrix}$$
where $a^0,\ldots,a^{s-1}$ are $N\times N$ matrices.
\end{definition}

The basic examples of such matrices come from the group $H_n^s$. Indeed, with $w=e^{2\pi i/s}$, each of the $N^2$ matrix coordinates $u_{ij}:H_N^s\to\mathbb C$ decomposes as follows:
$$u_{ij}=\sum_{r=0}^{s-1}w^ra^r_{ij}$$

To be more precise, here each of the functions $a^r_{ij}$ is a function taking values in $\{0,1\}$, and so is a projection in the $C^*$-algebra sense, and it follows from definitions that these projections form indeed a sudoku matrix, in the sense of Definition 11.19. 

\bigskip

Now with the above notion in hand, we have the following result:

\begin{theorem}
The following happen:
\begin{enumerate}
\item The algebra $C(H_N^s)$ is isomorphic to the universal commutative $C^*$-algebra generated by the entries of a $(s,N)$-sudoku matrix.

\item The algebra $C(H_N^{s+})$ is isomorphic to the universal $C^*$-algebra generated by the entries of a $(s,N)$-sudoku matrix.
\end{enumerate}
\end{theorem}

\begin{proof}
The first assertion follows from the second one, via Theorem 11.18. In order to prove the second assertion, consider the universal algebra in the statement:
$$A=C^*\left(a_{ij}^p\ \Big\vert \left(a^{q-p}_{ij}\right)_{pi,qj}=(s,N)-\mbox{sudoku }\right)$$

Consider also the algebra $C(H_N^{s+})$. According to Definition 11.17, this is presented by certain relations $R$, that we will call here level $s$ cubic conditions:
$$C(H_N^{s+})=C^*\left(u_{ij}\ \Big\vert\  u=N\times N\mbox{ level $s$ cubic }\right)$$

We can construct a pair of inverse morphisms between these algebras, as follows:

\medskip

(1) Our first claim is that $U_{ij}=\sum_pw^{-p}a^p_{ij}$ is a level $s$ cubic unitary. Indeed, by using the sudoku condition, the verification of (1-4) in Definition 11.17 is routine.

\medskip

(2) Our second claim is that the elements $A^p_{ij}=\frac{1}{s}\sum_rw^{rp}u^r_{ij}$, with the convention $u_{ij}^0=p_{ij}$, form a level $s$ sudoku unitary. Once again, the proof here is routine.

\medskip

(3) According to the above, we can define a morphism $\Phi:C(H_N^{s+})\to A$ by the formula $\Phi(u_{ij})=U_{ij}$, and a morphism $\Psi:A\to C(H_N^{s+})$ by the formula $\Psi(a^p_{ij})=A^p_{ij}$.

\medskip

(4) We check now the fact that $\Phi,\Psi$ are indeed inverse morphisms:
$$\Psi\Phi(u_{ij})
=\sum_pw^{-p}A^p_{ij}
=\frac{1}{s}\sum_pw^{-p}\sum_rw^{rp}u_{ij}^r
=\frac{1}{s}\sum_{pr}w^{(r-1)p}u_{ij}^r
=u_{ij}$$

As for the other composition, we have the following computation:
$$\Phi\Psi(a^p_{ij})
=\frac{1}{s}\sum_rw^{rp}U_{ij}^r
=\frac{1}{s}\sum_rw^{rp}\sum_qw^{-rq}a_{ij}^q
=\frac{1}{s}\sum_qa_{ij}^q\sum_rw^{r(p-q)}
=a^p_{ij}$$

Thus we have an identification $C(H_N^{s+})=A$, as claimed.
\end{proof}

We can say more about the quantum permutation group meaning of $H_N^{s+}$, in analogy with what we know from before regarding the group $H_N^s$, as follows:

\begin{theorem}
The disjoint union of $N$ copies of the oriented $s$-cycle
$$\xymatrix@R=15pt@C=17pt{
&\bullet\ar[r]&\bullet\ar[dr]\\
\bullet\ar[ur]&&&\bullet\ar[d]\\
\bullet\ar[u]&&&\bullet\ar[dl]\\
&\bullet\ar[ul]&\bullet\ar[l]}\qquad
\xymatrix@R=13pt@C=17pt{\\ \\ \ldots\ldots}
\qquad
\xymatrix@R=15pt@C=17pt{
&\bullet\ar[r]&\bullet\ar[dr]\\
\bullet\ar[ur]&&&\bullet\ar[d]\\
\bullet\ar[u]&&&\bullet\ar[dl]\\
&\bullet\ar[ul]&\bullet\ar[l]}$$
has $H_N^s$ as symmetry group, and $H_N^{s+}$ as quantum symmetry group.
\end{theorem}

\begin{proof}
In what regards the $H_N^s$ assertion, this is something which is obvious. Regarding now the $H_N^{s+}$ assertion, our precise claim here is that we have the following formula, with $d$ being the adjacency matrix of the graph in the statement:
$$C(H_N^{s+})=C(S_{Ns}^+)\big/\big<du=ud\Big>$$

In order to prove this, observe that the adjacency matrix $d$ is given by:
$$d=
\begin{pmatrix}
0&1_N&0&\ldots&0\\
0&0&1_N&\ldots&0\\
\vdots&\vdots&&\ddots&\\
0&0&0&\ldots&1_N\\
1_N&0&0&\ldots&0
\end{pmatrix}$$

Now given a $sN\times sN$ magic unitary, commutation with $d$ means that the matrix is circulant. Thus, we obtain the sudoku relations from Definition 11.19, as desired.
\end{proof}

Finally, we can talk about wreath products, the result here being as follows:

\begin{theorem}
We have the following results:
\begin{enumerate}
\item $H_N^s=\mathbb Z_s\wr S_N$, with $\wr$ being a wreath product.

\item $H_N^{s+}=\mathbb Z_s\wr_*S_N^+$, with $\wr_*$ being a free wreath product.
\end{enumerate}
\end{theorem}

\begin{proof}
In what regards the first result, this is something that we know well. Regarding the second result, the free wreath product of two quantum permutation groups $G\subset S_s^+$, $H\subset S_N^+$, with fundamental corepresentations $u,v$, can be constructed via:
$$C(G\wr_*H)=(C(G)^{*N}*C(H))\big/\big<[u_{ij}^{(a)},v_{ab}]=0\big>$$

But in our case, it is routine to show that $G=\mathbb Z_s$ and $H=S_N^+$ produce in this way the quantum group $G\wr_*H=H_N^{s+}$. Alternatively, all this comes as well from Theorem 11.21, via the following formulae, valid for any connected graph $X$:
$$G(X\ldots X)=G(X)\wr S_N\quad,\quad 
G^+(X\ldots X)=G^+(X)\wr_*S_N^+$$

Thus, one way or another, we are led to the conclusions in the statement. See \cite{bb+}.
\end{proof}

\section*{11c. Laws of characters} 

In order to further study the free quantum groups that we have, we can use representation theory. We have indeed the following version of Tannakian duality:

\begin{theorem}
The following operations are inverse to each other:
\begin{enumerate}
\item The construction $A\to C$, which associates to any Woronowicz algebra $A$ the tensor category formed by the intertwiner spaces $C_{kl}=Hom(u^{\otimes k},u^{\otimes l})$.

\item The construction $C\to A$, which associates to a tensor category $C$ the Woronowicz algebra $A$ presented by the relations $T\in Hom(u^{\otimes k},u^{\otimes l})$, with $T\in C_{kl}$.
\end{enumerate}
\end{theorem}

\begin{proof}
As before with such general results, this is something that we already know for the classical groups, from chapter 4, and the extension to the quantum groups is straightforward. For details here we refer to Woronowicz \cite{wo2} and Malacarne \cite{mal}.
\end{proof}

Next comes the notion of easiness. Let us start with the following definition:

\begin{definition}
A category of partitions is a collection $D=\bigsqcup_{k,l}D(k,l)$ of subsets $D(k,l)\subset P(k,l)$, having the following properties:
\begin{enumerate}
\item Stability under the horizontal concatenation, $(\pi,\sigma)\to[\pi\sigma]$.

\item Stability under vertical concatenation $(\pi,\sigma)\to[^\sigma_\pi]$, with matching middle symbols.

\item Stability under the upside-down turning $*$, with switching of colors, $\circ\leftrightarrow\bullet$.

\item Each set $P(k,k)$ contains the identity partition $||\ldots||$.

\item The sets $P(\emptyset,\circ\bullet)$ and $P(\emptyset,\bullet\circ)$ both contain the semicircle $\cap$.
\end{enumerate}
\end{definition} 

In other words, what we have here is our previous definition for the categories of partitions, from chapter 4, with the last axiom there, stating that $P(k,\bar{k})$ with $|k|=2$ must contain the crossing partition $\slash\hskip-2.0mm\backslash\,$, removed. With this, following \cite{bsp}, we have:

\begin{theorem}
Any category of partitions $D\subset P$ produces a series of compact quantum groups $G=(G_N)$, with $G_N\subset U_N^+$ for any $N\in\mathbb N$, via the formula
$$Hom(u^{\otimes k},u^{\otimes l})=span\left(T_\pi\Big|\pi\in D(k,l)\right)$$
for any $k,l$, and Tannakian duality. We call such quantum groups easy.
\end{theorem}

\begin{proof}
Indeed, the axioms in Definition 11.24 show, via the standard properties of $\pi\to T_\pi$, that the following spaces form a Tannakian category, in the sense of \cite{wo2}:
$$C_{kl}=span\left(T_\pi\Big|\pi\in D(k,l)\right)$$

Thus, Tannakian duality applies, and provides us with a closed subgroup $G_N\subset U_N^+$ such that $Hom(u^{\otimes k},u^{\otimes l})=C_{kl}$, for any colored integers $k,l$, as desired.
\end{proof}

At the level of the examples now, we certainly have the easy groups from chapter 4, and we have as well their free versions constructed above, the result being:

\begin{theorem}
We have easy quantum groups as follows, with calligraphic standing for matching, meaning satisfying $\#\circ=\#\bullet$, as a weighted equality, in each block: 
$$\xymatrix@R=16pt@C=15pt{
&K_N^+\ar[rr]&&U_N^+\\
H_N^+\ar[rr]\ar[ur]&&O_N^+\ar[ur]\\
&S_N^+\ar[rr]\ar[uu]&&C_N^+\ar[uu]\\
S_N^+\ar[uu]\ar[ur]\ar[rr]&&B_N^+\ar[uu]\ar[ur]}
\qquad\xymatrix@R=18pt@C=14pt{\\ \\ :}
\qquad\xymatrix@R=19.2pt@C=0pt{
&\mathcal{NC}_{even}\ar[dl]\ar[dd]&&\mathcal{NC}_2\ar[dl]\ar[ll]\ar[dd]\\
NC_{even}\ar[dd]&&NC_2\ar[dd]\ar[ll]\\
&NC\ar[dl]&&\mathcal{NC}_{12}\ar[dl]\ar[ll]\\
NC&&NC_{12}\ar[ll]}$$
Moreover, the quantum reflection groups $H_N^{s+}$ with $s\in\{1,2,\ldots,\infty\}$ fit on the left diagonal of the left cube, corresponding to the categories $NC^s$, fitting on the right cube. 
\end{theorem}

\begin{proof}
This is identical to the proof of the classical result, from chapter 4, with the diagrams, intertwiners and everything being identical, save for the very end in each case, with the categories generated by the relevant partitions being noncrossing, due to the fact that we removed the crossing axiom, in our updated Definition 11.24.
\end{proof}

Getting now to integration questions, again following chapter 4, we first have:

\begin{theorem}
The asymptotic laws of main characters are as follows,
$$\xymatrix@R=16pt@C=15pt{
&K_N^+\ar[rr]&&U_N^+\\
H_N^+\ar[rr]\ar[ur]&&O_N^+\ar[ur]\\
&S_N^+\ar[rr]\ar[uu]&&C_N^+\ar[uu]\\
S_N^+\ar[uu]\ar[ur]\ar[rr]&&B_N^+\ar[uu]\ar[ur]}
\qquad\xymatrix@R=18pt@C=14pt{\\ \\ :}
\qquad\xymatrix@R=17pt@C=20pt{
&\Pi_1\ar@{-}[rr]\ar@{-}[dd]&&\Gamma_1\ar@{-}[dd]\\
\pi^2_1\ar@{-}[rr]\ar@{-}[dd]\ar@{-}[ur]&&\gamma_1\ar@{-}[dd]\ar@{-}[ur]\\
&\pi_1\ar@{-}[rr]\ar@{-}[uu]&&\Gamma_1^1\ar@{.}[uu]\\
\pi_1\ar@{-}[uu]\ar@{-}[ur]\ar@{-}[rr]&&\gamma_1^1\ar@{-}[uu]\ar@{-}[ur]
}$$
with $H_N^{s+}\to\pi^s_1$ inserted, and each convergence is stationary, starting from $N=4$.
\end{theorem}

\begin{proof}
This is something quite tricky, the idea being as follows:

\medskip

(1) According to Theorem 11.8 and to easiness, we have, in each case:
$$\int_G\chi_1^k=\dim(Fix(u^{\otimes k}))
=\dim\left[span\left(T_\pi\Big|\pi\in D(k)\right)\right]$$

(2) Normally at this point we can use, as in chapter 4, the Lindst\"om formula, telling us that the vectors $T_\pi$ are linearly independent at $N\geq k$, and we obtain the laws in the statement, with the moment convergence being stationary, starting from $N=k$.

\medskip

(3) However, and here comes the point, since we are dealing with noncrossing partitions, instead of Lindst\"om we can use the Di Francesco formula from chapter 10, telling us that  that the vectors $T_\pi$ are linearly independent at $N\geq 4$. Thus, we obtain the laws in the statement, with the convergence being stationary, starting from $N=4$.

\medskip

(4) Let us also mention that, with a bit more care, namely using the appropriate Di Francesco determinant formula, for each of the quantum groups involved, we can see that for $O_N^+,U_N^+$ the convergence is actually stationary starting from $N=2$.

\medskip

(5) Finally, there is as well a purely algebraic proof, for all this, not using Gram determinants, obtained by computing the fusion rules for the various quantum groups involved. We will discuss this later, towards the end of this chapter.
\end{proof}

The above result is quite interesting, because save for having Bercovici-Pata, which was something that we were expecting, we have the stationarity of the convergence, which is something truly amazing, in stark contrast with what happens for the classical groups. So, are the free quantum groups simpler than their classical counterparts?

\bigskip

Which sounds like a good and deep question, and time I guess to ask the cat:

\begin{cat}
Meow.
\end{cat}

Humm, not very helpful all this, usually cat answers with a basic meow when the question is totally trivial. Guess we will have to think some more, at all this.

\bigskip

Moving on, still following the material from chapter 4, we can talk about Weingarten integration and application to the truncated characters, the result being as follows:

\begin{theorem}
The asymptotic laws of truncated characters are as follows,
$$\xymatrix@R=16pt@C=15pt{
&K_N^+\ar[rr]&&U_N^+\\
H_N^+\ar[rr]\ar[ur]&&O_N^+\ar[ur]\\
&S_N^+\ar[rr]\ar[uu]&&C_N^+\ar[uu]\\
S_N^+\ar[uu]\ar[ur]\ar[rr]&&B_N^+\ar[uu]\ar[ur]}
\qquad\xymatrix@R=18pt@C=14pt{\\ \\ :}
\qquad\xymatrix@R=17pt@C=20pt{
&\Pi_t\ar@{-}[rr]\ar@{-}[dd]&&\Gamma_t\ar@{-}[dd]\\
\pi^2_t\ar@{-}[rr]\ar@{-}[dd]\ar@{-}[ur]&&\gamma_t\ar@{-}[dd]\ar@{-}[ur]\\
&\pi_t\ar@{-}[rr]\ar@{-}[uu]&&\Gamma_t^t\ar@{.}[uu]\\
\pi_t\ar@{-}[uu]\ar@{-}[ur]\ar@{-}[rr]&&\gamma_t^t\ar@{-}[uu]\ar@{-}[ur]
}$$
with $H_N^{s+}\to\pi^s_t$ inserted, and with the convergence being no longer stationary.
\end{theorem}

\begin{proof}
As before in the classical case, we have a Weingarten integration formula for the easy quantum groups, as follows, with $W_{kN}=G_{kN}^{-1}$, where $G_{kN}(\pi,\nu)=N^{|\pi\vee\nu|}$:
$$\int_{G_N}u_{i_1j_1}\ldots u_{i_kj_k}=\sum_{\pi,\nu\in D(k)}\delta_\pi(i)\delta_\nu(j)W_{kN}(\pi,\nu)$$ 

To be more precise, this formula holds indeed, exactly as in the classical case, and with the remark that $*$ exponents can be added too, in the obvious way. Next, by using this formula, we have the following computation, for the truncated characters:
\begin{eqnarray*}
\int_{G_N}(u_{11}+\ldots +u_{ss})^k
&=&\sum_{i_1=1}^{s}\ldots\sum_{i_k=1}^s\int_{G_N}u_{i_1i_1}\ldots u_{i_ki_k}\\
&=&\sum_{\pi,\nu\in D(k)}W_{kN}(\pi,\nu)\sum_{i_1=1}^{s}\ldots\sum_{i_k=1}^s\delta_\pi(i)\delta_\nu(i)\\
&=&\sum_{\pi,\nu\in D(k)}W_{kN}(\pi,\nu)G_{ks}(\nu,\pi)\\
&=&Tr(W_{kN}G_{ks})
\end{eqnarray*}

The point now is that in the uniform case the Gram matrix is asymptotically diagonal, $G_{kN}\simeq diag(N^{|\pi|})$, with this being something clear for all the examples, and true in general too, either via an abstract argument, or via the classification results in \cite{twe}. Now the inverse of $G_{kN}$ being asymptotically diagonal too, $W_{kN}\simeq diag(N^{-|\pi|})$, we obtain:
$$\lim_{N\to\infty}\int_{G_N}\chi_t^k=\sum_{\pi\in D(k)}t^{|\pi|}$$

Thus, we are led to the conclusions in the statement. And with the convergence being no longer stationary, which brings peace, guess cat was wrong with his meow.
\end{proof}

As a next task, observe that our cube in Theorem 11.29 does not look very good, with $S_N^+$ appearing twice, and with this being something similar to what we met in chapter 4, with $S_N$ appearing in the classical cube there twice. So the question is, can we have soemthing better, based on this? In answer, in the orthogonal case we have:

\begin{theorem}
The basic orthogonal easy quantum groups are
$$\xymatrix@R=16pt@C=15pt{
&H_N^+\ar[rr]&&O_N^+\\
S_N^+\ar[rr]\ar[ur]&&B_N^+\ar[ur]\\
&H_N\ar[rr]\ar[uu]&&O_N\ar[uu]\\
S_N\ar[uu]\ar[ur]\ar[rr]&&B_N\ar[uu]\ar[ur]
}\qquad\xymatrix@R=18pt@C=14pt{\\ \\ :}
\qquad\xymatrix@R=15pt@C=20pt{
&\pi_t^2\ar@{-}[rr]\ar@{-}[dd]&&\gamma_t\ar@{-}[dd]\\
\pi_t\ar@{-}[rr]\ar@{-}[dd]\ar@{-}[ur]&&\gamma_t^t\ar@{-}[dd]\ar@{-}[ur]\\
&p_t^2\ar@{-}[rr]\ar@{-}[uu]&&g_t\ar@{.}[uu]\\
p_t\ar@{-}[uu]\ar@{-}[ur]\ar@{-}[rr]&&g_t^t\ar@{-}[uu]\ar@{-}[ur]
}$$
with on the right being the asymptotic laws of truncated characters.
\end{theorem}

\begin{proof}
This comes indeed from our results from chapter 4, and from here.
\end{proof}

In the unitary case now, we have a similar result, as follows:

\begin{theorem}
The basic unitary easy quantum groups are
$$\xymatrix@R=16pt@C=15pt{
&K_N^+\ar[rr]&&U_N^+\\
S_N^+\ar[rr]\ar[ur]&&C_N^+\ar[ur]\\
&K_N\ar[rr]\ar[uu]&&U_N\ar[uu]\\
S_N\ar[uu]\ar[ur]\ar[rr]&&C_N\ar[uu]\ar[ur]
}\qquad\xymatrix@R=18pt@C=14pt{\\ \\ :}
\qquad\xymatrix@R=16pt@C=20pt{
&\Pi_t\ar@{-}[rr]\ar@{-}[dd]&&\Gamma_t\ar@{-}[dd]\\
\pi_t\ar@{-}[rr]\ar@{-}[dd]\ar@{-}[ur]&&\Gamma_t^t\ar@{-}[dd]\ar@{-}[ur]\\
&P_t\ar@{-}[rr]\ar@{-}[uu]&&G_t\ar@{.}[uu]\\
p_t\ar@{-}[uu]\ar@{-}[ur]\ar@{-}[rr]&&G_t^t\ar@{-}[uu]\ar@{-}[ur]
}$$
with on the right being the asymptotic laws of truncated characters.
\end{theorem}

\begin{proof}
This comes again from our results from chapter 4, and from here.
\end{proof}

Next, the point is that we cannot really merge our pairs of 3D cubes into 4D cubes. However, we can get rid of some of the objects, and we are led in this way to:

\begin{theorem}
The very basic easy quantum groups are
$$\xymatrix@R=16.1pt@C=15pt{
&K_N^+\ar[rr]&&U_N^+\\
H_N^+\ar[rr]\ar[ur]&&O_N^+\ar[ur]\\
&K_N\ar[rr]\ar[uu]&&U_N\ar[uu]\\
H_N\ar[uu]\ar[ur]\ar[rr]&&O_N\ar[uu]\ar[ur]
}\qquad\xymatrix@R=18pt@C=14pt{\\ \\ :}
\qquad\xymatrix@R=16pt@C=20pt{
&\Pi_t\ar@{-}[rr]\ar@{-}[dd]&&\Gamma_t\ar@{-}[dd]\\
\pi_t^2\ar@{-}[rr]\ar@{-}[dd]\ar@{-}[ur]&&\gamma_t\ar@{-}[dd]\ar@{-}[ur]\\
&P_t\ar@{-}[rr]\ar@{-}[uu]&&G_t\ar@{.}[uu]\\
p_t^2\ar@{-}[uu]\ar@{-}[ur]\ar@{-}[rr]&&g_t\ar@{-}[uu]\ar@{-}[ur]
}$$
with on the right being the asymptotic laws of truncated characters.
\end{theorem}

\begin{proof}
This follows indeed from Theorems 11.30 and 11.31. Observe also that our cube, called ``standard cube'', is something quite subtle, featuring:

\medskip

(1) The classical/free correspondence on the vertical.

\medskip

(2) An interesting discrete/continuous correspondence on the horizontal.

\medskip

(3) Also, an interesting real/complex correspondence on the depth.

\medskip

Summarizing, what we have here are true 3D frames, for both quantum algebra and probability theory, taken in a large sense. Which is really nice, and conceptual.
\end{proof}

Getting now to classification matters, many things can be said here. To start with, in the orthogonal case we have the following result, going back to \cite{bsp}:

\begin{theorem}
The classical and free uniform orthogonal easy quantum groups, with inclusions between them, are as follows:
$$\xymatrix@R=16pt@C=15pt{
&H_N^+\ar[rr]&&O_N^+\\
S_N^+\ar[rr]\ar[ur]&&B_N^+\ar[ur]\\
&H_N\ar[rr]\ar[uu]&&O_N\ar[uu]\\
S_N\ar[uu]\ar[ur]\ar[rr]&&B_N\ar[uu]\ar[ur]
}$$
Moreover, this is an intersection/easy generation diagram, in the sense that for any of its square subdiagrams $P\subset Q,R\subset S$ we have $P=Q\cap R$ and $\{Q,R\}=S$.
\end{theorem}

\begin{proof}
There are several things to be proved, the idea being as follows:

\medskip

(1) We know that the quantum groups in the statement are indeed easy and uniform, the corresponding categories of partitions being as follows:
$$\xymatrix@R=16pt@C=3pt{
&NC_{even}\ar[dl]\ar[dd]&&NC_2\ar[dl]\ar[ll]\ar[dd]\\
NC\ar[dd]&&NC_{12}\ar[dd]\ar[ll]\\
&P_{even}\ar[dl]&&P_2\ar[dl]\ar[ll]\\
P&&P_{12}\ar[ll]
}$$

Since this latter diagram is an intersection and generation diagram, we conclude that we have an intersection and easy generation diagram of quantum groups, as stated.

\medskip

(2) Regarding now the classification, consider first an easy group $S_N\subset G_N\subset O_N$. This must come from a certain category $P_2\subset D\subset P$, and if we assume $G=(G_N)$ to be uniform, then $D$ is uniquely determined by the subset $L\subset\mathbb N$ consisting of the sizes of the blocks of the partitions in $D$. Our claim is that the admissible sets are as follows:

\medskip

-- $L=\{2\}$, producing $O_N$.

\medskip

-- $L=\{1,2\}$, producing $B_N$.

\medskip

-- $L=\{2,4,6,\ldots\}$, producing $H_N$.

\medskip

-- $L=\{1,2,3,\ldots\}$, producing $S_N$.

\medskip

(3) Indeed, in one sense, this follows from our easiness results for $O_N,B_N,H_N,S_N$. In the other sense now, assume that $L\subset\mathbb N$ is such that the set $P_L$ consisting of partitions whose sizes of the blocks belong to $L$ is a category of partitions. We know from the axioms of the categories of partitions that the semicircle $\cap$ must be in the category, so we have $2\in L$. We claim that the following conditions must be satisfied as well:
$$k,l\in L,\,k>l\implies k-l\in L$$
$$k\in L,\,k\geq 2\implies 2k-2\in L$$

(4) Indeed, we will prove that both conditions follow from the axioms of the categories of
partitions. Let us denote by $b_k\in P(0,k)$ the one-block partition:
$$b_k=\left\{\begin{matrix}\sqcap\hskip-0.7mm \sqcap&\ldots&\sqcap\\
1\hskip2mm 2&\ldots&k\end{matrix} \right\}$$

For $k>l$, we can write $b_{k-l}$ in the following way:
$$b_{k-l}=\left\{\begin{matrix}\sqcap\hskip-0.7mm
\sqcap&\ldots&\ldots&\ldots&\ldots&\sqcap\\ 1\hskip2mm 2&\ldots&l&l+1&\ldots&k\\
\sqcup\hskip-0.7mm \sqcup&\ldots&\sqcup&|&\ldots&|\\ &&&1&\ldots&k-l\end{matrix}\right\}$$

In other words, we have the following formula:
$$b_{k-l}=(b_l^*\otimes |^{\otimes k-l})b_k$$

Since all the terms of this composition are in $P_L$, we have $b_{k-l}\in P_L$, and this proves our first claim. As for the second claim, this can be proved in a similar way, by capping two adjacent $k$-blocks with a $2$-block, in the middle.

\medskip

(5) With these conditions in hand, we can conclude in the following way:

\medskip

\underline{Case 1}. Assume $1\in L$. By using the first condition with $l=1$ we get:
$$k\in L\implies k-1\in L$$

This condition shows that we must have $L=\{1,2,\ldots,m\}$, for a certain number $m\in\{1,2,\ldots,\infty\}$. On the other hand, by using the second condition we get:
\begin{eqnarray*}
m\in L
&\implies&2m-2\in L\\
&\implies&2m-2\leq m\\
&\implies&m\in\{1,2,\infty\}
\end{eqnarray*}

The case $m=1$ being excluded by the condition $2\in L$, we reach to one of the two sets producing the groups $S_N,B_N$.

\medskip

\underline{Case 2}. Assume $1\notin L$. By using the first condition with $l=2$ we get:
$$k\in L\implies k-2\in L$$

This condition shows that we must have $L=\{2,4,\ldots,2p\}$, for a certain number $p\in\{1,2,\ldots,\infty\}$. On the other hand, by using the second condition we get:
\begin{eqnarray*}
2p\in L
&\implies&4p-2\in L\\
&\implies&4p-2\leq 2p\\
&\implies&p\in\{1,\infty\}
\end{eqnarray*}

Thus $L$ must be one of the two sets producing $O_N,H_N$, and we are done. 

\medskip

(6) In the free case, $S_N^+\subset G_N\subset O_N^+$, the situation is quite similar, the admissible sets being once again the above ones, producing this time the free quantum groups $O_N^+,B_N^+,H_N^+,S_N^+$. For details regarding all this, and for more, we refer to \cite{bsp}.
\end{proof}

As a continuation of this, we can try to lift the various conditions in Theorem 11.33, except of course for easiness, and things here become quite technical, as follows:

\begin{fact}
The following happen, regarding the easy quantum groups:
\begin{enumerate}
\item In the classical orthogonal case, not necessarily uniform, we have two more solutions, namely $S_N'=S_N\times\mathbb Z_2$ and $B_N'=B_N\times\mathbb Z_2$.

\item In the free orthogonal case, not necessarily uniform, we have three more solutions, namely $S_N^{\prime+}=S_N^+\times\mathbb Z_2$ and the two possible liberations of $B_N'$.

\item In the general orthogonal case, not necessarily classical or free, we have as main examples the half-liberations $H_N^*,O_N^*$, coming via $abc=cba$, both non-uniform.

\item The classification is in fact possible in the general orthogonal case, notably featuring many interesting intermediate objects $H_N\subset G\subset H_N^+$.

\item In the unitary case, classical or free, the classification is possible too, by making your way through the jungle of products, namely $G\times\mathbb Z_s$, $G\cdot\mathbb Z_s$ and $\mathbb Z_s\wr G$.

\item The classification of intermediate quantum groups $U_N\subset G\subset U_N^+$ is possible too, having interesting similarities with the classification for $H_N\subset G\subset H_N^+$.
\end{enumerate}
\end{fact}

And we will end our discussion of easiness with this. For more on the above, you can have a look at the papers of Weber and his collaborators \cite{mwe}, \cite{rwe}, \cite{twe}. And with the remark that these papers usually focus on combinatorics and classification, with the probability theory, in relation with characters and their laws, still in need to be done.

\bigskip

Finally, in relation with the standard cube from Theorem 11.32, we have:

\begin{fact}[Ground Zero]
Under very strong combinatorial axioms
$$\xymatrix@R=16.1pt@C=15pt{
&K_N^+\ar[rr]&&U_N^+\\
H_N^+\ar[rr]\ar[ur]&&O_N^+\ar[ur]\\
&K_N\ar[rr]\ar[uu]&&U_N\ar[uu]\\
H_N\ar[uu]\ar[ur]\ar[rr]&&O_N\ar[uu]\ar[ur]
}$$
are the only quantum groups surviving.
\end{fact}

To be more precise, this is something which comes from the above-mentioned classification results of Weber and his collaborators, by piling enough axioms, as to kill everything, except for the very strong. Philosophically, this is of course something quite interesting, and for more on all this, you can have a look at any of my quantum group books.

\section*{11d. Clebsch-Gordan}

As a last topic for this chapter, let us go back to Cat 11.28, and try to better understand that comment. So, the question was whether free quantum groups are simpler than classical groups, and with cat suggesting that the answer to this is obvious.

\bigskip

So, is this correct or not? As a first further piece of evidence for this, we have:

\begin{theorem}
Unlike the passage $O_N\to U_N$, which is something non-trivial, coming from $\mathfrak o_N\to\mathfrak u_N$ at Lie algebra level, the passage $O_N^+\to U_N^+$ is a free complexification
$$U_N^+=\widetilde{O_N^+}$$
with this meaning that we have an embedding of algebras as follows,
$$C(U_N^+)\subset C(\mathbb T)*C(O_N^+)\quad,\quad u=zv$$
up to identifying the full and reduced versions of these algebras.
\end{theorem}

\begin{proof}
There are several things going on here, the idea being as follows:

\medskip

(1) The orthogonal and unitary groups $O_N,N_N$ are both Lie groups, and the corresponding Lie algebras $\mathfrak o_N,\mathfrak u_N$ are easy to compute, their formulae being as follows:
$$\mathfrak o_N=\left\{ A\in M_N(\mathbb R)\Big|A^t=-A\right\}\quad,\quad 
\mathfrak u_N=\left\{ B\in M_N(\mathbb C)\Big|B^*=-B\right\}$$

By doing the linear algebra, we conclude that we have a direct sum decomposition of real linear spaces as follows, with $\Delta\subset M_N(\mathbb R)$ being the diagonal matrices:
$$\mathfrak u_N\simeq\mathfrak o_N\oplus\Delta\oplus\mathfrak o_N$$

Now based on this, we can say that $\mathfrak o_N\to\mathfrak u_N$ is some sort of ``complexification'', and so that the passage $O_N\to U_N$ appears too as a complexification, as follows:
$$O_N\to\mathfrak o_N\to\mathfrak u_n\to U_N$$

(2) Getting now to quantum groups, given an arbitrary pair $(G,v)$, we can associate to it its free complexification $(\widetilde{G},\tilde{v})$, via the formula in the statement, namely:
$$C(\widetilde{G})=<\tilde{v}_{ij}>\subset C(\mathbb T)*C(G)\quad,\quad \tilde{v}=zv$$

Now let $v,zv,u$ the fundamental representations of the following quantum groups:
$$O_N^+\subset\widetilde{O_N^+}\subset U_N^+$$

At the level of the associated Hom spaces we obtain reverse inclusions, as follows:
$$Hom(v^{\otimes k},v^{\otimes l})
\supset Hom((zv)^{\otimes k},(zv)^{\otimes l})
\supset Hom(u^{\otimes k},u^{\otimes l})$$

But the spaces on the left and on the right are known from Theorem 11.26, with the result there stating that these spaces are as follows:
$$span\left(T_\pi\Big|\pi\in NC_2(k,l)\right)\supset span\left(T_\pi\Big|\pi\in\mathcal{NC}_2(k,l)\right)$$

Regarding now the spaces in the middle, these are obtained from those on the left by some sort of ``coloring'', so we obtain the same spaces as those on the right. Thus, our embedding $\widetilde{O_N^+}\subset U_N^+$ must be an isomorphism, modulo the usual equivalence relation.
\end{proof}

As an interesting theoretical consequence of the above result, which is something quite fascinating, bringing even more evidence for Cat 11.28, we have:

\begin{theorem}
We have an identification of quantum groups, as follows:
$$PO_N^+=PU_N^+$$
Thus, free projective geometry is scalarless, at the same time real and complex.
\end{theorem}

\begin{proof}
Again, we have several things going on here, as follows:

\medskip

(1) To start with, in the general context of the Woronowicz quantum groups, given a pair $(G,v)$, we can talk about its projective version $(PG,u)$, constructed as follows:
$$C(PG)=<u_{ij}>\subset C(G)\quad,\quad u=v\otimes\bar{v}$$

(2) Now in the case of $G=U_N^+$, by using Theorem 11.36, we obtain right away:
$$PU_N^+=P\widetilde{O_N^+}=PO_N^+$$

(3) We can deduce this as well directly. Indeed, with notations as before, we have:
$$Hom\left((v\otimes v)^k,(v\otimes v)^l\right)=span\left(T_\pi\Big|\pi\in NC_2((\circ\bullet)^k,(\circ\bullet)^l)\right)$$
$$Hom\left((u\otimes\bar{u})^k,(u\otimes\bar{u})^l\right)=span\left(T_\pi\Big|\pi\in \mathcal{NC}_2((\circ\bullet)^k,(\circ\bullet)^l)\right)$$

Now the sets on the right being equal, $PO_N^+\subset PU_N^+$ must be an isomorphism.
\end{proof}

Moving on, even more good news come by looking at the fusion rules:

\begin{theorem}
The irreducible representations of $O_N^+$  with $N\geq2$ can be labeled by positive integers, $r_k$ with $k\in\mathbb N$, the fusion rules for these representations are
$$r_k\otimes r_l=r_{|k-l|}+r_{|k-l|+2}+\ldots+r_{k+l}$$ 
exactly as the Clebsch-Gordan rules for $SU_2$, and the corresponding dimensions are 
$$\dim r_k=(q^{k+1}-q^{-k-1})/(q-q^{-1})$$
with $q+q^{-1}=N$. And with this being something far simpler than for $O_N$.
\end{theorem}

\begin{proof}
The proof is nearly identical to the proof for $SU_2$ from chapter 4, with the needed estimate there, involving the Catalan numbers, coming for $O_N^+$ as follows:
\begin{eqnarray*}
\int_{O_N^+}\chi^{2k}
&=&\dim(Fix(u^{\otimes 2k}))\\
&=&\dim\left(span\left(T_\pi|\pi\in NC_2(2k)\right)\right)\\
&\leq&|NC_2(2k)|\\
&=&C_k
\end{eqnarray*}

Thus, by reasoning as in chapter 4, we can construct the irreducible representations by recurrence, with the conclusions in the statement.
\end{proof}

Regarding now the quantum group $S_N^+$, we have a similar result for it, as follows:

\begin{theorem}
The irreducible representations of $S_N^+$ with $N\geq4$ can be labeled by positive integers, $r_k$ with $k\in\mathbb N$, the fusion rules for these representations are
$$r_k\otimes r_l=r_{|k-l|}+r_{|k-l|+1}+\ldots+r_{k+l}$$ 
exactly as the Clebsch-Gordan rules for $SO_3$, and the corresponding dimensions are 
$$\dim r_k=(q^{k+1}-q^{-k})/(q-1)$$
with $q+q^{-1}=N-2$. And with this being something far simpler than for $S_N$.
\end{theorem}

\begin{proof}
Similar story here, with however a subtlety coming at the technical level, due to the fact, explained in the case of $SO_3$ in chapter 4, that we need here an equality involving the Catalan numbers, instead of just an estimate. But this equality can be deduced as follows, by using easiness, and the Di Francesco results from chapter 10:
\begin{eqnarray*}
\int_{S_N^+}\chi^k
&=&\dim(Fix(u^{\otimes k}))\\
&=&\dim\left(span\left(T_\pi|\pi\in NC(k)\right)\right)\\
&=&|NC(k)|\\
&=&C_k
\end{eqnarray*}

As for the last assertion, as before with $O_N$ vs $O_N^+$, prove me wrong here.
\end{proof}

Next, regarding the quantum group $U_N^+$, the result here is as follows:

\begin{theorem}
The irreducible representations of the quantum group $U_N^+$ with $N\geq2$ can be indexed by $\mathbb N*\mathbb N$, with the fusion rules being as follows,
$$r_k\otimes r_l=\sum_{k=xy,l=\bar{y}z}r_{xz}$$
and the corresponding dimensions $\dim r_k$ can be computed by recurrence. And with this being something far simpler than for $U_N$.
\end{theorem}

\begin{proof}
Quite similar story here, the idea being as follows:

\medskip

(1) This can be proved directly, a bit as for $O_N^+$ before, with the needed estimate, when constructing the representations $r_k$ by recurrence, coming from easiness.

\medskip

(2) Alternatively, this comes from Theorems 11.36 and 11.38, with the above fusion rules being a ``free complexification'' of the Clebsch-Gordan rules for $O_N^+$. 

\medskip

(3) Finally, let me mention that this result is particularly sweet for me, because this was the main theorem in my PhD thesis, done at Paris 7 some 30 years ago, under the guidance of Georges Skandalis. Thanks Georges for everything.
\end{proof}

Finally, regarding the quantum groups $H_N^{s+}$, the result here is as follows:

\begin{theorem}
The irreducible representations of $H_N^{s+}$ with $N\geq4$ can be labeled $r_x$ with $x\in<\mathbb Z_s>$, the fusion rules being as follows,
$$r_x\otimes r_y=\sum_{x=vz,y=\bar{z}w}r_{vw}+r_{v\cdot w}$$
with the dot product on the right being given by 
$$(i_1\ldots i_k)\cdot (j_1\ldots j_l)=i_1\ldots i_{k-1}(i_k+j_1)j_2\ldots j_l$$
and with this being far simpler than for the group $H_N^s$ itself.
\end{theorem}

\begin{proof}
As before with other results of this type, this is something of Clebsch-Gordan nature, that we already know, from Theorem 11.39, to hold at $s=1$. Observe also the similarity with Theorem 11.40. As for the proof, this is via the standard methods, namely recurrence, and with the needed estimates coming from easiness. See \cite{bve}.
\end{proof}

Well, looks like cat was right, the free world seems infinitely simpler than the classical world, no question about this. So, what do you think cat, did I get it right?

\begin{cat}
Well, now that you understood, get back to work, and rewrite this book, with free probability coming before classical probability.
\end{cat}

Humm, makes sense, but in practice, how in the world to do this, assume that the reader is familiar with all that quantum mechanics algebra. We're not all cats.

\section*{11e. Exercises}

Welcome to free quantum groups, and free geometry, and as exercises, we have:

\begin{exercise}
Fight a bit with $m_t\leftrightarrow\mu_t$, via transforms or cumulants.
\end{exercise}

\begin{exercise}
Learn from Woronowicz all the compact quantum group basics.
\end{exercise}

\begin{exercise}
Clarify all the details, in relation with $B_N^+$ and $C_N^+$.
\end{exercise}

\begin{exercise}
Clarify as well all the details, in relation with $H_N^{s+}$.
\end{exercise}

\begin{exercise}
Learn the proof of Tannakian duality for compact quantum groups.
\end{exercise}

\begin{exercise}
Learn more about easy quantum groups, and their classification.
\end{exercise}

\begin{exercise}
Further meditate at the uniqueness of free projective geometry.
\end{exercise}

\begin{exercise}
Compute the fusion rules for all classical easy groups.
\end{exercise}

As bonus exercise, get to know reggae music. Can be of help, with quantum mechanics.

\chapter{Free geometry}

\section*{12a. Spheres, tori}

With free quantum groups discussed, time to get into the real thing, free geometry. This will open up a nearly infinite wonderland, where we can do lots of probability computations. We will explore a bit all this, in this chapter, with the conclusion that some basic computations that can be done fit with Bercovici-Pata, while some other computations that can be done, still very basic, substantially deviate from Bercovici-Pata.

\bigskip

To start with, following \cite{bgo}, we can talk about free spheres, as follows:

\index{free sphere}

\begin{definition}
We have free real and complex spheres, defined via
$$C(S^{N-1}_{\mathbb R,+})=C^*\left(x_1,\ldots,x_N\Big|x_i=x_i^*,\sum_ix_i^2=1\right)$$
$$C(S^{N-1}_{\mathbb C,+})=C^*\left(x_1,\ldots,x_N\Big|\sum_ix_ix_i^*=\sum_ix_i^*x_i=1\right)$$
with the symbol $C^*$ standing for universal enveloping $C^*$-algebra.
\end{definition}

Here the fact that these algebras are indeed well-defined comes from the following estimate, which shows that the biggest $C^*$-norms on these $*$-algebras are bounded:
$$||x_i||^2
=||x_ix_i^*||
\leq\Big|\Big|\sum_ix_ix_i^*\Big|\Big|
=1$$

As a first result now, regarding the above free spheres, we have:

\begin{theorem}
We have embeddings of compact quantum spaces, as follows,
$$\xymatrix@R=15mm@C=15mm{
S^{N-1}_{\mathbb R,+}\ar[r]&S^{N-1}_{\mathbb C,+}\\
S^{N-1}_\mathbb R\ar[r]\ar[u]&S^{N-1}_\mathbb C\ar[u]
}$$
and the spaces on top appear as liberations of the spaces on the bottom.
\end{theorem}

\begin{proof}
In practice, we must establish the following isomorphisms, where the symbol $C^*_{comm}$ stands for ``universal commutative $C^*$-algebra generated by'':
$$C(S^{N-1}_\mathbb R)=C^*_{comm}\left(x_1,\ldots,x_N\Big|x_i=x_i^*,\sum_ix_i^2=1\right)$$
$$C(S^{N-1}_\mathbb C)=C^*_{comm}\left(x_1,\ldots,x_N\Big|\sum_ix_ix_i^*=\sum_ix_i^*x_i=1\right)$$

It is enough to establish the second isomorphism. So, consider the second universal commutative $C^*$-algebra $A$ constructed above. Since the standard coordinates on $S^{N-1}_\mathbb C$ satisfy the defining relations for $A$, we have a quotient map of as follows:
$$A\to C(S^{N-1}_\mathbb C)$$

Conversely, let us write $A=C(S)$, by using the Gelfand theorem. The variables $x_1,\ldots,x_N$ become in this way true coordinates, providing us with an embedding $S\subset\mathbb C^N$. Also, the quadratic relations become $\sum_i|x_i|^2=1$, so we have $S\subset S^{N-1}_\mathbb C$. Thus, we have a quotient map $C(S^{N-1}_\mathbb C)\to A$, as desired, and this gives the results.
\end{proof}

By using the free spheres constructed above, we can now formulate:

\begin{definition}
A real algebraic manifold $X\subset S^{N-1}_{\mathbb C,+}$ is a closed quantum subspace defined, at the level of the corresponding $C^*$-algebra, by a formula of type
$$C(X)=C(S^{N-1}_{\mathbb C,+})\Big/\Big<f_i(x_1,\ldots,x_N)=0\Big>$$
for certain family of noncommutative polynomials, as follows:
$$f_i\in\mathbb C<x_1,\ldots,x_N>$$
We denote by $\mathcal C(X)$ the $*$-subalgebra of $C(X)$ generated by the coordinates $x_1,\ldots,x_N$. 
\end{definition}

As a basic example here, we have the free real sphere $S^{N-1}_{\mathbb R,+}$. The classical spheres $S^{N-1}_\mathbb C,S^{N-1}_\mathbb R$, and their real submanifolds, are covered as well by this formalism. We have as well as examples the compact quantum groups, embedded as follows:
$$G\subset U_N^+\ \implies \ G\subset U_N^+\subset S^{N^2-1}_{\mathbb C,+}\quad,\quad x_{ij}=\frac{u_{ij}}{\sqrt{N}}$$

At the level of the general theory, we have the following version of the Gelfand theorem:

\begin{theorem}
If $X\subset S^{N-1}_{\mathbb C,+}$ is an algebraic manifold, as above, we have
$$X_{class}=\left\{x\in S^{N-1}_\mathbb C\Big|f_i(x_1,\ldots,x_N)=0\right\}$$
and $X$ appears as a liberation of $X_{class}$.
\end{theorem}

\begin{proof}
This is something that we already met, in the context of the free spheres. In general, the proof is similar, by using the Gelfand theorem. Indeed, if we denote by $X_{class}'$ the manifold constructed in the statement, then we have a quotient map of $C^*$-algebras as follows, mapping standard coordinates to standard coordinates:
$$C(X_{class})\to C(X_{class}')$$

Conversely now, from $X\subset S^{N-1}_{\mathbb C,+}$ we obtain $X_{class}\subset S^{N-1}_\mathbb C$. Now since the relations defining $X_{class}'$ are satisfied by $X_{class}$, we obtain an inclusion $X_{class}\subset X_{class}'$. Thus, at the level of algebras of continuous functions, we have a quotient map of $C^*$-algebras as follows, mapping standard coordinates to standard coordinates:
$$C(X_{class}')\to C(X_{class})$$

Thus, we have constructed a pair of inverse morphisms, and we are done.
\end{proof}

Finally, once again at the level of the general theory, we have:

\begin{definition}
We agree to identify two real algebraic submanifolds $X,Y\subset S^{N-1}_{\mathbb C,+}$ when we have a $*$-algebra isomorphism between $*$-algebras of coordinates
$$f:\mathcal C(Y)\to\mathcal C(X)$$
mapping standard coordinates to standard coordinates.
\end{definition}

Observe that this convention, which is something quite technical, related to amenability issues, fits with the similar convention from chapter 11, via the embeddings $G\subset S^{N^2-1}_{\mathbb C,+}$ mentioned after Definition 12.3. Observe also that our whole formalism, Definitions 12.3 and 12.5, fixes some chronic bugs of Gelfand duality, that we preferred not to really talk about, in chapters 2 and 5. Well, now you know, and very nice all this.

\bigskip

With this discussed, let us go back now to the spheres, and try to integrate over them. For this purpose, it is convenient to deal at the same time with all 4 spheres in Theorem 12.2, that we will call ``basic spheres'', and generically denote $S$. With the convention that the corresponding unitary quantum groups are denoted $U$, we first have:

\begin{proposition}
For the basic spheres, we have a diagram as follows,
$$\xymatrix@R=50pt@C=50pt{
C(S)\ar[r]^\alpha\ar[d]^\gamma&C(S)\otimes C(U)\ar[d]^{\gamma\otimes id}\\
C(U)\ar[r]^\Delta&C(U)\otimes C(U)
}$$
where on top $\alpha(x_i)=\sum_jx_j\otimes u_{ji}$, and on the left $\gamma(x_i)=u_{1i}$.
\end{proposition}

\begin{proof}
The diagram in the statement commutes indeed on the standard coordinates $x_i$. Thus by linearity and multiplicativity, the whole the diagram commutes.
\end{proof}

In practice now, the map $\alpha$ constructed above is what is called a coaction, and skipping some discussion here, the conclusion is that we have a good relationship $S\leftrightarrow U$, a bit like in the classical case. Now based on this, in analogy with what happens in the classical case, we can introduce the uniform integration over our spheres $S$, as follows:

\index{integration over spheres}

\begin{definition}
We endow each of the algebras $C(S)$ with the functional
$$\int_S:C(S)\to C(U)\to\mathbb C$$
obtained by composing the morphism $x_i\to u_{1i}$ with the Haar integration of $C(U)$.
\end{definition}

In order to efficiently integrate now over our generic spheres $S$, in the lack of some trick like spherical coordinates, we can use the Weingarten formula:

\index{integration over spheres}
 
\begin{theorem}
The integration over the basic spheres is given by
$$\int_Sx_{i_1}^{e_1}\ldots x_{i_k}^{e_k}=\sum_{\pi,\nu\in D(k)}\delta_\nu(i)W_{kN}(\pi,\nu)$$
where $D\subset P$ are the partitions for $U$, and $W_{kN}=G_{kN}^{-1}$, with $G_{kN}(\pi,\nu)=N^{|\pi\vee\nu|}$. 
\end{theorem}

\begin{proof}
According to our conventions, the integration over $S$ is a particular case of the integration over $U$, via $x_i=u_{1i}$. By using the Weingarten formula for $U$, we get:
\begin{eqnarray*}
\int_Sx_{i_1}^{e_1}\ldots x_{i_k}^{e_k}
&=&\int_Uu_{1i_1}^{e_1}\ldots u_{1i_k}^{e_k}\\
&=&\sum_{\pi,\nu\in D(k)}\delta_\pi(1)\delta_\nu(i)W_{kN}(\pi,\nu)\\
&=&\sum_{\pi,\nu\in D(k)}\delta_\nu(i)W_{kN}(\pi,\nu)
\end{eqnarray*}

Thus, we are led to the formula in the statement.
\end{proof}

Still following \cite{bgo}, we have the following key result:

\index{ergodicity}

\begin{theorem}
The integration functional of $S$ has the ergodicity property
$$\left(id\otimes\int_U\right)\alpha(x)=\int_Sx$$
where $\alpha:C(S)\to C(S)\otimes C(U)$ is the universal affine coaction map.
\end{theorem}

\begin{proof}
In the real case, $x_i=x_i^*$, it is enough to check the equality in the statement on an arbitrary product of coordinates, $x_{i_1}\ldots x_{i_k}$. The left term is as follows:
\begin{eqnarray*}
\left(id\otimes\int_U\right)\alpha(x_{i_1}\ldots x_{i_k})
&=&\sum_{j_1\ldots j_k}x_{j_1}\ldots x_{j_k}\int_Uu_{j_1i_1}\ldots u_{j_ki_k}\\
&=&\sum_{j_1\ldots j_k}\ \sum_{\pi,\nu\in D(k)}\delta_\pi(j)\delta_\nu(i)W_{kN}(\pi,\nu)x_{j_1}\ldots x_{j_k}\\
&=&\sum_{\pi,\nu\in D(k)}\delta_\nu(i)W_{kN}(\pi,\nu)\sum_{j_1\ldots j_k}\delta_\pi(j)x_{j_1}\ldots x_{j_k}\\
&=&\sum_{\pi,\nu\in D(k)}\delta_\nu(i)W_{kN}(\pi,\nu)
\end{eqnarray*}

Now by comparing with the formula in Theorem 12.8, we obtain the result. As for the proof in the complex case, this is similar, by adding exponents everywhere.
\end{proof} 

Still following \cite{bgo}, we have as well the following useful abstract result:

\index{integration over spheres}

\begin{theorem}
There is a unique positive unital trace $tr:C(S)\to\mathbb C$ satisfying
$$(tr\otimes id)\alpha(x)=tr(x)1$$
where $\alpha$ is the coaction map of the corresponding quantum isometry group,
$$\alpha:C(S)\to C(S)\otimes C(U)$$
and this is the canonical integration, as constructed in Definition 12.7.
\end{theorem}

\begin{proof}
First of all, it follows from the Haar integral invariance condition for $U$ that the canonical integration has indeed the invariance property in the statement, namely:
$$(tr\otimes id)\alpha(x)=tr(x)1$$

In order to prove now the uniqueness, let $tr$ be as in the statement. We have:
$$tr\left(id\otimes\int_U\right)\alpha(x)
=\int_U(tr\otimes id)\alpha(x)
=\int_U(tr(x)1)
=tr(x)$$

On the other hand, according to Theorem 12.9, we have as well:
$$
tr\left(id\otimes\int_U\right)\alpha(x)
=tr\left(\int_Sx\right)
=\int_Sx$$

We therefore conclude that $tr$ equals the standard integration, as claimed.
\end{proof}

Getting now to probabilistic aspects, we know from chapter 1 that for the classical spheres the laws of the standard coordinates, which are called hyperspherical laws, become normal and independent, with $N\to\infty$. Here is the free analogue of this result:

\begin{theorem}
For the free spheres $S^{N-1}_{\mathbb R,+}$, $S^{N-1}_{\mathbb C,+}$, the rescaled coordinates 
$$y_i=\sqrt{N}x_i$$
become semicircular/circular and free, in the $N\to\infty$ limit.
\end{theorem}

\begin{proof}
In the real case, the Weingarten formula from Theorem 12.8, coupled with the standard fact that the Gram and Weingarten matrices are asymptotically diagonal, $G_{kN}\simeq diag(N^{k/2})$ and $W_{kN}\simeq diag(N^{-k/2})$, gives the following estimate:
$$\int_{S^{N-1}_{\mathbb R,+}}x_{i_1}\ldots x_{i_k}
=\sum_{\pi,\nu\in NC_2(k)}\delta_\nu(i)W_{kN}(\pi,\nu)
\simeq N^{-k/2}\sum_{\nu\in NC_2(k)}\delta_\nu(i)$$

With this formula in hand, we can compute the asymptotic moments of each coordinate $x_i$. Indeed, by setting $i_1=\ldots=i_k=i$, all Kronecker symbols are 1, and we obtain:
$$\int_{S^{N-1}_{\mathbb R,+}}x_i^k\,dx\simeq N^{-k/2}|NC_2(k)|=N^{-k/2}M_k(\gamma_1)$$

Thus the rescaled coordinates $y_i=\sqrt{N}x_i$ become semicircular in the $N\to\infty$ limit, and the asymptotic freeness result follows as well, from the above general joint moment estimate. As for the proof in the complex case, this is similar. See \cite{bgo}.
\end{proof}

With spheres discussed, what is next? Tori I guess, these being certainly the simplest manifolds coming after spheres, or perhaps even before them. But here, we are right away into a difficulty, because the tori are compact quantum groups, which are not easy, so we can expect the probability computations for them to deviate from easiness, and Bercovici-Pata. And this is indeed the case, the situation being as follows:

\index{Meixner laws}
\index{free Meixner laws}
\index{free torus}

\begin{theorem}
In the context of the classical and free, real and complex tori, which are all compact quantum groups, appearing as duals of discrete groups,
$$\xymatrix@R=16.7mm@C=16.7mm{
T_N^+\ar[r]&\mathbb T_N^+\\
T_N\ar[r]\ar[u]&\mathbb T_N\ar[u]}
\qquad\xymatrix@R=8.8mm@C=15mm{\\ =\\ }\qquad
\xymatrix@R=15mm@C=15mm{
\widehat{\mathbb Z_2^{*N}}\ar[r]&\widehat{\mathbb Z^{*N}}\\
\widehat{\mathbb Z_2^N}\ar[r]\ar[u]&\widehat{\mathbb Z^N}\ar[u]
}$$
when looking at laws of truncated characters, the correspondence on the vertical is not Bercovici-Pata, but Meixner/free Meixner on the right, and a version of it, on the left.
\end{theorem}

\begin{proof}
This is something quite self-explanatory, the idea being as follows:

\medskip

(1) Regarding our first claim, that the 4 main tori are the group duals in the statement, this is something quite clear, only requiring a bit of geometric thinking.

\medskip

(2) Getting now to probability computations, the truncated characters for the tori $T=\widehat{\Gamma}$, with $\Gamma=<g_1,\ldots,g_N>$ being a discrete group, are as follows:
$$\chi_t=g_1+\ldots+g_{[tN]}$$

(3) But this leads us to the Meixner laws for $\Gamma=\mathbb Z^N$, to the free Meixner laws for $\Gamma=\mathbb Z^{*N}$, and to some versions of these laws for $\Gamma=\mathbb Z_2^N$ and $\Gamma=\mathbb Z_2^{*N}$.

\medskip

(4) The continuation of the story, however, is more complicated, the idea being that the Meixner/free Meixner correspondence can be axiomatized, and with this substantially deviating from Bercovici-Pata. Things here are quite technical, and we refer to the literature, with a good reference, to start with, being the paper of Anshelevich \cite{ans}.
\end{proof}

\section*{12b. Free manifolds}

As a next objective, in order to get deeper into free geometry, we would like to find a suitable collection of ``free homogeneous spaces'', generalizing at the same time the free spheres $S$, and the free unitary groups $U$. This can be done at several levels of generality, and central here is the construction of the free spaces of partial isometries, which can be done in fact for any easy quantum group. In order to explain this, let us start with:

\begin{definition}
Associated to any integers $L\leq M,N$ are the spaces
$$O_{MN}^L=\left\{T:E\to F\ {\rm isometry}\Big|E\subset\mathbb R^N,F\subset\mathbb R^M,\dim_\mathbb RE=L\right\}$$
$$U_{MN}^L=\left\{T:E\to F\ {\rm isometry}\Big|E\subset\mathbb C^N,F\subset\mathbb C^M,\dim_\mathbb CE=L\right\}$$
where the notion of isometry is with respect to the usual real/complex scalar products.
\end{definition}

As a first observation, at $L=M=N$ we obtain the groups $O_N,U_N$:
$$O_{NN}^N=O_N\quad,\quad 
U_{NN}^N=U_N$$ 

Another interesting specialization is $L=M=1$. Here the elements of $O_{1N}^1$ are the isometries $T:E\to\mathbb R$, with $E\subset\mathbb R^N$ one-dimensional. But such an isometry is uniquely determined by $T^{-1}(1)\in\mathbb R^N$, which must belong to $S^{N-1}_\mathbb R$. Thus, we have $O_{1N}^1=S^{N-1}_\mathbb R$. Similarly, in the complex case we have $U_{1N}^1=S^{N-1}_\mathbb C$, and so our results here are:
$$O_{1N}^1=S^{N-1}_\mathbb R\quad,\quad 
U_{1N}^1=S^{N-1}_\mathbb C$$

Yet another interesting specialization is $L=N=1$. Here the elements of $O_{1N}^1$ are the isometries $T:\mathbb R\to F$, with $F\subset\mathbb R^M$ one-dimensional. But such an isometry is uniquely determined by $T(1)\in\mathbb R^M$, which must belong to $S^{M-1}_\mathbb R$. Thus, we have $O_{M1}^1=S^{M-1}_\mathbb R$. Similarly, in the complex case we have $U_{M1}^1=S^{M-1}_\mathbb C$, and so our results here are:
$$O_{M1}^1=S^{M-1}_\mathbb R\quad,\quad
U_{M1}^1=S^{M-1}_\mathbb C$$

In general, the most convenient is to view the elements of $O_{MN}^L,U_{MN}^L$ as rectangular matrices, and to use matrix calculus for their study. We have indeed:

\begin{proposition}
We have identifications of compact spaces
$$O_{MN}^L\simeq\left\{U\in M_{M\times N}(\mathbb R)\Big|UU^t={\rm projection\ of\ trace}\ L\right\}$$
$$U_{MN}^L\simeq\left\{U\in M_{M\times N}(\mathbb C)\Big|UU^*={\rm projection\ of\ trace}\ L\right\}$$
with each partial isometry being identified with the corresponding rectangular matrix.
\end{proposition}

\begin{proof}
We can indeed identify the partial isometries $T:E\to F$ with their corresponding extensions $U:\mathbb R^N\to\mathbb R^M$, $U:\mathbb C^N\to\mathbb C^M$, obtained by setting $U_{E^\perp}=0$. Then, we can identify these latter maps $U$ with the corresponding rectangular matrices.
\end{proof}

In order to advance, observe now that the isometries $T:E\to F$, or rather their extensions $U:\mathbb K^N\to\mathbb K^M$, with $\mathbb K=\mathbb R,\mathbb C$, obtained by setting $U_{E^\perp}=0$, can be composed with the isometries of $\mathbb K^M,\mathbb K^N$, according to the following scheme:
$$\xymatrix@R=17mm@C=17mm{
\mathbb K^N\ar[r]^{B^*}&\mathbb K^N\ar@.[r]^U&\mathbb K^M\ar[r]^A&\mathbb K^M\\
B(E)\ar@.[r]\ar[u]&E\ar[r]^T\ar[u]&F\ar@.[r]\ar[u]&A(F)\ar[u]
}$$

With the identifications in Proposition 12.14 made, the precise statement here is:

\begin{proposition}
We have action maps as follows, which are both transitive,
$$O_M\times O_N\curvearrowright O_{MN}^L\quad,\quad 
(A,B)U=AUB^t$$
$$U_M\times U_N\curvearrowright U_{MN}^L\quad,\quad 
(A,B)U=AUB^*$$
whose stabilizers are respectively $O_L\times O_{M-L}\times O_{N-L}$ and $U_L\times U_{M-L}\times U_{N-L}$.
\end{proposition}

\begin{proof}
We have indeed action maps as in the statement, which are transitive. Let us compute now the stabilizer $G$ of the following point:
$$U=\begin{pmatrix}1&0\\0&0\end{pmatrix}$$

Since $(A,B)\in G$ satisfy $AU=UB$, their components must be of the following form:
$$A=\begin{pmatrix}x&*\\0&a\end{pmatrix}\quad,\quad 
B=\begin{pmatrix}x&0\\ *&b\end{pmatrix}$$

Now since $A,B$ are unitaries, these matrices follow to be block-diagonal, and so:
$$G=\left\{(A,B)\Big|A=\begin{pmatrix}x&0\\0&a\end{pmatrix},B=\begin{pmatrix}x&0\\ 0&b\end{pmatrix}\right\}$$

The stabilizer of $U$ is parametrized by triples $(x,a,b)$ belonging to $O_L\times O_{M-L}\times O_{N-L}$ and $U_L\times U_{M-L}\times U_{N-L}$, and we are led to the conclusion in the statement.
\end{proof}

Finally, let us work out the quotient space description of $O_{MN}^L,U_{MN}^L$. We have here:

\begin{theorem}
We have isomorphisms of homogeneous spaces as follows,
\begin{eqnarray*}
O_{MN}^L&=&(O_M\times O_N)/(O_L\times O_{M-L}\times O_{N-L})\\
U_{MN}^L&=&(U_M\times U_N)/(U_L\times U_{M-L}\times U_{N-L})
\end{eqnarray*}
with the quotient maps being given by $(A,B)\to AUB^*$, where $U=(^1_0{\ }^0_0)$.
\end{theorem}

\begin{proof}
This is just a reformulation of Proposition 12.15, by taking into account the fact that the fixed point used in the proof there was $U=(^1_0{\ }^0_0)$.
\end{proof}

Summarizing, we have here some basic homogeneous spaces, unifying the spheres with the rotation groups. The point now is that we can liberate these spaces, as follows:

\begin{definition}
Associated to any integers $L\leq M,N$ are the algebras
\begin{eqnarray*}
C(O_{MN}^{L+})&=&C^*\left((v_{ij})_{i=1,\ldots,M,j=1,\ldots,N}\Big|v=\bar{v},vv^t={\rm projection\ of\ trace}\ L\right)\\
C(U_{MN}^{L+})&=&C^*\left((v_{ij})_{i=1,\ldots,M,j=1,\ldots,N}\Big|vv^*,\bar{v}v^t={\rm projections\ of\ trace}\ L\right)
\end{eqnarray*}
with the trace being by definition the sum of the diagonal entries.
\end{definition}

Observe that we have inclusions between our various spaces, as follows:
$$\xymatrix@R=15mm@C=15mm{
O_{MN}^{L+}\ar[r]&U_{MN}^{L+}\\
O_{MN}^L\ar[r]\ar[u]&U_{MN}^L\ar[u]}$$

At the level of basic examples now, at $L=M=1$ and at $L=N=1$ we obtain the following diagrams, showing that our formalism covers indeed the free spheres:
$$\xymatrix@R=15mm@C=15mm{
S^{N-1}_{\mathbb R,+}\ar[r]&S^{N-1}_{\mathbb C,+}\\
S^{N-1}_\mathbb R\ar[r]\ar[u]&S^{N-1}_\mathbb C\ar[u]}
\qquad\qquad 
\xymatrix@R=15mm@C=15mm{
S^{M-1}_{\mathbb R,+}\ar[r]&S^{M-1}_{\mathbb C,+}\\
S^{M-1}_\mathbb R\ar[r]\ar[u]&S^{M-1}_\mathbb C\ar[u]}$$

We have as well the following result, in relation with the free rotation groups:

\begin{proposition}
At $L=M=N$ we obtain the diagram
$$\xymatrix@R=15mm@C=15mm{
O_N^+\ar[r]&U_N^+\\
O_N\ar[r]\ar[u]&U_N\ar[u]}$$
consisting of the groups $O_N,U_N$, and their liberations.
\end{proposition}

\begin{proof}
According to the above, we have the following presentation results:
\begin{eqnarray*}
C(O_{NN}^{N\times})&=&C^*_\times\left((v_{ij})_{i,j=1,\ldots,N}\Big|v=\bar{v},vv^t={\rm projection\ of\ trace}\ N\right)\\
C(U_{NN}^{N\times})&=&C^*_\times\left((v_{ij})_{i,j=1,\ldots,N}\Big|vv^*,\bar{v}v^t={\rm projections\ of\ trace}\ N\right)
\end{eqnarray*}

We use now the standard fact that if $p=aa^*$ is a projection then $q=a^*a$ is a projection too. We use as well the following formulae:
$$Tr(vv^*)=Tr(v^t\bar{v})\quad,\quad 
Tr(\bar{v}v^t)=Tr(v^*v)$$

We therefore obtain the following formulae:
\begin{eqnarray*}
C(O_{NN}^{N\times})&=&C^*_\times\left((v_{ij})_{i,j=1,\ldots,N}\Big|v=\bar{v},\ vv^t,v^tv={\rm projections\ of\ trace}\ N\right)\\
C(U_{NN}^{N\times})&=&C^*_\times\left((v_{ij})_{i,j=1,\ldots,N}\Big|vv^*,v^*v,\bar{v}v^t,v^t\bar{u}={\rm projections\ of\ trace}\ N\right)
\end{eqnarray*}

Now observe that, in tensor product notation, the conditions at right are all of the form $(tr\otimes id)p=1$. Thus, $p$ must be follows, for the above conditions:
$$p=vv^*,v^*v,\bar{v}v^t,v^t\bar{v}$$

We therefore obtain that, for any faithful state $\varphi$, we have $(tr\otimes\varphi)(1-p)=0$. It follows from this that the following projections must be all equal to the identity:
$$p=vv^*,v^*v,\bar{v}v^t,v^t\bar{v}$$

But this leads to the conclusion in the statement.
\end{proof}

Regarding now the homogeneous space structure of $O_{MN}^{L\times},U_{MN}^{L\times}$, the situation here is a bit more complicated in the free case than in the classical case, due to a number of algebraic and analytic issues. We first have the following result:

\begin{proposition}
The spaces $U_{MN}^{L\times}$ have the following properties:
\begin{enumerate}
\item We have an action $U_M^\times\times U_N^\times\curvearrowright U_{MN}^{L\times}$, given by $v_{ij}\to\sum_{kl}v_{kl}\otimes a_{ki}\otimes b_{lj}^*$.

\item We have a map $U_M^\times\times U_N^\times\to U_{MN}^{L\times}$, given by $v_{ij}\to\sum_{r\leq L}a_{ri}\otimes b_{rj}^*$.
\end{enumerate}
Similar results hold for the spaces $O_{MN}^{L\times}$, with all the $*$ exponents removed.
\end{proposition}

\begin{proof}
In the classical case, consider the following action and quotient maps:
$$U_M\times U_N\curvearrowright U_{MN}^L\quad,\quad 
U_M\times U_N\to U_{MN}^L$$

The transposes of these two maps are as follows, where $J=(^1_0{\ }^0_0)$:
\begin{eqnarray*}
\varphi&\to&((U,A,B)\to\varphi(AUB^*))\\
\varphi&\to&((A,B)\to\varphi(AJB^*))
\end{eqnarray*}

But with $\varphi=v_{ij}$ we obtain precisely the formulae in the statement. The proof in the orthogonal case is similar. Regarding now the free case, the proof goes as follows:

\medskip

(1) Assuming $vv^*v=v$, let us set $U_{ij}=\sum_{kl}v_{kl}\otimes a_{ki}\otimes b_{lj}^*$. We have then:
\begin{eqnarray*}
(UU^*U)_{ij}
&=&\sum_{pq}\sum_{klmnst}v_{kl}v_{mn}^*v_{st}\otimes a_{ki}a_{mq}^*a_{sq}\otimes b_{lp}^*b_{np}b_{tj}^*\\
&=&\sum_{klmt}v_{kl}v_{ml}^*v_{mt}\otimes a_{ki}\otimes b_{tj}^*\\
&=&\sum_{kt}v_{kt}\otimes a_{ki}\otimes b_{tj}^*\\
&=&U_{ij}
\end{eqnarray*}

Also, assuming that we have $\sum_{ij}v_{ij}v_{ij}^*=L$, we obtain:
$$\sum_{ij}U_{ij}U_{ij}^*
=\sum_{ij}\sum_{klst}v_{kl}v_{st}^*\otimes a_{ki}a_{si}^*\otimes b_{lj}^*b_{tj}
=\sum_{kl}v_{kl}v_{kl}^*\otimes1\otimes1
=L$$

(2) Assuming $vv^*v=v$, let us set $V_{ij}=\sum_{r\leq L}a_{ri}\otimes b_{rj}^*$. We have then:
$$(VV^*V)_{ij}
=\sum_{pq}\sum_{x,y,z\leq L}a_{xi}a_{yq}^*a_{zq}\otimes b_{xp}^*b_{yp}b_{zj}^*
=\sum_{x\leq L}a_{xi}\otimes b_{xj}^*
=V_{ij}$$

Finally, assuming that we have $\sum_{ij}u_{ij}u_{ij}^*=L$, we obtain:
$$\sum_{ij}V_{ij}V_{ij}^*
=\sum_{ij}\sum_{r,s\leq L}a_{ri}a_{si}^*\otimes b_{rj}^*b_{sj}
=\sum_{l\leq L}1
=L$$

By removing all the $*$ exponents, we obtain as well the orthogonal results.
\end{proof}

Let us examine now the relation between the above maps. In the classical case, given a quotient space $X=G/H$, the associated action and quotient maps are given by:
$$\begin{cases}
a:X\times G\to X&:\quad (Hg,h)\to Hgh\\
p:G\to X&:\quad g\to Hg
\end{cases}$$

Thus we have $a(p(g),h)=p(gh)$. In our context, a similar result holds: 

\begin{theorem}
With $G=G_M\times G_N$ and $X=G_{MN}^L$, where $G_N=O_N^\times,U_N^\times$, we have
$$\xymatrix@R=15mm@C=30mm{
G\times G\ar[r]^m\ar[d]_{p\times id}&G\ar[d]^p\\
X\times G\ar[r]^a&X
}$$
where $a,p$ are the action map and the map constructed in Proposition 12.19.
\end{theorem}

\begin{proof}
At the level of the associated algebras of functions, we must prove that the following diagram commutes, where $\Phi,\alpha$ are morphisms of algebras induced by $a,p$:
$$\xymatrix@R=15mm@C=25mm{
C(X)\ar[r]^\Phi\ar[d]_\alpha&C(X\times G)\ar[d]^{\alpha\otimes id}\\
C(G)\ar[r]^\Delta&C(G\times G)
}$$

When going right, and then down, the composition is as follows:
$$(\alpha\otimes id)\Phi(u_{ij})
=(\alpha\otimes id)\sum_{kl}v_{kl}\otimes a_{ki}\otimes b_{lj}^*
=\sum_{kl}\sum_{r\leq L}a_{rk}\otimes b_{rl}^*\otimes a_{ki}\otimes b_{lj}^*$$

On the other hand, when going down, and then right, the composition is as follows, where $F_{23}$ is the flip between the second and the third components:
$$\Delta\pi(u_{ij})
=F_{23}(\Delta\otimes\Delta)\sum_{r\leq L}a_{ri}\otimes b_{rj}^*
=F_{23}\left(\sum_{r\leq L}\sum_{kl}a_{rk}\otimes a_{ki}\otimes b_{rl}^*\otimes b_{lj}^*\right)$$

Thus the above diagram commutes indeed, and this gives the result.
\end{proof}

Let us discuss now the integration over the above spaces $G_{MN}^L$. We first have:

\begin{definition}
The integration functional of $G_{MN}^L$ is the composition
$$\int_{G_{MN}^L}:C(G_{MN}^L)\to C(G_M\times G_N)\to\mathbb C$$
of the representation $v_{ij}\to\sum_{r\leq L}a_{ri}\otimes b_{rj}^*$ with the Haar functional of $G_M\times G_N$.
\end{definition}

As an illustration here, observe that in the case $L=M=N$ we obtain the integration over $G_N$. Also, at $L=M=1$, or at $L=N=1$, we obtain the integration over the sphere. In the general case now, we first have the following result:

\begin{proposition}
The integration functional of $G_{MN}^L$ has the invariance property 
$$\left(\int_{G_{MN}^L}\!\otimes\ id\right)\Phi(x)=\int_{G_{MN}^L}x$$
with respect to the coaction map $\Phi(v_{ij})=\sum_{kl}v_{kl}\otimes a_{ki}\otimes b_{lj}^*$.
\end{proposition}

\begin{proof}
This comes indeed from a straightforward computation, based on the invariance property of the Haar functionals of $G_M,G_N$, that we will leave as an exercise.
\end{proof}

Next, in analogy with what we did for the spheres, we have a Weingarten formula:

\index{Weingarten formula}

\begin{theorem}
We have the Weingarten type formula
$$\int_{G_{MN}^L}v_{i_1j_1}\ldots v_{i_sj_s}=\sum_{\pi\sigma\tau\nu}L^{|\pi\vee\tau|}\delta_\sigma(i)\delta_\nu(j)W_{sM}(\pi,\sigma)W_{sN}(\tau,\nu)$$
where the matrices on the right are given by $W_{sM}=G_{sM}^{-1}$, with $G_{sM}(\pi,\sigma)=M^{|\pi\vee\sigma|}$.
\end{theorem}

\begin{proof}
By using the Weingarten formula for $G_M,G_N$, we obtain:
\begin{eqnarray*}
\int_{G_{MN}^L}v_{i_1j_1}\ldots v_{i_sj_s}
&=&\sum_{l_1\ldots l_s\leq L}\int_{G_M}a_{l_1i_1}\ldots a_{l_si_s}\int_{G_N}b_{l_1j_1}^*\ldots b_{l_sj_s}^*\\
&=&\sum_{l_1\ldots l_s\leq L}\sum_{\pi\sigma}\delta_\pi(l)\delta_\sigma(i)W_{sM}(\pi,\sigma)\sum_{\tau\nu}\delta_\tau(l)\delta_\nu(j)W_{sN}(\tau,\nu)\\
&=&\sum_{\pi\sigma\tau\nu}\left(\sum_{l_1\ldots l_s\leq L}\delta_\pi(l)\delta_\tau(l)\right)\delta_\sigma(i)\delta_\nu(j)W_{sM}(\pi,\sigma)W_{sN}(\tau,\nu)
\end{eqnarray*}

The coefficient being $L^{|\pi\vee\tau|}$, we obtain the formula in the statement.
\end{proof}

We can now derive an abstract characterization of the integration, as follows:

\begin{theorem}
The integration of $G_{MN}^L$ is the unique positive unital trace 
$$C(G_{MN}^L)\to\mathbb C$$
which is invariant under the action of the quantum group $G_M\times G_N$.
\end{theorem}

\begin{proof}
This is indeed very standard, exactly as before for the spheres.
\end{proof}

As a main application of the above results, we have the following quite conceptual statement, making the link with the Bercovici-Pata bijection \cite{bpa}:

\index{Bercovici-Pata bijection}
\index{non-overlapping coordinates}

\begin{theorem}
In the context of the liberation operations $G_{MN}^L\to G_{MN}^{L+}$, the laws of the sums of non-overlapping coordinates,
$$\chi_E=\sum_{(ij)\in E}u_{ij}$$
are in Bercovici-Pata bijection, in the $|E|=\kappa N,L=\lambda N,M=\mu N$ and $N\to\infty$ regime.
\end{theorem}

\begin{proof}
With $K=|E|$, the moments of the variables in the statement are:
\begin{eqnarray*}
M_s
&=&\sum_{\pi\sigma\tau\nu}K^{|\pi\vee\tau|}L^{|\sigma\vee\nu|}W_{sM}(\pi,\sigma)W_{sN}(\tau,\nu)\\
&\simeq&\sum_{\pi\tau}K^{|\pi\vee\tau|}L^{|\pi\vee\tau|}M^{-|\pi|}N^{-|\tau|}\\
&\simeq&\sum_\pi K^{|\pi|}L^{|\pi|}M^{-|\pi|}N^{-|\pi|}\\
&=&\sum_\pi\left(\frac{\kappa\lambda}{\mu}\right)^{|\pi|}
\end{eqnarray*}

But this gives the result, via our standard moment formulae for the main laws.
\end{proof}

And we will end our discussion of free manifolds with this. It is possible to substantially build on the above, and for more on this, have a look at my quantum algebra books.

\section*{12c. Hyperspherical laws}

Getting back to spheres, an interesting problem is that of computing the free hyperspherical laws. Let us begin with a discussion in the classical case. Here we certainly know about hyperspherical laws from chapter 1, and as a continuation, we have:

\begin{theorem}
The moments and density of the hyperspherical laws, that of the coordinates $x_i:S^{N-1}_\mathbb R\to[-1,1]$, are as follows, with $e_N=2,\pi$ for $N$ odd, even,
$$M_{2r}=\frac{(N-1)!!(2r)!!}{(N+2r-1)!!}\quad,\quad\mu_N=\frac{(N-1)!!}{e_N(N-2)!!}(1-x^2)^{(N-3)/2}dx$$
with the corresponding expectation, variance, skewness and kurtosis being given by:
$$E=0\quad,\quad V=\frac{1}{N}\quad,\quad\gamma=0\quad,\quad\kappa=\frac{3N}{N+2}$$
At $N=1,2,3,4$ we obtain the following well-known laws, all supported by $[-1,1]$:
$$\mu_1=\frac{\delta_{-1}+\delta_1}{2}\quad,\quad\mu_2=\frac{1}{\pi\sqrt{1-x^2}}\,dx
\quad,\quad\mu_3=\frac{dx}{2}\quad,\quad\mu_4=\frac{2}{\pi}\sqrt{1-x^2}\,dx$$
Also, the joint moments of the variables $x_i:S^{N-1}_\mathbb R\to[-1,1]$ are given by
$$E(x_1^{2r_1}\ldots x_N^{2r_N})=\frac{(N-1)!!(2r_1)!!\ldots(2r_N)!!}{(N+\sum 2r_i-1)!!}$$
and $y_i=\sqrt{N}x_i$ become standard normal, following $g_1$, and independent with $N\to\infty$.
\end{theorem}

\begin{proof}
We already know the joint moment formula, from chapter 1, and all the rest comes from this, via standard calculus. For details, you can check my book \cite{ba2}.
\end{proof}

Next, we have the following result, referring to the beta laws from chapter 10:

\index{beta laws}
\index{chi squared}

\begin{theorem}[version]
The squared coordinates $x_i^2:S^{N-1}_\mathbb R\to[0,1]$ follow beta laws, having moments and density as follows, with $e_N=2,\pi$ for $N$ odd, even,
$$M_r=\frac{(N-1)!!(2r)!!}{(N+2r-1)!!}\quad,\quad\nu_N=\frac{(N-1)!!}{e_N(N-2)!!}\,x^{-1/2}(1-x)^{(N-3)/2}dx$$
the expectation being $E=1/N$, and the variance, skewness and kurtosis being:
$$V=\frac{2N-2}{N^2(N+2)}\ ,\quad 
\gamma=\frac{N-2}{N+4}\sqrt{\frac{8(N+2)}{N-1}}\ ,\quad 
\kappa=\frac{3(N+2)(5N^2-13N+12)}{(N-1)(N+4)(N+6)}$$
At $N=1,2,3,4$ we obtain the following well-known laws, all supported by $[0,1]$:
$$\nu_1=\delta_1\quad,\quad\nu_2=\frac{1}{\pi\sqrt{x-x^2}}\,dx
\quad,\quad\nu_3=\frac{dx}{2\sqrt{x}}\quad,\quad\nu_4=\frac{2}{\pi}\sqrt{x^{-1}-1}\,dx$$
Also, the joint moments of the variables $x_i^2:S^{N-1}_\mathbb R\to[0,1]$ are given by
$$E(x_1^{2r_1}\ldots x_N^{2r_N})=\frac{(N-1)!!(2r_1)!!\ldots(2r_N)!!}{(N+\sum 2r_i-1)!!}$$
and $y_i=Nx_i^2$ become chi squared, following $\chi_1^2$, and independent with $N\to\infty$.
\end{theorem}

\begin{proof}
This is a straightforward remake of Theorem 12.26, see \cite{ba2}.
\end{proof}

Finally, our love for calculus leads us to a third result, as follows:

\index{chi variable}

\begin{theorem}[version]
The moments and density of the absolute values of the coordinates $|x_i|:S^{N-1}_\mathbb R\to[0,1]$ are as follows, with $e_N=2,\pi$ for $N$ odd, even,
$$M_k=\left(\frac{2}{\pi}\right)^{\delta_{2|N}\delta_{2\not\:|\,k}}\frac{(N-1)!!k!!}{(N+k-1)!!}\quad,\quad\eta_N=\frac{2(N-1)!!}{e_N(N-2)!!}(1-x^2)^{(N-3)/2}dx$$
and at $N=1,2,3,4$ we obtain the following well-known laws, all supported by $[0,1]$:
$$\mu_1=\delta_1\quad,\quad\mu_2=\frac{2}{\pi\sqrt{1-x^2}}\,dx
\quad,\quad\mu_3=dx\quad,\quad\mu_4=\frac{4}{\pi}\sqrt{1-x^2}\,dx$$
The joint moments are as follows, with $\Sigma=[odds/2]$ if $N$ is odd and $\Sigma=[(odds+1)/2]$ if $N$ is even, where ``odds'' is the number of odd numbers in the sequence $k_1,\ldots,k_N$,
$$E(|x_1|^{k_1}\ldots|x_N|^{k_N})=\left(\frac{2}{\pi}\right)^{\Sigma(k_1,\ldots,k_N)}\frac{(N-1)!!k_1!!\ldots k_N!!}{(N+\sum k_i-1)!!}$$
and $y_i=\sqrt{N}|x_i|$ become standard chi, following $\chi_1$, and independent with $N\to\infty$.
\end{theorem}

\begin{proof}
This is something more tricky, needing a remake of the Wallis formulae from chapter 1, using absolute values. For details on this, have a look a my book \cite{ba2}. 
\end{proof}

Getting now to the free real spheres, the result here, from \cite{bcz}, based on advanced quantum group and calculus techniques, is something quite tricky, as follows:

\index{free hyperspherical law}
\index{special functions}
\index{deformed quantum groups}

\begin{theorem}
The moments of the free hyperspherical law are given by
$$\int_{S^{N-1}_{\mathbb R,+}}x_1^{2l}=\frac{1}{(N+2)^l}\cdot\frac{q+1}{q-1}\cdot\frac{1}{l+1}\sum_{r=-l-1}^{l+1}(-1)^r\binom{2l+2}{l+r+1}\frac{r}{1+q^r}$$
where $q\in [-1,0)$ is such that $q+q^{-1}=-N$.
\end{theorem}

\begin{proof}
Not a trivial result that we have here, the idea being as follows:

\medskip

(1) In order to get into the mood, let us first see what our formula tells us, at $l=1$. And here, with the dots standing for some computations, what we get is:
$$M_2
=\frac{1}{N+2}\cdot\frac{q+1}{q-1}\cdot\frac{1}{2}\sum_{r=-2}^2(-1)^r\binom{4}{r+2}\frac{r}{1+q^r}=\ldots\ldots=\frac{1}{N}$$

Which is not particularly fun, hope you agree with me. Note by the way that this is the correct answer, because we have $x_1^2+\ldots+x_N^2=1$, uniformly.

\medskip

(2) At $l=2$ now, the computation becomes particularly tough, but still can be done by hand. Then at $l=3$, things complicated. In any case, here are the formulae:
$$M_4=\frac{2}{N(N+1)}\quad,\quad 
M_6=\frac{5N-7}{N(N+1)(N^2-2)}$$

And, are these correct. In answer, the Weingarten formula in Theorem 12.8 gives:
$$M_4=Tr\left[\begin{pmatrix}1&1\\1&1\end{pmatrix}
\begin{pmatrix}N^2&N\\N&N^2\end{pmatrix}^{-1}\right]
=\frac{2}{N(N+1)}$$
$$M_6=Tr\left[\begin{pmatrix}
1&1&1&1&1\\
1&1&1&1&1\\
1&1&1&1&1\\
1&1&1&1&1\\
1&1&1&1&1
\end{pmatrix}
\begin{pmatrix}
N^3&N^2&N^2&N&N^2\\
N^2&N^3&N&N^2&N\\
N^2&N&N^3&N^2&N\\
N&N^2&N^2&N^3&N^2\\
N^2&N&N&N^2&N^3
\end{pmatrix}^{-1}\right]
=\frac{5N-7}{N(N+1)(N^2-2)}$$

(3) Summarizing, our formula looks good. Observe also that, as a byproduct of our computations, the expectation, variance, skewness and kurtosis are:
$$E=0\quad,\quad V=\frac{1}{N}\quad,\quad\gamma=0\quad,\quad\kappa=\frac{2N}{N+1}$$

As a comment here, $E,V,\gamma$ equal those for the classical hyperspherical laws, from Theorem 12.26, but the kurtosis, which was $\kappa=3N/(N+2)$, has changed. Which is not surprising, because we should have $\kappa\simeq3$ for classical laws, and $\kappa\simeq2$ for free laws.

\medskip

(4) Finally, before getting into the proof, let us mention that \cite{bcz} contains as well a number of useful reformulations and corollaries of the formula in the statement, notably with the fact that the support is as follows, with the density being analytic:
$$S=\left[-\frac{2}{\sqrt{N+2}}\ ,\ \frac{2}{\sqrt{N+2}}\right]$$ 

(5) Getting now to the proof, the idea will be that of transporting our computation from space to space, no less than 4 times, up to reaching a tough computation on the unit circle $\mathbb T$, and solving that with some help from advanced orthogonal polynomials:
$$S^{N-1}_{\mathbb R,+}
\ \longrightarrow\ O_N^+ 
\ \longrightarrow\ O_F^+=SU_2^q
\ \longrightarrow\ B(l^2(\mathbb N))
\ \longrightarrow\ \mathbb T$$

(6) Regarding the first transport, via $x_1=u_{11}$, we know this, from the previous section. Next, to any $F\in GL_N(\mathbb R)$ satisfying $F^2=1$ let us associate the following algebra:
$$C(O_F^+)=C^*\left((u_{ij})_{i,j=1,\ldots,N}\Big|u=F\bar{u}F={\rm unitary}\right)$$

Observe that we have $O_{I_N}^+=O_N^+$. In general, the above algebra satisfies Woronowicz's generalized axioms in \cite{wo1}, which do not include the strong antipode axiom $S^2=id$. 

\medskip

(7) At $N=2$, up to a trivial equivalence relation on the matrices $F$, and on the quantum groups $O_F^+$, we can assume that $F$ is as follows, with $q\in [-1,0)$:
$$F=\begin{pmatrix}0&\sqrt{-q}\\
1/\sqrt{-q}&0\end{pmatrix}$$

Our claim is that for this matrix $F$ we have an identification $O_F^+=SU^q_2$. Indeed, the relations $u=F\bar{u}F$ tell us that $u$ must be of the following form:
$$u=\begin{pmatrix}\alpha&-q\gamma^*\cr \gamma&\alpha^*\end{pmatrix}$$

Thus $C(O_F^+)$ is the universal algebra generated by two elements $\alpha,\gamma$, with the relations making the above matrix $u$ unitary. But these unitarity conditions are as follows:
$$\alpha\gamma=q\gamma\alpha\quad,\quad
\alpha\gamma^*=q\gamma^*\alpha\quad,\quad 
\gamma\gamma^*=\gamma^*\gamma$$
$$\alpha^*\alpha+\gamma^*\gamma=1\quad,\quad 
\alpha\alpha^*+q^2\gamma\gamma^*=1$$

We recognize here the relations in \cite{wo1} defining the algebra $C(SU^q_2)$, and it follows that we have an isomorphism of Hopf $C^*$-algebras, as follows:
$$C(O_F^+)\simeq C(SU^q_2)$$

(8) Now back to the general case, let us try to understand the integration over $O_F^+$. Given a noncrossing pairing $\pi\in NC_2(2k)$ and an index $i=(i_1,\ldots,i_{2k})$, we set:
$$\delta_\pi^F(i)=\prod_{s\in\pi}F_{i_{s_l}i_{s_r}}$$

Here the product is over all strings $s=\{s_l\curvearrowright s_r\}$ of $\pi$. Our claim is that the following family of vectors, with $\pi\in NC_2(2k)$, spans the space of fixed vectors of $u^{\otimes 2k}$:
$$\xi_\pi=\sum_i\delta_\pi^F(i)e_{i_1}\otimes\ldots\otimes e_{i_{2k}}$$ 

Indeed, having $\xi_\cap$ fixed by $u^{\otimes 2}$ is equivalent to assuming that $u=F\bar{u}F$ is unitary. By using now the above vectors, we obtain the following Weingarten formula:
$$\int_{O_F^+}u_{i_1j_1}\ldots u_{i_{2k}j_{2k}}=\sum_{\pi\sigma}\delta_\pi^F(i)\delta_\sigma^F(j)W_{kN}(\pi,\sigma)$$

(9) With these preliminaries in hand, let us start the computation. Let $q\in [-1,0)$ be such that $q+q^{-1}=-N$. Our claim is that we have the following formula:
$$\int_{O_N^+}\varphi(\sqrt{N+2}\,u_{ij})=\int_{SU^q_2}\varphi(\alpha+\alpha^*+\gamma-q\gamma^*)$$

Indeed, the moments of the variable on the left are given by:
$$\int_{O_N^+}u_{ij}^{2k}=\sum_{\pi\sigma}W_{kN}(\pi,\sigma)$$

On the other hand, the moments of the variable on the right, which in terms of the fundamental corepresentation $v=(v_{ij})$ is given by $w=\sum_{ij}v_{ij}$, are given by:
$$\int_{SU^q_2}w^{2k}=\sum_{ij}\sum_{\pi\sigma}\delta_\pi^F(i)\delta_\sigma^F(j)W_{kN}(\pi,\sigma)$$

We deduce that $w/\sqrt{N+2}$ has the same moments as $u_{ij}$, which proves our claim.

\medskip

(10) Next, before jumping into computations over $SU^q_2$, it is convenient to introduce a useful degree of flexibility, coming from an extra parameter. The point indeed is that the following variables over $SU^q_2$ have the same law, whenever $AB=-q$:
$$\alpha+\alpha^*+\gamma-q\gamma^*\sim \alpha+\alpha^*+A\gamma+B\gamma^*$$

Indeed, this comes via Weingarten, applied as in (9). Thus, as a conclusion, the law that we want to compute is that of the following variable, with  $AB=-q$:
$$\sqrt{N+2}\,u_{ij}\sim \alpha+\alpha^*+A\gamma+B\gamma^*$$

(11) In order to do now the computation over $SU^q_2$, we can use a matrix model due to Woronowicz \cite{wo1}, where the standard generators $\alpha,\gamma$ are mapped as follows:
$$\pi_u(\alpha)e_k=\sqrt{1-q^{2k}}e_{k-1}\quad,\quad 
\pi_u(\gamma)e_k=uq^k e_k$$

Here $u\in\mathbb T$ is a parameter, and $(e_k)$ is the standard basis of $l^2(\mathbb N)$. The point with this representation is that it allows the computation of the Haar functional. Indeed, if $D$ is the diagonal operator given by $D(e_k)=q^{2k}e_k$, then the formula is as follows:
$$\int _{SU^q_2}x=(1-q^2)\int_{\mathbb T}tr(D\pi_u(x))\frac{du}{2\pi iu}$$

(12) Thus, the law of the variable that we are interested in is of the following form:
$$\int_{SU^q_2}\varphi(\alpha+\alpha^*+A\gamma+B\gamma^*)=(1-q^2)\int_{\mathbb T}tr(D\varphi(M))\frac{du}{2\pi iu}$$

To be more precise, this formula holds indeed, with $M$ being as follows:
$$M(e_k)=\sqrt{1-q^{2k+2}}\,e_{k+1}+q^k(Au+Bu^{-1})e_k+\sqrt{1-q^{2k}}\,e_{k-1}$$

(13) Next, let us perform a change of basis. Consider another copy of $l^2(\mathbb N)$, with orthonormal basis denoted $\{e_n'\}$, and define a linear map $J:l^2(\mathbb N)\to l^2(\mathbb N)$ as follows:
$$J(e_k')=\sqrt{1-q^2}\sqrt{1-q^4}\ldots\sqrt{1-q^{2k}}\,e_k$$

In terms of $D'=J^{-1}DJ$ and $M'=J^{-1}MJ$, our formula in (12) becomes:
$$\int\varphi(\alpha+\alpha^*+A\gamma+B\gamma^*)=(1-q^2)\int_{\mathbb T}tr(D'\varphi(M'))\frac{du}{2\pi iu}$$

Since $J$ is diagonal, we have $D'(e_k')=q^{2k}e_k'$. As for $M'$, this is given by:
\begin{eqnarray*}
M'(e_k)
&=&J^{-1}MJ(e_k')\\
&=&\sqrt{1-q^2}\sqrt{1-q^4}\ldots\sqrt{1-q^{2k}}J^{-1}M(e_k)\\
&=&e_{k+1}'+q^k(Au+Bu^{-1})e_k'+(1-q^{2k})e_{k-1}'
\end{eqnarray*}

(14) As a conclusion to all this, by removing all prime signs, and by remembering as well from (10) what we wanted to do, we have a formula as follows:
$$\int_{O_N^+}\varphi(\sqrt{N+2}\,u_{ij})=(1-q^2)\int_{\mathbb T}tr(D\varphi(M))\frac{du}{2\pi iu}$$

To be more precise, this formula holds indeed, with $M$ being as follows:
$$M(e_k)=e_{k+1}+q^k(Au+Bu^{-1})e_k+(1-q^{2k})e_{k-1}$$

(15) This was for the easy part, basically algebra. The hard part, analysis, is that of computing the above integral. First, the $q$-shifted factorial is given by:
$$(a;q)_k=(1-a)(1-qa)\ldots(1-q^{k-1}a)$$

Next, the multivariable version of this $q$-shifted factorial is given by:
$$(a_1,\ldots,a_r;q)_k=(a_1;q)_k\ldots(a_r;q)_k$$

Finally, the $q$-hypergeometric series is given by the following formula:
$${}_r\varphi_s\left(\begin{matrix}{a_1,\ldots,a_r}\cr {b_1,\ldots,b_s}\end{matrix}\,;q,z\right)
=\sum_{k=0}^\infty\left( (-1)^k q^{k(k-1)/2}\right)^{s-r+1}\frac{(a_1,\ldots,a_r;q)_k}{(q,b_1,\ldots,b_s;q)_k}\,z^k$$

(16) With this, we can now introduce the Al-Salam-Chihara polynomials, depending on two parameters $a,b$, which are defined by the following formula:
$$Q_k(x)=\frac{(ab;q)_k}{a^k}\,{}_3\varphi_2\left(\begin{matrix}q^{-k},az,az^{-1}\cr ab,0\end{matrix}\,;q,q\right)$$

Here we use the convenient parameterization $2x=z+z^{-1}$. These polynomials are known to satisfy the following recurrence relation, and we refer to \cite{awi} for this:
$$2xQ_k(x)=Q_{k+1}(x)+q^k(a+b)Q_k(x)+(1-q^k)(1-abq^{k-1})Q_{k-1}(x)$$

(17) Now by getting back to (14), let us set $a=Au$ and $b=Bu^{-1}$. In terms of these new parameters $a,b$, the formula there for our tridiagonal operator $M$ reads:
$$M(e_k)=e_{k+1}+q^k(a+b)e_k+(1-q^{2k})e_{k-1}$$

Now recall from (14) that we have $AB=-q$, and so $ab=-q$. Thus the recurrence formula for the corresponding Al-Salam-Chihara polynomials is:
$$2xQ_k(x)=Q_{k+1}(x)+q^k(a+b)Q_k(x)+(1-q^{2k})Q_{k-1}(x)$$

We deduce that the collection $\{Q_k(x)\}$ of Al-Salam-Chihara polynomials evaluated at $x$ play the role of eigenvectors for the operator $M$, with corresponding eigenvalue $2x$.

\medskip

(18) And with this, we are now into orthogonal polynomials, and some 3 pages of calculations, based on the Poisson kernel for the Al-Salam-Chihara polynomials, and on some remarkable formulae from \cite{awi}, at black belt level ${}_8\varphi_7$, lead to the formula in the statement. For details on all this, and for more, we refer to \cite{bcz}. 
\end{proof}

\section*{12d. Hypergeometric laws}

Now that we studied the free hyperspherical laws, let us look at well at the free hypergeometric ones. Following \cite{bbc}, we have the following definition:

\index{free hypergeometric law}

\begin{definition}
The free hypergeometric law of parameters $(N,m,p)$ is the law of
$$\sigma_{mp}=\sum_{i=1}^m\sum_{j=1}^pu_{ij}$$
over the quantum permutation group $S_N^+$.
\end{definition}

Observe the similarity with the well-known formula for the classical hypergeometric laws, involving the usual permutation group $S_N$. We can actually exploit this similarity for getting some useful pieces of information about our laws, as follows:

\begin{theorem}
The moments of the free hypergeometric laws are given by
$$M_k=\sum_{\pi,\nu\in NC(k)}W_{kN}(\pi,\nu)m^{|\pi|}p^{|\nu|}$$
with $W_{kN}$ being the Weingarten matrix for $S_N^+$. In particular we have
$$E=\frac{mp}{N}\quad,\quad V=\frac{mp(N-m)(N-p)}{N^2(N-1)}\quad,\quad 
\gamma=\frac{(N-2m)(N-2p)\sqrt{N-1}}{(N-2)\sqrt{mp(N-m)(N-p)}}$$
exactly as for the usual hypergeometric laws. As for the kurtosis, this is
$$\kappa=\frac{P(N,m,p)}{N^{10}(N-1)^{11}(N-2)^6(N^2-3N+1)m^2p^2(N-m)^2(N-p)^2}$$
with $P(N,m,p)$ being a certain polynomial, which might absorb some bottom terms.
\end{theorem}

\begin{proof}
According to the Weingarten formula for $S_N^+$, we have:
\begin{eqnarray*}
M_k
&=&\int_{S_N^+}\sum_{i_1=1}^m\ldots\sum_{i_k=1}^m
\sum_{j_1=1}^p\ldots\sum_{j_k=1}^pu_{i_1j_1}\ldots u_{i_pj_p}\\
&=&\sum_{i_1=1}^m\ldots\sum_{i_k=1}^m
\sum_{j_1=1}^p\ldots\sum_{j_k=1}^p
\sum_{\pi,\nu\in NC(k)}\delta_\pi(i)\delta_\nu(j)W_{kN}(\pi,\nu)\\
&=&\sum_{\pi,\nu\in NC(k)}W_{kN}(\pi,\nu)
\sum_{i_1=1}^m\ldots\sum_{i_k=1}^m\delta_\pi(i)
\sum_{j_1=1}^p\ldots\sum_{j_k=1}^p\delta_\nu(j)\\
&=&\sum_{\pi,\nu\in NC(k)}W_{kN}(\pi,\nu)m^{|\pi|}p^{|\nu|}
\end{eqnarray*}

Now the point is that, since we have $P(k)=NC(k)$ at $k=1,2,3$, the Weingarten matrices for $S_N$ and $S_N^+$ coincide at $k=1,2,3$, so we have that $E,V,\gamma$ figures. Finally, regarding the kurtosis, yes I know that's not professional, that denominator is the Di Francesco determinant of $G_{4N}$, divided by the square of the variance $V$. See \cite{ba2}.
\end{proof}

Still following \cite{bbc}, we have the following basic asymptotic study:

\begin{theorem}
We have the following asymptotic results,
\begin{enumerate}
\item With $N,m,p\to\infty$ and $mp/N\to t\in(0,\infty)$, we have:
$$\sigma_{mp}\sim\pi_t$$

\item With $M,m,p\to\infty$ and $m/N\to\nu\in (0,1)$ and $p/N\to 0$, we have:
$$\frac{\sigma_{mp}-p\nu}{\sqrt{p\nu(1-\nu)}}\sim\gamma_1$$
\end{enumerate}
with $\sigma_{mp}$ following the free hypergeometric law of parameters $(N,m,p)$.
\end{theorem}

\begin{proof}
This is quite standard, by using the Weingarten formula for the moments of the free hypergeometric laws, from Theorem 12.31, namely:
$$M_k=\sum_{\pi,\nu\in NC(k)}W_{kN}(\pi,\nu)m^{|\pi|}p^{|\nu|}$$

Indeed, we have the following standard estimate, for the Weingarten matrix:
$$W_{kN}(\pi,\nu)= 
\begin{cases}N^{-|\pi|}+O(N^{-|\pi|-1})&{\rm if}\ \pi=\nu\\ 
O(N^{|\pi\vee\nu|-|\pi|-|\nu|})&{\rm if}\ \pi\neq\nu
\end{cases}$$

Now in the regime $N,m,p\to\infty$ and $mp/N\to t\in(0,\infty)$, it follows that we have:
$$W_{kN}(\pi,\nu)m^{|\pi|}p^{|\nu|} \to 
\begin{cases}
t^{|\pi|}&{\rm if}\ \pi=\nu\\ 
0&{\rm if}\ \pi\neq\nu 
\end{cases}$$

Thus, in this regime we obtain indeed the moments of the Marchenko-Pastur law:
$$M_k=\sum_{\pi\in NC(k)}t^{|\pi|}$$

As for the second assertion, the proof is a bit more technical. See \cite{bbc}.
\end{proof}

With this discussed, time now for abstract algebra to strike. We have indeed the following result from \cite{bbc}, which was at that time something quite unexpected:

\begin{theorem}
The projective orthogonal and unitary quantum groups coincide, $PO_n^+=PU_n^+$, and appear as a twist of the quantum permutation group $S_{n^2}^+$,
$$PO_n^+=PU_n^+\sim S_{n^2}^+$$
and with this being something having no classical counterpart.
\end{theorem}
 
\begin{proof}
As mentioned, this is something quite surprising, the idea being as follows:

\medskip

(1) To start with, as explained in chapter 11, the inclusion $O_n^+\subset U_n^+$ produces an inclusion $PO_n^+\subset PU_n^+$, and since the main character on both sides follows $\pi_1$, we have:
$$PO_n^+=PU_n^+$$

Observe that this isomorphism has no classical counterpart, because $O_n,U_n$ have nothing much to do with each other, and so do their projective versions $PO_n,PU_n$.

\medskip

(2) Next, the above isomorphism is even better understood in the context of the quantum permutation groups $S_X^+$ of the finite quantum spaces $X$. Indeed, with the convention that $M_n$ is the quantum space given by $C(M_n)=M_n(\mathbb C)$, we have:
$$PO_n^+=PU_n^+=S_{M_n}^+$$

(3) Next, the fact that we have $M_n\sim\{1,\ldots,n^2\}$, twist of quantum spaces, raises the possibility for the corresponding quantum permutation groups to be related by a twisting operation. So, putting everything together, we would reach to a scheme as follows:
$$PO_n^+=PU_n^+=S_{M_n}^+\sim S_{n^2}^+$$ 

(4) And, is is true or not. At $n=2$ things are quite special, with the groups $SU_2$, $SO_3$ and their magic being involved, and it was known, basically since the probabilistic computations in \cite{bco}, that we have isomorphims as follows, confirming the above:
$$PO_2^+=PU_2^+=PU_2=SO_3\sim SO_3^{-1}\simeq S_4^+$$

(5) Now the point is that this holds in general too, as explained in \cite{bbc}, the reason behind this being the fact that the quantum automorphism group of a twisted group algebra is the twist of the quantum automorphism group of the original group algebra.

\medskip

(6) Finally, let us mention that all this is of potential key important in physics, showing that in the context of the Standard Model, in Chamseddine-Connes formulation \cite{cco}, as analysed by Bhowmick, Andrea and Dabrowski in \cite{bdd}, acting on the QED and QCD parts are twists of $S_4^+,S_9^+$. Which is quite interesting, suggesting that QED and QCD, suitably twisted, might be some sort of Yang-Mills theories based on $S_4^+,S_9^+$, respectively.
\end{proof}

Quite interesting all this, we are now deep into quarks. Still following \cite{bbc}, we have the following remarkable probabilistic consequence of the above physics trickery:

\begin{theorem}
The following families of variables have the same joint law,
$$\Big\{v_{ij}^2\Big|i,j=1,\ldots,n\Big\}\in C(O_n^+)\ \sim\ \Big\{\rho_{ij}=\frac{1}{n}\sum_{ab}u_{ia,jb}\Big|i,j=1,\ldots,n\Big\}\in C(S_{n^2}^+)$$
so in particular, the squared free hyperspherical law of parameter $n$ coincides with the $1/n$ rescaling of the free hypergeometric law of parameters $(n^2,n,n)$.
\end{theorem}

\begin{proof}
This comes from Theorem 12.33, but as pointed out in \cite{bbc}, thinking a bit retrospectively, this can be recovered as well directly, via Weingarten, as follows:
$$\int_{O_n^+}v_{ij}^{2k}=\sum_{\pi,\nu\in NC_2(2k)}W_{2k,n}(\pi,\nu)$$
$$\int_{S_{n^2}^+}\rho_{ij}^k=\sum_{\pi,\nu\in NC_2(2k)}n^{|\pi'|+|\sigma'|-k}W_{k,n^2}(\pi',\nu')$$

To be more precise, here $\pi\to\pi'$ is the shrinking operation for the pairings, appearing in chapter 10 in the context of Gram determinants, and this gives the result.
\end{proof}

Observe now that the squared free hyperspherical law of parameter $n$ is something that we know well, having struggled with it for pages in a row, in the previous section. What a comeback. In other words, we are dealing here with a generalization of that.

\bigskip

As a puzzling consequence of the above remarkable coincidence, we have:

\begin{theorem}
The squared hypespherical law of parameter $n$, having
$$E=\frac{1}{n}\quad,\quad V=\frac{2n-2}{n^2(n+2)}\quad,\quad\gamma=\frac{n-2}{n+4}\sqrt{\frac{8(n+2)}{n-1}}$$
and the $1/n$ rescaling of the hypergeometric law of parameters $(n^2,n,n)$, having
$$E=\frac{1}{n}\quad,\quad V=\frac{n-1}{n^2(n+1)}\quad,\quad\gamma=\frac{(n-2)^2}{n^3(n^2-2)}\sqrt{\frac{n+1}{n-1}}$$
have as common free analogue the free law discussed above, having
$$E=\frac{1}{n}\quad,\quad V=\frac{n-1}{n^2(n+1)}\quad,\quad\gamma=\frac{(n-2)^2}{n^3(n^2-2)}\sqrt{\frac{n+1}{n-1}}$$
and with the corresponding $3$ kurtoses being distinct.
\end{theorem}

\begin{proof}
The classical formulae are indeed well-known, and the free formulae come from Theorem 12.29 and its proof, or from Theorem 12.31.
\end{proof}

In relation with this, which is quite headaching, wondering what cat would have to say, about it. But cat is gone outside, chasing mice, guess we will never know.

\section*{12e. Exercises}

Free geometry rules, and here are some exercises, in relation with this:

\begin{exercise}
Further meditate on algebraic manifolds, and Gelfand functoriality.
\end{exercise}

\begin{exercise}
Learn more about the classical and free Meixner laws.
\end{exercise}

\begin{exercise}
Construct partial isometry spaces, for any easy quantum group.
\end{exercise}

\begin{exercise}
Learn about the quantum affine homogeneous spaces.
\end{exercise}

\begin{exercise}
Do the calculus, for the classical hyperspherical laws.
\end{exercise}

\begin{exercise}
Do the calculus, for the free hyperspherical laws.
\end{exercise}

\begin{exercise}
Do the calculus, for the classical hypergeometric laws.
\end{exercise}

\begin{exercise}
Do the calculus, for the free hypergeometric laws.
\end{exercise}

As bonus exercise, difficult, convert what we said above into particle physics.

\part{Matrix models}

\ \vskip50mm

\begin{center}
{\em Strangers in the night

Exchanging glances

Wandering in the night

What were the chances}
\end{center}

\chapter{Wishart matrices}

\section*{13a. Block transposition}

With the free probability basics learned, and a few good questions to work on learned too, time I guess to end this book. But of course, not without giving a few pieces of advice, regarding the more precise continuation of your learning:

\bigskip

(1) If you're more of an algebraist, you might enjoy relearning free probability, in an axiomatic way, with full details, and with due attention to the various algebraic tools available. In this case, go first with the book of Voiculescu et al. \cite{vdn}, a quick, historical must-read, and then the book of Nica and Speicher \cite{vdn} is for you.

\bigskip

(2) If you're also an algebraist, but of the functional analysis type, you might want to know more about operator algebras, which are something quite fascinating, and about what free probability tells us, about them. In this case carefully read Voiculescu et al. \cite{vdn}, followed by a more recent, complete book, such as Hiai and Petz \cite{hpe}.

\bigskip

(3) On the opposite, if you're more of an analyst, a bit tired by all the algebra that you had to endure, from the present book, that is a healthy attitude too, and you can go with the many papers available, on the analytic aspects of free probability. Some good starters here are the papers of Bo\.zejko and collaborators \cite{bbl}, \cite{brs}, \cite{bwy}.

\bigskip

(4) Finally, if you're a hardline analyst, why bothering with algebras and abstract variables and other mathematical inventions, random matrices are what you need. Go here with the historical book of Mehta \cite{meh}, followed by \cite{agz}. With the remark tough that for freeness aspects of the random matrices, you will have to get into (1).

\bigskip

With this said, you might wonder, what about geometers? Well, stay with me guys, with our algebraist and analytic colleagues put on a logical loop, the world is ours. We will be interested in the following question, which is somehow of geometric flavor:

\begin{question}
Can we find models for the more specialized measures $\mu$ that we met? What about models for variables, such as the coordinates $x_i$ of free manifolds?
\end{question}

So, this will be the question that we will discuss, in the present Part IV of this book. With the comment of course that this will be something a bit specialized.

\bigskip

Getting started now, let us quickly review the free measures $\mu$ that we met. These include the Wigner law $\gamma_t$, majestically modeled by the Wigner matrices, and then the Marchenko-Pastur law $\pi_t$, majestically modeled by the Wishart matrices. Next, we have the free Bessel laws $\pi^s_t$, with several matrix models for them being available at $t=1$. And then we have more general compound free Poisson laws $\pi_\nu$, as well as other interesting free measures, such as the arcsine laws $\alpha_t,\alpha_t'$, or the mysterious uniform law $u_t$.

\bigskip

In view of all this, I would say that modeling $\pi_\nu$ would be a good question. Indeed, substantial family of measures that we have here, more or less at the center of the free galaxy. And in addition, with $p_\nu\to\pi_\nu$ being a liberation operation, further studying $\pi_\nu$ will certainly teach us more, about how free probability works, in general.

\bigskip

Which sounds good, so let us summarize this discussion as follows:

\begin{problem}
How to manipulate the Wishart matrices, and the Marchenko-Pastur theorem, as to suitably model the general compound free Poisson laws $\pi_\nu$?
\end{problem}

In order to discuss this question, let us first recall what we previously know about it, from the previous chapters. At the beginning of everything, we have the Marchenko-Pastur theorem, that we know well since chapter 3, which is as follows:

\index{asymptotic freeness}

\begin{theorem}
Given a sequence of complex Wishart matrices
$$W_N=Y_NY_N^*\in M_N(L^\infty(X))$$
with $Y_N$ being $N\times M$ complex Gaussian of parameter $1$, we have
$$\frac{W_N}{N}\sim\max(1-t,0)\delta_0+\frac{\sqrt{4t-(x-1-t)^2}}{2\pi x}\,dx$$
with $M=tN\to\infty$, with the limiting measure being the Marchenko-Pastur law $\pi_t$.
\end{theorem}

\begin{proof}
This is the Marchenko-Pastur theorem \cite{mpa}, that we know well from chapter 3. Let us also mention, at the level of the basics, that, as explained in chapter 8, as a consequence of the results of Voiculescu in \cite{vo3}, when looking at a family of sequences of such matrices, with $t\geq1$, in the $N\to\infty$ limit these become asymptotically free.
\end{proof}

Regarding now the possible manipulations on the Wishart matrices, there are many of them, gradually explored over the time, and quite systematically in the last 20 years, starting with the key 2003 paper of Graczyk, Letac and Massam \cite{glm}. We have already talked about this, and more specifically about various product techniques, in chapter 8, following \cite{bb+}, and many other things can be said, as a continuation of this.

\bigskip

However, it is not about products, which remain something quite technical, that we would like to talk about, in this chapter. As something which is perhaps more fruitful, or at least looks more basic, we have the idea of block-transposing or block-modifying a Wishart matrix, according to the following suprising finding of Aubrun \cite{aub}:

\begin{fact}
When suitably block-transposing the entries of a complex Wishart matrix, we obtain as asymptotic distribution a shifted version of Wigner's semicircle law. 
\end{fact}

So, let us first try to understand this, and for generalizations, we can see later. Following \cite{aub}, and some more recent work in \cite{bn1}, where the Aubrun finding was further studied, and better understood, let us start with the following definition:

\begin{definition}
The partial transpose of a complex Wishart matrix $W$ of parameters $(N,M)=(dn,dm)$ is the matrix
$$\widetilde{W}=(id\otimes t)W$$
where $id$ is the identity of $M_d(\mathbb C)$, and $t$ is the transposition of $M_n(\mathbb C)$. 
\end{definition}

In more familiar terms of bases and indices, the standard decomposition $\mathbb C^{dn}=\mathbb C^d\otimes\mathbb C^n$ induces an algebra decomposition $M_{dn}(\mathbb C)=M_d(\mathbb C)\otimes M_n(\mathbb C)$, and with this convention made, the partial transpose matrix $\widetilde{W}$ constructed above has entries as follows:
$$\widetilde{W}_{ia,jb}=W_{ib,ja}$$

Our goal in what follows will be that of computing the law of $\widetilde{W}$, first when $d,n,m$ are fixed, and then in the $d\to\infty$ regime. For this purpose, we will need a number of standard facts regarding the noncrossing partitions. Let us start with:

\begin{proposition}
For a permutation $\sigma\in S_p$, we have the formula
$$|\sigma|+\#\sigma=p$$
where $|\sigma|$ is the number of cycles of $\sigma$, and $\#\sigma$ is the minimal $k\in\mathbb N$ such that $\sigma$ is a product of $k$ transpositions. Also, the following formula defines a distance on $S_p$, 
$$(\sigma,\pi)\to\#(\sigma^{-1}\pi)$$
and the set of permutations $\sigma\in S_p$ which saturate the triangular inequality 
$$\#\sigma+\#(\sigma^{-1}\gamma)=\#\gamma=p-1$$
where $\gamma\in S_p$ is a full cycle, is in bijection with the set $NC(p)$. 
\end{proposition}

\begin{proof}
All this is standard combinatorics, that can be found for instance, along with many other things, in the paper of Biane \cite{bia}, and that we will leave as an exercise.
\end{proof}

We use the standard bijection $NC(p)\simeq NC_2(2p)$, denoted $\pi\to\widetilde{\pi}$, obtained by fattening the partitions. We have the following formula,  where $\vee$ is the join operation on $NC_2(2p)$, and $\rho_{12} = (12)(34)\ldots(2p-1,2p)$ is the fattened identity permutation:
$$|\pi|=|\widetilde{\pi}\vee\rho_{12}|$$

Similarly, we have the formula $|\pi\gamma|=|\widetilde{\pi}\vee\rho_{14}|$, where $\rho_{14}$ is the pairing corresponding to the fattening of the inverse full cycle $\gamma^{-1}(i) = i-1$, which pairs an element $2i$ with $2(i-1)-1 = 2i-3$, or, equivalently, an element $i \in \{1, \ldots, 2p\}$ with $i+(-1)^{i+1}3$. 

\bigskip

We will need the following well-known result, which is something standard:

\begin{proposition}
The number $||\pi||$ of blocks having even size is given by
$$1+||\pi||=|\pi\gamma|$$
for every noncrossing partition $\pi \in NC(p)$.
\end{proposition}

\begin{proof}
We use a recurrence over the number of blocks of $\pi$. If $\pi$ has just one block, its associated geodesic permutation is $\gamma$ and we have, as desired:
$$|\gamma^2|=1+\delta_{2|p}$$

For the partitions $\pi$ having more than one block, we can assume without loss of generality that $\pi = \hat 1_k \sqcup \pi'$, where $\hat 1_k$ is a contiguous block of size $k$. Recall now that the number of blocks of the permutation $\pi\gamma$ is given by the following formula, where $\rho_{14} \in P_2(2p)$ is the pair partition which pairs an element $i$ with $i+(-1)^{i+1}3$:
$$|\pi\gamma|=|\widetilde{\pi}\vee\rho_{14}|$$

If $k$ is an even number, $k=2r$, in order to prove the result, consider the following partition, which contains the block $(1 \, 4 \, 5 \, 8 \, \ldots 4r-3 \, 4r)$, along with the blocks coming from the elements of the form $4i+2, 4i+3$ from $\{1, \ldots, 4r\}$ and from $\pi'$:
$$\sigma=\widetilde{\hat 1_{2r} \sqcup \pi'}\vee \rho_{14}$$

Now recall that we can count the blocks of the join of two partitions by drawing them one beneath the other, and then counting the number of connected components of the curve, without taking into account the possible crossings. We conclude that we have the following formula, where $\rho'_{14}$ is $\rho_{14}$ restricted to the set $\{2k+1, 2k+2 \ldots, 2p\}$:
$$|\widetilde{\pi}\vee\rho_{14}|=1+|\widetilde{\pi'}\vee\rho'_{14}|$$

As for the case where $k$ is odd, here there is no extra block appearing, and so:
$$|\widetilde{\pi}\vee\rho_{14}|=|\widetilde{\pi'}\vee\rho'_{14}|$$

Summarizing, we are led to the conclusion in the statement.
\end{proof}

We can now investigate the block-transposed Wishart matrices, and we have:

\begin{theorem}
For any $p\geq 1$ we have the formula
$$\lim_{d\to\infty}(E\circ tr)\big(m\widetilde{W}\big)^p
=\sum_{\pi\in NC(p)}m^{|\pi|}n^{||\pi||}$$
where $|.|$ and $||.||$ are the number of blocks, and the number of blocks of even size. 
\end{theorem}

\begin{proof}
The matrix elements of the partial transpose matrix are given by:
$$\widetilde{W}_{ia,jb}=W_{ib,ja}=(dm)^{-1}\sum_{k=1}^d\sum_{c=1}^mG_{ib,kc}\bar{G}_{ja,kc}$$

This gives the following formula, for the corresponding normalized trace:
\begin{eqnarray*}
tr(\widetilde{W}^p)
&=&(dn)^{-1}(dm)^{-p}\sum_{i_1,\ldots,i_p=1}^d\sum_{a_1,\ldots,a_p=1}^n\prod_{s=1}^p W_{i_sa_{s+1},i_{s+1}a_s} \\
&=&(dn)^{-1}(dm)^{-p}\sum_{i_1,\ldots,i_p=1}^d\sum_{a_1,\ldots,a_p=1}^n\prod_{s=1}^p \sum_{j_1,\ldots,j_p=1}^d\sum_{b_1,\ldots,b_p=1}^mG_{i_sa_{s+1},j_sb_s}\bar{G}_{i_{s+1}a_s,j_sb_s}
\end{eqnarray*}

After interchanging the product with the last two sums, the average of the general term can be computed by the Wick rule, namely:
$$E\left(\prod_{s=1}^pG_{i_sa_{s+1},j_sb_s}\bar{G}_{i_{s+1}a_s,j_sb_s}\right)
=\sum_{\pi\in S_p}\prod_{s=1}^p\delta_{i_s,i_{\pi(s)+1}}\delta_{a_{s+1},a_{\pi(s)}}\delta_{j_s,j_{\pi(s)}}\delta_{b_s,b_{\pi(s)}}$$

Let $\gamma\in S_p$ be the full cycle $\gamma=(1 \, 2 \, \ldots \, p)^{-1}$. The general factor in the above product is 1 if and only if the following four conditions are simultaneously satisfied:
$$\gamma^{-1}\pi\leq \ker i\quad,\quad
\pi\gamma \leq \ker a\quad,\quad
\pi \leq \ker j\quad,\quad
\pi \leq \ker b$$

Counting the number of free parameters in the above equation, we obtain:
\begin{eqnarray*}
(E\circ tr)(\widetilde{W}^p)
&=&(dn)^{-1}(dm)^{-p}\sum_{\pi\in S_p}d^{|\pi|+|\gamma^{-1}\pi|}m^{|\pi|}n^{|\pi\gamma|}\\
&=&\sum_{\pi\in S_p}d^{|\pi|+|\gamma^{-1}\pi|-p-1}m^{|\pi|-p}n^{|\pi\gamma|-1}
\end{eqnarray*}

The exponent of $d$ in the last expression on the right is as follows:
\begin{eqnarray*}
N(\pi)
&=&|\pi|+|\gamma^{-1}\pi|-p-1\\
&=&p-1-(\#\pi+\#(\gamma^{-1}\pi))\\
&=&p-1-(\#\pi+\#(\pi^{-1}\gamma))
\end{eqnarray*}

As explained in the beginning of this section, this quantity is known to be $\leq 0$, with equality iff $\pi$ is geodesic, hence associated to a noncrossing partition. Thus:
$$(E\circ tr)(\widetilde{W}^p)=(1+O(d^{-1}))m^{-p}n^{-1}\sum_{\pi\in NC(p)}m^{|\pi|} n^{|\pi\gamma|}$$

But together with $|\pi\gamma|=||\pi||+1$, this gives the result.
\end{proof}

We would like now to find an equation for the moment generating function of the asymptotic law of $m\widetilde{W}$. This moment generating function is defined as:
$$F(z)=\lim_{d\to\infty}(E\circ tr)\left(\frac{1}{1-zm\widetilde{W}}\right)$$

We have the following result, regarding this moment generating function:

\begin{theorem}
The moment generating function of $m\widetilde{W}$ satisfies the equation
$$(F-1)(1-z^2F^2)=mzF(1+nzF)$$
in the $d\to\infty$ limit.
\end{theorem}

\begin{proof}
We use the formula in Theorem 13.8. If we denote by $N(p,b,e)$ the number of partitions in $NC(p)$ having $b$ blocks and $e$ even blocks, we have:
\begin{eqnarray*}
F
&=&1+\sum_{p=1}^\infty\sum_{\pi\in NC(p)} z^pm^{|\pi|}n^{||\pi||}\\
&=&1+\sum_{p=1}^\infty\sum_{b=0}^\infty\sum_{e=0}^\infty z^pm^bn^eN(p,b,e)
\end{eqnarray*}

Let us try to find a recurrence formula for the numbers $N(p,b,e)$. If we look at the block containing $1$, this block must have $r\geq 0$ other legs, and we get:
\begin{eqnarray*}
N(p,b,e) 
&=&\sum_{r\in 2\mathbb N}\sum_{p=\Sigma p_i+r+1}\sum_{b=\Sigma b_i+1}\sum_{e=\Sigma e_i}N(p_1,b_1,e_1)\ldots N(p_{r+1},b_{r+1},e_{r+1})\\
&+&\sum_{r\in 2\mathbb N+1}\sum_{p=\Sigma p_i+r+1}\sum_{b=\Sigma b_i+1}\sum_{e=\Sigma e_i+1}N(p_1,b_1,e_1)\ldots N(p_{r+1},b_{r+1},e_{r+1})
\end{eqnarray*}

Here $p_1,\ldots,p_{r+1}$ are the number of points between the legs of the block containing 1, so that we have $p=(p_1+\ldots+p_{r+1})+r+1$, and the whole sum is split over two cases, $r$ even or odd, because the parity of $r$ affects the number of even blocks of our partition. Now by multiplying everything by a $z^pm^bn^e$ factor, and by carefully distributing the various powers of $z,m,b$ on the right, we obtain the following formula:
\begin{eqnarray*}
z^pm^bn^eN(p,b,e)
&=&m\sum_{r\in 2\mathbb N}z^{r+1}\sum_{p=\Sigma p_i+r+1}\sum_{b=\Sigma b_i+1}\sum_{e=\Sigma e_i}\prod_{i=1}^{r+1}z^{p_i}m^{b_i}n^{e_i}N(p_i,b_i,e_i)\\
&+&mn\sum_{r\in 2\mathbb N+1}z^{r+1}\sum_{p=\Sigma p_i+r+1}\sum_{b=\Sigma b_i+1}\sum_{e=\Sigma e_i+1}\prod_{i=1}^{r+1}z^{p_i}m^{b_i}n^{e_i}N(p_i,b_i,e_i)
\end{eqnarray*}

Let us sum now all these equalities, over all $p\geq 1$ and over all $b,e\geq 0$. According to the definition of $F$, at left we obtain $F-1$. As for the two sums appearing on the right, that is, at right of the two $z^{r+1}$ factors, when summing them over all $p\geq 1$ and over all $b,e\geq 0$, we obtain in both cases $F^{r+1}$. So, we have the following formula:
\begin{eqnarray*}
F-1
&=&m\sum_{r\in 2\mathbb N}(zF)^{r+1}+mn\sum_{r\in 2\mathbb N+1}(zF)^{r+1}\\
&=&m\,\frac{zF}{1-z^2F^2}+mn\,\frac{z^2F^2}{1-z^2F^2}\\
&=&mzF\,\frac{1+nzF}{1-z^2F^2}
\end{eqnarray*}

But this gives the formula in the statement, and we are done.
\end{proof}

Our goal now will be that of further processing the formula in Theorem 13.9, as to reach to a formula for the density of the corresponding law. We first have:

\begin{theorem}
The Cauchy transform of $m\widetilde{W}$ satisfies the equation
$$(\xi G-1)(1-G^2)=mG(1+nG)$$
in the $d\to\infty$ limit. Moreover, this equation simply reads
$$R=\frac{m}{2}\left(\frac{n+1}{1-z}-\frac{n-1}{1+z}\right)$$
with the usual $R$-transform substitutions $G\to z$ and $\xi\to R+z^{-1}$.
\end{theorem}

\begin{proof}
We have two assertions to be proved, the first one being standard, and the second one being something quite magic, the idea being as follows:

\medskip

(1) Consider the equation of $F$, found in Theorem 13.9, namely:
$$(F-1)(1-z^2F^2)=mzF(1+nzF)$$

With $z\to\xi^{-1}$ and $F\to\xi G$, so that $zF\to G$, we obtain, as desired:
$$(\xi G-1)(1-G^2)=mG(1+nG)$$

(2) Next, let us look at the equation of the Cauchy transform that we have. With the substitutions $\xi \to K$ and $G\to z$, this equation becomes:
$$(zK-1)(1-z^2)=mz(1+nz)$$

The point now is that with $K\to R+z^{-1}$ this latter equation becomes:
$$zR(1-z^2)=mz(1+nz)$$

But the solution of this latter equation is trivial to compute, given by:
$$R=m\,\frac{1+nz}{1-z^2}=\frac{m}{2}\left(\frac{n+1}{1-z}-\frac{n-1}{1+z}\right)$$

Thus, we are led to the conclusion in the statement.
\end{proof}

We can now answer our original question, regarding the limiting law, as follows:

\begin{theorem}
With $d\to\infty$ we have the following convergence in law,
$$m\widetilde{W}\sim\mu_{mn}$$
where $\mu_{mn}$ is the free difference of free Poisson laws of parameters $m(n\pm 1)/2$.
\end{theorem}

\begin{proof}
Consider two free variables $a,b$, following free Poisson laws of parameters $s,t$. According to our formulae from chapter 8, the corresponding $R$-transforms are:
$$R_a=\frac{s}{1-z}\quad,\quad R_b=\frac{t}{1-z}$$

By using the general dilation formula $R_{qb}(z)=qR_b(qz)$ at $q=-1$, we deduce that we have $R_{-b}=-t/(1+z)$. Thus, our two formulae above can be written as follows:
$$R_a=\frac{s}{1-z}\quad,\quad R_{-b}=-\frac{t}{1+z}$$

Now when making the sum, the corresponding $R$-transform will be given by:
$$R_{a-b}=\frac{s}{1-z}-\frac{t}{1+z}$$

But with $s=m(n+1)/2$ and $t=m(n-1)/2$, this latter formula reads:
$$R=\frac{m}{2}\left(\frac{n+1}{1-z}-\frac{n-1}{1+z}\right)$$

And since this is the formula found in Theorem 13.10, we obtain the result.
\end{proof}

\section*{13b. Shifted semicircles}

Still following \cite{bn1}, let us have now a closer look at the asymptotic laws that we found in Theorem 13.11, among others with the goal of recovering, as a particular case of that result, the original finding by Aubrun, regarding shifted semicircles \cite{aub}.

\bigskip

It is convenient to enlarge our study to all possible real values of $m,n$, as follows:

\begin{definition}
The real measure $\mu_{mn}$ of parameters $m>0,n\geq 1$ is given by
$$a\sim\pi_{m(n+1)/2}\ ,\ b\sim\pi_{m(n-1)/2}\ ,\ {\rm free}\quad\implies\quad a-b\sim\pi_{mn}$$
that is, is the free difference of free Poisson laws of parameters $m(n\pm 1)/2$.
\end{definition}

As a first remark, the various formulae that we found before, namely moment formula, and various equations for the functional transforms, are valid in this more general setting. Let us collect these formulae, along with a bit more, in a single statement:

\begin{theorem}
The measure $\mu_{mn}$ is a compound free Poisson law,
$$\mu_{mn}=\pi_\nu\quad,\quad\nu=\frac{m(n-1)}{2}\,\delta_{-1}+\frac{m(n+1)}{2}\,\delta_1$$
and has the following properties:
\begin{enumerate}
\item The moments are $M_k=\sum_{\pi\in NC(k)}m^{|\pi|}n^{||\pi||}$.

\item The mean is $E=m$, the variance is $V=mn$.

\item $\gamma=1/(n\sqrt{mn})$ and $\kappa=2+1/(mn)$.

\item $(M-1)(1-z^2M^2)=mzM(1+nzM)$.

\item $(\xi G-1)(1-G^2)=mG(1+nG)$.

\item $R(z)=m(1+nz)/(1-z^2)$.
\end{enumerate}
\end{theorem}

\begin{proof}
These are things that we actually discussed already, in chapter 8, but a few reminders won't hurt. Regarding the first assertion, recall that we have:
$$\nu=\sum_ic_i\delta_{y_i}\quad\implies\quad R_{\pi_{\nu}}(z)=\sum_i\frac{c_iy_i}{1-y_iz}$$

For the measure $\nu$ in the statement, we obtain the following $R$-transform:
$$R(z)=\frac{m(n-1)}{2}\cdot\frac{-1}{1+z}+\frac{m(n+1)}{2}\cdot\frac{1}{1-z}
=\frac{m}{2}\left(\frac{n+1}{1-z}-\frac{n-1}{1+z}\right)$$

But this is the $R$-transform of $\mu_{mn}$ that we computed before. As for the rest, we already explained in chapter 8 how to perform the following passages:
\begin{eqnarray*}
R(z)=m(1+nz)/(1-z^2)
&\implies&(\xi G-1)(1-G^2)=mG(1+nG)\\
&\implies&(M-1)(1-z^2M^2)=mzM(1+nzM)
\end{eqnarray*}

Thus, it remains to compute $E,V,\gamma,\kappa$. The low order moments are as follows:
$$M_1=m$$
$$M_2=m^2+mn$$
$$M_3=m^3+3m^2n+m$$
$$M_4=m^4+6m^3n+4m^2+2m^2n^2+mn$$

Thus the mean is $E=m$, and the variance is $V=mn$. Next, we have:
$$M_3'=M_3-3EM_2+2E^3=m$$
$$M_4'=M_4-4EM_3+6E^2M_2-3E^4=mn(2mn+1)$$ 

Thus the skewness and kurtosis are given by the following formulae:
$$\gamma=\frac{m}{mn\sqrt{mn}}
=\frac{1}{n\sqrt{mn}}\quad,\quad  
\kappa=\frac{mn(2mn+1)}{m^2n^2}
=2+\frac{1}{mn}$$

Finally, observe that we have $\kappa\simeq 2$ as we should, our law being of free nature. 
\end{proof}

Regarding now the main particular cases of measures of type $\mu_{mn}$, corresponding to some previously known random matrix computations, the situation is as follows:

\begin{theorem}
The measures $\mu_{mn}$ are as follows:
\begin{enumerate}
\item At $n=1$ we have the Marchenko-Pastur laws, $\mu_{m1}=\pi_m$.

\item With $n=sm\to\infty$ we obtain Aubrun's shifted semicircles.

\item With $m=t/n\to 0$ we obtain the free Bessel law $\pi^2_t$.
\end{enumerate}
\end{theorem}

\begin{proof}
These results come either from Definition 13.12, or from Theorem 13.13:

\medskip

(1) This follows from Definition 13.12, because at $n=1$ we have $m(n\pm 1)/2=m,0$. Observe that this is in tune with the Theorem 13.3, because at $n=1$ the partial transposition is trivial, so Theorem 13.11 computes the asymptotic law of $mW$.

\medskip

(2) In order to explain what is going on here, recall from Theorem 13.11 that, with the notations there, we have $m\widetilde{W}\sim\mu_{mn}$ in the $d\to\infty$ limit. Thus, with the convention that $D$ is the dilation operation, given by $a\sim\mu\implies qa\sim D_q(\mu)$, we have:
$$\widetilde{W}\sim D_{m^{-1}}(\mu_{mn})$$

Now the point is that, what Aubrun discovered in \cite{aub} is the following formula, valid in the regime in the statement, $n=sm\to\infty$, and with $d\to\infty$ too, as usual:
$$\widetilde{W}\sim\gamma_s^1$$

So, let us prove this. In view of what we have, this amounts in proving that:
$$\lim_{m\to\infty}D_{m^{-1}}(\mu_{m,sm})=\gamma_s^1$$

For this purpose, let $\Gamma=\Gamma(\xi)$ be the Cauchy transform of $\nu=D_{m^{-1}}(\mu_{mn})$. The Cauchy transform of $\mu_{mn}=D_m(\nu)$ is then $G(m\xi)=m^{-1}\Gamma$, and by making the replacements $\xi\to m\xi$ and $G\to m^{-1}\Gamma$ in the equation of $G$ in Theorem 13.13, we obtain: 
$$(\xi\Gamma-1)(1-m^{-2}\Gamma^2)=\Gamma(1+nm^{-1}\Gamma)$$

In the limit $n=sm\to\infty$, as indicated above, this equation becomes:
\begin{eqnarray*}
\xi\Gamma-1=\Gamma(1+s\Gamma)
&\implies&s\Gamma^2+(1-\xi)\Gamma+1=0\\
&\implies&\Gamma=\frac{\xi-1\pm\sqrt{(1-\xi)^2-4s}}{2s}
\end{eqnarray*}

By applying now the Stieltjes inversion formula, we obtain the following density:
$$\varphi(x)=\frac{\sqrt{4s-(x-1)^2}}{2s\pi}$$

But this is exactly the density of the semicircle law $\gamma_s^1$, and we are done.

\medskip

(3) This comes from Theorem 13.13, because with $m=t/n\to 0$, the measure $\nu$ appearing there, producing the compound free Poisson law $\mu_{mn}=\pi_\nu$, is given by: 
\begin{eqnarray*}
\nu
&=&\frac{m(n-1)}{2}\,\delta_{-1}+\frac{m(n+1)}{2}\,\delta_1\\
&=&\frac{t-m}{2}\,\delta_{-1}+\frac{t+m}{2}\,\delta_1\\
&\simeq&t\cdot\frac{\delta_{-1}+\delta_1}{2}
\end{eqnarray*}

Thus, we obtain indeed the free Bessel law $\pi^2_t$, as constructed in chapter 8.
\end{proof}

As a comment now on (3), we were talking before the statement of the theorem about ``previous random matrix computations'', and in the case of the free Bessel law $\pi^2_t$, this refers to the computations in \cite{bb+} using products of Gaussian matrices, producing matrix models for $\pi^2_1$, or rather for its modified version $\tilde{\pi}^2_1$, in the sense of chapter 8.

\bigskip

Now the point is that these latter computations have nothing to do with our current block-transposition business for Wishart matrices. Thus, in short, what we have here are different approaches to the modeling on $\pi^2_1$. And in relation with this, an interesting question appears. Indeed, with the methods in \cite{bb+} working for any $\pi^s_1$, we have:

\begin{question}
Can our block-transposition business for Wishart matrices be suitably modified, as to cover the general free Bessel laws $\pi^s_t$?
\end{question}

And good question this is, very natural because, as explained earlier in this book, the free Bessel laws $\pi^s_t$ are at the core of free probability of discrete type, so any candidate for a serious subtheory in free probability should have these laws $\pi^s_t$ well understood.

\bigskip

Well, so this is at least how the general philosophy goes. In practice now, Question 13.15 has led to a lot of work, first with my second paper with Nechita, laying the foundations for a possible answer, and we will talk about this in a moment, and then with two follow-ups to that paper, one by myself not solving the problem, but working on that was fun, and one by Nechita with Arizmendi and Vargas \cite{anv}, basically solving the problem, but in a quite complicated way. Somehow, Question 13.15 still stands.

\bigskip

Nevermind. As a last result now, that we also know from chapter 8, we have:

\begin{theorem}
The density of $\mu_{mn}$ has the following properties:
\begin{enumerate}
\item It has at most one atom, located at $0$, of mass $\max(1-mn,0)$.

\item The support is positive precisely when $n\leq m/4+1/m$ and $m \geq 2$.
\end{enumerate}
\end{theorem}

\begin{proof}
This is something that we discussed in chapter 8, the idea being that (1) follows from the general results of Bercovici-Voiculescu in \cite{bvo}, and that (2) follows from a long study, featuring several discriminants, based on the Cauchy transform equation from Theorem 13.13, namely $(\xi G-1)(1-G^2)=mG(1+nG)$. See \cite{bn1}.
\end{proof}

\section*{13c. Block modifications}

As a continuation of the above, let us discuss now some more general block modifications of the complex Wishart matrices, following the various more recent continuations of \cite{bn1}, with the goal of solving Problem 13.2 and Question 13.15, and also, importantly, of providing tools for quantum information theory, where such things can help.

\bigskip

We have the following construction, that we will be interested in:

\begin{definition}
Given a complex Wishart matrix $W=YY^*\in M_{dn}(L^\infty(X))$ as befoe, with $Y$ being a complex Gaussian $dn\times dm$ matrix, and a linear map
$$\varphi:M_n(\mathbb C)\to M_n(\mathbb C)$$
we consider the following matrix, obtained by applying $\varphi$ to the $n\times n$ blocks of $W$,
$$\widetilde{W}=(id\otimes\varphi)W\in M_{dn}(L^\infty(X))$$
and call it block-modified Wishart matrix.
\end{definition}

Here we are using some standard tensor product identifications, the details being as follows. Let $Y$ be a complex Gaussian $dn\times dm$ matrix, as above:
$$Y\in M_{dn\times dm}(L^\infty(X))$$

We can then form the corresponding complex Wishart matrix, as follows: 
$$W=YY^*\in M_{dn}(L^\infty(X))$$

The size of this matrix being a composite number, $N=dn$, we can regard this matrix as being a $n\times n$ matrix, with random $d\times d$ matrices as entries. Equivalently, by using standard tensor product notations, this amounts in regarding $W$ as follows:
$$W\in M_d(L^\infty(X))\otimes M_n(\mathbb C)$$

With this done, we can come up with our linear map, namely:
$$\varphi:M_n(\mathbb C)\to M_n(\mathbb C)$$

We can apply $\varphi$ to the tensors on the right, and we obtain a matrix as follows:
$$\widetilde{W}=(id\otimes\varphi)W\in M_d(L^\infty(X))\otimes M_n(\mathbb C)$$

Finally, we can forget now about tensors, and as a conclusion to all this, we have constructed a matrix as follows, that we can call block-modified Wishart matrix:
$$\widetilde{W}\in M_{dn}(L^\infty(X))$$

In practice now, what we mostly need for fully understanding Definition 13.17 are examples. Following Aubrun \cite{aub}, then \cite{bn1}, and the paper by Collins and Nechita \cite{cne} too, we have the following basic examples, for our general construction:

\begin{proposition}
We have the following examples of block-modified Wishart matrices $\widetilde{W}=(id\otimes\varphi)W$, coming from various linear maps $\varphi:M_n(\mathbb C)\to M_n(\mathbb C)$:
\begin{enumerate}
\item Wishart matrices: $\widetilde{W}=W$, obtained via $\varphi=id$.

\item Aubrun matrices: $\widetilde{W}=(id\otimes t)W$, with $t$ being the transposition.

\item Collins-Nechita one: $\widetilde{W}=(id\otimes\varphi)W$, with $\varphi=tr(.)1$.

\item Collins-Nechita two: $\widetilde{W}=(id\otimes\varphi)W$, with $\varphi$ erasing the off-diagonal part.
\end{enumerate}
\end{proposition}

\begin{proof}
This is indeed something self-explanatory, standing more as a definition.
\end{proof}

Getting back now to the general case, that of Definition 13.17, the linear map there $\varphi:M_n(\mathbb C)\to M_n(\mathbb C)$ is certainly useful for understanding the construction of the block-modified Wishart matrix $\widetilde{W}=(id\otimes\varphi)W$, as illustrated by the above examples. In practice, however, we would like to have as block-modification ``data'' something more concrete, such as a usual matrix. To be more precise, we would like to use:

\begin{proposition}
We have a correspondence between linear maps 
$$\varphi:M_n(\mathbb C)\to M_n(\mathbb C)$$
and square matrices $\Lambda\in M_n(\mathbb C)\otimes M_n(\mathbb C)$, given by the formula
$$\Lambda_{ab,cd}=\varphi(e_{ac})_{bd}$$
where $e_{ab}\in M_n(\mathbb C)$ are the standard generators of the matrix algebra $M_n(\mathbb C)$, given by the formula $e_{ab}:e_b\to e_a$, with $\{e_1,\ldots,e_n\}$ being the standard basis of $\mathbb C^n$.
\end{proposition}

\begin{proof}
This is standard linear algebra. Given a linear map $\varphi:M_n(\mathbb C)\to M_n(\mathbb C)$, we can associated to it numbers $\Lambda_{ab,cd}\in\mathbb C$ by the formula in the statement, namely:
$$\Lambda_{ab,cd}=\varphi(e_{ac})_{bd}$$

Now by using these $n^4$ numbers, we can construct a $n^2\times n^2$ matrix, as follows:
$$\Lambda=\sum_{abcd}\Lambda_{ab,cd}e_{ac}\otimes e_{bd}\in M_n(\mathbb C)\otimes M_n(\mathbb C)$$

Thus, we have constructed a correspondence $\varphi\to\Lambda$, and since this correspondence is injective, and the dimensions match, this correspondence is bijective, as claimed.
\end{proof}

Now by getting back to the block-modified Wishart matrices, we have:

\begin{proposition}
Given a Wishart $dn\times dn$ matrix $W=YY^*$, and a linear map
$$\varphi:M_n(\mathbb C)\to M_n(\mathbb C)$$
the entries of the corresponding block-modified matrix $\widetilde{W}=(id\otimes\varphi)W$ are given by
$$\widetilde{W}_{ia,jb}=\sum_{cd}\Lambda_{ca,db}W_{ic,jd}$$
where $\Lambda\in M_n(\mathbb C)\otimes M_n(\mathbb C)$ is the square matrix associated to $\varphi$, as above.
\end{proposition}

\begin{proof}
Again, this is linear algebra, coming from the following computation:
$$\widetilde{W}_{ia,jb}
=\sum_{cd}W_{ic,jd}\varphi(e_{cd})_{ab}
=\sum_{cd}\Lambda_{ca,db}W_{ic,jd}$$

We are therefore led to the conclusion in the statement.
\end{proof}

At the level of the main examples, from Proposition 13.18, the very basic linear maps $\varphi:M_n(\mathbb C)\to M_n(\mathbb C)$ used there can only correspond to some basic examples of matrices $\Lambda\in M_n(\mathbb C)\otimes M_n(\mathbb C)$, via the correspondence in Proposition 13.19. This is indeed the case, and in order to clarify this, and at a rather conceptual level, let us formulate, inspired by the representation theory material from chapter 4, the following definition:

\begin{definition}
Let $P(k,l)$ be the set of partitions between an upper row of $k$ points, and a lower row of $l$ points. Associated to any $\pi\in P(k,l)$ is the linear map
$$T_\pi(e_{i_1}\otimes\ldots\otimes e_{i_k})=\sum_{j_1\ldots j_l}\delta_\pi\begin{pmatrix}i_1&\ldots&i_k\\ j_1&\ldots&j_l\end{pmatrix}e_{j_1}\otimes\ldots\otimes e_{j_l}$$
between tensor powers of $\mathbb C^N$, called ``easy'', with the Kronecker type symbol on the right being given by $\delta_\pi=1$ when the indices fit, and $\delta_\pi=0$ otherwise.
\end{definition}

Observe the obvious connection with notion of easy group, from chapter 4, the point being that a closed subgroup $G\subset U_N$ is easy precisely when its Tannakian category $C_G=(C_G(k,l))$ with $C_G(k,l)\subset\mathcal L((\mathbb C^N)^k,(\mathbb C^N)^l)$ is spanned by easy maps. 

\bigskip

For our purposes, we will need a slight modification of Definition 13.21, as follows:

\begin{definition}
Associated to any partition $\pi\in P(2s,2s)$ is the linear map
$$\varphi_\pi(e_{a_1\ldots a_s,c_1\ldots c_s})=\sum_{b_1\ldots b_s}\sum_{d_1\ldots d_s}\delta_\pi\begin{pmatrix}a_1&\ldots&a_s&c_1&\ldots&c_s\\ b_1&\ldots&b_s&d_1&\ldots&d_s\end{pmatrix}e_{b_1\ldots b_s,d_1\ldots d_s}$$
obtained from $T_\pi$ by contracting all the tensors, via the operation
$$e_{i_1}\otimes\ldots\otimes e_{i_{2s}}\to e_{i_1\ldots i_s,i_{s+1}\ldots i_{2s}}$$
with $\{e_1,\ldots,e_N\}$ standing as usual for the standard basis of $\mathbb C^N$.
\end{definition}

In relation with our Wishart matrix considerations, the point is that the above linear map $\varphi_\pi$ can be viewed as a ``block-modification'' map, as follows:
$$\varphi_\pi:M_{N^s}(\mathbb C)\to M_{N^s}(\mathbb C)$$

As an illustration, let us discuss the case $s=1$. There are 15 partitions $\pi\in P(2,2)$, and among them, the most ``basic'' are the 4 partitions $\pi\in P_{even}(2,2)$. We have:

\begin{theorem}
The partitions $\pi\in P_{even}(2,2)$ are as follows,
$$\pi_1=\begin{bmatrix}\circ&\bullet\\ \circ&\bullet\end{bmatrix}\quad,\quad
\pi_2=\begin{bmatrix}\circ&\bullet\\ \bullet&\circ\end{bmatrix}\quad,\quad
\pi_3=\begin{bmatrix}\circ&\circ\\ \bullet&\bullet\end{bmatrix}\quad,\quad
\pi_4=\begin{bmatrix}\circ&\circ\\ \circ&\circ\end{bmatrix}$$
with the associated linear maps $\varphi_\pi:M_n(\mathbb C)\to M_n(\mathbb C)$ being as follows,
$$\varphi_1(A)=A\quad,\quad
\varphi_2(A)=A^t\quad,\quad
\varphi_3(A)=Tr(A)1\quad,\quad
\varphi_4(A)=A^\delta$$
and the associated square matrices $\Lambda_\pi\in M_n(\mathbb C)\otimes M_n(\mathbb C)$ being as follows,
$$\Lambda^1_{ab,cd}=\delta_{ab}\delta_{cd}\quad,\quad 
\Lambda^2_{ab,cd}=\delta_{ad}\delta_{bc}\quad,\quad
\Lambda^3_{ab,cd}=\delta_{ac}\delta_{bd}\quad,\quad
\Lambda^4_{ab,cd}=\delta_{abcd}$$
producing the main examples of block-modified matrices, from Proposition 13.18.
\end{theorem}

\begin{proof}
This is something elementary, coming from the formula in Definition 13.21. Indeed, in the case $s=1$, that we are interested in here, that formula becomes:
$$\varphi_\pi(e_{ac})=\sum_{bd}\delta_\pi\begin{pmatrix}a&c\\ b&d\end{pmatrix}e_{bd}$$

Now in the case of the 4 partitions in the statement, such maps are given by:
$$\varphi_1(e_{ac})=e_{ac}\quad,\quad 
\varphi_2(e_{ac})=e_{ca}\quad,\quad 
\varphi_3(e_{ac})=\delta_{ac}\sum_be_{bb}\quad,\quad
\varphi_4(e_{ac})=\delta_{ac}e_{aa}$$

Thus, we obtain the formulae in the statement. Regarding now the associated square matrices, appearing via $\Lambda_{ab,cd}=\varphi(e_{ac})_{bd}$, these are given by:
$$\Lambda^1_{ab,cd}=\delta_{ab}\delta_{cd}\quad,\quad 
\Lambda^2_{ab,cd}=\delta_{ad}\delta_{bc}\quad,\quad
\Lambda^3_{ab,cd}=\delta_{ac}\delta_{bd}\quad,\quad
\Lambda^4_{ab,cd}=\delta_{abcd}$$

Thus, we are led to the conclusions in the statement.
\end{proof}

Moving ahead now, we would like to study the distribution of the arbitrary block-modified Wishart matrices $\widetilde{W}=(id\otimes\varphi)W$. We will use as before the moment method. However, things will be more tricky in the present setting, and we will need:

\begin{definition}
The generalized colored moments of a random matrix 
$$W\in M_N(L^\infty(X))$$
with respect to a colored integer $e=e_1\ldots e_p$, and a permutation $\sigma\in S_p$, are the numbers
$$M^\sigma_e(W)=\frac{1}{N^{|\sigma|}}\,E\left(\sum_{i_1,\ldots,i_p}W^{e_1}_{i_1i_{\sigma(1)}}\ldots W^{e_p}_{i_pi_{\sigma(p)}}\right)$$
where $|\sigma|$ is the number of cycles of $\sigma$.
\end{definition}

This is something quite technical, in the spirit of the free probability and free cumulant theory in \cite{nsp}, that we will need in what follows. In order to understand how these generalized moments work, consider the standard cycle in $S_p$, namely:
$$\gamma=(1\to2\to\ldots\to p\to 1)$$

If we use this cycle $\gamma\in S_p$ as our permutation $\sigma\in S_p$ in the above definition, the corresponding generalized moment of a random matrix $W$ is then the usual moment:
$$M^\gamma_e(W)
=\frac{1}{N}\,E\left(\sum_{i_1,\ldots,i_p}W^{e_1}_{i_1i_2}\ldots W^{e_p}_{i_pi_1}\right)
=(E\circ tr)(W^{e_1}\ldots W^{e_p})$$

In general, we can decompose the computation of $M^\sigma_e(W)$ over the cycles of $\sigma$, and we obtain in this way a certain product of moments of $W$. See \cite{nsp}.

\bigskip

As a second illustration now, in relation with the usual square matrices, and more specifically with the square matrices $\Lambda\in M_n(\mathbb C)\otimes M_n(\mathbb C)$ as in Proposition 13.19, we have the following formula, that we will use many times in what follows:

\begin{proposition}
Given a usual square matrix, of composed size,
$$\Lambda\in M_n(\mathbb C)\otimes M_n(\mathbb C)$$
we have the following generalized moment formula,
$$(M^\sigma_e\otimes M^\tau_e)(\Lambda)=\frac{1}{n^{|\sigma|+|\tau|}}\sum_{i_1,\ldots, i_p}\sum_{j_1,\ldots,j_p}\Lambda_{i_1j_1,i_{\sigma(1)}j_{\tau(1)}}^{e_1}\ldots\ldots\Lambda_{i_pj_p,i_{\sigma(p)}j_{\tau(p)}}^{e_p}$$
valid for any two permutations $\sigma,\tau\in S_p$, and any colored integer $e=e_1\ldots e_p$.
\end{proposition}

\begin{proof}
This is something obvious, applying the construction in Definition 13.24 with $N=n^2$, $X=\{.\}$, $W=\Lambda$, and then making a tensor product of the corresponding moments $M^\sigma_e$, $M^\tau_e$, regarded as linear functionals on $M_n(\mathbb C)\otimes M_n(\mathbb C)$.
\end{proof}

Consider now the embedding $NC(p)\subset S_p$ obtained by ``cycling inside each block''. That is, each block $b=\{b_1,\ldots,b_k\}$ with $b_1<\ldots<b_k$ of a given noncrossing partition $\sigma\in NC(p)$ produces by definition the cycle $(b_1\ldots b_k)$ of the corresponding permutation $\sigma\in S_p$. Observe that the one-block partition $\gamma\in NC(p)$ corresponds in this way to the standard cycle $\gamma\in S_p$. Also, the number of blocks $|\sigma|$ of a partition $\sigma\in NC(p)$ corresponds to the number of cycles $|\sigma|$ of the corresponding permutation $\sigma\in S_p$. 

\bigskip

With these conventions, we have the following result, which generalizes our various Wishart matrix moment computations, that we did so far in this book:

\index{block-modified matrix}

\begin{theorem}
The asymptotic moments of a block-modified Wishart matrix 
$$\widetilde{W}=(id\otimes\varphi)W$$
with parameters $d,m,n\in\mathbb N$ as before, are given by the formula
$$\lim_{d\to\infty}M_e\left(\frac{\widetilde{W}}{d}\right)=\sum_{\sigma\in NC(p)}(mn)^{|\sigma|}(M^\sigma_e\otimes M^\gamma_e)(\Lambda)$$
where $\Lambda\in M_n(\mathbb C)\otimes M_n(\mathbb C)$ is the square matrix associated to $\varphi:M_n(\mathbb C)\to M_n(\mathbb C)$.
\end{theorem}

\begin{proof}
We use the formula for the matrix entries of $\widetilde{W}$, directly in terms of the matrix $\Lambda$ associated to the map $\varphi$, from Proposition 13.20, namely:
$$\widetilde{W}_{ia,jb}=\sum_{cd}\Lambda_{ca,db}W_{ic,jd}$$

By conjugating this formula, we obtain the following formula for the entries of the adjoint matrix $\widetilde{W}^*$, that we will need as well, in what follows:
$$\widetilde{W}_{ia,jb}^*
=\sum_{cd}\bar{\Lambda}_{db,ca}\bar{W}_{jd,ic}
=\sum_{cd}\Lambda^*_{ca,db}W_{ic,jd}$$

Thus, we have the following global formula, valid for any exponent $e\in\{1,*\}$:
$$\widetilde{W}_{ia,jb}^e=\sum_{cd}\Lambda^e_{ca,db}W_{ic,jd}$$

In order to compute the moments of $\widetilde{W}$, observe first that we have:
\begin{eqnarray*}
tr(\widetilde{W}^{e_1}\ldots\widetilde{W}^{e_p})
&=&\frac{1}{dn}\sum_{i_ra_r}\prod_s\widetilde{W}_{i_sa_s,i_{s+1}a_{s+1}}^{e_s}\\
&=&\frac{1}{dn}\sum_{i_ra_rc_rd_r}\prod_s\Lambda_{c_sa_s,d_sa_{s+1}}^{e_s}W_{i_sc_s,i_{s+1}d_s}\\
&=&\frac{1}{dn}\sum_{i_ra_rc_rd_rj_rb_r}\prod_s\Lambda_{c_sa_s,d_sa_{s+1}}^{e_s}Y_{i_sc_s,j_sb_s}\bar{Y}_{i_{s+1}d_s,j_sb_s}
\end{eqnarray*}

The average of the general term can be computed by the Wick rule, which gives:
$$E\left(\prod_sY_{i_sc_s,j_sb_s}\bar{Y}_{i_{s+1}d_s,j_sb_s}\right)
=\#\left\{\sigma\in S_p\Big|i_{\sigma(s)}=i_{s+1},c_{\sigma(s)}=d_s,j_{\sigma(s)}=j_s,b_{\sigma(s)}=b_s\right\}$$

Let us look now at the above sum. The $i,j,b$ indices range over sets having respectively $d,d,m$ elements, and they have to be constant under the action of $\sigma\gamma^{-1},\sigma,\sigma$. Thus when summing over these $i,j,b$ indices we simply obtain a factor as follows:
$$f=d^{|\sigma\gamma^{-1}|}d^{|\sigma|}m^{|\sigma|}$$

Thus, we obtain the following moment formula:
$$(E\circ tr)(\widetilde{W}^{e_1}\ldots\widetilde{W}^{e_p})
=\frac{1}{dn}\sum_{\sigma\in S_p}d^{|\sigma\gamma^{-1}|}(dm)^{|\sigma|}\sum_{a_rc_r}\prod_s\Lambda_{c_sa_s,c_{\sigma(s)}a_{s+1}}^{e_s}$$

On the other hand, we know from Proposition 13.25 that the generalized moments of the matrix $\Lambda\in M_n(\mathbb C)\otimes M_n(\mathbb C)$ are given by the following formula:
$$(M^\sigma_e\otimes M^\tau_e)(\Lambda)=\frac{1}{n^{|\sigma|+|\tau|}}\sum_{i_1\ldots i_p}\sum_{j_1\ldots j_p}\Lambda_{i_1j_1,i_{\sigma(1)}j_{\tau(1)}}^{e_1}\ldots\ldots\Lambda_{i_pj_p,i_{\sigma(p)}j_{\tau(p)}}^{e_p}$$

By combining the above two formulae, we obtain the following moment formula:
$$(E\circ tr)(\widetilde{W}^{e_1}\ldots\widetilde{W}^{e_p})
=\sum_{\sigma\in S_p}d^{|\sigma|+|\sigma\gamma^{-1}|-1}(mn)^{|\sigma|}(M^\sigma_e\otimes M^\gamma_e)(\Lambda)$$

We use now the standard fact, that we know well from before, that for $\sigma\in S_p$ we have an inequality as follows, with equality precisely when $\sigma\in NC(p)$:
$$|\sigma|+|\sigma\gamma^{-1}|\leq p+1$$

Thus with $d\to\infty$ the sum restricts over the partitions $\sigma\in NC(p)$, and we get:
$$\lim_{d\to\infty}M_e\big(\widetilde{W}\big)=d^p\sum_{\sigma\in NC(p)}(mn)^{|\sigma|}(M^\sigma_e\otimes M^\gamma_e)(\Lambda)$$

Thus, we are led to the conclusion in the statement.
\end{proof}

\section*{13d. Compound Poisson} 

In order to interpret the asymptotic moment formula that we found in Theorem 13.26, things are quite tricky, because at that level of generality, the resulting measure can vary a lot, generically escaping from the class of measures that we can efficiently study.

\bigskip

Still following \cite{bn1} and its continuations, our purpose in what follows will be that of identifying the situations where Theorem 13.26 produces compound free Poisson laws, or at least the obvious situations when this happens. Let us introduce:

\index{multiplicative matrix}

\begin{definition}
We call a matrix $\Lambda\in M_n(\mathbb C)\otimes M_n(\mathbb C)$ multiplicative when
$$(M^\sigma_e\otimes M^\gamma_e)(\Lambda)=(M^\sigma_e\otimes M^\sigma_e)(\Lambda)$$
holds for any $p\in\mathbb N$, any exponents $e_1,\ldots,e_p\in\{1,*\}$, and any $\sigma\in NC(p)$.
\end{definition}

Obviously, this is something quite technical, but we will see in a moment where the above condition comes from, the idea being that the combinatorics in the formula from Theorem 13.26 will heavily simplify, and we will be led to compound Poisson laws. Also, in what regards the examples, there are many of them, and we will discuss this later.

\bigskip

Getting back now to the block-modified Wishart matrices, and to the asymptotic moment formula that we found in Theorem 13.26, as just mentioned, the above notion of multiplicative matrix is exactly what we need, for further improving all this. Indeed, we have the following result, substantially building on Theorem 13.26:

\begin{theorem}
Consider a block-modified Wishart matrix 
$$\widetilde{W}=(id\otimes\varphi)W$$
and assume that the matrix $\Lambda\in M_n(\mathbb C)\otimes M_n(\mathbb C)$ associated to $\varphi$ is multiplicative. Then
$$\frac{\widetilde{W}}{d}\sim\pi_{mn\rho}$$
holds, in moments, in the $d\to\infty$ limit, where $\rho=law(\Lambda)$.
\end{theorem}

\begin{proof}
This is something quite tricky, using all the above:

\medskip

(1) Our starting point is the asymptotic moment formula found in Theorem 13.26, for an arbitrary block-modified Wishart matrix, namely:
$$\lim_{d\to\infty}M_e\left(\frac{\widetilde{W}}{d}\right)=\sum_{\sigma\in NC_p}(mn)^{|\sigma|}(M^\sigma_e\otimes M^\gamma_e)(\Lambda)$$

(2) Now since our modification matrix $\Lambda\in M_n(\mathbb C)\otimes M_n(\mathbb C)$ was assumed to be multiplicative, in the sense of Definition 13.27, this formula reads:
$$\lim_{d\to\infty}M_e\left(\frac{\widetilde{W}}{d}\right)=\sum_{\sigma\in NC_p}(mn)^{|\sigma|}(M^\sigma_e\otimes M^\sigma_e)(\Lambda)$$

(3) Our claim now, which will lead right away to the result, is that given an arbitrary square matrix $\Lambda\in M_n(\mathbb C)\otimes M_n(\mathbb C)$, multiplicative or not, having eigenvalue distribution $\rho=law(\Lambda)$, the moments of the corresponding compound free Poisson law $\pi_{mn\rho}$ are given by the following formula, for any choice of an extra parameter $m\in\mathbb N$:
$$M_e(\pi_{mn\rho})=\sum_{\sigma\in NC_p}(mn)^{|\sigma|}(M^\sigma_e\otimes M^\sigma_e)(\Lambda)$$

(4) Indeed, we know from chapter 9 that in the real case, the free cumulants of $\pi_{mn\rho}$ are the moments of $mn\rho$. But, as explained in \cite{nsp}, the same happens in the general complex case, in the sense that the free $*$-cumulants of $\pi_{mn\rho}$ are the $*$-moments of $mn\rho$. We conclude that these free $*$-cumulants are given by the following formula:
\begin{eqnarray*}
\kappa_p^e(\pi_{mn\rho})
&=&M^p_e(mn\rho)\\
&=&mn\cdot M^p_e(\Lambda)\\
&=&mn\cdot (M^\gamma_e\otimes M^\gamma_e)(\Lambda)
\end{eqnarray*}

But with this in hand, by using now Speicher's moment-cumulant formula, that we discussed in chapter 9 in the real case, and which holds in the general complex case too, as explained for instance in \cite{nsp}, we are led to the moment formula in (3).

\medskip

(5) The point now is that with that formula from (3) in hand, our previous asymptotic moment formula for the block-modified Wishart matrix $\widetilde{W}$ simply reads:
$$\lim_{d\to\infty}M_e\left(\frac{\widetilde{W}}{d}\right)=M_e(\pi_{mn\rho})$$

Thus we have indeed $\widetilde{W}/d\sim\pi_{mn\rho}$, in the $d\to\infty$ limit, as stated.
\end{proof}

Let us work out now some explicit consequences of Theorem 13.28, by using the modified easy linear maps from Definition 13.22. We recall from there that any modified easy linear map $\varphi_\pi$ can be viewed as a ``block-modification'' map, as follows:
$$\varphi_\pi:M_{N^s}(\mathbb C)\to M_{N^s}(\mathbb C)$$

In order to verify that the corresponding matrices $\Lambda_\pi$ are multiplicative, we will need to check that all the functions $\varphi(\sigma,\tau)=(M_\sigma^e\otimes M_\tau^e)(\Lambda_\pi)$ have the following property:
$$\varphi(\sigma,\gamma)=\varphi(\sigma,\sigma)$$

For this purpose, we can use the following elementary result:

\begin{proposition}
The following functions $\varphi:NC(p)\times NC(p)\to\mathbb R$ are multiplicative, in the sense that they satisfy the condition $\varphi(\sigma,\gamma)=\varphi(\sigma,\sigma)$:
\begin{enumerate}
\item $\varphi(\sigma,\tau)=|\sigma\tau^{-1}|-|\tau|$.

\item $\varphi(\sigma,\tau)=|\sigma\tau|-|\tau|$.

\item $\varphi(\sigma,\tau)=|\sigma\wedge\tau|-|\tau|$.
\end{enumerate}
\end{proposition}

\begin{proof}
All this is elementary, and can be proved as follows:

\medskip

(1) This follows indeed from the following computation:
$$\varphi_1(\sigma,\gamma)
=|\sigma\gamma^{-1}|-1
=p-|\sigma|
=\varphi_1(\sigma,\sigma)$$

(2) This follows indeed from the following computation:
$$\varphi_2(\sigma,\gamma)
=|\sigma\gamma|-1
=|\sigma^2|-|\sigma|
=\varphi_2(\sigma,\sigma)$$

(3) This follows indeed from the following computation:
$$\varphi_3(\sigma,\gamma)
=|\gamma|-|\gamma|
=0
=|\sigma|-|\sigma|
=\varphi_3(\sigma,\sigma)$$

Thus, we are led to the conclusions in the statement.
\end{proof}

We can get back now to the easy modification maps, and we have:

\begin{proposition}
The partitions $\pi\in P_{even}(2,2)$ are as follows,
$$\pi_1=\begin{bmatrix}\circ&\bullet\\ \circ&\bullet\end{bmatrix}\quad,\quad
\pi_2=\begin{bmatrix}\circ&\bullet\\ \bullet&\circ\end{bmatrix}\quad,\quad
\pi_3=\begin{bmatrix}\circ&\circ\\ \bullet&\bullet\end{bmatrix}\quad,\quad
\pi_4=\begin{bmatrix}\circ&\circ\\ \circ&\circ\end{bmatrix}$$
with the associated linear maps $\varphi_\pi:M_n(\mathbb C)\to M_N(\mathbb C)$ being as follows:
$$\varphi_1(A)=A\quad,\quad
\varphi_2(A)=A^t\quad,\quad
\varphi_3(A)=Tr(A)1\quad,\quad
\varphi_4(A)=A^\delta$$
The corresponding matrices $\Lambda_\pi$ are all multiplicative, in the sense of Definition 13.27.
\end{proposition}

\begin{proof}
In order to prove this, recall from Theorem 13.23 that the associated square matrices, appearing via $\Lambda_{ab,cd}=\varphi(e_{ac})_{bd}$, are given by:
$$\Lambda^1_{ab,cd}=\delta_{ab}\delta_{cd}\quad,\quad 
\Lambda^2_{ab,cd}=\delta_{ad}\delta_{bc}\quad,\quad
\Lambda^3_{ab,cd}=\delta_{ac}\delta_{bd}\quad,\quad
\Lambda^4_{ab,cd}=\delta_{abcd}$$

Since these matrices are all self-adjoint, we can assume that all the exponents are 1 in Definition 13.27, and the multiplicativity condition there becomes:
$$(M_\sigma\otimes M_\gamma)(\Lambda)=(M_\sigma\otimes M_\sigma)(\Lambda)$$

In order to check this condition, observe that for the above 4 matrices, we have:
\begin{eqnarray*}
(M^\sigma\otimes M^\tau)(\Lambda_1)&=&\frac{1}{n^{|\sigma|+|\tau|}}\sum_{i_1\ldots i_p}\delta_{i_{\sigma(1)}i_{\tau(1)}}\ldots\delta_{i_{\sigma(p)}i_{\tau(p)}}=n^{|\sigma\tau^{-1}|-|\sigma|-|\tau|}\\
(M^\sigma\otimes M^\tau)(\Lambda_2)&=&\frac{1}{n^{|\sigma|+|\tau|}}\sum_{i_1\ldots i_p}\delta_{i_1i_{\sigma\tau(1)}}\ldots\delta_{i_pi_{\sigma\tau(p)}}=n^{|\sigma\tau|-|\sigma|-|\tau|}\\
(M^\sigma\otimes M^\tau)(\Lambda_3)&=&\frac{1}{n^{|\sigma|+|\tau|}}\sum_{i_1\ldots i_p}\sum_{j_1\ldots j_p}\delta_{i_1i_{\sigma(1)}}\delta_{j_1j_{\tau(1)}}\ldots\delta_{i_pi_{\sigma(p)}}\delta_{j_pj_{\tau(p)}}=1\\
(M^\sigma\otimes M^\tau)(\Lambda_4)&=&\frac{1}{n^{|\sigma|+|\tau|}}\sum_{i_1\ldots i_p}\delta_{i_1i_{\sigma(1)}i_{\tau(1)}}\ldots\delta_{i_pi_{\sigma(p)}i_{\tau(p)}}=n^{|\sigma\wedge\tau|-|\sigma|-|\tau|}
\end{eqnarray*}

By using now the results in Proposition 13.29, this gives the result.
\end{proof}

Summarizing, the partitions $\pi\in P_{even}(2,2)$ provide us with some concrete input for Theorem 13.28. The point now is that, when using this input, we obtain the main known computations for the block-modified Wishart matrices, from \cite{aub}, \cite{cne}, \cite{mpa}:

\begin{theorem}
The asymptotic distribution results for the block-modified Wishart matrices coming from the partitions $\pi_1,\pi_2,\pi_3,\pi_4\in P_{even}(2,2)$ are as follows:
\begin{enumerate}
\item Marchenko-Pastur: $\frac{1}{d}W\sim\pi_t$, where $t=m/n$.

\item Aubrun type: $\frac{1}{d}(id\otimes t)W\sim\pi_\nu$, with $\nu=\frac{m(n-1)}{2}\,\delta_{-1}+\frac{m(n+1)}{2}\,\delta_1$. 

\item Collins-Nechita one: $n(id\otimes tr(.)1)W\sim\pi_t$, where $t=mn$.

\item Collins-Nechita two: $\frac{1}{d}(id\otimes(.)^\delta)W\sim\pi_m$.
\end{enumerate}
\end{theorem}

\begin{proof}
All these results follow from Theorem 13.28, with the maps $\varphi_1,\varphi_2,\varphi_3,\varphi_4$ in Proposition 13.30 producing the 4 matrices in the statement, modulo some rescalings:

\medskip

(1) Here $\Lambda=\sum_{ac}e_{ac}\otimes e_{ac}$, and so $\Lambda=nP$, where $P$ is the rank one projection on $\sum_ae_a\otimes e_a\in\mathbb C^n\otimes\mathbb C^n$. Thus we have the following formula, which gives the result:
$$\rho=\frac{n^2-1}{n^2}\,\delta_0+\frac{1}{n^2}\,\delta_n$$

(2) Here $\Lambda=\sum_{ac}e_{ac}\otimes e_{ca}$ is the flip operator, $\Lambda(e_c\otimes e_a)=e_a\otimes e_c$. Thus $\rho=\frac{n-1}{2n}\,\delta_{-1}+\frac{n+1}{2n}\,\delta_1$, and so we have the following formula, which gives the result:
$$mn\rho=\frac{m(n-1)}{2}\,\delta_{-1}+\frac{m(n+1)}{2}\,\delta_1$$

(3) Here $\Lambda=\sum_{ab}e_{aa}\otimes e_{bb}$ is the identity matrix, $\Lambda=1$. Thus in this case we have the following formula, which gives $\pi_{mn\rho}=\pi_{mn}$, and so $n\widetilde{W}\sim\pi_{mn}$, as claimed:
$$\rho=\delta_1$$

(4) Here $\Lambda=\sum_ae_{aa}\otimes e_{aa}$ is the orthogonal projection on $span(e_a\otimes e_a)\subset\mathbb C^n\otimes\mathbb C^n$. Thus we have the following formula, which gives the result:
$$\rho=\frac{n-1}{n}\,\delta_0+\frac{1}{n}\,\delta_1$$

Summarizing, we have proved all the assertions in the statement.
\end{proof}

And we will stop our discussion here. Of course, all this was just a beginning.

\section*{13e. Exercises}

Welcome to matrix models, good to have you here, and as exercises, we have:

\begin{exercise}
Carefully read Voiculescu et al. \cite{vdn}.
\end{exercise}

\begin{exercise}
Get initiated to the combinatorial yoga of Nica-Speicher \cite{nsp}.
\end{exercise}

\begin{exercise}
Learn some analysis MMA, from Bo\.zejko et al. \cite{bbl}, \cite{brs}, \cite{bwy}.
\end{exercise}

\begin{exercise}
Terminator training, Anderson, Guionnet and Zeitouni \cite{agz}.
\end{exercise}

\begin{exercise}
Learn more operator algebras, from Voiculescu \cite{vo5}.
\end{exercise}

\begin{exercise}
Master asymptotic freeness, from Mingo-Speicher \cite{msp}.
\end{exercise}

\begin{exercise}
Read \cite{aub}, \cite{bb+}, \cite{bn1} and their various follow-ups.
\end{exercise}

\begin{exercise}
Advance on the various questions left open above.
\end{exercise}

As bonus exercise, random matrices of any type, this is always good learning.

\chapter{Weyl models}

\section*{14a. Stationary models}

Moving on, let us get now to the second modeling question mentioned in the beginning of chapter 13, namely finding suitable matrix models for the coordinates $x_i$ of the free manifolds $X\subset S^{N-1}_{\mathbb C,+}$. We will restrict the attention to the quantum group case, where the basic theory takes place, and have as starting point the following definition:

\begin{definition}
A matrix model for $G\subset U_N^+$ is a morphism of $C^*$-algebras
$$\pi:C(G)\to M_K(C(T))\quad,\quad u_{ij}\to U_{ij}$$
where $K\geq1$ is an integer, and $T$ is a compact space.
\end{definition}

In other words, we are looking for random matrices $U_{ij}\in M_K(C(T))$ suitably modeling the quantum group coordinates $u_{ij}\in C(G)$. Which is something of obvious theoretical interest, mathematically speaking, abstract things becoming concrete, being the idea. As for the physics, do not worry, in analogy with what we did in chapters 11-12, having as underlying theory, and source of inspiration, the work of Connes on the Standard Model \cite{cco}, \cite{co3}, what we will be doing here will have as underlying theory, and source of inspiration, the work of Jones on statistical mechanics \cite{jo1}, \cite{jo2}, \cite{jo3}, \cite{jo4}, \cite{jo5}.

\bigskip

Getting now to Definition 14.1 as stated, we would like of course $\pi$ to provide us with useful information about $G$. And obviously, the simplest such case is when $\pi$ is faithful in the usual sense. However, this is something quite restrictive, because in this case the algebra $C(G)$ must be quite small, admitting an embedding as follows:
$$\pi:C(G)\subset M_K(C(T))$$

Technically, this means that $C(G)$ must be of type I, as an operator algebra, and we will discuss this in a moment, with the comment that this is indeed something quite restrictive. However, there are many interesting examples here, and all this is worth a detailed look. First, we have the following result, providing us with basic examples:

\begin{proposition}
The following closed subgroups $G\subset U_N^+$ have faithful models:
\begin{enumerate}
\item The compact Lie groups $G\subset U_N$.

\item The finite quantum groups $G\subset U_N^+$.
\end{enumerate}
In both cases, we can arrange for $\int_G$ to be restriction of the random matrix trace.
\end{proposition}

\begin{proof}
These assertions are all elementary, the proofs being as follows:

\medskip

(1) This is clear, because we can simply use here the identity map:
$$id:C(G)\to M_1(C(G))$$

(2) Here we can use the left regular representation $\lambda:C(G)\to M_{|G|}(\mathbb C)$. Indeed, let us endow the linear space $H=C(G)$ with the scalar product $<a,b>=\int_Gab^*$. We have then a representation of $*$-algebras, as follows:
$$\lambda:C(G)\to B(H)\quad,\quad 
a\to[b\to ab]$$

Now since we have $H\simeq\mathbb C^{|G|}$, we can view $\lambda$ as a matrix model map, as above. 

\medskip

(3) Finally, our claim is that we can choose our model as for the following formula to hold, where $\int_T$ is the integration with respect to a given probability measure on $T$:
$$\int_G=\left(tr\otimes\int_T\right)\pi$$

But this is clear for the model in (1), by definition, and is clear as well for the model in (2), by using the basic properties of the left regular representation.
\end{proof}

In the above result, the last assertion is quite interesting, and suggests formulating the following definition, somewhat independently on the notion of faithfulness:

\begin{definition}
A matrix model $\pi:C(G)\to M_K(C(T))$ is called stationary when
$$\int_G=\left(tr\otimes\int_T\right)\pi$$
where $\int_T$ is the integration with respect to a given probability measure on $T$.
\end{definition}

Here the term ``stationary'' comes from a functional analytic interpretation of this, with a certain Ces\`aro limit needed to be stationary, and we will explain this later. Getting now to the  examples, besides those in Proposition 14.2, we can look at group duals. So, consider a discrete group $\Gamma$, and a model for the corresponding group algebra:
$$\pi:C^*(\Gamma)\to M_K(C(T))$$

Since a representation of a group algebra must come from a unitary representation of the group, such a matrix model must come from a representation as follows:
$$\rho:\Gamma\to C(T,U_K)$$

With this identification made, we have the following result:

\begin{proposition}
An matrix model $\rho:\Gamma\subset C(T,U_K)$ is stationary when:
$$\int_Ttr(g^x)dx=0,\forall g\neq1$$
Moreover, the examples include all abelian groups, and all finite groups.
\end{proposition}

\begin{proof}
Consider indeed a group embedding $\rho:\Gamma\subset C(T,U_K)$, which produces by linearity a matrix model, as follows:
$$\pi:C^*(\Gamma)\to M_K(C(T))$$

It is enough to formulate the stationarity condition on the group elements $g\in C^*(\Gamma)$. Let us set $\rho(g)=(x\to g^x)$. With this notation, the stationarity condition reads:
$$\int_Ttr(g^x)dx=\delta_{g,1}$$

Since this equality is trivially satisfied at $g=1$, where by unitality of our representation we must have $g^x=1$ for any $x\in T$, we are led to the condition in the statement. Regarding now the examples, these are both clear. More precisely:

\medskip

(1) When $\Gamma$ is abelian we can use the following trivial embedding:
$$\Gamma\subset C(\widehat{\Gamma},U_1)\quad,\quad 
g\to[\chi\to\chi(g)]$$

(2) When $\Gamma$ is finite we can use the left regular representation:
$$\Gamma\subset\mathcal L(\mathbb C\Gamma)\quad,\quad 
g\to[h\to gh]$$

Indeed, in both cases, the stationarity condition is trivially satisfied.
\end{proof}

In order to discuss now certain analytic aspects of the matrix models, let us go back to the von Neumann algebras. We recall that we have the following result:

\begin{theorem}
Given a von Neumann algebra $A\subset B(H)$, if we write its center as 
$$Z(A)=L^\infty(X)$$
then we have a decomposition as follows, with the fibers $A_x$ having trivial center:
$$A=\int_XA_x\,dx$$
Moreover, the factors, $Z(A)=\mathbb C$, can be basically classified in terms of the ${\rm II}_1$ factors, which are those satisfying $\dim A=\infty$, and having a faithful trace $tr:A\to\mathbb C$.
\end{theorem}

\begin{proof}
This is something which is clear in finite dimensions, and in the commutative case too. In general, this is something heavy, the idea being as follows:

\medskip

(1) The first assertion, regarding the decomposition into factors, is von Neumann's reduction theory main result, which is actually one of the heaviest results in fundamental mathematics, and whose proof uses advanced functional analysis techniques. 

\medskip

(2) The classification of factors, due to Murray-von Neumann and Connes, is again something heavy, the idea being that the ${\rm II}_1$ factors are the ``building blocks'', with the other factors basically appearing from them via crossed product type constructions.
\end{proof}

Back now to matrix models, as a first general result, which is something which is not exactly trivial, and whose proof requires some functional analysis, we have:

\begin{theorem}
Assuming that a closed subgroup $G\subset U_N^+$ has a stationary model
$$\pi:C(G)\to M_K(C(T))$$
it follows that $G$ must be coamenable, and that the model is faithful. Moreover, $\pi$ extends into an embedding of von Neumann algebras, as follows,
$$L^\infty(G)\subset M_K(L^\infty(T))$$
which commutes with the canonical integration functionals.
\end{theorem}

\begin{proof}
Assume that we have a stationary model, as in the statement. By performing the GNS construction with respect to $\int_G$, we obtain a factorization as follows, which commutes with the respective canonical integration functionals:
$$\pi:C(G)\to C(G)_{red}\subset M_K(C(T))$$

Thus, in what regards the coamenability question, we can assume that $\pi$ is faithful. With this assumption made, we have an embedding as follows:
$$C(G)\subset M_K(C(T))$$

By performing the GNS construction we obtain a better embedding, as follows:
$$L^\infty(G)\subset M_K(L^\infty(T))$$

Now since the von Neumann algebra on the right is of type I, so must be its subalgebra $A=L^\infty(G)$. But this means that, when writing the center of this latter algebra as  $Z(A)=L^\infty(X)$, the whole algebra decomposes over $X$, as an integral of type I factors:
$$L^\infty(G)=\int_XM_{K_x}(\mathbb C)\,dx$$

In particular, we can see from this that $C(G)\subset L^\infty(G)$ has a unique $C^*$-norm, and so $G$ is coamenable. Thus we have proved our first assertion, and the second assertion follows as well, because our factorization of $\pi$ consists of the identity, and of an inclusion.
\end{proof}

In relation with the above, we have the following well-known result of Thoma:

\begin{theorem}
For a discrete group $\Gamma$, the following are equivalent:
\begin{enumerate}
\item $C^*(\Gamma)$ is of type I, so that we have an embedding $\pi:C^*(\Gamma)\subset M_K(C(X))$, with $X$ being a compact space.

\item $C^*(\Gamma)$ has a stationary model of type $\pi:C^*(\Gamma)\to M_F(C(L))$, with $F$ being a finite group, and $L$ being a compact abelian group.

\item $\Gamma$ is virtually abelian, in the sense that we have an abelian subgroup $\Lambda\triangleleft\Gamma$ such that the quotient group $F=\Gamma/\Lambda$ is finite. 

\item $\Gamma$ has an abelian subgroup $\Lambda\subset\Gamma$ whose index $K=[\Gamma:\Lambda]$ is finite.
\end{enumerate}
\end{theorem}

\begin{proof}
There are several proofs for this fact, the idea being as follows:

\medskip

$(1)\implies(4)$ This is the non-trivial implication, that we will not prove here. We refer instead to the literature, either Thoma's orignal paper, or books like those of Dixmier, mixing advanced group theory and advanced operator algebra theory.

\medskip

$(4)\implies(3)$ We choose coset representatives $g_i\in\Gamma$, and we set:
$$\Lambda'=\bigcap_ig_i\Gamma g_i^{-1}$$

Then $\Lambda'\subset\Lambda$ has finite index, and we have $\Lambda'\triangleleft\Gamma$, as desired.

\medskip

$(3)\implies(2)$ This follows by using the theory of induced representations. We can define a model $\pi:C^*(\Gamma)\to M_F(C(\widehat{\Lambda}))$ by setting:
$$\pi(g)(\chi)=Ind_\Lambda^\Gamma(\chi)(g)$$

Indeed, any character $\chi\in\widehat{\Lambda}$ is a 1-dimensional representation of $\Lambda$, and we can therefore consider the induced representation $Ind_\Lambda^\Gamma(\chi)$ of the group $\Gamma$. This representation is $|F|$-dimensional, and so maps the group elements $g\in\Gamma$ into order $|F|$ matrices $Ind_\Lambda^\Gamma(\chi)(g)$. Thus the above map $\pi$ is well-defined, and the fact that it is a representation is clear as well. In order to check now the stationarity property of this representation, we can use the following well-known character formula, due to Frobenius:
$$Tr\left(Ind_\Lambda^\Gamma(\chi)(g)\right)=\sum_{x\in F}\delta_{x^{-1}gx\in\Lambda}\chi(x^{-1}gx)$$

By integrating with respect to $\chi\in\widehat{\Lambda}$, we deduce from this that we have:
\begin{eqnarray*}
\left(Tr\otimes \int_{\widehat{\Lambda}}\right)\pi(g)
&=&\sum_{x\in F}\delta_{x^{-1}gx\in \Lambda}\int_{\widehat{\Lambda}}\chi(x^{-1}gx)d\chi\\
&=&\sum_{x\in F}\delta_{x^{-1}gx\in \Lambda}\delta_{g,1}\\
&=&|F|\cdot\delta_{g,1}
\end{eqnarray*}

Now by dividing by $|F|$ we conclude that the model is stationary, as claimed.

\medskip

$(2)\implies(1)$ This is the trivial implication, with the faithfulness of $\pi$ following from the abstract functional analysis arguments from the proof of Theorem 14.6.
\end{proof}

In the quantum group case the extension of Theorem 14.7 is something quite tricky, and for more on this, you can have a look at my papers with Chirvasitu, and with Bichon too, and their follow-ups. Alternatively, check my book on quantum permutations.

\bigskip

Getting back now to Definition 14.1, more generally, we can model in that way the standard coordinates $x_i\in C(X)$ of various algebraic manifolds $X\subset S^{N-1}_{\mathbb C,+}$. Indeed, these manifolds generalize the compact matrix quantum groups, which appear as:
$$G\subset U_N^+\subset S^{N^2-1}_{\mathbb C,+}$$

Moreover, we have many other interesting examples of such manifolds, such as the homogeneous spaces discussed in chapter 12. However, at this level of generality, not much general theory is available. It is elementary to show that, under the technical assumption $X^{class}\neq\emptyset$, there exists a universal $K\times K$ model for the algebra $C(X)$, which factorizes as follows, with $X^{(K)}\subset X$ being a certain algebraic submanifold: 
$$\pi_K:C(X)\to C(X^{(K)})\subset M_K(C(T_K))$$

To be more precise, the universal $K\times K$ model space $T_K$ appears by imposing to the complex $K\times K$ matrices the relations defining $X$, and the algebra $C(X^{(K)})$ is then by definition the image of $\pi_K$. In relation with this, we can set as well:
$$X^{(\infty)}=\bigcup_{K\in\mathbb N}X^{(K)}$$

We are led in this way to a filtration of $X$, as follows:
$$X^{class}= X^{(1)}\subset X^{(2)}\subset X^{(3)}\subset\ldots\ldots\subset X^{(\infty)}\subset X$$

It is possible to say a few non-trivial things about these manifolds $X^{(K)}$. In the compact quantum group case, however, that we are mainly interested in here, the matrix truncations $G^{(K)}\subset G$ are generically not quantum subgroups at $K\geq2$, and so this theory is a priori not very useful, at least in its basic form presented above.

\bigskip

For more on this, and notably for some concrete results regarding the half-liberations, which are related to the second matrix truncations $X^{(2)}$ constructed above, have a look at the paper of Bichon and Dubois-Violette \cite{bdu}, and its various follow-ups.

\section*{14b. Inner faithfulness}

Let us discuss now the general, non-coamenable case, with the aim of finding a weaker notion of faithfulness, which still does the job, namely that of ``reminding'' the quantum group. The idea comes by looking at the group duals $G=\widehat{\Gamma}$. Consider indeed a general model for the associated group algebra, which can be written as follows:
$$\pi:C^*(\Gamma)\to M_K(C(T))$$

The point is that such a representation of the group algebra must come by linearization from a unitary group representation, as follows: 
$$\rho:\Gamma\to C(T,U_K)$$

Now observe that when this group representation $\rho$ is faithful, the representation $\pi$ is in general not faithful, for instance because when $T=\{.\}$ its target algebra is finite dimensional. On the other hand, this representation ``reminds'' $\Gamma$, so can be used in order to fully understand $\Gamma$. Thus, we have an idea here, basically saying that, for practical purposes, the faithfuless property can be replaced with something much weaker. 

\bigskip

This weaker notion, which will be of great interest for us, is called ``inner faithfulness''. The general theory here starts with the following definition:

\begin{definition}
Let $\pi:C(G)\to M_K(C(T))$ be a matrix model. 
\begin{enumerate}
\item The Hopf image of $\pi$ is the smallest quotient Hopf $C^*$-algebra $C(G)\to C(H)$ producing a factorization as follows:
$$\pi:C(G)\to C(H)\to M_K(C(T))$$

\item When the inclusion $H\subset G$ is an isomorphism, i.e. when there is no non-trivial factorization as above, we say that $\pi$ is inner faithful.
\end{enumerate}
\end{definition}

The above notions are quite tricky, and having them well understood will take us some time. As a first example, motivated by the above discussion, in the case where $G=\widehat{\Gamma}$ is a group dual, $\pi$ must come from a group representation, as follows:
$$\rho:\Gamma\to C(T,U_K)$$

Thus the minimal factorization in (1) is obtained by taking the image:
$$\rho:\Gamma\to\Lambda\subset C(T,U_K)$$

Thus, as a conclusion, in this case $\pi$ is inner faithful precisely when we have:
$$\Gamma\subset C(T,U_K)$$

Dually now, given a compact Lie group $G$, and elements $g_1,\ldots,g_K\in G$, we have a diagonal representation $\pi:C(G)\to M_K(\mathbb C)$, appearing as follows:
$$f\to\begin{pmatrix}
f(g_1)\\
&\ddots\\
&&f(g_K)
\end{pmatrix}$$

The minimal factorization of this representation $\pi$, as in Definition 16.34 (1), is then via the algebra $C(H)$, with $H$ being the following closed subgroup of $G$:
$$H=\overline{<g_1,\ldots,g_K>}$$

Thus, as a conclusion, $\pi$ is inner faithful precisely when we have:
$$G=H$$

There are many other examples of inner faithful representations, which are however substantially more technically advanced, and we will discuss them later. 

\bigskip

Back to general theory now, in the framework of Definition 14.8, the existence and uniqueness of the Hopf image come by dividing $C(G)$ by a suitable ideal, with this being something standard. Alternatively, in Tannakian terms, we have:

\begin{theorem}
Assuming $G\subset U_N^+$, with fundamental corepresentation $u=(u_{ij})$, the Hopf image of a model $\pi:C(G)\to M_K(C(T))$ comes from the Tannakian category
$$C_{kl}=Hom(U^{\otimes k},U^{\otimes l})$$
where $U_{ij}=\pi(u_{ij})$, and where the spaces on the right are taken in a formal sense.
\end{theorem}

\begin{proof}
Since the morphisms increase the intertwining spaces, when defined either in a representation theory sense, or just formally, we have inclusions as follows:
$$Hom(u^{\otimes k},u^{\otimes l})\subset Hom(U^{\otimes k},U^{\otimes l})$$

More generally, we have such inclusions when replacing $(G,u)$ with any pair producing a factorization of $\pi$. Thus, by Tannakian duality, the Hopf image must be given by the fact that the intertwining spaces must be the biggest, subject to the above inclusions. On the other hand, since $u$ is biunitary, so is $U$, and it follows that the spaces on the right form a Tannakian category. Thus, we have a quantum group $(H,v)$ given by:
$$Hom(v^{\otimes k},v^{\otimes l})=Hom(U^{\otimes k},U^{\otimes l})$$

By the above discussion, $C(H)$ follows to be the Hopf image of $\pi$, as claimed.
\end{proof}

Regarding now the study of the inner faithful models, a key problem is that of computing the Haar integration functional. The result here is as follows:

\begin{theorem}
Given an inner faithful model $\pi:C(G)\to M_K(C(T))$, we have
$$\int_G=\lim_{k\to\infty}\frac{1}{k}\sum_{r=1}^k\int_G^r\quad,\quad 
\int_G^r=(\varphi\circ\pi)^{*r}$$
where $\phi*\psi=(\phi\otimes\psi)\Delta$, and with $\varphi=tr\otimes\int_T$ being the random matrix trace.
\end{theorem}

\begin{proof}
This is something quite tricky, the idea being as follows:

\medskip

(1) As a first observation, there is an obvious similarity here with the Woronowicz construction of the Haar measure, from chapter 11. In fact, the above result holds more generally for any model $\pi:C(G)\to B$, with $\varphi\in B^*$ being a faithful trace.

\medskip

(2) In order to prove now the result, inspired by this, if we denote by $\int_G'$ the limit in the statement, we must prove that this limit converges, and that:
$$\int_G'=\int_G$$

It is enough to check this on the coefficients of the Peter-Weyl corepresentations, and if we let $v=u^{\otimes k}$ be one of these corepresentations, we must prove that we have:
$$\left(id\otimes\int_G'\right)v=\left(id\otimes\int_G\right)v$$

(3) In order to prove this, we already know, from the general Haar measure theory, that the matrix on the right is the orthogonal projection onto $Fix(v)$:
$$\left(id\otimes\int_G\right)v=Proj\Big[Fix(v)\Big]$$

Regarding now the matrix on the left, the usual trick applied to the linear form $\varphi\pi$ tells us that this is the orthogonal projection onto the $1$-eigenspace of $(id\otimes\varphi\pi)v$:
$$\left(id\otimes\int_G'\right)v=Proj\Big[1\in (id\otimes\varphi\pi)v\Big]$$

(4) Now observe that, if we set $V_{ij}=\pi(v_{ij})$, we have the following formula:
$$(id\otimes\varphi\pi)v=(id\otimes\varphi)V$$

Thus, we can apply the usual trick, and we conclude that the $1$-eigenspace that we are interested in equals $Fix(V)$. But, according to Theorem 14.9, we have:
$$Fix(V)=Fix(v)$$

Thus, we have proved that we have $\int_G'=\int_G$, as desired.
\end{proof}

In practice, Theorem 14.10 is something quite powerful. As an illustration, regarding the law of the main character, we obtain here the following result:

\begin{proposition}
Given $\pi:C(G)\to M_K(C(T))$ inner faithful, let $\mu=law(\chi)$, and let $\mu^r$ be the law of $\chi$ with respect to $\int_G^r=(\varphi\circ\pi)^{*r}$, where $\varphi=tr\otimes\int_T$. Then
$$\mu=\lim_{k\to\infty}\frac{1}{k}\sum_{r=0}^k\mu^r$$
with this convergence being in moments, and with the notation
$$(T_\varepsilon)_{i_1\ldots i_p,j_1\ldots j_p}=\left(tr\otimes\int_T\right)(U_{i_1j_1}^{\varepsilon_1}\ldots U_{i_pj_p}^{\varepsilon_p})$$
the moments of $\mu^r$ are the numbers $c_\varepsilon^r=Tr(T_\varepsilon^r)$.
\end{proposition}

\begin{proof}
The above two formulae are both elementary, by using the general convergence result established in Theorem 14.10. Indeed, the first formula comes from the limiting formula in Theorem 14.10, by applying the linear forms there to the main character $\chi$. As for the second formula, this comes from the definitions of the measure $\mu^r$ and of the matrix $T_e$, by summing the entries of $T_e$ over equal indices, $i_r=j_r$.
\end{proof}

Interestingly, the above results regarding inner faithfulness have applications as well to the notion of stationarity introduced before, clarifying among others the use of the word ``stationary''. To be more precise, in order to detect the stationary models, we have the following useful criterion, mixing linear algebra and analysis:

\begin{theorem}
For a model $\pi:C(G)\to M_K(C(T))$, the following are equivalent:
\begin{enumerate}
\item $Im(\pi)$ is a Hopf algebra, and the Haar integration on it is:
$$\psi=\left(tr\otimes\int_T\right)\pi$$

\item The linear form $\psi=(tr\otimes\int_T)\pi$ satisfies the idempotent state property:
$$\psi*\psi=\psi$$

\item We have $T_e^2=T_e$, $\forall p\in\mathbb N$, $\forall e\in\{1,*\}^p$, where:
$$(T_e)_{i_1\ldots i_p,j_1\ldots j_p}=\left(tr\otimes\int_T\right)(U_{i_1j_1}^{e_1}\ldots U_{i_pj_p}^{e_p})$$
\end{enumerate}
If these conditions are satisfied, we say that $\pi$ is stationary on its image.
\end{theorem}

\begin{proof}
Given a matrix model $\pi:C(G)\to M_K(C(T))$ as in the statement, we can factorize it via its Hopf image, as in Definition 14.8:
$$\pi:C(G)\to C(H)\to M_K(C(T))$$

Now observe that (1,2,3) above depend only on the factorized representation:
$$\nu:C(H)\to M_K(C(T))$$

Thus, we can assume in practice that we have $G=H$, which means that we can assume that $\pi$ is inner faithful. With this assumption made, the formula in Theorem 14.10 applies to our situation, and the proof of the equivalences goes as follows:

\medskip

$(1)\implies(2)$ This is clear from definitions, because the Haar integration on any compact quantum group satisfies the idempotent state equation:
$$\psi*\psi=\psi$$

$(2)\implies(1)$ Assuming $\psi*\psi=\psi$, we have $\psi^{*r}=\psi$ for any $r\in\mathbb N$, and Theorem 14.10 gives $\int_G=\psi$. By using now Theorem 14.6, we obtain the result.

\medskip

In order to establish now $(2)\Longleftrightarrow(3)$, we can use the following elementary formula, which comes from the definition of the convolution operation:
$$\psi^{*r}(u_{i_1j_1}^{e_1}\ldots u_{i_pj_p}^{e_p})=(T_e^r)_{i_1\ldots i_p,j_1\ldots j_p}$$

\medskip

$(2)\implies(3)$ Assuming $\psi*\psi=\psi$, by using the above formula at $r=1,2$ we obtain that the matrices $T_e$ and $T_e^2$ have the same coefficients, and so they are equal.

\medskip

$(3)\implies(2)$ Assuming $T_e^2=T_e$, by using the above formula at $r=1,2$ we obtain that the linear forms $\psi$ and $\psi*\psi$ coincide on any product of coefficients $u_{i_1j_1}^{e_1}\ldots u_{i_pj_p}^{e_p}$. Now since these coefficients span a dense subalgebra of $C(G)$, this gives the result.
\end{proof}

As an illustration, let us go back to the half-liberations $O_N^*,U_N^*$ briefly discussed in chapter 11, appearing via the relations $abc=cba$. We first have the following result:

\index{half-liberation}

\begin{proposition}
We have a matrix model as follows, 
$$C(O_N^*)\to M_2(C(U_N))\quad,\quad 
u_{ij}\to\begin{pmatrix}0&v_{ij}\\ \bar{v}_{ij}&0\end{pmatrix}$$
where $v$ is the fundamental corepresentation of $C(U_N)$, as well as a model as follows,
$$C(U_N^*)\to M_2(C(U_N\times U_N))\quad,\quad 
u_{ij}\to\begin{pmatrix}0&v_{ij}\\ w_{ij}&0\end{pmatrix}$$
where $v,w$ are the fundamental corepresentations of the two copies of $C(U_N)$.
\end{proposition}

\begin{proof}
It is routine to check that the matrices on the right are indeed biunitaries, and since the first matrix is also self-adjoint, we obtain models as follows:
$$C(O_N^+)\to M_2(C(U_N))\quad,\quad 
C(U_N^+)\to M_2(C(U_N\times U_N))$$

Consider now antidiagonal $2\times2$ matrices, with commuting entries, as follows:
$$X_i=\begin{pmatrix}0&x_i\\ y_i&0\end{pmatrix}$$

We have then $X_iX_jX_k=X_kX_jX_i$, so our models above factorize, as claimed.
\end{proof}

Following \cite{bn2}, \cite{bdu}, we can now formulate a first concrete theorem, as follows:

\begin{theorem}
The standard antidiagonal models for half-liberations,
$$C(O_N^*)\to M_2(C(U_N))\quad,\quad 
C(U_N^*)\to M_2(C(U_N\times U_N))$$
are both stationary, and in particular they are faithful.
\end{theorem}

\begin{proof}
Let us first discuss the case of $O_N^*$. We use Theorem 14.12 (3). Since the fundamental representation is self-adjoint, the various matrices $T_e$ with $e\in\{1,*\}^p$ are all equal. We denote this common matrix by $T_p$. We have, with $H=U_N$:
$$(T_p)_{i_1\ldots i_p,j_1\ldots j_p}
=\left(tr\otimes\int_H\right)\left[\begin{pmatrix}0&v_{i_1j_1}\\\bar{v}_{i_1j_1}&0\end{pmatrix}\ldots\ldots\begin{pmatrix}0&v_{i_pj_p}\\\bar{v}_{i_pj_p}&0\end{pmatrix}\right]$$

Since when multipliying an odd number of antidiagonal matrices we obtain an atidiagonal matrix, we have $T_p=0$ for $p$ odd. Also, when $p$ is even, we have:
\begin{eqnarray*}
(T_p)_{i_1\ldots i_p,j_1\ldots j_p}
&=&\left(tr\otimes\int_H\right)\begin{pmatrix}v_{i_1j_1}\ldots\bar{v}_{i_pj_p}&0\\0&\bar{v}_{i_1j_1}\ldots v_{i_pj_p}\end{pmatrix}\\
&=&\frac{1}{2}\left(\int_Hv_{i_1j_1}\ldots\bar{v}_{i_pj_p}+\int_H\bar{v}_{i_1j_1}\ldots v_{i_pj_p}\right)\\
&=&\int_HRe(v_{i_1j_1}\ldots\bar{v}_{i_pj_p})
\end{eqnarray*}

We have $T_p^2=T_p=0$ when $p$ is odd, so we are left with proving that for $p$ even we have $T_p^2=T_p$. For this purpose, we use the following formula:
$$Re(x)Re(y)=\frac{1}{2}\left(Re(xy)+Re(x\bar{y})\right)$$

By using this identity for each of the terms which appear in the product, and multi-index notations in order to simplify the writing, we obtain:
\begin{eqnarray*}
(T_p^2)_{ij}
&=&\sum_{k_1\ldots k_p}(T_p)_{i_1\ldots i_p,k_1\ldots k_p}(T_p)_{k_1\ldots k_p,j_1\ldots j_p}\\
&=&\int_H\int_H\sum_{k_1\ldots k_p}Re(v_{i_1k_1}\ldots\bar{v}_{i_pk_p})Re(w_{k_1j_1}\ldots\bar{w}_{k_pj_p})dvdw\\
&=&\frac{1}{2}\int_H\int_H\sum_{k_1\ldots k_p}Re(v_{i_1k_1}w_{k_1j_1}\ldots\bar{v}_{i_pk_p}\bar{w}_{k_pj_p})+Re(v_{i_1k_1}\bar{w}_{k_1j_1}\ldots\bar{v}_{i_pk_p}w_{k_pj_p})dvdw\\
&=&\frac{1}{2}\int_H\int_HRe((vw)_{i_1j_1}\ldots(\bar{v}\bar{w})_{i_pj_p})+Re((v\bar{w})_{i_1j_1}\ldots(\bar{v}w)_{i_pj_p})dvdw
\end{eqnarray*}

Now since $vw\in H$ is uniformly distributed when $v,w\in H$ are uniformly distributed, the quantity on the left integrates up to $(T_p)_{ij}$. Also, since $H$ is conjugation-stable, $\bar{w}\in H$ is uniformly distributed when $w\in H$ is uniformly distributed, so the quantity on the right integrates up to the same quantity, namely $(T_p)_{ij}$. Thus, we have:
$$(T_p^2)_{ij}
=\frac{1}{2}\Big((T_p)_{ij}+(T_p)_{ij}\Big)
=(T_p)_{ij}$$

Summarizing, we have obtained that for any $p$, we have $T_p^2=T_p$. Thus Theorem 14.12 applies, and shows that our model is stationary, as claimed. As for the proof of the stationarity for the standard antidiagonal model for $U_N^*$, this is similar.
\end{proof}

Many other things can be said about half-liberation, as a continuation of the above, and we refer here to the paper \cite{bdu}, and its various follow-ups.

\section*{14c. Pauli matrices}

Following \cite{bco}, we would like to discuss now a key matrix model for the quantum group $S_4^+$, coming from the Pauli matrices. These matrices are as follows:
$$c_1=\begin{pmatrix}1&0\\ 0&1\end{pmatrix}\quad,\quad
c_2=\begin{pmatrix}i&0\\ 0&-i\end{pmatrix}\quad,\quad
c_3=\begin{pmatrix}0&1\\ -1&0\end{pmatrix}\quad,\quad
c_4=\begin{pmatrix}0&i\\ i&0\end{pmatrix}$$

Observe that these matrices are elements of $SU_2$. In fact, $SU_2$ consists of linear combinations of Pauli matrices, with points on the real sphere $S^3$ as coefficients:
$$SU_2=\left\{\sum x_ic_i\,\Big| \sum|x_i|^2=1\right\}$$

The starting remark is that the Pauli matrices form an orthonormal basis of $M_2(\mathbb C)$, with respect to the scalar product $<a,b>=tr(ab^*)$. Moreover, the same is true if we multiply them to the left or to the right by an element of $SU_2$. Thus, we have:

\begin{proposition}
For any $x\in SU_2$ the elements 
$$\xi_{ij}=c_ixc_j$$
form a magic basis of $M_2({\mathbb C})$, with respect to the scalar product $<a,b>=tr(ab^*)$.
\end{proposition}

\begin{proof}
This follows indeed from the above discussion.
\end{proof}

We fix a Hilbert space isomorphism $M_2({\mathbb C})\simeq {\mathbb C}^4$, and we use the corresponding identification of operator algebras $B(M_2({\mathbb C}))\simeq M_4({\mathbb C})$. Associated to each element $x\in SU_2$ is the representation of $C(S_4^+)$ mapping $u_{ij}$ to the rank one projection on $c_ixc_j$:
$$\pi_x:C(S_4^+)\to M_4(\mathbb C)\quad,\quad u_{ij}\to Proj(c_ixc_j)$$

This representation depends on $x$. For getting a faithful representation, the idea is to regard all these representations as fibers of a single representation:

\begin{definition}
The Pauli representation of $C(S_4^+)$ is the map
$$\pi:C(S_4^+)\to C(SU_2, M_4({\mathbb C}))$$
mapping each standard coodinate $u_{ij}$ to the function $x\to Proj(c_ixc_j)$.
\end{definition}

Observe that we can replace here $SU_2$ by its quotient $PU_2=SO_3$. For reasons that will become clear later on, we prefer to use $SU_2$. We will prove now that $\pi$ is faithful, in the usual sense. Consider the standard linear form on the algebra on the right:
$$\int\varphi=\int_{SU_2}tr\left(\varphi(x)\right)\, dx$$

We use the following analytic formulation of faithfulness:

\begin{proposition}
The representation $\pi$ is faithful provided that
$$\int u_{i_1j_1}\ldots u_{i_kj_k}=
\int \pi_{i_1j_1}\ldots \pi_{i_kj_k}$$ for any choice of $k$ and
of various $i,j$ indices, where $\pi_{ij}=\pi(u_{ij})$.
\end{proposition}

\begin{proof}
This is indeed something obvious, coming from definitions.
\end{proof}

Our first goal will be to compute the integral on the right in the above statement. For this purpose, consider the canonical action of $SU_2$ on the algebra $M_2(\mathbb C)^{\otimes k}$, obtained as $k$-th tensor power of the adjoint action on $M_2(\mathbb C)$, as follows:
$$\alpha_x(a_1\otimes\ldots\otimes a_k)=xa_1x^*\otimes\ldots\otimes xa_kx^*$$

With this convention, we have the following key result:

\begin{proposition}
With $c_i=c_{i_1}\otimes\ldots\otimes c_{i_k}$, we have the following formula,
$$\int \pi_{i_1j_1}\ldots \pi_{i_kj_k}=<c_i,R^*ER(c_j)>$$ 
where $E$ is the expectation under the canonical action of $SU_2$, and where
$$R(c_{i_1}\otimes\ldots\otimes  c_{i_k})=\frac{1}{2}(c_{i_1}c_{i_2}^*
\otimes c_{i_2}c_{i_3}^*\otimes\ldots\otimes c_{i_k}c_{i_1}^*)$$
with the convention that at $k=1$ we have $R(c_i)=c_ic_i^*/2=1/2$.
\end{proposition}

\begin{proof}
We have indeed the following computation:
\begin{eqnarray*}
\int \pi_{i_1j_1}\ldots \pi_{i_kj_k} 
&=& \int\pi(u_{i_1j_1})\ldots \pi(u_{i_kj_k})\\
&=&\int Proj(c_{i_1}xc_{j_1})\ldots Proj(c_{i_k}xc_{j_k})\\
&=&\int_{SU_2}tr\left( Proj(c_{i_1}xc_{j_1})\ldots Proj(c_{i_k}xc_{j_k})\right)\,dx\\
&=&\frac{1}{4}\int_{SU_2}<c_{i_1}xc_{j_1},c_{i_k}xc_{j_k}>\ldots
<c_{i_2}xc_{j_2},c_{i_1}xc_{j_1}>\,dx\\
&=&\frac{1}{4}\int_{SU_2}tr(c_{i_1}xc_{j_1}c_{j_k}^*x^*c_{i_k}^*)\ldots tr(c_{i_2}xc_{j_2}c_{j_1}^*x^*c_{i_1}^*)\,dx \\
&=&\frac{1}{4}\int_{SU_2}tr(c_{i_k}^*c_{i_1}xc_{j_1}c_{j_k}^*x^*)\ldots tr(c_{i_1}^*c_{i_2}xc_{j_2}c_{j_1}^*x^*)\,dx\\
&=&\frac{1}{4}\int_{SU_2}tr(c_{i_k}c_{i_1}^*xc_{j_1}c_{j_k}^*x^*)\ldots 
tr(c_{i_1}c_{i_2}^*xc_{j_2}c_{j_1}^*x^*)\,dx\\
&=&\frac{1}{4}\int_{SU_2}tr\left(c_{i_1}c_{i_2}^*xc_{j_2}c_{j_1}^*x^*
\otimes\ldots\otimes c_{i_k}c_{i_1}^*xc_{j_1}c_{j_k}^*x^*\right)\,dx\\
&=&\int_{SU_2}tr\left( R(c_i)\alpha_x(R(c_j)^*)\right)\,dx\\
&=&tr\left(\int_{SU_2}R(c_i)\alpha_x(R(c_j)^*)\, dx\right)\\
&=&tr\left(R(c_i)\int_{SU_2} \alpha_x(R(c_j)^*)\,dx\right)\\
&=&tr\left(R(c_i)ER(c_j)^*\right)\\
&=&<R(c_i),ER(c_j)>\\
&=&<c_i,R^*ER(c_j)>
\end{eqnarray*}

Thus, we are led to the conclusion in the statement.
\end{proof}

Next, we have the following result, which provides the key to faithfulness:

\begin{proposition}
Given $\nu\in NC(k)$, the following element of $M_2(\mathbb C)^k$,
$$c_\nu=\sum_j\delta_\nu(j)c_j$$
satisfies the formula $R^*ER(c_\nu)=c_\nu$.
\end{proposition}

\begin{proof}
This is something quite long and routine, the idea being as follows:

\medskip

(1) Our first claim is that $f=\sum_ic_i\otimes c_i^*$ is invariant under $SU_2$. Indeed, we have $f=2(c_1\otimes c_1)-\sum_ic_i\otimes c_i$, and since $c_1\otimes c_1$ is invariant under $SU_2$, what is left to prove is that the sum on the right is invariant as well. This sum, viewed as a matrix, is:
$$\sum_ic_i\otimes c_i
=\begin{pmatrix}
 0&0&0&0\\
0&2&-2&0\\
0&-2&2&0\\
0&0&0&0
\end{pmatrix}$$

But this matrix is $4$ times the projection onto the determinant subspace of 
$\mathbb C^{2}\otimes\mathbb C^2$, which is invariant under the action of $SU_2$, and this concludes the proof of our claim.

\medskip

(2) Next, consider the trivial partition, $0_k=\{\{1\},\ldots ,\{k\}\}$. We have then:
\begin{eqnarray*}
R(c_{0_k})
&=&R\left(\sum_jc_j\right)\\
&=&R\left(\sum_{j_1\ldots j_k}c_{j_1}\otimes\ldots\otimes c_{j_k}\right)\\
&=&\frac{1}{2}\sum_{j_1\ldots j_k} c_{j_1}c^*_{j_2}\otimes
c_{j_2}c^*_{j_3}\otimes\ldots\otimes c_{j_k}c_{j_1}^*\\
&=&\frac{4}{2}\sum_{i_1\ldots i_{k-1}}c_{i_1}\otimes c_{i_1}^*c_{i_2}\otimes \ldots \otimes c_{i_{k-2}}^*c_{i_{k-1}}\otimes c_{i_{k-1}}^*\\
&=&2\sum_{i_1\ldots i_{k-1}}(c_{i_1}\otimes c_{i_1}^*)_{12}(c_{i_2}\otimes c_{i_2}^*)_{23}\ldots (c_{i_{k-1}}\otimes c_{i_{k-1}}^*)_{k-1,k}\\
&=&2\left(\sum_{i_1}c_{i_1}\otimes c_{i_1}^*\right)_{12}\left(\sum_{i_2}c_{i_2}\otimes c_{i_2}^*\right)_{23}\ldots \left(\sum_{i_{k-1}}c_{i_{k-1}}\otimes c_{i_{k-1}}^*\right)_{k-1,k}\cr
&=&2f_{12}f_{23}\ldots f_{k-1,k}
\end{eqnarray*}

Now $f$ being invariant under $SU_2$, the same happens for each $f_{i,i+1}$. Since the invariants form an algebra, the above product is invariant under $SU_2$ too.

\medskip

(3) Our next claim is that $SU_2$ leaves in fact invariant any $R(c_\nu)$. But this is something that can be proved by suitably modifying (2), and we refer to \cite{bco} for this.

\medskip

(4) Getting now to $R^*ER(c_\nu)=c_\nu$, we first have the following computation:
\begin{eqnarray*}
<R^*ER(c_\nu),c_j>
&=&<R^*R(c_\nu),c_j>\\
&=&<R(c_\nu),R(c_j)>\\
&=&\sum_i\delta_\nu(i)<R(c_i),R(c_j)>
\end{eqnarray*}

On the other hand, we have from definitions the following formula:
$$<c_\nu,c_j>=\delta_\nu(j)$$

Since the elements $c_j$ span the ambient space, what is left to prove is:
$$\sum_i\delta_\nu(i)<R(c_i),R(c_j)>=\delta_\nu(j)$$

But this can be checked by direct computation, as explained in \cite{bco}.
\end{proof}

Still following \cite{bco}, we can now prove our main result, as follows:

\begin{theorem}
The Pauli representation of $C(S_4^+)$, namely
$$\pi:C(S_4^+)\to C(SU_2, M_4({\mathbb C}))\quad,\quad u_{ij}\to[x\to Proj(c_ixc_j)]$$
is stationary, and so faithful.
\end{theorem}

\begin{proof}
We denote as usual by $c_1,\ldots,c_4$ the Pauli matrices, and we let $e_1,\ldots,e_4$ be the standard basis of $\mathbb C^4$. For a multi-index $i=(i_1,\ldots ,i_k)$, we set:
$$e_i=e_{i_1}\otimes\ldots\otimes e_{i_k}\quad,\quad 
c_i=c_{i_1}\otimes\ldots\otimes c_{i_k}$$

Each partition $\nu\in NC(k)$ creates then two tensors, in the following way:
$$e_\nu=\sum_i\delta_\nu(i)e_i\quad,\quad 
c_\nu=\sum_i\delta_\nu(i)c_i$$

Consider the following $4^k\times 4^k$ matrices, with entries labeled by multi-indices $i,j$:
$$P_{ij}=tr\left(\int \pi_{i_1j_1}\ldots \pi_{i_kj_k}\right)\quad,\quad 
U_{ij}=\left(\int u_{i_1j_1}\ldots u_{i_kj_k}\right)$$

According to Proposition 14.18, we have the following formula:
$$P_{ij}=<c_i,R^*ER(c_j)>$$

Now let $\Phi$ be the linear map $(\mathbb{C}^4)^{\otimes k}\to M_2(\mathbb{C})^{\otimes k}$
given by $\Phi (e_i)=c_i$. In terms of this map, the above equation can be rephrased as follows:
$$P_{ij}=<e_i,\Phi^*R^*ER\Phi(e_j)>$$

According to Proposition 14.17, it is enough to prove that we have $P=U$. But $U$ is the orthogonal projection of $(\mathbb{C}^4)^{\otimes k}$ onto the following space:
$$S_e=span\left\{e_\nu\Big|\nu\in NC(k)\right\}$$

Thus what we have to prove is that $\Phi^*R^*ER\Phi$ is the projection onto $S_e$. But this is equivalent to proving that $R^*ER$ is the projection onto the following space: 
$$S_c=span\left\{ c_\nu\Big|\nu\in NC(k)\right\}$$

We know from Proposition 14.19 that $R^*ER(c_\nu)=c_\nu$, for all $\nu\in NC(k)$. But this implies that $R^*ER$ restricted to $S_c$ is the identity. Now the vector space having as dimension the Catalan number $C_k$, which is the rank of the operator $E$, we conclude that  $S_c$ has to be the image of $R^*ER$. But since $R^*ER$ is self-adjoint, and is the identity on its image, it is the orthogonal projection onto $S_c$, and this finishes the proof.
\end{proof}

As a continuation of this, as explained in chapter 12, we have in fact an isomorphism $S_4^+\simeq SO_3^{-1}$, and in view of this, Theorem 14.20 suggests the following conjecture:

\begin{conjecture}
Given a quantum permutation group of $4$ points, 
$$G\subset S_4^+\simeq SO_3^{-1}$$
coming by twisting a usual ADE subgroup of the group $SO_3$,
$$H\subset SO_3$$
the restriction of the Pauli model for $C(S_4^+)$, with fibers coming from the elements of $H\subset SO_3$, has the algebra $C(G)$ as Hopf image.
\end{conjecture}

To be more precise, Theorem 14.20 tells us that the conjecture holds for $G=S_4^+$ itself. Indeed, here we have $H=SO_3$, so the corresponding restriction of the Pauli model for $C(S_4^+)$ is the Pauli model itself, and this model being stationary, its Hopf image is the algebra $C(S_4^+)$ itself, as stated. The problem is to generalize this. See \cite{bco}.

\section*{14d. Weyl matrices}

Following the more recent paper \cite{bn2}, we would like to discuss now some generalizations of the above, in relation with the Weyl matrices. We first have:

\index{Weyl matrix}
\index{Pauli matrix}

\begin{definition}
Given a finite abelian group $H$, the associated Weyl matrices are
$$W_{ia}:e_b\to<i,b>e_{a+b}$$
where $i\in H$, $a,b\in\widehat{H}$, and where $(i,b)\to<i,b>$ is the Fourier coupling $H\times\widehat{H}\to\mathbb T$.
\end{definition}

As a basic example, consider the cyclic group $H=\mathbb Z_2=\{0,1\}$. Here the Fourier coupling is given by $<i,b>=(-1)^{ib}$, and so the Weyl matrices act via:
$$W_{00}:e_b\to e_b\quad,\quad 
W_{10}:e_b\to(-1)^be_b$$
$$W_{11}:e_b\to(-1)^be_{b+1}\quad,\quad 
W_{01}:e_b\to e_{b+1}$$

Thus, we have the following formulae:
$$W_{00}=\begin{pmatrix}1&0\\0&1\end{pmatrix}\quad,\quad 
W_{10}=\begin{pmatrix}1&0\\0&-1\end{pmatrix}$$
$$W_{11}=\begin{pmatrix}0&-1\\1&0\end{pmatrix}\quad,\quad
W_{01}=\begin{pmatrix}0&1\\1&0\end{pmatrix}$$

We recognize here, up to some multiplicative factors, the four Pauli matrices. Now back to the general case, we have the following well-known result:

\begin{proposition}
The Weyl matrices are unitaries, and satisfy:
\begin{enumerate}
\item $W_{ia}^*=<i,a>W_{-i,-a}$.

\item $W_{ia}W_{jb}=<i,b>W_{i+j,a+b}$.

\item $W_{ia}W_{jb}^*=<j-i,b>W_{i-j,a-b}$.

\item $W_{ia}^*W_{jb}=<i,a-b>W_{j-i,b-a}$.
\end{enumerate}
\end{proposition}

\begin{proof}
The unitarity follows from (3,4), and the rest of the proof goes as follows:

\medskip

(1) We have indeed the following computation:
\begin{eqnarray*}
W_{ia}^*
&=&\left(\sum_b<i,b>E_{a+b,b}\right)^*\\
&=&\sum_b<-i,b>E_{b,a+b}\\
&=&\sum_b<-i,b-a>E_{b-a,b}\\
&=&<i,a>W_{-i,-a}
\end{eqnarray*}

(2) Here the verification goes as follows:
\begin{eqnarray*}
W_{ia}W_{jb}
&=&\left(\sum_d<i,b+d>E_{a+b+d,b+d}\right)\left(\sum_d<j,d>E_{b+d,d}\right)\\
&=&\sum_d<i,b><i+j,d>E_{a+b+d,d}\\
&=&<i,b>W_{i+j,a+b}
\end{eqnarray*}

(3,4) By combining the above two formulae, we obtain:
\begin{eqnarray*}
W_{ia}W_{jb}^*
&=&<j,b>W_{ia}W_{-j,-b}\\
&=&<j,b><i,-b>W_{i-j,a-b}
\end{eqnarray*}

We obtain as well the following formula:
\begin{eqnarray*}
W_{ia}^*W_{jb}
&=&<i,a>W_{-i,-a}W_{jb}\\
&=&<i,a><-i,b>W_{j-i,b-a}
\end{eqnarray*}

But this gives the formulae in the statement, and we are done.
\end{proof}

With $n=|H|$, we can use an isomorphism $l^2(\widehat{H})\simeq\mathbb C^n$ as to view each $W_{ia}$ as a usual matrix, $W_{ia}\in M_n(\mathbb C)$, and hence as a usual unitary, $W_{ia}\in U_n$. Also, given a vector $\xi$, we denote by $Proj(\xi)$ the orthogonal projection onto $\mathbb C\xi$. Following \cite{bn2}, we have:

\begin{proposition}
Given a closed subgroup $E\subset U_n$, we have a representation
$$\pi_H:C(S_N^+)\to M_N(C(E))$$
$$w_{ia,jb}\to[U\to Proj(W_{ia}UW_{jb}^*)]$$
where $n=|H|,N=n^2$, and where $W_{ia}$ are the Weyl matrices associated to $H$.
\end{proposition}

\begin{proof}
The Weyl matrices being given by $W_{ia}:e_b\to<i,b>e_{a+b}$, we have:
$$tr(W_{ia})=\begin{cases}
1&{\rm if}\ (i,a)=(0,0)\\
0&{\rm if}\ (i,a)\neq(0,0)
\end{cases}$$

Together with the formulae in Proposition 14.23, this shows that the Weyl matrices are pairwise orthogonal with respect to the following scalar product on $M_n(\mathbb C)$:
$$<x,y>=tr(xy^*)$$

Thus, these matrices form an orthogonal basis of $M_n(\mathbb C)$, consisting of unitaries:
$$W=\left\{W_{ia}\Big|i\in H,a\in\widehat{H}\right\}$$

Thus, each row and each column of the matrix $\xi_{ia,jb}=W_{ia}UW_{jb}^*$ is an orthogonal basis of $M_n(\mathbb C)$, and so the corresponding projections form a magic unitary, as claimed.
\end{proof}

We will need the following well-known result:

\begin{proposition}
With $T=Proj(x_1)\ldots Proj(x_p)$ and $||x_i||=1$ we have 
$$<T\xi,\eta>=<\xi,x_p><x_p,x_{p-1}>\ldots<x_2,x_1><x_1,\eta>$$
for any $\xi,\eta$. In particular, $Tr(T)=<x_1,x_p><x_p,x_{p-1}>\ldots<x_2,x_1>$.
\end{proposition}

\begin{proof}
For $||x||=1$ we have $Proj(x)\xi=<\xi,x>x$. This gives:
\begin{eqnarray*}
T\xi
&=&Proj(x_1)\ldots Proj(x_p)\xi\\
&=&Proj(x_1)\ldots Proj(x_{p-1})<\xi,x_p>x_p\\
&=&Proj(x_1)\ldots Proj(x_{p-2})<\xi,x_p><x_p,x_{p-1}>x_{p-1}\\
&=&\ldots\\
&=&<\xi,x_p><x_p,x_{p-1}>\ldots<x_2,x_1>x_1
\end{eqnarray*}

Now by taking the scalar product with $\eta$, this gives the first assertion. As for the second assertion, this follows from the first assertion, by summing over $\xi=\eta=e_i$.
\end{proof}

Now back to the Weyl matrix models, let us first compute $T_p$. We have:

\begin{proposition}
We have the formula
\begin{eqnarray*}
(T_p)_{ia,jb}
&=&\frac{1}{N}<i_1,a_1-a_p>\ldots<i_p,a_p-a_{p-1}><j_1,b_1-b_2>\ldots<j_p,b_p-b_1>\\
&&\int_Etr(W_{i_1-i_2,a_1-a_2}UW_{j_2-j_1,b_2-b_1}U^*)\ldots tr(W_{i_p-i_1,a_p-a_1}UW_{j_1-j_p,b_1-b_p}U^*)dU
\end{eqnarray*}
with all the indices varying in a cyclic way.
\end{proposition}

\begin{proof}
By using the trace formula in Proposition 14.25, we obtain:
\begin{eqnarray*}
&&(T_p)_{ia,jb}\\
&=&\left(tr\otimes\int_E\right)\left(Proj(W_{i_1a_1}UW_{j_1b_1}^*)\ldots Proj(W_{i_pa_p}UW_{j_pb_p}^*)\right)\\
&=&\frac{1}{N}\int_E<W_{i_1a_1}UW_{j_1b_1}^*,W_{i_pa_p}UW_{j_pb_p}^*>\ldots<W_{i_2a_2}UW_{j_2b_2}^*,W_{i_1a_1}UW_{j_1b_1}^*>dU
\end{eqnarray*}

In order to compute now the scalar products, observe that we have:
\begin{eqnarray*}
<W_{ia}UW_{jb}^*,W_{kc}UW_{ld}^*>
&=&tr(W_{jb}U^*W_{ia}^*W_{kc}UW_{ld}^*)\\
&=&tr(W_{ia}^*W_{kc}UW_{ld}^*W_{jb}U^*)\\
&=&<i,a-c><l,d-b>tr(W_{k-i,c-a}UW_{j-l,b-d}U^*)
\end{eqnarray*}

By plugging these quantities into the formula of $T_p$, we obtain the result.
\end{proof}

Consider now the Weyl group $W=\{W_{ia}\}\subset U_n$, that we already met in the proof of Proposition 14.24. We have the following result, from \cite{bn2}:

\index{Weyl matrix model}

\begin{theorem}
For any compact group $W\subset E\subset U_n$, the model
$$\pi_H:C(S_N^+)\to M_N(C(E))$$
$$w_{ia,jb}\to[U\to Proj(W_{ia}UW_{jb}^*)]$$
constructed above is stationary on its image.
\end{theorem}

\begin{proof}
We must prove that we have $T_p^2=T_p$. We have:
\begin{eqnarray*}
&&(T_p^2)_{ia,jb}\\
&=&\sum_{kc}(T_p)_{ia,kc}(T_p)_{kc,jb}\\
&=&\frac{1}{N^2}\sum_{kc}<i_1,a_1-a_p>\ldots<i_p,a_p-a_{p-1}><k_1,c_1-c_2>\ldots<k_p,c_p-c_1>\\
&&<k_1,c_1-c_p>\ldots<k_p,c_p-c_{p-1}><j_1,b_1-b_2>\ldots<j_p,b_p-b_1>\\
&&\int_Etr(W_{i_1-i_2,a_1-a_2}UW_{k_2-k_1,c_2-c_1}U^*)\ldots tr(W_{i_p-i_1,a_p-a_1}UW_{k_1-k_p,c_1-c_p}U^*)dU\\
&&\int_Etr(W_{k_1-k_2,c_1-c_2}VW_{j_2-j_1,b_2-b_1}V^*)\ldots tr(W_{k_p-k_1,c_p-c_1}VW_{j_1-j_p,b_1-b_p}V^*)dV
\end{eqnarray*}

By rearranging the terms, this formula becomes:
\begin{eqnarray*}
&&(T_p^2)_{ia,jb}\\
&=&\frac{1}{N^2}<i_1,a_1-a_p>\ldots<i_p,a_p-a_{p-1}><j_1,b_1-b_2>\ldots<j_p,b_p-b_1>\\
&&\int_E\int_E\sum_{kc}<k_1-k_p,c_1-c_p>\ldots<k_p-k_{p-1},c_p-c_{p-1}>\\
&&tr(W_{i_1-i_2,a_1-a_2}UW_{k_2-k_1,c_2-c_1}U^*)tr(W_{k_1-k_2,c_1-c_2}VW_{j_2-j_1,b_2-b_1}V^*)\\
&&\hskip50mm\ldots\ldots\\
&&tr(W_{i_p-i_1,a_p-a_1}UW_{k_1-k_p,c_1-c_p}U^*)tr(W_{k_p-k_1,c_p-c_1}VW_{j_1-j_p,b_1-b_p}V^*)dUdV
\end{eqnarray*}

Let us denote by $I$ the above double integral. By using $W_{kc}^*=<k,c>W_{-k,-c}$ for each of the couplings, and by moving as well all the $U^*$ variables to the left, we obtain:
\begin{eqnarray*}
I
&=&\int_E\int_E\sum_{kc}tr(U^*W_{i_1-i_2,a_1-a_2}UW_{k_2-k_1,c_2-c_1})tr(W_{k_2-k_1,c_2-c_1}^*VW_{j_2-j_1,b_2-b_1}V^*)\\
&&\hskip50mm\ldots\ldots\\
&&tr(U^*W_{i_p-i_1,a_p-a_1}UW_{k_1-k_p,c_1-c_p})tr(W_{k_1-k_p,c_1-c_p}^*VW_{j_1-j_p,b_1-b_p}V^*)dUdV
\end{eqnarray*}

In order to perform now the sums, we use the following formula:
\begin{eqnarray*}
tr(AW_{kc})tr(W_{kc}^*B)
&=&\frac{1}{N}\sum_{qrst}A_{qr}(W_{kc})_{rq}(W^*_{kc})_{st}B_{ts}\\
&=&\frac{1}{N}\sum_{qrst}A_{qr}<k,q>\delta_{r-q,c}<k,-s>\delta_{t-s,c}B_{ts}\\
&=&\frac{1}{N}\sum_{qs}<k,q-s>A_{q,q+c}B_{s+c,s}
\end{eqnarray*}

If we denote by $A_x,B_x$ the variables which appear in the formula of $I$, we have:
\begin{eqnarray*}
&&I\\
&=&\frac{1}{N^p}\int_E\int_E\sum_{kcqs}<k_2-k_1,q_1-s_1>\ldots<k_1-k_p,q_p-s_p>\\
&&(A_1)_{q_1,q_1+c_2-c_1}(B_1)_{s_1+c_2-c_1,s_1}\ldots (A_p)_{q_p,q_p+c_1-c_p}(B_p)_{s_p+c_1-c_p,s_p}\\
&=&\frac{1}{N^p}\int_E\int_E\sum_{kcqs}<k_1,q_p-s_p-q_1+s_1>\ldots<k_p,q_{p-1}-s_{p-1}-q_p+s_p>\\
&&(A_1)_{q_1,q_1+c_2-c_1}(B_1)_{s_1+c_2-c_1,s_1}\ldots (A_p)_{q_p,q_p+c_1-c_p}(B_p)_{s_p+c_1-c_p,s_p}
\end{eqnarray*}

Now observe that we can perform the sums over $k_1,\ldots,k_p$. We obtain in this way a multiplicative factor $n^p$, along with the condition:
$$q_1-s_1=\ldots=q_p-s_p$$

Thus we must have $q_x=s_x+a$ for a certain $a$, and the above formula becomes:
$$I=\frac{1}{n^p}\int_E\int_E\sum_{csa}(A_1)_{s_1+a,s_1+c_2-c_1+a}(B_1)_{s_1+c_2-c_1,s_1}\ldots(A_p)_{s_p+a,s_p+c_1-c_p+a}(B_p)_{s_p+c_1-c_p,s_p}$$

Consider now the variables $r_x=c_{x+1}-c_x$, which altogether range over the set $Z$ of multi-indices having sum 0. By replacing the sum over $c_x$ with the sum over $r_x$, which creates a multiplicative $n$ factor, we obtain the following formula:
$$I=\frac{1}{n^{p-1}}\int_E\int_E\sum_{r\in Z}\sum_{sa}(A_1)_{s_1+a,s_1+r_1+a}(B_1)_{s_1+r_1,s_1}\ldots(A_p)_{s_p+a,s_p+r_p+a}(B_p)_{s_p+r_p,s_p}$$

For an arbitrary multi-index $r$, we have the following formula:
$$\delta_{\sum_ir_i,0}=\frac{1}{n}\sum_i<i,r_1>\ldots<i,r_p>$$

Thus, we can replace the sum over $r\in Z$ by a full sum, as follows:
\begin{eqnarray*}
I
&=&\frac{1}{n^p}\int_E\int_E\sum_{rsia}<i,r_1>(A_1)_{s_1+a,s_1+r_1+a}(B_1)_{s_1+r_1,s_1}\\
&&\hskip40mm\ldots\ldots\\
&&\hskip20mm<i,r_p>(A_p)_{s_p+a,s_p+r_p+a}(B_p)_{s_p+r_p,s_p}
\end{eqnarray*}

In order to ``absorb'' now the indices $i,a$, we can use the following formula:
\begin{eqnarray*}
&&W_{ia}^*AW_{ia}\\
&=&\left(\sum_b<i,-b>E_{b,a+b}\right)\left(\sum_{bc}E_{a+b,a+c}A_{a+b,a+c}\right)\left(\sum_c<i,c>E_{a+c,c}\right)\\
&=&\sum_{bc}<i,c-b>E_{bc}A_{a+b,a+c}
\end{eqnarray*}

Thus we have the following formula:
$$(W_{ia}^*AW_{ia})_{bc}=<i,c-b>A_{a+b,a+c}$$

With this in hand, our formula becomes:
\begin{eqnarray*}
&&I\\
&=&\frac{1}{n^p}\int_E\int_E\sum_{rsia}(W_{ia}^*A_1W_{ia})_{s_1,s_1+r_1}(B_1)_{s_1+r_1,s_1}\ldots(W_{ia}^*A_pW_{ia})_{s_p,s_p+r_p}(B_p)_{s_p+r_p,s_p}\\
&=&\int_E\int_E\sum_{ia}tr(W_{ia}^*A_1W_{ia}B_1)\ldots\ldots tr(W_{ia}^*A_pW_{ia}B_p)
\end{eqnarray*}

Now by replacing $A_x,B_x$ with their respective values, we obtain:
\begin{eqnarray*}
I
&=&\int_E\int_E\sum_{ia}tr(W_{ia}^*U^*W_{i_1-i_2,a_1-a_2}UW_{ia}VW_{j_2-j_1,b_2-b_1}V^*)\\
&&\hskip30mm\ldots\ldots\\
&&tr(W_{ia}^*U^*W_{i_p-i_1,a_p-a_1}UW_{ia}VW_{j_1-j_p,b_1-b_p}V^*)dUdV
\end{eqnarray*}

By moving the $W_{ia}^*U^*$ variables at right, we obtain, with $S_{ia}=UW_{ia}V$:
\begin{eqnarray*}
I
&=&\sum_{ia}\int_E\int_Etr(W_{i_1-i_2,a_1-a_2}S_{ia}W_{j_2-j_1,b_2-b_1}S_{ia}^*)\\
&&\hskip30mm\ldots\ldots\\
&&tr(W_{i_p-i_1,a_p-a_1}S_{ia}W_{j_1-j_p,b_1-b_p}S_{ia}^*)dUdV
\end{eqnarray*}

Now since $S_{ia}$ is Haar distributed when $U,V$ are Haar distributed, we obtain:
$$I=N\int_E\int_Etr(W_{i_1-i_2,a_1-a_2}UW_{j_2-j_1,b_2-b_1}U^*)\ldots tr(W_{i_p-i_1,a_p-a_1}UW_{j_1-j_p,b_1-b_p}U^*)dU$$

But this is exactly $N$ times the integral in the formula of $(T_p)_{ia,jb}$, from Proposition 14.26. Since the $N$ factor cancels with one of the two $N$ factors that we found in the beginning of the proof, when first computing $(T_p^2)_{ia,jb}$, we are done.
\end{proof}

As an illustration for the above result, we can recover Theorem 14.20:

\index{Pauli model}

\begin{theorem}
We have a stationary matrix model as follows,
$$\pi:C(S_4^+)\subset M_4(C(SU_2))$$
given on the standard coordinates by the following formula,
$$\pi(u_{ij})=[x\to Proj(c_ixc_j)]$$
where $x\in SU_2$, and $c_1,c_2,c_3,c_4$ are the Pauli matrices.
\end{theorem}

\begin{proof}
As already explained in the comments following Definition 14.22, the Pauli matrices appear as particular cases of the Weyl matrices:
$$W_{00}=\begin{pmatrix}1&0\\0&1\end{pmatrix}\quad,\quad 
W_{10}=\begin{pmatrix}1&0\\0&-1\end{pmatrix}$$
$$W_{11}=\begin{pmatrix}0&-1\\1&0\end{pmatrix}\quad,\quad
W_{01}=\begin{pmatrix}0&1\\1&0\end{pmatrix}$$

Thus, Theorem 14.27 produces in this case the model in the statement.
\end{proof}

Observe that, since the projection $Proj(c_ixc_j)$ depends only on the image of $x$ in the quotient $SU_2\to SO_3$, we can replace the model space $SU_2$ by the smaller space $SO_3$.  This can be used in conjunction with the isomorphism $S_4^+\simeq SO_3^{-1}$, and as explained in \cite{bbc}, our model becomes in this way something more conceptual, as follows:
$$\pi:C(SO_3^{-1})\subset M_4(C(SO_3))$$

As a philosophical conclusion, to this and to some previous findings as well, no matter what we do, we always end up getting back to $SU_2,SO_3$. Thus, we are probably doing some physics here. This is indeed the case, the above computations being closely related to the standard computations for the Ising and Potts models. The general relation, however, between quantum permutations and lattice models, is not axiomatixed yet.

\section*{14e. Exercises}

This was a quite technical chapter, and as exercises on this, we have:

\begin{exercise}
Get to know more about stationarity, and the Thoma theorem.
\end{exercise}

\begin{exercise}
Learn about modeling questions for arbitrary free manifolds.
\end{exercise}

\begin{exercise}
Learn, from Franz, Skalski and others, about idempotent states.
\end{exercise}

\begin{exercise}
Learn more, from Bichon and Dubois-Violette, about half-liberation.
\end{exercise}

\begin{exercise}
Get to know about the Pauli matrices, and their use in physics.
\end{exercise}

\begin{exercise}
Meditate at the modeling of the ADE subgroups $G\subset S_4^+$.
\end{exercise}

\begin{exercise}
Learn about the Weyl matrices, and their various uses.
\end{exercise}

\begin{exercise}
Learn about cocyclic Weyl models, generalizing the above.
\end{exercise}

As bonus exercise, advanced, convert what we said above into statistical mechanics.

\chapter{Subfactor theory}

\section*{15a. Factors, subfactors}

We have been doing some quite advanced mathematics lately, but there is certainly more advanced than that, too. Following Jones \cite{jo1}, \cite{jo2}, \cite{jo3}, \cite{jo4}, \cite{jo5}, we will be talking here about subfactors, with our official motivations being as follows:

\begin{principle}
When looking at an inclusion of ${\rm II}_1$ factors $A_0\subset A_1$, also called subfactor, some interesting combinatorics and probability appears, as follows:
\begin{enumerate}
\item The symmetry group of such an inclusion is a quantum group type object $G$, more general than the usual quantum groups, and this due to $|G|\notin\mathbb N$, in general.

\item We can do some interesting probability over such generalized quantum groups $G$, notably by studying various arithmetic versions of $\gamma_t$ and $\pi_t$.

\item And with all this being relevant to modern theoretical physics, and more specifically, to advanced statistical mechanics, and conformal field theory.
\end{enumerate}
\end{principle}

So this will be our guiding principle, and in what follows we will explain (1), to some reasonable extent, with theory and examples, and then do some computations for (2), which will be something quite fun, completing our free probability knowledge. 

\bigskip

As for (3), things here are more complicated. Jones was originally motivated by basic statistical mechanics, namely Ising and Potts models, and by algebraic quantum field theory. Later, the key paper of Witten \cite{wit}, providing a direct geometric explanation for his knot polynomial, using more advanced quantum field theory, joined the picture too. And afterwards, many things happened, and from a modern perspective, subfactors and the related notion of planar algebras can be somehow considered as a ``theory of everything'', or rather as the tools of that theory, which still remains to be found.

\bigskip

Well, not quite sure that all this is understandable, but somehow, this is the idea. For more, purchase a copy of Guionnet et al. \cite{agz}, a copy of Di Francesco et al. \cite{dms} too, these two books are good for training for such things, and in addition they are both quite big, coming in nice gray-orange livery, making them perfect as construction bricks too, or for smashing intruders, read, meditate, levitate, and if you get the point, congratulations.

\bigskip

Finally, still talking motivations, as an addendum to Principle 15.1, putting us in tune with what we officially wanted to do in this present Part IV, matrix models, we have:

\begin{addendum}
Subfactors $A_0\subset A_1$ are something quite concrete, and even of matrix type, when looking at the hyperfinite ones, $A_0\simeq A_1\simeq R$, and so can be thought of as being ``matrix models'' for the various laws emerging from their study.
\end{addendum}

But probably too much talking, let us get to work. We already know about von Neumann algebras from chapter 5, and as explained there, the study of such algebras naturally leads to the study of the ${\rm II}_1$ factors, which are defined as follows:

\begin{definition}
A ${\rm II}_1$ factor is a von Neumann algebra $A\subset B(H)$ which is infinite dimensional, has trivial center, and has a trace $tr:A\to\mathbb C$.
\end{definition}

There are many things that can be said about the ${\rm II}_1$ factors, and of particular interest is the following key result of Murray and von Neumann \cite{mvo}, which clarifies the situation with the various Hilbert space representations $A\subset B(H)$ of a given ${\rm II}_1$ factor $A$:

\index{coupling constant}

\begin{theorem}
Given a representation of a ${\rm II}_1$ factor $A\subset B(H)$, we can talk about the corresponding coupling constant
$$\dim_AH\in(0,\infty]$$
which for the standard form, where $H=L^2(A)$, takes the value $1$, and which in general mesures how far is $A\subset B(H)$ from the standard form.
\end{theorem}

\begin{proof}
There are several proofs for this fact, the idea being as follows:

\medskip

(1) We can amplify the standard representation of $A$, on the Hilbert space $L^2(A)$, into a representation on $L^2(A)\otimes l^2(\mathbb N)$, and then cut it down with a projection. We obtain in this way a whole family of embeddings $A\subset B(H)$, which are quite explicit.

\medskip

(2) The point now is that of proving, via a technical $2\times2$ matrix trick, that any representation $A\subset B(H)$ appears in this way. In this picture, the coupling constant appears as the trace of the projection used to cut down $L^2(A)\otimes l^2(\mathbb N)$.

\medskip

(3) Thus, we are led to the above conclusion. Alternatively, the coupling constant can be defined as follows, with the number on the right being independent of the choice on a nonzero vector $x\in H$, and with this being the original definition from \cite{mvo}:
$$\dim_AH=\frac{tr_A(P_{A'x})}{tr_{A'}(P_{Ax})}$$

We refer to \cite{mvo}, or to any operator algebra book, for more details here.
\end{proof}

Following Jones \cite{jo1}, given a ${\rm II}_1$ factor $A_0$, let us discuss now the representations $A_0\subset A_1$, with $A_1$ being another ${\rm II}_1$ factor. This is a quite natural notion too, and perhaps even more natural than the representations $A_0\subset B(H)$, because we have decided in chapter 5 that the ${\rm II}_1$ factors $A_1$, and not the full operator algebras $B(H)$, are the correct infinite dimensional generalization of the usual matrix algebras $M_N(\mathbb C)$.

\bigskip

Given an inclusion of ${\rm II}_1$ factors $A_0\subset A_1$, a first question is that of defining its index, measuring how big is $A_1$, when compared to $A_0$. This can be done as follows:

\index{subfactor}
\index{index of subfactor}

\begin{theorem}
Given an inclusion of ${\rm II}_1$ factors $A_0\subset A_1$, the number
$$N=\frac{\dim_{A_0}H}{\dim_{A_1}H}$$
is independent of the ambient Hilbert space $H$, and is called index.
\end{theorem}

\begin{proof}
This is something very standard, with the fact that the index as defined by the above formula is indeed independent of the ambient Hilbert space $H$ coming from the various basic properties of the coupling constant, from Theorem 15.4 and its proof.
\end{proof}

In order to study the subfactors, let us start with the following standard result:

\index{conditional expectation}

\begin{proposition}
Given $A_0\subset A_1$, there is a unique linear map $E:A_1\to A_0$ which is positive, unital, trace-preserving and such that
$$E(a_1ba_2)=a_1E(b)a_2$$
for any $a_1,a_2\in A_0$. This map is called conditional expectation from $A_1$ onto $A_0$.
\end{proposition}

\begin{proof}
Consider the standard representation by left multiplication $A_1\subset L^2(A_1)$. If we denote by $\Omega$ the standard cyclic and separating vector of $L^2(A_1)$, we have an identification of vector spaces $A_0\Omega=L^2(A_0)$. Consider now the following projection:
$$e:L^2(A_1)\to L^2(A_0)$$

It follows from definitions that we have an inclusion $e(A_1\Omega)\subset A_0\Omega$, so our projection $e$ induces by restriction a certain linear map $E:A_1\to A_0$, related to $e$ by:
$$exe=E(x)e$$

But this shows that the linear map $E$ satisfies the various conditions in the statement, namely positivity, unitality, trace preservation and bimodule property. As for the uniqueness assertion, this follows by using the same argument, applied backwards.
\end{proof}

We will be interested in what follows in the orthogonal projection $e:L^2(A_1)\to L^2(A_0)$ producing the expectation $E:A_1\to A_0$, rather than in $E$ itself:

\index{Jones projection}

\begin{definition}
Associated to any subfactor $A_0\subset A_1$ is the orthogonal projection
$$e:L^2(A_1)\to L^2(A_0)$$
producing the conditional expectation $E:A_1\to A_0$ via the following formula:
$$exe=E(x)e$$
This projection is called Jones projection for the subfactor $A_0\subset A_1$.
\end{definition}

Quite remarkably, the subfactor $A_0\subset A_1$, as well as its commutant, can be recovered from the knowledge of this projection, in the following way:

\begin{proposition}
Given a subfactor $A_0\subset A_1$, with Jones projection $e$, we have
$$A_0=A_1\cap\{e\}'\quad,\quad 
A_0'=(A_1'\cap\{e\})''$$
as equalities of von Neumann algebras, acting on the space $L^2(A_1)$.
\end{proposition}

\begin{proof}
These formulae basically follow from $exe=E(x)e$, as follows:

\medskip

(1) Let us first prove that we have $A_0\subset A_1\cap\{e\}'$. Given $x\in A_0$, we have:
$$xe=E(x)e=exe\quad,\quad 
x^*e=E(x^*)e=ex^*e$$

Thus, we obtain, as desired, $ex
=(x^*e)^*
=(ex^*e)^*
=exe
=xe$.

\medskip

(2) Let us prove now that $A_1\cap\{e\}'\subset A_0$. Assuming $ex=xe$, we have:
$$E(x)e=exe=xe^2=xe$$

We conclude from this that we have the following equality:
$$(E(x)-x)\Omega=(E(x)-x)e\Omega=0$$

Now since $\Omega$ is separating for $A_1$ we have, as desired, $x=E(x)\in A_0$.

\medskip

(3) In order to prove now $A_0'=<A_1',e>$, observe that we have:
$$A_0
=A_1\cap\{e\}'
=A_1''\cap\{e\}'
=(A_1'\cap\{e\})'$$

Now by taking the commutant, we obtain $A_0'=(A_1'\cap\{e\})''$, as desired.
\end{proof}

We are now ready to formulate a key definition, as follows:

\index{basic construction}

\begin{definition}
Associated to any subfactor $A_0\subset A_1$ is the basic construction
$$A_0\subset_eA_1\subset A_2$$
with $A_2=<A_1,e>$ being the algebra generated by $A_1$ and by $e:L^2(A_1)\to L^2(A_0)$.
\end{definition}

The idea now, following as before Jones \cite{jo1}, will be that the inclusion $A_1\subset A_2$ appears as a kind of ``reflection'' of the original inclusion $A_0\subset A_1$, and also that the basic construction can be iterated, with this leading to non-trivial results. We first have:

\begin{proposition}
Given a subfactor $A_0\subset A_1$ having finite index, 
$$[A_1:A_0]<\infty$$
the basic construction $A_0\subset_eA_1\subset A_2$ has the following properties:
\begin{enumerate}
\item $A_2=JA_0'J$, with $J(T)=T^*$.

\item $A_2=\overline{A_1+A_1eb}$.

\item $A_2$ is a ${\rm II}_1$ factor.

\item $[A_2:A_1]=[A_1:A_0]$.

\item $eA_2e=A_0e$.

\item $tr(e)=[A_1:A_0]^{-1}$.

\item $tr(xe)=tr(x)[A_1:A_0]^{-1}$, for any $x\in A_1$.
\end{enumerate}
\end{proposition}

\begin{proof}
All this is very standard, the idea being as follows:

\medskip

(1) We have $JA_1'J=A_1$ and $JeJ=e$, which gives, as desired:
$$JA_0'J
=J<A_1',e>J
=<JA_1'J,JeJ>
=<A_1,e>
=A_2$$

(2) This follows from the fact that the vector space $A_1+A_1eA_1$ is closed under multiplication, and from the fact that we have $exe=E(x)e$.

\medskip

(3) This comes from the fact that $[A_1:A_0]<\infty$ is equivalent to the fact that $A_0'$ is a factor. Indeed, this is in turn equivalent to the fact that $A_2=JA_0'J$ is a factor.

\medskip

(4) We have indeed the following computation:
\begin{eqnarray*}
[A_2:A_1]
&=&\frac{\dim_{A_1}L^2(A_1)}{\dim_{A_2}L^2(A_1)}\\
&=&\frac{1}{\dim_{A_2}L^2(A_1)}\\
&=&\frac{1}{\dim_{JA_0'J}L^2(A_1)}\\
&=&\frac{1}{\dim_{A_0'}L^2(A_1)}\\
&=&\dim_{A_0}L^2(A_1)\\
&=&[A_1:A_0]
\end{eqnarray*}

(5) This follows indeed from (2) and from the formula $exe=E(x)e$.

\medskip

(6) We have indeed the following computation, to start with:
\begin{eqnarray*}
1
&=&\dim_{A_0}L^2(A_0)\\
&=&\dim_{A_0}(eL^2(A_1))\\
&=&tr_{A_0'}(e)\dim_{A_0}(L^2(A_1))\\
&=&tr_{A_0'}(a)[A_1:A_0]
\end{eqnarray*}

Now since $A_2=JA_0'J$ and $JeJ=e$, we obtain from this, as desired:
$$tr(e)
=tr_{JA_0'J}(JeJ)
=tr_{A_0'}(e)
=[A_1:A_0]^{-1}$$

(7) We already know from (6) that the formula in the statement holds for $x=1$. In order to discuss the general case, observe first that for $x,y\in A_0$ we have:
$$tr(xye)=tr(yex)=tr(yxe)$$

Thus the linear map $x\to tr(xe)$ is a trace on $A_0$, and by uniqueness of the trace on $A_0$, we must have $tr(xe)=c\cdot tr(x)$ for a certain constant $c>0$. Now by using (6) we obtain $c=[A_1:A_0]^{-1}$, so we have proved the formula in the statement for $x\in A_0$. The passage to the general case $x\in A_1$ can be done as follows:
$$tr(xe)
=tr(exe)
=tr(E(x)e)
=tr(E(x))c
=tr(x)c$$

Thus, we have proved the formula in the statement, in general.
\end{proof}

Let us perform now twice the basic construction. We have here:

\begin{proposition}
In the context of a double basic construction
$$A_0\subset_eA_1\subset_fA_2\subset A_3$$
we have the formulae $fef=[A_1:A_0]^{-1}f$ and $efe=[A_1:A_0]^{-1}e$.
\end{proposition}

\begin{proof}
The first formula is clear from definitions, because we have:
$$fef
=E(e)f
=tr(e)f
=[A_1:A_0]^{-1}f$$

Regarding the second formula, it is enough to check it on $(A_1+A_1eA_1)\Omega$, which is a dense subset. Thus, we must show that for any $x,y,z\in A_1$, we have:
$$efe(x+yez)\Omega=[A_1:A_0]^{-1}e(x+yez)\Omega$$

For this purpose, we will prove that we have, for any $x,y,z\in A_1$:
$$efex\Omega=[A_1:A_0]^{-1}ex\Omega\quad,\quad 
efeyez\Omega=[A_1:A_0]^{-1}eyez\Omega$$

The first formula can be established as follows:
\begin{eqnarray*}
efex\Omega
&=&efexf\Omega\\
&=&eE(ex)f\Omega\\
&=&eE(e)xf\Omega\\
&=&[A_1:A_0]^{-1}exf\Omega\\
&=&[A_1:A_0]^{-1}ex\Omega
\end{eqnarray*}

As for the second formula, this can be established as follows:
\begin{eqnarray*}
efeyez\Omega
&=&efeyezf\Omega\\
&=&eE(eyez)f\Omega\\
&=&eE(eye)zf\Omega\\
&=&eE(E(y)e)zf\Omega\\
&=&eE(y)E(e)zf\Omega\\
&=&[A_1:A_0]^{-1}eE(y)zf\Omega\\
&=&[A_1:A_0]^{-1}eyezf\Omega\\
&=&[A_1:A_0]^{-1}eyez\Omega
\end{eqnarray*}

Thus, we are led to the conclusion in the statement.
\end{proof}

We can in fact perform the basic construction by recurrence, and we obtain:

\index{Jones tower}
\index{Jones projections}

\begin{theorem}
Associated to any subfactor $A_0\subset A_1$ is the Jones tower
$$A_0\subset_{e_1}A_1\subset_{e_2}A_2\subset_{e_3}A_3\subset\ldots\ldots$$
with the Jones projections having the following properties:
\begin{enumerate}
\item $e_i^2=e_i=e_i^*$.

\item $e_ie_j=e_je_i$ for $|i-j|\geq2$.

\item $e_ie_{i\pm1}e_i=[A_1:A_0]^{-1}e_i$.

\item $tr(we_{n+1})=[A_1:A_0]^{-1}tr(w)$, for any word $w\in<e_1,\ldots,e_n>$.
\end{enumerate}
\end{theorem}

\begin{proof}
This follows indeed from Proposition 15.10 and Proposition 15.11, because the triple basic construction does not need any further study. See \cite{jo1}.
\end{proof}

The relations found in Theorem 15.12 are in fact well-known, from the standard theory of the Temperley-Lieb algebra. In order to discuss this, we will need:

\index{Jones projection}

\begin{proposition}
The Temperley-Lieb algebra $TL_N(k)$ is generated by the diagrams
$$\varepsilon_1={\ }^\cup_\cap\quad,\quad 
\varepsilon_2=|\!{\ }^\cup_\cap\quad,\quad
\varepsilon_3=||\!{\ }^\cup_\cap\quad,\quad 
\ldots$$
which are all multiples of projections, in the sense that their rescaled versions
$$e_i=N^{-1}\varepsilon_i$$
satisfy the abstract projection relations $e_i^2=e_i=e_i^*$.
\end{proposition}

\begin{proof}
We have two assertions here, the idea being as follows:

\medskip

(1) The fact that the Temperley-Lieb algebra $TL_N(k)$ is indeed generated by the sequence $\varepsilon_1,\varepsilon_2,\ldots$ follows by drawing pictures, and more specifically by decomposing each basis element $\pi\in NC_2(k,k)$ as a product of such elements $\varepsilon_i$.

\medskip

(2) Regarding the projection assertion, when composing $\varepsilon_i$ with itself we obtain $\varepsilon_i$ itself, times a circle. Thus $\varepsilon_i^2=N\varepsilon_i$, so  $e_i=N^{-1}\varepsilon_i$ satisfies $e_i^2=e_i^*=e_i$, as desired.
\end{proof}

As a second result now, making the link with Theorem 15.12, we have:

\begin{proposition}
The standard generators $e_i=N^{-1}\varepsilon_i$ of the Temperley-Lieb algebra $TL_N(k)$ have the following properties, where $tr$ is the trace obtained by closing:
\begin{enumerate}
\item $e_ie_j=e_je_i$ for $|i-j|\geq2$.

\item $e_ie_{i\pm1}e_i=N^{-1}e_i$.

\item $tr(we_{n+1})=N^{-1}tr(w)$, for any word $w\in<e_1,\ldots,e_n>$.
\end{enumerate}
\end{proposition}

\begin{proof}
This follows indeed by doing some elementary computations with diagrams, in the spirit of those performed in the proof of Proposition 15.13.
\end{proof}

With the above results in hand, and still following Jones' paper \cite{jo1}, we can now reformulate Theorem 15.12 into something more conceptual, as follows:

\index{subfactor}
\index{Temperley-Lieb algebra}

\begin{theorem}
Given a subfactor $A_0\subset A_1$, construct its the Jones tower:
$$A_0\subset_{e_1}A_1\subset_{e_2}A_2\subset_{e_3}A_3\subset\ldots\ldots$$
The rescaled sequence of projections $e_1,e_2,e_3,\ldots\in B(H)$ produces then a representation 
$$TL_N\subset B(H)$$
of the Temperley-Lieb algebra of index $N=[A_1:A_0]$.
\end{theorem}

\begin{proof}
We know from Theorem 15.12 that the rescaled sequence of Jones projections $e_1,e_2,e_3,\ldots\in B(H)$ behaves algebrically exactly as the following $TL_N$ diagrams:
$$\varepsilon_1={\ }^\cup_\cap\quad,\quad 
\varepsilon_2=|\!{\ }^\cup_\cap\quad,\quad
\varepsilon_3=||\!{\ }^\cup_\cap\quad,\quad 
\ldots$$

But these diagrams generate $TL_N$, and so we have an embedding $TL_N\subset B(H)$, where $H$ is the Hilbert space where our subfactor $A_0\subset A_1$ lives, as claimed.
\end{proof}

Let us make as well the following key observation, also from \cite{jo1}:

\index{relative commutant}
\index{higher commutant}
\index{planar algebra}

\begin{theorem} 
Given a finite index subfactor $A_0\subset A_1$, the graded algebra $P=(P_k)$ formed by the sequence of higher relative commutants
$$P_k=A_0'\cap A_k$$
contains the copy of the Temperley-Lieb algebra constructed above, $TL_N\subset P$. This graded algebra $P=(P_k)$ is called ``planar algebra'' of the subfactor.
\end{theorem}

\begin{proof}
As a first observation, since the Jones projection $e_1:A_1\to A_0$ commutes with $A_0$, we have $e_1\in P_2$. By translation we obtain, for any $k\in\mathbb N$:
$$e_1,\ldots,e_{k-1}\in P_k$$

Thus we have indeed an inclusion of graded algebras $TL_N\subset P$, as claimed.
\end{proof}

As an interesting consequence of the above results, somehow contradicting the ``continuous geometry'' philosophy that has being going on so far, in relation with the ${\rm II}_1$ factors, we have the following surprising result, also from Jones' original paper \cite{jo1}:

\index{index theorem}
\index{quantized index}

\begin{theorem}
The index of subfactors $A_0\subset A_1$ is ``quantized'' in the $[1,4]$ range,
$$N\in\left\{4\cos^2\left(\frac{\pi}{n}\right)\Big|n\geq3\right\}\cup[4,\infty]$$
with the obstruction coming from the existence of the representation $TL_N\subset B(H)$.
\end{theorem}

\begin{proof}
In order to best comment on what happens, when iterating the basic construction, let us record the first few values of the numbers in the statement:
$$4\cos^2\left(\frac{\pi}{3}\right)=1\quad,\quad 
4\cos^2\left(\frac{\pi}{4}\right)=2$$
$$4\cos^2\left(\frac{\pi}{5}\right)=\frac{3+\sqrt{5}}{2}\quad,\quad 
4\cos^2\left(\frac{\pi}{6}\right)=3$$

When performing a basic construction, we get, by trace manipulations on $e_1$:
$$N\notin(1,2)$$

With a double basic construction, we get, by trace manipulations on $<e_1,e_2>$:
$$N\notin\left(2\ ,\,\frac{3+\sqrt{5}}{2}\right)$$

With a triple basic construction, we get, by trace manipulations on $<e_1,e_2,e_3>$:
$$N\notin\left(\frac{3+\sqrt{5}}{2}\ ,\,3\right)$$

And so on. In practice now, the most elegant way of proving the result is by using the fundamental fact, from Theorem 15.15, that that sequence of Jones projections $e_1,e_2,e_3,\ldots\subset B(H)$ generate a copy of the Temperley-Lieb algebra of index $N$:
$$TL_N\subset B(H)$$

Indeed, with this result in hand, we must prove that such a representation cannot exist in index $N<4$, unless we are in the following special situation:
$$N=4\cos^2\left(\frac{\pi}{n}\right)$$

But this can be proved by using some suitable trace and positivity manipulations on $TL_N$, as in the situations discussed above. For full details here, we refer to \cite{jo1}.
\end{proof}

\section*{15b. Planar algebras}

Quite remarkably, the planar algebra structure of $TL_N$, taken in an intuitive sense, that of composing diagrams, in various possible ways, extends to a planar algebra structure on the sequence of higher relative commutants $P_n=A_0'\cap A_n$. In order to discuss this, let us start with axioms for the planar algebras. Following Jones \cite{jo3}, we have:

\index{planar tangle}
\index{planar algebra}

\begin{definition}
The planar algebras are defined as follows:
\begin{enumerate}
\item We consider rectangles in the plane, with the sides parallel to the coordinate axes, and taken up to planar isotopy, and we call such rectangles boxes.

\item A labeled box is a box with $2n$ marked points on its boundary, $n$ on its upper side, and $n$ on its lower side, for some integer $n\in\mathbb N$.

\item A tangle is labeled box, containing a number of labeled boxes, with all marked points, on the big and small boxes, being connected by noncrossing strings.

\item A planar algebra is a sequence of finite dimensional vector spaces $P=(P_n)$, together with linear maps $P_{n_1}\otimes\ldots\otimes P_{n_k}\to P_n$, one for each tangle, such that the gluing of tangles corresponds to the composition of linear maps.
\end{enumerate}
\end{definition}

In this definition we are using rectangles, but everything being up to isotopy, we could have used instead circles with marked points, as in \cite{jo3}. Our choice for using rectangles comes from the main examples that we have in mind, to be discussed below, where the planar algebra structure is best viewed by using rectangles, as above.

\bigskip

Let us also mention that Definition 15.18 is something quite simplified, based on \cite{jo3}. As explained in \cite{jo3}, in order for subfactors to produce planar algebras and vice versa, there are quite a number of supplementary axioms that must be added, and in view of this, it is perhaps better to start with something stronger than Definition 15.18, as basic axioms. However, as before with rectangles vs circles, our axiomatic choices here are mainly motivated by the concrete examples that we have in mind. More on this later.

\bigskip

As a basic example of a planar algebra, we have the Temperley-Lieb algebra:

\begin{theorem}
The Temperley-Lieb algebra $TL_N$, viewed as graded algebra
$$TL_N=(TL_N(n))_{n\in\mathbb N}$$
is a planar algebra, with the corresponding linear maps associated to the planar tangles
$$TL_N(n_1)\otimes\ldots\otimes TL_N(n_k)\to TL_N(n)$$
appearing by putting the various $TL_N(n_i)$ diagrams into the small boxes of the given tangle, which produces a $TL_N(n)$ diagram.
\end{theorem}

\begin{proof}
This is something trivial, which follows from definitions:

\medskip

(1) Assume indeed that we are given a planar tangle $\pi$, as in Definition 15.18, consisting of a box having $2n$ marked points on its boundary, and containing $k$ small boxes, having respectively $2n_1,\ldots,2n_k$ marked points on their boundaries, and then a total of $n+\Sigma n_i$ noncrossing strings, connecting the various $2n+\Sigma 2n_i$ marked points.

\medskip

(2) We want to associate to this tangle $\pi$ a linear map as follows:
$$T_\pi:TL_N(n_1)\otimes\ldots\otimes TL_N(n_k)\to TL_N(n)$$

For this purpose, by linearity, it is enough to construct elements as follows, for any choice of Temperley-Lieb diagrams $\sigma_i\in TL_N(n_i)$, with $i=1,\ldots,k$:
$$T_\pi(\sigma_1\otimes\ldots\otimes\sigma_k)\in TL_N(n)$$

(3) But constructing such an element is obvious, just by putting the various diagrams $\sigma_i\in TL_N(n_i)$ into the small boxes the given tangle $\pi$. Indeed, this procedure produces a certain diagram in $TL_N(n)$, that we can call $T_\pi(\sigma_1\otimes\ldots\otimes\sigma_k)$, as above.

\medskip

(4) Finally, we have to check that everything is well-defined up to planar isotopy, and that the gluing of tangles corresponds to the composition of linear maps. But both these checks are trivial, coming from the definition of $TL_N$, and we are done.
\end{proof}

As a conclusion to all this, $P=TL_N$ is indeed a planar algebra, but of somewhat ``trivial'' type, with the triviality coming from the fact that, in this case, the elements of $P$ are planar diagrams themselves, and so the planar structure appears trivially.

\bigskip

Moving ahead now, here is our second basic example of a planar algebra, which is also ``trivial'' in the above sense, with the elements of the planar algebra being planar diagrams themselves, but which appears in a bit more complicated way:

\index{Fuss-Catalan algebra}

\begin{theorem}
The Fuss-Catalan algebra $FC_{N,M}$, which appears by coloring the Temperley-Lieb diagrams with black/white colors, clockwise, as follows 
$$\circ\bullet\bullet\circ\circ\bullet\bullet\circ\ldots\ldots\ldots\circ\bullet\bullet\circ$$
and keeping those diagrams whose strings connect either $\circ-\circ$ or $\bullet-\bullet$, is a planar algebra, with again the corresponding linear maps associated to the planar tangles
$$FC_{N,M}(n_1)\otimes\ldots\otimes FC_{N,M}(n_k)\to FC_{N,M}(n)$$
appearing by putting the various $FC_{N,M}(n_i)$ diagrams into the small boxes of the given tangle, which produces a $FC_{N,M}(n)$ diagram.
\end{theorem}

\begin{proof}
The proof here is nearly identical to the proof of Theorem 15.19, with the only change appearing at the level of the colors. To be more precise:

\medskip

(1) Forgetting about upper and lower sequences of points, which must be joined by strings, a Temperley-Lieb diagram can be thought of as being a collection of strings, say black strings, which compose in the obvious way, with the rule that the value of the circle, which is now a black circle, is $N$. And it is this obvious composition rule that gives the planar algebra structure, as explained in the proof of Theorem 15.19. 

\medskip

(2) Similarly, forgetting about points, a Fuss-Catalan diagram can be thought of as being a collection of strings, which come now in two colors, black and white. These Fuss-Catalan diagrams compose then in the obvious way, with the rule that the value of the black circle is $N$, and the value of the white circle is $M$. And it is this obvious composition rule that gives the planar algebra structure, as before for $TL_N$.
\end{proof}

In general, the planar algebras can be more complicated than this, and we will see some further examples in a moment. However, the idea is very simple, namely ``the elements of a planar algebra are not necessarily diagrams, but they behave like diagrams".

\bigskip

In relation now with subfactors, the result, which extends Theorem 15.16, and which was found by Jones in \cite{jo3}, almost 20 years after \cite{jo1}, is as follows:

\index{higher commutant}
\index{planar algebra}

\begin{theorem} 
Given a subfactor $A_0\subset A_1$, the collection $P=(P_n)$ of linear spaces 
$$P_n=A_0'\cap A_n$$
has a planar algebra structure, extending the planar algebra structure of $TL_N$.
\end{theorem}

\begin{proof}
We know from Theorem 15.16 that we have an inclusion as follows, coming from the basic construction, and with $TL_N$ itself being a planar algebra:
$$TL_N\subset P$$

Thus, the whole point is that of proving that the trivial planar algebra structure of $TL_N$ extends into a planar algebra structure of $P$. But this can be done via a long algebraic study, and for the full computation here, we refer to Jones' paper \cite{jo3}.
\end{proof}

As a first illustration for the above result, we have:

\index{Temperley-Lieb}
\index{Fuss-Catalan algebra}

\begin{theorem}
We have the following universality results:
\begin{enumerate}
\item The Temperley-Lieb algebra $TL_N$ appears inside the planar algebra of any subfactor $A_0\subset A_1$ having index $N$.

\item The Fuss-Catalan algebra $FC_{N,M}$ appears inside the planar algebra of any subfactor $A_0\subset A_1$, in the presence of an intermediate subfactor $A_0\subset B\subset A_1$.
\end{enumerate}
\end{theorem}

\begin{proof}
Here the first assertion is something that we already know, from Theorem 15.16, and the second assertion is something quite standard as well, by carefully working out the basic construction for $A_0\subset A_1$, in the presence of an intermediate subfactor $A_0\subset B\subset A_1$. For details here, we refer to the paper of Bisch and Jones \cite{bjo}.
\end{proof}

As a free probability comment here, the Temperley-Lieb algebra, which appears by definition as the span of $NC_2$, is certainly a free probability object, and one way of being more concrete here is by saying that suitable fixed point subfactors associated to $S_N^+,O_N^+,U_N^+$ have planar algebra equal to $TL_N$. As in what regards the Fuss-Catalan algebra, this is related to the bicolored partitions appearing in the study of $H_N^+$, and more generally of $H_N^{s+}$, and again, the precise subfactor statement about this concerns fixed point subfactors associated to the quantum groups $H_N^{s+}$. More on this later.

\bigskip

The above results raise the question on whether any planar algebra produces a subfactor. The answer here is yes, but with many subtleties, and in order to talk about this, we first need to introduce a certain distinguished ${\rm II}_1$ factor $R$, as follows:

\index{R}

\begin{definition}
The Murray-von Neumann hyperfinite ${\rm II}_1$ factor is
$$R=\overline{\bigcup_iM_{n_i}(\mathbb C)}^{\,w}$$
independently of the choice of the algebras $M_{n_i}(\mathbb C)$, and of the embeddings between them.
\end{definition}

To be more precise, all this is based on two theorems of Murray and von Neumann \cite{mvo}, stating on one hand that when performing the above inductive limit construction we obtain, after taking the weak closure, a certain ${\rm II}_1$ factor, and on the other hand, that the factor that we obtain is independent on the choice of the algebras $M_{n_i}(\mathbb C)$, and of the embeddings between them. All this is certainly non-trivial, and even less trivial is the following theorem, coming as a continuation of the work in \cite{mvo}, due to Connes \cite{co2}:

\begin{theorem}
The Murray-von Neumann ${\rm II}_1$ factor $R$ is the unique ${\rm II}_1$ factor which is amenable, in the sense that we have a conditional expectation as follows:
$$E:B(H)\to R$$
In particular, for a discrete group $\Gamma$ we have $L(\Gamma)=R$ precisely when $\Gamma\neq\{1\}$ has the infinite conjugacy class (ICC) property, and is amenable.
\end{theorem}

\begin{proof}
This is something fairly complicated, to the point of causing troubles not only to mathematicians, and no surprise here, but to physicists as well. In case you know a good physicist, the best is to ask that physicist, but there is no guarantee here, guy might be well clueless on all this. So, read from time to time operator algebras, and once ready go through \cite{co2}. And in the meantime do not hesitate to ask around, this being a good test for distinguishing good physicists from first-class physicists.
\end{proof}

Jokes left aside now, what is difficult in the above is the proof of ``amenability implies hyperfiniteness''. Indeed, the converse can only be something standard, namely proving that a certain concrete algebra, $R$ from Definition 15.23, has a certain concrete property. As for the last assertion, this cannot be complicated either, because one of the possible definitions of the amenability of $\Gamma$ is in terms of an invariant mean $m:l^\infty(\Gamma)\to\mathbb C$, and this makes the connection with the expectation $E:B(l^2(\Gamma))\to L(\Gamma)$. See \cite{co2}.

\bigskip

Getting back now to subfactors, and to our questions regarding the correspondence between subfactors and planar algebras, these are difficult questions too, and the various answers to these questions can be summarized, a bit informally, as follows:

\begin{theorem}
We have the following results:
\begin{enumerate}
\item Any planar algebra with positivity produces a subfactor.

\item In particular, we have $TL$ and $FC$ type subfactors.

\item In the amenable case, and with $A_1=R$, the correspondence is bijective.

\item In general, we must take $A_1=L(F_\infty)$, and we do not have bijectivity.

\item The axiomatization of $P$, in the case $A_1=R$, is not known.
\end{enumerate}
\end{theorem}

\begin{proof}
All this is quite heavy, basically coming from the work of Popa in the 90s, using heavy functional analysis, the idea being as follows:

\medskip

(1) As already mentioned after Definition 15.18, our planar algebra axioms here are something quite simplified, based on \cite{jo3}. However, when getting back to Theorem 15.21, the conclusion is that the subfactor planar algebras there satisfy a number of supplementary ``positivity'' conditions, basically coming from the positivity of the ${\rm II}_1$ factor trace. And the point is that, with these positivity conditions axiomatized, we reach to something which is equivalent to Popa's axiomatization of the lattice of higher relative commutants $A_i'\cap A_j$ of the finite index subfactors, obtained in the 90s via heavy functional analysis. For the full story here, and details, we refer to Jones' paper \cite{jo3}.

\medskip

(2) The existence of the $TL_N$ subfactors, also known as ``$A_\infty$ subfactors'' in the literature, is something which was known for some time, since some early work of Popa on the subject. As for the existence of the $FC_{N,M}$ subfactors, this can be shown by using the intermediate subfactor picture, $A_0\subset B\subset A_1$, by composing two $A_\infty$ subfactors of suitable indices, $A_0\subset B$ and $B\subset A_1$. For the full story here, we refer to \cite{bjo}, \cite{jo3}.

\medskip

(3) This is something fairly heavy, as it is always the case with operator algebra results regarding hyperfiniteness and amenability, due to Popa \cite{pop}.

\medskip

(4) This is something a bit more recent, obtained by further building on the above-mentioned constructions of Popa, and we refer here to \cite{psh} and related work.

\medskip

(5) This is the big open question in subfactors. The story here goes back to Jones' original  paper \cite{jo1}, which contains at the end the question, due to Connes, of finding the possible values of the index for the irreducible subfactors of $R$. This question, which certainly looks much easier than (5) in the statement, is in fact still open, now 40 years after its formulation, and with no one having any valuable idea in dealing with it.
\end{proof}

\section*{15c. Basic examples} 

Let us discuss now some basic examples of subfactors, with concrete illustrations for all the above notions, constructions, and general theory. These examples will all come from group actions $G\curvearrowright Q$, which are assumed to be minimal, in the sense that:
$$(Q^G)'\cap Q=\mathbb C$$ 

As a starting point, we have the following result, due to Jones \cite{jo1}:

\index{Jones subfactor}

\begin{proposition}
Assuming that $G$ is a compact group, acting minimally on a ${\rm II}_1$ factor $Q$, and that $H\subset G$ is a subgroup of finite index, we have a subfactor
$$Q^G\subset Q^H$$
having index $N=[G:H]$, called Jones subfactor.
\end{proposition}

\begin{proof}
This is something standard, the idea being that the factoriality of $Q^G,Q^H$ comes from the minimality of the action, and that the index formula is clear.
\end{proof}

Along the same lines, we have the following result, due to Ocneanu:

\index{Ocneanu subfactor}

\begin{proposition}
Assuming that $G$ is a finite group, acting minimally on a ${\rm II}_1$ factor $Q$, we have a subfactor as follows,
$$Q\subset Q\rtimes G$$
having index $N=|G|$, called Ocneanu subfactor.
\end{proposition}

\begin{proof}
This is standard as well, the idea being that the factoriality of $Q\rtimes G$ comes from the minimality of the action, and that the index formula is clear.
\end{proof}

We have as well a third result of the same type, due to Wassermann \cite{was}, namely:

\index{Wassermann subfactor}

\begin{proposition}
Assuming that $G$ is a compact group, acting minimally on a ${\rm II}_1$ factor $Q$, and that $G\to PU_n$ is a projective representation, we have a subfactor
$$Q^G\subset (M_n(\mathbb C)\otimes Q)^G$$
having index $N=n^2$, called Wassermann subfactor.
\end{proposition}

\begin{proof}
As before, the idea is that the factoriality of $Q^G,(M_n(\mathbb C)\otimes Q)^G$ comes from the minimality of the action, and the index formula is clear.
\end{proof}

The above subfactors look quite related, and indeed they are, due to:

\begin{theorem}
The Jones, Ocneanu and Wassermann subfactors are all of the same nature, and can be written as follows,
$$\left( Q^G\subset Q^H\right)\,\simeq\, \left( ({\mathbb C}\otimes Q)^G\subset (l^\infty(G/H)\otimes Q)^G\right)$$
$$\left( Q\subset Q\rtimes G\right)\,\simeq\,  \left( (l^\infty (G)\otimes Q)^G\subset ({\mathcal L} (l^2(G))\otimes Q)^G\right)$$
$$\left( Q^G\subset (M_n(\mathbb C) \otimes Q)^G\right)\,\simeq\, \left( ({\mathbb C}\otimes Q)^G\subset (M_n(\mathbb C)\otimes Q)^G\right)$$
with standard identifications for the various tensor products and fixed point algebras.
\end{theorem}

\begin{proof}
This is something standard, modulo several standard identifications. We will explain all this more in detail later, after unifying these subfactors.
\end{proof}

In order to unify now the above constructions of subfactors, the idea is quite clear. Given a compact group $G$, acting minimally on a ${\rm II}_1$ factor $Q$, and an inclusion of finite dimensional algebras $B_0\subset B_1$, endowed as well with an action of $G$, we would like to construct a kind of generalized Wassermann subfactor, as follows:
$$(B_0\otimes Q)^G\subset (B_1\otimes Q)^G$$

In order to do this, we must talk first about the finite dimensional algebras $B$, and about inclusions of such algebras $B_0\subset B_1$. Let us start with the following definition:

\begin{definition}
Associated to any finite dimensional algebra $B$ is its canonical trace, obtained by composing the left regular representation with the trace of $\mathcal L(B)$:
$$tr:B\subset\mathcal L(B)\to\mathbb C$$
We say that an inclusion of finite dimensional algebras $B_0\subset B_1$ is Markov if it commmutes with the canonical traces of $B_0,B_1$.
\end{definition}

As a basic illustration for this, any inclusion of type $\mathbb C\subset B$ is Markov. In general, if we write $B_0=C(X_0)$ and $B_1=C(X_1)$, then the inclusion $B_0\subset B_1$ must come from a certain fibration $X_1\to X_0$, and the inclusion $B_0\subset B_1$ is Markov precisely when the fibration $X_1\to X_0$ commutes with the respective counting measures.

\bigskip

We will be back to Markov inclusions and their various properties on several occasions, in what follows. For our next purposes here, we just need the following result:

\begin{proposition}
Given a Markov inclusion of finite dimensional algebras $B_0\subset B_1$ we can perform to it the basic construction, as to obtain a Jones tower
$$B_0\subset_{e_1}B_1\subset_{e_2}B_2\subset_{e_3}B_3\subset\ldots\ldots$$
exactly as we did in the above for the inclusions of ${\rm II}_1$ factors.
\end{proposition}

\begin{proof}
This is something standard, from \cite{jo1}, by following the computations in the above, from the case of the ${\rm II}_1$ factors, and with everything extending well. It is of course possible to do something more general here, unifying the constructions for the inclusions of ${\rm II}_1$ factors $A_0\subset A_1$, and for the inclusions of Markov inclusions of finite dimensional algebras $B_0\subset B_1$, but we will not need this degree of generality, in what follows.
\end{proof}

With these ingredients in hand, getting back now to the Jones, Ocneanu and Wassermann subfactors, from Theorem 15.29, the point is that these constructions can be unified, and then further studied, the final result on the subject being as follows:

\index{fixed point subfactor}
\index{Markov inclusion}
\index{Wassermann subfactor}

\begin{theorem}
Let $G$ be a compact group, and $G\to Aut(Q)$ be a minimal action on a ${\rm II}_1$ factor. Consider a Markov inclusion of finite dimensional algebras
$$B_0\subset B_1$$
and let $G\to Aut(B_1)$ be an action which leaves invariant $B_0$, and which is such that its restrictions to the centers of $B_0$ and $B_1$ are ergodic. We have then a subfactor
$$(B_0\otimes Q)^G\subset (B_1\otimes Q)^G$$
of index $N=[B_1:B_0]$, called generalized Wassermann subfactor, whose Jones tower is 
$$(B_1\otimes Q)^G\subset(B_2\otimes Q)^G\subset(B_3\otimes Q)^G\subset\ldots$$
where $\{ B_i\}_{i\geq 1}$ are the algebras in the Jones tower for $B_0\subset B_1$, with the canonical actions of $G$ coming from the action $G\to Aut(B_1)$, and whose planar algebra is given by:
$$P_k=(B_0'\cap B_k)^G$$
These subfactors generalize the Jones, Ocneanu and Wassermann subfactors.
\end{theorem}

\begin{proof}
This is something which is routine, following Wassermann \cite{was}, and we will be back to this in a moment, with details, directly in a more general setting.
\end{proof}

In addition to the Jones, Ocneanu and Wassermann subfactors, discussed and unified in the above, we have the Popa subfactors, which are constructed as follows:

\index{Popa subfactor}

\begin{proposition}
Given a discrete group $\Gamma=<g_1,\ldots,g_n>$, acting faithfully via outer automorphisms on a ${\rm II}_1$ factor $P$, we have the following ``diagonal'' subfactor
$$\left\{ \begin{pmatrix}
g_1(q)\\
&\ddots\\
&& g_n(q)
\end{pmatrix} \Big| q\in P\right\} \subset M_n(P)$$
having index $N=n^2$, called Popa subfactor.
\end{proposition}

\begin{proof}
This is something standard, a bit as for the Jones, Ocneanu and Wassermann subfactors, with the result basically coming from the work of Popa \cite{pop}.
\end{proof}

In order to unify now Theorem 15.32 and Proposition 15.33, we can use:

\begin{definition}
A coaction of a Woronowicz algebra $A$ on a finite von Neumann algebra $Q$ is an injective morphism $\Phi:Q\to Q\otimes A''$ satisfying the following conditions:
\begin{enumerate}
\item Coassociativity: $(\Phi\otimes id)\Phi=(id\otimes\Delta)\Phi$.

\item Trace equivariance: $(tr\otimes id)\Phi=tr(.)1$.

\item Smoothness: $\overline{\mathcal Q}^{\,w}=Q$, where $\mathcal Q=\Phi^{-1}(Q\otimes_{alg}\mathcal A)$.
\end{enumerate}
\end{definition}

These conditions come from what happens in the commutative case, $A=C(G)$, where they correspond to the usual associativity, trace equivariance and smoothness of the corresponding action $G\curvearrowright Q$. Along the same lines, we have as well:

\index{ergodic coaction}
\index{faithful coaction}
\index{minimal coaction}

\begin{definition}
A coaction $\Phi:Q\to Q\otimes A''$ as above is called:
\begin{enumerate}
\item Ergodic, if the algebra $Q^\Phi=\left\{p\in Q\big|\Phi(p)=p\otimes1\right\}$ reduces to $\mathbb C$.

\item Faithful, if the span of $\left\{(f\otimes id)\Phi(Q)\big|f\in Q_*\right\}$ is dense in $A''$.

\item Minimal, if it is faithful, and satisfies $(Q^\Phi)'\cap Q=\mathbb C$.
\end{enumerate} 
\end{definition}

Observe that the minimality of the action implies in particular that the fixed point algebra $Q^\Phi$ is a factor. Thus, we are getting here to the case that we are interested in, actions producing factors, via their fixed point algebras. Next, we have:

\index{fixed point algebra}

\begin{proposition}
Consider a Woronowicz algebra $A=(A,\Delta,S)$, and denote by $A_\sigma$ the Woronowicz algebra $(A,\sigma\Delta ,S)$, where $\sigma$ is the flip. Given two coactions
$$\beta:B\to B\otimes A\quad,\quad 
\pi:Q\to Q\otimes A_\sigma$$
with $B$ being finite dimensional, the following linear map, while not being multiplicative in general, is coassociative with respect to the comultiplication $\sigma\Delta$ of $A_\sigma$,
$$\beta\odot\pi:B\otimes Q\to B\otimes Q\otimes A_\sigma$$
$$b\otimes p\to \pi (p)_{23}((id\otimes S)\beta(b))_{13}$$
and its fixed point space, which is by definition the following linear space,
$$(B\otimes Q)^{\beta\odot\pi}=\left\{x\in B\otimes Q\Big|(\beta\odot\pi )x=x\otimes 1\right\}$$
is then a von Neumann subalgebra of $B\otimes Q$. 
\end{proposition}

\begin{proof}
This is something standard, which follows from a straightforward algebraic verification. As mentioned in the statement, to be noted is that the tensor product coaction $\beta\odot\pi$ is not multiplicative in general. See my old subfactor papers.
\end{proof}

Our first task is to investigate the factoriality of such algebras, and we have here:

\begin{theorem}
If $\beta:B\to B\otimes A$ is a coaction and $\pi:Q\to Q\otimes A_\sigma$ is a minimal coaction, then the following conditions are equivalent:
\begin{enumerate}
\item The von Neumann algebra $(B\otimes Q)^{\beta\odot\pi}$ is a factor.

\item The coaction $\beta$ is centrally ergodic, $Z(B)\cap B^\beta=\mathbb C$.
\end{enumerate}
\end{theorem}

\begin{proof}
This is something standard, the idea being as follows:

\medskip

(1) Our first claim, which is something whose proof is a routine verification, is that the following diagram is a non-degenerate commuting square:
$$\begin{matrix}
Q&\subset&B\otimes Q\\ 
\cup &\ &\cup \\
Q^\pi&\subset&(B\otimes Q)^{\beta\odot\pi}
\end{matrix}$$

(2) In order to prove now the result, it is enough to check the following equality, between subalgebras of the von Neumann algebra $B\otimes Q$:
$$Z((B\otimes Q)^{\beta\odot\pi})=(Z(B)\cap B^\beta)\otimes 1$$

But this follows from the non-degeneracy of the above commuting square.
\end{proof}
 
With the above results in hand, we can now formulate our main theorem regarding the fixed point subfactors, of the most possible general type, as follows:

\index{fixed point subfactor}
\index{Wassermann subfactor}

\begin{theorem}
Let $G$ be a compact quantum group, and $G\to Aut(Q)$ be a minimal action on a ${\rm II}_1$ factor. Consider a Markov inclusion of finite dimensional algebras
$$B_0\subset B_1$$
and let $G\to Aut(B_1)$ be an action which leaves invariant $B_0$ and which is such that its restrictions to the centers of $B_0$ and $B_1$ are ergodic. We have then a subfactor
$$(B_0\otimes Q)^G\subset (B_1\otimes Q)^G$$
of index $N=[B_1:B_0]$, called generalized Wassermann subfactor, whose Jones tower is 
$$(B_1\otimes Q)^G\subset(B_2\otimes Q)^G\subset(B_3\otimes Q)^G\subset\ldots$$
where $\{ B_i\}_{i\geq 1}$ are the algebras in the Jones tower for $B_0\subset B_1$, with the canonical actions of $G$ coming from the action $G\to Aut(B_1)$, and whose planar algebra is given by:
$$P_k=(B_0'\cap B_k)^G$$
These subfactors generalize the Jones, Ocneanu, Wassermann and Popa subfactors.
\end{theorem}

\begin{proof}
This is indeed routine, based on the above general theory and results.
\end{proof}

The above result is important in connection with probability questions, because our usual character computations for $G$, for instance in the case where $G\subset U_N^+$ is easy, take place in the associated planar algebra $P_k=(B_0'\cap B_k)^G$. More on this later.

\bigskip

This was for the basic theory of the fixed point subfactors. Many more things can be said about them, notably with an axiomatization of the planar algebras that we can obtain in this way, as being the subalgebras of Jones' bipartite graph planar algebras from \cite{jo4}, and also with a number of results and open questions regarding amenability. For more on all this, have a look at my old subfactor papers, and their follow-ups.

\section*{15d. Circular measures}

We would like to discuss now various structure and classification questions for subfactors, all interesting questions, related to physics, regarded from a probabilistic viewpoint. In order to get started, we need invariants for our subfactors. We have the choice here between algebraic and analytic invariants, the situation being as follows:

\index{principal graph}
\index{fusion algebra}
\index{Poincar\'e series}
\index{spectral measure}

\begin{definition}
Associated to any finite index subfactor $A_0\subset A_1$, having planar algebra $P=(P_n)$, are the following invariants:
\begin{enumerate}
\item Its principal graph $\Gamma$, which describes the inclusions $P_0\subset P_1\subset P_2\subset\ldots\,$, with the reflections coming from basic constructions removed.

\item Its fusion algebra $F$, which describes the fusion rules for the various types of bimodules that can appear, namely $A_0-A_0$, $A_0-A_1$, $A_1-A_0$, $A_1-A_1$.

\item Its Poincar\'e series $f$, which is the generating series of the graded components of the planar algebra, $f(z)=\sum_n\dim(P_n)z^n$.

\item Its spectral measure $\mu$, which is the probability measure having as moments the dimensions of the planar algebra components, $\int x^nd\mu(x)=\dim(P_n)$.
\end{enumerate}
\end{definition}

This definition is of course something a bit informal, and there is certainly some work to be done, in order to fully define all these invariants $\Gamma,F,f,\mu$, and to work out the precise relation between them. We will be back to this later, but for the moment, let us keep in mind the fact that associated to a given subfactor $A_0\subset A_1$ are several invariants, which are not exactly equivalent, but are definitely versions of the same thing, the ``combinatorics of the subfactor'', coming in algebraic or analytic flavors.

\bigskip

More in detail now, let us begin by explaining how the principal graph $\Gamma$ is constructed. Consider a finite index irreducible subfactor $A_0\subset A_1$, with associated planar algebra $P_n=A_0'\cap A_n$, and let us look at the following system of inclusions:
$$P_0\subset P_1\subset P_2\subset\ldots$$

By taking the Bratelli diagram of this system of inclusions, and then deleting the reflections coming from basic constructions, we obtain a certain graph $\Gamma$, called principal graph of $A_0\subset A_1$. The main properties of $\Gamma$ can be summarized as follows:

\index{amenable subfactor}

\begin{proposition}
The principal graph $\Gamma$ has the following properties:
\begin{enumerate}
\item The higher relative commutant $P_n=A_0'\cap A_n$ is isomorphic to the abstract vector space spanned by the $2n$-loops on $\Gamma$ based at the root.

\item In the amenable case, where $A_1=R$ and when the subfactor is ``amenable'', the index of $A_0\subset A_1$ is given by $N=||\Gamma||^2$.
\end{enumerate} 
\end{proposition}

\begin{proof}
This is something standard, the idea being as follows:

\medskip

(1) The statement here, which explains among others the relation between the principal graph $\Gamma$, and the other subfactor invariants, from Definition 15.39, comes from the definition of the principal graph, as a Bratelli diagram, with the reflections removed.

\medskip

(2) This is actually a quite subtle statement, but for our purposes here, we can take the equality $N=||\Gamma||^2$, which reminds the Kesten amenability condition for discrete groups, or discrete quantum groups, as a definition for the amenability of the subfactor.
\end{proof}

As an illustration for all this, let us discuss the case of the small index subfactors, $N\in[1,4]$. Following Jones \cite{jo1} and related work, we first have the following result:

\index{ADE}
\index{Coxeter-Dynkin}

\begin{theorem}
The index of subfactors is subject to the condition
$$N\in\left\{4\cos^2\left(\frac{\pi}{n}\right)\Big|n\geq3\right\}\cup[4,\infty]$$
and at $N\leq4$, the principal graph must be one of the Coxeter-Dynkin ADE graphs.
\end{theorem}

\begin{proof}
Here the first assertion is something that we already know, and the second assertion comes from this, the key ingredient being the fact that in index $N\leq4$, and more generally in the amenable case, we must have $N=||\Gamma||^2$. See \cite{jo1}.
\end{proof}

More in detail now, the usual Coxeter-Dynkin ADE graphs are as follows:
$$A_k=\bullet-\circ-\circ\cdots\circ-\circ-\circ\hskip20mm A_{\infty}=\bullet-\circ-\circ-\circ\cdots\hskip7mm$$
\vskip-7mm
$$D_k=\bullet-\circ-\circ\dots\circ-
\begin{matrix}\ \circ\cr\ |\cr\ \circ \cr\ \cr\  \end{matrix}-\circ\hskip70mm$$
\vskip-7mm
$$\ \ \ \ \ \ \ \tilde{A}_{2k}=
\begin{matrix}
\circ&\!\!\!\!-\circ-\circ\cdots\circ-\circ-&\!\!\!\!\circ\cr
|&&\!\!\!\!|\cr
\bullet&\!\!\!\!-\circ-\circ-\circ-\circ-&\!\!\!\!\circ\cr\cr\cr\end{matrix}\hskip20mm \tilde{A}_\infty=
\begin{matrix}
\circ&\!\!\!\!-\circ-\circ-\circ\cdots\cr
|&\cr
\bullet&\!\!\!\!-\circ-\circ-\circ\cdots\cr\cr\cr\end{matrix}
\hskip15mm$$
\vskip-9mm
$$\;\tilde{D}_k=\bullet-
\begin{matrix}\circ\cr|\cr\circ\cr\ \cr\ \end{matrix}-\circ\dots\circ-
\begin{matrix}\ \circ\cr\ |\cr\ \circ \cr\ \cr\  \end{matrix}-\circ \hskip20mm D_\infty=\bullet-
\begin{matrix}\circ\cr|\cr\circ\cr\ \cr\ \end{matrix}-\circ-\circ\cdots\hskip7mm$$
\vskip-7mm

There are as well a number of exceptional Coxeter-Dynkin graphs. First we have:
$$E_6=\bullet-\circ-
\begin{matrix}\circ\cr|\cr\circ\cr\ \cr\ \end{matrix}-
\circ-\circ\hskip71mm$$
\vskip-13mm
$$E_7=\bullet-\circ-\circ-
\begin{matrix}\circ\cr|\cr\circ\cr\ \cr\ \end{matrix}-
\circ-\circ\hskip18mm$$
\vskip-15mm
$$\hskip30mm E_8=\bullet-\circ-\circ-\circ-
\begin{matrix}\circ\cr|\cr\circ\cr\ \cr\ \end{matrix}-
\circ-\circ$$
\vskip-5mm

Finally, we have index 4 versions of the above exceptional graphs, as follows:
$$\tilde{E}_6=\bullet-\circ-\begin{matrix}
\circ\cr|
\cr\circ\cr|&\cr\circ&\!\!\!\!-\ \circ\cr\ \cr\   \cr\ \cr\ \end{matrix}-\circ\hskip71mm$$
\vskip-22mm
$$\tilde{E}_7=\bullet-\circ-\circ-
\begin{matrix}\circ\cr|\cr\circ\cr\ \cr\ \end{matrix}-
\circ-\circ-\circ\hskip18mm$$
\vskip-15mm
$$\hskip30mm \tilde{E}_8=\bullet-\circ-\circ-\circ-\circ-
\begin{matrix}\circ\cr|\cr\circ\cr\ \cr\ \end{matrix}-
\circ-\circ$$
\vskip-5mm

The subfactors of index $\leq 4$ were intensively studied in the 80s and early 90s, and about 10 years after Jones' paper \cite{jo1}, a complete classification result was found, with contributions by many authors. A simplified form of this result is as follows:

\index{principal graph}

\begin{theorem}
The principal graphs of subfactors of index $\leq 4$ are:
\begin{enumerate}
\item Index $<4$ graphs: $A_k$, $D_{even}$, $E_6$, $E_8$. 

\item Index $4$ finite graphs: $\tilde{A}_{2k}$, $\tilde{D}_k$, $\tilde{E}_6$, $\tilde{E}_7$, $\tilde{E}_8$.

\item Index $4$ infinite graphs: $A_\infty$, $\tilde{A}_\infty$, $D_\infty$.
\end{enumerate}
\end{theorem}

\begin{proof}
As already mentioned, this is something quite heavy, with contributions by many authors, and notably Ocneanu. Observe that, as a subtlety of subfactor theory, the graphs $D_{odd}$ and $E_7$ don't appear in the above list. For a discussion, see \cite{jo3}.
\end{proof}

Getting now to probability, as a first observation, we know that $A_\infty$ and $\tilde{A}_\infty$ are the graphs that we previously called $\mathbb N$ and $\mathbb Z$, and studied before in some detail. So, based on our previous work for these graphs, let us formulate the following definition:

\begin{definition}
The positive spectral measure $\nu$ of a rooted bipartite graph $X$ is the real probability measure having as moments
$$\int_\mathbb Rx^kd\nu(x)=L_{2k}$$
where $L_{2k}$ is the number of $2k$-loops based at the root. Equivalently, we must have
$$\int_\mathbb R\frac{1}{1-xz}\,d\nu(x)=f(z)$$
where $f(z)=\sum_{k=0}^\infty L_{2k}z^k$ is the Poincar\'e series of $X$.
\end{definition}

Here the existence of $\nu$, and the fact that this is indeed a positive measure, meaning a measure supported on $[0,\infty)$, comes from the following simple fact:

\begin{proposition}
The positive spectral measure of a rooted bipartite graph $X$ is given by the following formula, with $d$ being the adjacency matrix of the graph,
$$\nu=law(d^2)$$
and with the probabilistic computation being with respect to the expectation 
$$A\to<A>$$
with $<A>$ being the $(*,*)$-entry of a matrix $A$, where $*$ is the root.
\end{proposition}

\begin{proof}
With the above conventions, we have the following computation:
$$f(z)
=\sum_{k=0}^\infty L_{2k}z^k
=\sum_{k=0}^\infty\left<d^{2k}\right>z^k
=\left<\frac{1}{1-d^2z}\right>$$

But this shows that we have $\nu=law(d^2)$, as desired.
\end{proof}

Let us introduce as well the following notion, coming from Jones' paper \cite{jo5}:

\index{circular measure}
\index{circular spectral measure}

\begin{definition}
The circular measure $\varepsilon$ of a rooted bipartite graph $X$ is given by
$$d\varepsilon(q)=d\nu((q+q^{-1})^2)$$
where $\nu$ is the associated positive spectral measure.
\end{definition}

As a basic illustration for this, assume that $\nu$ is a discrete measure, supported by $n$ positive numbers $x_1<\ldots<x_n$, with corresponding densities $p_1,\ldots,p_n$:
$$\nu=\sum_ip_i\delta_{x_i}$$

For each $i\in\{1,\ldots,n\}$ the equation $(q+q^{-1})^2=x_i$ has then four solutions, that we can denote $q_i,q_i^{-1},-q_i,-q_i^{-1}$. And with this notation, we have:
$$\varepsilon=\frac{1}{4}\sum_ip_i\left(\delta_{q_i}+\delta_{q_i^{-1}}+\delta_{-q_i}+\delta_{-q_i^{-1}}\right)$$

In general, the basic properties of $\varepsilon$ can be summarized as follows:

\begin{theorem}
The circular measure has the following properties:
\begin{enumerate}
\item $\varepsilon$ has equal density at $q,q^{-1},-q,-q^{-1}$.

\item The odd moments of $\varepsilon$ are $0$.

\item The even moments of $\varepsilon$ are half-integers.

\item When $X$ has norm $\leq 2$, $\varepsilon$ is supported by the unit circle.

\item When $X$ is finite, $\varepsilon$ is discrete.

\item If $K$ is a solution of $d=K+K^{-1}$, then $\varepsilon=law(K)$. 
\end{enumerate}
\end{theorem}

\begin{proof}
These results can be deduced from definitions, the idea being that (1-5) are trivial, and that (6) follows from the formula of $\nu$ from Proposition 15.44.
\end{proof}

Getting now to computations, we first have the following result:

\begin{theorem}
The circular measure of the basic index $4$ graph, namely 
$$\begin{matrix}
&\circ&\!\!\!\!-\circ-\circ\cdots\circ-\circ-&\!\!\!\!\circ\cr
\tilde{A}_{2n}=&|&&\!\!\!\!|\cr
&\bullet&\!\!\!\!-\circ-\circ-\circ-\circ-&\!\!\!\!\circ\cr\cr\cr\end{matrix}$$
\vskip-10mm

\noindent is the uniform measure on the $2n$-roots of unity.
\end{theorem}

\begin{proof}
Let us identify the vertices of $X=\tilde{A}_{2n}$ with the group $\{w^k\}$ formed by the $2n$-th roots of unity in the complex plane, where $w=e^{\pi i/n}$. The adjacency matrix of $X$ acts then on the functions $f\in C(X)$ in the following way:
$$df(w^s)=f(w^{s-1})+f(w^{s+1})$$

But this shows that we have $d=K+K^{-1}$, where $K$ is given by:
$$Kf(w^s)=f(w^{s+1})$$

Thus we can use Proposition 15.44 and Theorem 15.46 (6), and we get:
$$\varepsilon=law(K)$$

But this is the uniform measure on the $2n$-roots of unity, as claimed.
\end{proof}

Observe that Theorem 15.47 is something quite interesting, I mean you can try to do a few loop computations for $\tilde{A}_{2n}$, with bare hands, and these will lead you nowhere. Moving on now, with this circular measure idea of Jones adopted, in order to deal with the other ADE graphs, let us introduce the following densities:
$$\alpha=Re(1-q^2)$$
$$\beta=Re(1-q^4)$$
$$\gamma=Re(1-q^6)$$

We have then the following result, $d_n$ being the uniform measure on the $2n$-th roots of unity, and $d_n'=2d_{2n}-d_n$ being the uniform measure on the odd $4n$-roots of unity:

\index{ADE graph}
\index{circular measure}
\index{cyclotomic measure}

\begin{theorem}
The circular measures of ADE graphs are given by:
\begin{enumerate}
\item $A_{n-1}\to\alpha_n$.

\item $\tilde{A}_{2n}\to d_n$.

\item $D_{n+1}\to\alpha_n'$.

\item $\tilde{D}_{n+2}\to (d_n+d_1')/2$.

\item $E_6\to\alpha_{12}+(d_{12}-d_6-d_4+d_3)/2$.

\item $E_7\to\beta_9'+(d_1'-d_3')/2$.

\item $E_8\to\alpha_{15}'+\gamma_{15}'-(d_5'+d_3')/2$.

\item $\tilde{E}_{n+3}\to (d_n+d_3+d_2-d_1)/2$.
\end{enumerate}
\end{theorem}

\begin{proof}
This is something which can be proved in three steps, as follows:

\medskip

(1) For the simplest graph, namely the circle $\tilde{A}_{2n}$, we already have the result, from Theorem 15.47, with the proof there being something elementary.

\medskip

(2) For the other non-exceptional graphs, that is, of type A and D, the same method works, namely direct loop counting, with some matrix tricks.

\medskip

(3) As for the exceptional graphs, of type E, basically the same method works, but the computations are more complicated, using some root of unity know-how. See \cite{bbi}.
\end{proof}

For more on all this, and conformal field theory implications, see \cite{epu}.

\section*{15e. Exercises}

Welcome to subfactors, and no way back down, and as exercises, we have:

\begin{exercise}
Learn from Murray-von Neumann about ${\rm II}_1$ factors.
\end{exercise}

\begin{exercise}
Learn from Temperley and Lieb about their algebra.
\end{exercise}

\begin{exercise}
Learn the correspondence between subfactors and planar algebras. 
\end{exercise}

\begin{exercise}
Learn about the hyperfinite ${\rm II}_1$ factor $R$, as much as you can.
\end{exercise}

\begin{exercise}
Read the papers of Wassermann, on subfactors and related topics.
\end{exercise}

\begin{exercise}
Have a look as well at my old subfactor papers, all good stuff.
\end{exercise}

\begin{exercise}
Learn, from Morrison and others, about small index subfactors.
\end{exercise}

\begin{exercise}
Develop Weingarten calculus for suitable subfactors.
\end{exercise}

As bonus exercise, mandatory, read Jones's papers \cite{jo1}, \cite{jo2}, \cite{jo3}, \cite{jo4}, \cite{jo5}.

\chapter{Fourier models}

\section*{16a. Hadamard matrices}

Time to finish this book, and we kept the best for the end. We will be talking about Hadamard matrices, which are something quite deep, virtually related to everything that we have been doing in this book, in every single chapter, yet very easy to define:

\index{Hadamard matrix}

\begin{definition}
A complex Hadamard matrix is a square matrix 
$$H\in M_N(\mathbb C)$$
whose entries are on the unit circle, and whose rows are pairwise orthogonal.
\end{definition} 

Observe that the orthogonality condition tells us that the rescaled matrix $U=H/\sqrt{N}$ must be unitary. Thus, these matrices form a real algebraic manifold, given by:
$$X_N=M_N(\mathbb T)\cap\sqrt{N}U_N$$

The basic example is the Fourier matrix, $F_N=(w^{ij})$ with $w=e^{2\pi i/N}$.  With indices $i,j\in\{0,1,\ldots,N-1\}$, this matrix is as follows:
$$F_N=\begin{pmatrix}
1&1&1&\ldots&1\\
1&w&w^2&\ldots&w^{N-1}\\
1&w^2&w^4&\ldots&w^{2(N-1)}\\
\vdots&\vdots&\vdots&&\vdots\\
1&w^{N-1}&w^{2(N-1)}&\ldots&w^{(N-1)^2}
\end{pmatrix}$$

More generally, we have as example the Fourier coupling of any finite abelian group $G$, regarded via the isomorphism $G\simeq\widehat{G}$ as a square matrix, as follows:

\index{Fourier matrix}

\begin{theorem}
Given a finite abelian group $G$, with dual group $\widehat{G}=\{\chi:G\to\mathbb T\}$, consider the Fourier coupling $\mathcal F_G:G\times\widehat{G}\to\mathbb T$, given by $(i,\chi)\to\chi(i)$.
\begin{enumerate}
\item Via the standard isomorphism $G\simeq\widehat{G}$, this Fourier coupling can be regarded as a square matrix, $F_G\in M_G(\mathbb T)$, which is a complex Hadamard matrix.

\item For the cyclic group $G=\mathbb Z_N$ we obtain in this way, via the standard identification $\mathbb Z_N=\{1,\ldots,N\}$, the standard Fourier matrix, $F_N=(w^{ij})$ with $w=e^{2\pi i/N}$.

\item In general, when using a decomposition $G=\mathbb Z_{N_1}\times\ldots\times\mathbb Z_{N_k}$, the corresponding Fourier matrix is given by $F_G=F_{N_1}\otimes\ldots\otimes F_{N_k}$.
\end{enumerate}
\end{theorem}

\begin{proof}
This follows indeed from some basic facts from group theory:

\medskip

(1) With the identification $G\simeq\widehat{G}$ made our matrix is given by $(F_G)_{i\chi}=\chi(i)$, and the scalar products between the rows are computed as follows:
$$<R_i,R_j>
=\sum_\chi\chi(i)\overline{\chi(j)}
=\sum_\chi\chi(i-j)
=|G|\cdot\delta_{ij}$$

Thus, we obtain indeed a complex Hadamard matrix.

\medskip

(2) This follows from the well-known and elementary fact that, via the identifications $\mathbb Z_N=\widehat{\mathbb Z_N}=\{1,\ldots,N\}$, the Fourier coupling here is as follows, with $w=e^{2\pi i/N}$:
$$(i,j)\to w^{ij}$$

(3) We use here the following well-known formula, for the duals of products: 
$$\widehat{H\times K}=\widehat{H}\times\widehat{K}$$

At the level of the corresponding Fourier couplings, we obtain from this:
$$F_{H\times K}=F_H\otimes F_K$$

Now by decomposing $G$ into cyclic groups, as in the statement, and by using (2) for the cyclic components, we obtain the formula in the statement.
\end{proof}

In relation now with quantum groups, which will eventually lead us into von Neumann algebras, subfactors, planar algebras, spin models and many more, and into probability too, of rather free type, no worries for that, the starting observation is as follows:

\begin{proposition}
If $H\in M_N(\mathbb C)$ is Hadamard, the rank one projections 
$$P_{ij}=Proj\left(\frac{H_i}{H_j}\right)$$
where $H_1,\ldots,H_N\in\mathbb T^N$ are the rows of $H$, form a magic unitary.
\end{proposition}

\begin{proof}
This is clear, the verification for the rows being as follows:
$$\left<\frac{H_i}{H_j},\frac{H_i}{H_k}\right>
=\sum_l\frac{H_{il}}{H_{jl}}\cdot\frac{H_{kl}}{H_{il}}
=\sum_l\frac{H_{kl}}{H_{jl}}
=N\delta_{jk}$$

The verification for the columns is similar, as follows:
$$\left<\frac{H_i}{H_j},\frac{H_k}{H_j}\right>
=\sum_l\frac{H_{il}}{H_{jl}}\cdot\frac{H_{jl}}{H_{kl}}
=\sum_l\frac{H_{il}}{H_{kl}}
=N\delta_{ik}$$

Thus, we have indeed a magic unitary, as claimed. 
\end{proof}

We can proceed now in the same way as we did with before the Weyl matrices, namely constructing a model of $C(S_N^+)$, and performing the Hopf image construction:

\begin{definition}
To any Hadamard matrix $H\in M_N(\mathbb C)$ we associate the quantum permutation group $G\subset S_N^+$ given by the fact that $C(G)$ is the Hopf image of
$$\pi:C(S_N^+)\to M_N(\mathbb C)\quad,\quad 
u_{ij}\to Proj\left(\frac{H_i}{H_j}\right)$$
where $H_1,\ldots,H_N\in\mathbb T^N$ are the rows of $H$.
\end{definition}

Summarizing, we have a construction $H\to G$, and our claim is that this construction is something really useful, with $G$ encoding the combinatorics of $H$. To be more precise, our claim is that ``$H$ can be thought of as being a kind of Fourier matrix for $G$''. 

\bigskip

This is of course quite interesting, philosophically speaking. There are several results supporting this, with the main evidence coming from the following result:

\begin{theorem}
The construction $H\to G$ has the following properties:
\begin{enumerate}
\item For a Fourier matrix $H=F_G$ we obtain the group $G$ itself, acting on itself.

\item For $H\not\in\{F_G\}$, the quantum group $G$ is not classical, nor a group dual.

\item For a tensor product $H=H'\otimes H''$ we obtain a product, $G=G'\times G''$.
\end{enumerate}
\end{theorem}

\begin{proof}
All this material is standard, and elementary, as follows:

\medskip

(1) Let us first discuss the cyclic group case, where our Hadamard matrix is a usual Fourier matrix, $H=F_N$. Here the rows of $H$ are given by $H_i=\rho^i$, where:
$$\rho=(1,w,w^2,\ldots,w^{N-1})$$

Thus, we have the following formula, for the associated magic basis:
$$\frac{H_i}{H_j}=\rho^{i-j}$$

It follows that the corresponding rank 1 projections $P_{ij}=Proj(H_i/H_j)$ form a circulant matrix, all whose entries commute. Since the entries commute, the corresponding quantum group must satisfy $G\subset S_N$. Now by taking into account the circulant property of $P=(P_{ij})$ as well, we are led to the conclusion that we have:
$$G=\mathbb Z_N$$

In the general case now, where $H=F_G$, with $G$ being an arbitrary finite abelian group, the result can be proved either by extending the above proof, of by decomposing $G=\mathbb Z_{N_1}\times\ldots\times\mathbb Z_{N_k}$ and using (3) below, whose proof is independent from the rest.

\medskip

(2) This is something more tricky, needing some general study of the representations whose Hopf images are commutative, or cocommutative. For the story here, see \cite{bjb}.

\medskip

(3) Assume that we have a tensor product $H=H'\otimes H''$, and let $G,G',G''$ be the associated quantum permutation groups. We have then a diagram as follows:
$$\xymatrix@R=45pt@C25pt{
C(S_{N'}^+)\otimes C(S_{N''}^+)\ar[r]&C(G')\otimes C(G'')\ar[r]&M_{N'}(\mathbb C)\otimes M_{N''}(\mathbb C)\ar[d]\\
C(S_N^+)\ar[u]\ar[r]&C(G)\ar[r]&M_N(\mathbb C)
}$$

Here all the maps are the canonical ones, with those on the left and on the right coming from $N=N'N''$. At the level of standard generators, the diagram is as follows:
$$\xymatrix@R=45pt@C65pt{
u_{ij}'\otimes u_{ab}''\ar[r]&w_{ij}'\otimes w_{ab}''\ar[r]&P_{ij}'\otimes P_{ab}''\ar[d]\\
u_{ia,jb}\ar[u]\ar[r]&w_{ia,jb}\ar[r]&P_{ia,jb}
}$$

Now observe that this diagram commutes. We conclude that the representation associated to $H$ factorizes indeed through $C(G')\otimes C(G'')$, and this gives the result.
\end{proof}

At a more abstract level, one interesting question is that of abstractly characterizing the magic matrices coming from the complex Hadamard matrices. We have here:

\begin{proposition}
Given an Hadamard matrix $H\in M_N(\mathbb C)$, the vectors 
$$\xi_{ij}=\frac{H_i}{H_j}$$
on which the magic unitary entries $P_{ij}$ project, have the following properties:
\begin{enumerate}
\item $\xi_{ii}=\xi$ is the all-one vector.

\item $\xi_{ij}\xi_{jk}=\xi_{ik}$, for any $i,j,k$.

\item $\xi_{ij}\xi_{kl}=\xi_{il}\xi_{kj}$, for any $i,j,k,l$.
\end{enumerate}
\end{proposition}

\begin{proof}
All these assertions are trivial, by using the formula $\xi_{ij}=H_i/H_j$.
\end{proof}

Let us call now magic basis of a given Hilbert space $H$ any square array of vectors $\xi\in M_N(H)$, all whose rows and columns are orthogonal bases of $H$. With this convention, the above observations lead to the following result, at the magic basis level:

\index{magic basis}

\begin{theorem}
The magic bases $\xi\in M_N(S^{N-1}_\mathbb C)$ coming from the complex Hadamard matrices are those having the following properties:
\begin{enumerate}
\item We have $\xi_{ij}\in\mathbb T^N$, after a suitable rescaling. 

\item The conditions in Proposition 16.6 are satisfied.
\end{enumerate}
\end{theorem}

\begin{proof}
By using the multiplicativity conditions (1,2,3) in Proposition 16.6, we conclude that, up to a rescaling, we must have $\xi_{ij}=\xi_i/\xi_j$, where $\xi_1,\ldots,\xi_N$ is the first row of the magic basis. Together with our assumption $\xi_{ij}\in\mathbb T^N$, this gives the result. 
\end{proof}

At the general level now, regarding the representation theory of the quantum groups associated to the Hadamard matrices, we have the following result:

\index{Tannakian category}

\begin{theorem}
The Tannakian category of the quantum group $G\subset S_N^+$ associated to a complex Hadamard matrix $H\in M_N(\mathbb C)$ is given by
$$T\in Hom(u^{\otimes k},u^{\otimes l})\iff T^\circ G^{k+2}=G^{l+2}T^\circ$$
where the objects on the right are constructed as follows:
\begin{enumerate}
\item $T^\circ=id\otimes T\otimes id$.

\item $G_{ia}^{jb}=\sum_kH_{ik}\bar{H}_{jk}\bar{H}_{ak}H_{bk}$.

\item $G^k_{i_1\ldots i_k,j_1\ldots j_k}=G_{i_ki_{k-1}}^{j_kj_{k-1}}\ldots G_{i_2i_1}^{j_2j_1}$.
\end{enumerate}
\end{theorem}

\begin{proof}
We use the Tannakian result for the Hopf image of a representation, discussed in chapter 14. With the notations here, we have the following formula:
$$Hom(u^{\otimes k},u^{\otimes l})=Hom(U^{\otimes k},U^{\otimes l})$$

Now recall that the vector space on the right consists by definition of the complex $N^l\times N^k$ matrices $T$ satisfying the following relation:
$$TU^{\otimes k}=U^{\otimes l}T$$ 

If we denote this equality by $L=R$, the left term $L$ is given by:
$$L_{ij}
=(TU^{\otimes k})_{ij}
=\sum_aT_{ia}U^{\otimes k}_{aj}
=\sum_aT_{ia}U_{a_1j_1}\ldots U_{a_kj_k}$$

As for the right term $R$, this is given by the following formula:
$$R_{ij}
=(U^{\otimes l}T)_{ij}
=\sum_bU^{\otimes l}_{ib}T_{bj}
=\sum_bU_{i_1b_1}\ldots U_{i_lb_l}T_{bj}$$

Consider now the vectors $\xi_{ij}=H_i/H_j$. Since these vectors span the ambient Hilbert space, the equality $L=R$ is equivalent to the following equality:
$$<L_{ij}\xi_{pq},\xi_{rs}>=<R_{ij}\xi_{pq},\xi_{rs}>$$

We use now the following well-known formula, expressing a product of rank one projections $P_1,\ldots,P_k$ in terms of the corresponding image vectors $\xi_1,\ldots,\xi_k$:
$$<P_1\ldots P_kx,y>=<x,\xi_k><\xi_k,\xi_{k-1}>\ldots\ldots<\xi_2,\xi_1><\xi_1,y>$$

This gives the following formula for the above left term $L$:
\begin{eqnarray*}
<L_{ij}\xi_{pq},\xi_{rs}>
&=&\sum_aT_{ia}<P_{a_1j_1}\ldots P_{a_kj_k}\xi_{pq},\xi_{rs}>\\
&=&\sum_aT_{ia}<\xi_{pq},\xi_{a_kj_k}>\ldots<\xi_{a_1j_1},\xi_{rs}>\\
&=&\sum_aT_{ia}G_{pa_k}^{qj_k}G_{a_ka_{k-1}}^{j_kj_{k-1}}\ldots G_{a_2a_1}^{j_2j_1}G_{a_1r}^{j_1s}\\
&=&\sum_aT_{ia}G^{k+2}_{rap,sjq}\\
&=&(T^\circ G^{k+2})_{rip,sjq}
\end{eqnarray*}

As for the right term $R$, this is given by the following formula:
\begin{eqnarray*}
<R_{ij}\xi_{pq},\xi_{rs}>
&=&\sum_b<P_{i_1b_1}\ldots P_{i_lb_l}\xi_{pq},\xi_{rs}>T_{bj}\\
&=&\sum_b<\xi_{pq},\xi_{i_lb_l}>\ldots<\xi_{i_1b_1},\xi_{rs}>T_{bj}\\
&=&\sum_bG_{pi_l}^{qb_l}G_{i_li_{l-1}}^{b_lb_{l-1}}\ldots G_{i_2i_1}^{b_2b_1}G_{i_1r}^{b_1s}T_{bj}\\
&=&\sum_bG^{l+2}_{rip,sbq}T_{bj}\\
&=&(G^{l+2}T^\circ)_{rip,sjq}
\end{eqnarray*}

Thus, we obtain the formula in the statement.
\end{proof}

Summarizing, we have a construction of quantum groups from Hadamard matrices, along with a Tannakian result for it, that we can use for various purposes. We will see in what follows that this will lead us into some very interesting mathematics, related to key work of Popa, Haagerup, Jones and others. And of course, into free probability too.

\section*{16b. Spin models}

Let us discuss now the relation with subfactor theory, and with planar algebras. As a starting point, we have the following basic observation of Popa:

\index{orthogonal MASA}

\begin{proposition}
Up to a conjugation by a unitary, the pairs of orthogonal MASA in the simplest factor, namely the matrix algebra $M_N(\mathbb C)$, are as follows,
$$A=\Delta\quad,\quad 
B=H\Delta H^*$$
with $\Delta\subset M_N(\mathbb C)$ being the diagonal matrices, and with $H\in M_N(\mathbb C)$ being Hadamard.
\end{proposition}

\begin{proof}
Any maximal abelian subalgebra (MASA) in $M_N(\mathbb C)$ being conjugated to $\Delta$, we can assume, up to conjugation by a unitary, that we have, with $U\in U_N$:
$$A=\Delta\quad,\quad 
B=U\Delta U^*$$  

Now observe that given two diagonal matrices $D,E\in\Delta$, we have:
\begin{eqnarray*}
tr(D\cdot UEU^*)
&=&\frac{1}{N}\sum_i(DUEU^*)_{ii}\\
&=&\frac{1}{N}\sum_{ij}D_{ii}U_{ij}E_{jj}\bar{U}_{ij}\\
&=&\frac{1}{N}\sum_{ij}D_{ii}E_{jj}|U_{ij}|^2
\end{eqnarray*}

Thus, the orthogonality condition $A\perp B$ reformulates as follows:
$$\frac{1}{N}\sum_{ij}D_{ii}E_{jj}|U_{ij}|^2=\frac{1}{N^2}\sum_{ij}D_{ii}E_{jj}$$

We conclude that we must have $|U_{ij}|=\frac{1}{\sqrt{N}}$, for any $i,j$. But this tells us precisely that the matrix $H=\sqrt{N}U$ must be Hadamard, as claimed.
\end{proof}

Along the same lines, but at a more advanced level, we have:

\index{commuting square}

\begin{theorem}
Given a complex Hadamard matrix $H\in M_N(\mathbb C)$, the diagram formed by the associated pair of orthogonal MASA, namely
$$\xymatrix@R=35pt@C35pt{
\Delta\ar[r]&M_N(\mathbb C)\\
\mathbb C\ar[u]\ar[r]&H\Delta H^*\ar[u] }$$ is a commuting square in the sense of subfactor theory, in the sense that the expectations onto $\Delta,H\Delta H^*$ commute, and their product is the expectation onto $\mathbb C$.
\end{theorem}

\begin{proof}
The expectation $E_\Delta:M_N(\mathbb C)\to\Delta$ is the operation $M\to M_\Delta$ which consists in keeping the diagonal, and erasing the rest. Consider now the other expectation:
$$E_{H\Delta H^*}:M_N(\mathbb C)\to H\Delta H^*$$

It is better to identify this with the following expectation, with $U=H/\sqrt{N}$: 
$$E_{U\Delta U^*}:M_N(\mathbb C)\to U\Delta U^*$$

This latter expectation must be of the form $M\to UX_\Delta U^*$, with $X$ satisfying:
$$<M,UDU^*>=<UX_\Delta U^*,UDU^*>\quad,\quad\forall D\in\Delta$$

The scalar products being given by $<a,b>=tr(ab^*)$, this condition reads:
$$tr(MUD^*U^*)=tr(X_\Delta D^*)\quad,\quad\forall D\in\Delta$$

Thus $X=U^*MU$, and the formulae of our two expectations are as follows:
\begin{eqnarray*}
E_\Delta(M)&=&M_\Delta\\
E_{U\Delta U^*}(M)&=&U(U^*MU)_\Delta U^*
\end{eqnarray*}

With these formulae in hand, an elementary computation gives the result.
\end{proof}

The point now is that any commuting square $C$ produces a subfactor of the Murray-von Neumann hyperfinite ${\rm II}_1$ factor $R$. Consider indeed such a commuting square:
$$\xymatrix@R=35pt@C35pt{
C_{01}\ar[r]&C_{11}\\
C_{00}\ar[u]\ar[r]&C_{10}\ar[u]}$$

What happens is that, under some suitable, mild assumptions on the above inclusions $C_{00}\subset C_{10},C_{01}\subset C_{11}$, we can perform the basic construction for them, in finite dimensions, and we obtain in this way a whole array of commuting squares, as follows:
$$\xymatrix@R=35pt@C35pt{
A_0&A_1&A_2&\\
C_{02}\ar[r]\ar@.[u]&C_{12}\ar[r]\ar@.[u]&C_{22}\ar@.[r]\ar@.[u]&B_2\\
C_{01}\ar[r]\ar[u]&C_{11}\ar[r]\ar[u]&C_{21}\ar@.[r]\ar[u]&B_1\\
C_{00}\ar[u]\ar[r]&C_{10}\ar[u]\ar[r]&C_{20}\ar[u]\ar@.[r]&B_0}$$

Here the various $A,B$ letters stand for the von Neumann algebras obtained in the limit, which are all isomorphic to the hyperfinite ${\rm II}_1$ factor $R$, and we have:

\index{Ocneanu compactness}

\begin{theorem}
In the context of the above diagram, the following happen:
\begin{enumerate}
\item $A_0\subset A_1$ is a subfactor, and $\{A_i\}$ is the Jones tower for it.

\item The corresponding planar algebra is given by $A_0'\cap A_k=C_{01}'\cap C_{k0}$.

\item A similar result holds for the ``horizontal'' subfactor $B_0\subset B_1$.
\end{enumerate}
\end{theorem}

\begin{proof}
This is something very standard, the idea being that (1) is something quite routine, then (2) is a subtle result, called Ocneanu compactness theorem, and then (3) follows from (1,2), by flipping the diagram. Have a look here at any subfactor book.
\end{proof}

Getting back now to the Hadamard matrices, we can extend our lineup of results, namely Proposition 16.9 and Theorem 16.10, with an advanced result, as follows:

\begin{theorem}
Given a complex Hadamard matrix $H\in M_N(\mathbb C)$, the diagram formed by the associated pair of orthogonal MASA is a commuting square in the sense of subfactor theory, and the associated planar algebra $P=(P_k)$ is given by the following formula,
$$T\in P_k\iff T^\circ G^2=G^{k+2}T^\circ$$
where the objects on the right are constructed as follows:
\begin{enumerate}
\item $T^\circ=id\otimes T\otimes id$.

\item $G_{ia}^{jb}=\sum_kH_{ik}\bar{H}_{jk}\bar{H}_{ak}H_{bk}$.

\item $G^k_{i_1\ldots i_k,j_1\ldots j_k}=G_{i_ki_{k-1}}^{j_kj_{k-1}}\ldots G_{i_2i_1}^{j_2j_1}$.
\end{enumerate}
\end{theorem}

\begin{proof}
We have two assertions here, the idea being as follows:

\medskip

(1) The fact that we have indeed a commuting square is something that we already know, from Theorem 16.10, reproduced here for convenience.

\medskip

(2) The computation of the associated planar algebra is possible thanks to the Ocneanu compactness theorem, corresponding to the formula in Theorem 16.11 (2). To be more precise, by doing some direct computations, which are quite similar to those in the proof of Theorem 16.8, we obtain the formula in the statement. See \cite{jo3}.
\end{proof}

The point now is that all the above is very similar to Theorem 16.8, and we have: 

\index{spin model}

\begin{theorem}
Let $H\in M_N(\mathbb C)$ be a complex Hadamard matrix.
\begin{enumerate}
\item The planar algebra associated to $H$ is given by $P_k=Fix(u^{\otimes k})$, where $G\subset S_N^+$ is the associated quantum permutation group.

\item The corresponding Poincar\'e series $f(z)=\sum_k\dim(P_k)z^k$ equals the moment generating function $M(z)=\int_G\frac{1}{1-z\chi}$ of the law of the main character $\chi=\sum_iu_{ii}$.
\end{enumerate}
\end{theorem}

\begin{proof}
This follows by comparing the quantum group and subfactor results:

\medskip

(1) As already mentioned above, this assertion simply follows by comparing Theorem 16.8 with the subfactor computation in Theorem 16.12, yielding the same results. 

\medskip

(2) This is a consequence of (1), and of the Peter-Weyl type results of Woronowicz \cite{wo1}, telling us that fixed points can be counted by integrating characters.
\end{proof}

Summarizing, we have now a clarification of the various quantum algebraic objects associated to a complex Hadamard matrix $H\in M_N(\mathbb C)$, the idea being that the central object, which best encodes the ``symmetries'' of the matrix, and which allows the computation of the other quantum algebraic objects as well, such as the associated planar algebra, is the associated quantum permutation group $G\subset S_N^+$. 

\bigskip

Regarding now the subfactor itself, the result here is as follows:

\begin{theorem}
The subfactor associated to $H\in M_N(\mathbb C)$ is of the form
$$A^G\subset(\mathbb C^N\otimes A)^G$$
with $A=R\rtimes\widehat{G}$, where $G\subset S_N^+$ is the associated quantum permutation group.
\end{theorem}

\begin{proof}
This is something more technical, the idea being that the basic construction procedure for the commuting squares, explained before Theorem 16.11, can be performed in an ``equivariant setting'', for commuting squares having components as follows:
$$D\otimes_GE=(D\otimes(E\rtimes\widehat{G}))^G$$

To be more precise, starting with a commuting square formed by such algebras, we obtain by basic construction a whole array of commuting squares as follows, with $\{D_i\},\{E_i\}$ being by definition Jones towers, and with $D_\infty,E_\infty$ being their inductive limits:
$$\xymatrix@R=35pt@C35pt{
D_0\otimes_GE_\infty&D_1\otimes_GE_\infty&D_2\otimes_GE_\infty\\
D_0\otimes_GE_2\ar@.[u]\ar[r]&D_1\otimes_GE_2\ar@.[u]\ar[r]&D_2\otimes_GE_2\ar@.[u]\ar@.[r]&D_\infty\otimes_GE_2\\
D_0\otimes_GE_1\ar[u]\ar[r]&D_1\otimes_GE_1\ar[u]\ar[r]&D_2\otimes_GE_1\ar[u]\ar@.[r]&D_\infty\otimes_GE_1\\
D_0\otimes_GE_0\ar[u]\ar[r]&D_1\otimes_GE_0\ar[u]\ar[r]&D_2\otimes_GE_0\ar[u]\ar@.[r]&D_\infty\otimes_GE_0}$$

The point now is that this quantum group picture works in fact for any commuting square having $\mathbb C$ in the lower left corner. In the Hadamard matrix case, that we are interested in here, the corresponding commuting square is as follows:
$$\xymatrix@R=35pt@C35pt{
\mathbb C\otimes_G\mathbb C^N\ar[r]&\mathbb C^N\otimes_G\mathbb C^N\\
\mathbb C\otimes_G\mathbb C\ar[u]\ar[r]&\mathbb C^N\otimes_G\mathbb C\ar[u] }$$ 

Thus, the subfactor obtained by vertical basic construction appears as follows:
$$\mathbb C\otimes_GE_\infty\subset\mathbb C^N\otimes_GE_\infty$$

But this gives the conclusion in the statement, with the ${\rm II}_1$ factor appearing there being by definition $A=E_\infty\rtimes\widehat{G}$, and with the remark that we have $E_\infty\simeq R$.
\end{proof}

There are many other things that can be said about Hadamard matrices and their relation with operator algebras, and about the subfactors associated to the Hadamard matrices, and their generalizations, and we refer here to Jones \cite{jo3} and related papers, including my old subfactor papers, where Theorems 16.13 and 16.14 were proved.

\section*{16c. Fourier models}

Getting back now to Theorem 16.5, going beyond it is a quite delicate task. The next simplest models appear by deforming the Fourier matrices, or rather the tensor products of such matrices, $F_{G\times H}=F_G\otimes F_H$, via the following construction, due to Di\c t\u a:

\index{deformed Fourier matrix}

\begin{definition}
Given two finite abelian groups $G,H$, we consider the corresponding deformed Fourier matrix, given by the formula
$$(F_G\otimes_Q F_H)_{ia,jb}=Q_{ib}(F_G)_{ij}(F_H)_{ab}$$
and we factorize the associated representation $\pi_Q$ of the algebra $C(S_{G\times H}^+)$,
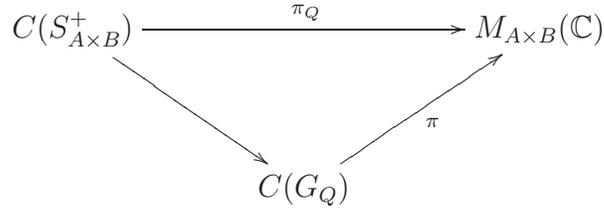
$$\xymatrix@R=38pt@C=40pt
{C(S_{G\times H}^+)\ar[rr]^{\pi_Q}\ar[rd]&&M_{G\times H}(\mathbb C)\\&C(G_Q)\ar[ur]_\pi&}$$
with $C(G_Q)$ being the Hopf image of this representation $\pi_Q$.
\end{definition}

Following \cite{bjb}, we will be interested in computing the above quantum permutation group $G_Q\subset S_{G\times H}^+$, as function of the parameter matrix $Q\in M_{G\times H}(\mathbb T)$. We will need:

\index{free wreath product}
\index{wreath product}

\begin{definition}
Gives two Hopf algebra quotients, as follows,
$$C(S_M^+)\to A\quad,\quad 
C(S_N^+)\to B$$
with fundamental corepresentations denoted $u,v$, we let
$$A*_wB=A^{*N}*B/<[u_{ab}^{(i)},v_{ij}]=0>$$
with the Hopf algebra structure making $w_{ia,jb}=u_{ab}^{(i)}v_{ij}$ a corepresentation.
\end{definition}

Here the fact that we have indeed a Hopf algebra follows from the fact that the above matrix $w$ is magic. In terms of quantum groups, we will write, as usual:
$$C(G)*_wC(H)=C(G\wr_*H)$$

The $\wr_*$ operation is then the free analogue of $\wr$, the usual wreath product. Now with this notion in hand, we can factorize representation $\pi_Q$ in Definition 16.15, as follows:

\begin{proposition}
We have a factorization as follows,
$$\xymatrix@R=38pt@C=35pt
{C(S_{G\times H}^+)\ar[rr]^{\pi_Q}\ar[rd]&&M_{G\times H}(\mathbb C)\\&C(H\wr_*G)\ar[ur]_\pi&}$$
given on the standard generators by the formulae
$$U_{ab}^{(i)}=\sum_jW_{ia,jb}\quad,\quad 
V_{ij}=\sum_aW_{ia,jb}$$
independently of $b$, where $W$ is the magic matrix producing $\pi_Q$.
\end{proposition}

\begin{proof}
With $K=F_G,L=F_H$ and $M=|G|,N=|H|$, the formula of the magic matrix $W\in M_{G\times H}(M_{G\times H}(\mathbb C))$ associated to $H=K\otimes_QL$ is as follows:
\begin{eqnarray*}
(W_{ia,jb})_{kc,ld}
&=&\frac{1}{MN}\cdot\frac{Q_{ic}Q_{jd}}{Q_{id}Q_{jc}}\cdot\frac{K_{ik}K_{jl}}{K_{il}K_{jk}}\cdot\frac{L_{ac}L_{bd}}{L_{ad}L_{bc}}\\
&=&\frac{1}{MN}\cdot\frac{Q_{ic}Q_{jd}}{Q_{id}Q_{jc}}\cdot K_{i-j,k-l}L_{a-b,c-d}
\end{eqnarray*}

Our claim now is that the representation $\pi_Q$ constructed in Definition 16.15 can be factorized in three steps, up to the factorization in the statement, as follows:
$$\xymatrix@R=60pt@C=50pt
{C(S_{G\times H}^+)\ar[rr]^{\pi_Q}\ar[d]&&M_{G\times H}(\mathbb C)\\
C(S_H^+\wr_*S_G^+)\ar[r]\ar@{.>}[rru]&C(S_H^+\wr_*G)\ar[r]\ar@{.>}[ur]&C(H\wr_*G)\ar@{.>}[u]}$$

Indeed, the construction of the map on the left is standard. Regarding the second factorization, this comes from the fact that since the elements $V_{ij}$ depend on $i-j$, they satisfy the defining relations for the quotient algebra $C(S_G^+)\to C(G)$. Finally, regarding the third factorization, observe that $W_{ia,jb}$ depends only on $i,j$ and on $a-b$. By summing over $j$ we obtain that the elements $U_{ab}^{(i)}$ depend only on $a-b$, and we are done.
\end{proof}

In order to further factorize the above representation, we will need:

\index{crossed product}
\index{crossed coproduct}

\begin{definition}
If $H\curvearrowright\Gamma$ is a finite group acting by automorphisms on a discrete group, the corresponding crossed coproduct Hopf algebra is
$$C^*(\Gamma)\rtimes C(H)=C^*(\Gamma)\otimes C(H)$$
with comultiplication given by the following formula,
$$\Delta(r\otimes\delta_k)=\sum_{h\in H}(r\otimes\delta_h)\otimes(h^{-1}\cdot r\otimes\delta_{h^{-1}k})$$
for $r\in\Gamma$ and $k\in H$. The corresponding quantum group is denoted $\widehat{\Gamma}\rtimes H$.
\end{definition}

Observe that $C(H)$ is a subcoalgebra, and that $C^*(\Gamma)$ is not a subcoalgebra. Now back to the factorization in Proposition 16.17, the point is that we have:

\begin{proposition}
With $L=F_H,N=|H|$ we have an isomorphism
$$C(H\wr_*G)\simeq C^*(H)^{*G}\rtimes C(G)$$
given by the following formulae
$$v_{ij}\to1\otimes v_{ij}\quad,\quad 
u_{ab}^{(i)}=\frac{1}{N}\sum_cL_{b-a,c}c^{(i)}\otimes 1$$
on the standard generators.
\end{proposition}

\begin{proof}
We know that $C(H\wr_*G)$ is the quotient of $C(H)^{*G}*C(G)$ by the relations $[u_{ab}^{(i)},v_{ij}]=0$. Now since the variable $v_{ij}$ depends only on $j-i$, we obtain:
$$[u_{ab}^{(i)},v_{kl}]
=[u_{ab}^{(i)},v_{i,l-k+i}]
=0$$

Thus, we are in a usual tensor product situation, and we have:
$$C(H\wr_*G)=C(H)^{*G}\otimes C(G)$$

Consider now the Fourier transform over $H$, which is a map as follows:
$$\Phi:C(H)\to C^*(H)$$

We can compose the above identification with $\Psi=\Phi^{*G}\otimes id$, and we obtain an isomorphism as in the statement. Now observe that we have:
$$\Phi(u_{ab})=\frac{1}{N}\sum_cL_{b-a,c}c$$

Thus the formula for the image of $u_{ab}^{(i)}$ is indeed the one in the statement.
\end{proof}

Here is now our key result, which will lead to further factorizations:

\begin{proposition}
With $c^{(i)}=\sum_aL_{ac}u_{a0}^{(i)}$ and $\varepsilon_{ke}=\sum_iK_{ik}e_{ie}$ we have:
$$\pi(c^{(i)})(\varepsilon_{ke})=\frac{Q_{i,e-c}Q_{i-k,e}}{Q_{ie}Q_{i-k,e-c}}\varepsilon_{k,e-c}$$
In particular if $c_1+\ldots+c_s=0$ then $\pi(c_1^{(i_1)}\ldots c_s^{(i_s)})$ is diagonal, for any $i_1,\ldots,i_s$.
\end{proposition}

\begin{proof}
With $c^{(i)}$ as in the statement, we have the following formula:
$$\pi(c^{(i)})
=\sum_aL_{ac}\pi(u_{a0}^{(i)})
=\sum_{aj}L_{ac}W_{ia,j0}$$

On the other hand, in terms of the basis in the statement, we have:
$$W_{ia,jb}(\varepsilon_{ke})=\frac{1}{N}\delta_{i-j,k}\sum_d\frac{Q_{id}Q_{je}}{Q_{ie}Q_{jd}}L_{a-b,d-e}\varepsilon_{kd}$$

We therefore obtain, as desired:
\begin{eqnarray*}
\pi(c^{(i)})(\varepsilon_{ke})
&=&\frac{1}{N}\sum_{ad}L_{ac}\frac{Q_{id}Q_{i-k,e}}{Q_{ie}Q_{i-k,d}}L_{a,d-e}\varepsilon_{kd}\\
&=&\frac{1}{N}\sum_d\frac{Q_{id}Q_{i-k,e}}{Q_{ie}Q_{i-k,d}}\varepsilon_{kd}\sum_aL_{a,d-e+c}\\
&=&\sum_d\frac{Q_{id}Q_{i-k,e}}{Q_{ie}Q_{i-k,d}}\varepsilon_{kd}\delta_{d,e-c}\\
&=&\frac{Q_{i,e-c}Q_{i-k,e}}{Q_{ie}Q_{i-k,e-c}}\varepsilon_{k,e-c}
\end{eqnarray*}

As for the last assertion, this is clear from this.
\end{proof}

We have now all needed ingredients for refining Proposition 16.17, as follows:

\begin{proposition}
We have a factorization as follows,
$$\xymatrix@R=38pt@C=30pt
{C(S_{G\times H}^+)\ar[rr]^{\pi_Q}\ar[rd]&&M_{G\times H}(\mathbb C)\\&C^*(\Gamma_{G,H})\rtimes C(G)\ar[ur]_\rho&}$$
where the group on the bottom is given by
$$\Gamma_{G,H}=H^{*G}\Big/\left<[c_1^{(i_1)}\ldots c_s^{(i_s)},d_1^{(j_1)}\ldots d_s^{(j_s)}]=1\Big|\sum_rc_r=\sum_rd_r=0\right>$$
with the above conventions and notations.
\end{proposition}

\begin{proof}
Assume that we have a representation, as follows:
$$\pi:C^*(\Gamma)\rtimes C(G)\to M_L(\mathbb C)$$

Let $\Lambda$ be a $G$-stable normal subgroup of $\Gamma$, so that $G$ acts on $\Gamma/\Lambda$, and we can form the product $C^*(\Gamma/\Lambda)\rtimes C(G)$, and assume that $\pi$ is trivial on $\Lambda$. Then $\pi$ factorizes as:
$$\xymatrix@R=40pt@C=30pt
{C^*(\Gamma)\rtimes C(G)\ar[rr]^\pi\ar[rd]&&M_L(\mathbb C)\\&C^*(\Gamma/\Lambda)\rtimes C(G)\ar[ur]_\rho}$$

With $\Gamma=H^{*G}$, and by using the above results, this gives the result.
\end{proof}

We have now all the needed ingredients for proving a main result, as follows:

\begin{theorem}
When $Q$ is generic, the minimal factorization for $\pi_Q$ is 
$$\xymatrix@R=40pt@C=40pt
{C(S_{G\times H}^+)\ar[rr]^{\pi_Q}\ar[rd]&&M_{G\times H}(\mathbb C)\\&C^*(\Gamma_{G,H})\rtimes C(G)\ar[ur]_\pi&}$$
where on the bottom $\Gamma_{G,H}\simeq\mathbb Z^{(|G|-1)(|H|-1)}\rtimes H$ is the discrete group constructed above.
\end{theorem}

\begin{proof}
Consider the factorization constructed in Proposition 16.21, which is as follows, with $L$ being the Hopf image of $\pi_Q$:
$$\theta : C^*(\Gamma_{G,H})\rtimes C(G)\to L$$

To be more precise, this morphism produces the following commutative diagram:
$$\xymatrix@R=30pt@C=40pt
{C(S_{G\times H}^+) \ar[rr]^{\pi_Q} \ar[dr]_{} \ar@/_/[ddr]_{}& & M_{G\times H}(\mathbb C) \\
& L \ar[ur]_{} & \\
& C^*(\Gamma_{G,H})\rtimes C(G) \ar@{-->}[u]_\theta \ar@/_/[uur]_{\pi}& 
}$$

The first observation is that the injectivity assumption on $C(G)$ holds by construction, and that for $f \in C(G)$, the matrix $\pi(f)$ is ``block scalar''. Now for $r \in \Gamma_{G,H}$ with $\theta(r\otimes 1)=\theta(1 \otimes f)$ for some $f \in C(G)$, we see, using the commutative diagram, that $\pi(r \otimes 1)$ is block scalar. Thus, modulo some standard algebra, we are done. See \cite{bjb}.
\end{proof}

Let us compute now the Kesten measure $\mu=law(\chi)$, in the case where the deformation matrix is generic, as before. Our results here will be a combinatorial moment formula, a geometric interpretation of it, and an asymptotic result. We first have:

\begin{theorem}
We have the moment formula
$$\int\chi^p
=\frac{1}{|G|\cdot|H|}\#\left\{\begin{matrix}i_1,\ldots,i_p\in G\\ d_1,\ldots,d_p\in H\end{matrix}\Big|\begin{matrix}[(i_1,d_1),(i_2,d_2),\ldots,(i_p,d_p)]\ \ \ \ \\=[(i_1,d_p),(i_2,d_1),\ldots,(i_p,d_{p-1})]\end{matrix}\right\}$$
where the sets between square brackets are by definition sets with repetition.
\end{theorem}

\begin{proof}
According to the various formulae above, the factorization found in Theorem 16.22 is, at the level of standard generators, as follows:
$$\begin{matrix}
C(S_{G\times H}^+)&\to&C^*(\Gamma_{G,H})\otimes C(G)&\to&M_{G\times H}(\mathbb C)\\
u_{ia,jb}&\to&\frac{1}{|H|}\sum_cF_{b-a,c}c^{(i)}\otimes v_{ij}&\to&W_{ia,jb}
\end{matrix}$$

Thus, the main character of the quantum permutation group that we found in Theorem 16.22 is given by the following formula:
$$\chi
=\frac{1}{|H|}\sum_{iac}c^{(i)}\otimes v_{ii}
=\sum_{ic}c^{(i)}\otimes v_{ii}
=\left(\sum_{ic}c^{(i)}\right)\otimes\delta_1$$

Now since the Haar functional of $C^*(\Gamma)\rtimes C(H)$ is the tensor product of the Haar functionals of $C^*(\Gamma),C(H)$, this gives the following formula, valid for any $p\geq1$:
$$\int\chi^p=\frac{1}{|G|}\int_{\widehat{\Gamma}_{G,H}}\left(\sum_{ic}c^{(i)}\right)^p$$

Consider the elements $S_i=\sum_cc^{(i)}$. With standard notations, we have:
$$S_i=\sum_c(b_{i0}-b_{ic},c)$$

Now observe that these elements multiply as follows:
$$S_{i_1}\ldots S_{i_p}=\sum_{c_1\ldots c_p}
\begin{pmatrix}
b_{i_10}-b_{i_1c_1}+b_{i_2c_1}-b_{i_2,c_1+c_2}&&\\
+b_{i_3,c_1+c_2}-b_{i_3,c_1+c_2+c_3}+\ldots\ldots&,&c_1+\ldots+c_p&\\
\ldots\ldots+b_{i_p,c_1+\ldots+c_{p-1}}-b_{i_p,c_1+\ldots+c_p}&&
\end{pmatrix}$$

In terms of the new indices $d_r=c_1+\ldots+c_r$, this formula becomes:
$$S_{i_1}\ldots S_{i_p}=\sum_{d_1\ldots d_p}
\begin{pmatrix}
b_{i_10}-b_{i_1d_1}+b_{i_2d_1}-b_{i_2d_2}&&\\
+b_{i_3d_2}-b_{i_3d_3}+\ldots\ldots&,&d_p&\\
\ldots\ldots+b_{i_pd_{p-1}}-b_{i_pd_p}&&
\end{pmatrix}$$

Now by integrating, we must have $d_p=0$ on one hand, and on the other hand:
$$[(i_1,0),(i_2,d_1),\ldots,(i_p,d_{p-1})]=[(i_1,d_1),(i_2,d_2),\ldots,(i_p,d_p)]$$

Equivalently, we must have $d_p=0$ on one hand, and on the other hand:
$$[(i_1,d_p),(i_2,d_1),\ldots,(i_p,d_{p-1})]=[(i_1,d_1),(i_2,d_2),\ldots,(i_p,d_p)]$$

Thus, by translation invariance with respect to $d_p$, we obtain:
$$\int_{\widehat{\Gamma}_{G,H}}S_{i_1}\ldots S_{i_p}
=\frac{1}{|H|}\#\left\{d_1,\ldots,d_p\in H\Big|\begin{matrix}[(i_1,d_1),(i_2,d_2),\ldots,(i_p,d_p)]\ \ \ \ \\=[(i_1,d_p),(i_2,d_1),\ldots,(i_p,d_{p-1})]\end{matrix}\right\}$$

It follows that we have the following moment formula:
$$\int_{\widehat{\Gamma}_{G,H}}\left(\sum_iS_i\right)^p
=\frac{1}{|H|}\#\left\{\begin{matrix}i_1,\ldots,i_p\in G\\ d_1,\ldots,d_p\in H\end{matrix}\Big|\begin{matrix}[(i_1,d_1),(i_2,d_2),\ldots,(i_p,d_p)]\ \ \ \ \\=[(i_1,d_p),(i_2,d_1),\ldots,(i_p,d_{p-1})]\end{matrix}\right\}$$

Now by dividing by $|G|$, we obtain the formula in the statement.
\end{proof} 

The formula in Theorem 16.23 can be further interpreted as follows:

\index{Gram matrix}

\begin{theorem}
With $M=|G|,N=|H|$ we have the formula
$$law(\chi)=\left(1-\frac{1}{N}\right)\delta_0+\frac{1}{N}law(A)$$
where $A\in C(\mathbb T^{MN},M_M(\mathbb C))$ is given by $A(q)=$ Gram matrix of the rows of $q$.
\end{theorem}

\begin{proof}
According to Theorem 16.23, we have the following formula:
\begin{eqnarray*}
\int\chi^p
&=&\frac{1}{MN}\sum_{i_1\ldots i_p}\sum_{d_1\ldots d_p}\delta_{[i_1d_1,\ldots,i_pd_p],[i_1d_p,\ldots,i_pd_{p-1}]}\\
&=&\frac{1}{MN}\int_{\mathbb T^{MN}}\sum_{i_1\ldots i_p}\sum_{d_1\ldots d_p}\frac{q_{i_1d_1}\ldots q_{i_pd_p}}{q_{i_1d_p}\ldots q_{i_pd_{p-1}}}\,dq\\
&=&\frac{1}{MN}\int_{\mathbb T^{MN}}\sum_{i_1\ldots i_p}\left(\sum_{d_1}\frac{q_{i_1d_1}}{q_{i_2d_1}}\right)\left(\sum_{d_2}\frac{q_{i_2d_2}}{q_{i_3d_2}}\right)\ldots\left(\sum_{d_p}\frac{q_{i_pd_p}}{q_{i_1d_p}}\right)dq
\end{eqnarray*}

Consider now the Gram matrix in the statement, namely:
$$A(q)_{ij}=<R_i,R_j>$$

To be more precise, $R_1,\ldots,R_M$ are the rows of the following matrix:
$$q\in \mathbb T^{MN}\simeq M_{M\times N}(\mathbb T)$$

We have then the following computation:
\begin{eqnarray*}
\int\chi^p
&=&\frac{1}{MN}\int_{\mathbb T^{MN}}<R_{i_1},R_{i_2}><R_{i_2},R_{i_3}>\ldots<R_{i_p},R_{i_1}>\\
&=&\frac{1}{MN}\int_{\mathbb T^{MN}}A(q)_{i_1i_2}A(q)_{i_2i_3}\ldots A(q)_{i_pi_1}\\
&=&\frac{1}{MN}\int_{\mathbb T^{MN}}Tr(A(q)^p)dq\\
&=&\frac{1}{N}\int_{\mathbb T^{MN}}tr(A(q)^p)dq
\end{eqnarray*}

But this gives the formula in the statement, and we are done.
\end{proof}

In general, the moments of the Gram matrix $A$ are given by a quite complicated formula, and we cannot expect to have a refinement of Theorem 16.24, with $A$ replaced by a plain, non-matricial random variable, say over a compact abelian group. However, this kind of interesting simplification does appear at $M=2$. We first have:

\begin{proposition}
For $F_2\otimes_QF_H$, with $Q\in M_{2\times N}(\mathbb T)$ generic, we have
$$N\int\left(\frac{\chi}{N}\right)^p=\int_{\mathbb T^N}\sum_{k\geq0}\binom{p}{2k}\left|\frac{a_1+\ldots+a_N}{N}\right|^{2k}da$$
where the integral on the right is with respect to the uniform measure on $\mathbb T^N$.
\end{proposition}

\begin{proof}
Consider the following quantity, from the proof of Theorem 16.24:
$$\Phi(q)=\sum_{i_1\ldots i_p}\sum_{d_1\ldots d_p}\frac{q_{i_1d_1}\ldots q_{i_pd_p}}{q_{i_1d_p}\ldots q_{i_pd_{p-1}}}$$

We can ``half-dephase'' the matrix $q\in M_{2\times N}(\mathbb T)$ if we want to, as follows:
$$q=\begin{pmatrix}1&\ldots&1\\ a_1&\ldots&a_N\end{pmatrix}$$

Let us compute now the above quantity $\Phi(q)$, in terms of the numbers $a_1,\ldots,a_N$. Our claim is that we have the following formula:
$$\Phi(q)=2\sum_{k\geq0}N^{p-2k}\binom{p}{2k}\left|\sum_ia_i\right|^{2k}$$

Indeed, in order to find the $N^{p-2k}|\sum_ia_i|^{2k}$ contribution, we have to count the circular configurations consisting of $p$ numbers $1,2$, such that the $1$ values are arranged into $k$ non-empty intervals, and the $2$ values are arranged into $k$ non-empty intervals as well. Now by looking at the endpoints of these $2k$ intervals, we have $2\binom{p}{2k}$ choices, and this gives the above formula. And by integrating, this gives the formula in the statement.
\end{proof}

Observe now that the integrals in Proposition 16.25 can be computed as follows:
\begin{eqnarray*}
\int_{\mathbb T^N}|a_1+\ldots+a_N|^{2k}da
&=&\int_{\mathbb T^N}\sum_{i_1\ldots i_k}\sum_{j_1\ldots j_k}\frac{a_{i_1}\ldots a_{i_k}}{a_{j_1}\ldots a_{j_k}}da\\
&=&\#\left\{i_1\ldots i_k,j_1\ldots j_k\Big|[i_1,\ldots,i_k]=[j_1,\ldots,j_k]\right\}\\
&=&\sum_{k=\sum r_i}\binom{k}{r_1,\ldots,r_N}^2
\end{eqnarray*}

We obtain in this way the following ``blowup'' result, for our measure:

\index{blowup}

\begin{proposition}
For $F_2\otimes_QF_H$, with $Q\in M_{2\times N}(\mathbb T)$ generic, we have
$$\mu=\left(1-\frac{1}{N}\right)\delta_0+\frac{1}{2N}\left(\Psi^+_*\varepsilon+\Psi^-_*\varepsilon\right)$$
where $\varepsilon$ is the uniform measure on $\mathbb T^N$, and where $\Psi^\pm(a)=N\pm\left|\sum_ia_i\right|$.
\end{proposition}

\begin{proof}
We use the formula found in Proposition 16.25, along with:
$$\sum_{k\geq0}\binom{p}{2k}x^{2k}=\frac{(1+x)^p+(1-x)^p}{2}$$

By using this identity, Proposition 16.25 reformulates as follows:
$$N\int\left(\frac{\chi}{N}\right)^p=\frac{1}{2}\int_{\mathbb T^N}\left(1+\left|\frac{\sum_ia_i}{N}\right|\right)^p+\left(1-\left|\frac{\sum_ia_i}{N}\right|\right)^p\,da$$

Now by multiplying by $N^{p-1}$, we obtain the following formula:
$$\int\chi^k=\frac{1}{2N}\int_{\mathbb T^N}\left(N+\left|\sum_ia_i\right|\right)^p+\left(N-\left|\sum_ia_i\right|\right)^p\,da$$

But this gives the formula in the statement, and we are done.
\end{proof}

We can further improve the above result, by reducing the maps $\Psi^\pm$ appearing there to a single one, and we are led to the following statement: 

\begin{theorem}
For $F_2\otimes_QF_H$, with $Q\in M_{2\times N}(\mathbb T)$ generic, we have
$$\mu=\left(1-\frac{1}{N}\right)\delta_0+\frac{1}{N}\Phi_*\varepsilon$$
where $\varepsilon$ is the uniform measure on $\mathbb Z_2\times\mathbb T^N$, and where $\Phi(e,a)=N+e\left|\sum_ia_i\right|$.
\end{theorem}

\begin{proof}
This is clear indeed from Proposition 16.26.
\end{proof}

As already mentioned, the above results at $M=2$ are something quite special. In the general case, $M\in\mathbb N$, it is not clear how to construct a nice blowup of the measure.

\section*{16d. Poisson laws}

Asymptotically, things are however quite simple, and we would like to explain this now. Let us go back indeed to the general case, where $M,N\in\mathbb N$ are arbitrary. The problem that we would like to solve is that of finding the good regime, of the following type, where the law in Theorem 16.23 converges, after some suitable manipulations:
$$M=f(K)\quad,\quad 
N=g(K)\quad,\quad 
K\to\infty$$

In order to do so, we have to do some combinatorics. Let $NC(p)$ be the set of noncrossing partitions of $\{1,\ldots,p\}$, and for $\pi\in P(p)$ let us denote by $|\pi|\in\{1,\ldots,p\}$ the number of blocks. With these conventions, we have the following result:

\begin{proposition}
With $M=\alpha K,N=\beta K$, $K\to\infty$ we have:
$$\frac{c_p}{K^{p-1}}\simeq\sum_{r=1}^p\#\left\{\pi\in NC(p)\Big||\pi|=r\right\}\alpha^{r-1}\beta^{p-r}$$
In particular, with $\alpha=\beta$ we have $c_p\simeq\frac{1}{p+1}\binom{2p}{p}(\alpha K)^{p-1}$.
\end{proposition}

\begin{proof}
We use the combinatorial formula in Theorem 16.23. Our claim is that, with $\pi=\ker(i_1,\ldots,i_p)$, the corresponding contribution to $c_p$ is:
$$C_\pi\simeq
\begin{cases}
\alpha^{|\pi|-1}\beta^{p-|\pi|}K^{p-1}&{\rm if}\ \pi\in NC(p)\\
O(K^{p-2})&{\rm if}\ \pi\notin NC(p)
\end{cases}$$

As a first observation, the number of choices for a multi-index $(i_1,\ldots,i_p)\in X^p$ satisfying $\ker i=\pi$ is given by the following formula:
$$M(M-1)\ldots (M-|\pi|+1)\simeq M^{|\pi|}$$

We conclude that we have the following estimate:
$$C_\pi\simeq M^{|\pi|-1}N^{-1}\#\left\{d_1,\ldots,d_p\in Y\Big|[d_\alpha|\alpha\in b]=[d_{\alpha-1}|\alpha\in b],\forall b\in\pi\right\}$$

The contribution of $\sigma=\ker d$ to the above quantity $C_\pi$ is then given by:
$$\Delta(\pi,\sigma)N(N-1)\ldots(N-|\sigma|+1)\simeq\Delta(\pi,\sigma)N^{|\sigma|}$$

To be more precise, the quantities on the right are as follows:
$$\Delta(\pi,\sigma)=\begin{cases}
1&{\rm if}\ |b\cap c|=|(b-1)\cap c|,\forall b\in\pi,\forall c\in\sigma\\
0&{\rm otherwise}
\end{cases}$$

We use now the standard fact that for $\pi,\sigma\in P(p)$ satisfying $\Delta(\pi,\sigma)=1$ we have:
$$|\pi|+|\sigma|\leq p+1$$

In addition, the equality case is well-known to happen when $\pi,\sigma\in NC(p)$ are inverse to each other, via Kreweras complementation. This shows that for $\pi\notin NC(p)$ we have:
$$C_\pi=O(K^{p-2})$$

On the other hand, this shows as well that for $\pi\in NC(p)$ we have:
$$C_\pi
\simeq M^{|\pi|-1}N^{-1}N^{p-|\pi|-1}
=\alpha^{|\pi|-1}\beta^{p-|\pi|}K^{p-1}$$

Thus, we are led to the connclusion in the statement.
\end{proof}

We denote by $D$ the dilation operation for measures, given by:
$$D_r(law(X))=law(rX)$$

With this convention, still following \cite{bjb}, we have the following result:

\begin{theorem}
With $M=\alpha K,N=\beta K$, $K\to\infty$ we have
$$\mu=\left(1-\frac{1}{\alpha\beta K^2}\right)\delta_0+\frac{1}{\alpha\beta K^2}D_{\frac{1}{\beta K}}(\pi_{\alpha/\beta})$$
and in particular with $\alpha=\beta$ we have the formula
$$\mu=\left(1-\frac{1}{\alpha^2K^2}\right)\delta_0+\frac{1}{\alpha^2K^2}D_{\frac{1}{\alpha K}}(\pi_1)$$
with $D$ being the dilation operation for measures, $D_r(law(X))=law(rX)$.
\end{theorem}

\begin{proof}
At $\alpha=\beta$, this follows from Proposition 16.28. In general now, we have:
\begin{eqnarray*}
\frac{c_p}{K^{p-1}}
&\simeq&\sum_{\pi\in NC(p)}\alpha^{|\pi|-1}\beta^{p-|\pi|}\\
&=&\frac{\beta^p}{\alpha}\sum_{\pi\in NC(p)}\left(\frac{\alpha}{\beta}\right)^{|\pi|}\\
&=&\frac{\beta^p}{\alpha}\int x^pd\pi_{\alpha/\beta}(x)
\end{eqnarray*}

When $\alpha\geq\beta$, where $d\pi_{\alpha/\beta}(x)=\varphi_{\alpha/\beta}(x)dx$ is continuous, we obtain:
$$c_p
=\frac{1}{\alpha K}\int(\beta Kx)^p\varphi_{\alpha/\beta}(x)dx
=\frac{1}{\alpha\beta K^2}\int x^p\varphi_{\alpha/\beta}\left(\frac{x}{\beta K}\right)dx$$

But this gives the formula in the statement. When $\alpha\leq\beta$ the computation is similar, with a Dirac mass as 0 dissapearing and reappearing, and gives the same result.
\end{proof}

Let us state as well an explicit result, regarding densities:

\begin{theorem}
With $M=\alpha K,N=\beta K$, $K\to\infty$ we have
$$\mu=\left(1-\frac{1}{\alpha\beta K^2}\right)\delta_0+\frac{1}{\alpha\beta K^2}\cdot\frac{\sqrt{4\alpha\beta K^2-(x-\alpha K-\beta K)^2}}{2\pi x}\,dx$$
and in particular with $\alpha=\beta$ we have the formula
$$\mu=\left(1-\frac{1}{\alpha^2K^2}\right)\delta_0+\frac{1}{\alpha^2K^2}\cdot\frac{\sqrt{\frac{4\alpha K}{x}-1}}{2\pi}$$
which at $\alpha=\beta=1$ gives $\mu=\left(1-\frac{1}{K^2}\right)\delta_0+\frac{1}{K^2}D_{\frac{1}{K}}(\pi_1)$.
\end{theorem}

\begin{proof}
According to the formula for the density of the free Poisson law, the density of the continuous part $D_{\frac{1}{\beta K}}(\pi_{\alpha/\beta})$ is indeed given by:
$$\frac{\sqrt{4\frac{\alpha}{\beta}-(\frac{x}{\beta K}-1-\frac{\alpha}{\beta})^2}}
{2\pi\cdot\frac{x}{\beta K}}=\frac{\sqrt{4\alpha\beta K^2-(x-\alpha K-\beta K)^2}}{2\pi x}$$

With $\alpha=\beta$ now, we obtain the second formula in the statement, and we are done.
\end{proof}

Observe that at $\alpha=\beta=1$ the above measure is supported by $[0,4K]$. On the other hand, since the groups $\Gamma_{M,N}$ are all amenable, the corresponding measures are supported on $[0,MN]$, and so on $[0,K^2]$ in the $M=N=K$ situation. The fact that we do not have a convergence of supports is not surprising, because our convergence is in moments.

\bigskip

The above results are of course not the end of the story, because we have now to understand what happens in the case of non-generic parameters. There has been some technical work here, by Bichon and by myself, and as a sample result, we have:

\index{free Poisson law}
\index{Marchenko-Pastur law}

\begin{theorem}
Given two finite abelian groups $G,H$, having cardinalities 
$$|G|=M\quad,\quad 
|H|=N$$ 
consider the main character $\chi$ of the quantum group associated to $\mathcal F_{G\times H}$. We have then
$$law\left(\frac{\chi}{N}\right)=\left(1-\frac{1}{M}\right)\delta_0+\frac{1}{M}\,\pi_t$$
in moments, with $M=tN\to\infty$, where $\pi_t$ is the free Poisson law of parameter $t>0$. In addition, this formula holds for any generic fiber of $\mathcal F_{G\times H}$.
\end{theorem}

\begin{proof}
We already know that the second assertion holds, as explained above. Regarding now the first assertion, our first claim is that for the representation coming from the parametric matrix $\mathcal F_{G\times H}$ we have the following formula, where $M=|G|,N=|H|$, and the sets between brackets are sets with repetitions:
$$c_p^r=\frac{1}{M^{r+1}N}\#\left\{\begin{matrix}i_1,\ldots,i_r,a_1,\ldots,a_p\in\{0,\ldots,M-1\},\\
b_1,\ldots,b_p\in\{0,\ldots,N-1\},\\
[(i_x+a_y,b_y),(i_{x+1}+a_y,b_{y+1})|y=1,\ldots,p]\\
=[(i_x+a_y,b_{y+1}),(i_{x+1}+a_y,b_y)|y=1,\ldots,p], \forall x
\end{matrix}\right\}$$

Indeed, by using the general moment formula with $K=F_G$, $L=F_H$, we have:
\begin{eqnarray*}
c_p^r
&=&\frac{1}{(MN)^r}\int_{T^r}\sum_{i_1^1\ldots i_p^r}\sum_{b_1^1\ldots b_p^r}\frac{Q^1_{i_1^1b_1^1}Q^1_{i_1^2b_2^1}}{Q^1_{i_1^1b_2^1}Q^1_{i_1^2b_1^1}}\ldots\frac{Q^1_{i_p^1b_p^1}Q^1_{i_p^2b_1^1}}{Q^1_{i_p^1b_1^1}Q^1_{i_p^2b_p^1}}\ldots\ldots\frac{Q^r_{i_1^rb_1^r}Q^r_{i_1^1b_2^r}}{Q^r_{i_1^rb_2^r}Q^r_{i_1^1b_1^r}}\ldots\frac{Q^r_{i_p^rb_p^r}Q^r_{i_p^1b_1^r}}{Q^r_{i_p^rb_1^r}Q^r_{i_p^1b_p^r}}\\
&&\hskip15mm\frac{1}{M^{pr}}\sum_{j_1^1\ldots j_p^r}\frac{K_{i_1^1j_1^1}K_{i_1^2j_2^1}}{K_{i_1^1j_2^1}K_{i_1^2j_1^1}}\ldots\frac{K_{i_p^1j_p^1}K_{i_p^2j_1^1}}{K_{i_p^1j_1^1}K_{i_p^2j_p^1}}\ldots\ldots\frac{K_{i_1^rj_1^r}K_{i_1^1j_2^r}}{K_{i_1^rj_2^r}K_{i_1^1j_1^r}}\ldots\frac{K_{i_p^rj_p^r}K_{i_p^1j_1^r}}{K_{i_p^rj_1^r}K_{i_p^1j_p^r}}\\
&&\hskip15mm\frac{1}{N^{pr}}\sum_{a_1^1\ldots a_p^r}\frac{L_{a_1^1b_1^1}L_{a_1^2b_2^1}}{L_{a_1^1b_2^1}L_{a_1^2b_1^1}}\ldots\frac{L_{a_p^1b_p^1}L_{a_p^2b_1^1}}{L_{a_p^1b_1^1}L_{a_p^2b_p^1}}\ldots\ldots\frac{L_{a_1^rb_1^r}L_{a_1^1b_2^r}}{L_{a_1^rb_2^r}L_{a_1^1b_1^r}}\ldots\frac{L_{a_p^rb_p^r}L_{a_p^1b_1^r}}{L_{a_p^rb_1^r}L_{a_p^1b_p^r}}\,dQ\\
\end{eqnarray*}

Since we are in the Fourier matrix case, $K=F_G,L=F_H$, we can perform the sums over $j,a$. To be more precise, the last two averages appearing above are respectively:
\begin{eqnarray*}
\Delta(i)&=&\prod_x\prod_y\delta(i^x_y+i^{x+1}_{y-1},i^{x+1}_y+i^x_{y-1})\\
\Delta(b)&=&\prod_x\prod_y\delta(b^x_y+b^{x+1}_{y-1},b^{x+1}_y+b^x_{y-1})
\end{eqnarray*}

We therefore obtain the following formula for the truncated moments of the main character, where $\Delta$ is the product of Kronecker symbols constructed above:
\begin{eqnarray*}
c_p^r
&=&\frac{1}{(MN)^r}\int_{T^r}\sum_{\Delta(i)=\Delta(b)=1}\frac{Q^1_{i_1^1b_1^1}Q^1_{i_1^2b_2^1}}{Q^1_{i_1^1b_2^1}Q^1_{i_1^2b_1^1}}\ldots\frac{Q^1_{i_p^1b_p^1}Q^1_{i_p^2b_1^1}}{Q^1_{i_p^1b_1^1}Q^1_{i_p^2b_p^1}}\ldots\ldots\frac{Q^r_{i_1^rb_1^r}Q^r_{i_1^1b_2^r}}{Q^r_{i_1^rb_2^r}Q^r_{i_1^1b_1^r}}\ldots\frac{Q^r_{i_p^rb_p^r}Q^r_{i_p^1b_1^r}}{Q^r_{i_p^rb_1^r}Q^r_{i_p^1b_p^r}}\,dQ
\end{eqnarray*}

Now by integrating with respect to $Q\in(\mathbb T^{G\times H})^r$, we are led to counting the multi-indices $i,b$ satisfying several conditions. First, we have the following condition:
$$\Delta(i)=\Delta(b)=1$$

We have as well the following conditions, where the sets between brackets are by definition sets with repetitions:
$$\begin{bmatrix}
i_1^1b_1^1&\ldots&i_p^1b_p^1&i_1^2b_2^1&\ldots&i_p^2b_1^1
\end{bmatrix}=
\begin{bmatrix}
i_1^1b_2^1&\ldots&i_p^1b_1^1&i_1^2b_1^1&\ldots&i_p^2b_p^1
\end{bmatrix}$$
$$\vdots$$
$$\begin{bmatrix}
i_1^rb_1^r&\ldots&i_p^rb_p^r&i_1^1b_2^r&\ldots&i_p^1b_1^r
\end{bmatrix}
=\begin{bmatrix}
i_1^rb_2^r&\ldots&i_p^rb_1^r&i_1^1b_1^r&\ldots&i_p^1b_p^r
\end{bmatrix}$$

In a more compact notation, the moment formula that we obtain in this way is therefore as follows:
$$c_p^r=\frac{1}{(MN)^r}\#\left\{i,b\Big|\Delta(i)=\Delta(b)=1,\ [i^x_yb^x_y,i^{x+1}_yb^x_{y+1}]=[i^x_yb^x_{y+1},i^{x+1}_yb^x_y],\forall x\right\}$$

Now observe that the above Kronecker type conditions $\Delta(i)=\Delta(b)=1$ tell us that the arrays of indices $i=(i^x_y),b=(b^x_y)$ must be of the following special form:
$$\begin{pmatrix}i^1_1&\ldots&i^1_p\\&\ldots\\ i^1_r&\ldots&i^r_p\end{pmatrix}=\begin{pmatrix}i_1+a_1&\ldots&i_1+a_p\\&\ldots\\ i_r+a_1&\ldots&i_r+a_p\end{pmatrix}$$

$$\begin{pmatrix}b^1_1&\ldots&b^1_p\\&\ldots\\ b^1_r&\ldots&b^r_p\end{pmatrix}=\begin{pmatrix}j_1+b_1&\ldots&j_1+b_p\\&\ldots\\ j_r+b_1&\ldots&j_r+b_p\end{pmatrix}$$

Here all the new indices $i_x,j_x,a_y,b_y$ are uniquely determined, up to a choice of $i_1,j_1$. Now by replacing $i^x_y,b^x_y$ with these new indices $i_x,j_x,a_y,b_y$, with a $MN$ factor added, which accounts for the choice of $i_1,j_1$, we obtain the following formula:
$$c_p^r=\frac{1}{(MN)^{r+1}}\#\left\{i,j,a,b\Big|\begin{matrix}[(i_x+a_y,j_x+b_y),(i_{x+1}+a_y,j_x+b_{y+1})]\\
=[(i_x+a_y,j_x+b_{y+1}),(i_{x+1}+a_y,j_x+b_y)],\forall x\end{matrix}\right\}$$

Now observe that we can delete if we want the $j_x$ indices, which are irrelevant. Thus, we obtain the announced formula. The continuation is via standard combinatorics.
\end{proof}

There are many interesting questions that are still open, regarding the computation of the spectral measure in the case where the parameter matrix $Q$ is not generic, and also regarding the computation for the deformations of the generalized Fourier matrices, which are not necessarily of Di\c t\u a type. We refer here to \cite{bjb} and related papers.

\section*{16e. Exercises}

Congratulations for having read this book, and no exercises for this final chapter. Instead, you can learn more free probability, from the various books and articles referenced below, and then have a look at various open problems, from here or elsewhere.

\baselineskip=14pt

\printindex


\begin{thebibliography}{99}

\baselineskip=11.4pt

\bibitem{afa}G.W. Anderson and B. Farrell, Asymptotically liberating sequences of random unitary matrices, {\em Adv. Math.} {\bf 255} (2014), 381--413.

\bibitem{agz}G.W. Anderson, A. Guionnet and O. Zeitouni, An introduction to random matrices, Cambridge Univ. Press (2010).

\bibitem{ans}M. Anshelevich, Free Meixner states, {\em Comm. Math. Phys.} {\bf 276} (2007), 863--899.

\bibitem{anv}O. Arizmendi, I. Nechita and C. Vargas, On the asymptotic distribution of block-modified random matrices, {\em J. Math. Phys.} {\bf 57} (2016), 1--27.

\bibitem{awi}R.A. Askey and J. Wilson, Some basic hypergeometric orthogonal polynomials that generalize Jacobi polynomials, AMS (1985).

\bibitem{aub}G. Aubrun, Partial transposition of random states and non-centered semicircular distributions, {\em Random Matrices Theory Appl.} {\bf 1} (2012), 125--145.

\bibitem{ba1}T. Banica, On the polar decomposition of circular variables, {\em Integral Equations Operator Theory} {\bf 24} (1996), 372--377. 

\bibitem{ba2}T. Banica, Normal random variables (2026).

\bibitem{bb+}T. Banica, S.T. Belinschi, M. Capitaine and B. Collins, Free Bessel laws, {\em Canad. J. Math.} {\bf 63} (2011), 3--37.

\bibitem{bjb}T. Banica and J. Bichon, Random walk questions for linear quantum groups, {\em Int. Math. Res. Not.} {\bf 24} (2015), 13406--13436.

\bibitem{bbc}T. Banica, J. Bichon and S. Curran, Quantum automorphisms of twisted group algebras and free hypergeometric laws, {\em Proc. Amer. Math. Soc.} {\bf 139} (2011), 3961--3971.

\bibitem{bbi}T. Banica and D. Bisch, Spectral measures of small index principal graphs, {\em Comm. Math. Phys.} {\bf 269} (2007), 259--281.

\bibitem{bco}T. Banica and B. Collins, Integration over the Pauli quantum group, {\em J. Geom. Phys.} {\bf 58} (2008), 942--961.

\bibitem{bcz}T. Banica, B. Collins and P. Zinn-Justin, Spectral analysis of the free orthogonal matrix, {\em Int. Math. Res. Not.} {\bf 17} (2009), 3286--3309.

\bibitem{bcu}T. Banica and S. Curran, Decomposition results for Gram matrix determinants, {\em J. Math. Phys.} {\bf 51} (2010), 1--14.

\bibitem{bcs}T. Banica, S. Curran and R. Speicher, De Finetti theorems for easy quantum groups, {\em Ann. Probab.} {\bf 40} (2012), 401--435.

\bibitem{bgo}T. Banica and D. Goswami, Quantum isometries and noncommutative spheres, {\em Comm. Math. Phys.} {\bf 298} (2010), 343--356.

\bibitem{bn1}T. Banica and I. Nechita, Asymptotic eigenvalue distributions of block-transposed Wishart matrices, {\em J. Theoret. Probab.} {\bf 26} (2013), 855--869.

\bibitem{bn2}T. Banica and I. Nechita, Flat matrix models for quantum permutation groups, {\em Adv. Appl. Math.} {\bf 83} (2017), 24--46.

\bibitem{bsp}T. Banica and R. Speicher, Liberation of orthogonal Lie groups, {\em Adv. Math.} {\bf 222} (2009), 1461--1501.

\bibitem{bve}T. Banica and R. Vergnioux, Fusion rules for quantum reflection groups, {\em J. Noncommut. Geom.} {\bf 3} (2009), 327--359.

\bibitem{bax}R.J. Baxter, Exactly solved models in statistical mechanics, Academic Press (1982).

\bibitem{bbl}S.T. Belinschi, M. Bo\.zejko, F. Lehner and R. Speicher, The normal distribution is $\boxplus$-infinitely divisible, {\em Adv. Math.} {\bf 226} (2011), 3677--3698.

\bibitem{bpa}H. Bercovici and V. Pata, Stable laws and domains of attraction in free probability theory, {\em Ann. of Math.} {\bf 149} (1999), 1023--1060.

\bibitem{bvo}H. Bercovici and D.V. Voiculescu, Regularity questions for free convolution, {\em Oper. Theory Adv. Appl.} {\bf 104} (1998), 37--47.

\bibitem{bdd}J. Bhowmick, F. D'Andrea and L. Dabrowski, Quantum isometries of the finite noncommutative geometry of the standard model, {\em Comm. Math. Phys.} {\bf 307} (2011), 101--131.

\bibitem{bia}P. Biane, Some properties of crossings and partitions, {\em Discrete Math.} {\bf 175} (1997), 41--53.

\bibitem{bcg}P. Biane, M. Capitaine and A. Guionnet, Large deviation bounds for matrix Brownian motion, {\em Invent. Math.} {\bf 152} (2003), 433--459.

\bibitem{bdu}J. Bichon and M. Dubois-Violette, Half-commutative orthogonal Hopf algebras, {\em Pacific J. Math.} {\bf 263} (2013), 13--28.

\bibitem{bjo}D. Bisch and V.F.R. Jones, Algebras associated to intermediate subfactors, {\em Invent. Math.} {\bf 128} (1997), 89--157.

\bibitem{bos}A. Bose, Random matrices and non-commutative probability, CRC Press (2021).

\bibitem{brs}M. Bo\.zejko and R. Speicher, An example of a generalized Brownian motion, {\em Comm. Math. Phys.} {\bf 137} (1991), 519--531.

\bibitem{bwy}M. Bo\.zejko and J. Wysocza\'nski, Remarks on $t$-transformations of measures and convolutions, {\em Ann. Inst. H. Poincar\'e Probab. Statist.} {\bf 37} (2001), 737--776.

\bibitem{cco}A.H. Chamseddine and A. Connes, The spectral action principle, {\em Comm. Math. Phys.} {\bf 186} (1997), 731--750.

\bibitem{col}B. Collins, Moments and cumulants of polynomial random variables on unitary groups, the Itzykson-Zuber integral, and free probability, {\em Int. Math. Res. Not.} {\bf 17} (2003), 953--982.

\bibitem{cma}B. Collins and S. Matsumoto, On some properties of orthogonal Weingarten functions, {\em J. Math. Phys.} {\bf 50} (2009), 1--18.

\bibitem{cne}B. Collins and I. Nechita, Random quantum channels I: graphical calculus and the Bell state phenomenon, {\em Comm. Math. Phys.} {\bf 297} (2010), 345--370.

\bibitem{csn}B. Collins and P. \'Sniady, Integration with respect to the Haar measure on unitary, orthogonal and symplectic groups, {\em Comm. Math. Phys.} {\bf 264} (2006), 773--795.

\bibitem{co1}A. Connes, Une classification des facteurs de type ${\rm III}$, {\em Ann. Sci. Ec. Norm. Sup.} {\bf 6} (1973), 133--252.

\bibitem{co2}A. Connes, Classification of injective factors. Cases ${\rm II}_1$, ${\rm II}_\infty$, ${\rm III}_\lambda$, $\lambda\neq1$, {\em Ann. of Math.} {\bf 104} (1976), 73--115. 

\bibitem{co3}A. Connes, Noncommutative geometry, Academic Press (1994).

\bibitem{csp}S. Curran and R. Speicher, Quantum invariant families of matrices in free probability, {\em J. Funct. Anal.} {\bf 261} (2011), 897--933.

\bibitem{dif}P. Di Francesco, Meander determinants, {\em Comm. Math. Phys.} {\bf 191} (1998), 543--583.

\bibitem{dms}P. Di Francesco, P. Mathieu and D. S\'en\'echal, Conformal field theory, Springer (1996).

\bibitem{ded}I. Dumitriu and A. Edelman, Matrix models for beta ensembles, {\em J. Math. Phys.} {\bf 43} (2002), 5830--5847.

\bibitem{dur}R. Durrett, Probability: theory and examples, Cambridge Univ. Press (1990).

\bibitem{dyk}K. Dykema, Free products of hyperfinite von Neumann algebras and free dimension, {\em Duke Math. J.} {\bf 69} (1993), 97--119.

\bibitem{epu}D.E. Evans and M. Pugh, Spectral measures and generating series for nimrep graphs in subfactor theory, {\em Comm. Math. Phys.} {\bf 295} (2010), 363--413.

\bibitem{fel}W. Feller, An introduction to probability theory and its applications, Wiley (1950).

\bibitem{fni}M. F\'evrier and A. Nica, Infinitesimal non-crossing cumulants and free probability of type B, {\em J. Funct. Anal.} {\bf 258} (2010), 2983--3023.

\bibitem{fsn}M. Fukuda and P. \'Sniady, Partial transpose of random quantum states: exact formulas and meanders, {\em J. Math. Phys.} {\bf 54} (2013), 1--31.

\bibitem{glm}P. Graczyk, G. Letac and H. Massam, The complex Wishart distribution and the symmetric group, {\em Ann. Statist.} {\bf 31} (2003), 287--309.

\bibitem{gjs}A. Guionnet, V.F.R. Jones and D. Shlyakhtenko, Random matrices, free probability, planar algebras and subfactors, {\em Quanta of maths} {\bf 11} (2010), 201--239.

\bibitem{gkz}A. Guionnet, M. Krishnapur and O. Zeitouni, The single ring theorem, {\em Ann. of Math.} {\bf 174} (2011), 1189--1217.

\bibitem{haa}U. Haagerup, On Voiculescu's R and S transforms for free non-commuting random variables, {\em Fields Inst. Comm.} {\bf 12} (1997), 127--148.

\bibitem{hth}U. Haagerup and S. Thorbj\o rnsen, Random matrices with complex Gaussian entries, {\em Exposition. Math.} {\bf 21} (2003), 293--337.

\bibitem{hpe}F. Hiai and D. Petz, The semicircle law, free random variables and entropy, AMS (2000).

\bibitem{jo1}V.F.R. Jones, Index for subfactors, {\em Invent. Math.} {\bf 72} (1983), 1--25.

\bibitem{jo2}V.F.R. Jones, On knot invariants related to some statistical mechanical models, {\em Pacific J. Math.} {\bf 137} (1989), 311--334.

\bibitem{jo3}V.F.R. Jones, Planar algebras I (1999).

\bibitem{jo4}V.F.R. Jones, The planar algebra of a bipartite graph, in ``Knots in Hellas '98'' (2000), 94--117.

\bibitem{jo5}V.F.R. Jones, The annular structure of subfactors, {\em Monogr. Enseign. Math.} {\bf 38} (2001), 401--463.

\bibitem{jun}K. Jung, Amenability, tubularity, and embeddings into $R^\omega$, {\em Math. Ann.} {\bf 338} (2007), 241--248.

\bibitem{ksp}C. K\"ostler, R. Speicher, A noncommutative de Finetti theorem: invariance under quantum permutations is equivalent to freeness with amalgamation, {\em Comm. Math. Phys.} {\bf 291} (2009), 473--490.

\bibitem{liu}W. Liu, General de Finetti type theorems in noncommutative probability, {\em  Comm. Math. Phys.} {\bf 369} (2019), 837--866.

\bibitem{mal}S. Malacarne, Woronowicz's Tannaka-Krein duality and free orthogonal quantum groups, {\em Math. Scand.} {\bf 122} (2018), 151--160.

\bibitem{mwe}A. Mang and M. Weber, Categories of two-colored pair partitions part II: categories indexed by semigroups, {\em J. Combin. Theory Ser. A} {\bf 180} (2021), 1--37.

\bibitem{mpa}V.A. Marchenko and L.A. Pastur, Distribution of eigenvalues in certain sets of random matrices, {\em Mat. Sb.} {\bf 72} (1967), 507--536.

\bibitem{meh}M.L. Mehta, Random matrices, Elsevier (1967).

\bibitem{mni}J.A. Mingo and A. Nica, Annular noncrossing permutations and partitions, and second-order asymptotics for random matrices, {\em Int. Math. Res. Not.} {\bf 28} (2004), 1413--1460.

\bibitem{mpo}J.A. Mingo and M. Popa, Freeness and the transposes of unitarily invariant random matrices, {\em J. Funct. Anal.} {\bf 271} (2016), 883--921.

\bibitem{msp}J.A. Mingo and R. Speicher, Free probability and random matrices, Springer (2017).

\bibitem{mvo}F.J. Murray and J. von Neumann, On rings of operators. IV,  {\em Ann. of Math.} {\bf 44} (1943), 716--808.

\bibitem{nsp}A. Nica and R. Speicher, Lectures on the combinatorics of free probability, Cambridge Univ. Press (2006).

\bibitem{pop}S. Popa, Classification of amenable subfactors of type II, {\em Acta Math.} {\bf 172} (1994), 163--255.

\bibitem{psh}S. Popa and D. Shlyakhtenko, Universal properties of $L(F_\infty)$ in subfactor theory, {\em Acta Math.} {\bf 191} (2004), 225--257.

\bibitem{rwe}S. Raum and M. Weber, The full classification of orthogonal easy quantum groups, {\em Comm. Math. Phys.} {\bf 341} (2016), 751--779.

\bibitem{sch}H. Schultz, Non-commutative polynomials of independent Gaussian random matrices. The real and symplectic cases, {\em Probab. Theory Related Fields} {\bf 131} (2005), 261--309. 

\bibitem{shl}D. Shlyakhtenko, Some applications of freeness with amalgamation, {\em J. Reine Angew. Math.} {\bf 500} (1998), 191--212.

\bibitem{sp1}R. Speicher, Multiplicative functions on the lattice of noncrossing partitions and free convolution, {\em Math. Ann.} {\bf 298} (1994), 611--628.

\bibitem{sp2}R. Speicher, Combinatorial theory of the free product with amalgamation and operator-valued free probability theory, {\em Mem. Amer. Math. Soc.} {\bf 132} (1998).

\bibitem{twa}P. Tarrago and J. Wahl, Free wreath product quantum groups and standard invariants of subfactors, {\em Adv. Math.} {\bf 331} (2018), 1--57.

\bibitem{twe}P. Tarrago and M. Weber, Unitary easy quantum groups: the free case and the group case, {\em Int. Math. Res. Not.} {\bf 18} (2017), 5710--5750.

\bibitem{twi}C.A. Tracy and H. Widom, Level-spacing distributions and the Airy kernel, {\em Comm. Math. Phys.} {\bf 159} (1994), 151--174.

\bibitem{vo1}D.V. Voiculescu, Addition of certain noncommuting random variables, {\em J. Funct. Anal.} {\bf 66} (1986), 323--346.

\bibitem{vo2}D.V. Voiculescu, Multiplication of certain noncommuting random variables, {\em J. Operator Theory} {\bf 18} (1987), 223--235.

\bibitem{vo3}D.V. Voiculescu, Circular and semicircular systems and free product factors, {\em Progr. Math.} {\bf 92} (1990), 45--60. 

\bibitem{vo4}D.V. Voiculescu, Limit laws for random matrices and free products, {\em Invent. Math.} {\bf 104} (1991), 201--220.

\bibitem{vo5}D.V. Voiculescu, The analogues of entropy and of Fisher's information measure in free probability theory, {\em Comm. Math. Phys.} {\bf 155} (1993), 71--92.

\bibitem{vdn}D.V. Voiculescu, K.J. Dykema and A. Nica, Free random variables, AMS (1992).

\bibitem{von}J. von Neumann, Mathematical foundations of quantum mechanics, Princeton Univ. Press (1955).

\bibitem{wan}S. Wang, Quantum symmetry groups of finite spaces, {\em Comm. Math. Phys.} {\bf 195} (1998), 195--211.

\bibitem{was}A. Wassermann, Coactions and Yang-Baxter equations for ergodic actions and subfactors, {\em London Math. Soc. Lect. Notes} {\bf 136} (1988), 203--236.

\bibitem{wei}D. Weingarten, Asymptotic behavior of group integrals in the limit of infinite rank, {\em J. Math. Phys.} {\bf 19} (1978), 999--1001.

\bibitem{wey}H. Weyl, The theory of groups and quantum mechanics, Princeton Univ. Press (1931).

\bibitem{wig}E. Wigner, Characteristic vectors of bordered matrices with infinite dimensions, {\em Ann. of Math.} {\bf 62} (1955), 548--564.

\bibitem{wit}E. Witten, Quantum field theory and the Jones polynomial, {\em Comm. Math. Phys.} {\bf 121} (1989), 351--399.

\bibitem{wo1}S.L. Woronowicz, Compact matrix pseudogroups, {\em Comm. Math. Phys.} {\bf 111} (1987), 613--665.

\bibitem{wo2}S.L. Woronowicz, Tannaka-Krein duality for compact matrix pseudogroups. Twisted SU(N) groups, {\em Invent. Math.} {\bf 93} (1988), 35--76.

\bibitem{zin}P. Zinn-Justin, Jucys-Murphy elements and Weingarten matrices, {\em Lett. Math. Phys.} {\bf 91} (2010), 119--127.

\end{thebibliography}
\end{document}